\documentclass[12pt]{article}
\usepackage[a4paper,margin=1in]{geometry}
\usepackage{amsmath,amssymb,amsthm,mathtools,bm}
\usepackage{enumerate}
\usepackage{enumitem}
\usepackage{mathrsfs}
\usepackage{tikz}
\usepackage{float}
\usepackage[hidelinks]{hyperref}
\usepackage{titlesec}
\usepackage{needspace}
\usepackage{booktabs}
\usepackage{makeidx}
\makeindex

\DeclareMathOperator{\Aut}{Aut}

\DeclareMathOperator{\PSL}{PSL}
\DeclareMathOperator{\SL}{SL}
\newcommand{\UHP}{\mathbb{H}}
\DeclareMathOperator{\dist}{dist}

\newcommand{\A}{\mathbb{A}}
\newcommand{\PP}{\mathbb{P}}

\newcommand{\OO}{\mathcal{O}}
\newcommand{\fm}{\mathfrak{m}}

\DeclareMathOperator{\Spec}{Spec}
\DeclareMathOperator{\Tor}{Tor}

\DeclareMathOperator{\ord}{ord}
\DeclareMathOperator{\Res}{Res}

\DeclareMathOperator{\length}{length}

\newcommand{\Div}{\mathrm{Div}}
\newcommand{\Prin}{\mathrm{Prin}}
\newcommand{\Cl}{\mathrm{Cl}}
\newcommand{\Pic}{\mathrm{Pic}}

\newtheorem{theorem}{Theorem}[section]
\newtheorem{proposition}[theorem]{Proposition}
\newtheorem{lemma}[theorem]{Lemma}
\theoremstyle{definition}
\newtheorem{definition}[theorem]{Definition}
\newtheorem{remark}[theorem]{Remark}
\newtheorem{corollary}[theorem]{Corollary}
\newtheorem{example}[theorem]{Example}
\newtheorem{exercise}[theorem]{Exercise}
\newcommand{\Hol}{\operatorname{Hol}}
\usetikzlibrary{calc,arrows.meta}

\newcommand{\C}{\mathbb{C}}
\newcommand{\R}{\mathbb{R}}
\newcommand{\D}{\mathbb{D}}

\newcommand{\Vol}{\operatorname{Vol}}

\newcommand{\CP}{\mathbb{CP}}
\newcommand{\Z}{\mathbb{Z}}
\newcommand{\Log}{\operatorname{Log}}
\newcommand{\Arg}{\operatorname{Arg}}
\newcommand{\Q}{\mathbb{Q}}

\newcommand{\Cstar}{\C^{\ast}}

\newcommand{\EssayByline}[1]{%
  \par\smallskip\noindent\textit{By #1}\par\medskip}

\newcommand{\EssayJunhyukBoo}{%
\clearpage
\subsection*{What Calculus Revealed to Me about Mathematics}
\EssayByline{Junhyuk Bu}

When one opens an advanced mathematics book, one often sees a dense sequence of
definitions, theorems, and proofs. I have always admired their logical completeness,
but logic alone did not make it easy for me to sustain the desire to study. When I
first encountered groups, I wondered why anyone would define such an object; when I
then saw the definition used in a proof, I wondered how anyone could have invented
the proof. For a long time, this distance made mathematics seem like an activity
reserved for people unlike me.

My experience changed through the \emph{Enjoying Math} community. Repeatedly working
through calculus and differential-form computations gave me the confidence that I,
too, could do mathematics. The decisive example was the winding form, which
measures how the angle around the origin changes along a path. Although it has
no global potential on the punctured plane, it has local
potentials. On suitable overlapping angular charts one may use local branches of
the argument. Subdivide a loop so that every subarc lies in one chart and
apply the fundamental theorem of line integrals on each subarc. The endpoint terms
cancel successively, except for the accumulated change of argument after one turn.
The integral is therefore a full turn multiplied by the loop's winding number.
For a positively oriented loop winding once around the origin, it is two pi.
I had known for years that integrating a derivative gives the difference between
its endpoint values, but seeing the same theorem recover
a topological integer was startling. The calculation explained both what fails
globally and what remains valid locally.

The same example suggested the shape of homotopy invariance. If a loop is deformed
without crossing the origin, the integral cannot change continuously from one
integer multiple of a full turn to another. Topological language appeared after the
calculation, as the language needed to describe a phenomenon already visible in
calculus.

Complex notation makes the same structure even more concentrated: the angular
form is the imaginary part of the logarithmic differential. Locally that
differential comes from a logarithm, but no single-valued logarithm exists
on the entire punctured complex plane. As the book proceeds, this elementary obstruction returns
in holomorphic line bundles, invertible sheaves, Cartier divisors, \v{C}ech
cohomology, and the first Chern class. What first looked like a collection of remote
definitions gradually became a series of increasingly abstract ways to record the
same failure of local data to glue globally.

This is why calculus is more than a preliminary tool for me. A computation can
reveal the criterion by which a theory is built. It can show which quantity survives,
which local choices are possible, and why a global obstruction must be introduced.
The result of a computation does not change with fashion; by returning to it, we can
retrace the route by which mathematicians discovered structure.
\clearpage
}

\newcommand{\EssayHyunminYoon}{%
\clearpage
\subsection*{Complex Analysis That Once Looked Like an Alien Language: Is the Square Root Really a Function?}
\EssayByline{Hyunmin Yoon}

I studied mathematics as an undergraduate, including complex analysis, but for a
long time I regarded the subject as multivariable calculus written in unfamiliar
symbols. A complex number is, after all, determined by two real coordinates. I also
learned about branch cuts, yet I did not understand what mathematical problem the
cutting and gluing were meant to solve.

The question of taking a square root changed that perspective. An equation and a function are not
the same thing. To define a function, one must assign exactly one value to every
point of the domain. In polar coordinates, one might take the square root of
the radius and halve the angle. But adding a full turn to the angle leaves
the original complex number unchanged while reversing the proposed square
root's sign. Thus this prescription does not define a single-valued function
on the punctured complex plane.

Locally, however, square roots can be built directly. Choose a disk whose
center is not the origin and whose radius is smaller than the center's
distance from the origin. The binomial series for a square root converges
throughout this disk after recentering and rescaling. Multiplying the series
by a chosen square root of the center gives a holomorphic function whose
square is the coordinate itself. On a connected overlap, two such
choices differ by a constant sign. Those signs are not arbitrary: analytic
continuation once around the origin changes a chosen square root to its negative.
The overlap signs therefore cannot all be removed by changing the local choices,
so the local roots do not glue to a global root on the punctured plane.

This calculation finally made a Riemann surface feel necessary rather than
decorative. Instead of forcing the two values onto the original punctured plane, we
place them on two sheets and glue the sheets so that analytic continuation becomes
single-valued. The resulting space is not an artificial picture attached to the
formula; it is the natural domain of the function.

I used to ask why complex numbers were needed at all. I now think they provide a
window through which calculus sees monodromy, local-to-global compatibility, and
geometry. The square-root example is elementary, but it contains the basic questions
that later reappear in covering spaces, line bundles, and sheaf cohomology: what can
be defined locally, what changes on overlaps, and what prevents the local objects
from becoming one global object?
\clearpage
}

\newcommand{\EssayDongsongBae}{%
\clearpage
\subsection*{The Road from Calculation to Concept}
\EssayByline{Dongsong Bae}

While revisiting calculus, I found myself seeing familiar formulas differently.
For a smooth function, integrating its differential along a piecewise smooth
curve gives the difference between the function's values at the endpoints.
I had once regarded this only as a formula for evaluating a line integral. Seen
from a little distance, however, it is striking: the left-hand side adds changes
along an entire curve, yet the right-hand side remembers only the endpoints. The
same phenomenon was already present in the one-variable fundamental theorem
of calculus.
From then on, I no longer regarded calculation merely as a way to obtain an
answer. Looking back at a completed calculation can reveal why a definition or a
theorem has the shape that it does. This has become one of the experiences I meet
most often in studying mathematics.

\paragraph{Constructing a function before invoking the theorem.}
The clearest instance arose while I was learning differential forms. Consider
a one-form in the plane with one coefficient in the horizontal direction
and another in the vertical direction. It is closed exactly when the
vertical derivative of the horizontal coefficient equals the horizontal
derivative of the vertical coefficient. Up to this point, I had only
unpacked a definition. The more important step was to resist quoting the
Poincar\'e lemma immediately and instead try to construct a function whose
differential is the given form.

Suppose, for this calculation, that the domain is a rectangle containing
the origin and the chosen point, so the coordinate paths stay inside it. I
integrate the horizontal coefficient along the first segment and the vertical
coefficient along the second. Differentiating the resulting function vertically
recovers the vertical coefficient. Differentiating horizontally first produces
an extra integral; the closedness condition changes its integrand into the
vertical derivative of the horizontal coefficient, so the fundamental theorem
of calculus reduces the result to the horizontal coefficient. The differential
of my constructed function is therefore the original form. I had known the
sentence ``on a star-shaped domain, every closed form is exact.'' After carrying
out this construction, I saw the statement differently. Closedness does not merely say that one
more differentiation produces zero; it supplies precisely the compatibility needed
to reconstruct the form as the differential of a function. Instead of memorizing a
theorem and then checking an example, I had built the function in an example and
recognized the Poincar\'e lemma as the result that guarantees this construction in
general.

\paragraph{Why take a quotient? A calculation on the circle.}
A similar change occurred when I studied de Rham cohomology. Its first group
on the circle identifies closed one-forms that differ by an exact one. At first, I did not feel why
this quotient was necessary.

Parametrize the circle by an angle and write a one-form as a periodic
coefficient times the angular differential. If the form is exact, its
coefficient is the derivative of a periodic function, and its integral over
one full turn vanishes. Conversely, a vanishing integral lets us integrate
the coefficient to obtain a periodic primitive. Thus exactness is precisely
the condition that the integral around the circle vanish.
Fourier series make the same fact especially transparent. The sine and cosine
modes integrate to zero over a full turn; only the constant mode remains. The
calculation shows exactly what survives after exact forms have been factored out.
The cohomological quotient no longer looked like a formal instruction to ``subtract
exact forms from closed forms.'' By watching the removable components disappear, I
could see what the quotient was designed to forget and what it was designed to
retain. The definition came first, but the calculation explained it afterward.

\paragraph{From determinants to permutations, and back through wedge products.}
I first met determinants as formulas for small square matrices. Later
I reconsidered them through multilinearity, alternation, and normalization. When
the columns are expanded by multilinearity, each term chooses one basis vector
from every column. A term in which a basis vector is chosen twice vanishes by
alternation. The surviving terms choose each basis vector exactly once,
and such a choice is a permutation. Restoring the basis vectors to their
original order contributes a sign. The determinant is therefore a sum over
permutations, each weighted by that sign and by the corresponding product
of matrix entries.
Following this expansion also leads naturally to cofactor expansion. The symmetric
group had not entered the determinant from nowhere: once multilinearity and
alternation were carried through to the end, permutations were unavoidable.

Much later, wedge products brought me back to the same structure. Evaluating
the wedge of two one-forms on two vectors gives the product of their values
in one order minus the product in the opposite order: exactly the pattern
of a two-by-two determinant.
The alternation I had met in determinants appeared in another setting: exchanging
the inputs reverses the sign, repeated directions give zero, and the value measures
oriented area or volume. Wedge product felt much less like a new operation to
memorize once I recognized the determinant inside it. The same recognition made it
natural that a top-degree form should measure volume.

This is one of the examples closest to what I mean by a road from calculation to
concept. A structure encountered through computation in one subject returned in
another, and only later did I see what the two definitions shared.

\paragraph{Meeting the same language again in physics.}
Recently this experience has continued in my study of physics. In Hamiltonian
mechanics, position and momentum coordinates combine into a symplectic
two-form that allows Hamilton's equations to be expressed through the
Hamiltonian vector field.
What impressed me more was changing from Cartesian to polar coordinates and doing
the calculation directly. The intermediate expressions became complicated, but
after choosing the corresponding radial and angular momenta correctly, the same
canonical structure reappeared: one contribution pairs the radius with its
momentum, and the other pairs the angle with its momentum.
At first, differential forms had seemed to me mainly a language for integration and
the general Stokes theorem. Here the same kind of two-form described the geometry
of the state space of a mechanical system. I am now studying Noether's theorem as
the point at which this geometric language connects symmetry with conservation
laws. That study is still in progress, but calculation has already shown me that
differential forms and mechanics, which I once learned separately, meet within the
same structure.

\paragraph{Calculation was not the opposite of definition.}
I used to think that the natural order of learning mathematics was to master the
definition, prove the theorem, and only then solve examples. Precise definitions and
proofs remain essential, but understanding does not always arise in that order.

Constructing a potential from a closed form showed me the Poincar\'e lemma. Computing
one-forms on the circle showed me what the cohomological quotient removes.
Expanding the determinant axioms showed me why permutations appear. Wedge products
and symplectic forms then let me recognize familiar calculations in unfamiliar
places.

I now try not to stop at memorizing a new definition. I calculate the simplest case
I can, and when the calculation is finished, I ask what repeated, what disappeared,
and what remained. Seen this way, abstraction is not the rejection of concrete
examples. It is the act of holding on to what remains common across many of them.

Calculation is not merely the final step in applying a concept that has already
been made. In my own study, looking closely enough at a calculation has often been
the way I discover where the concept begins.
\clearpage
}

\newcommand{\EssayInhoSong}{%
\clearpage
\subsection*{Complex Analysis and Differential Geometry as a Bridge to Mathematical Physics}
\EssayByline{Inho Song}

Before participating in the second Riemann complex-analysis program, I thought of
complex analysis mainly as the study of complex functions, contour integration, and
Riemann surfaces. The course made me see how harmonic functions, de Rham
cohomology, Hodge theory, Cartan's structure equations, characteristic classes,
and Riemann--Roch speak to one another. My own research on topological insulators,
Bogoliubov Fermi surfaces, and flux excitations in the toric code made the
connection to physics especially immediate.

I first understood the topological role of a differential equation through Fourier
series on the circle. A smooth one-form can be decomposed into Fourier modes.
Every nonconstant mode is the derivative of a periodic function; the constant
mode is not, because its integral around the circle need not vanish. The
angular differential therefore represents the single independent class in
the circle's first de Rham cohomology. This small calculation gave me a way to
recognize what Hodge theory says in general \cite{Warner1983,BottTu1982}.
On a compact oriented Riemannian
manifold without boundary, each de Rham class has a unique harmonic representative.
Here harmonic means annihilated by the Hodge Laplacian, formed from the
exterior derivative and its formal adjoint. This Laplacian is nonnegative. The
Poisson kernel from elementary harmonic analysis solves a different problem:
it extends boundary data harmonically into a disk. It does not, by itself,
choose harmonic representatives of arbitrary de Rham classes. Both constructions
use the Laplacian, but I find it useful to keep their purposes separate.

The sphere made this distinction concrete. Its spherical harmonics have
strictly positive Laplacian eigenvalues except in the constant mode, so its
harmonic functions are constant. On a two-sphere every closed one-form is
exact. If such a form is also harmonic, the energy identity forces both its
exterior derivative and its adjoint derivative to vanish. Its potential is
then a harmonic function and hence constant; the one-form itself vanishes.
The area form is harmonic and spans the harmonic two-forms. The sphere
therefore has one cohomology class in degree zero, none in degree one, and
one in degree two; their alternating count is two.
One can also expand nonharmonic one-forms using derivatives of positive-eigenvalue
spherical harmonics and their Hodge duals. I initially tried to deduce
the absence of harmonic one-forms from the failure of the usual polar
coordinate differentials to extend across the poles. That argument is insufficient:
a general one-form has variable coefficients. Exactness, followed by the
energy identity, supplies the missing step.

Cartan's structure equations gave me another calculation I could check directly.
They relate an orthonormal coframe to its connection and obtain curvature by
differentiating that connection and adding a quadratic correction. On the
unit sphere, away from the poles, a polar coframe gives a connection whose
curvature is exactly the area form. The quadratic correction vanishes here
because rotations in an oriented tangent plane form an abelian group, not
simply because the surface has complex dimension one. Integrating this
curvature over the sphere gives four pi; after the Gauss--Bonnet normalization,
the result is two, the sphere's Euler characteristic.
The same computation distinguishes curvature from holonomy. The standard flat
torus obtained by identifying opposite sides of a parallelogram has a globally parallel frame, so its
Levi-Civita holonomy is trivial. A different flat connection on a line bundle
over that torus can nevertheless have nontrivial holonomy around a noncontractible
loop. The coordinate differentials detect those loops: they are locally exact
but have nonzero periods and therefore are not globally exact. For a nonflat
metric on an oriented torus, the curvature may be nonzero at some points even
though its integral is zero, as Gauss--Bonnet requires. I learned not to infer
trivial holonomy merely from a zero total curvature integral.

Riemann--Roch brought these threads back to complex analysis. For a holomorphic
line bundle on a compact Riemann surface, the number of independent holomorphic
sections, corrected by the dimension of a first cohomology group, equals the
bundle's degree plus one minus the genus. The correction cannot generally be
ignored. With the appropriate convention, the degree is the integral of the
normalized Chern curvature. The higher-dimensional Hirzebruch--Riemann--Roch
formula explains why both the degree and the genus occur here
\cite{GriffithsHarris1994}: on a curve only degree-two characteristic-class
terms contribute, and the tangent bundle's Chern number is the Euler
characteristic. This is a useful interpretation of the curve formula,
not a replacement for its proof in this book.

In condensed-matter physics I meet a related construction in momentum space.
For a smooth isolated band, local normalized eigenstates
give a Berry connection. On overlaps the states differ by a phase,
so the local potentials need not combine into one global one-form, whereas
the curvature does. On a two-dimensional Brillouin torus, its integral,
normalized by a full turn, is an integer Chern number \cite{Haldane1988}.
The formula describes an
isolated band's line bundle; it does not apply unchanged to a group of
degenerate bands, which needs a nonabelian connection. This example helped me
see why transition functions, rather than a particular gauge choice, carry
the geometry.

I came to the course from questions about topological materials and left with
a much more concrete view of the mathematics behind them. Fourier modes made
cohomology calculable; the spherical coframe made curvature visible;
Riemann--Roch related a global section count to degree and genus. The value
of the connection between mathematics and physics, for me, is that each of
these calculations tells me exactly which local choices can disappear from
the final answer and which global information remains.
\clearpage
}

\newcommand{\EssaySuminKim}{%
\clearpage
\subsection*{What a Calculation on the Sphere Taught Me}
\EssayByline{Sumin Kim}

One calculation from the second Riemann complex-analysis program remains especially
vivid. We had been working in familiar Euclidean spaces, using Green's theorem and
the winding form, when the setting abruptly changed to the sphere. I could not
understand why a one-form on the sphere had to be discussed by splitting the sphere
into open sets and studying their overlap.

The area form cannot be exact on the whole sphere, because its integral
is four pi,
whereas Stokes' theorem would make the integral of an exact two-form vanish on a
closed surface. Yet after removing either pole, the form is exact on the
remaining contractible patch.

I wanted to find the potentials by hand. In the stereographic coordinates from the
north pole used in my calculation, the planar orientation is opposite to the chosen
orientation of the sphere. That reversal gives a minus sign in the planar
expression for area. I found a one-form there whose derivative recovered
the area form. Back in spherical coordinates, its angular coefficient is
minus one minus the cosine of the polar angle. On the patch obtained by
removing the south pole, the corresponding coefficient is one minus the
cosine. These corrections are exactly what cancel the apparent singularity
of the angular differential at the pole included in each patch.

On the overlap, the second potential exceeds the first by twice the angular
differential. This difference is closed but not exact on the annular overlap,
and its integral around a full turn is four pi. The same number appears both as the integral of area over the
whole sphere and as the integral of the mismatch around the overlap. That was the
moment when Mayer--Vietoris ceased to feel like a formal diagram. The sphere is split
because the desired potential exists locally but not globally; the overlap records
the precise obstruction to gluing the local answers.

I remember carrying out this calculation late at night, dissatisfied with spherical
coordinates and returning stubbornly to Cartesian coordinates. The emotional force
of the result mattered. Computation was not merely a way of checking a theorem I
already accepted. It allowed me to discover why the theorem needed its particular
language, and it gave me a concrete answer to the question, ``Why do we use two
patches?''
\clearpage
}

\newcommand{\EssayNamhoKim}{%
\clearpage
\subsection*{How Should Mathematics Be Learned?}
\EssayByline{Namho Kim}

The question ``How should mathematics be learned?'' often troubles me not when I am
completely lost, but when I seem to be following everything. I can read the
definitions and trace the arrows in a diagram, yet after closing the book I may be
unable to say what I actually know. This happened with de Rham cohomology and the
Mayer--Vietoris sequence. I searched repeatedly for an easier explanation, but the
number of explanations grew faster than my understanding.

An elementary formula finally gave me something to hold on to: Stokes'
theorem equates the integral of a derivative over a region with the
integral of its primitive around the boundary.
If the region has no boundary, every exact top-degree form integrates to
zero. On the unit sphere, however, the curvature form equals the area form,
whose integral is four pi. Therefore it is closed but not exact on the
whole sphere.

Trying to construct a potential directly, I took an angular differential
with a coefficient depending only on the polar angle. Differentiation showed that
the coefficient must be a constant minus the cosine of the polar angle.
Smoothness at the north pole fixes that constant one way; smoothness at the
south pole fixes it the opposite way. No single choice works globally, but
each choice gives a valid local potential on its own patch.
Splitting the sphere into northern and southern patches and applying Stokes' theorem
to both shows that the integral of curvature over the sphere equals the
integral of the potentials' difference around their common boundary.
That difference measures twice the angular change, giving four pi.

The calculation changed the diagram into a story. A global primitive fails to exist;
two local primitives exist; their difference on the overlap is closed; and the
integral of that difference carries the global information. Mayer--Vietoris names and
organizes precisely these steps.

I learned that understanding does not always come from finding a more polished
explanation. Sometimes it begins by asking a familiar formula a small question and
following the resulting computation until a structure becomes unavoidable. The
calculation on the sphere did not make every cohomological construction easy, but it
gave me a place from which I could reconstruct the meaning of the definitions for
myself.
\clearpage
}

\newcommand{\EssayYonghyunJi}{%
\clearpage
\subsection*{Counting the Dimension of Global Function Spaces on the Riemann Sphere}
\EssayByline{Yonghyun Ji}

Let a path join two points and let a holomorphic function be defined near it.
Integrating the function's derivative along the path gives the difference
of its values at the endpoints. For a closed path this is zero. Yet the
logarithmic differential has a nonzero integral around the unit circle:
one full turn gives two pi times the imaginary unit. The reason is that
the reciprocal function has local logarithmic primitives but no single-valued
logarithm on a neighborhood of the whole circle. The integral records the change of
the logarithm after one turn.

Now let a nonzero meromorphic function live on the Riemann sphere. Near any
point, it is a power of a local coordinate multiplied by a holomorphic
factor that does not vanish there. The exponent is positive for a zero
and negative for a pole. Integrating the logarithmic derivative around a
small loop reads off this exponent after normalization by a full complex turn.
Removing small disks around all zeros and poles and applying Stokes' theorem to the
remaining surface shows that the sum of all those exponents is zero.
Thus the principal divisor has degree zero: zeros and
poles balance globally.

This balance explains why poles are not merely defects. Holomorphic functions on the
compact sphere are constant, so they cannot distinguish points. Meromorphic functions
can. The standard coordinate has a simple zero at the origin and a simple
pole at infinity; it separates points once the value
at infinity is included.

Allow a pole of order at most a chosen nonnegative integer at infinity and
no other poles. The resulting functions are exactly the polynomials of
degree at most that integer, so the dimension of their space is one more
than the permitted pole order. Riemann--Roch generalizes this elementary
count. It says that the dimension of the meromorphic functions allowed
by a divisor, after subtracting a canonical correction, equals the
divisor's degree plus one minus the genus. On the Riemann sphere, the genus
is zero and the canonical divisor has degree minus two. For a nonnegative
pole allowance at infinity, the correction vanishes and the theorem
returns the polynomial count.

What I find most illuminating is the path from a local integral to a global dimension.
The logarithmic derivative reads the order at one point; Stokes' theorem forces a
global balance; divisors store the local integers; and Riemann--Roch converts their
degree into the dimension of a space of global functions. The abstract theorem is
therefore rooted in computations that can be performed explicitly on the sphere.
\clearpage
}

\newcommand{\EssayHojaeJung}{%
\clearpage
\subsection*{Why Riemann--Roch and Mayer--Vietoris Initially Felt Unfamiliar}
\EssayByline{Hojae Jung}

Riemann--Roch initially felt strange to me. For a holomorphic line bundle
on a compact Riemann surface, it says that the dimension of its holomorphic
sections minus the dimension of its first cohomology equals one minus the
genus plus the bundle's degree. The genus term resembles half the Euler characteristic, but I did not understand why
the degree of a line bundle belonged in the same formula.

Gauss--Bonnet provides a way to see the connection. A Hermitian metric on
a holomorphic line bundle produces Chern curvature; its normalized integral
is the bundle's degree. For the tangent bundle, Gauss--Bonnet identifies
the normalized total curvature with the Euler characteristic, twice one
minus the genus.
Thus curvature, degree, and Euler characteristic are different languages for closely
related global information. Choosing a special structure and calculating within it
can reveal a relation that is invisible at the level of definitions alone.

Mayer--Vietoris created a similar difficulty. A chain of cohomology groups connected
by arrows seemed technical rather than motivated. I therefore returned to the circle.
Cover it by two contractible open sets and choose an angular function on
each. On every component of their overlap the two angles differ by a
constant, so their derivatives agree. The local derivatives glue to
the global winding form, even though the
local potentials do not glue to a global function. Conversely, a locally constant
difference on the overlap can be converted, with a partition of unity, into a global
closed one-form. This is the connecting homomorphism in a concrete form.

The example showed me what the exact sequence is designed to measure. It tracks how
data defined on an overlap produces a global object one degree higher, and it records
whether local primitives can be made compatible. The winding form is both a source of
confusion and an ideal guide: its local potentials are elementary, its global
obstruction is visible, and its integral is nonzero.

These experiences changed my view of abstraction. Mathematics does not simply move
from a particular statement to a more general one. Often we deliberately choose a
special metric, a cover, or a model example, calculate carefully, and observe which
quantities are forced to agree. The abstract theorem is then a compressed record of
those observations. Spending a long time with one good example can therefore reveal
more than learning many definitions quickly.
\clearpage
}

\titleformat{\section}{\Large\bfseries}{\thesection.}{0.6em}{}
\titleformat{\subsection}{\large\bfseries}{\thesubsection.}{0.5em}{}
\titleformat{\subsubsection}{\normalsize\bfseries}{\thesubsubsection.}{0.5em}{}
\makeatletter
\renewcommand*\l@subsection{\@dottedtocline{2}{1.5em}{3.4em}}
\makeatother

\title{\Large\bf
Complex Analysis and Riemann Surfaces: \\ A Graduate Path to Algebraic Geometry}
\author{
	Gunhee Cho \\[0.8em]
	Department of Mathematics, Texas State University
}

\hypersetup{
	bookmarksdepth=2,
	pdftitle={Complex Analysis and Riemann Surfaces: A Graduate Path to Algebraic Geometry},
	pdfauthor={Gunhee Cho},
	pdfsubject={Graduate textbook on complex analysis, Riemann surfaces, and algebraic geometry},
	pdfkeywords={complex analysis, Riemann surfaces, differential forms, sheaf cohomology, Riemann--Roch, Jacobians}
}

\date{}

\begin{document}
	\maketitle
	\thispagestyle{empty}
	\begin{abstract}
		This book develops a calculation-led path from one-variable complex analysis
		to the geometry and cohomology of compact Riemann surfaces. Cauchy theory,
		power and Laurent series, residues, harmonic functions, conformal mapping,
		the Riemann mapping theorem, and invariant metrics supply the analytic
		foundation. Differential forms and Stokes' theorem then lead to curvature,
		Gauss--Bonnet, de~Rham cohomology, Hodge theory, and Poisson equations. The
		sphere and the complex torus recur as test cases on which signs, periods,
		curvature, harmonic representatives, and solvability conditions can be checked
		directly.

		The local-to-global transition is developed through complex manifolds,
		holomorphic line bundles, divisors, sheaves, \v{C}ech and derived-functor
		cohomology, exact sequences, and Dolbeault resolutions. Applications include
		Chern classes and curvature, Stein Riemann surfaces, the Mittag--Leffler
		problem, branched coverings, and canonical bundles. The Calabi--Yau condition
		is made concrete by proving the equivalence between triviality of the
		canonical sheaf and the existence of a nowhere-vanishing holomorphic top form;
		on a compact Riemann surface this singles out genus one. Poincar\'e and Serre
		duality prepare several complementary proofs of Riemann--Roch. The curve
		theorem is followed by adjunction, local and cohomological intersection
		pairings, surface Riemann--Roch, and Euler-characteristic calculations from
		the Euler sequence.

		The later chapters identify the Jacobian simultaneously as a period torus and
		as $\mathrm{Pic}^0$, interpret Abel's theorem and Jacobi inversion through
		linear systems, and pass from compact Riemann surfaces to affine and
		projective algebraic curves. Detailed divisor arithmetic on $\PP^1$, local
		degree and ramification calculations, and differentials on hyperelliptic curves
		serve as computational laboratories for the general theory. Hilbert's
		Nullstellensatz, very ample divisors,
		plane models, normalization, local intersection multiplicity, and B\'ezout's
		theorem make that passage explicit. Theta and Weierstrass functions, Mellin
		transforms and the zeta functional equation, Galois coverings, and modular
		curves provide further analytic, geometric, and arithmetic directions.
		Throughout, definitions and proofs are paired with representative calculations,
		motivating explanations, and exercises. The intended readers are advanced
		undergraduates and beginning graduate students seeking a self-contained bridge
		from computational complex analysis to algebraic curves and their cohomology.
	\end{abstract}
	\clearpage
	\tableofcontents
	\newpage

\section*{Introduction}
\addcontentsline{toc}{section}{Introduction}

This version is not the publisher's typeset edition.
Any later published edition should be cited using its own bibliographic record.

This book is intended for advanced undergraduates and beginning graduate students.
It begins with calculations that can be checked directly and uses them to motivate
the global structures developed later. It grew out of the \emph{Enjoying Math} YouTube channel
and a series of group reading courses led by Junhyuk Bu and other participants. In those
groups, participants met regularly (often online) to work through derivations,
estimates, normal forms, and local models before stepping back collectively to
distill the definitions and proofs used here.

Gunhee Cho is the sole author and is responsible for the book's mathematical
content, exposition, selection, and organization. The book also contains signed
reflections by members of the study community. Their writers are credited at
the point of inclusion; these contributions do not confer coauthorship.

The essays and reflections are by Dongsong Bae, Junhyuk Bu,
Byungjoo Jeon, Yonghyun Ji, Sumin Kim, Namho Kim, Minseung Kwak, Hojae Jung,
Inho Song, and Hyunmin Yoon. The placement of these contributions and the
surrounding mathematical exposition remain the author's responsibility.

The material began as lecture notes. The present version supplies the definitions,
examples, exercises, references, and transitions needed for independent reading.

\subsection*{Study groups, YouTube, and collaborative computation}

The book retains the intellectual history of an ongoing learning community. Many
examples and proof-checks were first developed
on the board or tablet during live sessions, refined in discussion, and only later
written up in a coherent form. Videos on the \emph{Enjoying Math} channel supplement
the text with demonstrations of computations, branch-cut gluings, and geometric
constructions that are expressed symbolically in the book.

The signed reflections appear after sections to which they are mathematically or
pedagogically related. Each is separated from the main exposition and attributed
to its writer.

\subsection*{Mathematical path}

The chapters follow a path already visible in Riemann's work:
begin with \emph{complex analysis} (Cauchy theory, harmonicity, residues),
use it to build and control \emph{Riemann surfaces} (branch points, branch cuts,
local models), then develop geometric structure via \emph{differential forms} and
\emph{Hodge theory}, and finally cross to \emph{algebraic geometry} (divisors, line
bundles, sheaves, cohomology), culminating in the \emph{Riemann--Roch theorem}
and the \emph{Jacobian}. At every stage, abstraction is anchored in explicit
computations: Taylor/Laurent expansions, Stokes/Cauchy--Green identities, residue
extraction, branch-cut gluings, curvature and area forms, Mayer--Vietoris patching,
exact-sequence calculations, local intersection lengths, and concrete
section-counting by Riemann--Roch. Canonical bundles are tested both by transition
functions and by global holomorphic volume forms, so that their triviality has a
geometric meaning before it is used in cohomology. The final algebraic and arithmetic chapters do
not abandon analysis: period integrals construct the Jacobian, local holomorphic
coordinates compute intersection multiplicities, and theta integrals produce the
functional equation of the zeta function.

\subsection*{Prerequisites and style}

A working knowledge of multivariable calculus and linear algebra is assumed.
Familiarity with basic topology and differential geometry (manifolds, forms) is helpful
but not required at the outset; we introduce what we need by computation and local models.
Proofs of the central results are included. When a major background theorem is used
without proof, the text identifies its role and points to the references. A displayed
calculation is meant to be reproducible: substitutions, orientations, convergence
conditions, and sign conventions are stated when they affect the result.

\subsection*{Organization of the book}

The six parts are arranged so that each major tool is first encountered in a
calculation, then formulated invariantly, and finally reused in a global or
algebraic argument. The section numbering below follows the present consolidated
structure.

\begin{itemize}
	\item \textbf{Part I: Foundations of Complex Analysis (Sections 1--4).}
	The first two sections develop holomorphicity, Cauchy theory, series, residues,
	analytic continuation, harmonic functions, conformal mapping, the Riemann mapping
	theorem, and invariant metrics. The next two introduce differential forms, Stokes'
	theorem, curvature, Gauss--Bonnet, de~Rham cohomology, Poisson equations, and Hodge
	theory. The sphere and torus are used repeatedly to verify the general formulas.

	\item \textbf{Part II: Sheaves and Cohomology on Complex Manifolds (Sections 5--9).}
	Line bundles, divisors, and the Picard group lead to sheaves, stalks, \v{C}ech
	cohomology, derived functors, and long exact sequences. Fine and Dolbeault
	resolutions connect abstract cohomology with differential forms. Canonical bundles
	are studied through transition functions, curvature, and holomorphic volume forms;
	in particular, the Calabi--Yau condition is related to the equivalence
	$K_X\simeq\mathcal O_X$. Stein curves, Mittag--Leffler problems, branched covers,
	and Riemann--Hurwitz provide the principal applications. Explicit calculations of
	divisors, local degrees, ramification, and hyperelliptic differentials keep these
	constructions accessible before they are used abstractly.

	\item \textbf{Part III: Duality (Section 10).}
	Poincar\'e duality, Hodge pairings, residues, and Serre duality are presented as
	compatible ways to turn cohomology classes into numerical pairings. Explicit checks
	on $\PP^1$ and a complex torus prepare the dimension calculations in Riemann--Roch.

	\item \textbf{Part IV: Riemann--Roch and Intersection Theory (Sections 11--13).}
	The curve theorem is developed through divisor calculations and a cohomological
	add-a-point proof. Its consequences include the degree of the canonical divisor,
	vanishing, and section-counting estimates. The surface extension introduces local
	and cohomological intersections, adjunction, surface Riemann--Roch, the Euler
	sequence, and Euler-characteristic computations on projective surfaces.

	\item \textbf{Part V: Algebraic Curves and Algebraic Geometry (Sections 14--15).}
	The Jacobian is constructed from periods and identified with $\Pic^0$; Abel's
	theorem, Jacobi inversion, and linear systems explain how divisors are encoded by
	the period torus. The following section introduces affine varieties and Hilbert's
	Nullstellensatz, embeds compact Riemann surfaces by very ample divisors, constructs
	plane models and their normalizations, and develops local intersection multiplicity
	through B\'ezout's theorem.

	\item \textbf{Part VI: Special Functions, Galois Theory, and Arithmetic
	(Sections 16--18).}
	Theta functions reconnect the Jacobian with elliptic and Weierstrass functions;
	the Mellin transform of theta yields the analytic continuation and functional
	equation of the zeta function. Holomorphic maps are then compared with function-field
	extensions and deck groups. Modular curves conclude the book by bringing analytic
	quotients, algebraic models, cohomology, and arithmetic Galois actions into a common
	framework.
\end{itemize}

\subsection*{How to use the book}

A first reading should distinguish three kinds of passage:
\begin{enumerate}
	\item In a computational passage, reproduce the calculation and check its domain,
	orientation, and limiting argument.
	\item In a proof, first identify the one step that carries the argument; then fill in
	the estimates or algebra surrounding it.
	\item For an abstract construction such as a sheaf, divisor, or Jacobian, keep one
	test case nearby: $\CP^1$, a complex torus, or a two-sheeted cover.
\end{enumerate}

There are (at least) three natural reading paths:
\begin{itemize}
	\item \textbf{Analysis-first path:}
	Sections 1--2 $\to$ Section 9 (branched covers) $\to$ Section 16
	(theta and Weierstrass functions).
	\item \textbf{Geometry-first path:}
	Sections 3--4 $\to$ Sections 7--10 (de~Rham cohomology, curvature, and
	duality) $\to$ Sections 11--13 (Riemann--Roch and intersections).
	\item \textbf{Algebraic-geometry-first path:}
	Sections 5--9 (bundles, divisors, sheaves, and exact sequences) $\to$
	Sections 11--12 (Riemann--Roch) $\to$ Sections 14--16 (Jacobians,
	algebraic curves, and theta functions).
\end{itemize}

The three paths share the same examples, so a reader can change emphasis without
starting again from the beginning.

\subsection*{References and further reading}

The complex-analysis chapters follow the classical development in
\cite{Ahlfors1979,Conway1978,Rudin1987,SteinShakarchi2003Complex}. For additional
worked problems and applications, compare
\cite{SaffSnider2003,CarrierKrookPearson2005,AblowitzFokas2003}; Needham's
geometric viewpoint~\cite{Needham1997} is especially useful alongside the
computational examples. The topological prerequisites can be reviewed in
\cite{GamelinGreene1999,Hatcher2002}. Differential forms, Stokes' theorem, and
de Rham theory may be compared with
\cite{Flanders1989,Spivak1965CalcManifolds,Warner1983,BottTu1982}; geometric and
elliptic-regularity background appears in
\cite{DoCarmo1992,Lee2018Riemannian,Morrey2008,Taylor2011PDE}. For Riemann
surfaces and Hodge theory, standard complements are
\cite{Forster1981,FarkasKra1992,Wells2008,Jost2006}.

The chapters on divisors, line bundles, sheaves, and algebraic curves use the
conventions developed in~\cite{Gunning1962,GunningRossi1965,GriffithsHarris1994,Hartshorne1977}.
For algebraic curves, Jacobians, and theta functions, see
\cite{Fulton2008Curves,Miranda1995,Fay1973,Mumford1983Tata1,Mumford1984Tata2}.
The modular and arithmetic chapters may be supplemented by
\cite{Apostol1990,Koblitz1993,Serre1973Course,DiamondShurman2005,DiamondIm1995,Milne1980Etale}.
For the function-field and Galois viewpoint, compare
\cite{Rotman1998Galois,FriedJarden2008,Stichtenoth2009}; the brief discussion of
Belyi maps and dessins is placed in the broader setting of
\cite{Schneps1994,LandoZvonkin2004}.
These references provide fuller proofs of major background theorems that are stated
without proof here and offer alternative conventions where signs or normalizations
vary across the literature.

\newcommand{\ReflectionMinseungKwak}{%
\subsection*{A Perspective on Rigidity and Local-to-Global Questions}
\noindent\textit{By Minseung Kwak}

The many equivalent formulations of holomorphicity can at first obscure their
strength. A basic comparison already shows the rigidity at stake: the smooth real
function that assigns each complex number its squared modulus is not holomorphic,
while the identity theorem says
that a holomorphic function vanishing locally must vanish on its entire connected
domain. Smooth functions are more flexible; partitions of unity allow local data to
be assembled globally.

That local--global comparison anticipates the later language of germs, stalks,
sheaves, and sheaf cohomology. The usefulness of undergraduate complex analysis is
therefore not confined to a collection of integral formulas: it supplies a concrete
model of the rigidity and gluing questions that return throughout geometry.

The same bridge appears later. The Riemann mapping theorem upgrades a topological
hypothesis to a biholomorphism with the disk, while normal families and harmonic
methods turn compactness and boundary data into global control. Runge approximation
and the Mittag--Leffler problem then point toward a cohomological question: when can
locally prescribed analytic data be realized globally?
}

\newcommand{\ReflectionJunhyukBoo}{%
\subsection*{Learning Mathematics through Computation}
\noindent\textit{By Junhyuk Bu}

When I study abstract mathematics written as a long list of definitions and theorems,
there are so many moments when I find myself asking:
``Why on earth are we defining \emph{this}?''  
Then, when I see an important proposition being proved using those very definitions,
I find myself wondering:
``Who are these people that they could even think of such things?''
From within that emotional distance, I think I gradually drifted away from mathematics.
I lived for a long time feeling that mathematics is a subject for geniuses,
and that I am simply not someone who can ``do math.''

But through \emph{Enjoying Math}, and especially through the calculus computations
we have been doing (and continue to do), I have had a series of precious experiences:
moments when I genuinely feel, ``Ah, even someone like me can actually do mathematics.''
I believe this will remain true going forward as well.
I have come to see calculation not only as a way to test a theory, but also as a
mathematical activity worth pursuing in its own right.

I still vividly remember the moment in our study when we locally computed potentials
of the winding form, applied the fundamental theorem for line integrals, and obtained
twice pi, the angle of one full turn in radians, for a positively oriented simple loop
winding once around the origin.
Up to that point, I had only thought of this as a situation where the fundamental theorem
of calculus (FTC) fails because there is a hole at the origin.
But through that first computation, I realized that locally we \emph{can} apply FTC.

Of course, the fundamental theorem of calculus---integrating a derivative
gives the change in the original function between the endpoints---was a
calculation I had done so often that it felt utterly familiar.
However, when I actually plugged the endpoints of the curve into the local expression
of the inverse tangent function and added everything up, seeing that this was \emph{really} the FTC that
I thought I knew, the FTC suddenly felt extremely unfamiliar.
For the first time while doing a computation, I experienced something like catharsis.

This is very different from simply remembering, in the context of Green's theorem,
that ``because the curl of the vector field is zero, the line integral along any boundary
curve is the same.''
Working through the calculation made the mechanism concrete: locally the form is a
gradient, its curl vanishes, and the fundamental theorem for line integrals forces
the local contributions to agree. That experience showed me how much mathematical
intuition can grow out of computation.

In complex analysis, the winding form becomes the differential obtained by
dividing the change in a complex variable by the variable itself. Integrating
it makes the number of turns of a curve around the origin visible.
Even though I myself have been looking at the winding form for almost a year now,
I still do not fully understand it.
This makes me feel that it is incredibly important to examine seemingly elementary objects
very thoroughly and meticulously.

The world of mathematics that I get to know through computation is truly fun,
and it feels as if I am having a conversation with mathematicians of the past.
The numerical values produced by computations do not change; through them we can observe
which patterns they discovered and how they formulated theories.
Looking at these values, I feel a kind of awe, in a sense different from what I wrote
at the beginning: ``Ah, so this is what it means for someone to do mathematics.''

Even so, the reason I still think I can continue doing mathematics is that I have come
to understand the criteria by which mathematical theories are built, and that the core
of those criteria lies in calculus computation.

}
\newcommand{\ReflectionByungjooJeon}{%
\subsection*{Why Calculus Computation Matters}
\noindent\textit{By Byungjoo Jeon}

Why does calculation in calculus matter? I cannot give a final answer, but I
can describe what it has done for my own understanding of mathematics.

For me, Riemann's complex analysis in \emph{Enjoying Math} provides a lens through
which ideas from abstract algebra, topology, and differential geometry can be traced
back to one-, two-, and three-dimensional calculus, often through
the fundamental theorem of calculus.

Practically speaking, all the computational tools we used were essentially just calculus
and linear algebra that we learned in our first year of university.

When I look at the logically complete proofs of countless theorems, I find that,
line by line, I can follow the argument, but I do not really gain intuition from it.
When I ask whether I could have discovered such a proof on my own, however, the
finished argument can feel remote: logically complete, but stripped of the
experiments that might have suggested it.

However, calculations in Euclidean space often suggest ideas
that later appear in such proofs; from carefully observing those computations,
I have experienced that it is possible to make new discoveries.

At first I knew nothing, but as my learning accumulated, I slowly began to realize that
many abstract objects become intelligible through concrete models and computations.
If I had to pick one mathematically concrete object, I would, without hesitation,
say differential forms.
Conceptually, but also in terms of their role in computations,
they have extremely wide applicability.
It feels as though differential forms have appeared in almost every area we have touched.
If I had not known differential forms, I think I would have lost my way,
trying to follow all the theories purely logically, without any computations.

One charm of computation is that as you deepen your understanding of mathematics,
even the calculations you performed in the past start to look different, and the same
object can be interpreted differently, giving you a new perspective.
The belief---or perhaps the working hope---that much of mathematics can be described
in the language of calculus and linear algebra is, I think, both a gateway
to understanding abstract objects and, at the same time,
a way of learning not to dismiss elementary objects as ``easy,''
but instead to see beauty and depth in them.

Short exact sequences, cohomology, homology, homotopy, groups, rings, fields, modules,
tensor products, projective spaces, and sheaves can look remote when one first sees
only their definitions. Concrete models from calculus and linear algebra do not
replace those definitions, but they give me a way to test them and to perform
representative computations.
That is what I now consider the main charm of complex analysis as I learned it
in \emph{Enjoying Math}.

I am still learning and still in the process of realizing new things, but each year,
I look forward to gaining a new viewpoint and a deeper mathematical literacy.

}
\newcommand{\ReflectionYonghyunJi}{%
\subsection*{From Structure Back to Computation}
\noindent\textit{By Yonghyun Ji}

For a long time, I practiced a selective kind of mathematics.
Rather than doing computations in calculus, I preferred to set up algebraic structures
and look for systems and patterns within them.
I felt much more at home in the world where structures built from sets and operations---
like groups, rings, and vector spaces---are defined and then studied.

So, concretely speaking, calculus computations were not really to my taste.
The grad, curl, and div operators that I had learned as a first-year student
seemed like unrelated formulas to memorize.
I solved a few problems for exams, but I never seriously tried to think about
why they had those particular forms, or how they are related to each other.
Computation was just a procedure for obtaining results and felt a bit removed
from the kind of ``structures'' I liked.

But after encountering differential forms, my thinking changed a lot.
I saw that grad, curl, and div do not exist separately but appear naturally unified
in a single language---and the thing that really allowed me to understand
this language was none other than the calculus computations that I had done
half-heartedly in the past.

When we carry this language of differential forms over to complex space,
the story extends in yet another direction.
Studying Riemann's complex analysis, I realized that, in a sense, I had been preparing
to cross over into the ``algebraic worldview'' via calculus computations.

One of the biggest shocks was that my view of polynomials completely changed.
Before studying Riemann's complex analysis, a polynomial was, to me, just an algebraic object
formed by multiplying and adding variables and coefficients.
Within that framework, I was mainly interested in structure---what properties
the polynomial ring has, what ideals look like, and so on.

But when I learned that the most basic data of a polynomial---its coefficients
and the order of its zeros---can be recovered from integrals, it was quite a shock.
By integrating along a suitable curve, you see the coefficients pop out, and the order
of poles or zeros becomes visible.
I realized that the information I thought of as ``algebraic'' was actually hidden
inside integrals.

At that moment, I began to feel that the dichotomy of ``algebra vs.\ analysis'' is,
in actual mathematics, a somewhat artificial division.
In complex analysis, computing integrals---for example, actually computing integrals
along specific paths---is not just a matter of polishing calculation techniques.
It plays the role of reading off algebraic information.
Without doing the integrals, you simply do not see certain things.
Through computation, structures that were invisible suddenly reappear in a new form.

So now I interpret the phrase ``the importance of calculus computation'' somewhat differently.
In the past, I thought of computation as just a short rite of passage before
entering the world of structures.
But through complex analysis and differential forms, I came to feel that computation is
what illuminates the very structures I want to see.

For me, calculus computations became a bridge connecting the algebraic mathematics
I liked with complex analysis.
Whereas I used to want to avoid computation and see only structure,
now I accept that we can approach structure precisely through computation.
That realization itself is something I value very highly.

}
\newcommand{\ReflectionHyunminYoon}{%
\subsection*{Calculus as a Way of Seeing}
\noindent\textit{By Hyunmin Yoon}

Before talking about the importance of calculus, I would like to organize
how I used to think about complex numbers.

\paragraph{1. Thoughts on complex numbers.}
I first encountered complex numbers when I was told that the equation
``a number squared, plus one, equals zero'' actually \emph{has} solutions.
At that time, I did not feel anything in particular; I just thought,
``Oh, this is part of the curriculum.''
I only thought of complex conjugates and computations with them as something that appears
in problems about the relation between roots and coefficients.

When I first saw Euler's formula, I was impressed for a moment:
the base of the natural logarithm, pi, and the imaginary unit appearing
together in a single equation!
What on earth is its true nature?
But then, when I saw it again in university, that feeling collapsed.
I thought, ``Oh, it's just polar coordinates.''
Even when complex numbers appeared together with the exponential function in analysis,
it did not move me.

So after that, complex numbers were really out of sight and out of mind.
Later, in an electrical engineering course on signal processing,
I encountered complex numbers again.
But hearing about them from the engineering side, I ended up thinking,
``Complex numbers are just a kind of seasoning you bring in when needed.''
This made me distance complex numbers even further from my mind.
And so, when I graduated from the math department, I felt as if I had graduated
from complex numbers as well.

\paragraph{2. The importance of calculus.}
Looking back on when I studied mathematics as an undergraduate,
I was obsessed with how quickly I could reach more advanced courses.
As a result, each course was merely a stepping stone to the next,
and in my last semester I took graduate-level functional analysis and Riemannian geometry.
At that time I was intoxicated with the feeling of
``I am finally studying such difficult mathematics!''
I had no concept of reflecting on or ``appreciating'' mathematics.

Even when I decided to pursue a career in AI, I thought that
``bringing mathematics into my career'' meant having been exposed
to as difficult mathematics as possible.

But as I began learning differential forms in \emph{Enjoying Math},
I realized I had been studying mathematics the wrong way.
Whenever I faced a difficult proposition, I used to chew over each word in the statement.
However, studying differential forms was nothing like that.
Once you study differential forms, you spend ages talking about what integration is.

Until then, integration to me was just
``Riemann integral: first--second year integration,''
``Lebesgue integral: third--fourth year integration.''
Even when solving exercise problems, I often skipped the actual computation
in integrals.

Differential forms grabbed me by the collar and asked:
``Do you really know what an integral is?''
I responded, ``Of course I know what an integral is---why are you asking?''
But the more I learned, the more I felt:
``I actually have no clue what integration really is.''
I realized I did not truly understand even the elementary integral of
the identity function over the unit interval.
I began to see that even the routine exam problems I had breezed through rested on
geometric ideas that differential forms make explicit.

Differential forms are also a good friend that encourages self-reflection.
As we study mathematics, we often begin to think that studying means
encountering many difficult concepts and hearing about them.
It seems I am not the only one.
But after studying like that, I found I could not extend my thinking by myself.

In those moments, differential forms keep poking me from the side:
``Friend, what can you actually do?''
At the same time, they remind me that the only thing I can truly do
is calculus.

Ironically, I think the importance of calculus is precisely that
``calculus is important in order to study calculus.''
There is an enormous amount to unpack even in a single integral problem,
and likewise in a single differentiation problem.
As I started to unpack them one by one, I was finally able to appreciate mathematics.
The mathematics I learned afterward was built by piling up such observations.

Difficult things are still difficult, of course, but the difficulty is not in
how alien or impressive something looks.
Rather, it lies in how densely calculus has been condensed into it.
I think of calculus as a lens through which we can appreciate mathematics.
What I truly have to do is computation, and those computations can only be carried out
through the lens (or tool) of calculus.

\paragraph{3. Then why should we look at complex numbers?}
I know differential forms only a little, and I know even less about complex numbers.
If I ask myself why I should care about complex numbers, the first thing that comes to mind is:
``I can't do this integral---complex numbers, please help.''
When I examined the winding form at the level of computations,
I came to understand that this calculation really measures rotation.
Complex notation makes the angular part of the calculation explicit.
Angles themselves are real numbers, so even though we bring in complex numbers,
it still feels as if we are ultimately revisiting the same calculus
we were looking at before.

So I do not feel that complex numbers ``flip the entire setup''---rather,
they serve as a window that allows us to carry calculus further.
To be honest, I still cannot confidently explain ``why we should look at complex numbers.''
But when I decide to study calculus and sit down with complex numbers at my side,
I find that I can see integrals from a broader perspective.

}
\newcommand{\ReflectionSuminKim}{%
\subsection*{What It Means to Know through Computation}
\noindent\textit{By Sumin Kim}

Perhaps this sounds overly highbrow, but I think we have to start with the question,
``Why do we study mathematics at all?''
The natural answer is that we study in order to understand.
People will have very different standards for what it means ``to know'' something,
but I personally think that you should at least be able to explain it to someone else
before you can say you know it.

However, as abstraction increases, explanations often start to self-reference:
you end up recycling your previous explanations, and they risk turning into something
like a religious text.
I do not think that is a good explanation.
A good explanation, in my view, is one that is easy for the listener to understand.

For that reason, I think computation is important.
No matter how difficult a concept is, I understand it better once I know what
it allows me to compute or compare.
If you show it through computation, at the very least you can persuade someone
through the numerical outcome.
(Of course, persuading someone why we suddenly perform a certain computation is also
a very important problem in itself.)

Also, even if, in the distant future, theory becomes ever more abstract,
the concrete results of computations will not change.
So I believe that computations themselves can become a kind of lantern that illuminates
hidden stories (perhaps future stories) behind the scenes.

It might be better if I share a concrete computation experience.
Among the computations I have valued while studying Riemann surfaces and complex
analysis, one example stands out.

At first, I could not understand at all why we needed to look at the Mayer--Vietoris
sequence.
When I asked myself, ``Why can't I accept this?'', my self-reflection led me
to the conclusion:
``The real issue is that I don't even understand why we should cover
the sphere by two open sets in the first place.''

From there, I came to realize that the fundamental problem was my inability to accept
that the integral over the sphere is not zero.
Since polar coordinates felt inconvenient, I tried to compute potentials
in Cartesian coordinates, and through this I came to accept that the
difference of the local potentials records twice the change in angle.
Then, after computing two different local potentials, I finally began to sympathize
with the intention behind the partition of unity---namely, to glue these together.

In this way, computation gives someone who does not yet know the theory a point of
entry. It also reveals why the theory was organized as it was.
In that sense, I think computation is important to the act of \emph{knowing} in mathematics.

}
\newcommand{\ReflectionHojaeJung}{%
\subsection*{Calculus, Geometry, and Pattern}
\noindent\textit{By Hojae Jung}

I like discovering patterns in numbers and seeing geometric information.
Perhaps because of that, as a middle and high school student I liked number theory,
where you find patterns in integers, and Euclidean geometry,
where you use proofs to derive information about lengths and angles.

Calculus, for me, was not like that.
More precisely, differentiating and integrating functions seemed useful for computing
areas under graphs or extrema, but I did not see much connection to Euclidean geometry
or number theory.

However, studying with people in \emph{Enjoying Math}, I came to see that geometry,
number theory, and calculus are all connected, and that treating them
as separate subjects can cause serious misunderstandings.

You can see this, for example, when computing the de Rham cohomology of a circle.
The first examples make a simple but striking point: different spaces can have
different de Rham cohomology.

Because a differential form is defined on a particular space, its global behavior
can depend on the topology of that space.
But de Rham cohomology is a story about calculus, and at the same time it contains
information about the space itself.
This is astonishing.

Moreover, this story is not limited to the circle.
It naturally extends to higher-dimensional spheres and tori.
The most characteristic visual difference between the ordinary sphere
and the two-dimensional torus is that every loop on the sphere is contractible,
whereas on the torus there are two distinct types of loops that are not contractible.
Even more surprising is that there are coordinate differential forms associated
to these two types of loops.

And de Rham cohomology is not something detached from geometric information
and our geometric intuition.
On the contrary, calculus, in a way that is deeply aligned with our geometric instincts,
speaks about those instincts in a clear and profound language.

For that reason, calculus is important.
When visual intuition is not available, calculus, in turn, provides us with
a calculable way of obtaining geometric intuition.

When complex numbers enter calculus, even more intriguing phenomena appear.
Complex numbers have many good properties.
Unlike real two-dimensional vectors, they come with a natural multiplication.
We introduced the imaginary unit to take a square root of negative one,
but we do not need a new number system to take a square root of that unit:
complex numbers already contain it.
Moreover, by Euler's formula, we get a natural polar form.

Thus complex numbers have an in-built computational convenience.
It is very natural, then, to consider complex-valued functions and to attempt
to differentiate and integrate them.

Complex functions have the very nice property that, once differentiable,
they admit power series expansions.
And in the process of integrating complex functions, we gain the perspective
that line integrals of vector fields can be viewed as integrals of complex functions.

Then, just as we study patterns in numbers, we begin to count zeros and poles with
multiplicity.  The argument principle turns that count into a winding number, connecting
the local order of a zero or pole with the global geometry of an image curve.

I believe we are all born with a bit of curiosity and ability to detect numerical patterns
and geometric structure.
Calculus gradually reveals those pieces of information in a way that is very closely
aligned with our instincts.

And I want to emphasize: calculus is a computational method that is very familiar
to us (even if not always easy).
As a tool that clarifies our thoughts in a computable way, calculus is something
we could study all our lives and still never exhaust.

}
\newcommand{\ReflectionNamhoKim}{%
\subsection*{Holding on to Calculus}
\noindent\textit{By Namho Kim}

I used to be someone who could not do complex analysis.
I disliked the word ``complex'' itself, and I disliked the word ``analysis''
even more.
Complex numbers felt like symbols that had nothing to do with my life,
and analysis seemed like a world reserved for math majors who walk around
with their shoulders full of tension.

That I came to find complex analysis interesting was, surprisingly, thanks to computation.
Computation was almost the only thing I could actually hold on to,
and one of the few ways I could verify something with my own hands
rather than relying on someone else's authority.

In fact, I had never heard terms like differential forms or Green's theorem.
I did not know even the names of Stokes's theorem or the divergence theorem.
It was the same for grad, curl, and div.
I could not properly integrate the tangent function.
So the only thing I could really do was perform the given calculations,
and with that one small tool called calculus, I trudged forward.

But strangely, through repeatedly doing such simple calculations, a door opened.
Though I knew almost no mathematics, by holding on to calculus alone I managed
to endure complex analysis, and along the way I directly saw moments when
algebraic structure seeped out of the computations.
Even though I did not know advanced tools like Mayer--Vietoris, by looking at the winding form
through the lens of calculus and checking with my own hands why certain values appear,
I began to see directly how calculus could lead to algebraic questions.

In hindsight, that was a kind of ``naive intuition.''
Because I did not know abstract terminology, I think I was able to directly touch
the structures hidden inside calculations.

Complex analysis showed me where it connects with algebra:
why boundary integrals carry topological information, and so on.
I learned this not through general theory but gradually through computation.

What remains in my hands even now is, in the end, calculus.
With that one small tool, I am able to grab hold of propositions that once felt
like rocks rolled in from Mars---strange, unfriendly claims with no visible origin---
and tear them apart to see what is inside.
That is my only consolation.

In the past, algebra felt like hell.
Symbols appeared and disappeared; proofs ended somewhere far away from where they began,
without any visible starting point; the more I read, the more only \emph{why} remained.
But when I revisited complex analysis holding only calculus in my hand,
those hellish sentences started to look a little more human.

When you approach something through computation, you can grope for its meaning
with your fingertips, even if you do not know the full background.

The winding form, for instance, was utterly unfamiliar at first, but as I followed it
with a calculus mindset, I experienced the quiet realization
of ``Ah, this is how it connects.''

I still do not feel at home in abstract algebra or fully understand its language.
But thanks to that one familiar technique called calculus, I can now dissect
propositions that once seemed to come from another world.
That experience, I think, is the reason I continue to hold on to complex analysis.

}
\newpage 

\clearpage
\part{Foundations of Complex Analysis}
\label{part:complex-foundations}

This part develops the analytic results used later for branched covers, Riemann
surfaces, Hodge theory, and duality. Each theorem is paired with either a proof
mechanism or a calculation that shows how its hypotheses enter.

\section{Core Complex Analysis}
\index{complex analysis}

\subsection{Holomorphic Functions, Cauchy--Riemann Equations, and Local Conformality}
\index{Cauchy--Riemann equations}\index{conformal maps!local}
\label{subsec:holomorphic-CR-English}

\subsubsection{Complex differentiability as a real-linear condition}
\begin{quote}
	\textbf{Goal.} We explain why complex differentiability is much more rigid than real differentiability.
	Concretely, we will prove:
	\begin{enumerate}
		\item A $C^1$ function $f=u+iv$ is holomorphic $\iff$ it satisfies the Cauchy--Riemann (CR) equations.
		\item At a point where the derivative is nonzero, a holomorphic map has a
		\emph{local rotation + scaling} geometry: the differential has the form
		$\begin{psmallmatrix} a & -b \\ b & a \end{psmallmatrix}$.
		\item Complex line integrals $\int_\gamma f(z)\,dz$ are just coupled real line integrals; the first
		unit circle computations already foreshadow residues.
	\end{enumerate}
\end{quote}

\subsubsection{Domains and Complex Differentiability}

\paragraph{The main distinction.}
In real analysis, differentiability at a point controls only the best linear approximation along real directions.
In complex analysis, the difference quotient must approach the same limit
through \emph{every} complex direction. This compatibility will manifest as
the CR equations.

\begin{definition}[Domain]
	An open set $U\subset\C$ is a \emph{domain} if it is connected.
\end{definition}

\begin{definition}[Complex differentiability and holomorphicity]
	Let $U\subset\C$ be a domain and $f:U\to\C$.
	We say that $f$ is \emph{complex differentiable} at $z_0\in U$ if the limit
	\[
	f'(z_0)=\lim_{z\to z_0}\frac{f(z)-f(z_0)}{z-z_0}
	\]
	exists.
	If $f$ is complex differentiable at every point of $U$, we say $f$ is \emph{holomorphic} on $U$,
	and write $f\in\Hol(U)$.
\end{definition}

\paragraph{A necessary two-direction test.}
To disprove complex differentiability, compare the difference quotient along two paths
(e.g.\ $z=z_0+h$ and $z=z_0+ih$). If the two limits disagree, $f'(z_0)$ cannot exist.

\begin{exercise}[Two-path test]
	\label{exr:selected-two-path}
	Let $f(z)=\overline z$ and fix $z_0\in\C$.
	Compute the difference quotient along $z=z_0+h$ and $z=z_0+ih$ ($h\in\R$, $h\to0$)
	and show the two limits are different. Conclude that $f$ is nowhere complex differentiable.
\end{exercise}

\subsubsection{Real Jacobian viewpoint: what should $Df$ look like?}

\paragraph{Real-linear interpretation.}
Write $z=x+iy$ and $f(z)=u(x,y)+iv(x,y)$.
Then $f$ is a map $F:\R^2\to\R^2$, $F(x,y)=(u(x,y),v(x,y))$.
If $f$ is complex differentiable at $z_0$, then the best real linear approximation $DF(z_0)$
must be equal to \emph{multiplication by a complex number}.

\[
DF(x_0,y_0)=
\begin{pmatrix}
	u_x & u_y\\
	v_x & v_y
\end{pmatrix}_{(x_0,y_0)}.
\]

\begin{lemma}[Real-linear maps that are complex multiplications]
	A real-linear map $L:\R^2\to\R^2$ is multiplication by a complex number $a+ib$
	(i.e.\ $L(h)= (a+ib)\,h$ under $\R^2\simeq\C$) if and only if its matrix has the form
	\[
	\begin{pmatrix}
		a & -b\\
		b & a
	\end{pmatrix}.
	\]
\end{lemma}

\begin{proof}
	If $L$ is multiplication by $a+ib$, then
	\[
	(a+ib)(h_1+ih_2)=(ah_1-bh_2)+i(bh_1+ah_2),
	\]
	so in coordinates $(h_1,h_2)$ the matrix is $\begin{psmallmatrix}a&-b\\ b&a\end{psmallmatrix}$.
	Conversely, a map with this matrix satisfies exactly the above identity, hence is multiplication by $a+ib$.
\end{proof}

\begin{exercise}[Recognizing complex multiplication]
	Show that a real matrix $\begin{psmallmatrix}\alpha&\beta\\ \gamma&\delta\end{psmallmatrix}$
	represents multiplication by some complex number if and only if $\alpha=\delta$ and $\gamma=-\beta$.
	Interpret this as: ``columns are obtained by a $90^\circ$ rotation and same scaling.''
\end{exercise}

\subsubsection{Wirtinger derivatives and the CR equations}

\paragraph{Wirtinger notation.}
The Wirtinger derivatives package the CR equations into the single condition $\partial f/\partial\overline z=0$,
which becomes the gateway to $\overline\partial$-methods, differential forms, and sheaf cohomology later.

\begin{definition}[Wirtinger derivatives]
	Write $z=x+iy$ and $\overline z=x-iy$. Define
	\[
	\frac{\partial}{\partial z}
	:=\frac12\Big(\frac{\partial}{\partial x}-i\frac{\partial}{\partial y}\Big),\qquad
	\frac{\partial}{\partial\overline z}
	:=\frac12\Big(\frac{\partial}{\partial x}+i\frac{\partial}{\partial y}\Big).
	\]
	For $f=u+iv$ set
	\[
	\frac{\partial f}{\partial z}=\frac12(f_x-i f_y),\qquad
	\frac{\partial f}{\partial\overline z}=\frac12(f_x+i f_y),
	\]
	where $f_x=u_x+i v_x$ and $f_y=u_y+i v_y$.
\end{definition}

\begin{lemma}[CR equations in Wirtinger form]
	\label{lem:CR-English-refactored}
	Let $f=u+iv\in C^1(U)$. Then
	\[
	\frac{\partial f}{\partial\overline z}=0
	\quad\Longleftrightarrow\quad
	u_x=v_y\ \text{ and }\ u_y=-v_x.
	\]
\end{lemma}

\begin{proof}
	Compute:
	\[
	\frac{\partial f}{\partial\overline z}
	=\frac12(f_x+i f_y)
	=\frac12\big((u_x+i v_x)+i(u_y+i v_y)\big)
	=\frac12\big((u_x-v_y)+i(v_x+u_y)\big).
	\]
	This equals $0$ iff both real and imaginary parts vanish:
	$u_x=v_y$ and $v_x=-u_y$.
\end{proof}

\paragraph{Two basic calculations.}
\begin{example}[A holomorphic and a non-holomorphic map]
	\leavevmode
	\begin{enumerate}
		\item $f(z)=z^2$. Write $z=x+iy$:
		\[
		f=(x+iy)^2=(x^2-y^2)+i(2xy),
		\quad u=x^2-y^2,\ v=2xy.
		\]
		Then $u_x=2x,\ u_y=-2y,\ v_x=2y,\ v_y=2x$, so CR holds everywhere.
		Hence $f$ is entire and $f'(z)=2z$.
		
		\item $f(z)=\overline z=x-iy$. Here $u=x,\ v=-y$ so
		$u_x=1$ while $v_y=-1$, violating CR everywhere.
		Thus $f$ is nowhere complex differentiable.
		Geometrically it is a reflection (orientation reversing), which cannot be holomorphic.
	\end{enumerate}
\end{example}

\begin{exercise}[Practice with Wirtinger]
	\leavevmode
	\begin{enumerate}
		\item Compute $\partial z/\partial z$, $\partial z/\partial\overline z$, $\partial\overline z/\partial z$,
		$\partial\overline z/\partial\overline z$.
		\item Show that for $f(z)=z^n$ ($n\in\Z_{\ge0}$), we have $\partial f/\partial\overline z=0$.
		\item Show that for $f(z)=|z|^2=z\overline z$, one has
		$\partial f/\partial\overline z=z$.  This derivative vanishes at $0$ but
		not on any neighborhood of $0$.  Conclude that $|z|^2$ is complex
		differentiable at $0$ but is not holomorphic on any domain.
	\end{enumerate}
\end{exercise}

\subsubsection{From the CR equations to complex differentiability}

\paragraph{Proof mechanism.}
Use real differentiability: a $C^1$ map admits a first-order approximation by its Jacobian.
Then CR forces the Jacobian to be a complex multiplication matrix, hence the complex difference quotient converges.

\begin{theorem}[CR equations imply complex differentiability]
	\label{thm:CR-implies-holomorphic-English-refactored}
	Let $U\subset\C$ be open and let $f=u+iv\in C^1(U)$.
	If the CR equations hold at $z_0=x_0+iy_0$, then $f$ is complex differentiable at $z_0$ and
	\[
	f'(z_0)=u_x(x_0,y_0)+i v_x(x_0,y_0)=v_y(x_0,y_0)-i u_y(x_0,y_0).
	\]
\end{theorem}

\begin{proof}
	Let $F:\R^2\to\R^2$, $F(x,y)=(u(x,y),v(x,y))$.
	Since $f\in C^1$, $F$ is Fr\'echet differentiable at $(x_0,y_0)$, so for $h=(h_1,h_2)\to 0$,
	\[
	F(z_0+h)-F(z_0)=DF(z_0)\,h+r(h),\qquad \frac{\|r(h)\|}{\|h\|}\to0.
	\]
	Here
	\[
	DF(z_0)=
	\begin{pmatrix}
		u_x & u_y\\
		v_x & v_y
	\end{pmatrix}_{(x_0,y_0)}.
	\]
	By CR, $u_x=v_y$ and $u_y=-v_x$, hence
	\[
	DF(z_0)=
	\begin{pmatrix}
		u_x & -v_x\\
		v_x & u_x
	\end{pmatrix}_{(x_0,y_0)}
	=
	\begin{pmatrix}
		a & -b\\
		b & a
	\end{pmatrix},
	\quad a:=u_x(x_0,y_0),\ b:=v_x(x_0,y_0).
	\]
	By the previous lemma, this matrix is multiplication by $a+ib$ under $\R^2\simeq\C$.
	Identifying $h=(h_1,h_2)$ with $h_1+i h_2\in\C$, we obtain
	\[
	f(z_0+h)-f(z_0)=(a+ib)\,h+\varepsilon(h),
	\qquad\text{where }\frac{|\varepsilon(h)|}{|h|}\to0.
	\]
	Divide by $h\neq 0$:
	\[
	\frac{f(z_0+h)-f(z_0)}{h}=(a+ib)+\frac{\varepsilon(h)}{h}\ \longrightarrow\ a+ib.
	\]
	Therefore $f'(z_0)=a+ib$.
\end{proof}

\begin{exercise}[Direct derivative computation]
	Let $f(z)=z^2$ and $z_0\neq 0$.
	Compute $f'(z_0)$ by the difference quotient definition.
	Then compute $u_x(x_0,y_0)+i v_x(x_0,y_0)$ and verify they match.
\end{exercise}

\subsubsection{Equivalent characterizations of holomorphicity}

\begin{theorem}[The $C^1$ criterion for holomorphicity]
	\label{thm:holomorphic-equivalences-English-refactored}
	Let $U\subset\C$ be a domain and let $f=u+iv\in C^1(U)$.
	The following are equivalent:
	\begin{enumerate}
		\item $f$ is holomorphic on $U$.
		\item $\dfrac{\partial f}{\partial\overline z}=0$ on $U$.
		\item The CR equations hold on $U$:
		\[
		u_x=v_y,\qquad u_y=-v_x.
		\]
	\end{enumerate}
\end{theorem}

\begin{proof}
	$(2)\Leftrightarrow(3)$ is Lemma~\ref{lem:CR-English-refactored}.
	That (3) implies (1) follows from Theorem~\ref{thm:CR-implies-holomorphic-English-refactored}
	applied pointwise. For the converse, fix $z_0=x_0+iy_0$. Real increments
	show that $f'(z_0)=u_x(z_0)+iv_x(z_0)$, while imaginary increments give
	$f'(z_0)=(u_y(z_0)+iv_y(z_0))/i=v_y(z_0)-iu_y(z_0)$.
	Equating real and imaginary parts yields $u_x=v_y$ and $v_x=-u_y$ at
	$z_0$. Since $z_0$ was arbitrary, (3) holds throughout $U$.
\end{proof}

\begin{exercise}[Holomorphic functions satisfy the CR equations]
	Assume $f$ is complex differentiable at $z_0=x_0+i y_0$.
	Show that the limits along $z=z_0+h$ and $z=z_0+ih$ force
	$u_x(z_0)=v_y(z_0)$ and $u_y(z_0)=-v_x(z_0)$.
\end{exercise}

\subsubsection{Local conformality: infinitesimal rotation and scaling}

\paragraph{Infinitesimal geometry.}
At a point where $f$ is holomorphic and $f'(z_0)\neq 0$, the Jacobian is complex multiplication:
\[
DF(z_0)=
\begin{pmatrix}
	a & -b\\
	b & a
\end{pmatrix},
\qquad a+ib=f'(z_0).
\]
Thus $DF(z_0)$ maps every infinitesimal vector by the same scaling factor $|f'(z_0)|$ and the same
rotation angle $\arg f'(z_0)$. The nonzero derivative makes $f$ locally
invertible, and its differential preserves angles at every point in a
small neighborhood of $z_0$; this is the sense in which $f$ is locally conformal.

\begin{exercise}[Angle preservation from matrices]
	Let $L:\R^2\to\R^2$ have matrix $\begin{psmallmatrix}a&-b\\ b&a\end{psmallmatrix}$ with $a^2+b^2\neq0$.
	Show that for any nonzero vectors $v,w$, the angle between $Lv$ and $Lw$ equals the angle between $v$ and $w$.
	(Hint: write $L=\lambda R$ where $\lambda=\sqrt{a^2+b^2}$ and $R$ is orthogonal.)
\end{exercise}

\subsubsection{Complex line integrals and first experiments}

\paragraph{Line integrals.}
Once holomorphic functions are viewed as geometric maps, integration along curves becomes natural.
This is also the entry point to Cauchy's theorem and residues later.

\begin{definition}[Piecewise $C^1$ curve and complex line integral]
	A continuous curve $\gamma:[a,b]\to U$ is \emph{piecewise $C^1$} if there is a finite partition of $[a,b]$
	such that $\gamma$ is $C^1$ on each subinterval.
	For such $\gamma$ and a continuous function $f:U\to\C$, define
	\[
	\int_\gamma f(z)\,dz:=\int_a^b f(\gamma(t))\,\gamma'(t)\,dt.
	\]
\end{definition}

\begin{remark}[Real and imaginary decomposition (as two coupled real integrals)]
	Write $\gamma(t)=x(t)+i y(t)$ and $f=u+iv$.
	Then $\gamma'(t)=x'(t)+i y'(t)$ and
	\[
	f(\gamma(t))\gamma'(t)
	=(u+iv)(x'+i y')
	=(u x'-v y')+i(v x'+u y').
	\]
	Hence
	\[
	\int_\gamma f(z)\,dz
	=\int_a^b (u x'-v y')\,dt
	+i\int_a^b (v x'+u y')\,dt.
	\]
	This shows complex line integrals are not mysterious: they are precisely two real line integrals.
\end{remark}

\begin{example}[Unit circle integral: the first residue glimpse]
	Let $\gamma(t)=e^{it}$ for $0\le t\le 2\pi$ and $f(z)=z^n$ with $n\in\Z$.
	Then $dz=i e^{it}\,dt$, so
	\[
	\int_\gamma z^n\,dz
	=\int_0^{2\pi} (e^{it})^n\, i e^{it}\,dt
	=i\int_0^{2\pi} e^{i(n+1)t}\,dt
	=
	\begin{cases}
		0,& n\neq -1,\\[0.3em]
		2\pi i,& n=-1.
	\end{cases}
	\]
\end{example}

\begin{exercise}[Path dependence vs path independence]
	\leavevmode
	\begin{enumerate}
		\item Compute $\displaystyle \int_\gamma z\,dz$ for $\gamma(t)=t(1+i)$, $0\le t\le1$.
		\item Compute the same integral along the broken path $0\to 1\to 1+i$.
		\item Explain why both answers agree (hint: find an antiderivative of $z$).
	\end{enumerate}
\end{exercise}

\begin{exercise}[More circle integrals]
	Let $\gamma(t)=Re^{it}$, $0\le t\le 2\pi$ with $R>0$.
	Compute $\int_\gamma \frac{1}{z}\,dz$ and $\int_\gamma z^n\,dz$ for $n\in\Z$.
\end{exercise}

\medskip

\subsection{Cauchy's Integral Theorem and Integral Formula}
\index{Cauchy's integral theorem}\index{Cauchy's integral formula}
\label{subsec:CIT-CIF}

\subsubsection{Cauchy's theorem and formula: statements and hypotheses}
\begin{quote}
	\textbf{Goal.} We prove the two central results:
	\begin{enumerate}
		\item \textbf{Cauchy Integral Theorem (CIT):} if $f$ is holomorphic
		near a simple closed curve $\gamma$ and its interior,
		then $\int_\gamma f(z)\,dz=0$.
		\item \textbf{Cauchy Integral Formula (CIF):} for $z_0$ inside $\gamma$,
		\[
		f(z_0)=\frac{1}{2\pi i}\int_\gamma \frac{f(z)}{z-z_0}\,dz.
		\]
	\end{enumerate}
	The second formula shows that the values of a holomorphic function on a boundary
	determine its values in the interior.
\end{quote}

\subsubsection{Geometry/topology input: simple closed curves and orientation}

\paragraph{Boundary convention.}
We want a curve $\gamma$ that bounds a region $D$ so that we can apply Green's theorem (a planar form of Stokes).
Positivity means ``counterclockwise around the interior.''

\begin{definition}[Simple closed curve]
	A piecewise $C^1$ curve $\gamma:[a,b]\to\C$ is \emph{simple closed} if $\gamma(a)=\gamma(b)$ and
	$\gamma$ is injective on $[a,b)$. Thus the endpoints are the only repeated point.
	It is \emph{positively oriented} if it winds counterclockwise around its interior.
\end{definition}

\begin{exercise}[Orientation check]
	Let $\gamma(t)=e^{it}$, $0\le t\le 2\pi$. Show it is positively oriented.
	Now consider $\widetilde\gamma(t)=e^{-it}$. Show $\int_{\widetilde\gamma} f(z)\,dz=-\int_\gamma f(z)\,dz$.
\end{exercise}

\subsubsection{Cauchy Integral Theorem on planar domains (CIT)}

The proof below uses Green's theorem and therefore keeps the $C^1$ hypothesis
visible.  Goursat's theorem, discussed after the proof, removes this auxiliary
regularity hypothesis and yields the usual theorem for holomorphic functions.

\begin{theorem}[Green--Cauchy theorem]
	\label{thm:CIT-refactored}
	Let $\gamma$ be a simple, positively oriented, piecewise $C^1$ closed curve in $\C$,
	and let $D$ be the bounded domain enclosed by $\gamma$.
	If $f=u+iv$ is $C^1$ on an open neighborhood of $\overline D$ and satisfies
	the Cauchy--Riemann equations there, then
	\[
	\int_{\gamma} f(z)\,dz=0.
	\]
\end{theorem}

\paragraph{Proof mechanism.}
Rewrite the complex integral as two real line integrals.
Apply Green's theorem to convert them into area integrals.
Holomorphicity (CR equations) kills the area integrands.

\begin{proof}
	Write $f=u+iv$ and $dz=dx+i\,dy$. Then
	\[
	\int_\gamma f(z)\,dz=\int_\gamma (u+iv)(dx+i\,dy)
	=\int_\gamma (u\,dx-v\,dy)+i\int_\gamma (v\,dx+u\,dy).
	\]
	Green's theorem states: for $P,Q\in C^1$,
	\[
	\int_{\partial D} P\,dx+Q\,dy
	=\iint_D\Big(\frac{\partial Q}{\partial x}-\frac{\partial P}{\partial y}\Big)\,dA.
	\]
	Apply it to $P=u$, $Q=-v$ to get
	\[
	\int_\gamma (u\,dx-v\,dy)=\iint_D\big((-v)_x-u_y\big)\,dA=\iint_D(-v_x-u_y)\,dA,
	\]
	and to $P=v$, $Q=u$ to get
	\[
	\int_\gamma (v\,dx+u\,dy)=\iint_D(u_x-v_y)\,dA.
	\]
	Since $f$ is holomorphic, the Cauchy--Riemann equations hold on $D$:
	\[
	u_x=v_y,\qquad u_y=-v_x.
	\]
	Hence $-v_x-u_y=0$ and $u_x-v_y=0$ pointwise, so both area integrals vanish.
	Therefore $\int_\gamma f(z)\,dz=0$.
\end{proof}

\begin{remark}[Stokes' theorem and $\overline\partial f=0$]
	CIT is a 2D Stokes theorem in disguise: the contour integral is the integral of a 1-form over $\partial D$,
	which becomes an integral of $d(\text{that 1-form})$ over $D$. Holomorphicity forces that exterior derivative
	to vanish. Later we will rewrite this cleanly as a $\overline\partial$-statement on forms.
\end{remark}

\paragraph{Comparison with the fundamental theorem of line integrals.}
If $f$ has a global primitive on a neighborhood of $\gamma$ (that is,
$f=g'$), then the fundamental theorem of calculus along curves gives
$\int_\gamma f(z)\,dz=0$. The Cauchy theorem obtains this vanishing
directly from holomorphicity on the enclosed region, without requiring
us to find a primitive first. Conversely, on a simply connected domain,
the vanishing of closed-curve integrals can be used to construct one.

\begin{exercise}[Warm-up: a primitive case]
	Let $f(z)=z^3$. Compute $\int_{|z|=2} z^3\,dz$ directly using the parameterization $z=2e^{it}$.
	Explain why the answer must be $0$ using primitives.
\end{exercise}

\subsubsection{The basic kernel $\boldsymbol{1/z}$}

\paragraph{The first singular kernel.}
CIT says holomorphic integrands give $0$. The first nontrivial phenomenon is when the integrand
fails to be holomorphic at one point inside the curve. The kernel $1/(z-z_0)$ is the simplest such singularity.

\begin{definition}[Disks and circles]
	For $z_0\in\C$ and $r>0$, define
	\[
	B(z_0,r)=\{z\in\C:\ |z-z_0|<r\},\qquad
	\partial B(z_0,r)=\{z\in\C:\ |z-z_0|=r\}.
	\]
\end{definition}

\begin{lemma}[The fundamental circle integral]
	\label{lem:1overz-refactored}
	For every $r>0$, with the circle oriented counterclockwise,
	\[
	\int_{\partial B(0,r)}\frac{1}{z}\,dz=2\pi i.
	\]
\end{lemma}

\begin{proof}
	Parameterize $z=re^{it}$, $0\le t\le 2\pi$. Then $dz=ire^{it}\,dt$ and
	\[
	\int_{\partial B(0,r)}\frac{1}{z}\,dz
	=\int_0^{2\pi}\frac{1}{re^{it}}(ire^{it})\,dt
	=i\int_0^{2\pi}dt
	=2\pi i.
	\]
\end{proof}

\begin{remark}[Topology in one line: winding number preview]
	The value $2\pi i$ is not an accident: it detects that the curve winds once around the origin.
	Later we will express $\int_\gamma \frac{1}{z}\,dz=2\pi i\cdot\mathrm{wind}(\gamma,0)$.
\end{remark}

\begin{exercise}[Translate the kernel]
	Let $z_0\in\C$ and $\gamma=\partial B(z_0,r)$ positively oriented.
	Show $\displaystyle \int_{\gamma}\frac{1}{z-z_0}\,dz=2\pi i$ by substituting $w=z-z_0$.
\end{exercise}

\subsubsection{Cauchy integral formula}

The formula reconstructs $f(z_0)$ from the boundary values of $f$.

\begin{theorem}[Cauchy Integral Formula (CIF)]
	\label{thm:CIF-refactored}
	Let $U\subset\C$ be open and let $f\in C^1(U)$ satisfy the Cauchy--Riemann equations.
	Let $\gamma$ be a simple, positively oriented, piecewise $C^1$ closed curve
	whose image and interior lie in $U$, and let $z_0$ lie inside $\gamma$.
	Then
	\[
	f(z_0)=\frac{1}{2\pi i}\int_{\gamma}\frac{f(z)}{z-z_0}\,dz.
	\]
\end{theorem}

\paragraph{Proof mechanism.}
Remove a small disk about $z_0$. Cauchy's theorem applies on the remaining
annular region, and the integral on the small circle tends to $2\pi i f(z_0)$.

\begin{proof}
	Choose $\varepsilon>0$ so that the closed disk
	$\overline{B(z_0,\varepsilon)}$ lies inside $\gamma$. Apply Green's
	theorem directly to $f(z)/(z-z_0)$ on the multiply connected region
	between $\gamma$ and the small circle
	$C_\varepsilon=\partial B(z_0,\varepsilon)$. One way to justify this
	version of Green's theorem is to cut the region along a piecewise smooth
	arc joining its two boundary components, apply the ordinary theorem to
	the resulting Jordan region, and cancel the integrals along the two
	oppositely oriented copies of the cut. The
	inner boundary has the opposite orientation, so
	\[
	\int_\gamma\frac{f(z)}{z-z_0}\,dz
	=\int_{C_\varepsilon}\frac{f(z)}{z-z_0}\,dz.
	\]
	Parametrize the small circle by $z=z_0+\varepsilon e^{it}$:
	\[
	\int_{C_\varepsilon}\frac{f(z)}{z-z_0}\,dz
	=i\int_0^{2\pi}f(z_0+\varepsilon e^{it})\,dt.
	\]
	Continuity of $f$ at $z_0$ makes the last integrand converge uniformly to
	$f(z_0)$. Letting $\varepsilon\to0$ gives $2\pi i f(z_0)$ and proves the formula.
\end{proof}

\begin{remark}[Boundary determination]
	CIF reconstructs the interior value $f(z_0)$ by a contour integral of
	boundary values. It is a boundary-determination result for holomorphic
	functions, not a general principle that arbitrary sheaf sections are
	determined by their boundary values.
\end{remark}

\paragraph{An integral evaluated by Cauchy's formula.}
\begin{example}[Compute a contour integral using CIF]
	Let $\gamma=\partial B(0,2)$ oriented counterclockwise. Compute
	\[
	\int_\gamma \frac{\cos z}{z}\,dz.
	\]
	\emph{Solution.} Here $f(z)=\cos z$ is entire. CIF with $z_0=0$ gives
	\[
	\cos(0)=\frac{1}{2\pi i}\int_\gamma \frac{\cos z}{z-0}\,dz
	\quad\mathrel{\text{and hence}}\quad
	\int_\gamma \frac{\cos z}{z}\,dz=2\pi i\cos(0)=2\pi i.
	\]
\end{example}

\begin{exercise}[CIF as an integral-evaluation machine]
	Let $\gamma=\partial B(1,3)$ oriented counterclockwise.
	Compute $\displaystyle \int_\gamma \frac{e^z}{z-1}\,dz$ and $\displaystyle \int_\gamma \frac{e^z}{(z-1)^2}\,dz$.
	(Use Theorem~\ref{thm:CIF-refactored} and the derivative formula in the next subsection.)
\end{exercise}

\subsubsection{Consequences: smoothness, derivative formulas, and power series}

\paragraph{Regularity gained from the integral formula.}
In real analysis, differentiability does not automatically imply analyticity.
In complex analysis, holomorphicity forces infinite differentiability and power series expansions.

\begin{corollary}[Cauchy estimates and $C^\infty$-regularity]
	\label{cor:Cinf-refactored}
	Assume initially that $f\in C^1(U)$ satisfies the Cauchy--Riemann equations.
	Then $f\in C^\infty(U)$.
	Moreover, for any integer $n\ge 0$ and any simple closed curve $\gamma$
	enclosing $z_0$, provided that $\gamma$ and its interior lie in $U$,
	\[
	f^{(n)}(z_0)=\frac{n!}{2\pi i}\int_\gamma \frac{f(z)}{(z-z_0)^{n+1}}\,dz.
	\]
	In particular, if $\overline{B(z_0,R)}\subset U$ and
	$M_R=\max_{|z-z_0|=R}|f(z)|$, then
	\[
	|f^{(n)}(z_0)|\le \frac{n!M_R}{R^n}.
	\]
	It follows that $f$ admits a convergent power series expansion around every $z_0\in U$.
\end{corollary}

\paragraph{Justification for differentiation under the integral.}
Differentiate under the integral sign. The only analytic input is that $f$ is bounded on $\gamma$ and
the kernels $(z-z_0)^{-(n+1)}$ behave uniformly on $\gamma$ away from $z_0$.

\begin{proof}
	Fix $z_0$ and choose $\gamma$ enclosing $z_0$.
	From CIF,
	\[
	f(z_0)=\frac{1}{2\pi i}\int_\gamma \frac{f(z)}{z-z_0}\,dz.
	\]
	Differentiate both sides with respect to $z_0$ (justified since the integrand is smooth in $z_0$ and uniformly
	bounded on $\gamma$ away from $z_0$), obtaining
	\[
	f'(z_0)=\frac{1}{2\pi i}\int_\gamma \frac{f(z)}{(z-z_0)^2}\,dz.
	\]
	Iterating yields the general formula with $n!$.
	On the circle $|z-z_0|=R$, its absolute value is bounded by
	\[
	\frac{n!}{2\pi}\,(2\pi R)\frac{M_R}{R^{n+1}}
	=\frac{n!M_R}{R^n},
	\]
	which is the Cauchy estimate.
	To obtain the power series, keep a circle $|w-z_0|=R$ inside $U$ and
	let $|z-z_0|<R$. Uniformly in $w$ on that circle,
	\[
	\frac{1}{w-z}
	=\frac{1}{w-z_0}\sum_{k=0}^{\infty}
	  \left(\frac{z-z_0}{w-z_0}\right)^k.
	\]
	Insert this expansion into Cauchy's formula for $f(z)$ and integrate
	term by term. The coefficient of $(z-z_0)^k$ is
	$(2\pi i)^{-1}\int_{|w-z_0|=R}f(w)(w-z_0)^{-k-1}\,dw
	=f^{(k)}(z_0)/k!$, as claimed.
\end{proof}

\begin{remark}[Removing the preliminary $C^1$ hypothesis]
	Goursat's theorem proves Cauchy's theorem for a merely complex-differentiable
	function by repeatedly subdividing a triangle. Applying the resulting Cauchy
	formula gives the regularity proved above, so every holomorphic function is
	automatically analytic. The $C^1$ version was presented first because its proof
	is the transparent Green--Stokes calculation; the general version has the same
	conclusions and will be used from this point onward. See
	\cite{Ahlfors1979,Conway1978} for the triangle-subdivision proof.
\end{remark}

\begin{exercise}[Derivatives of $e^z$ from the integral formula]
	Use Corollary~\ref{cor:Cinf-refactored} with $f(z)=e^z$ to compute $f^{(n)}(z_0)$.
	Conclude $e^{(n)}(z_0)=e^{z_0}$ for all $n$.
\end{exercise}

\subsubsection{Morera's theorem (one direction): integrals detect holomorphicity}

\paragraph{Converse direction.}
CIT says that holomorphicity makes all closed integrals vanish.
Morera says (under continuity) the converse holds: vanishing integrals force holomorphicity.
This is a bridge between integral conditions and differential conditions ($\overline\partial f=0$).

\begin{corollary}[Morera's theorem, one direction]
	\label{cor:Morera-refactored}
	If $f$ is continuous on $U$ and $\int_\gamma f(z)\,dz=0$ for every simple closed curve $\gamma\subset U$,
	then $f$ is holomorphic on $U$.
\end{corollary}

\paragraph{Construction of a local primitive.}
On a disk, vanishing integrals around triangles makes the integral from a fixed
base point independent of the polygonal path. This produces a primitive of $f$.

\begin{proof}
	Fix a disk $B=B(z_0,r)\subset U$. For $z\in B$, define $F(z)$ by
	integrating $f(w)\,dw$ along the straight segment from $z_0$ to $z$.
	If $z,z+h\in B$, the two straight segments from $z_0$ together with the
	segment from $z$ to $z+h$ form the boundary of a triangle. When its
	vertices are distinct and noncollinear, its integral is zero by
	hypothesis; in the collinear case, the directed segments cancel and
	the same conclusion holds. Thus
	\[
	F(z+h)-F(z)=\int_{[z,z+h]}f(w)\,dw
	=h\int_0^1 f(z+th)\,dt.
	\]
	Consequently
	\[
	\frac{F(z+h)-F(z)}{h}=\int_0^1f(z+th)\,dt\longrightarrow f(z)
	\]
	by continuity. Thus $F'=f$ on $B$. Since derivatives of holomorphic
	functions are holomorphic, $f$ is holomorphic on $B$. Such disks cover $U$.
\end{proof}

\begin{exercise}[Morera's theorem using polygonal paths]
Rework the proof of Corollary~\ref{cor:Morera-refactored} by:
	\begin{enumerate}
		\item proving path-independence on a disk from the vanishing of closed integrals,
		\item defining a primitive $F$ and proving $F'(z)=f(z)$ using the definition of complex derivative,
		\item concluding holomorphicity.
	\end{enumerate}
\end{exercise}

\subsubsection{Examples and first applications}

\begin{example}[Integrals of monomials around a circle]
	Let $\gamma=\partial B(z_0,r)$ oriented counterclockwise. For $n\in\Z$,
	\[
	\int_\gamma (z-z_0)^n\,dz
	=
	\begin{cases}
		0,& n\neq -1,\\[0.3em]
		2\pi i,& n=-1.
	\end{cases}
	\]
	\emph{Reason.} If $n\ge 0$, the integrand is holomorphic, so CIT gives $0$.
	If $n\le -2$, the function $(z-z_0)^{n+1}/(n+1)$ is a single-valued
	primitive on $\C\setminus\{z_0\}$. The special case $n=-1$ is
	Lemma~\ref{lem:1overz-refactored}.
\end{example}

\begin{example}[Branch cut and topology: $z^{3/2}$]
	As a preview of branches, define $f(z)=z^{3/2}$ on
	$\C\setminus(-\infty,0]$ using the principal argument. It is holomorphic
	on this slit plane. If one analytically continues its value once around
	$0$ in the punctured plane, the argument gains $2\pi$ and the value
	changes sign. Thus the branch cannot extend to a single-valued
	holomorphic function on all of $\C\setminus\{0\}$.
	Cauchy's integral formula still applies to every simple contour whose
	interior lies in the slit plane. What fails is the existence of one
	branch on the entire punctured plane; we construct branches carefully
	in Subsection~\ref{subsec:identity-log} below.
\end{example}

\subsubsection{Exercises}

\begin{exercise}[Core: Green's theorem computation]
	Let $D$ be a bounded piecewise smooth Jordan domain and let $f=u+iv$
	be $C^1$ on a neighborhood of $\overline D$ and satisfy the
	Cauchy--Riemann equations there. Apply Green's theorem separately to
	$u\,dx-v\,dy$ and $v\,dx+u\,dy$ on $\partial D$. Identify the two area
	integrands and explain why each vanishes pointwise. This reproduces the
	mechanism of the Cauchy integral theorem.
\end{exercise}

\begin{exercise}[Core: CIF on a square]
	Let $\gamma$ be the positively oriented square centered at $0$ with side length $2$.
	Compute $\displaystyle \int_\gamma \frac{\cos z}{z}\,dz$ using CIF and interpret the answer.
\end{exercise}

\begin{exercise}[Optional: independence of curve for CIF]
	Let $f\in\operatorname{Hol}(U)$, and let $\gamma_1,\gamma_2$ be positively oriented
	simple closed curves in $U\setminus\{z_0\}$ that both enclose $z_0$
	and are homotopic in $U\setminus\{z_0\}$.
	Show
	\[
	\int_{\gamma_1}\frac{f(z)}{z-z_0}\,dz=\int_{\gamma_2}\frac{f(z)}{z-z_0}\,dz.
	\]
	(Hint: apply homotopy invariance to the holomorphic integrand
	$f(z)/(z-z_0)$ on $U\setminus\{z_0\}$. An annular-region argument
	also works when the curves are nested and the intervening region lies in $U$.)
\end{exercise}

\begin{exercise}[Challenge: winding number preview]
	Let $\gamma$ be a simple closed positively oriented curve and assume $0$ lies inside.
	Show $\int_\gamma \frac{1}{z}\,dz=2\pi i$ by mapping the interior conformally to a disk (or by a homotopy argument
	once you establish that the integral is invariant under deformation avoiding $0$).
\end{exercise}

\subsection{Complex Analytic Functions and Power Series Representations}
\index{power series!holomorphic functions}
\label{subsec:complex-analytic-functions}

\subsubsection{Power series as rigid local models}
\begin{quote}
	\textbf{Goal.} We make precise the rigidity principle:
	\begin{enumerate}
		\item A power series expansion near $p$ determines \emph{everything} about the function near $p$.
		\item Complex-analytic functions are automatically $C^\infty$, with derivatives obtained by termwise differentiation.
		\item The coefficients are unique and satisfy $a_k=\dfrac{f^{(k)}(p)}{k!}$.
		\item \textbf{Holomorphic functions are analytic:} Cauchy's integral formula forces a convergent Taylor series.
	\end{enumerate}
	\textbf{Algebraic-geometry preview.} A convergent power series is a first model of a \emph{germ} and a \emph{local coordinate}:
	near $p$, a holomorphic function has a convergent power series in $(z-p)$.
\end{quote}

\subsubsection{Analyticity and radius of convergence}

\paragraph{Why define analyticity separately?}
In real analysis, differentiability (even smoothness) does \emph{not} imply analytic expandability.
In complex analysis, holomorphicity will imply analyticity; but it is instructive to isolate the notion of
``power series representation'' first.

\begin{definition}[Complex analytic function]
	Let $U\subset\C$ be open. A function $f:U\to\C$ is \emph{complex analytic at $p\in U$} if there exists $r>0$
	such that $B(p,r)\subset U$ and
	\[
	f(z)=\sum_{n=0}^\infty a_n (z-p)^n,\qquad z\in B(p,r),
	\]
	where the series converges uniformly on every compact subset of $B(p,r)$.
	Such an $r$ gives a disk on which the expansion represents $f$. The
	\emph{radius of convergence} of the coefficient series is instead the
	supremum of the radii on which that series converges, whether or not
	those larger disks lie in the originally chosen domain $U$.
\end{definition}

\begin{remark}[What ``uniform on compacts'' buys you]
	For a power series, convergence on compact subsets of its disk of convergence gives:
	\begin{itemize}
		\item termwise integration and, after the derivative series has been
		shown to converge locally uniformly, termwise differentiation,
		\item uniform control of tails,
		\item continuity of the sum, and smoothness by the differentiation argument below.
	\end{itemize}
	Uniform convergence of arbitrary differentiable functions alone does not justify
	passing a derivative through the limit.
	Thus once $f$ is analytic at $p$, its behavior on $B(p,r)$ is encoded in the coefficient sequence $\{a_n\}$.
\end{remark}

\begin{remark}[Radius of convergence and singularities (preview)]
	Later we will relate the maximal radius of convergence of the Taylor series at $p$ to the nearest
	obstruction to holomorphic continuation of that series. A boundary point of a chosen
	domain is not necessarily such an obstruction. For example, the function
	$1/(1-z)$ restricted to $|z|<1/2$ still has the Taylor series
	$\sum_{n\ge0}z^n$ at $0$, whose radius of convergence is $1$;
	the obstruction is at $z=1$, not at the boundary of the restricted domain.
	This connects local expansions to global geometry on
	Riemann surfaces.
\end{remark}

\begin{exercise}[Warm-up: geometric series and radius]
	Let $R>0$. Show that for $|z-p|<R$,
	\[
	\frac{1}{1-\frac{z-p}{R}}=\sum_{n=0}^\infty \frac{(z-p)^n}{R^n}.
	\]
	What is the radius of convergence at $p$? Where is the singularity?
\end{exercise}

\subsubsection{Analytic functions are smooth: termwise differentiation and estimates}

\paragraph{Key mechanism.}
On a slightly smaller disk, the power series behaves like a uniformly convergent sum of smooth functions,
and the derivative series converges uniformly by the Weierstrass $M$-test.
This justifies differentiating term-by-term.

\begin{lemma}[Smoothness and termwise differentiation of analytic functions]
	\label{lem:analytic-smooth-refactored}
	Let $f(z)=\sum_{n=0}^\infty a_n (z-p)^n$ on $B(p,r)$.
	Then $f\in C^\infty(B(p,r))$, and for every integer $k\ge 0$,
	\[
	f^{(k)}(z)=\sum_{n=k}^\infty n(n-1)\cdots(n-k+1)\,a_n (z-p)^{n-k},
	\qquad z\in B(p,r),
	\]
	with uniform convergence on every compact subset of $B(p,r)$.
\end{lemma}

\begin{proof}
	Fix $z_0\in B(p,r)$. Choose $0<\rho<r-|z_0-p|$, so that $\overline{B(z_0,\rho)}\subset B(p,r)$.
	Set
	\[
	R:=\sup_{z\in\overline{B(z_0,\rho)}}|z-p| \le |z_0-p|+\rho < r.
	\]
	Choose $R_1$ with $R<R_1<r$. Since the power series converges
	absolutely at every point of modulus $R_1$, its terms are bounded:
	there is a constant $M>0$ such that
	\[
	|a_n|R_1^n\le M\qquad(n\ge0).
	\]
	Thus, with $q:=R/R_1<1$,
	\[
	|a_n|R^n=(|a_n|R_1^n)q^n\le Mq^n.
	\]
	
	Now fix $k\ge 1$ and consider the $k$-th derivative term:
	\[
	n(n-1)\cdots(n-k+1)\,a_n (z-p)^{n-k}.
	\]
	For $z\in\overline{B(z_0,\rho)}$ we have $|z-p|\le R$, so
	\[
	\Big|n(n-1)\cdots(n-k+1)\,a_n (z-p)^{n-k}\Big|
	\le n^k |a_n| R^{n-k}.
	\]
	Using the preceding geometric bound, we get
	\[
	n^k |a_n| R^{n-k}
	= n^k (|a_n|R^n)\,R^{-k}
	\le (M R^{-k})\, n^k q^n.
	\]
	The numerical series $\sum_{n\ge0} n^k q^n$ converges, so by the Weierstrass $M$-test,
	the derivative series converges uniformly on $\overline{B(z_0,\rho)}$.
	
	Uniform convergence of the derivative series implies we may differentiate termwise on $\overline{B(z_0,\rho)}$,
	hence on a neighborhood of $z_0$. Since $z_0$ was arbitrary, $f\in C^\infty(B(p,r))$ and the displayed formula holds.
\end{proof}

\begin{remark}[What the estimates are really doing]
	We shrink from $B(p,r)$ to a smaller closed disk where $|z-p|\le R<r$.
	On that closed disk, the tails of the power series are controlled by a geometric sequence,
	which is precisely what the $M$-test needs to justify exchanging $\sum$ and $\frac{d}{dz}$.
\end{remark}

\begin{exercise}[Do it once by hand: differentiate a power series]
	Let $f(z)=\sum_{n=0}^\infty a_n (z-p)^n$ converge for $|z-p|<r$.
	Fix $0<R<r$ and show directly (using $M$-test) that
	\[
	\sum_{n=1}^\infty n a_n (z-p)^{n-1}
	\]
	converges uniformly on $|z-p|\le R$.
\end{exercise}

\subsubsection{Uniqueness of coefficients: Taylor formula and rigidity}

\paragraph{Why uniqueness matters.}
If a function has a power series expansion, the coefficients are not ``choices''---they are forced.
This is the first hard rigidity statement: the infinite jet at one point determines the function near that point.

\begin{theorem}[Uniqueness of coefficients and Taylor formula]
	\label{thm:coeff-unique-refactored}
	Let $f(z)=\sum_{n=0}^\infty a_n (z-p)^n$ on $B(p,r)$.
	Then for every $k\ge 0$,
	\[
	a_k=\frac{f^{(k)}(p)}{k!}.
	\]
	In particular, the analytic expansion of $f$ at $p$ is unique.
\end{theorem}

\begin{proof}
	By Lemma~\ref{lem:analytic-smooth-refactored}, we may differentiate termwise:
	\[
	f^{(k)}(z)=\sum_{n=k}^\infty n(n-1)\cdots(n-k+1)\,a_n (z-p)^{n-k}.
	\]
	Evaluate at $z=p$. Then $(p-p)^{n-k}=0$ for all $n>k$, and only the $n=k$ term survives:
	\[
	f^{(k)}(p)=k!\,a_k.
	\]
	Hence $a_k=f^{(k)}(p)/k!$.
\end{proof}

\begin{remark}[Rigidity rephrased: infinite jets determine germs]
	For analytic functions, the infinite jet
	\[
	\big(f(p),f'(p),f''(p),\dots\big)
	\]
	determines the entire function on some disk around $p$.
	In real analysis, there exist $C^\infty$ functions whose Taylor series at a point does \emph{not} recover the function.
	Thus complex differentiability imposes a rigidity that smoothness alone does not provide in real analysis.
\end{remark}

\begin{exercise}[Two analytic functions with same jet]
	Assume $f$ and $g$ are analytic near $p$ and satisfy $f^{(k)}(p)=g^{(k)}(p)$ for all $k\ge 0$.
	Show that their power series coefficients at $p$ agree, hence $f=g$ on a neighborhood of $p$.
\end{exercise}

\subsubsection{Holomorphic functions are analytic: Cauchy expansion and coefficients}

\paragraph{Why this is the main theorem here.}
So far, analyticity was a separate definition.
Now we prove holomorphicity \emph{forces} analyticity using Cauchy's integral formula.
This is the central rigidity theorem: complex differentiability implies power series.

\begin{proposition}[Holomorphic functions are analytic (Cauchy expansion)]
	\label{prop:holomorphic-analytic-refactored}
	Let $U\subset\C$ be open and let $f:U\to\C$ be holomorphic.
	Then for each $p\in U$ there exists $r>0$ such that for all $z\in B(p,r)$,
	\[
	f(z)=\sum_{k=0}^\infty a_k (z-p)^k,
	\]
	where for any $0<\rho<r$,
	\[
	a_k=\frac{1}{2\pi i}\int_{|\zeta-p|=\rho}\frac{f(\zeta)}{(\zeta-p)^{k+1}}\,d\zeta.
	\]
	Equivalently, $a_k=\dfrac{f^{(k)}(p)}{k!}$.
\end{proposition}

\paragraph{Proof strategy.}
CIF gives an integral representation involving the kernel $\frac{1}{\zeta-z}$.
When $|z-p|<|\zeta-p|=\rho$, expand $\frac{1}{\zeta-z}$ as a geometric series in $(z-p)$,
then interchange sum and integral (uniform convergence on $|z-p|\le r_0<\rho$).

\begin{proof}
	Translate so that $p=0$ (otherwise replace $z$ by $z-p$ and $\zeta$ by $\zeta-p$).
	Choose $r>0$ such that $\overline{B(0,r)}\subset U$.
	Fix $0<\rho<r$. For $|z|<\rho$, Cauchy's integral formula gives
	\[
	f(z)=\frac{1}{2\pi i}\int_{|\zeta|=\rho}\frac{f(\zeta)}{\zeta-z}\,d\zeta.
	\]
	For $|\zeta|=\rho$ and $|z|<\rho$,
	\[
	\frac{1}{\zeta-z}
	=\frac{1}{\zeta}\cdot \frac{1}{1-\frac{z}{\zeta}}
	=\sum_{k=0}^\infty \frac{z^k}{\zeta^{k+1}},
	\]
	and the series converges uniformly for $|z|\le r_0<\rho$ since $\big|\frac{z}{\zeta}\big|\le \frac{r_0}{\rho}<1$.
	Hence we may substitute into the integral and exchange sum and integral:
	\begin{align*}
		f(z)
		&=\frac{1}{2\pi i}\int_{|\zeta|=\rho} f(\zeta)\sum_{k=0}^\infty \frac{z^k}{\zeta^{k+1}}\,d\zeta\\
		&=\sum_{k=0}^\infty \left(\frac{1}{2\pi i}\int_{|\zeta|=\rho}\frac{f(\zeta)}{\zeta^{k+1}}\,d\zeta\right) z^k.
	\end{align*}
	Define
	\[
	a_k:=\frac{1}{2\pi i}\int_{|\zeta|=\rho}\frac{f(\zeta)}{\zeta^{k+1}}\,d\zeta.
	\]
	This yields a power series expansion valid for $|z|<\rho$. Since $\rho$ can be chosen with $0<\rho<r$,
	we obtain analyticity on $B(0,r)$.
	
	Finally, by Cauchy's differentiation formula (Corollary~\ref{cor:Cinf-refactored} in the previous subsection),
	\[
	f^{(k)}(0)=\frac{k!}{2\pi i}\int_{|\zeta|=\rho}\frac{f(\zeta)}{\zeta^{k+1}}\,d\zeta
	\quad\mathrel{\text{and hence}}\quad
	a_k=\frac{f^{(k)}(0)}{k!}.
	\]
\end{proof}

\begin{remark}[Holomorphic $\iff$ analytic near a point]
	Combining Proposition~\ref{prop:holomorphic-analytic-refactored} with Theorem~\ref{thm:coeff-unique-refactored}, we get
	These results identify three equivalent local conditions: $f$ is
	holomorphic near $p$; it is analytic there; and it admits a
	convergent Taylor series centered at $p$.
	The corresponding implication from smoothness to analyticity fails for real-valued functions.
\end{remark}

\begin{exercise}[Uniform convergence justification]
	Fix $0<r_0<\rho$.
	Show that for $|\zeta|=\rho$,
	\[
	\sum_{k=0}^\infty \left|\frac{z^k}{\zeta^{k+1}}\right|
	\le \frac{1}{\rho}\sum_{k=0}^\infty \left(\frac{r_0}{\rho}\right)^k
	=\frac{1}{\rho}\cdot\frac{1}{1-r_0/\rho}
	\]
	for all $|z|\le r_0$.
	Use this to justify exchanging $\sum$ and $\int$ in the proof above.
\end{exercise}

\subsubsection{Examples: radius of convergence as a detector of singularities}

\paragraph{Example A: a pole determines the radius.}
\begin{example}[Geometric series and a boundary pole]
	For $|z|<1$,
	\[
	\frac{1}{1-z}=\sum_{n=0}^\infty z^n.
	\]
	The radius of convergence at $0$ is $1$, and the nearest singularity is the pole at $z=1$ (distance $1$ from $0$).
\end{example}

\paragraph{Example B: entire functions have infinite radius.}
\begin{example}[Exponential function]
	For all $z\in\C$,
	\[
	e^{z}=\sum_{n=0}^\infty \frac{z^n}{n!}.
	\]
	The radius of convergence is $\infty$ since $\limsup_{n\to\infty}|1/n!|^{1/n}=0$.
\end{example}

\paragraph{Example C: branch points and local expansions.}
\begin{example}[A local branch of $\log z$ around $z=1$]
	On $|z-1|<1$, define a holomorphic branch of $\log z$ and expand
	\[
	\log z=\log\big(1+(z-1)\big)=\sum_{n=1}^\infty (-1)^{n+1}\frac{(z-1)^n}{n}.
	\]
	The radius of convergence is $1$, and the obstruction comes from the branch point at $z=0$.
\end{example}

\begin{exercise}[Compute a radius from coefficients]
	Use the root test to show that the power series $\sum_{n=0}^\infty a_n (z-p)^n$ has radius $r$ satisfying
	\[
	\limsup_{n\to\infty}|a_n|^{1/n}=\frac{1}{r}.
	\]
	Use the extended-real conventions $1/0=\infty$ and $1/\infty=0$;
	the coefficient $a_0$ does not affect the limit.
\end{exercise}

\subsubsection{Exercises (compute first, then abstract)}

\medskip
\noindent\textbf{Type A (Computational).}

\begin{exercise}
	Let $f$ be holomorphic on $B(0,1)$ and assume $|f(z)|\le M$ on $|z|<1$.
	Using Cauchy's estimate, prove
	\[
	\left|\frac{f^{(k)}(0)}{k!}\right|\le \frac{M}{\rho^k}
	\qquad (0<\rho<1).
	\]
	Deduce that if $|f(z)|\le M$ on $|z|<1$, then the Taylor coefficients satisfy $|a_k|\le M$.
\end{exercise}

\begin{exercise}
	Compute the power series expansion on $|z|<1$:
	\[
	\frac{1}{1-z^2}=\sum_{n=0}^\infty z^{2n}.
	\]
	Find its radius of convergence and relate it to the poles at $z=\pm 1$.
\end{exercise}

\begin{exercise}
	Find the Taylor series of $\frac{1}{(1-z)^2}$ on $|z|<1$ by differentiating the geometric series.
	Then use Theorem~\ref{thm:coeff-unique-refactored} to identify $f^{(n)}(0)$.
\end{exercise}

\medskip
\noindent\textbf{Type B (Conceptual / geometric).}

\begin{exercise}[Removable singularity and radius]
	Define $f$ on $U=B(0,2)$ by $f(z)=\frac{\sin z}{z}$ for $z\neq 0$ and $f(0)=1$.
	\begin{enumerate}
		\item Show that $f$ is holomorphic on $U$ (removable singularity).
		\item Compute the Taylor series at $0$ and determine its radius of convergence.
		\item Explain why the removed function is entire and hence its Taylor
		series has infinite radius of convergence.  The point at infinity does
		not impose a finite Euclidean radius on a Taylor series centered at $0$.
	\end{enumerate}
\end{exercise}

\begin{exercise}[Germs and jets]
	Assume $f$ and $g$ are holomorphic near $p$ and satisfy $f^{(k)}(p)=g^{(k)}(p)$ for all $k\ge 0$.
	Show $f=g$ near $p$. Interpret this as: \emph{a holomorphic germ is determined by its infinite jet.}
\end{exercise}

\begin{exercise}[Real-variable contrast]
	Give an example of a smooth function $h:\R\to\R$ whose Taylor series at $0$ is identically $0$ but $h\not\equiv 0$.
	Explain why this cannot happen for holomorphic functions (use the uniqueness of coefficients and the identity principle).
\end{exercise}

\subsection{Cauchy Estimates, Liouville's Theorem, and Applications}
\label{subsec:cauchy-liouville}

\subsubsection{From boundary estimates to global rigidity}
Cauchy's integral formula turns boundary bounds into quantitative control of
the function in the interior. Three consequences will be used repeatedly:

\begin{itemize}
	\item \textbf{Local quantitative rigidity (Cauchy estimates).}
	Boundary size on a circle controls every derivative at the center.
	\item \textbf{Global qualitative rigidity (Liouville).}
	Every bounded entire function is constant.
	\item \textbf{From growth to algebra (polynomial growth forces a polynomial).}
	If an entire function grows like $|z|^n$, then all derivatives beyond order $n$ vanish.
\end{itemize}

These results are prototypes for later complex-geometry statements:
growth conditions at infinity (or at divisors) cut down holomorphic objects to finite-dimensional spaces.

\subsubsection{Cauchy estimates}

\paragraph{Geometric intuition.}
Fix a disk $\overline{B(p,r)}\subset U$. Once $|f|$ is controlled on the boundary circle, holomorphicity forces
the interior Taylor coefficients to be controlled, hence controls every derivative at $p$.
Thus boundary estimates control every derivative at an interior point.

\begin{proposition}[Cauchy estimates: boundary controls derivatives]
	\label{prop:cauchy-estimates}
	Let $U\subset\C$ be open, $p\in U$, and assume $\overline{B(p,r)}\subset U$.
	Let $f\in\Hol(U)$ and define
	\[
	M_r:=\max_{|\zeta-p|=r}|f(\zeta)|.
	\]
	Then for every integer $k\ge 0$,
	\[
	|f^{(k)}(p)|\le \frac{k!}{r^k}\,M_r.
	\]
\end{proposition}

\begin{proof}
	By Cauchy's integral formula for derivatives,
	\[
	f^{(k)}(p)=\frac{k!}{2\pi i}\int_{|\zeta-p|=r}\frac{f(\zeta)}{(\zeta-p)^{k+1}}\,d\zeta.
	\]
	Take absolute values:
	\[
	|f^{(k)}(p)|
	\le \frac{k!}{2\pi}\int_{|\zeta-p|=r}\frac{|f(\zeta)|}{|\zeta-p|^{k+1}}\,|d\zeta|
	\le \frac{k!}{2\pi}\int_{|\zeta-p|=r}\frac{M_r}{r^{k+1}}\,|d\zeta|.
	\]
	Since the circle has length $2\pi r$, we obtain
	\[
	|f^{(k)}(p)|\le \frac{k!}{2\pi}\cdot \frac{M_r}{r^{k+1}}\cdot 2\pi r=\frac{k!}{r^k}M_r.
	\]
\end{proof}

\begin{remark}[A concrete reading of the inequality]
	For $k=1$ the estimate says:
	\[
	|f'(p)|\le \frac{1}{r}\max_{|\zeta-p|=r}|f(\zeta)|.
	\]
	So if $f$ is small on a circle, it cannot change rapidly at the center.
	Higher $k$ bounds control the corresponding Taylor coefficients and hence
	the higher-order local variation of $f$ at $p$.
\end{remark}

\paragraph{Sharpness and extremizers.}
The estimate is sharp (as an inequality with the constant $\frac{k!}{r^k}$).

\begin{lemma}[Sharpness of Cauchy estimates]
	\label{lem:cauchy-sharp}
	Fix $k\ge 0$ and $r>0$. The entire function $f(\zeta)=(\zeta-p)^k$ satisfies
	\[
	|f^{(k)}(p)|=\frac{k!}{r^k}M_r.
	\]
\end{lemma}

\begin{proof}
	We have $f^{(k)}(p)=k!$ and $M_r=\max_{|\zeta-p|=r}|\zeta-p|^k=r^k$.
\end{proof}

\paragraph{Worked example}
\begin{example}[Using Cauchy estimates for $e^z$]
	Let $f(z)=e^z$ and $p=0$. On $|z|=r$ we have $|e^z|=e^{\Re z}\le e^{|z|}=e^r$.
	Hence $M_r\le e^r$ and
	\[
	|f^{(k)}(0)|\le \frac{k!}{r^k}e^r.
	\]
	Since $f^{(k)}(0)=1$, we get the inequality
	\[
	1\le \frac{k!}{r^k}e^r\qquad\text{for all }r>0.
	\]
	For $k\ge1$, choosing $r=k$ gives the elementary lower bound
	$k!\ge k^k e^{-k}$.
\end{example}

\begin{exercise}[Cauchy estimate on a disk]
	Let $f$ be holomorphic on a neighborhood of $\overline{B(0,2)}$ and assume
	$|f(z)|\le 5$ on $|z|=2$.
	Show that $|f^{(3)}(0)|\le \frac{3!}{2^3}\cdot 5$.
\end{exercise}

\begin{exercise}[Coefficient bounds]
	If $f(z)=\sum_{n\ge 0} a_n (z-p)^n$ is holomorphic on $B(p,r)$, fix
	$0<\rho<r$ and show that
	\[
	|a_n|=\left|\frac{f^{(n)}(p)}{n!}\right|\le \frac{M_\rho}{\rho^n},
	\qquad M_\rho:=\max_{|z-p|=\rho}|f(z)|.
	\]
	Interpret this as: boundary bounds control Taylor coefficients.
\end{exercise}

\subsubsection{Liouville's theorem and immediate consequences}

\paragraph{What changes from local to global?}
On $\C$, you may take circles of arbitrarily large radius. Cauchy estimates then force derivatives to vanish.

\begin{theorem}[Liouville: bounded entire functions are constant]
	\label{thm:liouville}
	If $f:\C\to\C$ is entire and bounded, then $f$ is constant.
\end{theorem}

\begin{proof}
	Assume $|f(z)|\le M$ for all $z\in\C$.
	Fix $p\in\C$ and apply Proposition~\ref{prop:cauchy-estimates} to the disk $\overline{B(p,R)}$ for any $R>0$:
	\[
	|f'(p)|\le \frac{1!}{R}\max_{|\zeta-p|=R}|f(\zeta)|\le \frac{M}{R}.
	\]
	Letting $R\to\infty$ yields $f'(p)=0$. Since $p$ was arbitrary, $f'\equiv 0$ and $f$ is constant.
\end{proof}

\begin{remark}[Compact Riemann surfaces: a preview]
	On a compact Riemann surface $X$, every holomorphic function is also
	constant, but this does not follow by applying Liouville directly:
	Liouville concerns functions on $\C$. On $X$, compactness gives a
	maximum of the modulus and the local maximum-modulus principle forces
	the function to be constant. We prove that principle in
	Subsection~\ref{subsec:max-min-modulus}. The common
	lesson is that global conditions can force strong holomorphic rigidity.
\end{remark}

\paragraph{A standard corollary: the global form of the maximum principle.}
\begin{corollary}[A nonconstant entire function is unbounded]
	If $f$ is entire and nonconstant, then $f(\C)$ is unbounded.
\end{corollary}

\begin{proof}
	Contrapositive of Theorem~\ref{thm:liouville}.
\end{proof}

\begin{exercise}[Trigonometric functions are unbounded on $\C$]
	Show that $\sin z$ and $\cos z$ are unbounded on $\C$.
	(Hint: evaluate at $z=iy$ and let $y\to\infty$.)
\end{exercise}

\subsubsection{Polynomial growth forces a polynomial}

\paragraph{Key idea.}
Polynomial growth beats the decay $R^{-k}$ in Cauchy estimates when $k$ is large.
So sufficiently high derivatives must vanish.

\begin{lemma}[Polynomial growth kills high derivatives (centered version)]
	\label{lem:poly-growth-0}
	Let $f$ be entire. Assume there exist constants $C>0$ and integer $n\ge 0$ such that
	\[
	|f(z)|\le C(1+|z|)^n \quad \forall z\in\C.
	\]
	Then for every integer $k>n$,
	\[
	f^{(k)}(0)=0.
	\]
\end{lemma}

\begin{proof}
	Fix $k>n$. On $|z|=R$ we have $|f(z)|\le C(1+R)^n$.
	By Cauchy estimates (Proposition~\ref{prop:cauchy-estimates} with $p=0$),
	\[
	|f^{(k)}(0)|\le \frac{k!}{R^k}\max_{|z|=R}|f(z)|\le \frac{k!}{R^k}C(1+R)^n.
	\]
	Let $R\to\infty$. Since $(1+R)^n/R^k\to 0$ for $k>n$, the right-hand side tends to $0$,
	so $f^{(k)}(0)=0$.
\end{proof}

\paragraph{Upgrade: vanishing at every center.}
The lemma at $0$ implies the same at every point by translating.

\begin{lemma}[Polynomial growth kills high derivatives (any center)]
	\label{lem:poly-growth-any}
	Under the assumptions of Lemma~\ref{lem:poly-growth-0}, for every $p\in\C$ and every $k>n$,
	\[
	f^{(k)}(p)=0.
	\]
\end{lemma}

\begin{proof}
	Fix $p$ and define $g(z)=f(z+p)$. Then $g$ is entire and satisfies the same type of growth bound:
	\[
	|g(z)|=|f(z+p)|\le C(1+|z+p|)^n \le C'(1+|z|)^n
	\]
	for a constant $C'$ (use $1+|z+p|\le 1+|z|+|p|\le (1+|p|)(1+|z|)$).
	Apply Lemma~\ref{lem:poly-growth-0} to $g$ to get $g^{(k)}(0)=0$ for $k>n$.
	But $g^{(k)}(0)=f^{(k)}(p)$, so $f^{(k)}(p)=0$.
\end{proof}

\begin{proposition}[Entire function of polynomial growth is a polynomial]
	\label{prop:poly-growth-poly}
	If $f$ is entire and $|f(z)|\le C(1+|z|)^n$, then $f$ is a polynomial of degree at most $n$.
\end{proposition}

\begin{proof}
	By Lemma~\ref{lem:poly-growth-0}, $f^{(k)}(0)=0$ for all $k>n$.
	Hence the Taylor series at $0$ truncates:
	\[
	f(z)=\sum_{k=0}^\infty \frac{f^{(k)}(0)}{k!}z^k=\sum_{k=0}^n \frac{f^{(k)}(0)}{k!}z^k,
	\]
	so $f$ is a polynomial of degree $\le n$.
\end{proof}

\begin{remark}[Bridge to Riemann--Roch language]
	On a compact Riemann surface $X$, controlling pole orders (growth near finitely many points)
	produces finite-dimensional spaces of meromorphic functions. Proposition~\ref{prop:poly-growth-poly}
	is the simplest flat-model analogue of that principle.
\end{remark}

\begin{exercise}[Linear growth implies affine]
	Let $f$ be entire and satisfy $|f(z)|\le A+B|z|$ for all $z\in\C$.
	Show that $f$ is a polynomial of degree at most $1$, i.e.\ $f(z)=\alpha z+\beta$.
\end{exercise}

\subsubsection{Fundamental Theorem of Algebra via Liouville}

\begin{theorem}[Fundamental Theorem of Algebra]
	\label{thm:FTA}
	Every nonconstant complex polynomial $P(z)$ has at least one root in $\C$.
\end{theorem}

\begin{proof}
	Assume $P$ is nonconstant and has no zeros. Then $f(z):=1/P(z)$ is entire.
	Write $P(z)=a_n z^n + a_{n-1}z^{n-1}+\cdots + a_0$ with $a_n\neq 0$.
	Factor:
	\[
	P(z)=a_n z^n\Bigl(1+\varepsilon(z)\Bigr),
	\qquad
	\varepsilon(z):=\frac{a_{n-1}}{a_n}\frac{1}{z}+\cdots+\frac{a_0}{a_n}\frac{1}{z^n}.
	\]
	Then $\varepsilon(z)\to 0$ as $|z|\to\infty$. Choose $R>0$ such that $|\varepsilon(z)|\le \frac12$ when $|z|\ge R$.
	For $|z|\ge R$,
	\[
	|P(z)|=|a_n||z|^n|1+\varepsilon(z)|\ge |a_n||z|^n\Bigl(1-\tfrac12\Bigr)=\frac12|a_n||z|^n,
	\]
	so
	\[
	|f(z)|=\frac{1}{|P(z)|}\le \frac{2}{|a_n|}\,|z|^{-n}\le \frac{2}{|a_n|}R^{-n}
	\quad\text{for }|z|\ge R.
	\]
	Thus $f$ is bounded on $\C\setminus B(0,R)$, and since $f$ is continuous on the compact disk $\overline{B(0,R)}$,
	it is bounded there as well. Hence $f$ is bounded on all of $\C$.
	By Liouville (Theorem~\ref{thm:liouville}), $f$ is constant, so $P$ is constant, contradiction.
	Therefore $P$ must have a root.
\end{proof}

\begin{remark}[Role of the growth hypothesis]
	FTA is a global rigidity theorem disguised as algebra:
	``a polynomial cannot avoid $0$ everywhere without forcing $1/P$ to be a bounded entire function.''
\end{remark}

\begin{exercise}[Concrete estimate for a specific polynomial]
	Let $n\ge1$ and $P(z)=z^n+1$. Show that for $|z|\ge 2$,
	\[
	\left|\frac{1}{P(z)}\right|\le \frac{2}{|z|^n}.
	\]
	Conclude $1/P$ is bounded outside $B(0,2)$ and complete the Liouville argument to deduce that $P$ has a root.
\end{exercise}

\subsubsection{Uniform limits, holomorphicity, and convergence of derivatives}

\paragraph{Uniform limits and holomorphicity.}
Many existence/compactness arguments in complex analysis (e.g.\ normal families, Montel) pass to limits of holomorphic functions.
Cauchy's formula makes holomorphicity stable under uniform convergence on compact subsets.

\begin{proposition}[Uniform limit on compacts of holomorphic functions is holomorphic]
	\label{prop:uniform-limit-holo}
	Let $U\subset\C$ be open and $f_n\in\Hol(U)$.
	Assume $f_n\to f$ uniformly on every compact subset of $U$.
	Then $f\in\Hol(U)$.
\end{proposition}

\begin{proof}
	Fix $p\in U$ and choose $r>0$ with $\overline{B(p,r)}\subset U$.
	For any $z$ with $|z-p|<r$, Cauchy's integral formula gives
	\[
	f_n(z)=\frac{1}{2\pi i}\int_{|\zeta-p|=r}\frac{f_n(\zeta)}{\zeta-z}\,d\zeta.
	\]
	Since the circle $|\zeta-p|=r$ is compact and $f_n\to f$ uniformly on it, we may pass to the limit inside the integral:
	\[
	f(z)=\frac{1}{2\pi i}\int_{|\zeta-p|=r}\frac{f(\zeta)}{\zeta-z}\,d\zeta.
	\]
	The right-hand side is a holomorphic function of $z$ on $B(p,r)$ (as an integral of a holomorphic kernel in $z$),
	so $f$ is holomorphic on $B(p,r)$. Since $p$ was arbitrary, $f\in\Hol(U)$.
\end{proof}

\begin{theorem}[Convergence of derivatives under uniform convergence on compacts]
	\label{thm:derivative-convergence}
	Under the assumptions of Proposition~\ref{prop:uniform-limit-holo}, for every integer $k\ge 0$,
	\[
	f_n^{(k)} \to f^{(k)}
	\quad\text{uniformly on every compact subset of }U.
	\]
\end{theorem}

\begin{proof}
	Let $K\subset U$ be compact. Compactness of $K$ and openness of $U$ give a
	$\rho>0$ such that
	\[
	K_\rho:=\{\zeta:\operatorname{dist}(\zeta,K)\le\rho\}\subset U.
	\]
	This is a compact subset of $U$. For each $z\in K$, Cauchy's formula on
	the circle $|\zeta-z|=\rho$ gives
	\[
	f_n^{(k)}(z)-f^{(k)}(z)
	=\frac{k!}{2\pi i}\int_{|\zeta-z|=\rho}
	\frac{f_n(\zeta)-f(\zeta)}{(\zeta-z)^{k+1}}\,d\zeta.
	\]
	The circle is contained in $K_\rho$, so the ML estimate yields, uniformly
	for $z\in K$,
	\[
	|f_n^{(k)}(z)-f^{(k)}(z)|
	\le \frac{k!}{\rho^k}\sup_{\zeta\in K_\rho}|f_n(\zeta)-f(\zeta)|.
	\]
	The right-hand side tends to $0$ by uniform convergence on the compact
	set $K_\rho$, proving uniform convergence on $K$.
\end{proof}

\begin{remark}[Normal families: preview]
	Montel's theorem will upgrade ``uniform boundedness on compact sets'' into subsequential compactness.
	Proposition~\ref{prop:uniform-limit-holo} and Theorem~\ref{thm:derivative-convergence}
	are the analytic engine behind such compactness arguments.
\end{remark}

\begin{exercise}[Local boundedness from uniform convergence]
	Let $f_n\to f$ uniformly on compact subsets of $U$.
	Prove: for each compact $K\subset U$ there exists $M_K$ such that $|f_n(z)|\le M_K$ for all $z\in K$ and all $n$.
\end{exercise}

\begin{exercise}[Pointwise convergence does not control the boundary]
	Let $f_n(z)=z^n$ on a neighborhood of the closed unit disk.  Determine the
	pointwise limit on $[0,1]$ and show that it is discontinuous at $1$.  Then prove that
	$f_n\to0$ uniformly on every compact subset of the open unit disk.  Explain why these
	two conclusions are consistent and why the compact-open hypothesis in the preceding
	theorems concerns compact subsets of the \emph{open domain}, not its boundary.
\end{exercise}

\subsubsection{Exercises}

\medskip
\noindent\textbf{Type A: Computational / estimate-based.}

\begin{exercise}
	Let $f$ be holomorphic on an open neighborhood of $\overline{B(0,3)}$ and assume $\max_{|z|=3}|f(z)|\le 10$.
	Show that $|f^{(5)}(0)|\le \dfrac{5!}{3^5}\cdot 10$.
\end{exercise}

\begin{exercise}
	Let $f$ be entire and satisfy $|f(z)|\le 7(1+|z|)^4$.
	Prove that $f$ is a polynomial of degree at most $4$, and explain why you cannot conclude degree $\le 3$ from this data.
\end{exercise}

\medskip
\noindent\textbf{Type B: Conceptual / geometric.}

\begin{exercise}
	Explain why Liouville's theorem can be read as:
	``if a holomorphic function does not grow anywhere on $\C$, then it cannot vary anywhere.''
	Compare this with the compact Riemann surface statement ``holomorphic implies constant'':
	what is the geometric reason boundedness becomes automatic on compact domains?
\end{exercise}

\begin{exercise}[Where Liouville stops working]
	Let $f$ be entire and assume $|f(z)|\le A+Be^{|z|}$ for all $z$.
	Explain why Cauchy estimates do not force $f'$ to vanish.
	Suggest a growth condition that \emph{would} still force $f$ to be a polynomial
	(you may refer to Proposition~\ref{prop:poly-growth-poly} as guidance).
\end{exercise}

\subsection{Identity Theorem and the Holomorphic Logarithm}
\index{identity theorem}\index{holomorphic logarithm}
\label{subsec:identity-log}

\subsubsection{Two rigidity principles}
This subsection contains two closely related ``no-room-to-wiggle'' phenomena.

\begin{itemize}
	\item \textbf{Identity theorem / zero-set rigidity.}
	A holomorphic function cannot have a ``thick'' zero set inside a domain unless it is the zero function.
	Equivalently, knowing a holomorphic function on any subset with an interior accumulation point determines it everywhere.
	\item \textbf{Holomorphic logarithm / topology and branches.}
	The obstruction to defining $\log z$ globally on $\C\setminus\{0\}$ is \emph{topological}:
	loops that wind around $0$ force a $2\pi i$ ambiguity. On simply connected domains avoiding $0$,
	one can choose a single-valued holomorphic branch.
\end{itemize}

\subsubsection{Identity theorem and zero sets}

\begin{definition}[Zero set]
	Let $U\subset \C$ be open and $f\in \Hol(U)$. The \emph{zero set} of $f$ is
	\[
	Z(f):=\{z\in U: f(z)=0\}.
	\]
\end{definition}

\paragraph{Key local fact: zeros are isolated unless $f\equiv 0$.}

\begin{lemma}[Zeros are isolated unless everything is zero]
	\label{lem:zeros-isolated}
	Let $U\subset\C$ be a domain and $f\in\Hol(U)$. Then either $f\equiv 0$ on $U$, or else every zero
	$z_0\in Z(f)$ is \emph{isolated}: there exists $r>0$ such that
	\[
	0<|z-z_0|<r \quad\Longrightarrow\quad f(z)\neq 0.
	\]
\end{lemma}

\begin{proof}
	Fix a zero $z_0\in Z(f)$. Since $f$ is holomorphic, it admits a power series expansion on some disk $B(z_0,r_0)\subset U$:
	\[
	f(z)=\sum_{n=0}^\infty a_n (z-z_0)^n \qquad (|z-z_0|<r_0).
	\]
	If $a_n=0$ for all $n$, then $f\equiv 0$ on $B(z_0,r_0)$. To avoid
	using the identity theorem before it has been proved, set
	\[
	S:=\{z\in U:\ f\text{ vanishes on a neighborhood of }z\}.
	\]
	The set $S$ is open. It is also closed: if $z_j\in S$ and $z_j\to z\in U$,
	then $f^{(k)}(z_j)=0$ for every $j$ and every $k\ge0$; continuity of the
	holomorphic derivatives gives $f^{(k)}(z)=0$ for all $k$, so the Taylor
	series at $z$ is zero and $z\in S$. Hence connectedness of $U$ implies
	that a nonempty $S$ equals $U$. In the all-coefficients-zero case,
	therefore, $f\equiv0$ on $U$.
	
	Otherwise, let $m\ge 1$ be the smallest index such that $a_m\neq 0$.
	Factor out $(z-z_0)^m$:
	\[
	f(z)=(z-z_0)^m\,g(z),
	\qquad
	g(z):=\sum_{n=0}^\infty a_{n+m}(z-z_0)^n.
	\]
	Then $g$ is holomorphic on $B(z_0,r_0)$ and satisfies $g(z_0)=a_m\neq 0$.
	By continuity, there exists $0<r\le r_0$ such that $g(z)\neq 0$ for all $|z-z_0|<r$.
	Therefore, in $B(z_0,r)$ the only zero of $f$ is $z_0$ itself (with multiplicity $m$),
	so $z_0$ is isolated.
\end{proof}

\paragraph{A refinement: order of vanishing.}

\begin{definition}[Order of a zero]
	\label{def:order-zero}
	Let $U\subset\C$ be a domain, $f\in\Hol(U)$, and $z_0\in U$ with
	$f(z_0)=0$ but $f\not\equiv 0$.
	The unique integer $m\ge 1$ such that
	\[
	f(z)=(z-z_0)^m g(z),\qquad g\in\Hol(U),\ g(z_0)\neq 0,
	\]
	is called the \emph{order} (or \emph{multiplicity}) of the zero of $f$ at $z_0$,
	and we write
	\[
	\operatorname{ord}_{z_0}(f)=m.
	\]
\end{definition}

\begin{theorem}[Identity theorem]
	\label{thm:identity}
	Let $U\subset\C$ be open and connected, and let $f,g\in\Hol(U)$.
	Assume the coincidence set
	\[
	E:=\{z\in U: f(z)=g(z)\}
	\]
	has a limit point in $U$.
	Then $f\equiv g$ on $U$.
	
	Equivalently: if $h\in\Hol(U)$ has zeros accumulating at a point of $U$, then $h\equiv 0$.
\end{theorem}

\begin{proof}
	Let $h:=f-g\in\Hol(U)$. Then $E=Z(h)$.
	If $h\not\equiv 0$, Lemma~\ref{lem:zeros-isolated} implies every zero of $h$ is isolated,
	so $Z(h)$ cannot have a limit point in $U$, contradicting the hypothesis.
	Therefore $h\equiv 0$, i.e.\ $f\equiv g$.
\end{proof}

\begin{corollary}[Uniqueness of holomorphic extension]
	\label{cor:unique-extension}
	Let $U\subset\C$ be connected and $f,g\in\Hol(U)$.
	If $f=g$ on a subset $A\subset U$ that has an accumulation point in $U$, then $f\equiv g$ on $U$.
\end{corollary}

\begin{example}[Multiplicity at $0$]
	Let $f(z)=z^m e^{z}$ with $m\ge 1$. Then $\operatorname{ord}_{0}(f)=m$ because
	\[
	f(z)=z^m\underbrace{e^z}_{g(z)},\qquad g(0)=1\neq 0.
	\]
	Equivalently, $f^{(k)}(0)=0$ for $k=0,1,\dots,m-1$ but $f^{(m)}(0)\neq 0$.
\end{example}

\begin{exercise}[Order of vanishing via derivatives]
	Suppose $f\in\Hol(U)$ and $z_0\in U$ with $f(z_0)=0$ and $f\not\equiv 0$.
	Show that $\operatorname{ord}_{z_0}(f)=m$ if and only if
	\[
	f(z_0)=f'(z_0)=\cdots=f^{(m-1)}(z_0)=0
	\quad\text{and}\quad
	f^{(m)}(z_0)\neq 0.
	\]
\end{exercise}

\subsubsection{The holomorphic logarithm and branches}

The exponential map $\exp:\C\to\C^{\*}$ satisfies $\exp(w+2\pi i)=\exp(w)$, hence is not injective.
Therefore it cannot have a global holomorphic inverse on $\C^{\*}$.
Instead, $\log z$ is naturally multi-valued.

\begin{definition}[Multi-valued logarithm]
	\label{def:multivalued-log}
	Write $z=re^{i\theta}$ with $r>0$ and $\theta\in\R$.
	Define the multi-valued logarithm
	\[
	\log z := \left\{\ln r + i(\theta+2\pi k)\mid k\in\Z\right\}.
	\]
\end{definition}

\begin{theorem}[Existence of a holomorphic logarithm]
	\label{thm:log-exists}
	Let $\Omega\subset\C^{\*}$ be simply connected, and fix $z_\ast\in\Omega$.
	Then there exists a holomorphic function $L:\Omega\to\C$ such that
	\[
	e^{L(z)}=z\quad\text{for all }z\in\Omega,
	\]
	and $L(z_\ast)$ can be chosen to be any element of the multi-valued set $\log(z_\ast)$.
	Any two such branches differ by an additive constant $2\pi i k$, $k\in\Z$.
\end{theorem}

\begin{proof}
	Since $\Omega\subset\C^{\*}$, the function $1/z$ is holomorphic on $\Omega$.
	Fix $z_\ast\in\Omega$ and define
	\[
	L(z):=\int_{\gamma_{z_\ast,z}}\frac{1}{\zeta}\,d\zeta + c,
	\]
	where $\gamma_{z_\ast,z}$ is any piecewise $C^1$ path in $\Omega$ from $z_\ast$ to $z$.
	
	\emph{Step 1: path-independence.}
	If $\gamma_1,\gamma_2$ are two paths from $z_\ast$ to $z$, then $\gamma_1\cdot\overline{\gamma_2}$ is a closed loop in $\Omega$.
	Because $\Omega$ is simply connected, the loop is null-homotopic in
	$\Omega$. The homotopy form of Cauchy's theorem applied to $1/\zeta$
	therefore gives
	\[
	\int_{\gamma_1}\frac{1}{\zeta}\,d\zeta-\int_{\gamma_2}\frac{1}{\zeta}\,d\zeta
	=\int_{\gamma_1\cdot\overline{\gamma_2}}\frac{1}{\zeta}\,d\zeta=0,
	\]
	so the integral is path-independent and $L$ is well-defined.
	
	\emph{Step 2: holomorphy and derivative.}
	For $z$ near $z_0$, take a short path from $z_0$ to $z$ inside $\Omega$:
	\[
	L(z)-L(z_0)=\int_{z_0}^{z}\frac{1}{\zeta}\,d\zeta.
	\]
	Dividing by $z-z_0$ and letting $z\to z_0$ gives $L'(z_0)=1/z_0$.
	Hence $L$ is holomorphic and satisfies $L'(z)=1/z$.
	
	\emph{Step 3: exponentiating.}
	Define $F(z):=e^{L(z)}/z$ on $\Omega$. Then $F$ is holomorphic and
	\[
	\frac{F'(z)}{F(z)}=L'(z)-\frac{1}{z}=0,
	\]
	so $F$ is constant on connected $\Omega$. Choose the constant $c$ so that $F(z_\ast)=1$,
	i.e.\ $e^{L(z_\ast)}=z_\ast$. Then $F\equiv 1$ and $e^{L(z)}=z$ for all $z\in\Omega$.
	
	Finally, if $L_1,L_2$ are two branches, then $e^{L_1-L_2}=1$, so $L_1-L_2$ takes values in $2\pi i\Z$.
	By connectedness it must be constant, hence equals $2\pi i k$ for some $k\in\Z$.
\end{proof}

\begin{definition}[Principal branch on the slit plane]
	\label{def:principal-Log}
	Let $\Omega_0=\C\setminus(-\infty,0]$.
	Define the \emph{principal argument} $\Arg z\in(-\pi,\pi)$ and the \emph{principal logarithm}
	\[
	\Log z:=\ln|z|+i\,\Arg z \qquad (z\in\Omega_0).
	\]
\end{definition}

\begin{proposition}[Holomorphy and derivative of the principal branch]
	\label{prop:principal-Log-holo}
	The function $\Log$ is holomorphic on $\Omega_0$ and satisfies
	\[
	(\Log z)'=\frac{1}{z}.
	\]
\end{proposition}

\begin{proof}
	This is a special case of Theorem~\ref{thm:log-exists} since $\Omega_0$ is simply connected and avoids $0$.
	The normalization $\Arg z\in(-\pi,\pi)$ pins down the branch uniquely.
\end{proof}

\begin{example}[No holomorphic logarithm on an annulus]
	Let $\Omega=\{z\in\C:1<|z|<2\}$. Assume for contradiction that $L\in\Hol(\Omega)$ satisfies $e^{L(z)}=z$.
	Differentiating gives $L'(z)=1/z$ on $\Omega$.
	Let $\gamma(t)=\frac{3}{2}e^{it}$, $0\le t\le 2\pi$. Then $\gamma$ is a closed loop in $\Omega$ and
	\[
	\int_{\gamma}L'(z)\,dz=\int_{\gamma}\frac{1}{z}\,dz.
	\]
	The left-hand side equals $0$ because $L$ is single-valued:
	\[
	\int_{\gamma}L'(z)\,dz = L(\gamma(2\pi))-L(\gamma(0))=0.
	\]
	But $\int_{\gamma}\frac{1}{z}\,dz=2\pi i$, a contradiction. Hence no holomorphic logarithm exists on $\Omega$.
\end{example}

\begin{exercise}[Rotated branch cuts]
	Fix $\alpha\in\R$ and set
	\[
	\Omega_\alpha:=\C\setminus\{-re^{i\alpha}: r\ge 0\}.
	\]
	\begin{enumerate}
		\item Show $\Omega_\alpha$ is simply connected.
		\item Construct a branch $\Log_\alpha$ on $\Omega_\alpha$ and prove $(\Log_\alpha)'(z)=1/z$.
		\item Compare $\Log_\alpha$ and $\Log$ on $\Omega_\alpha\cap\Omega_0$.
	\end{enumerate}
\end{exercise}

\begin{exercise}[Winding and $2\pi i$ jumps]
	Let $\gamma:[0,1]\to\C^{\*}$ be a closed piecewise $C^1$ curve.
	Show that
	\[
	\frac{1}{2\pi i}\int_\gamma \frac{1}{z}\,dz \in \Z.
	\]
	(Preview: this integer will later be identified with the winding number of $\gamma$ about $0$.)
\end{exercise}

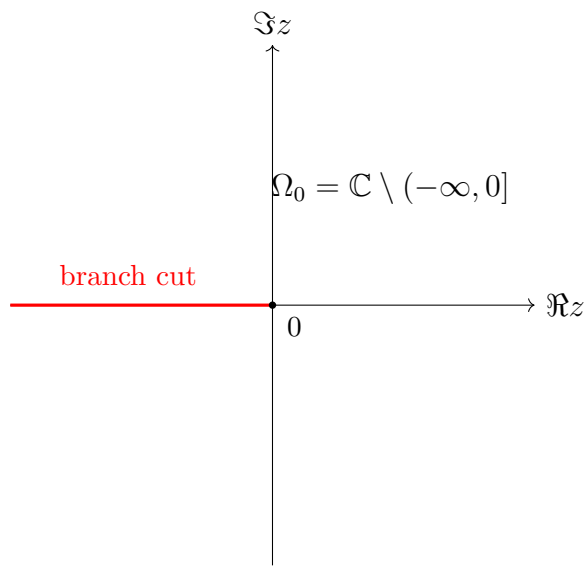
\begin{figure}[ht]
	\centering
	\begin{tikzpicture}[scale=1.15]
		\draw[->] (-3,0)--(3,0) node[right]{$\Re z$};
		\draw[->] (0,-3)--(0,3) node[above]{$\Im z$};
		\draw[very thick,red] (-3,0)--(0,0);
		\node[red] at (-1.65,0.35){\small branch cut};
		\node at (1.65,1.35){$\Omega_0=\C\setminus(-\infty,0]$};
		\fill (0,0) circle (1.2pt);
		\node at (0.25,-0.25){\small $0$};
	\end{tikzpicture}
	\caption{The principal logarithm is defined after removing the nonpositive
	real axis; crossing the cut changes the argument by $2\pi$.}
	\label{fig:principal-log-cut}
\end{figure}

\noindent\textbf{Summary.}
Identity theorem: zeros do not accumulate in the interior unless $f\equiv 0$.
Holomorphic logarithm: $\log$ exists precisely when the domain has no winding around $0$
(equivalently, every closed curve in the domain has winding number zero about
$0$). Simple connectivity is sufficient, but it is not necessary: a domain may
have holes that do not wind around the origin.

\subsection{Multivalued Analytic Functions and Riemann Surfaces}
\label{subsec:multivalued-riemann}

In this subsection we revisit \emph{multivalued} analytic functions such as
$\log z$ and $z^\alpha$, and explain a standard resolution principle:

\begin{center}
	\emph{Multivalued on the $z$--plane \quad$\leadsto$\quad single-valued on a suitable Riemann surface.}
\end{center}

Rather than forcing the function to live on the flat $z$--plane, we build a new
surface by cutting and gluing \emph{sheets}.  On that surface the function becomes
an honest single-valued holomorphic function.

\begin{center}
	\emph{Cut the plane $\longrightarrow$ glue sheets along the cut $\longrightarrow$ obtain a new surface.}
\end{center}

\begin{definition}[Complex manifold and Riemann surface]
	\label{def:complex-manifold}
	\index{complex manifold!definition}\index{Riemann surface!definition}
	A \emph{complex manifold of dimension $n$} is a Hausdorff,
	second-countable space $X$ covered by open sets $U_\alpha$ with
	homeomorphisms $z_\alpha:U_\alpha\to V_\alpha\subset\C^n$ onto
	open subsets such that every change of coordinates
	\[
	z_\beta\circ z_\alpha^{-1}:
	z_\alpha(U_\alpha\cap U_\beta)
	\longrightarrow z_\beta(U_\alpha\cap U_\beta)
	\]
	is holomorphic. Here a map between open subsets of complex Euclidean
	spaces is holomorphic when its component functions have convergent
	local power series in the complex coordinates. The charts form a
	\emph{holomorphic atlas}; their transition maps are biholomorphic
	because the reverse changes of coordinates belong to the same atlas.
	A connected complex manifold of dimension one is a \emph{Riemann
	surface}. Unless specified otherwise, manifolds in this book have no
	boundary.
\end{definition}

A function or map on complex manifolds is holomorphic when its expression
in these charts is holomorphic. This condition does not depend on the
charts: changing coordinates composes with biholomorphic maps. For
example, the Riemann sphere $\PP^1$ has the two coordinates $z$ and
$w=1/z$; the transition $w=1/z$ is holomorphic on $\C^*$, although
neither coordinate covers the whole sphere. Later, quotients of $\C$
by lattices give further Riemann surfaces whose chart changes are
translations.

\subsubsection{From branches to monodromy}

The preceding subsection constructed logarithm branches on suitable planar
domains.  We now ask what happens when a branch is analytically continued along
a path.  Continuing once around the origin increases a chosen argument by
$2\pi$ and hence changes a logarithm value by $2\pi i$.  More generally, the
change along a closed path $\gamma\subset\Cstar$ is
\[
	2\pi i\,\operatorname{Ind}(\gamma,0).
\]
This dependence of the continued value on the homotopy class of the path is the
basic example of \emph{monodromy}.

\begin{definition}[Holomorphic branch]
	Let $F$ be a multivalued analytic function on a domain $\Omega\subset\C$.
	A \emph{holomorphic branch} of $F$ on $\Omega$ is a holomorphic function
	$f\in\Hol(\Omega)$ such that $f(z)\in F(z)$ for every $z\in\Omega$.
	In particular, a holomorphic branch of $\log z$ on $\Omega\subset\Cstar$ is a holomorphic
	function $L\in\Hol(\Omega)$ satisfying $e^{L(z)}=z$ for all $z\in\Omega$.
\end{definition}

The definition of a holomorphic branch below applies equally to logarithms,
roots, and the algebraic functions considered later.  For logarithms, the
existence theorem and the principal branch were established in
Theorem~\ref{thm:log-exists} and Proposition~\ref{prop:principal-Log-holo}.

\begin{example}[Walking around the origin]
	Let $\gamma(t)=re^{it}$ for $0\le t\le 2\pi$. Analytic continuation of $\log z$ along $\gamma$
	returns to the same point $z=r$ but changes the value by $2\pi i$.
	This is the obstruction to having a global single-valued logarithm on $\Cstar$.
\end{example}

\subsubsection{Complex powers $z^\alpha$ and (multi)valuedness}

Fix $\alpha\in\C$.  Using the multivalued logarithm, define the multivalued power
\[
z^\alpha := e^{\alpha\log z}.
\]
Concretely, if $z=re^{i\theta}$ then the set of values is
\[
z^\alpha
= \bigl\{\, r^\alpha\, e^{\,i\alpha(\theta+2\pi k)} \ \bigm|\ k\in\Z\,\bigr\}.
\]

\begin{theorem}[Single-, finite-, and infinite-valuedness of $z^\alpha$]
	\label{thm:z-alpha-values}
	Let $\alpha\in\C$.
	\begin{enumerate}[label=\textup{(\arabic*)}]
		\item If $\alpha\in\Z$, then $z^\alpha$ is single-valued on $\Cstar$
		(and extends holomorphically across $0$ if $\alpha\ge 0$).
		\item If $\alpha\in\Q\setminus\Z$ and $\alpha=\frac{m}{n}$ with $\gcd(m,n)=1$ and $n\ge2$,
		then $z^\alpha$ has exactly $n$ distinct values at each $z\neq 0$.
		\item If $\alpha\notin\Q$, then $z^\alpha$ is infinitely multivalued:
		different $k\in\Z$ give infinitely many distinct values.
	\end{enumerate}
\end{theorem}

\begin{proof}
	Write $z=re^{i\theta}$.
	
	\noindent\textup{(1)}
	If $\alpha=n\in\Z$, then
	$e^{i n(\theta+2\pi k)}=e^{in\theta}$ for all $k\in\Z$, so the value is independent of $k$.
	
	\noindent\textup{(2)}
	If $\alpha=\frac{m}{n}$ with $\gcd(m,n)=1$, then increasing $k$ by $n$ adds
	$2\pi i m$ to the exponent, hence does not change the value.
	Thus there are at most $n$ values.
	They are distinct for $k=0,1,\dots,n-1$ since $\gcd(m,n)=1$ forces the angles
	$\frac{m}{n}(\theta+2\pi k)$ to be distinct modulo $2\pi$.
	
	\noindent\textup{(3)}
	For two indices $k,\ell$, the quotient of the corresponding values is
	$e^{2\pi i\alpha(k-\ell)}$. Equality for $k\ne\ell$ would imply
	$\alpha(k-\ell)\in\Z$, hence $\alpha\in\Q$. Thus if
	$\alpha\notin\Q$ (including nonreal $\alpha$), all values are distinct.
\end{proof}

\begin{example}[Rational powers]
	\begin{enumerate}
		\item $\alpha=\tfrac12$: $z^{1/2}$ has two values (the two square roots).
		\item $\alpha=\tfrac23$: $z^{2/3}$ has three values.
	\end{enumerate}
\end{example}

\subsubsection{Geometric picture: the Riemann surface of $z^{1/2}$}

The function $z^{1/2}$ cannot be single-valued on $\Cstar$:
going once around $0$ changes the sign of the square root.
We resolve this by passing to a two-sheeted surface built from two copies of a cut plane.

\medskip
\noindent\textbf{Cut-and-glue model.}
Cut the $z$--plane along the ray $[0,\infty)$ and take two copies of
\[
\C\setminus[0,\infty).
\]
On each sheet we have a continuous argument $\vartheta$.
Glue the sheets so that crossing the cut moves you to the other sheet.
The argument can be continued locally across the glued edges; after a full turn
on the two-sheeted surface it changes by $4\pi$. Thus there is no global
real-valued argument on the punctured surface, but the expression
\[
\sqrt{z}=\sqrt{r}\,e^{i\vartheta/2}
\]
is unchanged when $\vartheta$ is replaced by $\vartheta+4\pi$.
It therefore defines a single-valued holomorphic function on the glued
surface, extending across the point over $z=0$; in the algebraic model below
this function is simply the coordinate $w$.

\medskip
\noindent\textbf{Algebraic model.}
Equivalently, consider the affine curve
\[
X := \{(z,w)\in\C^2 : w^2=z\}.
\]
The projection $\pi:X\to\C$, $\pi(z,w)=z$, is $2$-to-$1$ away from $z=0$.
On $X$ the function $w$ is a single-valued holomorphic square root of $z$.

\begin{proposition}[Two-sheeted square-root surface]
	\label{prop:sqrt-surface}
	On $X=\{(z,w)\in\C^2: w^2=z\}$, the function $w$ is holomorphic and satisfies $w^2=z$.
	The projection $\pi:X\to\C$ is a $2$-to-$1$ holomorphic map away from $z=0$,
	and it is branched at $z=0$ (and also at $z=\infty$ after compactification to $\CP^1$).
\end{proposition}

\begin{proof}
	The map $\C\to X$, $w\mapsto(w^2,w)$, is a biholomorphic
	parameterization with inverse $(z,w)\mapsto w$. In this coordinate the
	function called $w$ is the identity and the projection is
	$w\mapsto w^2$. Its derivative is $2w$, so it is locally biholomorphic
	and two-to-one away from $w=0$, while it has ramification index $2$ at
	$w=0$. In the coordinate $u=1/w$ near the added point at infinity, the
	target coordinate is $1/z=u^2$, which gives the second branch point.
\end{proof}

\begin{figure}[ht]
	\centering
	\begin{tikzpicture}[scale=1.2,>=Stealth]
		\draw[rounded corners] (-4,-1.15) rectangle (-0.7,1.15);
		\draw[rounded corners] (0.7,-1.15) rectangle (4,1.15);
		\node at (-2.35,0.75) {sheet $+$};
		\node at (2.35,0.75) {sheet $-$};
		\draw[thick] (-3.65,0.08)--(-1.0,0.08);
		\draw[thick] (-3.65,-0.08)--(-1.0,-0.08);
		\draw[thick] (1.0,0.08)--(3.65,0.08);
		\draw[thick] (1.0,-0.08)--(3.65,-0.08);
		\node[font=\scriptsize] at (-2.5,0.30) {upper bank};
		\node[font=\scriptsize] at (-2.5,-0.33) {lower bank};
		\node[font=\scriptsize] at (2.5,0.30) {upper bank};
		\node[font=\scriptsize] at (2.5,-0.33) {lower bank};
		\draw[->,red!70!black,thick] (-1.0,0.08)--(1.0,-0.08);
		\draw[->,blue!70!black,thick] (-1.0,-0.08)--(1.0,0.08);
		\node[font=\small] at (0,-0.70) {crosswise identification};
	\end{tikzpicture}
	\caption{The upper bank on each slit sheet is identified with the lower
	bank on the other. This crosswise gluing makes the square root single-valued
	on the resulting surface; the algebraic model $w^2=z$ supplies its global coordinate $w$.}
	\label{fig:square-root-gluing}
\end{figure}

\subsubsection{Three sheets: the Riemann surface of $z^{1/3}$}

For the cube root we obtain a three-sheeted surface. Consider
\[
\mathcal{C}:=\{(z,w)\in\C^2 : w^3=z\}.
\]
Then $\pi(z,w)=z$ is $3$-to-$1$ away from $z=0$ (and, after compactification, away from $\infty$).
Writing $z=re^{i\theta}$, the three values
\[
w=r^{1/3}e^{i(\theta+2\pi k)/3},\qquad k=0,1,2,
\]
are realized as three sheets. Cutting along a ray and gluing three copies cyclically
produces the same surface as the algebraic model $\mathcal{C}$.

\begin{exercise}[Three-sheeted gluing]
	Construct the Riemann surface of $z^{1/3}$ by:
	\begin{enumerate}[label=\textup{(\arabic*)}]
		\item taking three copies of $\C\setminus[0,\infty)$,
		\item gluing the edges of the cut cyclically (sheet $1\to 2\to 3\to 1$),
		\item adding the points over $0$ and $\infty$ and explaining why the resulting compact surface is simply connected,
		and why the projection to the $z$--sphere is branched at $0$ and $\infty$.
	\end{enumerate}
\end{exercise}

\begin{exercise}[$n$ sheets for $\alpha=\frac{m}{n}$]
	Let $\alpha=\frac{m}{n}\in\Q$ with $\gcd(m,n)=1$ and $n\ge2$.
	Construct an $n$-sheeted surface by gluing $n$ copies of a slit plane so that
	a local argument returns to its starting value modulo $2\pi n$ after $n$ turns.
	Add the points over $0$ and $\infty$. Explain why the projection to $\CP^1$
	is an $n$-sheeted branched covering with branch points at $0$ and $\infty$.
	As a check, parametrize the compact surface by $z=t^n$ and
	$z^{m/n}=t^m$; interpret $t^m$ meromorphically when $m<0$.
\end{exercise}

\begin{exercise}[Where can $z^\alpha$ be single-valued?]
	Let $\alpha\in\C$. Determine the domains $\Omega\subset\Cstar$ on which one can choose a
	single-valued holomorphic branch of $z^\alpha$. Phrase your answer in terms of winding numbers of
	closed curves in $\Omega$ around $0$. Hint: continuation around a loop of
	winding number $k$ multiplies a local branch by $e^{2\pi i\alpha k}$.
\end{exercise}

\subsection{Meromorphic Functions, Singularities, and the Residue}
\label{subsec:mero-sing-res}

This subsection introduces \emph{meromorphic} functions and explains how
their behavior near isolated points is encoded by the Laurent expansion.
Conceptually, this is about \emph{classifying the local geometry of a function
	near a puncture}:
\begin{itemize}
	\item removable singularities behave like ordinary points after filling in a hole,
	\item poles look like ``finite order blow-ups'',
	\item essential singularities exhibit chaotic behavior (Casorati--Weierstrass).
\end{itemize}
The residue will appear as the coefficient of the $(z-p)^{-1}$ term and will
serve as a basic counting tool in residue theory.

\subsubsection{Meromorphic functions and Laurent expansions}

\begin{definition}[Meromorphic function]
	Let $\Omega\subset\C$ be open.
	A function $f:\Omega\to\widehat{\C}$ is \emph{meromorphic on $\Omega$} if there exists a
	discrete set $S\subset\Omega$ such that
	\begin{enumerate}[label=\textup{(\arabic*)}]
		\item $f$ is holomorphic on $\Omega\setminus S$, and
		\item for every $p\in S$, the point $p$ is a \emph{pole} of $f$.
	\end{enumerate}
	Equivalently, $f$ is holomorphic on $\Omega\setminus S$ and has only isolated poles
	(no essential singularities) in $\Omega$.
\end{definition}

\begin{remark}
	Meromorphic functions are to holomorphic functions what rational functions are to
	polynomials: ``holomorphic except for controlled poles''.
\end{remark}

\begin{theorem}[Laurent expansion on a punctured disk]
	\label{thm:laurent}
	Let $p\in\Omega$, and suppose $f$ is holomorphic on $0<|z-p|<r$ for some $r>0$.
	Then there exist unique coefficients $\{a_k\}_{k\in\Z}$ such that
	\[
	f(z)=\sum_{k=-\infty}^{\infty} a_k (z-p)^k,
	\]
	and the series converges absolutely and uniformly on every annulus
	$\rho\le |z-p|\le R$ with $0<\rho<R<r$.
\end{theorem}

\begin{proof}
	Fix numbers $0<\rho<R<r$. Choose two positively oriented circles
	\[
	C_R:\ \zeta=p+Re^{i\theta},\qquad C_\rho:\ \zeta=p+\rho e^{i\theta},
	\qquad 0\le \theta\le 2\pi.
	\]
	Let $A_{\rho,R}:=\{z:\rho<|z-p|<R\}$ be the annulus. Since $f$ is holomorphic on
	$A_{\rho,R}$ and on a neighborhood of its closure, Cauchy's integral formula on
	the annulus (equivalently, applying Cauchy on the boundary $\partial A_{\rho,R}=C_R- C_\rho$)
	gives, for every $z\in A_{\rho,R}$,
	\[
	f(z)
	= \frac{1}{2\pi i}\int_{C_R}\frac{f(\zeta)}{\zeta-z}\,d\zeta
	-\frac{1}{2\pi i}\int_{C_\rho}\frac{f(\zeta)}{\zeta-z}\,d\zeta.
	\]
	We expand each kernel as a geometric series.
	
	\smallskip
	\noindent\textbf{(Outer circle contribution).}
	On $C_R$ we have $|\zeta-p|=R$ and for $z\in A_{\rho,R}$ we have $|z-p|<R$, hence
	\[
	\frac{1}{\zeta-z}
	=\frac{1}{(\zeta-p)-(z-p)}
	=\frac{1}{\zeta-p}\cdot \frac{1}{1-\frac{z-p}{\zeta-p}}
	=\frac{1}{\zeta-p}\sum_{n=0}^\infty \Bigl(\frac{z-p}{\zeta-p}\Bigr)^n,
	\]
	where the series converges absolutely and uniformly because
	$\bigl|\frac{z-p}{\zeta-p}\bigr|\le \frac{R'}{R}<1$ on $\{ |z-p|\le R'\}$ for any fixed
	$R'<R$.
	
	Thus
	\[
	\frac{1}{2\pi i}\int_{C_R}\frac{f(\zeta)}{\zeta-z}\,d\zeta
	=\sum_{n=0}^\infty
	\left(\frac{1}{2\pi i}\int_{C_R} \frac{f(\zeta)}{(\zeta-p)^{n+1}}\,d\zeta\right)(z-p)^n.
	\]
	Define for $n\ge 0$
	\[
	a_n:=\frac{1}{2\pi i}\int_{C_R} \frac{f(\zeta)}{(\zeta-p)^{n+1}}\,d\zeta.
	\]
	These do not depend on the chosen $R$ as long as $\rho<R<r$, by deforming contours inside
	the holomorphic region (Cauchy's theorem).
	
	\smallskip
	\noindent\textbf{(Inner circle contribution).}
	On $C_\rho$ we have $|\zeta-p|=\rho$ and for $z\in A_{\rho,R}$ we have $|z-p|>\rho$, hence
	\[
	\frac{1}{\zeta-z}
	=-\frac{1}{z-\zeta}
	=-\frac{1}{(z-p)-(\zeta-p)}
	=-\frac{1}{z-p}\cdot\frac{1}{1-\frac{\zeta-p}{z-p}}
	=-\frac{1}{z-p}\sum_{n=0}^\infty \Bigl(\frac{\zeta-p}{z-p}\Bigr)^n,
	\]
	with absolute and uniform convergence on $\{|z-p|\ge \rho'\}$ for any $\rho'>\rho$.
	
	Therefore
	\[
	-\frac{1}{2\pi i}\int_{C_\rho}\frac{f(\zeta)}{\zeta-z}\,d\zeta
	=
	\sum_{n=0}^\infty
	\left(\frac{1}{2\pi i}\int_{C_\rho} f(\zeta)(\zeta-p)^n\,d\zeta\right)(z-p)^{-(n+1)}.
	\]
	Define for $m\ge 1$ (set $m=n+1$)
	\[
	a_{-m}:=\frac{1}{2\pi i}\int_{C_\rho} f(\zeta)(\zeta-p)^{m-1}\,d\zeta.
	\]
	Again $a_{-m}$ is independent of $\rho$ by contour deformation.
	
	\smallskip
	\noindent\textbf{(Assembling).}
	Combining the two expansions yields, for every $z\in A_{\rho,R}$,
	\[
	f(z)=\sum_{n=0}^\infty a_n(z-p)^n + \sum_{m=1}^\infty a_{-m}(z-p)^{-m}
	=\sum_{k=-\infty}^{\infty} a_k (z-p)^k.
	\]
	Uniform convergence on compact subannuli follows from uniform convergence of the geometric
	series and boundedness of $f$ on each circle. For a prescribed closed annulus
	$\rho\le |z-p|\le R$ with $0<\rho<R<r$, choose the two contour radii
	strictly below $\rho$ and strictly above $R$; the same argument then gives
	uniform convergence on the prescribed closed annulus.
	
	\smallskip
	\noindent\textbf{Uniqueness.}
	If two Laurent series represent $f$ on the annulus, their difference
	$\sum_{k\in\Z}c_k(z-p)^k$ represents $0$.  For a fixed integer $n$,
	multiply by $(z-p)^{-n-1}$ and integrate term-by-term on a circle.  Since
	\[
	\frac{1}{2\pi i}\int_{|\zeta-p|=\rho} (\zeta-p)^k\,d\zeta
	=
	\begin{cases}
		1,& k=-1,\\
		0,& k\neq -1,
	\end{cases}
	\]
	the only surviving term after division by $2\pi i$ is $c_n$. Thus $c_n=0$ for every
	$n\in\Z$, and the coefficients are unique.
\end{proof}

\begin{definition}[Residue]
	\label{def:residue}
	Under the assumptions of Theorem~\ref{thm:laurent}, the coefficient $a_{-1}$
	is called the \emph{residue} of $f$ at $p$ and is denoted by
	\[
	\operatorname{Res}(f,p):=a_{-1}.
	\]
\end{definition}

\begin{theorem}[Cauchy integral formula for residues]
	\label{thm:cauchy-residue}
	Let $p\in\Omega$ and let $\gamma$ be a positively oriented piecewise $C^1$ simple closed curve
	contained in a punctured disk $0<|z-p|<r$ and winding once around $p$.
	If $f$ is holomorphic on and inside $\gamma$ except possibly at $p$, then
	\[
	\frac{1}{2\pi i}\int_\gamma f(z)\,dz = \operatorname{Res}(f,p).
	\]
\end{theorem}

\begin{proof}
	Choose $\rho>0$ so small that the closed disk
	$\overline{B(p,\rho)}$ lies inside $\gamma$ and contains no singularity
	other than the possible one at $p$. Apply Cauchy's theorem to $f$ on the
	region between $\gamma$ and the positively oriented circle
	$C_\rho=\{|z-p|=\rho\}$. Since the inner boundary carries the negative
	boundary orientation,
	\[
	\int_\gamma f(z)\,dz=\int_{C_\rho}f(z)\,dz.
	\]
	The Laurent series at $p$ converges uniformly on the circle $C_\rho$, so
	it may be integrated term-by-term there. Every term except the
	coefficient of $(z-p)^{-1}$ integrates to zero, while
	\[
	\int_{|z-p|=\rho}\frac{1}{z-p}\,dz
	=
	\int_0^{2\pi}\frac{1}{\rho e^{it}}\, i\rho e^{it}\,dt
	= i\int_0^{2\pi} dt
	=2\pi i.
	\]
	Therefore
	\[
	\int_\gamma f(z)\,dz = a_{-1}\cdot 2\pi i,
	\]
	and dividing by $2\pi i$ yields the claim.
\end{proof}

\subsubsection{Classification of isolated singularities}

\begin{definition}[Isolated singularities]
	Let $f$ be holomorphic on $0<|z-p|<r$.
	\begin{itemize}
		\item $p$ is a \emph{removable singularity} if $f$ extends holomorphically to $p$.
		\item $p$ is a \emph{pole} if $|f(z)|\to\infty$ as $z\to p$.
		\item $p$ is an \emph{essential singularity} if it is neither removable nor a pole.
	\end{itemize}
	A pole has \emph{order $m$} when $(z-p)^m f(z)$ extends holomorphically
	to $p$ with a nonzero value there. The growth criterion below shows
	that every pole in the preceding sense has such a finite order.
\end{definition}

\begin{proposition}[Laurent-series criteria]
	\label{prop:sing-class-laurent}
	Let
	\[
	f(z)=\sum_{k=-\infty}^{\infty} a_k (z-p)^k
	\]
	be the Laurent expansion on $0<|z-p|<r$.
	\begin{enumerate}[label=\textup{(\arabic*)}]
		\item $p$ is removable $\iff a_k=0$ for all $k<0$.
		\item $p$ is a pole of order $m\ge1$ $\iff a_{-m}\neq 0$ and $a_k=0$ for all $k<-m$.
		\item $p$ is essential $\iff a_k\neq 0$ for infinitely many $k<0$.
	\end{enumerate}
\end{proposition}

\begin{proof}
	\textup{(1)} If $a_k=0$ for all $k<0$, then $f(z)=\sum_{k\ge 0}a_k(z-p)^k$ on $0<|z-p|<r$,
	and the series converges at $z=p$ as well, giving a holomorphic extension.
	
	Conversely, if $f$ extends holomorphically to $p$, then $f$ is holomorphic on $|z-p|<r$ and
	has a Taylor series there, so all negative Laurent coefficients must be $0$ by uniqueness.
	
	\textup{(2)} If $a_k=0$ for all $k<-m$ and $a_{-m}\neq 0$, define
	\[
	g(z):=(z-p)^m f(z)=\sum_{k=-m}^{\infty} a_k (z-p)^{k+m}
	=\sum_{j=0}^{\infty} a_{j-m}(z-p)^j.
	\]
	Then $g$ extends holomorphically to $p$ with $g(p)=a_{-m}\neq 0$, hence is nonzero near $p$.
	Thus $f(z)=(z-p)^{-m}g(z)$ blows up like $|z-p|^{-m}$, so $p$ is a pole of order $m$.
	
	Conversely, if $p$ is a pole of order $m$, then $(z-p)^m f(z)$ extends holomorphically and
	nonvanishingly to $p$, forcing the Laurent series of $f$ to have precisely $m$ negative powers.
	
	\textup{(3)} If infinitely many negative powers occur, then for every $m$ the product
	$(z-p)^m f(z)$ still has a negative term in its Laurent series, hence cannot extend holomorphically
	to $p$. In particular $f$ is not removable, and also not a pole of finite order. Thus $p$ is essential.
	
	Conversely, if $p$ is essential, it is neither removable nor a pole, so the principal part cannot be
	finite and cannot be empty; hence there must be infinitely many negative terms.
\end{proof}

\begin{remark}[Order notation]
	If $p$ is a pole of order $m$, write $\ord_p(f)=-m$.
	If $f$ has a zero of order $m\ge 1$ at $p$, write $\ord_p(f)=m$.
\end{remark}

\subsubsection{Examples}

\begin{example}[Removable singularity]
	Consider $f(z)=\dfrac{\sin z}{z}$ on $\Cstar$.
	Since
	\[
	\sin z = z - \frac{z^3}{6} + \frac{z^5}{120} - \cdots,
	\]
	we obtain
	\[
	\frac{\sin z}{z}=1-\frac{z^2}{6}+\frac{z^4}{120}-\cdots,
	\]
	which converges at $z=0$. Defining $f(0)=1$ fills in the hole, so $z=0$ is removable.
\end{example}

\begin{example}[Pole and residue]
	Let $f(z)=\dfrac{1}{(z-p)^3}$. Then $p$ is a pole of order $3$ and the residue is $0$.
	Indeed the Laurent series is exactly $(z-p)^{-3}$, with no $(z-p)^{-1}$ term.
\end{example}

\begin{example}[Essential singularity]
	The function $f(z)=e^{1/z}$ has an essential singularity at $z=0$, since
	\[
	e^{1/z}=\sum_{n=0}^{\infty}\frac{1}{n!}\,z^{-n}
	\]
	has infinitely many negative powers.
\end{example}

\subsubsection{Removable singularities and extension}

\begin{theorem}[Riemann removable singularity theorem]
	\label{thm:removable}
	Let $f$ be holomorphic on $0<|z-p|<r$ and bounded near $p$.
	Then $f$ extends uniquely to a holomorphic function on $|z-p|<r$.
\end{theorem}

\begin{proof}
	Assume $|f(z)|\le M$ for all $0<|z-p|<r_0$ for some $r_0\le r$.
	Fix a radius $\rho$ with $0<\rho<r_0$. By Theorem~\ref{thm:laurent}, $f$ has a Laurent
	expansion on $0<|z-p|<\rho$:
	\[
	f(z)=\sum_{k=-\infty}^{\infty} a_k (z-p)^k,
	\qquad
	a_k=\frac{1}{2\pi i}\int_{|\zeta-p|=\rho}\frac{f(\zeta)}{(\zeta-p)^{k+1}}\,d\zeta.
	\]
	We show $a_k=0$ for all $k<0$. Let $k=-m$ with $m\ge 1$. Then
	\[
	a_{-m}
	=\frac{1}{2\pi i}\int_{|\zeta-p|=\rho} f(\zeta)\,(\zeta-p)^{m-1}\,d\zeta.
	\]
	Estimate its absolute value:
	\[
	|a_{-m}|
	\le \frac{1}{2\pi}\int_{|\zeta-p|=\rho} |f(\zeta)|\,|\zeta-p|^{m-1}\,|d\zeta|
	\le \frac{1}{2\pi}\int_{|\zeta-p|=\rho} M\,\rho^{m-1}\,|d\zeta|
	= M\rho^{m}.
	\]
	Since this holds for every $0<\rho<r_0$, letting $\rho\to 0$ forces $|a_{-m}|\le M\rho^{m}\to 0$,
	hence $a_{-m}=0$ for each $m\ge 1$. Therefore all negative coefficients vanish and
	\[
	f(z)=\sum_{k=0}^{\infty} a_k (z-p)^k
	\]
	on $0<|z-p|<r_0$. The right-hand side converges at $z=p$, so defining $f(p):=a_0$ extends
	$f$ holomorphically to $|z-p|<r_0$ (and hence to $|z-p|<r$ by analyticity on overlaps).
	
	Uniqueness: if two holomorphic extensions existed, their difference would vanish on
	$0<|z-p|<r_0$, hence everywhere by the Identity Theorem.
\end{proof}

\subsubsection{Essential singularities and Casorati--Weierstrass}

\begin{theorem}[Casorati--Weierstrass]
	\label{thm:casorati}
	Suppose $f$ is holomorphic on $0<|z-p|<r_0$ and has an essential
	singularity at $p$. Then for every $r$ with $0<r<r_0$ the image
	\[
	f\bigl(\{z:0<|z-p|<r\}\bigr)
	\]
	is dense in $\C$.
\end{theorem}

\begin{proof}
	Fix $r>0$ and suppose the image is not dense. Then there exist $a\in\C$ and $\varepsilon>0$
	such that
	\[
	f\bigl(\{z:0<|z-p|<r\}\bigr)\cap D(a,\varepsilon)=\emptyset,
	\]
	equivalently,
	\[
	|f(z)-a|\ge \varepsilon\qquad (0<|z-p|<r).
	\]
	Define
	\[
	g(z):=\frac{1}{f(z)-a}.
	\]
	Then $g$ is holomorphic on $0<|z-p|<r$ and satisfies $|g(z)|\le 1/\varepsilon$ there,
	so $g$ is bounded near $p$. By Theorem~\ref{thm:removable}, $g$ extends holomorphically
	to $|z-p|<r$; denote the extension again by $g$.
	
	Now consider the value $g(p)$.
	
	\smallskip
	\noindent\textbf{Case 1: $g(p)\neq 0$.}
	Then $1/g$ is holomorphic near $p$, and hence
	\[
	f(z)=a+\frac{1}{g(z)}
	\]
	extends holomorphically to $p$. This means $p$ is removable for $f$.
	
	\smallskip
	\noindent\textbf{Case 2: $g(p)=0$.}
	Then $g$ has a zero at $p$ of some order $m\ge 1$, so we can write
	\[
	g(z)=(z-p)^m h(z)
	\]
	with $h$ holomorphic and $h(p)\neq 0$. Hence $1/g$ has a pole of order $m$ at $p$, and thus
	\[
	f(z)=a+\frac{1}{g(z)}
	\]
	has a pole at $p$.
	
	\smallskip
	In both cases, $p$ is not an essential singularity of $f$, contradicting the hypothesis.
	Therefore the image must be dense in $\C$.
\end{proof}

\subsubsection{Exercises}

\begin{exercise}[Principal part and residue]
	Compute the principal part of
	\[
	f(z)=\frac{z^2+1}{(z-2)^3(z+1)}
	\]
	at $z=2$ and determine $\operatorname{Res}(f,2)$.
	Hint: write the partial fraction expansion near $z=2$ and read off the $(z-2)^{-1}$ term.
\end{exercise}

\begin{exercise}[Density near $z=0$]
	Show directly that for every $w\in\C$, $\varepsilon>0$, and $\delta>0$,
	there exists $z$ with $0<|z|<\delta$
	such that
	\[
	|e^{1/z}-w|<\varepsilon.
	\]
	For $w\ne0$, choose one logarithm $L$ of $w$ and use
	$z_n=(L+2\pi i n)^{-1}$. For $w=0$, take a negative real sequence
	tending to zero. Relate this behavior to the infinitely long principal part in
	the Laurent expansion of $e^{1/z}$.
\end{exercise}

\begin{exercise}[Meromorphic functions on the Riemann sphere]
	A meromorphic function on the Riemann sphere $\widehat{\C}\simeq \CP^1$ is the same thing
	as a rational function $P(z)/Q(z)$. Explain why:
	\begin{enumerate}[label=\textup{(\arabic*)}]
		\item a meromorphic function on $\widehat{\C}$ has only finitely many poles
		(including possibly a pole at $\infty$),
		\item using partial fractions (and the behavior at $\infty$), one can write it as a ratio
		of polynomials.
	\end{enumerate}
	This is a first step toward viewing meromorphic functions as morphisms of compact
	Riemann surfaces.
\end{exercise}

\subsection{Laurent Series and Cauchy Theory on Annuli}
\label{subsec:laurent-series-annulus}

The analysis of Laurent series on annuli requires the following tools:
\begin{itemize}
	\item the Weierstrass $M$-test for uniform convergence,
	\item uniform convergence of power series on disks and of principal parts on punctured disks,
	\item simultaneous convergence of both parts on an annulus,
	\item and the Cauchy theory (coefficient formulas and contour deformation) adapted to annular regions.
\end{itemize}
Algebraically, a Laurent series is the sum of
\emph{(i)} an ordinary power series near a finite center $p$ and
\emph{(ii)} a power series at $\infty$ in the variable $w=1/(z-p)$,
glued on a common annulus.

\subsubsection{Weierstrass $M$-test and uniform convergence}

\begin{theorem}[Weierstrass $M$-test]
	\label{thm:M-test}
	Let $\{f_k\}$ be complex-valued functions on a set $A$. Suppose there exist
	numbers $M_k\ge 0$ such that
	\[
	|f_k(z)|\le M_k\qquad\text{for all }z\in A,
	\]
	and $\sum_{k=0}^\infty M_k<\infty$. Then $\sum_{k=0}^\infty f_k(z)$ converges
	absolutely and uniformly on $A$.
\end{theorem}

\begin{proof}
	For $N>m$ and $z\in A$, the triangle inequality gives
	\[
	 \left|\sum_{k=m+1}^{N}f_k(z)\right|
	 \le \sum_{k=m+1}^{N}M_k
	 \le \sum_{k=m+1}^{\infty}M_k.
	\]
	The last tail tends to zero as $m\to\infty$, independently of $z$.
	Thus the partial sums are uniformly Cauchy and converge uniformly on $A$.
	Applying the same estimate to $\sum |f_k(z)|$ proves absolute convergence at
	every point.
\end{proof}

\begin{remark}[What uniform convergence permits]
	Uniform convergence permits term-by-term integration along a fixed rectifiable
	curve. For holomorphic summands, locally uniform convergence also makes the sum
	holomorphic, and Cauchy's formula then gives convergence of the derivative series
	on smaller compact subsets. Differentiation does not follow from pointwise or
	uniform convergence alone. These distinctions will be used in Laurent theory on
	annuli.
\end{remark}

\begin{lemma}[Uniform convergence on smaller disks]
	\label{lem:power-unif-smaller}
	Let $\sum_{k=0}^{\infty} a_k (z - p)^k$ be a power series that converges at some
	point $z_0\neq p$. Set $r_0:=|z_0-p|$. Then the series converges absolutely and
	uniformly on every closed disk $\overline{B(p,r)}$ with $0<r<r_0$; hence it
	converges uniformly on every compact subset of $B(p,r_0)$.
\end{lemma}

\begin{proof}
	Since the series converges at $z_0$, the terms $a_k(z_0-p)^k$ tend to $0$, hence
	are bounded: there exists $C>0$ such that
	\[
	|a_k(z_0-p)^k|\le C\qquad\text{for all }k\ge 0.
	\]
	Fix $r$ with $0<r<r_0$, and take $z$ with $|z-p|\le r$. Then
	\[
	|a_k(z-p)^k|
	=
	|a_k(z_0-p)^k|\left|\frac{z-p}{z_0-p}\right|^k
	\le
	C\left(\frac{r}{r_0}\right)^k.
	\]
	The geometric series $\sum_{k=0}^\infty C(r/r_0)^k$ converges. By the
	Weierstrass $M$-test (Theorem~\ref{thm:M-test}), $\sum_{k=0}^\infty a_k(z-p)^k$
	converges absolutely and uniformly on $\overline{B(p,r)}$. Since every compact
	$K\Subset B(p,r_0)$ is contained in some $\overline{B(p,r)}$ with $r<r_0$, the
	convergence is uniform on $K$.
\end{proof}

\subsubsection{Inner series (principal parts) as power series in $1/(z-p)$}

A principal part
\[
\sum_{k=1}^\infty a_{-k}(z-p)^{-k}
\]
is a genuine power series in the variable $w:=1/(z-p)$:
\[
\sum_{k=1}^\infty a_{-k}(z-p)^{-k}=\sum_{k=1}^\infty a_{-k} w^k.
\]
Thus it behaves like an ordinary power series near $w=0$, i.e.\ near
$z=\infty$ relative to the center $p$.

\begin{lemma}[Uniform convergence of principal parts outside a disk]
	\label{lem:principal-unif-outside}
	Suppose $\sum_{k=1}^\infty a_{-k}(z-p)^{-k}$ converges at some point $z_1\neq p$.
	Let $r_1:=|z_1-p|$. Then for every $r>r_1$, the series converges absolutely and
	uniformly on the exterior region $\{z:|z-p|\ge r\}$ (and hence uniformly on every
	compact subset of $\{z:|z-p|>r_1\}$).
\end{lemma}

\begin{proof}
	Convergence at $z_1$ implies the terms $a_{-k}(z_1-p)^{-k}$ are bounded:
	there exists $C>0$ such that
	\[
	|a_{-k}(z_1-p)^{-k}|\le C\qquad\text{for all }k\ge 1.
	\]
	Fix $r>r_1$ and take $z$ with $|z-p|\ge r$. Then
	\[
	|a_{-k}(z-p)^{-k}|
	=
	|a_{-k}(z_1-p)^{-k}|\left|\frac{z_1-p}{z-p}\right|^k
	\le
	C\left(\frac{r_1}{r}\right)^k.
	\]
	Again $\sum_{k=1}^\infty C(r_1/r)^k$ is geometric, so by the $M$-test the
	principal-part series converges absolutely and uniformly on $\{ |z-p|\ge r\}$.
	As before, any compact subset of $\{|z-p|>r_1\}$ lies in some $\{|z-p|\ge r\}$
	with $r>r_1$, giving uniform convergence on compacts.
\end{proof}

\subsubsection{Simultaneous convergence on an annulus}

A Laurent series about $p$ is the formal sum of an ``outer'' (nonnegative)
and an ``inner'' (negative) part:
\[
\sum_{k=-\infty}^\infty a_k(z-p)^k
=
\underbrace{\sum_{k=0}^\infty a_k(z-p)^k}_{\text{outer part}}
\;+\;
\underbrace{\sum_{k=1}^\infty a_{-k}(z-p)^{-k}}_{\text{inner part}}.
\]

\begin{lemma}[Convergence on an annulus]
	\label{lem:annulus-conv-refined}
	Assume the outer series $\sum_{k=0}^\infty a_k(z-p)^k$ converges at some point
	$z_0$ with $r_0:=|z_0-p|>0$, and the inner series $\sum_{k=1}^\infty a_{-k}(z-p)^{-k}$
	converges at some point $z_1$ with $r_1:=|z_1-p|>0$.
	Assume $r_1<r_0$. Then both series converge absolutely and uniformly on every compact subset of the annulus
	\[
	A_{r_1,r_0}(p):=\{z\in\C: r_1<|z-p|<r_0\}.
	\]
	In particular, the full Laurent series converges absolutely and uniformly on every
	compact subset of $A_{r_1,r_0}(p)$.
\end{lemma}

\begin{proof}
	By Lemma~\ref{lem:power-unif-smaller}, the outer series converges uniformly on every
	compact subset of $\{|z-p|<r_0\}$. By Lemma~\ref{lem:principal-unif-outside}, the inner
	series converges uniformly on every compact subset of $\{|z-p|>r_1\}$.
	Intersecting these regions yields uniform convergence of both parts on every compact
	subset of $A_{r_1,r_0}(p)$. Summing two uniformly convergent series gives uniform
	convergence of the full Laurent series on the same compact sets.
\end{proof}

\begin{corollary}[Uniform convergence on closed subannuli]
	\label{cor:closed-subannulus}
	Under the hypotheses and notation of Lemma~\ref{lem:annulus-conv-refined},
	fix $\varepsilon>0$ with $r_1+\varepsilon<r_0-\varepsilon$. Then on the closed subannulus
	\[
	\{z: r_1+\varepsilon \le |z-p|\le r_0-\varepsilon\},
	\]
	both the inner and outer series converge absolutely and uniformly, hence the Laurent
	series may be integrated term-by-term. It may also be differentiated
	term-by-term: choose a slightly larger closed subannulus still contained
	in $A_{r_1,r_0}(p)$ and apply the termwise differentiation theorem
	separately to the power series in $z-p$ and the power series in
	$w=1/(z-p)$, followed by the chain rule.
\end{corollary}

\subsubsection{Laurent coefficients: uniqueness and contour formulas}

\begin{proposition}[Coefficient formula and uniqueness]
	\label{prop:laurent-coeff-formula}
	Let $f$ be holomorphic on the annulus $A_{r_1,r_0}(p)$ and suppose
	\[
	f(z)=\sum_{k=-\infty}^{\infty} a_k(z-p)^k
	\]
	converges absolutely and uniformly on each compact subset of the annulus.
	Then for every integer $n$ and every circle $C_\rho:=\{|\zeta-p|=\rho\}$ with
	$r_1<\rho<r_0$,
	\[
	a_n=\frac{1}{2\pi i}\int_{C_\rho}\frac{f(\zeta)}{(\zeta-p)^{n+1}}\,d\zeta.
	\]
	Moreover, the right-hand side is independent of $\rho$, so the coefficients are uniquely
	determined by $f$. In particular, the Laurent expansion of $f$ about $p$ is unique.
\end{proposition}

\begin{proof}
	Fix $n$ and $\rho$. On $C_\rho$ the Laurent series converges uniformly, hence we may
	divide term-by-term:
	\[
	\frac{f(\zeta)}{(\zeta-p)^{n+1}}
	=
	\sum_{k=-\infty}^{\infty} a_k(\zeta-p)^{k-n-1},
	\qquad \zeta\in C_\rho.
	\]
	Uniform convergence on $C_\rho$ justifies term-by-term integration:
	\[
	\frac{1}{2\pi i}\int_{C_\rho}\frac{f(\zeta)}{(\zeta-p)^{n+1}}\,d\zeta
	=
	\sum_{k=-\infty}^{\infty} a_k\left(\frac{1}{2\pi i}\int_{C_\rho}(\zeta-p)^{k-n-1}\,d\zeta\right).
	\]
	Parameterize $C_\rho$ by $\zeta-p=\rho e^{it}$ ($0\le t\le 2\pi$). Then
	\[
	\int_{C_\rho}(\zeta-p)^m\,d\zeta
	=
	\int_0^{2\pi}(\rho e^{it})^m \, i\rho e^{it}\,dt
	=
	i\rho^{m+1}\int_0^{2\pi}e^{i(m+1)t}\,dt
	=
	\begin{cases}
		2\pi i, & m=-1,\\
		0, & m\neq -1.
	\end{cases}
	\]
	Hence the only surviving term is the one with $k-n-1=-1$, i.e.\ $k=n$, and the sum
	equals $a_n$.
	
	To see independence of $\rho$, fix $r_1<\rho_1<\rho_2<r_0$ and consider
	\[
	h(\zeta):=\frac{f(\zeta)}{(\zeta-p)^{n+1}},
	\]
	which is holomorphic on $A_{r_1,r_0}(p)$. By contour deformation on an annulus
	(Theorem~\ref{thm:annulus-deform} below), we have
	\[
	\int_{C_{\rho_2}} h(\zeta)\,d\zeta=\int_{C_{\rho_1}} h(\zeta)\,d\zeta,
	\]
	so the coefficient formula does not depend on the chosen radius.
\end{proof}

\subsubsection{Cauchy theory on annuli}

\begin{theorem}[Contour deformation on an annulus]
	\label{thm:annulus-deform}
	Let $f$ be holomorphic on the annulus $A_{r_1,r_0}(p)$.
	For any radii $r_1<\rho_1<\rho_2<r_0$ we have
	\[
	\int_{C_{\rho_2}} f(\zeta)\,d\zeta
	=
	\int_{C_{\rho_1}} f(\zeta)\,d\zeta,
	\]
	where each $C_\rho$ is positively oriented.
\end{theorem}

\begin{proof}
	Consider the closed annular region
	\[
	\overline{A_{\rho_1,\rho_2}}(p):=\{\zeta:\rho_1\le |\zeta-p|\le \rho_2\}.
	\]
	Its boundary is $\partial \overline{A_{\rho_1,\rho_2}}(p)=C_{\rho_2}-C_{\rho_1}$,
	since the inner boundary is negatively oriented.
	By Cauchy's integral theorem applied to $f$ on a neighborhood of this closed region,
	\[
	0=\int_{\partial \overline{A_{\rho_1,\rho_2}}(p)} f(\zeta)\,d\zeta
	=\int_{C_{\rho_2}} f(\zeta)\,d\zeta-\int_{C_{\rho_1}} f(\zeta)\,d\zeta,
	\]
	which gives the claim.
\end{proof}

\begin{theorem}[Cauchy integral formula on an annulus]
	\label{thm:cauchy-annulus}
	Let $f$ be holomorphic on $A_{r_1,r_0}(p)$ and choose radii
	$r_1<\rho_1<\rho_2<r_0$.
	Then for every $z$ with $\rho_1<|z-p|<\rho_2$,
	\[
	f(z)
	=
	\frac{1}{2\pi i}\left(
	\int_{C_{\rho_2}} \frac{f(\zeta)}{\zeta-z}\,d\zeta
	-
	\int_{C_{\rho_1}} \frac{f(\zeta)}{\zeta-z}\,d\zeta
	\right).
	\]
	Equivalently,
	\[
	\int_{C_{\rho_2}} \frac{f(\zeta)}{\zeta-z}\,d\zeta
	-
	\int_{C_{\rho_1}} \frac{f(\zeta)}{\zeta-z}\,d\zeta
	=
	2\pi i\, f(z).
	\]
\end{theorem}

\begin{proof}
	Fix $z$ with $\rho_1<|z-p|<\rho_2$. Choose $\varepsilon>0$ so small that
	$\overline{B(z,\varepsilon)}\subset \{\rho_1<|\zeta-p|<\rho_2\}$, and set
	\[
	R:=\{\zeta:\rho_1<|\zeta-p|<\rho_2\}\setminus \overline{B(z,\varepsilon)}.
	\]
	Define $g(\zeta):=f(\zeta)/(\zeta-z)$. Then $g$ is holomorphic on $R$.
	The boundary of $R$ has three components:
	\[
	\partial R = C_{\rho_2}-C_{\rho_1}-\partial B(z,\varepsilon),
	\]
	where $\partial B(z,\varepsilon)$ is positively oriented as a circle, and appears
	with a minus sign because it is an inner boundary of $R$.
	
	To justify Cauchy's theorem on this multiply connected region, cut $R$
	along two disjoint arcs joining its inner boundary circles to the outer
	boundary. Apply the theorem to the resulting simply connected region;
	the integrals along the two sides of each cut cancel. Thus
	\[
	0=\int_{\partial R} g(\zeta)\,d\zeta
	=
	\int_{C_{\rho_2}} \frac{f(\zeta)}{\zeta-z}\,d\zeta
	-\int_{C_{\rho_1}} \frac{f(\zeta)}{\zeta-z}\,d\zeta
	-\int_{\partial B(z,\varepsilon)} \frac{f(\zeta)}{\zeta-z}\,d\zeta.
	\]
	Hence
	\[
	\int_{C_{\rho_2}} \frac{f(\zeta)}{\zeta-z}\,d\zeta
	-\int_{C_{\rho_1}} \frac{f(\zeta)}{\zeta-z}\,d\zeta
	=
	\int_{\partial B(z,\varepsilon)} \frac{f(\zeta)}{\zeta-z}\,d\zeta.
	\]
	Since $f$ is holomorphic on and inside $\partial B(z,\varepsilon)$, the usual
	Cauchy integral formula on a disk gives
	\[
	\int_{\partial B(z,\varepsilon)} \frac{f(\zeta)}{\zeta-z}\,d\zeta
	=
	2\pi i\, f(z).
	\]
	Substituting yields the claimed identity, and dividing by $2\pi i$ gives the formula for $f(z)$.
\end{proof}

\begin{remark}[Two boundary components]
	On a disk there is only one boundary component, so Cauchy's formula uses one contour.
	On an annulus, the correct formula is ``outer minus inner''. This is exactly the
	analytic reflection of why Laurent series (with positive and negative powers) are the
	natural expansions on annular domains.
\end{remark}

\subsubsection{Visualization of the annulus}

\begin{figure}[H]
	\centering
	\begin{tikzpicture}[scale=1.1]
		\draw[thick] (0,0) circle (1.2);
		\draw[thick] (0,0) circle (2.4);
		\draw[-{Stealth[length=2.5mm]}] (2.4,0)
			arc[start angle=0,end angle=50,radius=2.4];
		\draw[-{Stealth[length=2.5mm]}] (1.2,0)
			arc[start angle=0,end angle=-50,radius=1.2];
		\node at (0,0) {$p$};
		\node at (1.2,0) [right] {$r_1$};
		\node at (2.4,0) [right] {$r_0$};
		\filldraw (1.8,0.8) circle (2pt) node[right] {$z$};
	\end{tikzpicture}
	\caption{The annular Cauchy formula uses the positively oriented outer
	circle and the negatively oriented inner boundary.}
	\label{fig:annular-cauchy}
\end{figure}

\subsection{Local Consequences of Laurent Theory}
\label{subsec:laurent-and-residue-no-bold}

We now put the Laurent expansion
\[
f(z)=\sum_{n=-\infty}^{\infty} a_n (z-p)^n
\]
to work near an isolated point $p$. The coefficient test and residue identity
were proved above; here we collect them in a form useful for calculations
and derive a growth test for poles. The derivative formula for residues
follows in the next subsection.

\subsubsection{Laurent coefficients from contour integrals}

Let $f$ be holomorphic on an annulus
\[
A_{r_1,r_2}(p):=\{z\in\C: r_1<|z-p|<r_2\}.
\]
Assume $f$ admits a Laurent expansion on $A_{r_1,r_2}(p)$ that converges
absolutely and uniformly on every compact subannulus:
\[
f(z)=\sum_{n=-\infty}^{\infty} a_n (z-p)^n.
\]

\begin{proposition}[Laurent coefficient formula]
	\label{prop:laurent-coeff-formula-local}
	For every integer $n$ and every circle $C_\rho:=\{|\zeta-p|=\rho\}$ with
	$r_1<\rho<r_2$,
	\[
	a_n=\frac{1}{2\pi i}\int_{C_\rho}\frac{f(\zeta)}{(\zeta-p)^{n+1}}\,d\zeta.
	\]
	The right-hand side is independent of $\rho$.
\end{proposition}

\begin{proof}
	This is Proposition~\ref{prop:laurent-coeff-formula}, applied with outer
	radius $r_0=r_2$. Its proof integrates the locally uniformly convergent
	series on $C_\rho$ and uses annular contour deformation for radius
	independence.
\end{proof}

\subsubsection{Classification of isolated singularities via Laurent coefficients}

Assume $f$ is holomorphic on $0<|z-p|<r$ and admits a Laurent expansion
\[
f(z)=\sum_{n=-\infty}^{\infty} a_n (z-p)^n,
\qquad 0<|z-p|<r.
\]

\begin{theorem}[Classification of isolated singularities]
	\label{thm:singularity-classification-laurent}
	Exactly one of the following cases occurs:
	\begin{itemize}
		\item The point $p$ is \emph{removable} exactly when $a_{-k}=0$ for every $k\ge 1$.
		\item The point $p$ is a \emph{pole of order $m$} exactly when $a_{-m}\neq 0$ and
		$a_{-k}=0$ for every $k>m$.
		\item The point $p$ is \emph{essential} exactly when infinitely many $a_{-k}$,
		$k\ge 1$, are nonzero.
	\end{itemize}
\end{theorem}

\begin{proof}
	These are the three criteria of Proposition~\ref{prop:sing-class-laurent},
	written with negative indices $-k$. In particular, a pole of order $m$
	is detected by the holomorphic extension of $(z-p)^m f(z)$ taking the
	nonzero value $a_{-m}$ at $p$; an infinite principal part cannot have
	this property for any $m$.
\end{proof}

\subsubsection{Detecting poles via growth}

\begin{proposition}[Blow-up implies finite pole order]
	\label{prop:pole-growth}
	Let $f$ be holomorphic on $0<|z-p|<r$. If $|f(z)|\to\infty$ as $z\to p$, then $p$
	is a pole of some finite order $m\ge1$. Equivalently, $(z-p)^m f(z)$ extends
	holomorphically across $p$ with a nonzero value there.
\end{proposition}

\begin{proof}
	Since $|f(z)|\to\infty$, choose $0<r_0<r$ so that
	$|f(z)|>1$, and in particular $f(z)\neq0$, throughout
	$0<|z-p|<r_0$.
	Define
	\[
	g(z):=\frac{1}{f(z)}
	\qquad (0<|z-p|<r_0).
	\]
	Then $g$ is holomorphic on the punctured disk and $g(z)\to 0$ as $z\to p$.
	In particular, $g$ is bounded near $p$: there exists $M>0$ and $0<r_1<r_0$ such that
	$|g(z)|\le M$ for $0<|z-p|<r_1$.
	
	By Riemann's removable singularity theorem, $g$ extends holomorphically to $p$.
	Denote the extension again by $g$, so $g(p)=0$. Since $g(z)=1/f(z)\ne0$
	on the punctured disk, its extension is not identically zero on $B(p,r_1)$.
	The zero of $g$ at $p$ has finite order: there exist an integer $m\ge 1$ and a
	holomorphic function $h$ on $B(p,r_1)$ with $h(p)\neq 0$ such that
	\[
	g(z)=(z-p)^m h(z).
	\]
	Therefore, on $0<|z-p|<r_1$,
	\[
	f(z)=\frac{1}{g(z)}=\frac{1}{(z-p)^m h(z)}=(z-p)^{-m}\cdot \frac{1}{h(z)}.
	\]
	Since $h(p)\neq 0$, the function $1/h$ is holomorphic near $p$, so $f$ has a pole
	of order $m$ at $p$.
\end{proof}

\subsubsection{Principal part and holomorphic part}

If $p$ is a pole of order $m$, then the Laurent series splits uniquely as
\[
f(z)
=
\underbrace{\sum_{k=-m}^{-1} a_k (z-p)^k}_{\text{principal part}}
\;+\;
\underbrace{\sum_{k=0}^{\infty} a_k (z-p)^k}_{\text{holomorphic part}}.
\]
The integer $m$ is intrinsic.

\begin{lemma}[Uniqueness of the pole order]
	\label{lem:pole-order-unique}
	Suppose the Laurent expansion of $f$ at $p$ satisfies
	\[
	a_{-m}\neq 0,\qquad a_{-k}=0\ \text{for all }k>m.
	\]
	Then $m$ is the unique pole order of $f$ at $p$.
\end{lemma}

\begin{proof}
	Consider $(z-p)^m f(z)$. Multiplying the Laurent expansion by $(z-p)^m$ yields
	\[
	(z-p)^m f(z)=a_{-m}+a_{-m+1}(z-p)+a_{-m+2}(z-p)^2+\cdots,
	\]
	which is a power series convergent near $p$ and takes the nonzero value $a_{-m}$
	at $p$. Hence $(z-p)^m f(z)$ extends holomorphically and is nonvanishing at $p$.
	
	If $N<m$, then $(z-p)^N f(z)$ still has a term $a_{-m}(z-p)^{N-m}$ with negative
	exponent $N-m<0$, so it cannot extend holomorphically to $p$. If $N>m$, then
	\[
	(z-p)^N f(z)=(z-p)^{N-m}\bigl((z-p)^m f(z)\bigr)
	\]
	extends holomorphically but vanishes at $p$ (since $N-m\ge 1$). Thus the only
	exponent producing a holomorphic \emph{nonvanishing} extension at $p$ is $N=m$,
	which is exactly the pole order.
\end{proof}

\subsubsection{Residue and the residue integral formula}

Recall from Definition~\ref{def:residue} that the residue is the
coefficient of $(z-p)^{-1}$. The next statement restates the local
integral formula at the point where it becomes our main computational tool.

\begin{theorem}[Residue integral formula]
	\label{thm:residue-integral-formula}
	Let $f$ be holomorphic on $0<|z-p|<r$. Let $\gamma$ be a positively oriented,
	piecewise $C^1$ simple closed curve contained in $\{0<|z-p|<r\}$
	that winds once around $p$. Then
	\[
	\int_\gamma f(z)\,dz=2\pi i\,\operatorname{Res}(f,p).
	\]
\end{theorem}

\begin{proof}
	This is the local residue formula of Theorem~\ref{thm:cauchy-residue}.
	Indeed, a sufficiently small circle about $p$ lies inside $\gamma$;
	deforming $\gamma$ to that circle crosses no singularity, and termwise
	integration of the Laurent series leaves only $a_{-1}$.
\end{proof}

\begin{figure}[ht]
	\centering
	\begin{tikzpicture}[scale=1.1]
		\draw[thick] (0,0) circle (2);
		\node at (0,0) {$p$};
		\draw[->,thick] (2,0) arc (0:300:2);
		\node at (-1.5,1.3) {$\gamma$};
		\node at (0,-2.4) {$\displaystyle \int_\gamma f(z)\,dz = 2\pi i\, \operatorname{Res}(f,p)$};
	\end{tikzpicture}
	\caption{A contour enclosing one isolated pole contributes
	$2\pi i$ times the residue.}
	\label{fig:single-residue-contour}
\end{figure}

\subsection{Poles, Residues, Holomorphic Simple Connectivity, and the Residue Theorem}
\label{subsec:poles-residues}

We now globalize the local residue formula to contours enclosing many poles,
and relate it to the notion of \emph{holomorphic simple connectivity} (also
called \emph{Cauchy domains} in some texts): domains where integrals of
holomorphic functions over closed curves always vanish.

\subsubsection{Coefficient formulas at a pole}

Let $f$ have a pole of order $k\ge 1$ at $p$. Then $f$ is holomorphic on
$0<|z-p|<r$ for some $r>0$ and has a Laurent expansion
\[
f(z)=\sum_{n=-k}^{\infty} a_n (z-p)^n,
\qquad a_{-k}\neq 0.
\]
Equivalently,
\[
(z-p)^k f(z)=a_{-k}+a_{-k+1}(z-p)+a_{-k+2}(z-p)^2+\cdots
\]
is holomorphic near $p$.

\begin{proposition}[Coefficient formulas at a pole]
	\label{prop:coeff-formulas-pole}
	Let $f$ have a pole of order $k$ at $p$. Then
	\[
	a_{-k+j}
	=
	\frac{1}{j!}
	\left[
	\frac{d^{\,j}}{dz^{\,j}}
	\left( (z-p)^k f(z) \right)
	\right]_{z=p}
	\qquad (j=0,1,2,\dots).
	\]
	In particular, the residue is obtained by taking $j=k-1$:
	\[
	\operatorname{Res}(f,p)=a_{-1}
	=
	\frac{1}{(k-1)!}
	\left[
	\frac{d^{\,k-1}}{dz^{\,k-1}}
	\left( (z-p)^k f(z) \right)
	\right]_{z=p}.
	\]
\end{proposition}

\begin{proof}
	Set
	\[
	g(z):=(z-p)^k f(z).
	\]
	By assumption $g$ is holomorphic in a neighborhood of $p$, so it admits a Taylor
	expansion
	\[
	g(z)=\sum_{j=0}^{\infty} b_j (z-p)^j
	\quad \text{for } |z-p|<r
	\]
	with coefficients
	\[
	b_j=\frac{1}{j!}g^{(j)}(p)
	=
	\frac{1}{j!}\left[
	\frac{d^{\,j}}{dz^{\,j}}
	\left((z-p)^k f(z)\right)\right]_{z=p}.
	\]
	But also
	\[
	g(z)=(z-p)^k f(z)=(z-p)^k\sum_{n=-k}^{\infty} a_n (z-p)^n
	=
	\sum_{n=-k}^{\infty} a_n (z-p)^{n+k}
	=
	\sum_{j=0}^{\infty} a_{-k+j}(z-p)^j,
	\]
	where we relabeled $j=n+k$. Comparing with the Taylor expansion of $g$ shows
	$b_j=a_{-k+j}$ for every $j\ge 0$, giving the stated coefficient formula.
	Taking $j=k-1$ yields $a_{-1}$, i.e.\ the residue formula.
\end{proof}

\subsubsection{Holomorphically simply connected domains}

\begin{definition}[Holomorphically simply connected]
	\label{def:hol-simply-connected}
	A domain $\Omega\subset\C$ is \emph{holomorphically simply connected} if it is
	connected and
	\[
	\int_\gamma f(z)\,dz=0
	\quad \text{for every closed piecewise $C^1$ curve }\gamma\subset\Omega
	\quad \text{and every } f\in \Hol(\Omega).
	\]
	Equivalently, every holomorphic function on $\Omega$ admits a global antiderivative.
\end{definition}

\begin{proposition}[Equivalence with existence of primitives]
	\label{prop:primitive-equivalence}
	For a domain $\Omega\subset\C$, the following are equivalent:
	\begin{enumerate}[label=(\roman*)]
		\item $\displaystyle \int_\gamma f(z)\,dz=0$ for every closed piecewise $C^1$ curve
		$\gamma\subset\Omega$ and every $f\in\Hol(\Omega)$.
		\item Every $f\in\Hol(\Omega)$ has an antiderivative $F\in\Hol(\Omega)$ with $F'=f$.
	\end{enumerate}
\end{proposition}

\begin{proof}
	Assume (i) and prove (ii).
	Fix $f\in\Hol(\Omega)$ and fix a base point $z_0\in\Omega$.
	For any $z\in\Omega$, choose a piecewise $C^1$ curve $\gamma_{z_0,z}$ in $\Omega$
	from $z_0$ to $z$, and define
	\[
	F(z):=\int_{\gamma_{z_0,z}} f(\zeta)\,d\zeta.
	\]
	We first show $F$ is well-defined (independent of the chosen path).
	If $\gamma_1,\gamma_2$ are two curves from $z_0$ to $z$, then $\gamma_1-\gamma_2$
	(the curve obtained by traversing $\gamma_1$ then $\gamma_2$ backwards) is a closed
	curve in $\Omega$, so by (i),
	\[
	0=\int_{\gamma_1-\gamma_2} f(\zeta)\,d\zeta
	=
	\int_{\gamma_1} f(\zeta)\,d\zeta-\int_{\gamma_2} f(\zeta)\,d\zeta.
	\]
	Hence the two path integrals agree, and $F$ is well-defined.
	
	Next we show $F'=f$. Fix $z\in\Omega$ and choose $\varepsilon>0$ so that the disk
	$B(z,\varepsilon)\subset\Omega$. For $h$ with $|h|<\varepsilon$, let $\sigma_h$
	be the straight line segment from $z$ to $z+h$. By additivity of path integrals,
	\[
	F(z+h)-F(z)=\int_{\sigma_h} f(\zeta)\,d\zeta.
	\]
	Parameterize $\sigma_h(t)=z+th$, $0\le t\le 1$. Then $d\zeta=h\,dt$ and
	\[
	\frac{F(z+h)-F(z)}{h}
	=
	\int_0^1 f(z+th)\,dt.
	\]
	As $h\to 0$, continuity of $f$ implies $f(z+th)\to f(z)$ uniformly in $t\in[0,1]$,
	so by dominated convergence,
	\[
	\lim_{h\to 0}\frac{F(z+h)-F(z)}{h}
	=
	\int_0^1 f(z)\,dt
	=f(z).
	\]
	Thus $F'(z)=f(z)$ for all $z\in\Omega$, and $F$ is holomorphic with derivative $f$.
	
	Now assume (ii) and prove (i).
	If $f=F'$ on $\Omega$, then for any closed piecewise $C^1$ curve $\gamma$,
	the fundamental theorem of calculus for complex line integrals gives
	\[
	\int_\gamma f(z)\,dz=\int_\gamma F'(z)\,dz = 0,
	\]
	because the integral of an exact differential over a closed curve vanishes.
\end{proof}

\begin{remark}
	If $\Omega$ is (topologically) simply connected, then $\Omega$ is holomorphically
	simply connected.  For plane domains the converse also holds: if every
	holomorphic function has a primitive, then the index of every closed curve about
	every point of $\C\setminus\Omega$ is zero (apply the hypothesis to
	$1/(z-a)$), which is a standard characterization of simple connectivity.
	Thus the analytic and topological notions agree here.
\end{remark}

\subsubsection{Contour decomposition: reducing a global integral to local circles}

Let $f$ be meromorphic on a domain $\Omega$, and let $\Gamma$ be a positively oriented,
simple closed contour in $\Omega$. Assume $f$ is holomorphic on a neighborhood of
$\Gamma$ and meromorphic in the interior of $\Gamma$, with only finitely many poles
$p_1,\dots,p_n$ inside $\Gamma$.

Choose disjoint disks $B(p_j,\varepsilon_j)\subset\Omega$ so small that:
\begin{enumerate}[label=(\alph*)]
	\item $f$ is holomorphic on each punctured disk $0<|z-p_j|<\varepsilon_j$,
	\item the closed disks $\overline{B(p_j,\varepsilon_j)}$ are pairwise disjoint, and
	\item all these disks lie strictly inside $\Gamma$.
\end{enumerate}
Let $C_j:=\partial B(p_j,\varepsilon_j)$, oriented positively (counterclockwise).

\begin{proposition}[Contour decomposition]
	\label{prop:contour-decomposition}
	With the notation above,
	\[
	\int_{\Gamma} f(z)\,dz
	=
	\sum_{j=1}^{n} \int_{C_j} f(z)\,dz.
	\]
\end{proposition}

\begin{proof}
	Let
	\[
	D:=\mathrm{Int}(\Gamma)\setminus \bigcup_{j=1}^n B(p_j,\varepsilon_j).
	\]
	Then $D$ is a region whose boundary consists of the outer curve $\Gamma$
	and the inner circles $C_j$ with the \emph{opposite} orientation (because each
	$C_j$ bounds a hole removed from the interior). More precisely,
	\[
	\partial D = \Gamma - \sum_{j=1}^n C_j,
	\]
	where ``$-$'' means reversing orientation.
	
	By construction, $f$ is holomorphic on an open neighborhood of $\overline{D}$.
	Cut $D$ along disjoint arcs joining each inner circle to the outer boundary.
	Cauchy's theorem applies to the resulting region, and the two sides of every
	cut contribute equal integrals with opposite orientations. Therefore
	\[
	0=\int_{\partial D} f(z)\,dz
	=
	\int_\Gamma f(z)\,dz-\sum_{j=1}^n \int_{C_j} f(z)\,dz.
	\]
	Rearranging gives the stated decomposition.
\end{proof}

\subsubsection{Residue theorem}

\begin{theorem}[Residue theorem]
	\label{thm:residue-theorem}
	Let $f$ be meromorphic on an open set containing a positively oriented,
	piecewise $C^1$ simple closed contour $\Gamma$ and its interior. Assume that
	$f$ has no pole on $\Gamma$, and that its poles inside $\Gamma$ are
	$p_1,\dots,p_n$ (finitely many). Then
	\[
	\int_\Gamma f(z)\,dz
	=
	2\pi i \sum_{j=1}^{n}\operatorname{Res}(f,p_j).
	\]
\end{theorem}

\begin{proof}
	Choose small circles $C_j$ around each pole as in Proposition
	\ref{prop:contour-decomposition}. By contour decomposition,
	\[
	\int_\Gamma f(z)\,dz=\sum_{j=1}^n \int_{C_j} f(z)\,dz.
	\]
	For each $j$, the function $f$ is holomorphic on a punctured disk about $p_j$, so
	the local residue formula (applied to $C_j$) yields
	\[
	\int_{C_j} f(z)\,dz = 2\pi i\,\operatorname{Res}(f,p_j).
	\]
	Summing over $j=1,\dots,n$ gives
	\[
	\int_\Gamma f(z)\,dz
	=
	\sum_{j=1}^n 2\pi i\,\operatorname{Res}(f,p_j)
	=
	2\pi i \sum_{j=1}^n \operatorname{Res}(f,p_j),
	\]
	which is the desired identity.
\end{proof}

\subsubsection{Winding number and the special case $f(z)=1/(z-p)$}

The integrand $\dfrac{1}{z-p}$ has residue $1$ at $p$ and no other singularities.
Its contour integral measures how many times the curve winds around $p$.

\begin{definition}[Winding number / index]
	\label{def:winding-number}
	Let $\gamma:[a,b]\to\C\setminus\{p\}$ be a closed piecewise $C^1$ curve.
	The winding number of $\gamma$ around $p$ is
	\[
	\operatorname{Ind}(\gamma,p)
	:=
	\frac{1}{2\pi i}\int_\gamma \frac{1}{z-p}\,dz.
	\]
\end{definition}

\begin{proposition}[Integer-valuedness of the winding number]
	\label{prop:winding-number-integer}
	For every closed piecewise $C^1$ curve $\gamma$ avoiding $p$,
	\[
	\operatorname{Ind}(\gamma,p)\in\Z.
	\]
	Moreover, if $\gamma$ is a positively oriented simple closed curve and $p$ lies
	inside $\gamma$, then $\operatorname{Ind}(\gamma,p)=1$, while if $p$ lies outside,
	then $\operatorname{Ind}(\gamma,p)=0$.
\end{proposition}

\begin{proof}
	Step 1: lift the curve to an argument function.
	Let
	\[
	\alpha(t):=\gamma(t)-p \in \Cstar.
	\]
	The map $\alpha:[a,b]\to\Cstar$ is continuous and closed: $\alpha(a)=\alpha(b)$.
	Write $\alpha(t)$ in polar form $\alpha(t)=r(t)e^{i\theta(t)}$.
	Because $[a,b]$ is simply connected and $\alpha(t)\neq 0$ for all $t$, one can choose
	a continuous argument function $\theta:[a,b]\to\R$ such that
	\[
	\alpha(t)=r(t)e^{i\theta(t)}\quad\text{for all }t,
	\]
	with $r(t)=|\alpha(t)|>0$. (Equivalently, we choose a continuous branch of $\Arg$
	along the curve.)
	
	Step 2: compute the integrand as a total derivative.
	Differentiate $\log\alpha(t)=\ln r(t)+i\theta(t)$:
	\[
	\frac{d}{dt}\log(\alpha(t))
	=
	\frac{\alpha'(t)}{\alpha(t)}
	=
	\frac{\gamma'(t)}{\gamma(t)-p}.
	\]
	Hence
	\[
	\int_\gamma \frac{1}{z-p}\,dz
	=
	\int_a^b \frac{\gamma'(t)}{\gamma(t)-p}\,dt
	=
	\int_a^b \frac{d}{dt}\log(\alpha(t))\,dt
	=
	\log(\alpha(b))-\log(\alpha(a)).
	\]
	Since $\alpha(a)=\alpha(b)$, the real parts $\ln r(b)-\ln r(a)$ cancel, leaving
	\[
	\int_\gamma \frac{1}{z-p}\,dz = i\bigl(\theta(b)-\theta(a)\bigr).
	\]
	
	Step 3: conclude integrality.
	Because $\alpha(a)=\alpha(b)$, we must have $r(a)=r(b)$ and
	$e^{i\theta(a)}=e^{i\theta(b)}$, so $\theta(b)-\theta(a)\in 2\pi\Z$.
	Therefore
	\[
	\operatorname{Ind}(\gamma,p)
	=
	\frac{1}{2\pi i}\int_\gamma \frac{1}{z-p}\,dz
	=
	\frac{1}{2\pi i}\cdot i\bigl(\theta(b)-\theta(a)\bigr)
	=
	\frac{\theta(b)-\theta(a)}{2\pi}
	\in \Z.
	\]
	The final statement is the index form of the Jordan curve theorem. A
	positively oriented Jordan curve is homotopic in its complement of $p$
	to a small positive circle when $p$ is in the bounded component, and
	to a constant loop when $p$ is in the unbounded component. Homotopy
	invariance of the integral therefore gives index $1$ and $0$,
	respectively.
\end{proof}

\begin{remark}
	The residue theorem admits an equivalent ``winding number'' form:
	if $f$ is meromorphic, $\Gamma$ is a null-homologous closed cycle in the
	domain avoiding the poles, and only finitely many poles have nonzero index, then
	\[
	\int_\Gamma f(z)\,dz
	=
	2\pi i \sum_{p\in\mathrm{Poles}(f)} \operatorname{Ind}(\Gamma,p)\,\operatorname{Res}(f,p),
	\]
	This version is useful when $\Gamma$ winds multiple times around some poles.  In
	particular, the finiteness condition is automatic when the region on which the
	index is nonzero has compact closure in the domain.
\end{remark}

\subsection{Meromorphic Functions at \texorpdfstring{$\infty$}{Infinity} and Rationality}
\label{subsec:meromorphic-infinity}

\subsubsection{Meromorphic functions on a domain}

Let $\Omega\subset\C$ be a domain. A function $f:\Omega\to\widehat\C$ is \emph{meromorphic}
on $\Omega$ if there exists a \emph{closed discrete} subset $S\subset\Omega$ such that
\begin{enumerate}[label=(\roman*)]
	\item $f$ is holomorphic on $\Omega\setminus S$, and
	\item for each $p\in S$, the point $p$ is a pole of $f$.
\end{enumerate}
We write $f\in M(\Omega)$.

\begin{remark}
	The condition ``closed discrete'' means: for every compact $K\subset\Omega$,
	the intersection $S\cap K$ is finite. Equivalently, $S$ has no accumulation point in $\Omega$.
\end{remark}

\begin{proposition}[The value at a pole]
	\label{prop:change-values-poles}
	Let $\Omega\subset\C$ be a domain and let $S\subset\Omega$ be the
	closed discrete pole set of $f\in M(\Omega)$.  As a map to the Riemann
	sphere, $f$ has a unique continuous---indeed holomorphic---extension across
	each $p\in S$, and that extension has value $f(p)=\infty$.  Assigning an
	arbitrary finite point value at $p$ preserves the punctured germ but does not
	produce a continuous sphere-valued map at $p$.
\end{proposition}

\begin{proof}
	Fix $p\in S$. Since $S$ is discrete, there exists $r>0$ such that
	$B_r(p)\cap S=\{p\}$. Because $f$ has a pole at $p$, there exists $m\ge 1$ and a
	holomorphic function $h$ on $B_r(p)$ with $h(p)\neq 0$ such that
	\[
	f(z)=\frac{h(z)}{(z-p)^m}\qquad (0<|z-p|<r).
	\]
	In the target coordinate $u=1/w$ near $w=\infty$, the map is represented by
	$1/f(z)=(z-p)^m/h(z)$, which extends holomorphically and vanishes at $p$.
	Thus setting $f(p)=\infty$ gives a holomorphic map to $\widehat\C$.
	Moreover $|f(z)|\to\infty$ as $z\to p$, so no finite assigned value can make
	the pointwise map continuous at $p$.
\end{proof}

\begin{example}[The reciprocal sine function]
	The function
	\[
	g(z)=\frac{1}{\sin z}
	\]
	is meromorphic on $\C$ with pole set $S=\{k\pi: k\in\Z\}$. Indeed, $\sin z$ is holomorphic
	and has simple zeros at each $k\pi$, so $1/\sin z$ has simple poles there and is holomorphic
	elsewhere.
\end{example}

\subsubsection{Meromorphicity and singularities at infinity}

Let $f$ be holomorphic on the full exterior region $\{z\in\C:|z|>R\}$
for some $R>0$. This hypothesis makes $\infty$ an isolated point of the
domain's complement; an arbitrary unbounded domain need not contain a
neighborhood of $\infty$.
To study the behavior ``near $\infty$'' we invert the coordinate:
\[
w=\frac{1}{z}.
\]
A neighborhood of $\infty$ in the $z$-plane corresponds to a punctured neighborhood
of $0$ in the $w$-plane. Accordingly, define the transformed function
\[
G(w):=f\!\left(\frac{1}{w}\right),
\]
defined on the entire punctured disk $0<|w|<1/R$.

\begin{definition}[Classification of the singularity at $\infty$]
	We say $f$ has:
	\begin{enumerate}[label=(\roman*)]
		\item a \emph{removable singularity at $\infty$} if $G$ extends holomorphically to $w=0$;
		\item a \emph{pole at $\infty$} if $G$ has a pole at $w=0$;
		\item an \emph{essential singularity at $\infty$} if $G$ has an essential singularity at $w=0$.
	\end{enumerate}
	Equivalently, the type of singularity of $f$ at $\infty$ is the same as the type of
	singularity of $G(w)=f(1/w)$ at $0$.
\end{definition}

\begin{remark}[Laurent expansions near $\infty$]
	If $G$ has a Laurent expansion with finite principal part near $0$,
	\[
	G(w)=\sum_{n=-m}^{\infty} a_n w^n \qquad (0<|w|<\varepsilon),
	\]
	then substituting $w=1/z$ gives an expansion for $f$ for $|z|>1/\varepsilon$:
	\[
	f(z)=G\!\left(\frac{1}{z}\right)=\sum_{n=-m}^{\infty} a_n z^{-n}.
	\]
	Thus $f$ is holomorphic at $\infty$ precisely when only nonnegative powers of $w$
	appear in $G$ (equivalently, only nonpositive powers of $z$ appear in the expansion of $f$),
	and $f$ has a pole at $\infty$ precisely when finitely many but at least
	one nonzero \emph{positive}-power terms in $z$ appear. Infinitely many
	such terms instead signal an essential singularity.
\end{remark}

\subsubsection{Entire functions and growth at infinity}

\begin{proposition}[Pole at $\infty$ forces a polynomial]
	\label{prop:entire-pole-infty-poly}
	If $f$ is entire and has a pole at $\infty$, then $f$ is a polynomial.
	Moreover, if the pole order of $f$ at $\infty$ is $m$, then $\deg f = m$.
\end{proposition}

\begin{proof}
	Assume $f$ is entire and has a pole of order $m\ge 1$ at $\infty$.
	Then $G(w)=f(1/w)$ has a pole of order $m$ at $w=0$, so there exists a Laurent expansion
	\[
	G(w)=\sum_{n=-m}^{\infty} a_n w^n
	\qquad (0<|w|<\varepsilon)
	\]
	with $a_{-m}\neq 0$ and $a_n=0$ for all $n<-m$.
	
	Substitute $w=1/z$ (so $|z|>1/\varepsilon$):
	\[
	f(z)=G\!\left(\frac{1}{z}\right)=\sum_{n=-m}^{\infty} a_n z^{-n}
	=
	\sum_{j=0}^{m} a_{-j} z^{j} \;+\; \sum_{n=1}^{\infty} a_n z^{-n}.
	\]
	We claim the negative-power tail $\sum_{n=1}^{\infty} a_n z^{-n}$ must vanish identically.
	Let
	\[
	P(z):=\sum_{j=0}^{m}a_{-j}z^j,\qquad H(z):=f(z)-P(z).
	\]
	The function $H$ is entire because both $f$ and $P$ are entire. On
	$|z|>1/\varepsilon$ its Laurent representation is
	\[
	H(z)=\sum_{n=1}^{\infty}a_nz^{-n},
	\]
	so $H(z)\to0$ as $z\to\infty$. It is therefore bounded outside a disk
	and, by continuity, bounded on the remaining compact disk. Liouville's
	theorem gives that $H$ is constant; its limit at infinity makes that
	constant zero. Uniqueness of Laurent coefficients then yields
	$a_n=0$ for all $n\ge1$.
	
	Therefore
	\[
	f(z)=\sum_{j=0}^{m} a_{-j} z^{j},
	\]
	a polynomial of degree exactly $m$ because $a_{-m}\neq 0$.
\end{proof}

\begin{corollary}[Nonpolynomial entire functions are essential at $\infty$]
	\label{cor:entire-nonpoly-essential-infty}
	If $f$ is entire and not a polynomial, then $f$ has an essential singularity at $\infty$.
\end{corollary}

\begin{proof}
	Let $f$ be entire and not a polynomial. By Proposition~\ref{prop:entire-pole-infty-poly},
	$f$ cannot have a pole at $\infty$.
	
	If $\infty$ were removable, then $f(z)$ would be bounded for sufficiently
	large $|z|$. It is also bounded on the complementary compact disk, so
	Liouville's theorem would make $f$ constant, contradicting that $f$ is
	nonpolynomial.
	Hence $\infty$ is neither removable nor a pole, so it must be essential.
\end{proof}

\subsubsection{Meromorphic functions on the Riemann sphere are rational}

Let $\widehat{\C}=\C\cup\{\infty\}$ be the Riemann sphere.

\begin{theorem}[Rationality criterion]
	\label{thm:meromorphic-sphere-rational}
	A function $f:\widehat{\C}\to\widehat{\C}$ is meromorphic on $\widehat{\C}$
	if and only if $f$ is a rational function.
\end{theorem}

\begin{proof}
	$(\Leftarrow)$ Suppose $f$ is rational: $f(z)=P(z)/Q(z)$ with coprime polynomials $P,Q$ and
	$Q\not\equiv 0$. Then $f$ is holomorphic on $\C\setminus\{Q=0\}$, and each zero of $Q$
	gives a pole of finite order. Since $Q$ has finitely many zeros, $f$ has finitely many poles
	in $\C$. At $\infty$, write $z=1/w$ and consider $f(1/w)=P(1/w)/Q(1/w)$; multiplying numerator
	and denominator by a suitable power of $w$ shows $f(1/w)$ has at worst a pole at $w=0$.
	Hence $f$ is meromorphic on $\widehat{\C}$.
	
	Conversely, suppose $f$ is meromorphic on $\widehat{\C}$.
	Then $f$ has only isolated poles on $\widehat{\C}$. Because $\widehat{\C}$ is compact,
	any infinite set of poles would have an accumulation point, contradicting that poles are isolated.
	Hence $f$ has only finitely many poles in $\widehat{\C}$.
	
	List the finite poles in $\C$ as $p_1,\dots,p_k$ with orders $m_1,\dots,m_k$.
	Define
	\[
	D(z):=\prod_{j=1}^k (z-p_j)^{m_j},
	\qquad
	F(z):=D(z)\,f(z).
	\]
	On $\C\setminus\{p_1,\dots,p_k\}$, the product is holomorphic.
	Near each $p_j$, the factor $(z-p_j)^{m_j}$ cancels the pole of order $m_j$, so $F$
	extends holomorphically across $p_j$. Therefore $F$ is entire on $\C$.
	
	We now analyze $F$ at $\infty$. Since $f$ is meromorphic on $\widehat{\C}$, it is meromorphic at
	$\infty$; multiplying by the polynomial $D$ cannot create an essential singularity at $\infty$.
	Hence $F$ is entire and has at worst a pole at $\infty$. If it has a
	pole, Proposition~\ref{prop:entire-pole-infty-poly} makes it a
	polynomial. If the singularity is removable, then $F$ is bounded near
	infinity and hence, by Liouville's theorem, constant. Thus in either
	case $F$ is a polynomial.
	
	Therefore,
	\[
	f(z)=\frac{F(z)}{D(z)}
	=
	\frac{F(z)}{\prod_{j=1}^k (z-p_j)^{m_j}},
	\]
	which is a rational function.
\end{proof}

\subsubsection{Exercises}

\begin{exercise}[Pole at $\infty$ and growth]
	Let $f$ be entire. Prove that $f$ has at worst a pole at $\infty$ if and only if there exist
	constants $C>0$ and $m\in\mathbb{N}$ such that
	\[
	|f(z)|\le C(1+|z|)^m
	\qquad\text{for all }z\in\C.
	\]
	Conclude that an entire function has at worst a pole at $\infty$ if and only if it is a
	polynomial, and that the singularity is an actual pole precisely for nonconstant polynomials.
\end{exercise}

\begin{exercise}[Laurent expansion at infinity]
	Show that $f$ is meromorphic at $\infty$ if and only if there exist $R>0$ and integers
	$m\ge 0$ such that for $|z|>R$,
	\[
	f(z)=\sum_{n=-m}^{\infty} a_n z^{-n}
	\]
	with convergence as a Laurent series in $1/z$.
\end{exercise}

\begin{exercise}[Degree and the value at $\infty$]
	Let $f$ be a nonzero meromorphic function on $\widehat{\C}$ and write
	$f(z)=P(z)/Q(z)$ in lowest terms.
	Prove that $f(\infty)=\infty$ if and only if $\deg P > \deg Q$, and that in this case
	the pole order of $f$ at $\infty$ equals $\deg P-\deg Q$.
\end{exercise}

\subsection{Concrete Examples with the Residue Theorem}
\label{subsec:concrete-residue-examples}

The examples below determine isolated singularities and residues by three
complementary methods: Laurent expansion, limit formulas, and algebraic
rearrangement. Comparing the methods makes it easier to choose the shortest
reliable computation in later applications.

\subsubsection{Local order language: zeros and poles}

Recall the order conventions from the identity-theorem and Laurent-series
sections. A zero of order $n$ is a factor $(z-z_0)^n$ times a holomorphic
nonvanishing factor; a pole of order $n$ is the reciprocal kind of local
factorization. Equivalently, the Laurent expansion begins with a nonzero
term of degree $-n$. The coefficient of degree $-1$ is the residue,
whether or not the pole is simple.

\begin{remark}[An order test]
	For nonzero holomorphic germs $u,v$ at $z_0$, set $\ord_{z_0}(u)=0$
	when $u(z_0)\ne0$. If $f=u/v$, then
	\[
	\ord_{z_0}(f)=\ord_{z_0}(u)-\ord_{z_0}(v).
	\]
	In particular, $f$ has a pole at $z_0$ if $\ord_{z_0}(v)>\ord_{z_0}(u)$.
\end{remark}

\subsubsection{Three standard residue extraction tools}

\begin{proposition}[Simple poles: limit formula]
	If $f$ has a simple pole at $z_0$, then
	\[
	\Res(f,z_0)=\lim_{z\to z_0}(z-z_0)f(z).
	\]
\end{proposition}

\begin{proof}
	Write the Laurent expansion at $z_0$:
	\[
	f(z)=\frac{a_{-1}}{z-z_0}+a_0+a_1(z-z_0)+\cdots.
	\]
	Multiplying by $(z-z_0)$ gives
	\[
	(z-z_0)f(z)=a_{-1}+a_0(z-z_0)+a_1(z-z_0)^2+\cdots,
	\]
	which is holomorphic near $z_0$ and tends to $a_{-1}$ as $z\to z_0$. Hence the limit equals
	$\Res(f,z_0)=a_{-1}$.
\end{proof}

\begin{proposition}[Higher-order poles: derivative formula]
	If $f$ has a pole of order $k$ at $z_0$, then
	\[
	\Res(f,z_0)=\frac{1}{(k-1)!}\left[\frac{d^{\,k-1}}{dz^{\,k-1}}\Bigl((z-z_0)^k f(z)\Bigr)\right]_{z=z_0}.
	\]
\end{proposition}

\begin{proof}
	Since $f$ has a pole of order $k$ at $z_0$, the product
	\[
	h(z):=(z-z_0)^k f(z)
	\]
	is holomorphic near $z_0$ and admits a Taylor expansion
	\[
	h(z)=\sum_{m=0}^{\infty} b_m (z-z_0)^m.
	\]
	Then
	\[
	f(z)=\frac{h(z)}{(z-z_0)^k}=\sum_{m=0}^{\infty} b_m (z-z_0)^{m-k}.
	\]
	The residue is the coefficient of $(z-z_0)^{-1}$, which corresponds to $m-k=-1$,
	i.e.\ $m=k-1$. Hence $\Res(f,z_0)=b_{k-1}$.
	On the other hand, $h^{(k-1)}(z_0)=(k-1)!\,b_{k-1}$, so
	\[
	\Res(f,z_0)=\frac{h^{(k-1)}(z_0)}{(k-1)!}
	=\frac{1}{(k-1)!}\left[\frac{d^{\,k-1}}{dz^{\,k-1}}\Bigl((z-z_0)^k f(z)\Bigr)\right]_{z=z_0}.
	\]
\end{proof}

\begin{proposition}[Reciprocal of a simple zero]
	Let $h$ be holomorphic near $z_0$ and suppose $h$ has a simple zero at $z_0$.
	Then $1/h$ has a simple pole at $z_0$ and
	\[
	\Res\!\left(\frac{1}{h},z_0\right)=\frac{1}{h'(z_0)}.
	\]
\end{proposition}

\begin{proof}
	Since $h$ has a simple zero at $z_0$, write $h(z)=(z-z_0)u(z)$ with $u$ holomorphic and
	$u(z_0)\neq 0$. Then
	\[
	\frac{1}{h(z)}=\frac{1}{z-z_0}\cdot\frac{1}{u(z)}.
	\]
	Thus the residue is $1/u(z_0)$. Differentiating $h(z)=(z-z_0)u(z)$ gives $h'(z_0)=u(z_0)$,
	so the residue equals $1/h'(z_0)$.
\end{proof}

\subsubsection{Worked examples}

\begin{example}[A rational function with a zero and poles]
	Consider
	\[
	f(z)=\frac{z+1}{z^3(z-1)}.
	\]
	The numerator vanishes at $z=-1$ while the denominator does not, hence $f$ has a simple zero at $z=-1$.
	
	At $z=0$, the factor $z^3$ in the denominator forces a pole of order $3$. Indeed,
	\[
	g(z):=z^3 f(z)=\frac{z+1}{z-1}
	\]
	is holomorphic near $0$ with $g(0)=-1\neq 0$.
	
	To compute $\Res(f,0)$, expand near $0$:
	\[
	\frac{1}{z-1}=-\frac{1}{1-z}=-\sum_{k=0}^{\infty} z^k,\qquad |z|<1.
	\]
	Hence
	\[
	g(z)=-(z+1)\sum_{k=0}^{\infty}z^k
	=-(1+2z+2z^2+2z^3+\cdots),
	\]
	so
	\[
	f(z)=\frac{g(z)}{z^3}
	=-\frac{1}{z^3}-\frac{2}{z^2}-\frac{2}{z}+\cdots,
	\qquad
	\Res(f,0)=-2.
	\]
	At $z=1$ we have a simple pole, and
	\[
	\Res(f,1)=\lim_{z\to 1}(z-1)\frac{z+1}{z^3(z-1)}=\frac{2}{1^3}=2.
	\]
\end{example}

\begin{example}[Residues at several simple poles]
	Let
	\[
	f(z)=\frac{z^2 + z^7}{(z-2)(z-3)(z-4)(z-5)}.
	\]
	Each pole is simple. For instance,
	\[
	\Res(f,2)=\lim_{z\to 2}(z-2)f(z)=\frac{2^2+2^7}{(2-3)(2-4)(2-5)}.
	\]
	Analogous formulas hold at $z=3,4,5$.
\end{example}

\begin{example}[Essential singularity]
	The function
	\[
	f(z)=\exp\!\left(\frac{1}{z}\right)
	\]
	has an essential singularity at $z=0$ since
	\[
	e^{1/z}=\sum_{n=0}^{\infty}\frac{1}{n!}\,z^{-n}
	=1+\frac{1}{z}+\frac{1}{2!z^2}+\cdots.
	\]
	Hence $\Res(f,0)=1$.
\end{example}

\begin{example}[Combining zeros and poles]
	Consider
	\[
	f(z)=\frac{\sin z}{z^2}.
	\]
	Since $\sin z = z-\frac{z^3}{3!}+\cdots$, we obtain
	\[
	\frac{\sin z}{z^2}=\frac{1}{z}-\frac{z}{3!}+\cdots,
	\]
	so $f$ has a simple pole at $z=0$ and $\Res(f,0)=1$.
\end{example}

\begin{example}[Higher-order pole]
	Let
	\[
	f(z)=\frac{1}{z^2}+\frac{7}{z}+6z.
	\]
	Then $f$ has a pole of order $2$ at $z=0$ and $\Res(f,0)=7$.
\end{example}

\begin{example}[Simple pole from algebraic manipulation]
	Let $f(z)=\dfrac{z}{z+1}$. Then
	\[
	f(z)=\frac{z+1-1}{z+1}=1-\frac{1}{z+1},
	\]
	so $f$ has a simple pole at $z=-1$ and $\Res(f,-1)=-1$.
\end{example}

\begin{example}[Convergence region and residues]
	Let
	\[
	f(z)=\frac{1}{z(1-z)}.
	\]
	On $0<|z|<1$ we have $\dfrac{1}{1-z}=\sum_{n=0}^{\infty}z^n$, hence
	\[
	f(z)=\frac{1}{z}+1+z+z^2+\cdots,
	\qquad
	\Res(f,0)=1.
	\]
\end{example}

\begin{example}[Branch point: no residue]
	The function $f(z)=\log z$ has a branch point at $z=0$, so the singularity is not isolated
	(as a single-valued holomorphic function on a punctured disk). Therefore a residue is not defined.
\end{example}

\begin{example}[Residue via reciprocal of a simple zero]
	Since $\sin z$ has simple zeros at $z=n\pi$ and $\sin'(n\pi)=\cos(n\pi)=(-1)^n\neq 0$,
	each pole of $\csc z=1/\sin z$ is simple and
	\[
	\Res(\csc z,n\pi)=\frac{1}{\cos(n\pi)}=(-1)^n.
	\]
\end{example}

\begin{example}[Residue of a rational function at $z=2$]
	Let
	\[
	f(z)=\frac{2z+7z^2}{(z-2)(z-3)(z-4)(z-5)}.
	\]
	Then
	\[
	\Res(f,2)=\lim_{z\to 2}(z-2)f(z)=\frac{2\cdot 2+7\cdot 2^2}{(2-3)(2-4)(2-5)}.
	\]
\end{example}

\subsubsection{Exercises}

\begin{exercise}
	Compute the residues of the following functions:
	\begin{enumerate}[label=(\alph*)]
		\item $\displaystyle f(z)=\frac{z^2}{(z-1)^3}$ at $z=1$.
		\item $\displaystyle f(z)=\frac{1}{z^2+4}$ at $z=2i$.
		\item $\displaystyle f(z)=\frac{e^{1/z}}{z^2}$ at $z=0$.
		\item $\displaystyle f(z)=\frac{\cos z}{z^4}$ at $z=0$.
	\end{enumerate}
\end{exercise}

\subsection{Residue Calculations on Multiple Contours}
\label{subsec:applications-residue-multicontour}

This subsection applies Theorem~\ref{thm:residue-theorem} to contour integrals.
The examples separate three tasks that should be checked in this order: determine
which singularities lie inside the contour, compute the corresponding residues,
and then form their sum. We also derive two formulas for a quotient $f/g$ near a
simple zero of $g$, including the case $f(a)=g(a)=0$.

\subsubsection{Using the global residue formula}

The global formula and its contour-deformation proof were established in
Theorem~\ref{thm:residue-theorem}. We now use that result without repeating its
hypotheses or proof.

\subsubsection{Example 1: A rational function with poles at $0$ and $1$}

Consider
\[
f(z)=\frac{z^2}{z(z-1)}=\frac{z}{z-1}.
\]
Let $C$ be a positively oriented simple closed contour enclosing both $0$
and $1$, with neither point on $C$.

At $z=0$, the singularity is removable because $f(z)=z/(z-1)$ is analytic at $0$:
indeed,
\[
\Res(f,0)=0.
\]
At $z=1$, the pole is simple, and
\[
\Res(f,1)=\lim_{z\to 1}(z-1)\frac{z}{z-1}=1.
\]
Therefore
\[
\int_C f(z)\,dz = 2\pi i\,[0+1]=2\pi i.
\]

\subsubsection{Example 2: A rational integrand (residues or partial fractions)}

Compute
\[
\int_C \frac{5z-2}{z(z-1)}\,dz,
\]
where $C$ is a positively oriented simple closed contour enclosing $0$ and
$1$, with neither point on $C$.

Both poles are simple. At $z=0$,
\[
\Res\!\left(\frac{5z-2}{z(z-1)},0\right)
=\lim_{z\to 0}\frac{5z-2}{z-1}=2.
\]
At $z=1$,
\[
\Res\!\left(\frac{5z-2}{z(z-1)},1\right)
=\lim_{z\to 1}\frac{5z-2}{z}=3.
\]
Hence
\[
\int_C \frac{5z-2}{z(z-1)}\,dz
=
2\pi i(2+3)=10\pi i.
\]

\subsubsection{Example 3: Laurent expansion at $0$ and detecting the residue}

Compute
\[
\int_{|z|=1} z^{-2}\sin\!\left(\frac{1}{z}\right)\, dz.
\]
On $0<|z|<\infty$ we have
\[
\sin\!\left(\frac{1}{z}\right)
=
\frac{1}{z}-\frac{1}{3!z^3}+\frac{1}{5!z^5}-\cdots,
\]
so
\[
z^{-2}\sin\!\left(\frac{1}{z}\right)
=
\frac{1}{z^3}-\frac{1}{3!z^5}+\frac{1}{5!z^7}-\cdots.
\]
There is no $(1/z)$ term, hence $\Res(f,0)=0$ and therefore
\[
\int_{|z|=1} z^{-2}\sin(1/z)\,dz = 0.
\]

\subsubsection{Example 4: A pole inside the contour and expansion around $z=2$}

Compute
\[
\int_{|z-2|=1} \frac{dz}{z^2(z-2)^4}.
\]
The contour encloses $z=2$ and excludes $z=0$. Thus only the pole at $z=2$ contributes.

Write
\[
\frac{1}{z^2}
=
\frac{1}{(2+(z-2))^2}
=
\frac{1}{4}\,\frac{1}{\bigl(1+\frac{z-2}{2}\bigr)^2}.
\]
For $|z-2|<2$ we may expand (binomial series)
\[
\frac{1}{(1+w)^2}=\sum_{m=0}^{\infty}(-1)^m(m+1)w^m
=1-2w+3w^2-4w^3+\cdots.
\]
With $w=\frac{z-2}{2}$ this gives
\[
\frac{1}{z^2}
=
\frac{1}{4}\sum_{m=0}^{\infty}(-1)^m(m+1)\left(\frac{z-2}{2}\right)^m.
\]
Therefore
\[
\frac{1}{z^2(z-2)^4}
=
\frac{1}{4}\sum_{m=0}^{\infty}(-1)^m(m+1)\left(\frac{1}{2}\right)^m (z-2)^{m-4}.
\]
The residue at $z=2$ is the coefficient of $(z-2)^{-1}$, i.e.\ the term with $m-4=-1$,
so $m=3$. Hence
\[
\Res\!\left(\frac{1}{z^2(z-2)^4},2\right)
=
\frac{1}{4}\,(-1)^3(3+1)\left(\frac{1}{2}\right)^3
=
-\frac{1}{8}.
\]
By the residue theorem (here: single pole inside),
\[
\int_{|z-2|=1} \frac{dz}{z^2(z-2)^4}
=
2\pi i\left(-\frac{1}{8}\right)
=
-\frac{\pi i}{4}.
\]

\subsubsection{Example 5: Residues of $\csc z$ at finitely many enclosed poles}

Consider $f(z)=\csc z = 1/\sin z$. The zeros of $\sin z$ are $k\pi$, and
\[
\sin'(k\pi)=\cos(k\pi)=(-1)^k\neq 0,
\]
so each pole of $\csc z$ is simple. If a contour $C$ encloses $0,\pi,2\pi$ and no other
poles, then
\[
\Res(\csc z,k\pi)=\frac{1}{\sin'(k\pi)}=\frac{1}{\cos(k\pi)}=(-1)^k,
\]
hence
\[
\int_C \frac{dz}{\sin z}
=
2\pi i\,\bigl[1-1+1\bigr]
=2\pi i.
\]

\subsubsection{Quotients, simple zeros, and residues}

The next two propositions are often the fastest way to compute residues of
quotients $f/g$.

\begin{proposition}[Residue at a simple zero of the denominator]
	\label{prop:residue-simple-zero-denominator}
	Let $f$ and $g$ be analytic near $a$, and assume $g$ has a simple zero at $a$:
	\[
	g(a)=0,\qquad g'(a)\neq 0,
	\]
	and assume $f(a)\ne0$.
	Then $\dfrac{f}{g}$ has a simple pole at $a$ with
	\[
	\Res\!\left(\frac{f}{g},a\right)=\frac{f(a)}{g'(a)}.
	\]
\end{proposition}

\begin{proof}
	Since $g(a)=0$ and $g'(a)\neq 0$, the Taylor expansion of $g$ at $a$ begins with a
	nonzero linear term. Define
	\[
	h(z):=
	\begin{cases}
		\dfrac{g(z)}{z-a}, & z\neq a,\\[4pt]
		g'(a), & z=a.
	\end{cases}
	\]
	We show $h$ is analytic near $a$ and satisfies $h(a)=g'(a)\neq 0$.
	Indeed, write the Taylor series
	\[
	g(z)=\sum_{m=1}^{\infty} \frac{g^{(m)}(a)}{m!}(z-a)^m,
	\]
	(valid near $a$), then for $z\neq a$,
	\[
	\frac{g(z)}{z-a}=\sum_{m=1}^{\infty}\frac{g^{(m)}(a)}{m!}(z-a)^{m-1}
	=\sum_{m=0}^{\infty}\frac{g^{(m+1)}(a)}{(m+1)!}(z-a)^m,
	\]
	which extends holomorphically to $z=a$ with value $g'(a)$. Thus $h$ is analytic and
	nonvanishing near $a$.
	
	Now write $g(z)=(z-a)h(z)$. Then
	\[
	\frac{f(z)}{g(z)}=\frac{1}{z-a}\cdot\frac{f(z)}{h(z)}.
	\]
	Since $f$ and $h$ are analytic and $h(a)\neq 0$, the quotient $f/h$ is analytic near $a$.
	Its value at $a$ equals $f(a)/h(a)=f(a)/g'(a)$. Therefore the Laurent expansion of $f/g$ at $a$
	has the form
	\[
	\frac{f(z)}{g(z)}=\frac{f(a)}{g'(a)}\cdot\frac{1}{z-a}+\text{holomorphic terms},
	\]
	so the residue is $f(a)/g'(a)$.
\end{proof}

\begin{proposition}[Limit of a quotient when both numerator and denominator vanish]
	\label{prop:limit-quotient-both-vanish}
	Let $f$ and $g$ be analytic near $a$ with
	\[
	f(a)=g(a)=0,\qquad g'(a)\neq 0.
	\]
	Then $f/g$ extends holomorphically across $a$ and
	\[
	\lim_{z\to a}\frac{f(z)}{g(z)}=\frac{f'(a)}{g'(a)}.
	\]
\end{proposition}

\begin{proof}
	As above, since $g$ has a simple zero at $a$, write $g(z)=(z-a)h(z)$ with $h$ analytic and
	$h(a)=g'(a)\neq 0$. Because $f(a)=0$, define
	\[
	k(z):=
	\begin{cases}
		\dfrac{f(z)}{z-a}, & z\neq a,\\[4pt]
		f'(a), & z=a.
	\end{cases}
	\]
	The same Taylor-series argument shows $k$ is analytic near $a$ and $k(a)=f'(a)$, and we have
	$f(z)=(z-a)k(z)$. Hence for $z\neq a$,
	\[
	\frac{f(z)}{g(z)}=\frac{(z-a)k(z)}{(z-a)h(z)}=\frac{k(z)}{h(z)}.
	\]
	The right-hand side is analytic near $a$ (since $h(a)\neq 0$), so $f/g$ extends holomorphically
	across $a$. Taking $z\to a$ gives
	\[
	\lim_{z\to a}\frac{f(z)}{g(z)}=\frac{k(a)}{h(a)}=\frac{f'(a)}{g'(a)}.
	\]
\end{proof}

\subsection{Applications of the Residue Theorem to Real Integrals}
\label{subsec:residue-real-integrals}

To evaluate a real integral by residues, one must choose a contour, locate the
enclosed poles, compute their residues, and estimate every auxiliary arc. The
examples carry out all four steps, including the indentation argument for
$\int_0^\infty \frac{\sin x}{x}\,dx$.

\subsubsection{Simple poles and residues}

At a simple pole the principal part has only one term:
$f(z)=a_{-1}/(z-z_0)+h(z)$ with $h$ holomorphic near $z_0$.
Its coefficient $a_{-1}$ is the residue, as in Definition~\ref{def:residue}.

\begin{lemma}[Residue at a simple pole via a limit]
	\label{lem:res-simple-limit}
	If $f$ has a simple pole at $z_0$, then
	\[
	\Res(f,z_0)=\lim_{z\to z_0}(z-z_0)f(z).
	\]
\end{lemma}

\begin{proof}
	By definition, write $f(z)=\frac{a_{-1}}{z-z_0}+h(z)$ with $h$ holomorphic near $z_0$.
	Then
	\[
	(z-z_0)f(z)=a_{-1}+(z-z_0)h(z)\xrightarrow[z\to z_0]{}a_{-1}.
	\]
\end{proof}

\subsubsection{A standard large--arc estimate (Jordan-type lemma)}

\begin{lemma}[Upper--half--plane exponential decay]
	\label{lem:jordan-basic}
	Fix $a>0$.  Let $R>0$ and parametrize the upper semicircle
	\[
	\Gamma_R:=\{z=Re^{it}:0\le t\le \pi\}.
	\]
	Then for every $z\in\Gamma_R$,
	\[
	|e^{iaz}| = e^{-a\,\Im z}=e^{-aR\sin t}\le 1.
	\]
	Moreover, if $H$ is meromorphic and satisfies $|H(z)|\le \frac{M}{R^{1+\delta}}$
	on $\Gamma_R$ for some constants $M,\delta>0$ (uniformly in $t$), then
	\[
	\left|\int_{\Gamma_R} e^{iaz}H(z)\,dz\right|
	\le \pi R\cdot \frac{M}{R^{1+\delta}}
	\longrightarrow 0\qquad(R\to\infty).
	\]
\end{lemma}

\begin{proof}
	For $z=Re^{it}$ we have $\Im z=R\sin t\ge 0$, hence
	\[
	|e^{iaz}|=\exp(\Re(iaz))=\exp(-a\,\Im z)=e^{-aR\sin t}\le 1.
	\]
	Also $|dz|=R\,dt$ on $\Gamma_R$.  Therefore,
	\[
	\left|\int_{\Gamma_R} e^{iaz}H(z)\,dz\right|
	\le \int_0^\pi |e^{iaRe^{it}}|\,|H(Re^{it})|\,|dz|
	\le \int_0^\pi 1\cdot \frac{M}{R^{1+\delta}}\cdot R\,dt
	= \frac{\pi M}{R^\delta}\xrightarrow[R\to\infty]{}0.
	\]
\end{proof}

\subsubsection{Example 1: $\displaystyle \int_{0}^{\infty} \frac{\cos x}{1+x^{4}}\, dx$}

Consider
\[
I:=\int_{0}^{\infty} \frac{\cos x}{1+x^{4}}\, dx.
\]
Let
\[
f(z)=\frac{e^{iz}}{1+z^{4}}.
\]
Let $C_R$ be the positively oriented contour consisting of the real segment
$[-R,R]$ and the upper semicircle $\Gamma_R$.

\paragraph{Step 1: poles in the upper half--plane.}
The equation $1+z^4=0$ is $z^4=e^{i\pi}$, so the roots are
\[
z_k=e^{i(\pi/4+k\pi/2)}\qquad(k=0,1,2,3).
\]
Those in the upper half--plane are
\[
z_1=e^{i\pi/4}=\frac{1+i}{\sqrt2},\qquad
z_2=e^{3i\pi/4}=\frac{-1+i}{\sqrt2}.
\]
Each is a simple pole because $(1+z^4)'=4z^3$ does not vanish at a root.

\paragraph{Step 2: residues.}
By Lemma~\ref{lem:res-simple-limit},
\[
\Res(f,z_0)=\frac{e^{iz_0}}{(1+z^4)'|_{z=z_0}}
=\frac{e^{iz_0}}{4z_0^3}.
\]
Hence
\[
\int_{C_R} f(z)\,dz
=
2\pi i\Bigl(\Res(f,z_1)+\Res(f,z_2)\Bigr).
\]

\paragraph{Step 3: large arc vanishes.}
On $\Gamma_R$ we have $|1+z^4|\ge ||z|^4-1|=R^4-1$, so
\[
|f(z)|=\left|\frac{e^{iz}}{1+z^4}\right|
\le \frac{|e^{iz}|}{R^4-1}
\le \frac{1}{R^4-1}.
\]
Therefore by Lemma~\ref{lem:jordan-basic} (with $\delta=3$),
\[
\left|\int_{\Gamma_R} f(z)\,dz\right|
\le \pi R\cdot \frac{1}{R^4-1}\longrightarrow 0\qquad(R\to\infty).
\]

\paragraph{Step 4: pass to the limit and extract the real part.}
We have
\[
\int_{-R}^{R}\frac{e^{ix}}{1+x^4}\,dx+\int_{\Gamma_R}f(z)\,dz
=2\pi i\Bigl(\Res(f,z_1)+\Res(f,z_2)\Bigr).
\]
Letting $R\to\infty$ and using $\int_{\Gamma_R}f\to 0$ gives
\[
\int_{-\infty}^{\infty}\frac{e^{ix}}{1+x^4}\,dx
=2\pi i\Bigl(\Res(f,z_1)+\Res(f,z_2)\Bigr).
\]
Taking real parts,
\[
\int_{-\infty}^{\infty}\frac{\cos x}{1+x^4}\,dx
=\Re\left(2\pi i\Bigl(\Res(f,z_1)+\Res(f,z_2)\Bigr)\right).
\]
Since the integrand $\frac{\cos x}{1+x^4}$ is even,
\[
2I=\int_{-\infty}^{\infty}\frac{\cos x}{1+x^4}\,dx.
\]
To make the final simplification explicit, put $t=1/\sqrt2$. Since
$z_1=t(1+i)$, $z_2=t(-1+i)$, $z_1^3=z_2$, and $z_2^3=z_1$, the sum of
residues is
\[
\frac{e^{-t}}4\left(\frac{e^{it}}{z_2}+\frac{e^{-it}}{z_1}\right)
=-\frac{it e^{-t}}2\bigl(\cos t+\sin t\bigr).
\]
Thus the classical closed form is
\[
\int_{0}^{\infty}\frac{\cos x}{1+x^4}\,dx
=
\frac{\pi}{2\sqrt2}\,e^{-1/\sqrt2}\Bigl(\cos(1/\sqrt2)+\sin(1/\sqrt2)\Bigr).
\]

\subsubsection{Example 2: $\displaystyle \int_{-\infty}^{\infty} \frac{\sin x}{1+x^{2}}\, dx$}

Let
\[
f(z)=\frac{e^{iz}}{1+z^2}.
\]
The poles are at $z=\pm i$ and are simple.  Only $z=i$ lies in the upper half--plane.

\paragraph{Step 1: residue at $z=i$.}
Since $(1+z^2)'=2z$, by Lemma~\ref{lem:res-simple-limit},
\[
\Res(f,i)=\frac{e^{ii}}{2i}=\frac{e^{-1}}{2i}.
\]

\paragraph{Step 2: apply the semicircle contour and let $R\to\infty$.}
Exactly as in Example 1, the arc contribution vanishes because
$|f(z)|\le \frac{1}{|z|^2-1}$ on $\Gamma_R$ (for $R$ large) and $|e^{iz}|\le 1$ in the upper half--plane.
Thus
\[
\int_{-\infty}^{\infty}\frac{e^{ix}}{1+x^2}\,dx
=2\pi i\,\Res(f,i)
=2\pi i\cdot \frac{e^{-1}}{2i}
=\pi e^{-1}.
\]
The right-hand side is real, so taking imaginary parts gives
\[
\int_{-\infty}^{\infty}\frac{\sin x}{1+x^2}\,dx=0.
\]
This conclusion also follows immediately from oddness; here it serves as a
consistency check on extracting real and imaginary parts of the contour integral.

\subsubsection{Example 3: $\displaystyle \int_{0}^{\infty} \frac{\sin x}{x}\, dx$ (indentation at the origin)}

Consider
\[
f(z)=\frac{e^{iz}}{z}.
\]
The point $z=0$ is a simple pole with residue $1$.

\paragraph{Contour.}
Fix $R>1$ and $0<\varepsilon<1$.  Let $C_{R,\varepsilon}$ be the positively oriented contour
consisting of:
\begin{itemize}
	\item the segment $[-R,-\varepsilon]$ on the real axis,
	\item the small semicircle $\gamma_\varepsilon$ of radius $\varepsilon$ in the upper half--plane
	joining $-\varepsilon$ to $+\varepsilon$, traversed clockwise,
	\item the segment $[\varepsilon,R]$ on the real axis,
	\item the large semicircle $\Gamma_R$ in the upper half--plane joining $R$ to $-R$.
\end{itemize}
The indentation excludes the pole at $0$.

\paragraph{Residue theorem on $C_{R,\varepsilon}$.}
The function is analytic on and inside the indented contour, so
\[
	\int_{C_{R,\varepsilon}} f(z)\,dz=0.
\]
Writing the contour integral as a sum gives
\begin{equation}
	\label{eq:keyhole-sinx-over-x}
	\int_{-R}^{-\varepsilon}\frac{e^{ix}}{x}\,dx
	+\int_{\gamma_\varepsilon}\frac{e^{iz}}{z}\,dz
	+\int_{\varepsilon}^{R}\frac{e^{ix}}{x}\,dx
	+\int_{\Gamma_R}\frac{e^{iz}}{z}\,dz
	=0.
\end{equation}

\paragraph{Large arc $\Gamma_R$ vanishes.}
On $\Gamma_R$, $|e^{iz}|\le 1$ and $|z|=R$, hence $|e^{iz}/z|\le 1/R$ and
\[
\left|\int_{\Gamma_R}\frac{e^{iz}}{z}\,dz\right|
\le \pi R\cdot \frac{1}{R}=\pi.
\]
To get vanishing, we use the sharper decay $|e^{iz}|=e^{-R\sin t}$:
\[
\left|\int_{\Gamma_R}\frac{e^{iz}}{z}\,dz\right|
\le \int_0^\pi \frac{e^{-R\sin t}}{R}\,R\,dt
=\int_0^\pi e^{-R\sin t}\,dt
\longrightarrow 0\qquad(R\to\infty),
\]
by dominated convergence (for $t\in(0,\pi)$, $\sin t>0$ so $e^{-R\sin t}\to 0$, and the integrand is bounded by $1$).

\paragraph{Small arc $\gamma_\varepsilon$.}
Parametrize $\gamma_\varepsilon$ by $z=\varepsilon e^{it}$,
$t\in[\pi,0]$. Then $dz=i\varepsilon e^{it}dt$ and
\[
\int_{\gamma_\varepsilon}\frac{e^{iz}}{z}\,dz
=
	\int_\pi^0 \frac{e^{i\varepsilon e^{it}}}{\varepsilon e^{it}}\, i\varepsilon e^{it}\,dt
=
	i\int_\pi^0 e^{i\varepsilon e^{it}}\,dt
	\longrightarrow i\int_\pi^0 1\,dt=-i\pi
\qquad(\varepsilon\to 0),
\]
again by dominated convergence (since $|e^{i\varepsilon e^{it}}|=e^{-\varepsilon\sin t}\le 1$).

\paragraph{Combine the real-axis pieces.}
Use the substitution $x\mapsto -x$ in the first integral:
\[
\int_{-R}^{-\varepsilon}\frac{e^{ix}}{x}\,dx
=
\,-\int_{\varepsilon}^{R}\frac{e^{-ix}}{x}\,dx.
\]
Hence the sum of the two real-axis pieces is
\[
\int_{\varepsilon}^{R}\frac{e^{ix}}{x}\,dx-\int_{\varepsilon}^{R}\frac{e^{-ix}}{x}\,dx
=
2i\int_{\varepsilon}^{R}\frac{\sin x}{x}\,dx.
\]
Equation~\eqref{eq:keyhole-sinx-over-x} therefore gives
\[
2i\int_{\varepsilon}^{R}\frac{\sin x}{x}\,dx
=
-\int_{\gamma_\varepsilon}\frac{e^{iz}}{z}\,dz
-\int_{\Gamma_R}\frac{e^{iz}}{z}\,dz.
\]
Now let $R\to\infty$ and then $\varepsilon\to0$. Since the large
arc tends to $0$ and the small arc tends to $-i\pi$, we obtain
\[
2i\int_{0}^{\infty}\frac{\sin x}{x}\,dx
=
i\pi,
\]
hence
\[
\int_{0}^{\infty}\frac{\sin x}{x}\,dx=\frac{\pi}{2}.
\]

\subsubsection{Example 4: A semicircle template for cosine integrals}

Let $a>0$ and $n=2m\ge2$ be even. Consider
\[
\int_{0}^{\infty}\frac{\cos(ax)}{1+x^{n}}\,dx.
\]
Set
\[
f(z)=\frac{e^{iaz}}{1+z^{n}}.
\]
The poles are the roots of $1+z^n=0$, i.e.
\[
z_k=e^{i( (2k+1)\pi/n )},\qquad k=0,1,\dots,n-1,
\]
and the relevant ones are those with $\Im z_k>0$. Each pole is simple since $(1+z^n)'=n z^{n-1}\neq 0$ at a root.
For such a root $z_0$,
\[
\Res\!\left(\frac{e^{iaz}}{1+z^n},z_0\right)
=
\frac{e^{iaz_0}}{n z_0^{\,n-1}}.
\]
On the upper semicircle $\Gamma_R$, the exponential factor satisfies $|e^{iaz}|\le 1$
(Lemma~\ref{lem:jordan-basic}), and the denominator has modulus at least $R^n-1$;
thus the large-arc contribution vanishes as $R\to\infty$.
There are no poles on the real axis when $n$ is even. Since
$\cos(ax)/(1+x^n)$ is even, taking the real part gives the concrete formula
\[
\int_0^\infty\frac{\cos(ax)}{1+x^n}\,dx
=\frac12\operatorname{Re}\!\left(
2\pi i\!\sum_{\Im z_k>0}
\frac{e^{iaz_k}}{n z_k^{\,n-1}}\right).
\]
The same semicircle integral does not by itself determine the half-line
sine integral: $\sin(ax)/(1+x^n)$ is odd, so its whole-line integral is
zero. For odd $n$, the denominator also has a pole at $-1$ on the real
axis, requiring an indentation or a different contour.

\subsubsection{Exercises}

\begin{exercise}
	Evaluate $\displaystyle \int_{0}^{\infty} \frac{\cos x}{x^{2}+4}\, dx$ using the residue theorem.
\end{exercise}

\begin{exercise}
	Show that $\displaystyle \int_{0}^{\infty} \frac{\sin(ax)}{x}\, dx = \frac{\pi}{2}$ for all $a>0$.
	(Hint: scale $x\mapsto ax$, or redo Example~3 with $e^{iaz}/z$.)
\end{exercise}

\begin{exercise}
	Compute $\displaystyle \int_{-\infty}^{\infty} \frac{x\sin x}{1+x^{4}}\, dx$.
	(Hint: consider $\int_{-\infty}^{\infty} \frac{z e^{iz}}{1+z^{4}}\,dz$ and take imaginary parts.)
\end{exercise}

\subsection{Jordan's Lemma and Branch-Cut Integrals}
\label{subsec:advanced-contours}

\paragraph{Conventions.}
Throughout, all contours are piecewise $C^1$ and positively oriented unless stated otherwise.
Residues are denoted by $\operatorname{Res}(f,p)$.
When using branch cuts, we specify the branch explicitly.

\subsubsection{Jordan's Lemma}

\begin{lemma}[Jordan's lemma]\label{lem:jordan}
	Let $a>0$, let $C_R:=\{Re^{it}:0\le t\le\pi\}$, and let $g$ be
	continuous on $C_R$ with $|g(z)|\le M$ there. Then
	\[
	\left|\int_{C_R} e^{iaz}g(z)\,dz\right|\le \frac{\pi M}{a}.
	\]
	Moreover, if for all sufficiently large $R$ one has
	$|g(z)|\le M/R$ on $C_R$, then
	\[
	\left|\int_{C_R} e^{iaz}g(z)\,dz\right|\le \frac{\pi M}{aR}\xrightarrow[R\to\infty]{}0.
	\]
\end{lemma}

\begin{proof}
	Parametrize $C_R$ by $z=Re^{it}$, $0\le t\le \pi$, so $dz=iRe^{it}\,dt$ and
	\[
	\left|\int_{C_R} e^{iaz}g(z)\,dz\right|
	\le \int_{0}^{\pi} \bigl|e^{iaRe^{it}}\bigr|\,|g(Re^{it})|\,R\,dt.
	\]
	Now $e^{iaRe^{it}} = e^{iaR(\cos t+i\sin t)}$, hence
	\[
	\bigl|e^{iaRe^{it}}\bigr| = e^{-aR\sin t}.
	\]
	Using $|g|\le M$ on $C_R$ gives
	\[
	\left|\int_{C_R} e^{iaz}g(z)\,dz\right|
	\le MR\int_{0}^{\pi} e^{-aR\sin t}\,dt.
	\]
	For $t\in[0,\pi]$ we have the elementary bound
	\[
	\sin t \ge \frac{2}{\pi}\,t \quad \text{for } 0\le t\le \frac{\pi}{2},
	\qquad
	\sin t \ge \frac{2}{\pi}\,(\pi-t)\quad \text{for } \frac{\pi}{2}\le t\le \pi.
	\]
	Splitting at $\pi/2$ and changing variables,
	\begin{align*}
		\int_{0}^{\pi} e^{-aR\sin t}\,dt
		&\le \int_{0}^{\pi/2} e^{-aR(2t/\pi)}\,dt + \int_{\pi/2}^{\pi} e^{-aR(2(\pi-t)/\pi)}\,dt\\
		&= 2\int_{0}^{\pi/2} e^{-2aRt/\pi}\,dt
		= 2\cdot \frac{\pi}{2aR}\left(1-e^{-aR}\right)
		\le \frac{\pi}{aR}.
	\end{align*}
	Therefore
	\[
	\left|\int_{C_R} e^{iaz}g(z)\,dz\right|\le MR\cdot \frac{\pi}{aR}=\frac{\pi M}{a}.
	\]
	If additionally $|g(z)|\le M/R$ on $C_R$, the same estimate yields
	\[
	\left|\int_{C_R} e^{iaz}g(z)\,dz\right|
	\le \frac{M}{R}\,R\int_{0}^{\pi} e^{-aR\sin t}\,dt
	\le \frac{M}{R}\,R\cdot \frac{\pi}{aR}=\frac{\pi M}{aR}\to 0.
	\]
\end{proof}

\subsubsection{A warning example: $\displaystyle \int_0^\infty \frac{x}{1+x^2}\,dx$}

The integral
\[
\int_0^\infty \frac{x}{1+x^2}\,dx
=
\left[\frac12\log(1+x^2)\right]_{0}^{\infty}
=\infty
\]
diverges (logarithmically).  Nevertheless,
\[
\operatorname{p.v.}\int_{-\infty}^{\infty}\frac{x}{1+x^2}\,dx=0
\]
by oddness.  Contour methods frequently compute \emph{principal values} when poles lie on (or symmetry is used across) the real axis.

\subsubsection{Branch-cut method: $\displaystyle \int_0^\infty \frac{x^{1/3}}{1+x^2}\,dx$}

We now compute the genuinely branch-cut-driven integral
\[
I:=\int_{0}^{\infty}\frac{x^{1/3}}{1+x^2}\,dx.
\]

\paragraph{Step 0: Choose a branch.}
Let $\Log z = \ln|z| + i\Arg z$ with $0<\Arg z<2\pi$ (branch cut along the \emph{positive} real axis),
and define
\[
z^{1/3}:=\exp\!\left(\frac13\Log z\right).
\]
With this choice, for $x>0$:
\[
x^{1/3}_{\text{upper}} = x^{1/3} \quad (\Arg x=0^+),
\qquad
x^{1/3}_{\text{lower}} = x^{1/3}e^{2\pi i/3}\quad (\Arg x=2\pi^-).
\]

Set
\[
f(z)=\frac{z^{1/3}}{1+z^2}.
\]
The poles of $f$ are at $z=i$ and $z=-i$, both simple.

\paragraph{Step 1: Keyhole contour and reduction to $(1-e^{2\pi i/3})I$.}
Let $\Gamma_{R,\varepsilon}$ be the keyhole contour around the positive real axis:
\begin{itemize}
	\item upper ray: $z=x$ from $\varepsilon$ to $R$ (just above the cut),
	\item outer circle: $|z|=R$ from angle $0$ to $2\pi$ avoiding the cut,
	\item lower ray: $z=x$ from $R$ back to $\varepsilon$ (just below the cut),
	\item inner circle: $|z|=\varepsilon$ closing the contour around the origin.
\end{itemize}
Because $f(z)=O(|z|^{-5/3})$ as $|z|\to\infty$ and $f(z)=O(|z|^{1/3})$ as $z\to 0$,
the outer and inner circular contributions vanish as $R\to\infty$ and $\varepsilon\to 0$:
\[
\int_{|z|=R} f(z)\,dz \to 0,
\qquad
\int_{|z|=\varepsilon} f(z)\,dz \to 0.
\]

Along the upper ray we have $z^{1/3}=x^{1/3}$ and $dz=dx$, so the contribution tends to $I$.
Along the lower ray, the contour runs from $R$ to $\varepsilon$ (reverse direction), and
$z^{1/3}=x^{1/3}e^{2\pi i/3}$, hence the contribution tends to
\[
\int_{R}^{\varepsilon}\frac{(x)^{1/3}e^{2\pi i/3}}{1+x^2}\,dx
= -e^{2\pi i/3}\int_{\varepsilon}^{R}\frac{x^{1/3}}{1+x^2}\,dx
\to -e^{2\pi i/3} I.
\]
Therefore,
\begin{equation}\label{eq:keyhole-factor}
	\lim_{R\to\infty,\ \varepsilon\to 0}\int_{\Gamma_{R,\varepsilon}} f(z)\,dz
	=
	(1-e^{2\pi i/3})\,I.
\end{equation}

\paragraph{Step 2: Compute residues at $\pm i$ on this branch.}
Since $1+z^2=(z-i)(z+i)$, for a simple pole at $z=z_0\in\{i,-i\}$,
\[
\operatorname{Res}(f,z_0)=\lim_{z\to z_0}(z-z_0)\frac{z^{1/3}}{(z-i)(z+i)}
=\frac{z_0^{1/3}}{2z_0}.
\]
On our branch $0<\Arg z<2\pi$, we have
\[
i = e^{i\pi/2}\ \mathrel{\text{and hence}}\ i^{1/3}=e^{i\pi/6},
\qquad
-i = e^{i3\pi/2}\ \mathrel{\text{and hence}}\ (-i)^{1/3}=e^{i\pi/2}=i.
\]
Hence
\[
\operatorname{Res}(f,i)=\frac{i^{1/3}}{2i}=\frac{e^{i\pi/6}}{2i},
\qquad
\operatorname{Res}(f,-i)=\frac{(-i)^{1/3}}{-2i}=\frac{i}{-2i}=-\frac12.
\]
Therefore,
\begin{equation}\label{eq:sum-res}
	\sum \operatorname{Res}(f,\pm i)=\frac{e^{i\pi/6}}{2i}-\frac12.
\end{equation}

\paragraph{Step 3: Apply the residue theorem and solve for $I$.}
The keyhole contour encloses both poles $\pm i$, so
\[
\int_{\Gamma_{R,\varepsilon}} f(z)\,dz
=2\pi i\left(\operatorname{Res}(f,i)+\operatorname{Res}(f,-i)\right).
\]
Taking limits and using \eqref{eq:keyhole-factor}--\eqref{eq:sum-res},
\[
(1-e^{2\pi i/3})\,I
=
2\pi i\left(\frac{e^{i\pi/6}}{2i}-\frac12\right).
\]
Compute the right-hand side:
\[
2\pi i\cdot \frac{e^{i\pi/6}}{2i}=\pi e^{i\pi/6},
\qquad
2\pi i\cdot\left(-\frac12\right)=-\pi i,
\]
so
\[
(1-e^{2\pi i/3})\,I = \pi e^{i\pi/6}-\pi i.
\]
Now simplify using $e^{2\pi i/3}=-\frac12+i\frac{\sqrt3}{2}$ and $e^{i\pi/6}=\frac{\sqrt3}{2}+i\frac12$.
Then
\[
1-e^{2\pi i/3}=\frac32-i\frac{\sqrt3}{2},
\qquad
\pi e^{i\pi/6}-\pi i=\pi\left(\frac{\sqrt3}{2}-i\frac12\right).
\]
Thus
\[
I=\pi\,\frac{\frac{\sqrt3}{2}-i\frac12}{\frac32-i\frac{\sqrt3}{2}}.
\]
Multiply numerator and denominator by the conjugate $\frac32+i\frac{\sqrt3}{2}$:
\[
I=\pi\,
\frac{\left(\frac{\sqrt3}{2}-i\frac12\right)\left(\frac32+i\frac{\sqrt3}{2}\right)}
{\left(\frac32\right)^2+\left(\frac{\sqrt3}{2}\right)^2}.
\]
The denominator equals $\frac{9}{4}+\frac{3}{4}=3$. The numerator expands to
\[
\left(\frac{\sqrt3}{2}\cdot\frac32 + \frac{\sqrt3}{2}\cdot i\frac{\sqrt3}{2}\right)
+
\left(-i\frac12\cdot\frac32 - i\frac12\cdot i\frac{\sqrt3}{2}\right)
=
\frac{3\sqrt3}{4} + i\frac{3}{4} - i\frac{3}{4} + \frac{\sqrt3}{4}
=
\sqrt3.
\]
Hence
\[
I=\pi\cdot \frac{\sqrt3}{3}=\frac{\pi}{\sqrt3}.
\]

\subsection{Zeros, Poles, and the Argument Principle}
\index{argument principle}\index{zeros and poles}
\label{subsec:zeros-poles-argument-principle}

In this subsection we summarize the structure of zeros and poles of holomorphic
and meromorphic functions. We discuss (i) the order of a zero, (ii) integral formulas
for computing multiplicities, (iii) the argument principle, and (iv) geometric intuition
behind local covering behavior.

\subsubsection{Order of a Zero}

Let $\Omega\subset\mathbb{C}$ be a domain and let $f\in\Hol(\Omega)$.
Assume $f\not\equiv 0$ and $p\in\Omega$ satisfies $f(p)=0$.
Then there exists a unique integer $n\ge 1$ and a function $g\in\Hol(B_r(p))$
for some $r>0$ such that
\[
f(z)=(z-p)^n g(z),\qquad g(p)\neq 0.
\]
The integer $n$ is called the \emph{order} (or \emph{multiplicity}) of the zero at $p$,
and we write
\[
\ord_p(f)=n.
\]
In particular, zeros of a nonzero holomorphic function are isolated; hence $f$ has only
finitely many zeros in any compact set.

\begin{remark}
	If $g(p)\neq 0$ then $g$ is nonvanishing on a sufficiently small neighborhood of $p$.
	Thus locally, all the vanishing of $f$ near $p$ is captured by the factor $(z-p)^n$.
\end{remark}

\subsubsection{Computing the Order via a Contour Integral}

Let $r>0$ be small enough so that $f$ has no zeros on the circle
\[
\partial B_r(p)=\{z:|z-p|=r\}.
\]
Then
\begin{equation}\label{eq:order-by-log-derivative}
	\frac{1}{2\pi i}\int_{\partial B_r(p)} \frac{f'(z)}{f(z)}\,dz=\ord_p(f).
\end{equation}

\begin{proof}
	Write $f(z)=(z-p)^n g(z)$ with $g\in\Hol(B_r(p))$ and $g(p)\neq 0$.
	Then $g$ has no zeros on $\overline{B_r(p)}$ for $r$ small, so $g'/g$ is holomorphic
	on $B_r(p)$. Differentiate and divide:
	\[
	\frac{f'(z)}{f(z)}=\frac{n}{z-p}+\frac{g'(z)}{g(z)}.
	\]
	Integrate over $\partial B_r(p)$:
	\[
	\frac{1}{2\pi i}\int_{\partial B_r(p)}\frac{f'}{f}\,dz
	=
	\frac{1}{2\pi i}\int_{\partial B_r(p)}\frac{n}{z-p}\,dz
	+
	\frac{1}{2\pi i}\int_{\partial B_r(p)}\frac{g'}{g}\,dz.
	\]
	The second integral is $0$ because $g'/g$ is holomorphic on $B_r(p)$.
	The first integral equals $n$ because $\int_{\partial B_r(p)}\frac{dz}{z-p}=2\pi i$.
\end{proof}

\subsubsection{Argument Principle (Holomorphic Case)}

Let $f\in\Hol(\Omega)$ be not identically zero, and let $\gamma$ be a positively
oriented piecewise $C^1$ simple closed curve whose image and interior lie in
$\Omega$. Assume $f$ has no zeros on $\gamma$.
Let $z_1,\dots,z_k$ be the zeros of $f$ in the interior of $\gamma$, and let
$m_j=\ord_{z_j}(f)$. Then
\begin{equation}\label{eq:arg-principle-holo}
	\frac{1}{2\pi i}\int_{\gamma}\frac{f'(z)}{f(z)}\,dz=\sum_{j=1}^{k} m_j.
\end{equation}

\begin{proof}
	The function $f'/f$ is meromorphic, and at a zero $z_j$ of order $m_j$ we have
	$f(z)=(z-z_j)^{m_j}g(z)$ with $g(z_j)\neq 0$, hence
	\[
	\frac{f'(z)}{f(z)}=\frac{m_j}{z-z_j}+\frac{g'(z)}{g(z)}.
	\]
	So $\frac{f'}{f}$ has a simple pole at $z_j$ with residue $m_j$.
	By the residue theorem,
	\[
	\int_{\gamma}\frac{f'}{f}\,dz = 2\pi i\sum_{j=1}^{k}\Res\!\left(\frac{f'}{f},z_j\right)
	=2\pi i\sum_{j=1}^{k} m_j,
	\]
	which gives \eqref{eq:arg-principle-holo}.
\end{proof}

\begin{remark}[Geometric meaning: change of argument]
	If $\gamma$ is parametrized by $z=\gamma(t)$ and $f(\gamma(t))\neq 0$, then
	$f(\gamma(t))$ traces a closed curve in $\mathbb{C}^*$. The integral
	$\frac{1}{2\pi i}\int_{\gamma} f'/f$ equals the winding number of $f(\gamma)$
	around $0$, i.e.\ the total change of $\arg f(\gamma(t))$ divided by $2\pi$.
\end{remark}

\subsubsection{Argument Principle (Meromorphic Case)}

Let $f\in M(\Omega)$ be meromorphic and not identically zero, and let $\gamma$
be a positively oriented piecewise $C^1$ simple closed curve whose image and
interior lie in $\Omega$ and whose image contains no zeros or poles of $f$.
Let $z_1,\dots,z_k$ be the zeros of $f$ inside $\gamma$ with multiplicities $m_j$,
and let $p_1,\dots,p_\ell$ be the poles of $f$ inside $\gamma$ with orders $\mu_i$.
Then
\begin{equation}\label{eq:arg-principle-meromorphic}
	\frac{1}{2\pi i}\int_{\gamma}\frac{f'(z)}{f(z)}\,dz
	=
	\sum_{j=1}^{k} m_j
	-
	\sum_{i=1}^{\ell} \mu_i.
\end{equation}

\begin{proof}
	The logarithmic derivative $f'/f$ is meromorphic. At a zero $z_j$ of order $m_j$,
	it has a simple pole with residue $m_j$ as above. At a pole $p_i$ of order $\mu_i$,
	write $f(z)=(z-p_i)^{-\mu_i}h(z)$ with $h(p_i)\neq 0$, so
	\[
	\frac{f'(z)}{f(z)}=-\frac{\mu_i}{z-p_i}+\frac{h'(z)}{h(z)},
	\]
	hence the residue is $-\mu_i$.
	Apply the residue theorem to $\frac{f'}{f}$ on $\gamma$ to obtain \eqref{eq:arg-principle-meromorphic}.
\end{proof}

\subsubsection{Local geometric intuition}

If $f$ is holomorphic and $\ord_p(f)=n$, then in a small punctured neighborhood of $p$,
the map behaves like $z\mapsto (z-p)^n$ up to multiplication by a nonzero holomorphic factor.
Geometrically, small circles around $p$ are mapped to loops winding $n$ times around $0$.
This is the local ``$n$-to-$1$ covering'' picture behind why $\frac{1}{2\pi i}\int f'/f$
counts multiplicities.

\subsubsection{Exercises}

\begin{exercise}
	Let $f(z) = (z-1)^3(z+2)^2$.
	\begin{enumerate}
		\item Find $\ord_{1}(f)$ and $\ord_{-2}(f)$.
		\item Compute
		\[
		\frac{1}{2\pi i}\int_{|z-1|=r}\frac{f'(z)}{f(z)}\,dz
		\]
		for sufficiently small $r>0$, and interpret the result via \eqref{eq:order-by-log-derivative}.
	\end{enumerate}
\end{exercise}

\begin{exercise}
	Let $f(z)=z^5-2z^3+z$.
	\begin{enumerate}
		\item Factor $f$ and list all zeros with multiplicities.
		\item Use \eqref{eq:arg-principle-holo} to compute
		\[
		\frac{1}{2\pi i}\int_{|z|=2} \frac{f'(z)}{f(z)}\,dz,
		\]
		and check that it matches the total multiplicity of zeros inside $|z|<2$.
	\end{enumerate}
\end{exercise}

\begin{exercise}
	Let $f(z)=\dfrac{z^2+1}{z^2(z-3)}$.
	Compute
	\[
	\frac{1}{2\pi i}\int_{|z|=2}\frac{f'(z)}{f(z)}\,dz.
	\]
	Identify explicitly the zeros and poles of $f$ inside $|z|<2$, and verify \eqref{eq:arg-principle-meromorphic}.
\end{exercise}

\par\medskip
\begin{exercise}
	Let $p(z)=z^4+3z^2+1$.
	Set $w=z^2$ and locate all four zeros explicitly.  Verify that they lie on the
	imaginary axis, so the open right half-plane contains no zero.  Then explain why a
	right-half-disk argument-principle contour cannot be used without indenting around
	the zeros on its boundary.  Carry out that indentation as an optional check.
\end{exercise}

\subsection{Rouché's Theorem and Zero Counting by Boundary Data}
\index{Rouch\'e's theorem}
\label{subsec:rouche}

\subsubsection{Motivation: zeros are stable under small boundary perturbations}

The argument principle shows that the integer
\[
N_\Gamma(f):=\frac{1}{2\pi i}\int_\Gamma \frac{f'(z)}{f(z)}\,dz
\]
counts the number of zeros of a holomorphic function $f$ inside $\Gamma$ (with multiplicity),
as long as $f$ has no zeros on $\Gamma$.
Rouch\'e's theorem makes this comparison precise: if two holomorphic functions are close
on the boundary, then they must have the same number of zeros inside.

\subsubsection{Rouché's Theorem}

\begin{theorem}[Rouché's Theorem]
	Let $\Omega\subset\mathbb{C}$ be a domain, and let $f,g\in\mathcal{H}(\Omega)$.
	Let $\Gamma$ be a positively oriented simple closed contour with
	$\Gamma$ and its interior contained in $\Omega$.
	If
	\[
	|g(z)|<|f(z)| \qquad \text{for all } z\in\Gamma,
	\]
	then $f$ and $f+g$ have the same number of zeros (counted with multiplicity)
	in the interior of $\Gamma$.
\end{theorem}

\begin{proof}
	For $t\in[0,1]$ define $h_t(z)=f(z)+t\,g(z)$.
	On $\Gamma$,
	\[
	|h_t(z)-f(z)|=t|g(z)|\le |g(z)|<|f(z)|.
	\]
	By the reverse triangle inequality,
	\[
	|h_t(z)|\ge |\,|f(z)|-|h_t(z)-f(z)|\,|>0,
	\]
	so $h_t$ has no zeros on $\Gamma$ for all $t\in[0,1]$.
	Hence $h_t'/h_t$ is holomorphic on a neighborhood of $\Gamma$ and
	\[
	N_\Gamma(h_t)=\frac{1}{2\pi i}\int_\Gamma \frac{h_t'(z)}{h_t(z)}\,dz
	\]
	is integer-valued and continuous in $t$, thus constant in $t$.
	Evaluating at $t=0$ and $t=1$ gives $N_\Gamma(f)=N_\Gamma(f+g)$.
\end{proof}

\begin{remark}
	The hypothesis $|g|<|f|$ is used only on the boundary.
	No comparison is needed in the interior; the conclusion is topological
	(stability of the winding number of $h_t(\Gamma)$ about the origin).
\end{remark}

\subsubsection{Example 1: counting zeros via a dominant term}

\begin{example}[All zeros of $z^4+4z+1$ lie in $|z|<2$]
	Let $p(z)=z^4+4z+1$.
	On the circle $|z|=2$ we estimate
	\[
	|z^4|=16,\qquad |4z+1|\le 8+1=9,
	\]
	so $|4z+1|<|z^4|$ on $|z|=2$.
	By Rouch\'e with $f(z)=z^4$ and $g(z)=4z+1$, the functions $p(z)=f(z)+g(z)$
	and $z^4$ have the same number of zeros in $|z|<2$.
	Since $z^4$ has four zeros (counted with multiplicity) at $0$, $p$ also has
	four zeros in $|z|<2$. In particular, \emph{all} zeros of $p$ lie in $|z|<2$.
\end{example}

\subsubsection{Example 2: a root count on a small circle}

\begin{example}[Zeros of $z^3+z^2$ in $|z|<\tfrac12$]
	Let $f(z)=z^3$ and $g(z)=z^2$.
	On $|z|=\tfrac12$,
	\[
	|f(z)|=|z|^3=\frac{1}{8},\qquad |g(z)|=|z|^2=\frac{1}{4}.
	\]
	Here $|g|>|f|$, so Rouch\'e \emph{does not} apply in the direction $(f,g)$.
	
	Instead write $z^3+z^2=z^2(z+1)$ and count directly:
	the zeros are $z=0$ (multiplicity $2$) and $z=-1$ (multiplicity $1$).
	Inside $|z|<\tfrac12$ only $z=0$ lies inside, so there are exactly $2$ zeros
	(counted with multiplicity) in $|z|<\tfrac12$.
\end{example}

\begin{remark}
	This example is pedagogically useful: it shows that one should not apply
	Rouch\'e blindly. If the inequality fails, either choose a different contour
	or compare different pairs (e.g.\ compare $z^2$ vs $z^3$ on a different radius).
\end{remark}

\subsubsection{Example 3: using Rouché to show ``exactly one zero''}

\begin{example}[Exactly one zero in a disk]
	Let $q(z)=z^5+10z-3$ and consider $|z|=1$.
	On $|z|=1$ we have
	\[
	|10z-3|\ge 10-3=7,\qquad |z^5|=1.
	\]
	Thus $|z^5|<|10z-3|$ on $|z|=1$.
	By Rouch\'e, $q(z)=(10z-3)+z^5$ has the same number of zeros in $|z|<1$
	as $10z-3$, which has exactly one zero at $z=3/10$.
	Hence $q$ has exactly one zero in the unit disk.
\end{example}

\subsubsection{Exercises}

\begin{exercise}[Rouché: choose the right dominant term]
	Let $p(z)=z^6+3z^2+2$.
	\begin{enumerate}
		\item Show that on $|z|=2$ one has $|3z^2+2|<|z^6|$.
		\item Deduce that $p$ has six zeros in $|z|<2$.
	\end{enumerate}
\end{exercise}

\begin{exercise}[Rouché on an annulus boundary]
	Let $p(z)=z^4+z^3+1$.
	\begin{enumerate}
		\item Show that $p$ has no zeros on $|z|=2$ and conclude that it has four zeros in $|z|<2$.
		\item Show that $p$ has no zeros on $|z|=\tfrac12$ and conclude that it has no zeros in $|z|<\tfrac12$.
		\item Conclude that all zeros lie in $\tfrac12<|z|<2$.
	\end{enumerate}
\end{exercise}

\begin{exercise}[Counting zeros in the right half-plane (outline)]
	Let $p(z)=z^4+1$.
	Verify first that $p$ has no zeros on the imaginary axis. Use the argument
	principle on a sufficiently large right half-disk to determine how many zeros
	lie in the open right half-plane, and check the answer against the explicit
	fourth roots of $-1$.
\end{exercise}

\begin{exercise}[A ``one zero'' disk problem]
	Show that $z^7-8z+1$ has exactly one zero in $|z|<1$.
	(Choose an appropriate dominant term on $|z|=1$.)
\end{exercise}

\ReflectionMinseungKwak
\EssayJunhyukBoo

\section{Advanced Topics in Complex Analysis}
\index{complex analysis!advanced topics}

The basic theorems established above turn boundary information into
statements about the interior of a domain.  We now use that principle in
two directions.  Maximum principles and Schwarz--Pick give rigidity:
they constrain what a holomorphic map can do.  Normal families turn
uniform bounds into convergent subsequences, making it possible to
construct maps, most notably in the Riemann mapping theorem.  The
intervening results on entire functions, harmonic functions, and
analytic continuation supply examples and tools for these two themes.
As the arguments become more global, it is useful to keep asking which
hypothesis supplies compactness and which supplies uniqueness.

\subsection{Maximum and Minimum Modulus Principles and Consequences}
\label{subsec:max-min-modulus}

\subsubsection{Maximum modulus principle (local and global forms)}

\begin{theorem}[Maximum Modulus Principle]
	Let $\Omega\subset\mathbb{C}$ be a connected open set and let $f\in\mathcal{H}(\Omega)$.
	If $|f|$ attains a local maximum at some $z_0\in\Omega$, then $f$ is constant on $\Omega$.
\end{theorem}

\begin{proof}
	Choose a disk $B(z_0,r)\Subset\Omega$ on which
	$|f(z)|\le|f(z_0)|$. If $f(z_0)=0$, this inequality makes $f$ vanish on
	the disk, and the identity theorem finishes the proof. Otherwise set
	$h=f/f(z_0)$. Then $|h|\le1$ on the disk and $h(z_0)=1$. For every
	$0<\rho<r$, Cauchy's mean-value formula gives
	\[
	1=\Re h(z_0)=\frac1{2\pi}\int_0^{2\pi}
	\Re h(z_0+\rho e^{it})\,dt.
	\]
	The continuous integrand is everywhere at most $1$, so it must equal
	$1$ on the whole circle. Since $|h|\le1$, this forces $h=1$ there.
	The identity theorem now gives $h\equiv1$ on $\Omega$, and hence $f$ is
	constant.  This proof uses only Cauchy's formula and does not presuppose
	the open mapping theorem proved later.
\end{proof}

\begin{corollary}[Boundary maximum on a bounded domain]
	Let $\Omega$ be bounded and $f\in\mathcal{H}(\Omega)\cap C^0(\overline{\Omega})$.
	Then
	\[
	\max_{\overline{\Omega}}|f|=\max_{\partial\Omega}|f|.
	\]
\end{corollary}

\begin{proof}
	The continuous function $|f|$ attains its maximum on compact $\overline{\Omega}$.
	If a maximizer lies in the interior, the maximum modulus principle makes
	$f$ constant on $\Omega$, hence on its closure by continuity; the same
	value is then attained on $\partial\Omega$. Otherwise a maximizer
	already lies on $\partial\Omega$.
\end{proof}

\subsubsection{Minimum modulus principle}

\begin{theorem}[Minimum Modulus Principle]
	Let $\Omega$ be a connected open set and let $f\in\mathcal{H}(\Omega)$ be nonconstant.
	If $f$ has no zeros in $\Omega$, then $|f|$ has no local minimum in the interior of $\Omega$.
\end{theorem}

\begin{proof}
	If $f$ has no zeros, then $1/f\in\mathcal{H}(\Omega)$.
	A local minimum of $|f|$ at $z_0$ is equivalent to a local maximum of $|1/f|$ at $z_0$.
	By the maximum modulus principle, $1/f$ must be constant, hence $f$ is constant,
	contradicting the hypothesis.
\end{proof}

\begin{remark}
	If $f$ has zeros, then $\min_{\overline{\Omega}}|f|=0$ can occur (at zeros),
	so the zero-free hypothesis is essential.
\end{remark}

\subsubsection{Worked estimates}

\begin{example}[Maximum modulus for $e^z$ on a disk]
	On $\Omega=B_R(0)$ let $f(z)=e^z$. Then $|f(z)|=e^{\Re z}$.
	On $|z|\le R$, the inequality $\Re z\le |z|\le R$ yields
	\[
	|e^z|\le e^R,
	\]
	with equality only at $z=R$.
	Thus $\max_{|z|\le R}|e^z|=e^R$ and it occurs only on the boundary.
\end{example}

\begin{example}[Minimum modulus for a zero-free function]
	On $\Omega=B_1(0)$ let $f(z)=e^z$.
	Since $e^z\neq 0$ on $\mathbb{C}$, the minimum modulus principle implies
	that $|e^z|$ cannot have an interior local minimum.
	Indeed, $|e^z|=e^{\Re z}$ has infimum $e^{-1}$ on the open disk,
	approached as $z\to-1$ from inside; on the closed disk that value is
	attained at the boundary point $-1$.
\end{example}

\begin{example}[A function with zeros: minimum equals $0$]
	Let $f(z)=z^2+1$ on $|z|<2$.
	The zeros $z=\pm i$ lie inside the disk, so
	\[
	\min_{|z|\le 2}|z^2+1|=0,
	\]
	and the minimum is attained exactly at $z=\pm i$.
\end{example}

\subsubsection{Strong maximum principle (identity-theorem viewpoint)}

\begin{theorem}[Strong Maximum Principle for holomorphic functions]
	Let $\Omega$ be a domain, $f\in\mathcal{H}(\Omega)$, and suppose that
	$|f(z)|$ is constant on a nonempty open set
	$U\subset\Omega$. Then $f$ is constant on $\Omega$.
\end{theorem}

\begin{proof}
	Every point of $U$ is an interior local maximum point of $|f|$.
	Applying the maximum modulus principle at any one such point shows
	directly that $f$ is constant on the connected domain $\Omega$.
\end{proof}

\subsubsection{A standard application: uniqueness for the Dirichlet problem}

\begin{example}[Uniqueness on the disk (harmonic maximum principle)]
	Let $u$ be harmonic on the unit disk and continuous on the closed disk.
	If $u$ attains a maximum at an interior point, then $u$ is constant.
	Indeed, because the disk is simply connected, take a harmonic conjugate
	$v$ on the entire disk and set
	$F=u+iv$ holomorphic. Then $e^{F}$ is holomorphic and
	\[
	|e^{F}|=e^{u}.
	\]
	A maximum of $u$ yields a maximum of $|e^{F}|$, hence $e^{F}$ is constant by the
	maximum modulus principle, so $u$ is constant.
	This implies uniqueness of solutions to the Dirichlet problem on the disk.
\end{example}

\subsubsection{Exercises}

\begin{exercise}[Liouville via maximum modulus on large disks]
	Let $f$ be entire and assume there exist $R,M>0$ such that $|f(z)|\le M$ for all $|z|\ge R$.
	Show that $f$ is constant.
	(Hint: apply the maximum modulus principle to $f$ on $|z|\le n$ and let $n\to\infty$.)
\end{exercise}

\begin{exercise}[Minimum modulus: interior minimum forces a zero]
	Let $\Omega$ be bounded, $f\in\mathcal{H}(\Omega)\cap C^0(\overline{\Omega})$,
	and assume $|f(z_0)|=\min_{\overline{\Omega}}|f|$ for some $z_0\in\Omega$.
	Show that either $f(z_0)=0$ or $f$ is constant.
\end{exercise}

\begin{exercise}[Coefficient bounds from boundary control]
	Let $f(z)=\sum_{n=0}^\infty a_n z^n$ be holomorphic in $|z|<1$ and continuous on $\overline{\mathbb{D}}$.
	Assume $\max_{|z|\le 1}|f(z)|\le 1$.
	\begin{enumerate}
		\item Use Cauchy's integral formula to show $|a_n|\le 1$ for all $n\ge 0$.
		\item Show that if $|a_k|=1$ for some $k\ge 1$, then
		$f(z)=e^{i\theta}z^k$ for some $\theta\in\mathbb{R}$.
		(Hint: inspect equality in the boundary coefficient integral.)
	\end{enumerate}
\end{exercise}

\begin{exercise}[A holomorphic function cannot have both a maximum and minimum inside]
	Let $f\in\mathcal{H}(\Omega)$ on a connected domain $\Omega$.
	Show that if $|f|$ attains both a local maximum and a local minimum at interior points,
	then $f$ is constant.
\end{exercise}

\begin{exercise}[Harmonic maximum principle from holomorphic maximum modulus]
	Let $u$ be harmonic on a bounded domain $\Omega$ and continuous on $\overline{\Omega}$.
	Assume $\Omega$ is simply connected.
	\begin{enumerate}
		\item Let $v$ be a harmonic conjugate of $u$ and set $F=u+iv$.
		Show that $|e^{F}|=e^{u}$.
		\item Deduce the maximum principle for $u$: $\max_{\overline{\Omega}}u=\max_{\partial\Omega}u$.
	\end{enumerate}
\end{exercise}

\begin{exercise}[Strong maximum principle via identity theorem]
	Let $\Omega$ be a domain, $f\in\mathcal{H}(\Omega)$, and assume there exists a sequence $z_n\to z_0\in\Omega$
	such that $|f(z_n)|=\sup_\Omega |f|$ for all $n$.
	Show that $f$ is constant on $\Omega$.
\end{exercise}

\subsection{Schwarz Lemma, Schwarz--Pick Lemma, and the Poincar\'e Metric on the Disk}
\index{Schwarz lemma}\index{Schwarz--Pick lemma}\index{Poincar\'e metric}
\label{subsec:schwarz-pick-poincare}

\paragraph{Goal.}
We develop the Schwarz lemma and its M\"obius-invariant extension, the Schwarz--Pick lemma.
We then introduce the Poincar\'e metric on $\mathbb{D}$ and explain how Schwarz--Pick is exactly
the statement that holomorphic self-maps of the disk are non-expanding for the hyperbolic metric.
As an application we show that disk automorphisms act by hyperbolic isometries and
record the explicit structure of $\Aut(\mathbb{D})$.

\subsubsection{Preliminaries: the unit disk and its automorphisms}

\begin{definition}[Unit disk]
	\[
	\mathbb{D}:=\{z\in\mathbb{C}:|z|<1\}.
	\]
\end{definition}

\begin{definition}[M\"obius automorphisms of the disk]
	For $a\in\mathbb{D}$ define
	\[
	\varphi_a(z):=\frac{a-z}{1-\overline{a}z}.
	\]
	Then $\varphi_a:\mathbb{D}\to\mathbb{D}$ is biholomorphic, satisfies $\varphi_a(a)=0$ and
	$\varphi_a(0)=a$, and $\varphi_a^{-1}=\varphi_a$.
\end{definition}

\begin{lemma}[Basic identities for $\varphi_a$]
	For $a\in\mathbb{D}$ and $z\in\mathbb{D}$,
	\begin{align}
		1-|\varphi_a(z)|^2
		&=\frac{(1-|a|^2)(1-|z|^2)}{|1-\overline{a}z|^2}, \label{eq:phi-identity}\\[2mm]
		\varphi_a'(z)
		&=-\,\frac{1-|a|^2}{(1-\overline{a}z)^2},
		\qquad
		|\varphi_a'(z)|
		=\frac{1-|a|^2}{|1-\overline{a}z|^2}. \label{eq:phi-derivative}
	\end{align}
\end{lemma}

\begin{proof}
	First compute
	\[
	|\varphi_a(z)|^2=\frac{|a-z|^2}{|1-\overline{a}z|^2}.
	\]
	Then
	\begin{align*}
		1-|\varphi_a(z)|^2
		&=\frac{|1-\overline{a}z|^2-|a-z|^2}{|1-\overline{a}z|^2}.
	\end{align*}
	Expand the numerator:
	\begin{align*}
		|1-\overline{a}z|^2
		&=(1-\overline{a}z)(1-a\overline{z})
		=1-a\overline{z}-\overline{a}z+|a|^2|z|^2,\\
		|a-z|^2
		&=(a-z)(\overline{a}-\overline{z})
		=|a|^2-a\overline{z}-\overline{a}z+|z|^2.
	\end{align*}
	Subtracting gives
	\[
	|1-\overline{a}z|^2-|a-z|^2
	=(1-|a|^2)(1-|z|^2),
	\]
	which proves \eqref{eq:phi-identity}.
	
	For the derivative,
	\[
	\varphi_a(z)=\frac{a-z}{1-\overline{a}z}
	\quad\mathrel{\text{and hence}}\quad
	\varphi_a'(z)=\frac{-(1-\overline{a}z)-(a-z)(-\overline{a})}{(1-\overline{a}z)^2}
	=-\,\frac{1-|a|^2}{(1-\overline{a}z)^2}.
	\]
	Taking absolute values yields \eqref{eq:phi-derivative}.
\end{proof}

\subsubsection{Schwarz lemma}

\begin{theorem}[Schwarz lemma]
	Let $f:\mathbb{D}\to\mathbb{D}$ be holomorphic and assume $f(0)=0$.
	Then
	\[
	|f(z)|\le |z|\quad(z\in\mathbb{D}),
	\qquad\text{and}\qquad
	|f'(0)|\le 1.
	\]
	Moreover, if equality holds at some $z\neq 0$ (i.e.\ $|f(z)|=|z|$) or if $|f'(0)|=1$,
	then there exists $\theta\in\mathbb{R}$ such that $f(z)=e^{i\theta}z$ for all $z$.
\end{theorem}

\begin{proof}
	Define
	\[
	g(z):=
	\begin{cases}
		\dfrac{f(z)}{z}, & z\neq 0,\\[2mm]
		f'(0), & z=0.
	\end{cases}
	\]
	Since $f(0)=0$, $f(z)=z\cdot g(z)$ for $z\neq 0$. We show $g$ is holomorphic on $\mathbb{D}$.
	Indeed, $f$ has a power series $f(z)=a_1 z+a_2 z^2+\cdots$ near $0$, hence
	\[
	\frac{f(z)}{z}=a_1+a_2 z+a_3 z^2+\cdots,
	\]
	which is holomorphic near $0$ and equals $a_1=f'(0)$ at $0$. Thus $g\in\mathcal{H}(\mathbb{D})$.
	
	Now fix $0<r<1$. On the circle $|z|=r$ we have $|f(z)|<1$, hence
	\[
	|g(z)|=\left|\frac{f(z)}{z}\right|
	\le \frac{|f(z)|}{r}
	<\frac{1}{r}.
	\]
	This is not yet the desired bound. Instead we use the maximum modulus principle on $g$ as follows.
	
	For $|z|=r$,
	\[
	|g(z)|
	=\left|\frac{f(z)}{z}\right|
	\le \frac{|f(z)|}{|z|}
	\le \frac{1}{r}.
	\]
	Therefore $\max_{|z|=r}|g(z)|\le 1/r$.
	Letting $r\uparrow 1$, we obtain $\sup_{|z|<1}|g(z)|\le 1$.
	More explicitly: for any fixed $z$ with $|z|<1$, choose $r$ with $|z|<r<1$.
	Then by the maximum modulus principle applied to $g$ on $|w|\le r$,
	\[
	|g(z)|\le \max_{|w|=r}|g(w)|
	=\max_{|w|=r}\left|\frac{f(w)}{w}\right|
	\le \frac{1}{r}.
	\]
	Now let $r\uparrow 1$ to conclude $|g(z)|\le 1$ for all $z\in\mathbb{D}$.
	
	Hence for $z\neq 0$,
	\[
	|f(z)|=|z|\,|g(z)|\le |z|.
	\]
	Also $|f'(0)|=|g(0)|\le 1$.
	
	For the rigidity: if $|f(z_0)|=|z_0|$ for some $z_0\neq 0$, then
	$|g(z_0)|=1$. Since $|g|\le 1$ on $\mathbb{D}$ and $g$ attains its maximum modulus
	at an interior point, the maximum modulus principle forces $g$ to be constant:
	$g(z)\equiv e^{i\theta}$ for some $\theta$. Thus $f(z)=e^{i\theta}z$.
	Similarly, if $|f'(0)|=|g(0)|=1$, then $|g|$ attains its maximum at $0$ and again
	$g$ is constant, giving the same conclusion.
\end{proof}

\subsubsection{Schwarz--Pick lemma}

\begin{theorem}[Schwarz--Pick lemma]
	Let $f:\mathbb{D}\to\mathbb{D}$ be holomorphic. Then for all $z,w\in\mathbb{D}$,
	\begin{equation}\label{eq:SP-crossratio}
		\left|\frac{f(z)-f(w)}{1-\overline{f(w)}\,f(z)}\right|
		\le
		\left|\frac{z-w}{1-\overline{w}\,z}\right|.
	\end{equation}
	Equivalently, for all $z\in\mathbb{D}$,
	\begin{equation}\label{eq:SP-derivative}
		\frac{|f'(z)|}{1-|f(z)|^2}\le \frac{1}{1-|z|^2}.
	\end{equation}
	Moreover, equality in \eqref{eq:SP-crossratio} for some distinct $z\neq w$
	(or equality in \eqref{eq:SP-derivative} at some point) holds iff $f\in\Aut(\mathbb{D})$.
\end{theorem}

\begin{proof}
	Fix $w\in\mathbb{D}$ and set $a:=f(w)\in\mathbb{D}$.
	Consider the holomorphic map
	\[
	F(z):=\bigl(\varphi_a\circ f\circ \varphi_w\bigr)(z).
	\]
	We have $\varphi_w(0)=w$ and $\varphi_a(a)=0$, so
	\[
	F(0)=\varphi_a(f(\varphi_w(0)))=\varphi_a(f(w))=\varphi_a(a)=0.
	\]
	Also $F:\mathbb{D}\to\mathbb{D}$ since each factor maps $\mathbb{D}$ to itself.
	By Schwarz lemma applied to $F$,
	\[
	|F(z)|\le |z|\qquad(z\in\mathbb{D}).
	\]
	Now substitute $z=\varphi_w(z_0)$ (i.e.\ evaluate at $z=\varphi_w(z_0)$ and rename):
	\[
	|\varphi_a(f(z_0))|
	=
	|(\varphi_a\circ f)(z_0)|
	=
	|F(\varphi_w(z_0))|
	\le
	|\varphi_w(z_0)|.
	\]
	Using the explicit formula for $\varphi_a$ and $\varphi_w$ gives exactly
	\eqref{eq:SP-crossratio}.
	
	To derive the differential form, set $z=w+h$ and let $h\to 0$ in \eqref{eq:SP-crossratio},
	or differentiate the identity obtained from Schwarz lemma for $F$ at $0$.
	From Schwarz lemma we also have $|F'(0)|\le 1$.
	Compute $F'(0)$ by the chain rule:
	\[
	F'(0)=\varphi_a'(f(w))\cdot f'(w)\cdot \varphi_w'(0).
	\]
	Using \eqref{eq:phi-derivative} and the fact $\varphi_w'(0)=-(1-|w|^2)$, we get
	\[
	|\varphi_a'(a)|
	=\frac{1-|a|^2}{|1-\overline{a}a|^2}
	=\frac{1-|a|^2}{(1-|a|^2)^2}
	=\frac{1}{1-|a|^2},
	\qquad
	|\varphi_w'(0)|=1-|w|^2.
	\]
	Hence
	\[
	|F'(0)|
	=
	\frac{1}{1-|a|^2}\,|f'(w)|\,(1-|w|^2)
	\le 1.
	\]
	Recalling $a=f(w)$ yields \eqref{eq:SP-derivative}.
	
	For rigidity: if equality holds in \eqref{eq:SP-crossratio} for some $z\neq w$,
	then equality holds in Schwarz lemma for $F$ at some point, forcing $F$ to be a rotation:
	$F(\zeta)=e^{i\theta}\zeta$. Then $f=\varphi_a^{-1}\circ F\circ \varphi_w^{-1}$
	is a disk automorphism. Conversely, if $f\in\Aut(\mathbb{D})$, then composing with
	automorphisms preserves equality, so \eqref{eq:SP-crossratio} and \eqref{eq:SP-derivative}
	hold with equality everywhere.
\end{proof}

\subsubsection{The Poincar\'e (hyperbolic) metric and curvature}

\begin{definition}[Poincar\'e metric on $\mathbb{D}$]
	The hyperbolic (Poincar\'e) metric on $\mathbb{D}$ is the conformal Riemannian metric
	\[
	ds_{\mathbb{D}}^2=\lambda_{\mathbb{D}}(z)^2\,|dz|^2,
	\qquad
	\lambda_{\mathbb{D}}(z):=\frac{2}{1-|z|^2}.
	\]
\end{definition}

\begin{remark}[Gauss curvature in conformal coordinates]
	If a metric on a planar domain is written as $ds^2=\lambda(z)^2 |dz|^2$ with $\lambda>0$,
	then its Gauss curvature is
	\[
	K(z)=-\frac{1}{\lambda(z)^2}\,\Delta\bigl(\log \lambda(z)\bigr),
	\qquad
	\Delta=\partial_x^2+\partial_y^2.
	\]
	For $\lambda_{\mathbb{D}}(z)=\dfrac{2}{1-|z|^2}$ one computes $K\equiv -1$.
\end{remark}

\subsubsection{Schwarz--Pick as hyperbolic non-expansion}

\begin{theorem}[Infinitesimal contraction]
	If $f:\mathbb{D}\to\mathbb{D}$ is holomorphic, then for all $z\in\mathbb{D}$,
	\[
	\lambda_{\mathbb{D}}(f(z))\,|f'(z)|
	\le
	\lambda_{\mathbb{D}}(z).
	\]
	Equivalently, $f$ does not increase hyperbolic lengths of tangent vectors.
\end{theorem}

\begin{proof}
	Rewrite \eqref{eq:SP-derivative} as
	\[
	\frac{2}{1-|f(z)|^2}\,|f'(z)|
	\le
	\frac{2}{1-|z|^2},
	\]
	which is exactly the stated inequality.
\end{proof}

\begin{corollary}[Distance contraction]
	Let $d_{\mathbb{D}}$ denote the path metric induced by $ds_{\mathbb{D}}^2$:
	\[
	d_{\mathbb{D}}(z,w)
	=
	\inf_{\gamma}\int_0^1 \lambda_{\mathbb{D}}(\gamma(t))\,|\gamma'(t)|\,dt,
	\]
	where the infimum is over all piecewise $C^1$ curves $\gamma$ joining $z$ to $w$.
	Then every holomorphic $f:\mathbb{D}\to\mathbb{D}$ satisfies
	\[
	d_{\mathbb{D}}(f(z),f(w))\le d_{\mathbb{D}}(z,w).
	\]
\end{corollary}

\begin{proof}
	Let $\gamma$ be any curve from $z$ to $w$. Then $f\circ\gamma$ joins $f(z)$ to $f(w)$, and
	\[
	\lambda_{\mathbb{D}}(f(\gamma(t)))\,|(f\circ\gamma)'(t)|
	=
	\lambda_{\mathbb{D}}(f(\gamma(t)))\,|f'(\gamma(t))|\,|\gamma'(t)|
	\le
	\lambda_{\mathbb{D}}(\gamma(t))\,|\gamma'(t)|
	\]
	by infinitesimal contraction. Integrating over $t\in[0,1]$ yields
	\[
	L_{\mathbb{D}}(f\circ\gamma)\le L_{\mathbb{D}}(\gamma).
	\]
	Taking infimum over $\gamma$ gives the distance inequality.
\end{proof}

\subsubsection{Disk automorphisms preserve the hyperbolic metric}

\begin{theorem}[Automorphisms preserve the hyperbolic metric]
	If $f\in\Aut(\mathbb{D})$, then
	\[
	\lambda_{\mathbb{D}}(f(z))\,|f'(z)|
	=
	\lambda_{\mathbb{D}}(z)
	\qquad(z\in\mathbb{D}),
	\]
	and consequently
	\[
	d_{\mathbb{D}}(f(z),f(w))=d_{\mathbb{D}}(z,w)
	\qquad(z,w\in\mathbb{D}).
	\]
\end{theorem}

\begin{proof}
	From Schwarz--Pick we already have
	$\lambda_{\mathbb{D}}(f(z))|f'(z)|\le \lambda_{\mathbb{D}}(z)$.
	Apply the same inequality to $f^{-1}\in\Aut(\mathbb{D})$ at the point $f(z)$:
	\[
	\lambda_{\mathbb{D}}(z)\,|(f^{-1})'(f(z))|
	\le
	\lambda_{\mathbb{D}}(f(z)).
	\]
	Using $(f^{-1})'(f(z))=1/f'(z)$ gives
	\[
	\lambda_{\mathbb{D}}(z)\le \lambda_{\mathbb{D}}(f(z))\,|f'(z)|.
	\]
	Combine both inequalities to get equality.
	
	For distances: for any curve $\gamma$, equality of densities implies
	$L_{\mathbb{D}}(f\circ\gamma)=L_{\mathbb{D}}(\gamma)$, hence the induced path
	distance is preserved.
\end{proof}

\begin{theorem}[Explicit form of $\Aut(\mathbb{D})$]
	Every disk automorphism is of the form
	\[
	f(z)=e^{i\theta}\,\frac{a-z}{1-\overline{a}z},
	\qquad a\in\mathbb{D},\ \theta\in\mathbb{R}.
	\]
	In particular, $\Aut(\mathbb{D})$ acts transitively on $\mathbb{D}$.
\end{theorem}

\begin{proof}
	Let $f\in\Aut(\mathbb{D})$ and set $a:=f(0)\in\mathbb{D}$.
	Consider
	\[
	F:=\varphi_a\circ f.
	\]
	Then $F$ is an automorphism and $F(0)=\varphi_a(f(0))=\varphi_a(a)=0$.
	Apply Schwarz lemma to $F$ and also to $F^{-1}$: since both fix $0$,
	\[
	|F(z)|\le |z|
	\quad\text{and}\quad
	|F^{-1}(z)|\le |z|.
	\]
	Apply the second inequality to $z=F(w)$ to get $|w|\le |F(w)|$. Hence $|F(w)|=|w|$ for all $w$.
	By the rigidity part of Schwarz lemma, $F$ must be a rotation:
	\[
	F(z)=e^{i\theta}z.
	\]
	Therefore
	\[
	f=\varphi_a^{-1}\circ F
	=
	\varphi_a\circ (e^{i\theta}\,\mathrm{id})
	\quad\mathrel{\text{and hence}}\quad
	f(z)=\varphi_a(e^{i\theta}z).
	\]
	Writing $b=e^{-i\theta}a$, this is
	\[
	f(z)=\frac{a-e^{i\theta}z}{1-\overline a e^{i\theta}z}
	=e^{i\theta}\frac{b-z}{1-\overline b z},
	\]
	which is the stated normal form after renaming $b$ as $a$.
	Transitivity follows by choosing $a$ so that $\varphi_a(w)=0$ and composing with a rotation.
\end{proof}

\begin{remark}[Group structure and stabilizer]
	The stabilizer of $0$ in $\Aut(\mathbb{D})$ consists of rotations $z\mapsto e^{i\theta}z$.
	Moreover, every automorphism decomposes as
	\[
	\Aut(\mathbb{D})
	=
	\{\text{translation to }0\}\cdot\{\text{rotation}\}\cdot\{\text{translation back}\}.
	\]
\end{remark}

\subsubsection{Examples and computations}

\begin{example}[Schwarz lemma from Schwarz--Pick]
	If $f(0)=0$, then \eqref{eq:SP-crossratio} with $w=0$ becomes
	\[
	|f(z)|\le |z|.
	\]
	The derivative form at $z=0$ gives $|f'(0)|\le 1$.
\end{example}

\begin{example}[A strict contraction]
	Let $f(z)=\tfrac12 z$. Then $|f'(z)|=\tfrac12$ and
	\[
	\lambda_{\mathbb{D}}(f(z))\,|f'(z)|
	=
	\frac{2}{1-|z|^2/4}\cdot \frac12
	<
	\frac{2}{1-|z|^2}
	=\lambda_{\mathbb{D}}(z),
	\]
	so $f$ strictly decreases hyperbolic lengths everywhere.
\end{example}

\begin{example}[Automorphism is an isometry: direct check for $\varphi_a$]
	Using \eqref{eq:phi-identity} and \eqref{eq:phi-derivative},
	The factors involving $1-|a|^2$ and $|1-\overline a z|^2$
	cancel when the two identities are substituted:
	\begin{align*}
	\lambda_{\mathbb{D}}(\varphi_a(z))\,|\varphi_a'(z)|
	&=\frac{2|1-\overline a z|^2}{(1-|a|^2)(1-|z|^2)}
	  \frac{1-|a|^2}{|1-\overline a z|^2}\\
	&=\frac{2}{1-|z|^2}=\lambda_{\mathbb{D}}(z).
	\end{align*}
\end{example}

\subsubsection{Exercises}

\begin{exercise}[Rigidity in Schwarz lemma]
	Let $f:\mathbb{D}\to\mathbb{D}$ be holomorphic with $f(0)=0$.
	Assume there exists $z_0\neq 0$ such that $|f(z_0)|=|z_0|$.
	Show that $f(z)=e^{i\theta}z$ for some $\theta\in\mathbb{R}$ by applying the
	maximum modulus principle to $g(z)=f(z)/z$.
\end{exercise}

\begin{exercise}[Derivative form from cross-ratio form]
	Assume \eqref{eq:SP-crossratio}. Fix $w$ and set $z=w+h$.
	Let $h\to 0$ and prove \eqref{eq:SP-derivative} by carefully expanding
	$f(w+h)=f(w)+f'(w)h+o(h)$ and comparing both sides.
\end{exercise}

\begin{exercise}[Hyperbolic distance to the origin]
	Show that
	\[
	d_{\mathbb{D}}(0,r)
	=
	\int_0^r \frac{2\,dt}{1-t^2}
	=
	\log\frac{1+r}{1-r},
	\qquad 0\le r<1.
	\]
\end{exercise}

\begin{exercise}[Distance formula via automorphisms]
	Let $z,w\in\mathbb{D}$ and set $\alpha=\varphi_w(z)$.
	Use invariance under automorphisms to show
	\[
	d_{\mathbb{D}}(z,w)=d_{\mathbb{D}}(\alpha,0)
	\]
	and deduce
	\[
	\tanh\!\Big(\frac12 d_{\mathbb{D}}(z,w)\Big)=|\alpha|
	=
	\left|\frac{z-w}{1-\overline{w}z}\right|.
	\]
\end{exercise}

\begin{exercise}[Transitivity of $\Aut(\mathbb{D})$]
	Given $u,v\in\mathbb{D}$, construct explicitly an $f\in\Aut(\mathbb{D})$ such that
	$f(u)=v$. (Hint: send $u$ to $0$ using $\varphi_u$, rotate, then send $0$ to $v$.)
\end{exercise}

\begin{exercise}[Fixed points and strict contraction]
	Let $f:\mathbb{D}\to\mathbb{D}$ be holomorphic with $f(0)=0$ and $|f'(0)|<1$.
	Show that $0$ is the unique fixed point of $f$ in $\mathbb{D}$ and that
	$f^{\circ n}(z)\to 0$ uniformly on compact sets. (Hint: for each $r<1$,
	apply the maximum principle to $f(z)/z$ on $|z|\le r$ to obtain a
	strict bound $q_r<1$, and note that Schwarz's lemma keeps every orbit
	starting in that disk inside it.)
\end{exercise}

\subsection{Invariant Metrics on the Unit Disk and the K\"ahler--Einstein Condition}
\index{invariant metric}\index{Carath\'eodory metric}\index{Kobayashi metric}\index{Bergman metric}\index{K\"ahler--Einstein metric}
\label{subsec:five-disk-metrics}

The metric $ds_{\mathbb D}^2=4|dz|^2/(1-|z|^2)^2$ arose above from
Schwarz--Pick.  There are other ways to obtain the same geometry: test
holomorphic maps \emph{out of} the disk, test maps \emph{into} it, or use
the Hilbert space of square-integrable holomorphic functions.  These
constructions matter on more general domains, where they need not agree
\cite{Kobayashi1998Hyperbolic,Krantz2013Bergman}.  We keep the factor $2$
in the Poincar\'e length $2|v|/(1-|z|^2)$ throughout; another common
convention omits it and gives curvature $-4$ instead of $-1$.

\begin{definition}[Invariant metric]
A family of tangent-vector lengths $F_\Omega(z;v)$ on complex domains
is called \emph{biholomorphically invariant} if every biholomorphism
$\Phi:\Omega\to\Omega'$ satisfies
\[
F_{\Omega'}\bigl(\Phi(z);\Phi'(z)v\bigr)=F_\Omega(z;v).
\]
For a single domain, this says in particular that its automorphisms
act by isometries.  An \emph{invariant metric} may be a genuine metric
or, for the extremal constructions below, only a pseudometric.  This
definition concerns behavior under biholomorphisms; the
K\"ahler--Einstein condition introduced below instead concerns
curvature.  The disk metric has both properties.
\end{definition}

\begin{definition}[Carath\'eodory and Kobayashi--Royden lengths]
For a complex domain $\Omega$, $z\in\Omega$, and $v\in\mathbb C$, set
\begin{align*}
C_\Omega(z;v)
 &=\sup_{f\in\operatorname{Hol}(\Omega,\mathbb D)}
       \frac{|f'(z)v|}{1-|f(z)|^2},\\
K_\Omega(z;v)
 &=\inf\left\{\alpha>0:\exists h\in\operatorname{Hol}(\mathbb D,\Omega),
       h(0)=z,\ h'(0)=v/\alpha\right\}.
\end{align*}
The first measures how strongly maps \emph{from} $\Omega$ can separate
nearby points; the second asks how large an analytic disk can realize a
given tangent vector.  They may be pseudometric lengths on a general
domain.  The corresponding distance constructions are especially useful
because holomorphic maps cannot increase them.
Indeed, composition in the defining supremum and infimum shows that
every holomorphic map contracts these lengths; applying this to a
biholomorphism and its inverse gives the equality in the definition
of invariance.  The Bergman metric is biholomorphically invariant as
well: under $\Phi$ the Bergman kernel transforms on the diagonal by
$B_\Omega(z,z)=|\Phi'(z)|^2B_{\Omega'}(\Phi(z),\Phi(z))$;
the extra $\log|\Phi'|^2$ is harmonic, so its mixed derivative
vanishes.  Unlike the two extremal lengths, the Bergman metric is not
in general contracted by arbitrary holomorphic maps.
\end{definition}

\begin{proposition}[Extremal maps on the disk]
For every $z\in\mathbb D$ and $v\in\mathbb C$,
\[
C_{\mathbb D}(z;v)=K_{\mathbb D}(z;v)
=\frac{|v|}{1-|z|^2}=\frac12\lambda_{\mathbb D}(z)|v|.
\]
Thus the real length elements $2C_{\mathbb D}(z;dz)$ and
$2K_{\mathbb D}(z;dz)$ are precisely the Poincar\'e length element.
\end{proposition}
\begin{proof}
Schwarz--Pick gives
$|f'(z)|/(1-|f(z)|^2)\leq 1/(1-|z|^2)$ for every
$f:\mathbb D\to\mathbb D$; the identity attains the bound in $C$.
For an admissible $h$ in the definition of $K$, Schwarz--Pick at $0$
gives $|h'(0)|\leq1-|h(0)|^2=1-|z|^2$.  Since
$h'(0)=v/\alpha$, this forces
$\alpha\geq |v|/(1-|z|^2)$.  If $v\ne0$, choose a disk
automorphism $h$ with $h(0)=z$ and derivative at $0$ of modulus
$1-|z|^2$, rotating its parameter so that $h'(0)$ has the direction
of $v$.  Then equality holds with
$\alpha=|v|/(1-|z|^2)$.  For $v=0$ both lengths vanish.
\end{proof}

\begin{definition}[Bergman metric and conventions]
Let $A^2(\Omega)$ be the holomorphic functions square-integrable for
Euclidean area $dA$, and let $B_\Omega(z,w)$ be their reproducing
kernel.  Wherever $B_\Omega(z,z)>0$, define the Hermitian coefficient
$g^B_{z\bar z}=\partial_z\partial_{\bar z}\log B_\Omega(z,z)$.
In this subsection the associated \emph{real} line element and
K\"ahler form are, respectively,
$ds_B^2=2g^B_{z\bar z}|dz|^2$ and
$\omega_B=i g^B_{z\bar z}\,dz\wedge d\bar z$.
Some authors instead call $g^B_{z\bar z}|dz|^2$ the Bergman line
element; their displayed metric is one half of ours.
\end{definition}

\begin{proposition}[Bergman kernel and metric of the disk]
For Euclidean area measure on $\mathbb D$,
\[
B_{\mathbb D}(z,w)=\frac{1}{\pi(1-z\bar w)^2},\qquad
g^B_{z\bar z}=\frac{2}{(1-|z|^2)^2},\qquad
ds_B^2=ds_{\mathbb D}^2.
\]
\end{proposition}
\begin{proof}
The monomials
$e_n(z)=\sqrt{(n+1)/\pi}\,z^n$ form an orthonormal basis of
$A^2(\mathbb D)$: polar integration gives
$\int_{\mathbb D}|z|^{2n}dA=\pi/(n+1)$, while angular integration
gives orthogonality.  Expanding any holomorphic $L^2$ function into
its Taylor series and integrating on concentric circles shows that
these monomials are complete.  Hence, for $|z\bar w|<1$,
\[
B_{\mathbb D}(z,w)=\sum_{n=0}^{\infty}e_n(z)\overline{e_n(w)}
=\frac1\pi\sum_{n=0}^{\infty}(n+1)(z\bar w)^n
=\frac{1}{\pi(1-z\bar w)^2}.
\]
On the diagonal, $\log B_{\mathbb D}(z,z)=-\log\pi
-2\log(1-z\bar z)$.  Since
$\partial_z\partial_{\bar z}[-\log(1-z\bar z)]
=(1-|z|^2)^{-2}$, the formulas for $g^B$ and $ds_B^2$ follow.
\end{proof}

\begin{definition}[K\"ahler--Einstein condition in complex dimension one]
For a positive Hermitian coefficient $g_{z\bar z}$, set
$\omega=i g_{z\bar z}\,dz\wedge d\bar z$ and define its Ricci form by
$\rho=-i\partial\bar\partial\log g_{z\bar z}$.
It is K\"ahler--Einstein if $\rho=\kappa\omega$ for a real constant
$\kappa$.  In complex dimension one every Hermitian metric is
K\"ahler: a $2$-form on a real surface is automatically closed.
\end{definition}

\begin{proposition}[The disk metric is K\"ahler--Einstein]
With the preceding normalization, $\omega_B$ satisfies
$\rho_B=-\omega_B$; equivalently, its Gaussian curvature is $-1$.
Thus the Poincar\'e real metric is also a complete negative
K\"ahler--Einstein metric on $\mathbb D$.
\end{proposition}
\begin{proof}
Write $g=g^B_{z\bar z}=2(1-|z|^2)^{-2}$.  Direct differentiation gives
\[
\partial_z\partial_{\bar z}\log g
=\frac{2}{(1-|z|^2)^2}=g.
\]
Therefore $\rho_B=-i g\,dz\wedge d\bar z=-\omega_B$.
The underlying real metric is the Poincar\'e metric, whose curvature
was computed above.  Completeness follows also from the distance
formula $d_{\mathbb D}(0,r)=\log((1+r)/(1-r))$, which diverges as
$r\to1$; disk automorphisms reduce any escaping path to this
boundary behavior.
\end{proof}

Their agreement is special to the disk, not a definition of these
metrics on arbitrary domains.  Carath\'eodory and Kobayashi lengths
encode extremal holomorphic maps, the Bergman metric encodes the
$L^2$ holomorphic function space, and the K\"ahler--Einstein equation
encodes curvature.  Their common formula here is a useful benchmark
when comparing more complicated domains.

\subsection{Montel's Theorem and Normal Families}
\index{Montel's theorem}\index{normal families}

We introduce normal families and prove Montel's theorem: local boundedness
implies subsequential convergence uniformly on compact subsets. This compactness
principle is used below in Hurwitz's theorem and the Riemann mapping theorem.

\subsubsection{Normal families and local boundedness}

\begin{definition}[Normal family]
	Let $\Omega\subset\mathbb{C}$ be a domain and let $\mathcal{F}\subset\mathcal{O}(\Omega)$.
	We say that $\mathcal{F}$ is a \emph{normal family} if every sequence
	$\{f_n\}\subset\mathcal{F}$ admits a subsequence that converges uniformly on
	compact subsets of $\Omega$ to a holomorphic limit function.
\end{definition}

\begin{definition}[Locally bounded family]
	A family $\mathcal{F}\subset\mathcal{O}(\Omega)$ is \emph{locally bounded} if for
	every compact set $K\Subset\Omega$ there exists a constant $M_K>0$ such that
	$\lvert f(z)\rvert\le M_K$ for all $z\in K$ and all $f\in\mathcal{F}$.
\end{definition}

\subsubsection{Cauchy estimates and equicontinuity}

\begin{lemma}[Uniform derivative bounds on compact sets]
	Let $\mathcal{F}\subset\mathcal{O}(\Omega)$ be locally bounded and let
	$K\Subset\Omega$.  Then there exists a constant $C_K>0$ such that
	$\lvert f'(z)\rvert\le C_K$ for all $z\in K$ and all $f\in\mathcal{F}$.
\end{lemma}

\begin{proof}
	Since $K$ is compact and $\Omega$ is open, there exists $r>0$ such that every disk
	$\overline{D(z,r)}$ with $z\in K$ is contained in $\Omega$.  By local boundedness,
	there exists $M>0$ such that $\lvert f(\zeta)\rvert\le M$ for all $\zeta$ with
	$\lvert \zeta-z\rvert\le r$ and all $f\in\mathcal{F}$.
	
	By Cauchy's integral formula for derivatives, for each $z\in K$ we have
	$f'(z)=\frac{1}{2\pi i}\int_{\lvert \zeta-z\rvert=r} \frac{f(\zeta)}{(\zeta-z)^2}\,d\zeta$.
	Taking absolute values yields
	$\lvert f'(z)\rvert\le \frac{M}{r}$.
	Setting $C_K=M/r$ proves the claim.
\end{proof}

\begin{lemma}[Equicontinuity on compact sets]
	If $\mathcal{F}\subset\mathcal{O}(\Omega)$ is locally bounded, then it is
	equicontinuous on every compact subset of $\Omega$.
\end{lemma}

\begin{proof}
	Fix $K\Subset\Omega$ and choose $\delta>0$ so that the closed
	$\delta$-neighborhood $K_\delta$ is compactly contained in $\Omega$.
	By the previous lemma, there exists $C>0$ such that
	$|f'|\le C$ on $K_\delta$ for all $f\in\mathcal F$.
	If $z,w\in K$ and $|z-w|<\delta$, the segment from $w$ to $z$
	lies in $K_\delta$, and therefore
	$f(z)-f(w)=\int_0^1 f'(w+t(z-w))(z-w)\,dt$,
	so $|f(z)-f(w)|\le C|z-w|$. Given $\varepsilon>0$, choosing
	$|z-w|<\min(\delta,\varepsilon/C)$ proves equicontinuity on $K$.
\end{proof}

\subsubsection{Montel's theorem}

\begin{theorem}[Montel]
	Let $\Omega\subset\mathbb{C}$ be a domain and let $\mathcal{F}\subset\mathcal{O}(\Omega)$
	be locally bounded.  Then $\mathcal{F}$ is a normal family.
\end{theorem}

\begin{proof}
	Fix a compact set $K\Subset\Omega$.  Local boundedness implies uniform boundedness
	of $\mathcal{F}$ on a slightly larger compact set, which in turn yields uniform
	bounds on derivatives and equicontinuity on $K$.  By the Arzel\`a--Ascoli theorem,
	every sequence in $\mathcal{F}$ has a subsequence converging uniformly on $K$.
	
	Applying this argument to an exhaustion of $\Omega$ by compact sets and using a
	diagonal subsequence, we obtain a subsequence converging uniformly on every
	compact subset of $\Omega$.  The limit function is holomorphic by Morera's theorem.
\end{proof}

\subsubsection{Examples}

\begin{example}[A normal family]
	The family $\{f_n(z)=z^n\}$ on the unit disk is normal.  Indeed,
	$\lvert z^n\rvert\le 1$ for all $z$ in the disk, so the family is locally bounded.
	Moreover, for every $r<1$ we have $\sup_{\lvert z\rvert\le r}\lvert z^n\rvert=r^n$,
	which converges to zero as $n\to\infty$.
\end{example}

\begin{example}[A non-normal family]
	The family $\{g_n(z)=nz\}$ on the unit disk is not normal.  Near the origin the
	family is not locally bounded. More directly, any subsequence converging uniformly
	on compact sets would converge at $z=1/2$, whereas $g_n(1/2)=n/2\to\infty$.
	Thus no subsequence has a finite holomorphic compact-open limit.
\end{example}

\begin{example}[Blaschke products]
	Finite Blaschke products satisfy $\lvert B(z)\rvert\le 1$ on the unit disk.
	Any family of such functions is locally bounded and therefore normal by Montel's
	theorem.
\end{example}

\subsubsection{Exercises}

\begin{exercise}
	Decide which of the families $\{z^n\}$ and $\{nz\}$ is normal on the unit disk and
	justify your answer using Montel's theorem.
\end{exercise}

\begin{exercise}
	Let $\mathcal{F}$ be a locally bounded family on a domain $\Omega$ and let
	$\{K_j\}$ be an exhaustion of $\Omega$ by compact sets.
	Explain carefully how to extract a diagonal subsequence converging uniformly on
	every $K_j$.
\end{exercise}

\begin{exercise}
	Assume $f_n\to f$ uniformly on compact subsets of $\Omega$, where each $f_n$ is
	holomorphic.  Prove that $f$ is holomorphic by applying Morera's theorem.
\end{exercise}

\subsection{Hurwitz's Theorem and the Open Mapping Theorem}

Uniform convergence on compact sets preserves zeros with their multiplicities,
except when the limit is identically zero. Hurwitz's theorem makes this precise;
the same zero-counting argument then yields the Open Mapping Theorem.

\subsubsection{Zeros and multiplicity}

\begin{definition}[Zero of multiplicity $m$]
	Let $f$ be holomorphic in a neighborhood of a point $z_0\in\mathbb{C}$.
	We say that $z_0$ is a zero of multiplicity $m$ if there exists a holomorphic
	function $g$ with $g(z_0)\neq 0$ such that
	$f(z)=(z-z_0)^m g(z)$ in a neighborhood of $z_0$.
\end{definition}

\subsubsection{Hurwitz's theorem}

\begin{theorem}[Hurwitz]
	Let $\Omega\subset\mathbb{C}$ be a domain and let $f_n\in\mathcal{O}(\Omega)$
	converge uniformly on compact subsets to a holomorphic function $f$.
	
	If each $f_n$ has no zeros in $\Omega$, then either $f$ has no zeros in $\Omega$
	or $f$ is identically zero.
	
	More generally, if $f\not\equiv0$ has a zero of multiplicity $m$ at a point $z_0$, then for
	every sufficiently small disk centered at $z_0$, all sufficiently large $n$
	satisfy that $f_n$ has exactly $m$ zeros in that disk, counted with multiplicity.
\end{theorem}

\begin{proof}
	Let $z_0\in\Omega$ be a zero of $f$ of multiplicity $m$.  Choose a radius $r>0$
	so small that $f(z)\neq 0$ for all $z$ with $\lvert z-z_0\rvert=r$.
	
	Uniform convergence of $f_n$ to $f$ on the closed disk implies that for all
	sufficiently large $n$ we have
	$\lvert f_n(z)-f(z)\rvert<\lvert f(z)\rvert$ for all $\lvert z-z_0\rvert=r$.
	By Rouch\'e's theorem, $f_n$ and $f$ have the same number of zeros inside the disk.
	Hence $f_n$ has exactly $m$ zeros in the disk for all large $n$.
	
	Now suppose each $f_n$ has no zeros. If $f\not\equiv0$ had any
	zero, the preceding conclusion would force all sufficiently large $f_n$
	to have a zero, a contradiction. Thus either $f$ is zero-free or
	$f\equiv0$.
\end{proof}

\begin{remark}
	Hurwitz's theorem formalizes the idea that zeros of holomorphic functions are
	\emph{stable under uniform convergence}.  Zeros may move continuously, but they
	cannot appear spontaneously in the limit unless the limit function collapses to
	the zero function.
\end{remark}

\subsubsection{Examples}

\begin{example}[Convergence of simple zeros]
	Let $f_n(z)=z+\frac{1}{n}$.  Then $f_n$ converges uniformly on compact sets to
	$f(z)=z$.  Each $f_n$ has a simple zero at $z=-\frac{1}{n}$, and these zeros
	converge to the zero of the limit function.
\end{example}

\begin{example}[Degenerate limit]
	Let $f_n(z)=\frac{z}{n}$.  Each $f_n$ has a zero at the origin, and $f_n$ converges
	uniformly on compact sets to the zero function.  This illustrates the exceptional
	case allowed by Hurwitz's theorem.
\end{example}

\subsubsection{The Open Mapping Theorem}

\begin{theorem}[Open Mapping Theorem]
	Let $\Omega\subset\mathbb{C}$ be a domain.  Every nonconstant holomorphic function
	$f:\Omega\to\mathbb{C}$ maps open sets to open sets.
\end{theorem}

\begin{proof}
	Fix a point $z_0\in\Omega$.  Since $f$ is holomorphic and nonconstant, its Taylor
	expansion at $z_0$ has the form
	$f(z)=f(z_0)+a_m(z-z_0)^m+\text{higher order terms}$ with $m\ge1$ and $a_m\neq 0$.
	
	Choose a closed disk $\overline{B(z_0,r)}\Subset\Omega$ so small that
	$f(z)\ne f(z_0)$ on its boundary. Set
	\[
	\eta:=\min_{|z-z_0|=r}|f(z)-f(z_0)|>0.
	\]
	If $|w-f(z_0)|<\eta$, then on the boundary
	\[
	|(f(z)-w)-(f(z)-f(z_0))|=|w-f(z_0)|<|f(z)-f(z_0)|.
	\]
	Rouché's theorem shows that $f-w$ and $f-f(z_0)$ have the same
	number of zeros in the disk. The latter has at least the zero $z_0$,
	so $w\in f(B(z_0,r))$. Thus the full disk
	$B(f(z_0),\eta)$ lies in the image, proving openness.
\end{proof}

\subsubsection{Further examples}

\begin{example}[Failure for constant functions]
	A constant function maps every open set to a single point, which is not open.
	Thus the nonconstancy assumption in the Open Mapping Theorem is essential.
\end{example}

\begin{example}[Complex versus real analysis]
	The smooth nonconstant map $(x,y)\mapsto(x,0)$ sends an open disk to a
	line segment, not an open subset of $\mathbb R^2$. The Open Mapping
	Theorem is therefore a consequence of holomorphicity, not mere smoothness.
\end{example}

\subsubsection{Exercises}

\begin{exercise}
	Let $f_n\in\mathcal{O}(\Omega)$ converge uniformly on compact sets to a nonconstant
	$f$.  Assume that all zeros of $f_n$ lie in a fixed compact subset of $\Omega$.
	Show that every zero of $f$ is the limit of a sequence of zeros of $f_n$.
\end{exercise}

\begin{exercise}
	Let $f_n:\Omega\to\mathbb{C}$ be injective and holomorphic, and assume that
	$f_n$ converges uniformly on compact subsets to a holomorphic function $f$.
	Show that $f$ is either injective or constant.
\end{exercise}

\begin{exercise}
	Use the Open Mapping Theorem to derive the Maximum Modulus Principle for a
	nonconstant holomorphic function on a bounded domain with continuous boundary
	values.
\end{exercise}

\begin{exercise}
	Let $f\in\mathcal{O}(\Omega)$ and suppose $f(z_0)=w_0$ with multiplicity $m\ge1$.
	Show that there exist neighborhoods of $z_0$ and $w_0$ such that every
	$w\neq w_0$ near $w_0$ has exactly $m$ preimages near $z_0$, counted with
	multiplicity.
\end{exercise}

\subsection{Picard Theorems and Essential Singularities}

Near an essential singularity, the image is far less constrained than near a
pole. Casorati--Weierstrass gives density, while the Little and Great Picard
theorems sharpen density to a value-distribution statement. The functions
$e^{1/z}$ and $\sin(1/z)$ provide concrete tests of each conclusion.

\subsubsection{Essential singularities}

\begin{definition}[Isolated singularity and essential singularity]
	Let $f$ be holomorphic on a punctured disk
	$\{z\in\mathbb{C}:0<\lvert z-z_0\rvert<r\}$.
	We say that $z_0$ is an \emph{isolated singularity} of $f$.
	
	Write the Laurent expansion of $f$ about $z_0$ as
	$f(z)=\sum_{n=-\infty}^{\infty} a_n (z-z_0)^n$ on $0<\lvert z-z_0\rvert<r$.
	Then:
	\begin{itemize}
		\item $z_0$ is \emph{removable} iff $a_n=0$ for all $n<0$;
		\item $z_0$ is a \emph{pole of order $m$} iff $a_{-m}\neq 0$ and $a_n=0$ for all $n<-m$;
		\item $z_0$ is \emph{essential} iff $a_n\neq 0$ for infinitely many $n<0$.
	\end{itemize}
\end{definition}

\begin{example}[A classical essential singularity]
	The function $f(z)=e^{1/z}$ has an essential singularity at $z=0$.
	Indeed, the Taylor series $e^w=\sum_{n=0}^\infty \frac{w^n}{n!}$ gives the Laurent
	expansion $e^{1/z}=\sum_{n=0}^\infty \frac{1}{n!}z^{-n}$, which contains infinitely
	many negative powers.
\end{example}

\subsubsection{Casorati--Weierstrass}

\begin{theorem}[Casorati--Weierstrass]
	Let $f$ be holomorphic on $0<\lvert z-z_0\rvert<r$ and suppose $z_0$ is an essential
	singularity.  Then for every $\rho\in(0,r)$ the set
	$f(\{z:0<\lvert z-z_0\rvert<\rho\})$ is dense in $\mathbb{C}$.
	
	Equivalently: for every $w\in\mathbb{C}$ and every $\varepsilon>0$, for every
	$\rho\in(0,r)$ there exists $z$ with $0<\lvert z-z_0\rvert<\rho$ such that
	$\lvert f(z)-w\rvert<\varepsilon$.
\end{theorem}

\begin{proof}
	Fix $\rho\in(0,r)$.  Suppose, toward a contradiction, that
	$f(\{z:0<\lvert z-z_0\rvert<\rho\})$ is \emph{not} dense in $\mathbb{C}$.
	Then there exist $w_0\in\mathbb{C}$ and $\varepsilon_0>0$ such that
	$\lvert f(z)-w_0\rvert\ge \varepsilon_0$ for all $z$ with
	$0<\lvert z-z_0\rvert<\rho$.
	
	Define $g(z)=\frac{1}{f(z)-w_0}$ on $0<\lvert z-z_0\rvert<\rho$.
	Then $g$ is holomorphic on the punctured disk and satisfies
	$\lvert g(z)\rvert\le \frac{1}{\varepsilon_0}$ there, hence $g$ is bounded near $z_0$.
	By the Removable Singularity Theorem, $g$ extends holomorphically to $z_0$.
	Let $\widetilde g$ denote the extension to $\lvert z-z_0\rvert<\rho$.
	
	Now consider two cases.
	
	\smallskip
	\noindent\textbf{Case 1: $\widetilde g(z_0)\neq 0$.}
	Then $\widetilde g$ is nonzero in a small neighborhood of $z_0$, so
	$\frac{1}{\widetilde g}$ is holomorphic near $z_0$.  On $0<\lvert z-z_0\rvert<\rho$
	we have $f(z)-w_0=\frac{1}{g(z)}=\frac{1}{\widetilde g(z)}$, hence
	$f(z)=w_0+\frac{1}{\widetilde g(z)}$ extends holomorphically across $z_0$.
	Thus $z_0$ is removable for $f$, contradiction.
	
	\smallskip
	\noindent\textbf{Case 2: $\widetilde g(z_0)=0$.}
	Then $z_0$ is a zero of $\widetilde g$ of some finite order $m\ge 1$, so
	$\widetilde g(z)=(z-z_0)^m h(z)$ with $h$ holomorphic and $h(z_0)\neq 0$.
	Hence $1/\widetilde g$ has a pole of order $m$ at $z_0$, and
	$f(z)=w_0+\frac{1}{\widetilde g(z)}$ extends meromorphically across $z_0$ with a pole.
	Thus $z_0$ is a pole for $f$, contradiction.
	
	\smallskip
	In all cases we contradict the assumption that $z_0$ is essential.  Therefore
	$f(\{z:0<\lvert z-z_0\rvert<\rho\})$ must be dense in $\mathbb{C}$.
\end{proof}

\subsubsection{Concrete value computations near an essential singularity}

\begin{example}[Explicit solutions for $e^{1/z}=w$ and accumulation at $0$]
	Fix $w\in\mathbb{C}$ with $w\neq 0$.  Solving $e^{1/z}=w$ is equivalent to
	$1/z\in \Log(w)$, where $\Log(w)=\{\log\lvert w\rvert+i(\arg w+2\pi k):k\in\mathbb{Z}\}$.
	Thus the solutions are
	$z_k=\frac{1}{\log\lvert w\rvert+i(\arg w+2\pi k)}$ for those
	$k\in\mathbb{Z}$ for which the denominator is nonzero. The excluded
	denominator occurs only when $w=1$, for one choice of $k$.
	As $\lvert k\rvert\to\infty$, the denominator has modulus asymptotic to $2\pi\lvert k\rvert$,
	so $\lvert z_k\rvert\to 0$.  Hence $e^{1/z}$ takes every nonzero value infinitely
	often in any punctured neighborhood of $0$.
	
	Note that $e^{1/z}$ never takes the value $0$, so this example also illustrates
	the possibility of omitting exactly one value.
\end{example}

\begin{example}[Zeros of $\sin(1/z)$ accumulate at the essential singularity]
	Let $f(z)=\sin(1/z)$.  Then $z=0$ is essential because
	$\sin w=\sum_{n=0}^{\infty}(-1)^n\frac{w^{2n+1}}{(2n+1)!}$ gives
	$\sin(1/z)=\sum_{n=0}^\infty (-1)^n\frac{1}{(2n+1)!}z^{-(2n+1)}$, with infinitely many
	negative powers.
	
	Moreover, $f(z)=0$ iff $1/z=\pi k$ for some $k\in\mathbb{Z}\setminus\{0\}$, i.e.
	$z_k=\frac{1}{\pi k}$.  Then $z_k\to 0$ as $\lvert k\rvert\to\infty$.
	So $0$ is taken infinitely often near the essential singularity.
\end{example}

\subsubsection{Picard theorems}

\begin{theorem}[Great Picard]
	\label{thm:great-picard}
	Let $f$ be holomorphic on $0<\lvert z-z_0\rvert<r$ with an essential singularity at $z_0$.
	Then in every punctured neighborhood of $z_0$, the function $f$ assumes every
	complex value infinitely often, with at most one possible exception.
\end{theorem}

The proof of Great Picard is not reproduced here; it requires a substantially
stronger normal-family or covering-space argument than Casorati--Weierstrass.
For proofs of this precise infinite-attainment statement, see
\cite{Ahlfors1979,Conway1978}.

\begin{remark}
	Great Picard strengthens Casorati--Weierstrass in two ways:
	(i) it upgrades density to actual attainment of values, and
	(ii) it asserts infinite repetition of values arbitrarily close to the singularity.
\end{remark}

\begin{theorem}[Little Picard]
	If $f:\mathbb{C}\to\mathbb{C}$ is a nonconstant entire function, then $f(\mathbb{C})$
	omits at most one complex value.
\end{theorem}

\begin{proof}
	If $f$ is a nonconstant polynomial, the fundamental theorem of algebra shows
	that it assumes every complex value. Otherwise $g(w)=f(1/w)$ has an essential
	singularity at $w=0$: a removable singularity would make $f$ bounded near
	infinity and hence constant by Liouville. A pole of order $m$ would give
	$|f(z)|=O(|z|^m)$ at infinity, so Cauchy's estimates force every
	derivative of order greater than $m$ to vanish; then $f$ is a polynomial.
	Apply Theorem~\ref{thm:great-picard} to $g$ on a punctured disk. At most one
	complex value can be absent there, and therefore at most one can be absent
	from $f(\mathbb C)$.
\end{proof}

\begin{remark}
	The example $f(z)=e^z$ shows that one omitted value is possible. The proof
	above establishes Little Picard as a consequence of Great Picard; it does not
	provide a new proof of the latter theorem.
\end{remark}

\subsubsection{A useful corollary: classification of isolated singularities}

\begin{proposition}[Picard-flavored characterization]
	Let $f$ be holomorphic on $0<\lvert z-z_0\rvert<r$.
	\begin{itemize}
		\item If $f$ is bounded near $z_0$, then $z_0$ is removable.
		\item If $\lvert f(z)\rvert\to\infty$ as $z\to z_0$, then $z_0$ is a pole.
		\item If $f$ omits two distinct complex values in some punctured neighborhood of $z_0$,
		then $z_0$ cannot be essential; it is removable or a pole.
	\end{itemize}
\end{proposition}

\begin{proof}
	The first two items are the classical classification of isolated singularities
	using the Riemann Removable Singularity Theorem and the definition of poles.
	
	For the third item, assume $f$ omits two distinct values $a\neq b$ near $z_0$.
	If $z_0$ were essential, Great Picard would imply that in every
	punctured neighborhood $f$ assumes every complex value, with at most
	one exception, infinitely often. It therefore cannot omit both $a$
	and $b$. Hence the singularity is removable or a pole.
\end{proof}

\subsubsection{Exercises}

\begin{exercise}[Essential singularity of $\sin(1/z)$ and explicit zeros]
	Let $f(z)=\sin(1/z)$.
	\begin{enumerate}
		\item Prove that $z=0$ is an essential singularity of $f$ by writing its Laurent expansion.
		\item Show that $f(z)=0$ has infinitely many solutions $z_k$ with $z_k\to 0$ and find them explicitly.
		\item Show that $f$ takes the values $\pm 1$ infinitely often near $0$ by solving $\sin(1/z)=\pm 1$.
	\end{enumerate}
\end{exercise}

\begin{exercise}[Casorati--Weierstrass versus removable/pole]
	Let $f$ be holomorphic on $0<\lvert z\rvert<1$ and assume $f(0<\lvert z\rvert<1)$
	is contained in a horizontal line $\{w:\Im w=c\}$ for some $c\in\mathbb{R}$.
	\begin{enumerate}
		\item Use the Open Mapping Theorem to show that $f$ is constant, and hence bounded near $0$.
		\item Conclude that $0$ cannot be an essential singularity of $f$.
		\item Explain how this is consistent with Casorati--Weierstrass.
	\end{enumerate}
\end{exercise}

\begin{exercise}[Essential singularity at infinity]
	Let $f$ be entire and nonconstant, and define $g(z)=f(1/z)$ on $0<\lvert z\rvert<1$.
	\begin{enumerate}
		\item Show that if $f$ is a polynomial of degree $d\ge 1$, then $z=0$ is a pole of order $d$ for $g$.
		\item Show that if $f$ is not a polynomial, then $z=0$ is an essential singularity for $g$.
		\item Interpret this as: every non-polynomial entire function has an essential singularity at $\infty$.
	\end{enumerate}
\end{exercise}

\begin{exercise}[Infinite valence of transcendental entire functions]
	Let $f$ be a transcendental entire function. Use Great Picard at infinity
	to show that, apart from at most one complex value, every value has
	infinitely many preimages with arbitrarily large modulus.
\end{exercise}

\begin{exercise}[A value-omission test near a suspected essential singularity]
	Let $f$ be holomorphic on $0<\lvert z\rvert<1$.
	Assume there exist two values $a\neq b$ such that $f(z)\neq a$ and $f(z)\neq b$
	for all $0<\lvert z\rvert<\rho$ for some $\rho\in(0,1)$.
	Show that $z=0$ is not an essential singularity of $f$.
\end{exercise}

\begin{exercise}[A sharpened computation for $e^{1/z}$]
	Fix $w\neq 0$. For each integer $k$ with nonzero denominator, let
	$z_k=\frac{1}{\log\lvert w\rvert+i(\arg w+2\pi k)}$.
	\begin{enumerate}
		\item Show that $e^{1/z_k}=w$ for every defined $z_k$.
		\item Show that $\lvert z_k\rvert\le \frac{1}{\pi\lvert k\rvert}$ for all sufficiently large $\lvert k\rvert$.
		\item Conclude that $z_k\to0$ with $|z_k|=O(1/|k|)$ (indeed,
		$|z_k|\sim(2\pi|k|)^{-1}$).
	\end{enumerate}
\end{exercise}

\begin{exercise}[Doubly periodic entire functions]
	Let $f$ be entire and doubly periodic with respect to a lattice $\Lambda\subset\mathbb{C}$.
	\begin{enumerate}
		\item Show that $f$ is bounded on a fundamental parallelogram for $\Lambda$.
		\item Conclude that $f$ is bounded on $\mathbb{C}$ and hence constant by Liouville's theorem.
		\item Explain why a nonconstant doubly periodic function cannot be entire.
		If it is meromorphic on $\mathbb C$, explain why it must have a pole
		in a fundamental parallelogram.
	\end{enumerate}
\end{exercise}

\subsection{Infinite Products and Weierstrass Factorization}
\index{Weierstrass factorization}\index{infinite products}

We develop infinite products of numbers and holomorphic functions and then prove
the Weierstrass factorization theorem. Infinite products construct entire
functions with prescribed zero sets, extending polynomial factorization. The
analytic issue is compact-uniform convergence of the partial products, which
ensures that the limit is holomorphic.

\subsubsection{Infinite products of complex numbers}

\begin{definition}[Infinite product]
	Let $\{a_n\}_{n\ge 1}$ be a sequence of complex numbers with $a_n\neq 0$.
	The infinite product $\prod_{n=1}^{\infty} a_n$ is said to \emph{converge} if the
	partial products $P_N:=\prod_{n=1}^{N} a_n$ converge to a \emph{nonzero} complex number.
	If $P_N\to 0$, the product is said to \emph{diverge to zero}.  Otherwise it is
	\emph{divergent}.
\end{definition}

\begin{remark}[A necessary condition]
	If $\prod_{n=1}^{\infty} a_n$ converges to a nonzero limit, then necessarily
	$a_n\to 1$.  Indeed, $a_N=P_N/P_{N-1}\to L/L=1$ if $P_N\to L\neq 0$.
\end{remark}

\begin{proposition}[Logarithmic criterion near $1$]
	Assume $a_n\to 1$.  Then there exists $N_0$ such that for all $n\ge N_0$,
	the principal branch $\Log(a_n)$ is well-defined (after restricting to a small
	disk around $1$ that avoids the branch cut).  In this situation,
	$\prod_{n=1}^{\infty} a_n$ converges to a nonzero limit if and only if
	$\sum_{n=N_0}^{\infty} \Log(a_n)$ converges in $\mathbb{C}$.
\end{proposition}

\begin{proof}
	Since $a_n\to 1$, choose $\delta\in(0,1)$ and $N_0$ so that $\lvert a_n-1\rvert<\delta$
	for all $n\ge N_0$ and the disk $\{w:\lvert w-1\rvert<\delta\}$ does not intersect the
	branch cut of $\Log$.  Then $\Log(a_n)$ is defined for $n\ge N_0$ and is continuous.
	
	If $\sum_{n=N_0}^{\infty}\Log(a_n)$ converges, then
	$S_N:=\sum_{n=N_0}^{N}\Log(a_n)$ converges to some $S$, hence
	$\prod_{n=N_0}^{N} a_n=\exp(S_N)\to \exp(S)\neq 0$.
	Multiplying by the finite constant $\prod_{n=1}^{N_0-1} a_n\neq 0$ yields convergence
	of $\prod_{n=1}^{\infty} a_n$.
	
	Conversely, let $P_N=\prod_{n=1}^N a_n\to L\ne0$. Choose a disk
	$V$ about $L$ avoiding $0$ and a holomorphic logarithm $\log_V$ on
	$V$. For all large $N$, both $P_N$ and $P_{N-1}$ lie in $V$. Then
	\[
	d_N:=\log_V(P_N)-\log_V(P_{N-1})
	\]
	is a logarithm of $a_N$. Both $d_N$ and the principal $\Log(a_N)$
	tend to $0$, while their difference lies in $2\pi i\Z$; hence they
	are equal for all sufficiently large $N$. The tail sums of
	$\Log(a_N)$ therefore telescope to $\log_V(P_N)$ up to a constant
	and converge.
\end{proof}

\begin{remark}[A practical sufficient condition]
	If $a_n=1+b_n$ and $\sum_{n=1}^{\infty}\lvert b_n\rvert<\infty$ with $\lvert b_n\rvert$
	eventually small, then $\sum \Log(1+b_n)$ converges and the product converges.
	This is often used with the estimate $\lvert\Log(1+b)\rvert\le C\lvert b\rvert$ for
	$\lvert b\rvert$ small.
\end{remark}

\subsubsection{Infinite products of holomorphic functions}

\begin{definition}[Compact-uniform convergence of products]
	Let $\Omega\subset\mathbb{C}$ be a domain and let $u_n$ be holomorphic on $\Omega$.
	We say the product $\prod_{n=1}^{\infty} u_n(z)$ converges \emph{uniformly on compact
		sets} to a nonzero holomorphic function if its partial products
	$P_N(z)=\prod_{n=1}^{N}u_n(z)$ converge uniformly on every compact subset to a
	holomorphic function $P$ that is not identically zero.  In the standard product
	criterion below, the tail is nonvanishing on each fixed compact set; zeros of $P$
	may still come from the finitely many earlier factors.
\end{definition}

\begin{proposition}[Weierstrass M-test for products]
	Let $\Omega$ be a domain and let $u_n$ be holomorphic functions on $\Omega$, none
	of which is identically zero.
	Assume that for every compact $K\Subset\Omega$ there exist numbers $M_n(K)\ge 0$ such that
	$\sup_{z\in K}\lvert u_n(z)-1\rvert\le M_n(K)$ and $\sum_{n=1}^{\infty} M_n(K)<\infty$.
	Then $\prod_{n=1}^{\infty} u_n(z)$ converges uniformly on compact subsets
	of $\Omega$ to a holomorphic function that is not identically zero.
	Moreover, near each point a sufficiently late tail product is nonvanishing.
	Consequently, the zeros of the limit are exactly the zeros of the factors
	$u_n$, with their multiplicities added when several factors vanish at the
	same point.
\end{proposition}

\begin{proof}
	Fix compact $K\Subset\Omega$, and choose a compact neighborhood
	$L\Subset\Omega$ of $K$. Since $\sum M_n(L)<\infty$, we have
	$M_n(L)\to 0$. Hence for $n$ large enough,
	$\sup_{L}\lvert u_n-1\rvert\le 1/2$, so $\lvert u_n\rvert\ge 1/2$
	on $L$ for large $n$. Thus the tail factors have no zeros near $K$.
	
	For these $n$, define $v_n(z)=\Log(u_n(z))$ on the interior of $L$ using the
	principal branch on $\{w:\lvert w-1\rvert<1/2\}$. We have the estimate
	$\lvert v_n(z)\rvert\le C\lvert u_n(z)-1\rvert$ uniformly on $L$
	for a universal constant $C$ (e.g.\ from the power series of $\Log(1+\xi)$ for $\lvert\xi\rvert<1/2$).
	Hence $\sum \sup_L\lvert v_n\rvert<\infty$.
	
	Therefore the series $\sum v_n(z)$ converges uniformly on $K$ to a holomorphic function
	$V(z)$ (uniform limit of holomorphic partial sums).  Exponentiating gives uniform convergence
	of the product of the tail: $\prod_{n\ge N}u_n(z)=\exp(\sum_{n\ge N}v_n(z))\to \exp(V(z))$.
	Multiplying by the finite initial product yields convergence of $\prod_{n=1}^{\infty}u_n(z)$.
	
	Finally, the tail limit is the exponential of a holomorphic function
	and is therefore nonvanishing. Near any prescribed point, all sufficiently
	late factors belong to such a nonvanishing tail. Only finitely many earlier
	factors can therefore contribute zeros there, and their orders add under
	multiplication. Since no finite product of the nonzero holomorphic factors is
	identically zero, neither is the limiting product.
\end{proof}

\subsubsection{Canonical factors}

\begin{definition}[Weierstrass canonical factors]
	For $p\in\mathbb{N}\cup\{0\}$ define the canonical factor of genus $p$ by
	$E_0(w)=1-w$ and, for $p\ge 1$,
	$E_p(w)=(1-w)\exp\!\bigl(w+\frac{w^2}{2}+\cdots+\frac{w^p}{p}\bigr)$.
\end{definition}

\begin{lemma}[Key cancellation estimate]
	Fix $p\in\mathbb{N}\cup\{0\}$.  There exist constants $C_p>0$ and $\delta_p\in(0,1)$ such that
	for all $w$ with $\lvert w\rvert\le \delta_p$,
	$\lvert \Log(E_p(w))\rvert \le C_p \lvert w\rvert^{p+1}$.
	In particular, $E_p(w)=1+O(\lvert w\rvert^{p+1})$ as $w\to 0$.
\end{lemma}

\begin{proof}
	For $\lvert w\rvert$ small, $\Log(1-w)=-\sum_{k=1}^{\infty}\frac{w^k}{k}$.
	By definition,
	$\Log(E_p(w))=\Log(1-w)+w+\frac{w^2}{2}+\cdots+\frac{w^p}{p}$.
	Substituting the series for $\Log(1-w)$, the terms $-w,-w^2/2,\dots,-w^p/p$ cancel,
	leaving
	$\Log(E_p(w))=-\sum_{k=p+1}^{\infty}\frac{w^k}{k}$.
	Hence for $\lvert w\rvert\le 1/2$,
	$\lvert \Log(E_p(w))\rvert \le \sum_{k=p+1}^{\infty}\frac{\lvert w\rvert^k}{k}
	\le \sum_{k=p+1}^{\infty}\lvert w\rvert^k
	= \frac{\lvert w\rvert^{p+1}}{1-\lvert w\rvert}
	\le 2\lvert w\rvert^{p+1}$.
	This gives the desired bound with $C_p=2$ and $\delta_p=1/2$.
	Exponentiating yields $E_p(w)=1+O(\lvert w\rvert^{p+1})$.
\end{proof}

\subsubsection{Convergence of canonical products}

\begin{proposition}[Canonical products converge on compacts]
	Let $\{a_n\}$ be a sequence of nonzero complex numbers with $\lvert a_n\rvert\to\infty$.
	Fix an integer $p\ge 0$ and assume $\sum_{n=1}^{\infty} \lvert a_n\rvert^{-(p+1)}<\infty$.
	Define
	$P_N(z)=\prod_{n=1}^{N} E_p\!\bigl(\frac{z}{a_n}\bigr)$.
	Then $P_N$ converges uniformly on compact subsets of $\mathbb{C}$ to an entire function
	$P(z)$.  Moreover, $P$ vanishes precisely at the points $z=a_n$ (with multiplicities
	matching repetitions in the sequence) and has no other zeros.
\end{proposition}

\begin{proof}
	Fix a compact set $K\Subset\mathbb{C}$ and set $R=\sup_{z\in K}\lvert z\rvert$.
	Since $\lvert a_n\rvert\to\infty$, choose $N_0$ so that $\lvert z/a_n\rvert\le 1/2$
	for all $z\in K$ and all $n\ge N_0$.  Then the lemma applies to $w=z/a_n$ and gives
	$\lvert \Log(E_p(z/a_n))\rvert \le C_p \lvert z/a_n\rvert^{p+1}\le C_p R^{p+1}\lvert a_n\rvert^{-(p+1)}$
	uniformly on $K$.
	
	Therefore the series $\sum_{n\ge N_0}\Log(E_p(z/a_n))$ converges uniformly on $K$ by the
	Weierstrass M-test, because $\sum \lvert a_n\rvert^{-(p+1)}<\infty$.
	Its sum is holomorphic on a neighborhood of $K$.  Exponentiating shows the tail product
	$\prod_{n\ge N_0}E_p(z/a_n)$ converges uniformly on $K$ to a holomorphic nonvanishing limit.
	Multiplying by the finite initial product $\prod_{n=1}^{N_0-1}E_p(z/a_n)$ gives uniform
	convergence of $P_N$ on $K$ to a holomorphic limit.  Since $K$ was arbitrary, the limit is entire.
	
	For zeros, fix $z_0\notin\{a_n\}$.  A small closed disk about $z_0$ avoids every $a_n$;
	the finite initial product is nonzero there after the disk is made smaller, and the uniformly
	convergent logarithmic tail is the exponential of a holomorphic function.  Hence $P(z_0)\ne0$.
	At a listed point $a$, isolate the finitely many factors for which $a_n=a$.  Each has a simple
	zero there, while all remaining factors give a nonzero holomorphic function near $a$ by the same
	tail argument.  Thus the order of the zero is exactly the number of repetitions of $a$.
\end{proof}

\subsubsection{Weierstrass Factorization Theorem}

\begin{theorem}[Weierstrass Factorization Theorem]
	Let $f$ be an entire function, not identically zero.  Let $m\ge 0$ be the multiplicity
	of the zero of $f$ at $0$, and let $\{a_n\}$ be the (finite or infinite) multiset of nonzero
	zeros of $f$, listed with multiplicities and with no accumulation point in $\mathbb{C}$.
	Then there exist integers $p_n\ge 0$ and an entire function $g$ such that
	$f(z)=z^{m}\exp(g(z))\prod_{n=1}^{\infty} E_{p_n}\!\bigl(\frac{z}{a_n}\bigr)$,
	where the product is interpreted as $1$ if there are finitely many zeros.
	Conversely, any expression of this form defines an entire function with the prescribed zeros.
\end{theorem}

\begin{proof}
	Step 1: Choose genera $p_n$ so that the canonical product converges.
	Because $|a_n|\to\infty$, the simple choice $p_n=n$ suffices.  Indeed, on a fixed disk
	$|z|\le R$, one has $|z/a_n|\le1/2$ for all sufficiently large $n$, and then
	$|z/a_n|^{n+1}\le2^{-(n+1)}$.
	
	Define the canonical product
	$P(z)=z^{m}\prod_{n=1}^{\infty}E_{p_n}\!\bigl(\frac{z}{a_n}\bigr)$.
	The cancellation estimate gives
	$|\Log E_n(z/a_n)|\le2|z/a_n|^{n+1}\le2^{-n}$ for all sufficiently large $n$,
	uniformly on that disk.  The M-test therefore gives compact-uniform convergence of the
	logarithmic tails.  The zero argument of the proposition
	shows that $P$ has precisely the prescribed zeros.
	
	Step 2: The quotient $f/P$ is entire and has no zeros.
	Define $h(z)=f(z)/P(z)$.  This is holomorphic on $\mathbb{C}\setminus\{0,a_1,a_2,\dots\}$.
	Near any zero $a_n$, both $f$ and $P$ vanish to the same order, so the singularity of $h$
	at $a_n$ is removable.  Similarly at $0$.  Hence $h$ extends holomorphically to all of $\mathbb{C}$,
	so $h$ is entire.  By construction, $h$ has no zeros.
	
	Step 3: A nonvanishing entire function is an exponential.
	Since $h$ is entire and never vanishes, the logarithmic derivative $h'(z)/h(z)$ is entire.
	Choose a base point $z_0\in\mathbb{C}$ and define
	$g(z)=\int_{z_0}^{z}\frac{h'(\zeta)}{h(\zeta)}\,d\zeta$, where the integral is taken along any
	piecewise $C^1$ path in $\mathbb{C}$.  This is well-defined because $\mathbb{C}$ is simply connected
	and $h'/h$ is entire, so the integral is path-independent.  Then $g$ is entire and
	$\exp(g(z))$ satisfies $(\exp g)'/(\exp g)=g'=h'/h$, hence $h/\exp(g)$ is constant.
	Adjusting $g$ by an additive constant, we obtain $h(z)=\exp(g(z))$.
	
	Combining steps yields $f(z)=P(z)\exp(g(z))$, which is the desired factorization.
	The converse direction is immediate from the convergence proposition and the fact that
	$\exp(g)$ has no zeros.
\end{proof}

\subsubsection{Classical examples}

\begin{example}[Sine product]
	The zeros of $\sin(\pi z)$ are the integers $\mathbb{Z}$, all simple.
	The classical sine product is
	\[
	\sin(\pi z)=\pi z \prod_{n=1}^{\infty}\left(1-\frac{z^2}{n^2}\right).
	\]
	This identity, including its normalization, is recorded in
	\cite[Eq.~4.22.1]{DLMFTrig2026}.
	The sum $\sum_{n\ge1}n^{-2}$ ensures compact-uniform convergence of the
	paired genus-zero factors. Convergence and the zero set alone do not fix
	the nonvanishing exponential factor or the constant $\pi$; these are
	part of the cited classical identity.
\end{example}

\begin{example}[Euler product for $1/\Gamma$]
	Euler's product can be written as
	$\frac{1}{\Gamma(z)}=z e^{\gamma z}\prod_{n=1}^{\infty}\bigl(1+\frac{z}{n}\bigr)e^{-z/n}$,
	where $\gamma$ is the Euler--Mascheroni constant.  This is a canonical product using
	genus $1$ factors since $\bigl(1+\frac{z}{n}\bigr)e^{-z/n}=E_1(-z/n)$;
	see \cite[Eq.~5.8.2]{DLMFGamma2026} for the normalization.
\end{example}

\subsubsection{Exercises}

\begin{exercise}[Asymptotics of $E_1$]
	Compute $E_1(z)=(1-z)e^z$ and show that $E_1(z)=1-\frac{z^2}{2}+O(\lvert z\rvert^3)$ as $z\to 0$.
	Deduce that if $\prod_{n=1}^{\infty}E_1(z/a_n)$ is to converge compact-uniformly, it is natural to
	require $\sum_{n=1}^{\infty}\lvert a_n\rvert^{-2}<\infty$.
\end{exercise}

\begin{exercise}[Entire function with finitely many zeros]
	Let $f$ be entire and assume it has only finitely many zeros (counted with multiplicity).
	Show that $f(z)=e^{g(z)}P(z)$ where $P$ is a polynomial and $g$ is entire.
\end{exercise}

\begin{exercise}[A simple canonical product]
	Show that $f(z)=\prod_{n=1}^{\infty}\bigl(1+\frac{z}{n^2}\bigr)$ defines an entire function.
	Hint: use the product M-test on compact sets with $u_n(z)=1+\frac{z}{n^2}$ and
	$\sup_{z\in K}\lvert u_n(z)-1\rvert\le R/n^2$.
\end{exercise}

\begin{exercise}[Uniqueness up to exponentials]
	Let $f$ and $g$ be nonzero entire functions with the same zeros (including multiplicity).
	Use Weierstrass factorization to show there exists an entire $h$ with $f(z)=e^{h(z)}g(z)$.
\end{exercise}

\begin{exercise}[Fixed genus for finite order: Hadamard style]
	Assume $f$ is an entire function of finite order $\rho<\infty$ and its zeros $\{a_n\}$ satisfy
	$\#\{n:\lvert a_n\rvert\le R\}=O(R^\rho)$ as $R\to\infty$.
	Explain why one may choose a single integer $p\ge \rho$ and write a factorization of $f$
	using canonical factors $E_p(z/a_n)$ of fixed genus $p$.
\end{exercise}

\subsection{Entire Functions: Order and Type}

Order and type quantify how fast an entire function grows at infinity. They are
central in Hadamard factorization and in the analytic study of special functions
(Gamma, Bessel, zeta, \dots).

\subsubsection{Maximum modulus and basic properties}

\begin{definition}[Maximum modulus]
	Let $f$ be entire.  For $r\ge 0$ define
	$M_f(r):=\max_{\lvert z\rvert=r}\lvert f(z)\rvert$.
	When $f$ is clear we simply write $M(r)$.
\end{definition}

\begin{proposition}[Monotonicity of $M(r)$]
	For an entire function $f$, the function $r\mapsto M(r)$ is nondecreasing.
	Moreover, $M(r)=\max_{\lvert z\rvert\le r}\lvert f(z)\rvert$.
\end{proposition}

\begin{proof}
	Fix $0\le r<R$.  By the maximum modulus principle applied to $f$ on the disk
	$\{z:\lvert z\rvert\le R\}$, the maximum of $\lvert f\rvert$ on that closed disk
	is achieved on the boundary circle $\lvert z\rvert=R$.  Hence
	$\max_{\lvert z\rvert\le r}\lvert f(z)\rvert\le \max_{\lvert z\rvert\le R}\lvert f(z)\rvert
	=\max_{\lvert z\rvert=R}\lvert f(z)\rvert=M(R)$.
	Taking maxima over $\lvert z\rvert=r$ shows $M(r)\le M(R)$ and also
	$M(r)=\max_{\lvert z\rvert\le r}\lvert f(z)\rvert$.
\end{proof}

\subsubsection{Order and equivalent growth tests}

\begin{definition}[Order]
	For a nonconstant entire function $f$, its \emph{order} $\rho$ is
	\[
	\rho:=\limsup_{r\to\infty}\frac{\log\log M(r)}{\log r},
	\]
	where $M(r)>1$ for all sufficiently large $r$.  Constant entire functions are assigned
	order $0$ by convention.
	If $\rho<\infty$, we say $f$ has \emph{finite order}.
\end{definition}

\begin{remark}[What $\log\log M(r)$ measures]
	Roughly, $\log M(r)$ measures exponential growth, and $\log\log M(r)$ measures
	\emph{power} growth of the exponent.  For instance:
	polynomials satisfy $\log M(r)\sim n\log r$ so $\log\log M(r)\sim \log\log r$,
	hence order $0$, while $e^{az}$ satisfies $\log M(r)\sim \lvert a\rvert r$ so
	$\log\log M(r)\sim \log r$, hence order $1$.
\end{remark}

\begin{proposition}[Equivalent characterization of order]
	Let $f$ be a nonconstant entire function and $\rho\in[0,\infty]$.  Then $\rho$ is the order of $f$ if and only if:
	\begin{enumerate}
		\item for every $\alpha>\rho$ there exists $R_\alpha$ such that
		$M(r)\le \exp(r^\alpha)$ for all $r\ge R_\alpha$;
		\item for every $\beta<\rho$ and every $C>0$ there exist arbitrarily large $r$ such that
		$M(r)\ge \exp(C r^\beta)$.
	\end{enumerate}
\end{proposition}

\begin{proof}
	Assume $\rho<\infty$ first.
	If $\alpha>\rho$, then by definition of limsup there exists $R_\alpha$ such that
	$\frac{\log\log M(r)}{\log r}\le \alpha$ for all $r\ge R_\alpha$.
	Equivalently $\log\log M(r)\le \alpha\log r$, so $\log M(r)\le r^\alpha$ and
	$M(r)\le \exp(r^\alpha)$ for $r\ge R_\alpha$.
	
	Conversely, if $\beta<\rho$, then by definition of limsup there exist arbitrarily large $r$
	for which $\frac{\log\log M(r)}{\log r}\ge \beta+\varepsilon$ for some $\varepsilon>0$ with
	$\beta+\varepsilon<\rho$.  Then $\log\log M(r)\ge (\beta+\varepsilon)\log r$, hence
	$\log M(r)\ge r^{\beta+\varepsilon}$.  Given $C>0$, for sufficiently large such $r$ we have
	$r^{\beta+\varepsilon}\ge C r^\beta$, so $M(r)\ge \exp(C r^\beta)$ along an unbounded sequence.
	The case $\rho=\infty$ is the same limsup argument.  When $\rho=0$, condition (i) is the
	corresponding upper-growth assertion; condition (ii), which only concerns $\beta<0$, follows
	because a nonconstant entire function has $M(r)\to\infty$, whereas
	$\exp(Cr^\beta)\to1$ for $\beta<0$.

	Conversely, suppose (i) and (ii) hold and let $\tau$ denote the order computed from
	the defining limsup.  Condition (i) gives $\tau\le\alpha$ for every $\alpha>\rho$,
	so $\tau\le\rho$.  Taking $C=1$ in (ii), applying $\log\log$, and passing to the
	limsup gives $\tau\ge\beta$ for every $\beta<\rho$; hence $\tau\ge\rho$.
	The same argument in the extended real line covers $\rho=\infty$.  Thus $\tau=\rho$.
\end{proof}

\begin{proposition}[A simple comparison criterion]
	Let $f$ be entire.  If there exists $\alpha>0$ and $R_0$ such that
	$\lvert f(z)\rvert\le \exp(\lvert z\rvert^\alpha)$ for all $\lvert z\rvert\ge R_0$,
	then the order of $f$ satisfies $\rho\le \alpha$.
\end{proposition}

\begin{proof}
	If $f$ is constant, the conclusion holds by convention.  Otherwise $M(r)>1$
	for all sufficiently large $r$, and the defining logarithms are available.
	For $r\ge R_0$, the hypothesis gives $M(r)\le \exp(r^\alpha)$.
	Hence $\log\log M(r)\le \log(r^\alpha)=\alpha\log r$, and dividing by $\log r$ and taking
	$\limsup$ yields $\rho\le \alpha$.
\end{proof}

\subsubsection{Type (for finite positive order)}

\begin{definition}[Type]
	Assume $f$ has finite order $\rho\in(0,\infty)$.
	The \emph{type} $\sigma$ of $f$ (with respect to $\rho$) is
	$\sigma:=\limsup_{r\to\infty}\frac{\log M(r)}{r^\rho}$.
	If $0<\sigma<\infty$, $f$ is of \emph{normal type}.  If $\sigma=0$ it is of \emph{minimal type},
	and if $\sigma=\infty$ it is of \emph{infinite type}.
\end{definition}

\begin{remark}[Type depends on the chosen order]
	Type is only defined after fixing the order $\rho$.  Order $1$ implies the upper
	bound $\log M(r)\le r^{1+\varepsilon}$ eventually for every $\varepsilon>0$,
	but it does not make $\log M(r)$ comparable to $r$: minimal and infinite type
	are both possible.  The type records the limsup scale at exactly $r^\rho$.
\end{remark}

\subsubsection{Detailed computations in basic examples}

\begin{example}[Exponentials $e^{az}$]
	Let $f(z)=e^{az}$ with $a\neq 0$.  For $\lvert z\rvert=r$ we have
	$\lvert e^{az}\rvert=e^{\Re(az)}\le e^{\lvert a\rvert \lvert z\rvert}=e^{\lvert a\rvert r}$,
	and equality occurs at $z=r\cdot \overline{a}/\lvert a\rvert$ (so that $az=\lvert a\rvert r$ is real).
	Thus $M(r)=e^{\lvert a\rvert r}$.
	Then $\log M(r)=\lvert a\rvert r$ and $\log\log M(r)=\log(\lvert a\rvert r)=\log r+O(1)$, so
	$\rho=\limsup_{r\to\infty}\frac{\log r+O(1)}{\log r}=1$.
	Moreover
	$\sigma=\limsup_{r\to\infty}\frac{\log M(r)}{r}= \lim_{r\to\infty}\frac{\lvert a\rvert r}{r}=\lvert a\rvert$.
\end{example}

\begin{example}[Polynomials and the ``order $0$'' phenomenon]
	Let $f(z)=a_n z^n+\cdots+a_0$ with $a_n\neq 0$ and $n\ge1$.
	For $\lvert z\rvert=r$,
	$\lvert f(z)\rvert\le \lvert a_n\rvert r^n+\cdots+\lvert a_0\rvert\le C r^n$ for $r\ge 1$,
	hence $M(r)\le C r^n$.
	For every $z$ with $|z|=r$, the reverse triangle inequality gives
	\[
	 |f(z)|\ge |a_n|r^n-\sum_{k<n}|a_k|r^k.
	\]
	The right-hand side is at least \(\frac12|a_n|r^n\) for all sufficiently
	large $r$, so the same lower bound holds for $M(r)$.
	Thus $M(r)\asymp r^n$ and $\log M(r)= n\log r+O(1)$.
	Hence $\log\log M(r)=\log(n\log r+O(1))=\log\log r+O(1)$ and therefore $\rho=0$.
	Together with the convention for constants, every polynomial therefore has order $0$.
	The converse is false: there are transcendental
	entire functions of order $0$.  What characterizes polynomials is the stronger bound
	$M(r)=O(r^A)$ for some fixed $A$, as shown in the exercises below.
\end{example}

\begin{example}[The function $e^{z^m}$]
	Let $f(z)=e^{z^m}$ with integer $m\ge 1$.
	If $\lvert z\rvert=r$, then $\lvert e^{z^m}\rvert=e^{\Re(z^m)}\le e^{\lvert z^m\rvert}=e^{r^m}$,
	and equality occurs when $z^m=r^m$ is real positive, i.e.\ for $z=r e^{2\pi i k/m}$.
	Thus $M(r)=e^{r^m}$, so $\log M(r)=r^m$ and $\log\log M(r)=\log(r^m)=m\log r$.
	Hence $\rho=\limsup \frac{m\log r}{\log r}=m$.
	The type (with respect to $\rho=m$) is
	$\sigma=\limsup_{r\to\infty}\frac{\log M(r)}{r^m}=\lim_{r\to\infty}\frac{r^m}{r^m}=1$.
\end{example}

\begin{example}[The reciprocal Gamma function: order $1$ and infinite type]
	The Gamma function itself is meromorphic, not entire.  Its reciprocal $1/\Gamma$ is entire.
	For a positive integer $n$, the reflection formula gives
	$|\Gamma(\tfrac12-n)|^{-1}=\Gamma(n+\tfrac12)/\pi$.
	Stirling's formula therefore yields
	$\log M_{1/\Gamma}(n-\tfrac12)\ge n\log n+O(n)$.
	Monotonicity of $M$ extends this lower bound, up to $O(r)$, to every
	large radius $r$.
	For the upper bound, use the product in the preceding example. On
	$|z|=r$, split its logarithm at $n=2r$. The terms with $n\le2r$
	contribute at most
	$\sum_{n\le2r}[\log(1+r/n)+r/n]=O(r\log r)$.
	For $n>2r$, the estimate
	$|\Log(1+z/n)-z/n|\le C r^2/n^2$ gives a tail of order $O(r)$.
	Together with $\log|z e^{\gamma z}|=O(r+\log r)$, this proves
	$\log M_{1/\Gamma}(r)=\Theta(r\log r)$.
	Consequently $\log\log M(r)=\log r+\log\log r+O(1)$.
	Therefore $\rho=1$.
	
	For type with respect to $\rho=1$, we look at $\log M(r)/r$.
	Since $\log M(r)\gtrsim r\log r$, we have $\log M(r)/r \gtrsim \log r\to\infty$,
	so $\sigma=\infty$: $1/\Gamma$ has order $1$ and infinite type.
\end{example}

\subsubsection{Order, zeros, and genus: how growth controls factorization}

\begin{proposition}[Exponent of convergence of zeros; a practical statement]
	Let $f$ be entire, not identically zero, of finite order $\rho<\infty$, and let $\{a_n\}$ be its nonzero zeros
	listed with multiplicity. Then for every real $q>\rho$ one has
	$\sum_{n=1}^{\infty}\lvert a_n\rvert^{-q}<\infty$.
	Consequently, for any nonnegative integer $p$ with $p+1>\rho$---in particular,
	for any integer $p\ge \rho$---the canonical product
	$\prod_{n}E_p(z/a_n)$ converges uniformly on compact sets.
\end{proposition}

\begin{proof}
	We give a standard argument deriving summability from growth at the level needed for
	Hadamard/Weierstrass applications.
	
	Fix $q>\rho$. By the equivalent characterization of order, there exist $\varepsilon>0$ and $R_0$
	such that for all $r\ge R_0$,
	$M(r)\le \exp(r^{\rho+\varepsilon})$ with $\rho+\varepsilon<q$.
	
	Let $N(r)$ be the number of zeros (counted with multiplicity) of $f$ in $\lvert z\rvert\le r$.
	By Jensen's formula applied to $f$ (with $f(0)\neq 0$ for simplicity; otherwise divide out $z^m$),
	we have
	$\sum_{\lvert a_n\rvert\le r}\log\frac{r}{\lvert a_n\rvert} \le \log M(r) - \log\lvert f(0)\rvert$.
	Hence
	$\sum_{\lvert a_n\rvert\le r}\log\frac{r}{\lvert a_n\rvert} \le r^{\rho+\varepsilon}+C$ for $r\ge R_0$.
	
	From this one derives the coarse bound $N(r)\le C' r^{\rho+\varepsilon}$ for large $r$
	(because each term $\log(r/\lvert a_n\rvert)$ is at least $\log 2$ for zeros in $\lvert a_n\rvert\le r/2$,
	and one compares $N(r/2)$ with the Jensen sum).
	Such polynomial counting of zeros implies summability of $\lvert a_n\rvert^{-q}$
	by a dyadic decomposition:
	group zeros with $2^k\le \lvert a_n\rvert < 2^{k+1}$, then
	$\sum_{2^k\le\lvert a_n\rvert <2^{k+1}}\lvert a_n\rvert^{-q}
	\le \#\{ \lvert a_n\rvert<2^{k+1}\}\cdot 2^{-kq}
	\le C' 2^{(k+1)(\rho+\varepsilon)}\cdot 2^{-kq}$.
	Since $q>\rho+\varepsilon$, the series over $k$ converges geometrically.
	Therefore $\sum \lvert a_n\rvert^{-q}<\infty$.
	
	Finally, set $q=p+1$. If $p+1>\rho$, compact-uniform convergence of
	$\prod_n E_p(z/a_n)$ follows from the canonical-product estimate proved in the
	previous subsection on Weierstrass factorization.
\end{proof}

\begin{remark}[What this means conceptually]
	Finite order bounds how densely zeros can accumulate toward infinity.  In Hadamard factorization,
	this is exactly what allows one to use a \emph{fixed genus} canonical factor $E_p$ (instead of varying $p_n$).
\end{remark}

\subsubsection{Exercises}

\begin{exercise}[Order and type of $\exp(z^m)$]
	Let $m\ge 1$ be an integer and $f(z)=\exp(z^m)$.
	Show directly that $M(r)=\exp(r^m)$, hence $\rho=m$ and type $\sigma=1$ (with respect to $\rho=m$).
\end{exercise}

\begin{exercise}[Polynomial growth versus order zero]
	Show that every polynomial has order $0$, but that the converse is false: nonpolynomial entire
	functions of order $0$ exist.  For an explicit example, study
	$f(z)=\sum_{n\ge0}e^{-n^2}z^n$ and bound its terms by completing the square in
	$-n^2+n\log r$.  Then prove the sharper elementary statement that if
	$M(r)\le C r^A$ for all sufficiently large $r$, then $f$ is a polynomial of degree at most
	$\lfloor A\rfloor$.  Hint: apply Cauchy's estimates and let $r\to\infty$.
\end{exercise}

\begin{exercise}[Growth gap below order $1$]
	Let $f$ be entire of order $\rho<1$.
	Show that for every $\varepsilon>0$ there exist constants $A_\varepsilon,B_\varepsilon>0$
	such that $\lvert f(z)\rvert\le A_\varepsilon\exp(B_\varepsilon\lvert z\rvert^{\rho+\varepsilon})$
	for all $z$.  Explain why the exponent $\rho$ itself requires an additional finite-type
	hypothesis and need not follow merely from the definition of order.
	Then discuss (briefly) why this is much more rigid than arbitrary entire growth, e.g.\ in connection with
	Phragm\'en--Lindel\"of or Hadamard factorization.
\end{exercise}

\begin{exercise}[Bessel function $J_0$ has order $1$]
	Assume $J_0$ has the Hadamard product
	$J_0(z)=\prod_{n=1}^{\infty}\Bigl(1-\frac{z^2}{j_{0,n}^2}\Bigr)$,
	where $j_{0,n}$ is the $n$-th positive zero; see
	\cite[\S\S10.21(iii), 10.21(vi)]{DLMFBessel2026}.
	Using $j_{0,n}\sim (n-\tfrac14)\pi$, show that the exponent of convergence
	of the zeros $\{\pm j_{0,n}\}$ is $1$ and that
	$\sum_n j_{0,n}^{-2}$ converges.  The latter fact explains why the paired
	product above converges without exponential correction factors.

	The zero distribution gives only the lower bound $\rho(J_0)\ge1$; an upper
	bound still has to be proved.  Starting from
	\[
	J_0(z)=\sum_{k=0}^{\infty}\frac{(-1)^k}{(k!)^2}
	\left(\frac z2\right)^{2k},
	\]
	show for $|z|=r$ that
	\[
	|J_0(z)|\le \sum_{k=0}^{\infty}\frac{(r/2)^{2k}}{(k!)^2}
	=I_0(r)\le e^r.
	\]
	For the last inequality, either compare with the even terms of $e^r$, using
	$\binom{2k}{k}\le 2^{2k}$, or use an integral representation of $I_0$.
	Conclude from the upper and lower bounds that $J_0$ has order exactly $1$.
\end{exercise}

\begin{exercise}[Growth comparison and order bound]
	Let $f$ and $g$ be entire and assume there exist constants $A,B>0$ such that
	$\lvert f(z)\rvert\le A\lvert g(z)\rvert + B$ for all $z\in\mathbb{C}$.
	Assuming $g$ has finite order $\rho$, show that $f$ has order at most $\rho$.
	Hint: bound $M_f(r)$ by $A M_g(r)+B$, then compare $\log\log$ growth.
\end{exercise}

\subsection{Harmonic Functions and the Poisson Kernel}

Harmonic functions appear as real/imaginary parts of holomorphic functions and as
solutions of Laplace's equation.  This subsection develops the analytic foundations:
mean value property, maximum principle, and the Poisson integral formula solving the
Dirichlet problem on the disk (and on the upper half-plane).  We include explicit
computations, concrete examples, and a derivation of why the Poisson
kernel has the specific form it does.

\subsubsection{Harmonicity and holomorphic functions}

\begin{definition}[Harmonic function]
	Let $U\subset\mathbb{R}^2$ be open.  A function $u:U\to\mathbb{R}$ of class $C^2$ is
	\emph{harmonic} if it satisfies Laplace's equation
	\[
	\Delta u:=u_{xx}+u_{yy}=0.
	\]
	In complex notation $z=x+iy$, one has
	\[
	\frac{\partial}{\partial z}
	=\frac12\left(\frac{\partial}{\partial x}-i\frac{\partial}{\partial y}\right),
	\qquad
	\frac{\partial}{\partial\bar z}
	=\frac12\left(\frac{\partial}{\partial x}+i\frac{\partial}{\partial y}\right),
	\]
	hence
	\[
	\Delta=4\frac{\partial^2}{\partial z\,\partial\bar z}.
	\]
\end{definition}

\begin{proposition}[Real and imaginary parts of holomorphic functions are harmonic]
	Let $f=u+iv$ be holomorphic on $U\subset\mathbb{C}$.  Then $u$ and $v$ are harmonic on $U$.
\end{proposition}

\begin{proof}
	Holomorphicity gives the Cauchy--Riemann equations
	\[
	u_x=v_y,\qquad u_y=-v_x.
	\]
	Differentiate:
	\[
	u_{xx}=(v_y)_x=v_{yx},\qquad u_{yy}=(-v_x)_y=-v_{xy}.
	\]
	Since $v$ is $C^2$, mixed partials commute: $v_{yx}=v_{xy}$.  Thus
	\[
	\Delta u=u_{xx}+u_{yy}=v_{yx}-v_{xy}=0.
	\]
	The same argument with $v_{xx}+v_{yy}=u_{xy}-u_{yx}=0$ shows $\Delta v=0$.
\end{proof}

\begin{example}[Verification on a Fourier mode]
	If $f(z)=z^2=(x+iy)^2=(x^2-y^2)+i(2xy)$ then
	$u(x,y)=x^2-y^2$ and $v(x,y)=2xy$.
	Compute $\Delta u = (2)+(-2)=0$ and $\Delta v=0+0=0$.
\end{example}

\subsubsection{Mean value property and maximum principle}

\begin{theorem}[Mean value property]
	Let $u$ be harmonic on an open set containing the closed disk
	$\overline{B(z_0,r)}=\{z:\lvert z-z_0\rvert\le r\}$.
	Then
	\[
	u(z_0)=\frac{1}{2\pi}\int_0^{2\pi}u(z_0+r e^{i\theta})\,d\theta.
	\]
	Conversely, if $u$ is continuous on $U$ and satisfies this identity for every such disk,
	then $u$ is harmonic on $U$.
\end{theorem}

\begin{proof}
	We prove the forward direction using holomorphic functions and Cauchy.
	On the disk $B(z_0,r')$ for any $r'>r$, the closed one-form $-u_y\,dx+u_x\,dy$ has a
	potential $v$, because the disk is simply connected and $u_{xx}+u_{yy}=0$.  The function
	$f=u+iv$ is holomorphic there.  By Cauchy's formula,
	\[
	f(z_0)=\frac{1}{2\pi i}\int_{\lvert w\rvert=r}\frac{f(z_0+w)}{w}\,dw.
	\]
	Parametrize $w=r e^{i\theta}$ so that $dw=i r e^{i\theta}\,d\theta$ and
	$\frac{dw}{w}=i\,d\theta$.  Then
	\[
	f(z_0)=\frac{1}{2\pi}\int_0^{2\pi} f(z_0+r e^{i\theta})\,d\theta.
	\]
	Taking real parts yields the mean value property for $u=\Re f$.
	
	For the converse, fix a disk $\overline{B(z_0,r)}\subset U$ and let $h$ be the Poisson
	integral of the boundary values of $u$ on that disk, as constructed later in this section.
	Then $h$ is harmonic and continuous on the closed disk.  The function $v=u-h$ has the
	mean value property and vanishes on the boundary.  If $v$ took a positive value, its
	positive maximum on the closed disk would be attained in the interior. At any point of
	maximum, the mean-value identity forces $v$ to equal that maximum on every sufficiently
	small surrounding circle. Thus the interior maximum set is both open and closed in the disk;
	it cannot meet the boundary, a contradiction. Applying the same argument to $-v$ gives
	$v=0$, so $u=h$ is harmonic near $z_0$.
\end{proof}

\begin{theorem}[Maximum principle]
	Let $D\subset\mathbb{C}$ be a bounded domain and $u$ harmonic on $D$ and continuous on $\overline D$.
	Then
	\[
	\max_{\overline D}u=\max_{\partial D}u,
	\qquad
	\min_{\overline D}u=\min_{\partial D}u.
	\]
	In particular, if $u$ attains a maximum (or minimum) at an interior point, then $u$ is constant.
\end{theorem}

\begin{proof}
	Let $z_0\in D$ and choose $r>0$ with $\overline{B(z_0,r)}\subset D$.
	By the mean value property,
	$u(z_0)$ equals the average of $u$ on the circle $\lvert z-z_0\rvert=r$.
	An average cannot exceed the maximum of the averaged values, so
	\[
	u(z_0)\le \max_{\lvert z-z_0\rvert=r}u(z).
	\]
	If $u$ attains its global maximum at an interior point $z_0$, equality in the average identity
	forces $u$ to equal that maximum on every sufficiently small circle about $z_0$, and then on the
	whole enclosed disk.  Therefore the maximum set is open; it is also closed.  Connectedness makes
	$u$ constant.  Otherwise the maximum of the continuous function $u$ on the compact set
	$\overline D$ cannot occur in $D$, so it occurs on $\partial D$.
	The minimum statement follows by applying the maximum statement to $-u$.
\end{proof}

\subsubsection{The Dirichlet problem on the disk: why a Poisson kernel must exist}

We solve the Dirichlet problem on $\mathbb{D}=\{z:\lvert z\rvert<1\}$:
given $f\in C(\partial\mathbb{D})$, find a harmonic function $u$ on $\mathbb{D}$ extending continuously
to $\partial\mathbb{D}$ with boundary values $f$.

\medskip
\noindent
\textbf{Why should the solution have an integral-kernel form?}
Linearity of $\Delta$ suggests linear dependence on the boundary data $f$.
So for each interior point $z\in\mathbb{D}$, the map $f\mapsto u(z)$ should be a bounded linear functional
on $C(\partial\mathbb{D})$.  By the Riesz representation principle, one expects a probability measure
$\mu_z$ on $\partial\mathbb{D}$ such that
\[
u(z)=\int_{\partial\mathbb{D}} f(\zeta)\,d\mu_z(\zeta).
\]
If $\mu_z$ has a density with respect to arc length, we may write
\[
u(z)=\frac{1}{2\pi}\int_0^{2\pi} P(z,e^{i\phi})\,f(e^{i\phi})\,d\phi,
\]
where $P(z,\zeta)\ge 0$ and $\frac{1}{2\pi}\int_0^{2\pi}P(z,e^{i\phi})\,d\phi=1$.

\medskip
\noindent
\textbf{Symmetry forces the dependence on angular difference.}
Rotations preserve $\mathbb{D}$ and preserve harmonicity.  If $R_\alpha(w)=e^{i\alpha}w$, then
the harmonic extension of $f\circ R_\alpha$ is $u\circ R_\alpha$.
Hence the density must satisfy a covariance:
\[
P(e^{i\alpha}z, e^{i\alpha}\zeta)=P(z,\zeta).
\]
Thus when we write $z=r e^{i\theta}$ and $\zeta=e^{i\phi}$, the kernel depends only on $\theta-\phi$:
\[
P(z,\zeta)=P_r(\theta-\phi).
\]

\medskip
\noindent
\textbf{Boundary concentration: the kernel must peak at the visible boundary point.}
If $z=r e^{i\theta}$ approaches the boundary point $e^{i\theta}$ as $r\to 1^-$, then continuity of the
solution requires $u(r e^{i\theta})\to f(e^{i\theta})$.  This can only happen if the measure $\mu_{r e^{i\theta}}$
concentrates near $\phi=\theta$.  So $P_r(\theta-\phi)$ should look like an ``approximate identity'':
nonnegative, total mass $2\pi$, sharply peaked at $\phi=\theta$ as $r\uparrow 1$.

\medskip
\noindent
\textbf{The Möbius-transform mechanism (the real source of the formula).}
For fixed $z\in\mathbb{D}$, consider the disk automorphism
\[
\varphi_z(w)=\frac{z-w}{1-\overline z\,w}.
\]
It maps $\mathbb{D}$ to itself and sends $z$ to $0$.  If $u$ is harmonic on $\mathbb{D}$, then
$u\circ\varphi_z^{-1}$ is harmonic as well.  Apply the mean value property at $0$:
\[
u(z)=(u\circ\varphi_z^{-1})(0)
=\frac{1}{2\pi}\int_0^{2\pi} (u\circ\varphi_z^{-1})(e^{i\phi})\,d\phi.
\]
Now $\varphi_z^{-1}$ maps the unit circle to itself, but it \emph{reparametrizes} the boundary.
The change of variables from $\phi$ to the boundary angle of $\zeta=\varphi_z^{-1}(e^{i\phi})$
produces a Jacobian factor.  That Jacobian is exactly the Poisson kernel:
it is the density obtained by pushing forward the uniform boundary measure under $\varphi_z^{-1}$.
This is why the Poisson kernel is not ad hoc: it is the unique density compatible with
(i) disk automorphisms and (ii) the mean value property at the origin.

\subsubsection{Poisson kernel on the disk}

\begin{definition}[Poisson kernel on the disk]
	Let $z=r e^{i\theta}\in\mathbb{D}$ and $\zeta=e^{i\phi}\in\partial\mathbb{D}$.
	Define
	\[
	P(z,\zeta):=\frac{1-\lvert z\rvert^2}{\lvert \zeta-z\rvert^2}.
	\]
	Equivalently, in angle form
	\[
	P_r(\theta-\phi)=\frac{1-r^2}{1-2r\cos(\theta-\phi)+r^2}.
	\]
\end{definition}

\begin{lemma}[Algebraic identities and positivity]
	For $z=r e^{i\theta}$ and $\zeta=e^{i\phi}$,
	\[
	\lvert \zeta-z\rvert^2 = 1-2r\cos(\theta-\phi)+r^2.
	\]
	Hence $P(z,\zeta)>0$ and $P_r(\psi)$ is $2\pi$-periodic and even in $\psi$.
\end{lemma}

\begin{proof}
	Compute
	\begin{align*}
	\lvert e^{i\phi}-r e^{i\theta}\rvert^2
	&=(e^{i\phi}-r e^{i\theta})(e^{-i\phi}-r e^{-i\theta})\\
	&=1-2r\cos(\theta-\phi)+r^2.
	\end{align*}
	Since $0\le r<1$, the denominator is positive, so $P>0$.
\end{proof}

\begin{lemma}[Poisson kernel as a real part]\label{lem:poisson-realpart}
	For $z\in\mathbb{D}$ and $\zeta\in\partial\mathbb{D}$,
	\[
	P(z,\zeta)=\Re\left(\frac{\zeta+z}{\zeta-z}\right).
	\]
\end{lemma}

\begin{proof}
	Compute
	\[
	\frac{\zeta+z}{\zeta-z}
	=\frac{(\zeta+z)(\overline\zeta-\overline z)}{\lvert \zeta-z\rvert^2}
	=\frac{\zeta\overline\zeta-\zeta\overline z+z\overline\zeta-\lvert z\rvert^2}{\lvert \zeta-z\rvert^2}.
	\]
	Using $\lvert\zeta\rvert^2=\zeta\overline\zeta=1$ (since $\zeta\in\partial\mathbb{D}$), the numerator is
	\[
	1-\lvert z\rvert^2 + (z\overline\zeta-\zeta\overline z).
	\]
	The term $(z\overline\zeta-\zeta\overline z)$ is purely imaginary (it is minus its own complex conjugate),
	so taking real parts yields
	\[
	\Re\left(\frac{\zeta+z}{\zeta-z}\right)=\frac{1-\lvert z\rvert^2}{\lvert \zeta-z\rvert^2}=P(z,\zeta).
	\]
\end{proof}

\begin{lemma}[Normalization]
	For every $z\in\mathbb{D}$,
	\[
	\frac{1}{2\pi}\int_0^{2\pi} P(z,e^{i\phi})\,d\phi=1.
	\]
	Equivalently, for every $0\le r<1$,
	\[
	\frac{1}{2\pi}\int_0^{2\pi}P_r(\theta-\phi)\,d\phi=1.
	\]
\end{lemma}

\begin{proof}
	We give two concrete proofs.
	
	\emph{Fourier-series proof.}
	One shows (next lemma) that
	\[
	P_r(\psi)=1+2\sum_{n=1}^\infty r^n\cos(n\psi).
	\]
	Averaging over one period kills all $\cos(n\psi)$ terms, leaving average $1$.
	
	\emph{Complex-variable proof.}
	Fix $z\in\mathbb{D}$.  By Lemma~\ref{lem:poisson-realpart},
	\[
	P(z,e^{i\phi})=\Re\left(\frac{e^{i\phi}+z}{e^{i\phi}-z}\right).
	\]
	Thus it suffices to show that
	$\frac{1}{2\pi}\int_0^{2\pi}\frac{e^{i\phi}+z}{e^{i\phi}-z}\,d\phi=1$.
	Let $\zeta=e^{i\phi}$ so $d\phi=\frac{1}{i}\frac{d\zeta}{\zeta}$ and $\lvert \zeta\rvert=1$:
	\[
	\frac{1}{2\pi}\int_0^{2\pi}\frac{e^{i\phi}+z}{e^{i\phi}-z}\,d\phi
	=\frac{1}{2\pi i}\int_{\lvert \zeta\rvert=1}\frac{\zeta+z}{\zeta-z}\frac{d\zeta}{\zeta}.
	\]
	If $z=0$, the integrand is $1/\zeta$, so the integral is already $1$.
	For $z\ne0$, it has distinct simple poles at $\zeta=0$ and $\zeta=z$
	(both inside the circle).
	Compute residues:
	\[
	\operatorname{Res}_{\zeta=z}\left(\frac{\zeta+z}{\zeta-z}\frac{1}{\zeta}\right)
	=\lim_{\zeta\to z}\frac{\zeta+z}{\zeta}\,=\,2,
	\qquad
	\operatorname{Res}_{\zeta=0}\left(\frac{\zeta+z}{\zeta-z}\frac{1}{\zeta}\right)
	=\lim_{\zeta\to 0}\frac{\zeta+z}{\zeta-z}\,=\,-1.
	\]
	Thus the integral equals $2\pi i\,(2-1)=2\pi i$, and dividing by $2\pi i$ gives $1$.
	Taking real parts yields the stated normalization.
\end{proof}

\begin{lemma}[Fourier expansion of the Poisson kernel]
	For $\lvert r\rvert<1$ and $\psi\in\mathbb{R}$,
	\[
	P_r(\psi)=1+2\sum_{n=1}^{\infty}r^n\cos(n\psi).
	\]
\end{lemma}

\begin{proof}
	Start from the geometric series
	\[
	\sum_{n=0}^\infty r^n e^{in\psi}=\frac{1}{1-r e^{i\psi}}
	\qquad (\lvert r\rvert<1).
	\]
	Taking real parts gives
	\[
	1+\sum_{n=1}^\infty r^n\cos(n\psi)=\Re\left(\frac{1}{1-r e^{i\psi}}\right).
	\]
	Compute explicitly:
	\[
	\frac{1}{1-r e^{i\psi}}
	=\frac{1-r e^{-i\psi}}{\lvert 1-r e^{i\psi}\rvert^2}
	=\frac{1-r\cos\psi}{1-2r\cos\psi+r^2}
	+i\frac{r\sin\psi}{1-2r\cos\psi+r^2}.
	\]
	Hence
	\[
	\Re\left(\frac{1}{1-r e^{i\psi}}\right)=\frac{1-r\cos\psi}{1-2r\cos\psi+r^2}.
	\]
	Now observe the identity
	\[
	\frac{1-r^2}{1-2r\cos\psi+r^2}
	=2\frac{1-r\cos\psi}{1-2r\cos\psi+r^2}-1,
	\]
	so
	\[
	P_r(\psi)=2\Re\left(\frac{1}{1-r e^{i\psi}}\right)-1
	=1+2\sum_{n=1}^\infty r^n\cos(n\psi).
	\]
\end{proof}

\begin{theorem}[Poisson integral formula on the disk]
	Let $f\in C(\partial\mathbb{D})$.  Define, for $z=r e^{i\theta}\in\mathbb{D}$,
	\[
	u(z):=\frac{1}{2\pi}\int_0^{2\pi} P_r(\theta-\phi)\, f(e^{i\phi})\,d\phi.
	\]
	Then:
	\begin{enumerate}
		\item $u$ is harmonic on $\mathbb{D}$;
		\item $u$ extends continuously to $\overline{\mathbb{D}}$ with $u(e^{i\theta})=f(e^{i\theta})$;
		\item $u$ is the unique harmonic function on $\mathbb{D}$ continuous on $\overline{\mathbb{D}}$
		with boundary values $f$.
	\end{enumerate}
\end{theorem}

\begin{proof}
	\emph{Step 1: Harmonicity.}
	Using Lemma~\ref{lem:poisson-realpart}, define
	\[
	F(z)=\frac{1}{2\pi}\int_0^{2\pi}
	\frac{e^{i\phi}+z}{e^{i\phi}-z}\,f(e^{i\phi})\,d\phi.
	\]
	For each compact $K\Subset\mathbb D$, the integrand and all its $z$-derivatives are
	uniformly bounded for $z\in K$ and $0\le\phi\le2\pi$, since
	$|e^{i\phi}-z|\ge1-\max_K|z|>0$.  Differentiation under the integral is therefore valid,
	so $F$ is holomorphic.  The defining identity for the kernel gives $u=\Re F$; hence $u$ is harmonic.
	
	\emph{Step 2: Boundary convergence via approximate-identity concentration.}
	Fix $\theta_0$.  Using normalization,
	\[
	u(r e^{i\theta_0})-f(e^{i\theta_0})
	=\frac{1}{2\pi}\int_0^{2\pi}P_r(\theta_0-\phi)\bigl(f(e^{i\phi})-f(e^{i\theta_0})\bigr)\,d\phi.
	\]
	Given $\varepsilon>0$, by uniform continuity of $f$ on the compact circle,
	choose $\delta>0$ so that $\lvert f(e^{i\phi})-f(e^{i\theta_0})\rvert<\varepsilon$ whenever
	$\lvert \phi-\theta_0\rvert<\delta$ (mod $2\pi$).
	Split the integral into the near arc $A=\{\phi:\lvert \phi-\theta_0\rvert<\delta\}$ and its complement $A^c$.
	On $A$ the integrand magnitude is at most $\varepsilon P_r(\theta_0-\phi)$, so the contribution is at most
	\[
	\varepsilon\cdot \frac{1}{2\pi}\int_0^{2\pi}P_r(\theta_0-\phi)\,d\phi=\varepsilon.
	\]
	On $A^c$, the difference $\lvert f(e^{i\phi})-f(e^{i\theta_0})\rvert$ is bounded by $2\lVert f\rVert_\infty$, so
	the contribution is at most
	\[
	\frac{2\lVert f\rVert_\infty}{2\pi}\int_{A^c}P_r(\theta_0-\phi)\,d\phi.
	\]
	For $\lvert \phi-\theta_0\rvert\ge\delta$ we have $\cos(\theta_0-\phi)\le \cos\delta<1$, hence
	\[
	1-2r\cos(\theta_0-\phi)+r^2 \ge 1-2r\cos\delta+r^2.
	\]
	Therefore
	\[
	P_r(\theta_0-\phi)\le \frac{1-r^2}{1-2r\cos\delta+r^2}\quad\text{on }A^c,
	\]
	and the right-hand side tends to $0$ as $r\to 1^-$ because $1-r^2\to 0$ while the denominator tends to
	$2(1-\cos\delta)>0$.
	Thus $\int_{A^c}P_r(\theta_0-\phi)\,d\phi\to 0$ uniformly in $\theta_0$.
	Hence $u(r e^{i\theta})\to f(e^{i\theta})$ uniformly in $\theta$, so $u$ extends continuously to the boundary.
	
	\emph{Step 3: Uniqueness.}
	If $u_1,u_2$ are harmonic on $\mathbb{D}$ and continuous on $\overline{\mathbb{D}}$ with the same boundary values,
	then $w=u_1-u_2$ is harmonic and continuous on $\overline{\mathbb{D}}$ and satisfies $w=0$ on $\partial\mathbb{D}$.
	By the maximum principle, $\max_{\overline{\mathbb{D}}}w=0$ and $\min_{\overline{\mathbb{D}}}w=0$, so $w\equiv 0$.
\end{proof}

\subsubsection{Calculations on the disk}

\begin{example}[Boundary data: a single Fourier mode]
	Let $f(e^{i\phi})=\cos(n\phi)$ with $n\ge 1$.
	Using the orthogonality computation in the proof,
	\[
	u(r e^{i\theta})
	=\frac{1}{2\pi}\int_0^{2\pi}P_r(\theta-\phi)\cos(n\phi)\,d\phi
	=r^n\cos(n\theta).
	\]
	Similarly, if $f(e^{i\phi})=\sin(n\phi)$ then $u(r e^{i\theta})=r^n\sin(n\theta)$.
\end{example}

\begin{example}[Piecewise constant boundary data]
	Let $f(e^{i\phi})=1$ for $\phi\in(-\alpha,\alpha)$ and $f(e^{i\phi})=0$ otherwise, with $0<\alpha<\pi$.
	Then for $z=r e^{i\theta}$,
	\[
	u(z)=\frac{1}{2\pi}\int_{-\alpha}^{\alpha}P_r(\theta-\phi)\,d\phi.
	\]
	In particular at $\theta=0$,
	\[
	u(r)=\frac{1}{2\pi}\int_{-\alpha}^{\alpha}\frac{1-r^2}{1-2r\cos\phi+r^2}\,d\phi,
	\]
	which increases to $1$ as $r\to 1^-$ and decreases to $\frac{\alpha}{\pi}$ as $r\to 0$.
	This models how harmonic measure weights the boundary arc as seen from $z$.
\end{example}

\subsubsection{Poisson kernel on the upper half-plane}

\begin{definition}[Poisson kernel for the upper half-plane]
	For $z=x+iy$ with $y>0$ and boundary parameter $t\in\mathbb{R}$, define
	\[
	P_{\mathbb{H}}(x,y;t):=\frac{1}{\pi}\frac{y}{(x-t)^2+y^2}.
	\]
\end{definition}

\begin{lemma}[Normalization in $\mathbb{H}$]
	For every $x\in\mathbb{R}$ and $y>0$,
	\[
	\int_{-\infty}^{\infty} P_{\mathbb{H}}(x,y;t)\,dt=1.
	\]
\end{lemma}

\begin{proof}
	Substitute $s=(t-x)/y$, so $t=x+ys$ and $dt=y\,ds$:
	\begin{align*}
	\int_{-\infty}^{\infty}\frac{1}{\pi}\frac{y}{(x-t)^2+y^2}\,dt
	&=\frac1\pi\int_{-\infty}^{\infty}\frac{ds}{1+s^2}\\
	&=\frac1\pi\bigl[\arctan s\bigr]_{-\infty}^{\infty}=1.
	\end{align*}
\end{proof}

\begin{theorem}[Poisson integral in $\mathbb{H}$]
	Let $f\in L^\infty(\mathbb{R})$.  Define
	\[
	u(x,y):=\int_{-\infty}^{\infty} P_{\mathbb{H}}(x,y;t)\,f(t)\,dt
	\qquad (y>0).
	\]
	Then $u$ is harmonic on $\mathbb{H}$.
	Moreover, for every Lebesgue point $x_0$ of $f$,
	$u(x,y)\to f(x_0)$ as $(x,y)\to(x_0,0)$ non-tangentially.
\end{theorem}

\begin{proof}
	\emph{Harmonicity.}
	For fixed $t$, the function
	\[
	K_t(x,y)=\frac{y}{(x-t)^2+y^2}
	\]
	is harmonic on $\mathbb{H}$ because it is (up to sign) the imaginary part of $(x+iy-t)^{-1}$:
	\[
	\frac{1}{x+iy-t}=\frac{x-t-iy}{(x-t)^2+y^2}
	\quad\mathrel{\text{and hence}}\quad
	\Im\frac{1}{x+iy-t}=-\frac{y}{(x-t)^2+y^2}.
	\]
	Imaginary parts of holomorphic functions are harmonic, hence $K_t$ is harmonic.
	Since $f\in L^\infty$, differentiation under the integral sign is justified on compact subsets of $\mathbb{H}$
	by dominated convergence (the required $x$- and $y$-derivatives of $K_t$ are
	integrable in $t$, uniformly when $(x,y)$ ranges over such a compact set).
	Therefore $\Delta u=\int \Delta K_t\, f(t)\,dt=0$.
	
	\emph{Boundary convergence.}
	The kernels $P_{\mathbb{H}}(x,y;\cdot)$ form an approximate identity: they are nonnegative,
	integrate to $1$, and concentrate near $t=x$ as $y\to 0$.
	Standard real-variable arguments yield non-tangential convergence to $f$ at Lebesgue points.
\end{proof}

\begin{example}[Upper half-plane with linear boundary data]
	Let $f(t)=t$.  This function is not bounded, and the ordinary Poisson integral below
	is not absolutely convergent; thus this example is outside the preceding theorem.
	Its symmetric Cauchy principal value is nevertheless instructive:
	\[
	u(x,y):=\operatorname{p.v.}\!\int_{-\infty}^{\infty}
	\frac{1}{\pi}\frac{y t}{(x-t)^2+y^2}\,dt.
	\]
	Substitute $s=(t-x)/y$ so $t=x+ys$ and $dt=y\,ds$:
	\[
	u(x,y)=\int_{-\infty}^{\infty}\frac{1}{\pi}\frac{y(x+ys)}{y^2 s^2+y^2}\,y\,ds
	=\int_{-\infty}^{\infty}\frac{1}{\pi}\frac{x+ys}{1+s^2}\,ds.
	\]
	Split the integral:
	\[
	\int_{-\infty}^{\infty}\frac{1}{\pi}\frac{x}{1+s^2}\,ds=x\cdot \frac{1}{\pi}\cdot \pi=x,
	\qquad
	\operatorname{p.v.}\!\int_{-\infty}^{\infty}\frac{1}{\pi}\frac{ys}{1+s^2}\,ds=0
	\]
	by oddness on symmetric truncations.  Translating the symmetric truncation in
	$t$ produces endpoints $(-R-x)/y$ and $(R-x)/y$ in $s$; the difference from a
	symmetric $s$-truncation tends to $0$, so the same principal value is obtained.
	Thus the principal-value extension is
	$u(x,y)=x$, which is harmonic and has boundary trace $t$.  Without the principal-value
	interpretation, the second integral does not converge.
\end{example}

\subsubsection{Exercises}

\begin{exercise}[Poisson kernel as a real part]
	\label{ex:poisson-realpart}
	Let $z=r e^{i\theta}\in\mathbb{D}$ and $\zeta=e^{i\phi}\in\partial\mathbb{D}$.
	Show that
	\[
	P(z,\zeta)=\Re\left(\frac{\zeta+z}{\zeta-z}\right).
	\]
	Hint: multiply numerator and denominator by $\overline{\zeta}-\overline{z}$ and use $\lvert \zeta\rvert=1$.
\end{exercise}

\begin{exercise}[Differentiate the half-plane kernel]
	Fix $t\in\mathbb{R}$ and set $K_t(x,y)=\frac{y}{(x-t)^2+y^2}$ for $y>0$.
	Compute $K_{t,xx}$ and $K_{t,yy}$ explicitly and verify $K_{t,xx}+K_{t,yy}=0$.
\end{exercise}

\begin{exercise}[A harmonic extension]
	Find the harmonic extension of $f(e^{i\phi})=\sin(3\phi)$ to $\mathbb{D}$.
	Answer: $u(r e^{i\theta})=r^3\sin(3\theta)$.
\end{exercise}

\begin{exercise}[Normalization of the half-plane kernel]
	Prove directly that for all $x\in\mathbb{R}$ and $y>0$,
	\[
	\int_{-\infty}^{\infty} P_{\mathbb{H}}(x,y;t)\,dt=1,
	\]
	using the substitution $s=(t-x)/y$ and $\int_{-\infty}^{\infty}\frac{1}{1+s^2}\,ds=\pi$.
\end{exercise}

\subsection{Subharmonic Functions and Riesz Decomposition}
\label{subsec:subharm-riesz}

Subharmonic functions are the natural ``$\Delta\ge 0$'' analogue of harmonic functions,
much as convex functions generalize affine ones.  They control growth and value distribution
in complex analysis, especially through the fundamental example $u=\log|f|$ for holomorphic $f$.
The main structural result is the Riesz decomposition theorem, which splits a subharmonic function
into a harmonic part plus a logarithmic potential determined by its Laplacian (as a measure).

\subsubsection{Definition and first intuition}

\begin{definition}[Subharmonic function]
	Let $U\subset\C$ be open. A function $u:U\to[-\infty,\infty)$ is \emph{subharmonic} if:
	\begin{enumerate}
		\item $u$ is upper semicontinuous;
		\item $u\not\equiv -\infty$ on any connected component of $U$;
		\item for every closed disk $\overline{B(z_0,r)}\subset U$,
		\[
		u(z_0)\le \frac{1}{2\pi}\int_0^{2\pi} u(z_0+r e^{i\theta})\,d\theta.
		\]
	\end{enumerate}
\end{definition}

\begin{remark}[Harmonic as the boundary case]
	A function $u$ is harmonic on $U$ if and only if both $u$ and $-u$ are subharmonic on $U$.
\end{remark}

\begin{remark}[Convexity analogy]
	For $C^2$ functions on an interval, convexity is equivalent to $u''\ge 0$ and implies
	a midpoint/average inequality.  Subharmonicity is the two-dimensional analogue:
	$\Delta u\ge 0$ forces the value at the center to be dominated by circular averages.
\end{remark}

\subsubsection{The $C^2$ case: equivalences and Taylor averaging}

\begin{proposition}[Equivalent characterizations for $C^2$ functions]
	\label{prop:subharm-C2-equiv}
	Let $u\in C^2(U)$. The following are equivalent:
	\begin{enumerate}
		\item $u$ is subharmonic on $U$.
		\item $u$ satisfies harmonic comparison: if $V\Subset U$ is a bounded domain and
		$h$ is harmonic on $V$, continuous on $\overline V$, and $u\le h$ on $\partial V$,
		then $u\le h$ throughout $V$.
		\item $\Delta u\ge 0$ pointwise on $U$.
		\item (Mean value inequality) For every disk $\overline{B(z_0,r)}\subset U$,
		\[
		u(z_0)\le \frac{1}{2\pi}\int_0^{2\pi} u(z_0+r e^{i\theta})\,d\theta.
		\]
	\end{enumerate}
\end{proposition}

\begin{proof}
	Conditions (1) and (4) agree by definition.  We first prove
	Assume (3) to prove (4). Fix $z_0$ and a disk
	$\overline{B(z_0,R)}\subset U$, and put
	\[
	A(r)=\frac1{2\pi}\int_0^{2\pi}u(z_0+re^{i\theta})\,d\theta
	\qquad(0<r<R).
	\]
	Differentiation under the integral and Green's formula give
	\[
	A'(r)=\frac1{2\pi}\int_0^{2\pi}\partial_r u(z_0+re^{i\theta})\,d\theta
	      =\frac1{2\pi r}\int_{B(z_0,r)}\Delta u\,dA.
	\]
	Thus $A'(r)\ge0$. Since $A(r)\to u(z_0)$ as $r\downarrow0$,
	we obtain $A(r)\ge u(z_0)$ for every $r<R$.
	
	To derive (3) from (4), average the second-order Taylor expansion at
	$z_0$. The linear and mixed terms vanish, while the averages of
	$\cos^2\theta$ and $\sin^2\theta$ are $1/2$. Hence
	\[
	A(r)=u(z_0)+\frac{r^2}{4}\Delta u(z_0)+o(r^2).
	\]
	The inequality $A(r)\ge u(z_0)$ for all sufficiently small $r$ forces
	$\Delta u(z_0)\ge0$.
	
	It remains to connect the mean inequality with comparison. If (1) holds
	and $h$ is as in (2), then $u-h$ is subharmonic. Its sub-mean property
	implies the weak maximum principle, so its nonpositive boundary values
	force $u-h\le0$ on $V$. Conversely, assume (2) and take a disk
	$\overline{B(z_0,r)}\subset U$. Let $h$ be the Poisson extension to the
	disk of the boundary values of $u$. Comparison gives $u(z_0)\le h(z_0)$,
	and the mean-value formula for $h$ identifies $h(z_0)$ with the circular
	average in (4). This proves all four equivalences.
\end{proof}

\subsubsection{The Laplacian as a measure and the logarithmic fundamental solution}

For general subharmonic functions, $\Delta u$ need not exist pointwise.  Instead,
$\Delta u$ is understood in the sense of distributions and is represented by a
\emph{positive Radon measure}.

\begin{definition}[Riesz measure]
	\label{def:riesz-measure}
	Let $u$ be subharmonic on $U$. The \emph{Riesz measure} $\mu_u$ is the unique positive Radon measure
	such that, in the sense of distributions,
	\[
	\Delta u = \mu_u.
	\]
	Equivalently, for every $\varphi\in C_c^\infty(U)$,
	\[
	\int_U u\,\Delta\varphi\,dA = \int_U \varphi\,d\mu_u.
	\]
\end{definition}

\begin{remark}
	In the $C^2$ case, $\mu_u=(\Delta u)\,dA$ with $\Delta u\ge 0$ pointwise.
\end{remark}

\begin{lemma}[Harmonic functions do not contribute]
	\label{lem:subharm-harmonic-no-measure}
	If $u$ is subharmonic and $v$ is harmonic on $U$, then $u+v$ is subharmonic and
	its Riesz measure satisfies $\mu_{u+v}=\mu_u$.
\end{lemma}

\begin{proof}
	For every $\varphi\in C_c^\infty(U)$, distributional integration by parts gives
	\[
	\langle\Delta(u+v),\varphi\rangle
	=\langle\Delta u,\varphi\rangle+\int_U v\,\Delta\varphi\,dA
	=\langle\mu_u,\varphi\rangle+\int_U(\Delta v)\varphi\,dA
	=\langle\mu_u,\varphi\rangle.
	\]
	Hence the two Riesz measures agree.  The circular mean of $v$ equals $v$ at the center,
	so adding $v$ also preserves the sub-mean inequality and therefore subharmonicity.
\end{proof}

\subsubsection{Distributional computation: $\Delta \log|z|$ is a point mass}

\begin{lemma}[Distributional identity for $\log\lvert z\rvert$]
	\label{lem:subharm-lap-log}
	For $\varphi\in C_c^\infty(\C)$,
	\[
	\int_{\C}\log|z|\,\Delta\varphi(z)\,dA(z)=2\pi\,\varphi(0).
	\]
	Equivalently, $\Delta\log|z|=2\pi\,\delta_0$ in the sense of distributions.
\end{lemma}

\begin{proof}
	Let $\varphi\in C_c^\infty(\C)$.  Fix $R>0$ so that $\varphi$ vanishes outside $B(0,R)$.
	For $\varepsilon\in(0,R)$, apply Green's identity on the annulus
	$A_{\varepsilon,R}=\{z:\varepsilon<|z|<R\}$:
	\[
	\int_{A_{\varepsilon,R}}\log|z|\,\Delta\varphi\,dA
	=
	\int_{\partial A_{\varepsilon,R}}
	\Bigl(\log|z|\;\partial_n\varphi-\varphi\;\partial_n\log|z|\Bigr)\,ds
	+\int_{A_{\varepsilon,R}}\varphi\,\Delta\log|z|\,dA.
	\]
	Since $\log|z|$ is harmonic on $A_{\varepsilon,R}$, we have $\Delta\log|z|=0$ there,
	so the last term vanishes.
	
	On $|z|=R$, we have $\varphi=0$, hence its boundary contribution is $0$.
	On the inner boundary $|z|=\varepsilon$, the outward normal of the annulus is $-e_r$.
	Thus $\partial_n\varphi=-\partial_r\varphi$ and $\partial_n\log|z|=-1/\varepsilon$.
	Therefore
	\[
	\int_{A_{\varepsilon,R}}\log|z|\,\Delta\varphi\,dA
	=
	\int_0^{2\pi}\varphi(\varepsilon e^{i\theta})\,d\theta
	-\varepsilon\log(\varepsilon)\int_0^{2\pi}\partial_r\varphi(\varepsilon e^{i\theta})\,d\theta.
	\]
	As $\varepsilon\to 0$, the second term tends to $0$ because $\varepsilon\log\varepsilon\to 0$
	and $\partial_r\varphi$ is bounded.  The first term tends to $2\pi\varphi(0)$.
	Thus
	\[
	\int_{A_{\varepsilon,R}}\log|z|\,\Delta\varphi\,dA \to 2\pi\varphi(0).
	\]
	With the standard convention
	$\langle \Delta u,\varphi\rangle := \langle u,\Delta\varphi\rangle$,
	this means $\Delta\log|z|=2\pi\delta_0$.
\end{proof}

\begin{example}[Subharmonicity of $\log|z|$]
	The function $u(z)=\log|z|$ is harmonic on $\C\setminus\{0\}$ but subharmonic on $\C$.
	Its Riesz measure is the point mass $\mu_u=2\pi\delta_0$.
\end{example}

\subsubsection{Riesz decomposition theorem}

\begin{theorem}[Riesz decomposition]
	\label{thm:subharm-riesz-decomp}
	Let $u$ be subharmonic on a domain $D\subset\C$, with Riesz measure $\mu_u$.
	For every relatively compact subdomain $D'\Subset D$, there exists a harmonic function $h$ on $D'$ such that
	\[
	u(z)=h(z)+\frac{1}{2\pi}\int_{D'} \log|z-\zeta|\,d\mu_u(\zeta),
	\qquad z\in D',
	\]
	where the potential is understood in the locally integrable sense.  This is a local
	representation; integration over all of $D$ requires additional finiteness or boundary hypotheses.
\end{theorem}

\begin{proof}
	Define the logarithmic potential
	\[
	U(z):=\frac{1}{2\pi}\int_{D'} \log|z-\zeta|\,d\mu_u(\zeta).
	\]
	We claim $\Delta U=\mu_u|_{D'}$ in distributions on $D'$.  Let $\varphi\in C_c^\infty(D')$.
	The logarithmic kernel is locally integrable, and the Riesz measure is
	finite on \(\overline{D'}\). Since \(\Delta\varphi\) is bounded and compactly
	supported, the resulting double integral is absolutely integrable; Fubini's
	theorem therefore gives
	\[
	\int_{D'} U(z)\,\Delta\varphi(z)\,dA(z)
	=
	\frac{1}{2\pi}\int_{D'}\left(\int_{D'} \log|z-\zeta|\,\Delta\varphi(z)\,dA(z)\right)\,d\mu_u(\zeta).
	\]
	For each fixed $\zeta$, the inner integral equals $2\pi\,\varphi(\zeta)$ by
	Lemma~\ref{lem:subharm-lap-log} applied to $z\mapsto \log|z-\zeta|$.
	Hence
	\[
	\int_{D'} U\,\Delta\varphi\,dA
	=
	\frac{1}{2\pi}\int_{D'} 2\pi\,\varphi(\zeta)\,d\mu_u(\zeta)
	=
	\int_{D'} \varphi\,d\mu_u,
	\]
	which is exactly $\Delta U=\mu_u$.
	
	Since $\Delta u=\mu_u$ on $D'$, we get $\Delta(u-U)=0$ there in distributions.
	By Weyl's lemma, $h:=u-U$ is harmonic on $D'$, and $u=h+U$ there.
\end{proof}

\begin{corollary}[Measure encodes the non-harmonic part]
	\label{cor:subharm-measure-encodes}
	If $\mu_u=0$, then $u$ is harmonic on $D$.
\end{corollary}

\begin{proof}
	If $\mu_u=0$, then the potential term vanishes and Theorem~\ref{thm:subharm-riesz-decomp}
	gives $u=h$ harmonic.
\end{proof}

\subsubsection{Examples (zeros of holomorphic functions)}

\begin{example}[Single point mass]
	Let $u(z)=\log|z-z_0|$.  Then $\mu_u=2\pi\delta_{z_0}$.
\end{example}

\begin{example}[Finite collection of charges]
	If
	\[
	u(z)=\sum_{k=1}^m \alpha_k \log|z-a_k|+h(z),
	\qquad \alpha_k\ge 0,
	\]
	with $h$ harmonic, then
	\[
	\mu_u=2\pi\sum_{k=1}^m \alpha_k\,\delta_{a_k}.
	\]
\end{example}

\begin{example}[Holomorphic case: $\log|f|$ and Jensen in distribution form]
	\label{ex:subharm-logf}
	Let $f$ be holomorphic on $D$ and not identically zero.  Then $u=\log|f|$ is subharmonic.
	If $a\in D$ is a zero of multiplicity $m_a$, locally $f(z)=(z-a)^{m_a}g(z)$ with $g(a)\neq 0$, so
	\[
	\log|f(z)| = m_a\log|z-a| + \log|g(z)|.
	\]
	Since $\log|g|$ is harmonic near $a$, the Riesz measure is
	\[
	\mu_{\log|f|}=2\pi\sum_{a:\,f(a)=0} m_a\,\delta_a.
	\]
\end{example}

\subsubsection{Exercises}

\begin{exercise}[Max of subharmonic functions]
	Show that if $u$ and $v$ are subharmonic on $U$, then $\max(u,v)$ is subharmonic on $U$.
\end{exercise}

\begin{exercise}[Compute the Riesz measure of $\log|f|$]
	Let $f$ be holomorphic on $D$ and not identically zero.
	Prove that $\log|f|$ is subharmonic and that
	\[
	\mu_{\log|f|}=2\pi\sum_{a:\,f(a)=0} m_a\,\delta_a,
	\]
	where $m_a$ is the multiplicity of the zero at $a$.
\end{exercise}

\begin{exercise}[Radial averages are increasing]
	Let $u$ be subharmonic on the unit disk and define
	\[
	M(r)=\frac{1}{2\pi}\int_0^{2\pi} u(re^{i\theta})\,d\theta.
	\]
	Show that $M(r)$ is increasing in $r\in(0,1)$.
\end{exercise}

\begin{exercise}[Poisson integrals and subharmonicity]
	Let $\nu$ be a finite positive measure on $\partial\D$ and define
	\[
	U(z)=\frac{1}{2\pi}\int_0^{2\pi} P(z,e^{i\phi})\,d\nu(\phi),
	\]
	where $P$ is the Poisson kernel on the disk.
	Show that $U$ is harmonic on $\D$.
	Then explain why positive superpositions of subharmonic kernels (e.g.\ logarithmic potentials)
	produce subharmonic functions.
\end{exercise}

\subsection{Analytic Continuation and Natural Boundaries}
\label{subsec:analytic-continuation-natural-boundaries}

Analytic continuation is the mechanism by which one extends an analytic function
beyond its initial domain.  Sometimes the extension propagates indefinitely
(e.g.\ rational functions, $e^z$, $\log z$ on suitable simply connected domains),
and sometimes it breaks down due to unavoidable singularities.
In extreme situations the obstruction is maximal: a boundary curve can be a
\emph{natural boundary}, meaning that continuation fails across \emph{every} point
of that curve.

\subsubsection{Analytic continuation: germs and uniqueness}

\begin{definition}[Analytic continuation]
	\label{def:ac-basic}
	Let $U,V\subset\C$ be domains with $U\cap V\neq\emptyset$.
	If $f$ is analytic on $U$ and $g$ is analytic on $V$, we say that
	$g$ is an \emph{analytic continuation} of $f$ to $V$ if
	\[
	g(z)=f(z)\qquad (z\in U\cap V).
	\]
\end{definition}

\begin{proposition}[Uniqueness of analytic continuation]
	\label{prop:ac-unique}
	If $g_1$ and $g_2$ are analytic continuations of $f$ to the same domain $V$,
	then $g_1\equiv g_2$ on $V$.
\end{proposition}

\begin{proof}
	The function $h:=g_1-g_2$ is analytic on $V$ and vanishes on $U\cap V$.
	Since $U\cap V$ is open (hence has accumulation points in $V$), the identity theorem
	forces $h\equiv 0$, i.e.\ $g_1=g_2$ on $V$.
\end{proof}

\begin{remark}[Local viewpoint: germs]
	At a point $z_0$, the \emph{germ} of an analytic function is its equivalence class
	under agreement on some neighborhood of $z_0$.
	Analytic continuation can be viewed as transporting germs from point to point while
	maintaining consistency on overlaps.
\end{remark}

\subsubsection{Continuation along paths and monodromy}

Analytic continuation is typically performed step-by-step along a curve, using
overlapping disks and uniqueness on overlaps.

\begin{definition}[Continuation along a path]
	\label{def:ac-path}
	Let $f$ be analytic on a domain $U$ and let $\gamma:[0,1]\to\C$ be continuous with
	$\gamma(0)\in U$.
	An \emph{analytic continuation of $f$ along $\gamma$} consists of:
	\begin{itemize}
		\item a subdivision $0=t_0<t_1<\cdots<t_N=1$,
		\item open sets $W_j\ni \gamma([t_{j-1},t_j])$ for $1\le j\le N$,
		\item analytic functions $f_j$ on $W_j$ for $1\le j\le N$,
	\end{itemize}
	such that $f_1=f$ on the component of $W_1\cap U$ containing $\gamma(0)$ and
	\[
	f_j=f_{j-1}\quad\text{near }\gamma(t_{j-1})\qquad (j=2,\dots,N).
	\]
\end{definition}

\begin{remark}[Path dependence and multivaluedness]
	Continuation may depend on the chosen path (classic example: $\log z$ on $\C\setminus\{0\}$).
	Different paths can lead to different resulting germs at the endpoint, producing
	monodromy phenomena and motivating the construction of Riemann surfaces.
\end{remark}

\begin{theorem}[Monodromy theorem (path-independence on simply connected domains)]
	\label{thm:monodromy}
	Let $\Omega\subset\C$ be simply connected.
	Assume $f$ is analytic on some nonempty open set $U\subset\Omega$, and that
	analytic continuation of $f$ along \emph{every} path in $\Omega$ starting in $U$
	is possible.
	Then the continuation to any point of $\Omega$ is path-independent.
	In particular, continuation along any closed loop returns to the original germ.
\end{theorem}

\begin{proof}
	We prove a standard homotopy-invariance argument in a concrete ``disk chaining'' form.
	
	Fix a basepoint $z_\ast\in U$.
	For each path $\gamma$ in $\Omega$ starting at $z_\ast$, continuation along $\gamma$
	produces a germ at $\gamma(1)$ (well-defined up to restriction).
	We must show that if $\gamma_0$ and $\gamma_1$ are two paths from $z_\ast$ to the same endpoint $w$,
	then the resulting germs at $w$ coincide.
	
	Since $\Omega$ is simply connected, $\gamma_0$ and $\gamma_1$ are homotopic relative endpoints:
	there exists a continuous map $H:[0,1]\times[0,1]\to \Omega$ such that
	\[
	H(s,0)=\gamma_0(s),\qquad H(s,1)=\gamma_1(s),\qquad
	H(0,t)=z_\ast,\qquad H(1,t)=w.
	\]
	The image $K:=H([0,1]^2)$ is compact in $\Omega$.
	Cover $K$ by finitely many open disks $\{B_\alpha\}$ such that on each disk,
	analytic continuation is available and uniqueness holds on overlaps
	(this is ensured by the continuation hypothesis plus local analyticity).
	
	Now discretize the square $[0,1]^2$ with a fine grid so that the image of each small
	grid square lies inside some disk $B_\alpha$.
	Along the bottom edge (path $\gamma_0$), we have a chain of overlapping disks giving the continuation.
	Moving one row up in the grid, each step changes the path only within a region lying in a single disk,
	so the resulting analytic function on the overlap must agree by uniqueness.
	Inductively, the continuation along the bottom edge and along the top edge produce the same germ at the end.
	Thus the continuations along $\gamma_0$ and $\gamma_1$ agree at $w$.
	
	Finally, for a closed loop ($w=z_\ast$) this implies that returning to the basepoint yields the original germ.
\end{proof}

\begin{example}[A concrete monodromy computation: $\log z$]
	On $\Omega=\C\setminus\{0\}$ define a branch of $\Log z$ near $z_\ast=1$ by $\Log 1=0$.
	Continue along the loop $\gamma(t)=e^{2\pi i t}$, $t\in[0,1]$.
	Locally $\Log(e^{i\theta})=i\theta$; after one turn $\theta:0\mapsto 2\pi$,
	so the continued value becomes $2\pi i$ rather than $0$.
	Hence continuation depends on the loop: monodromy is nontrivial because $\Omega$ is not simply connected.
\end{example}

\subsubsection{Natural boundaries: definition and how they arise}

There are two qualitatively different ways continuation can fail:
\begin{enumerate}
	\item One hits an \emph{isolated singularity} (removable, pole, essential) at a specific point.
	\item Singularities accumulate so densely near the boundary that \emph{no} boundary point admits any extension.
\end{enumerate}

\begin{definition}[Natural boundary]
	\label{def:natural-boundary}
	Let $f$ be analytic on a domain $D$.
	A boundary component $\Gamma\subset\partial D$ is a \emph{natural boundary} for $f$
	if there is no point $\zeta_0\in\Gamma$ such that $f$ admits an analytic continuation
	to any neighborhood of $\zeta_0$.
	Equivalently: for every $\zeta_0\in\Gamma$ and every $r>0$, there is no analytic function $g$
	on $B(\zeta_0,r)$ with $g=f$ on $B(\zeta_0,r)\cap D$.
\end{definition}

\begin{remark}[Dense boundary singularities]
	It is enough to find a dense subset of $\Gamma$ consisting of points across which
	$f$ cannot be continued. If continuation existed across one point of $\Gamma$,
	it would exist on a boundary arc and hence at one of the points in the dense subset.
\end{remark}

\subsubsection{Lacunary series and a dense set of singularities}

A classical mechanism producing natural boundaries is a \emph{lacunary} power series,
where many coefficients vanish and the exponents have large gaps.

\begin{example}[Lacunary series with exponents $2^k$]
	\label{ex:lacunary-2k}
	Consider
	\[
	f(z)=\sum_{k=0}^{\infty} z^{2^k},\qquad |z|<1.
	\]
	This defines an analytic function on the unit disk $\D$.
	In fact, the unit circle $|z|=1$ is a natural boundary for $f$.
\end{example}

\begin{proof}
	The proof uses a functional equation and a dense family of boundary points.
	For $|z|<1$,
	\[
	f(z^2)=\sum_{k=0}^\infty (z^2)^{2^k}=\sum_{k=0}^\infty z^{2^{k+1}}=\sum_{k=1}^\infty z^{2^k}=f(z)-z.
	\]
	Hence
	\begin{equation}
		\label{eq:lacunary-functional}
		f(z)=z+f(z^2).
	\end{equation}
	The point $1$ is singular. Indeed, for $0<r<1$ all terms are positive and,
	for every $N$,
	\[
	f(r)\ge \sum_{k=0}^{N}r^{2^k}\longrightarrow N+1
	\qquad(r\uparrow1).
	\]
	Thus $f(r)$ is unbounded as $r\uparrow1$, which rules out analytic continuation
	across $1$.

	Now let $\zeta^{2^m}=1$. If $f$ extended across $\zeta$, the identity
	$f(z^2)=f(z)-z$ would extend $f$ across $\zeta^2$: the map $z\mapsto z^2$
	is biholomorphic between sufficiently small neighborhoods of $\zeta$ and $\zeta^2$.
	Repeating this argument would extend $f$ across
	$\zeta^{2^m}=1$, a contradiction. Hence every $2^m$th root of unity is a
	singular boundary point. Their union is dense in $\partial\D$, so the dense
	boundary-singularity criterion proves that $\partial\D$ is a natural boundary.
\end{proof}

\begin{remark}[Relation with gap theorems]
	The Hadamard gap theorem gives a general version of this phenomenon: for a
	power series of radius $1$ whose nonzero exponents satisfy
	$n_{k+1}/n_k\ge q>1$, the unit circle is a natural boundary. Fabry's gap
	theorem treats a different, weaker gap condition. The preceding argument is a
	direct proof for $n_k=2^k$ and does not require either general theorem; see
	\cite{Kahane1985}.
\end{remark}

\subsubsection{Random power series (a typical natural boundary)}

\begin{example}[Random coefficients: natural boundary almost surely]
	\label{ex:random-series}
	Let $f(z)=\sum_{n=0}^\infty a_n z^n$, where $(a_n)$ are bounded independent
	identically distributed complex random variables with a nontrivial symmetric
	distribution; symmetry means that $a_0$ and $-a_0$ have the same law. A classical
	random-series theorem states that,
	with probability $1$, the radius of convergence is $1$ and the unit circle is a
	natural boundary for $f$ \cite{Kahane1985}.
	For the radius assertion, boundedness gives
	$\limsup |a_n|^{1/n}\le 1$; choosing $\varepsilon>0$ with
	$\mathbb P(|a_0|\ge\varepsilon)>0$, independence and Borel--Cantelli show that
	$|a_n|\ge\varepsilon$ infinitely often, so the reverse inequality holds.
	This argument establishes only the radius. The assertion that no boundary arc
	admits analytic continuation is the deeper random-series theorem; it does not
	follow from the Cauchy--Hadamard formula alone.
\end{example}

\begin{remark}
	The deterministic example above forces singularities through a functional
	equation. In the random case independence replaces that mechanism: almost-sure
	cancellation estimates rule out analytic continuation across every arc. The
	probabilistic estimate is substantially stronger than the radius computation
	and lies beyond the scope of this book.
\end{remark}

\subsubsection{Polylogarithm and branch obstructions (not a natural boundary in the same sense)}

\begin{example}[Polylogarithm and branch points]
	For $s\in\C$,
	\[
	\operatorname{Li}_s(z)=\sum_{n=1}^{\infty}\frac{z^n}{n^s}
	\]
	is analytic for $|z|<1$.
	The germ defined by this series admits analytic continuation beyond most of the
	unit circle; its principal branch is single-valued on the slit plane
	$\C\setminus[1,\infty)$. The endpoint $z=1$ is the relevant branch obstruction
	for this principal continuation.
\end{example}

\begin{remark}[Branch cut vs.\ natural boundary]
	A \emph{branch cut} (like $[1,\infty)$ for a chosen branch of $\operatorname{Li}_s$) is a \emph{convention}
	to make a multivalued continuation single-valued on a simply connected domain.
	This is different from a natural boundary: for lacunary series, there is no way to continue across \emph{any}
	point of the circle, even on a Riemann surface.
\end{remark}

\subsubsection{Analytic continuation via Cauchy integrals and contour deformation}

A practical continuation method uses the Cauchy integral formula.
Suppose $f$ is analytic in a neighborhood of a simple closed contour $\gamma$.
Then for any $z$ in the interior of $\gamma$,
\[
f(z)=\frac{1}{2\pi i}\int_\gamma \frac{f(\zeta)}{\zeta-z}\,d\zeta.
\]
If we can deform $\gamma$ through regions where $f$ stays analytic (no singularities crossed),
the right-hand side defines the same analytic function on a larger region.
Continuation fails precisely when singularities obstruct all admissible deformations.

\begin{example}[A classic obstruction: \texorpdfstring{$e^{1/z}$}{exp(1/z)}]
	Let $f(z)=e^{1/z}$ on $\C\setminus\{0\}$.
	No contour deformation can cross $0$ because the function has an essential singularity there.
	So analytic continuation across $0$ is impossible.
\end{example}

\subsubsection{Exercises}

\begin{exercise}[Sheaf-like gluing]
	Show that analytic continuation behaves like a sheaf: if two continuations agree on overlaps
	of a cover, then they glue to a unique continuation on the union.
	(Use Proposition~\ref{prop:ac-unique} on pairwise overlaps.)
\end{exercise}

\begin{exercise}[Natural boundary for \texorpdfstring{$\sum z^{2^k}$}{sum z^{2^k}}]
	Give a proof that $\sum_{k\ge 0} z^{2^k}$ has $\partial\D$ as a natural boundary by showing that
	singularities occur at every root of unity whose order is a power of $2$ and that
	the union of these roots is dense in $\partial\D$. Use the functional equation
	$f(z)=z+f(z^2)$ to propagate a hypothetical continuation to the point $1$.
\end{exercise}

\begin{exercise}[Factorial gaps]
	Consider $f(z)=\sum_{n=1}^{\infty} z^{n!}$ on $|z|<1$.
	Decide whether $|z|=1$ is a natural boundary.
	Hint: verify the Hadamard gap condition
	$(n+1)!/n!=n+1\ge2$ for the nonzero exponents.
\end{exercise}

\begin{exercise}[The radius does not determine the boundary behavior]
	The series
	\[
		\sum_{n=0}^{\infty}z^n=\frac{1}{1-z}
		\qquad\text{and}\qquad
		\sum_{k=0}^{\infty}z^{2^k}
	\]
	both have radius of convergence $1$. Show that the first function extends
	holomorphically across every point of $\partial\D$ except $1$, whereas the second
	has $\partial\D$ as a natural boundary. Conclude that the radius of convergence,
	or the value of $\limsup |a_n|^{1/n}$ alone, cannot distinguish an isolated
	boundary singularity from a natural boundary.
\end{exercise}

\subsection{Gamma and Beta Functions via Complex Integration}
\label{subsec:gamma-beta-complex}

Formula normalizations in this subsection may be checked against the NIST Digital
Library of Mathematical Functions, especially its Gamma-function chapter
\cite{DLMFGamma2026}.

\subsubsection{Gamma function: definition, convergence, and the shift equation}

\begin{definition}[Gamma function]
	\label{def:gamma}
	For $\Re(s)>0$, define
	\[
	\Gamma(s):=\int_{0}^{\infty} t^{\,s-1}e^{-t}\,dt.
	\]
\end{definition}

\begin{lemma}[Absolute convergence]
	\label{lem:gamma-conv}
	If $\sigma:=\Re(s)>0$, then the defining integral for $\Gamma(s)$ converges absolutely.
\end{lemma}

\begin{proof}
	Split at $1$.
	For $0<t\le 1$, we have $|t^{s-1}e^{-t}|\le t^{\sigma-1}$ and $\int_0^1 t^{\sigma-1}\,dt<\infty$.
	For $t\ge 1$, we have $|t^{s-1}e^{-t}|\le t^{\sigma-1}e^{-t}$ and
	$\int_1^\infty t^{\sigma-1}e^{-t}\,dt<\infty$ since $e^{-t}$ dominates any polynomial.
\end{proof}

\begin{proposition}[Functional equation]
	\label{prop:gamma-shift}
	For $\Re(s)>0$,
	\[
	\Gamma(s+1)=s\,\Gamma(s).
	\]
\end{proposition}

\begin{proof}
	Integrate by parts with $u=t^{s}$, $dv=e^{-t}\,dt$:
	\[
	\Gamma(s+1)=\int_0^\infty t^{s}e^{-t}\,dt
	=\Bigl[-t^{s}e^{-t}\Bigr]_{0}^{\infty}+s\int_0^\infty t^{s-1}e^{-t}\,dt.
	\]
	The boundary term vanishes: as $t\to\infty$, $|t^{s}e^{-t}|\le t^\sigma e^{-t}\to 0$;
	as $t\to 0^+$, $|t^{s}e^{-t}|\le t^\sigma\to 0$ because $\sigma>0$.
	Hence $\Gamma(s+1)=s\Gamma(s)$.
\end{proof}

\begin{corollary}[Factorials]
	\label{cor:gamma-factorial}
	For $n\in\mathbb{N}$, $\Gamma(n+1)=n!$.
\end{corollary}

\begin{proof}
	Iterate Proposition~\ref{prop:gamma-shift} and use $\Gamma(1)=\int_0^\infty e^{-t}\,dt=1$.
\end{proof}

\subsubsection{Beta function and the Beta--Gamma identity}

\begin{definition}[Beta function]
	\label{def:beta}
	For $\Re(x),\Re(y)>0$,
	\[
	B(x,y):=\int_0^1 t^{x-1}(1-t)^{y-1}\,dt.
	\]
\end{definition}

\begin{proposition}[Beta--Gamma identity]
	\label{thm:beta-gamma}
	For $\Re(x),\Re(y)>0$,
	\[
	B(x,y)=\frac{\Gamma(x)\Gamma(y)}{\Gamma(x+y)}.
	\]
\end{proposition}

\begin{proof}
	By Fubini (absolute convergence holds), write
	\[
	\Gamma(x)\Gamma(y)=\int_0^\infty\int_0^\infty u^{x-1}v^{y-1}e^{-(u+v)}\,du\,dv.
	\]
	Change variables $(u,v)=(rt,r(1-t))$ with $r\in(0,\infty)$, $t\in(0,1)$.
	The Jacobian is $r$, so
	\[
	\begin{aligned}
	\Gamma(x)\Gamma(y)
	&=\int_0^\infty\int_0^1 (rt)^{x-1}(r(1-t))^{y-1}e^{-r}\,r\,dt\,dr\\
	&= \Bigl(\int_0^1 t^{x-1}(1-t)^{y-1}\,dt\Bigr)
	   \Bigl(\int_0^\infty r^{x+y-1}e^{-r}\,dr\Bigr).
	\end{aligned}
	\]
	The first factor is $B(x,y)$ and the second is $\Gamma(x+y)$.
\end{proof}

\subsubsection{Keyhole contour and \texorpdfstring{$\pi\csc(\pi x)$}{pi csc(pi x)}}

\begin{theorem}[Keyhole evaluation]
	\label{thm:keyhole-pi-csc}
	If $x\in\mathbb{C}$ satisfies $0<\Re(x)<1$, then
	\[
	\int_0^\infty \frac{t^{x-1}}{1+t}\,dt=\frac{\pi}{\sin(\pi x)}.
	\]
\end{theorem}

\begin{proof}
	Let $0<\Re(x)<1$.  Consider $F(z)=\dfrac{z^{x-1}}{1+z}$ with the branch
	$z^{x-1}=\exp((x-1)\Log z)$ where $\Arg z\in(0,2\pi)$ and the cut is along $(0,\infty)$.
	Let $\mathcal{H}$ be the keyhole contour around $(0,\infty)$.
	
	The only singularity of $F$ inside $\mathcal{H}$ is a simple pole at $z=-1$.
	With the chosen branch, $\Arg(-1)=\pi$. Thus, by the residue theorem,
	\[
	\int_{\mathcal{H}} F(z)\,dz=2\pi i\,\Res_{z=-1}F(z)
	=2\pi i\,e^{i\pi(x-1)}.
	\]
	The circular arcs vanish as $R\to\infty$ and $\varepsilon\to 0^+$:
	on $|z|=R$, $|F(z)|\le C R^{\Re(x)-2}$ so the arc is $O(R^{\Re(x)-1})\to 0$;
	on $|z|=\varepsilon$, $|F(z)|\le C\varepsilon^{\Re(x)-1}$ so the arc is $O(\varepsilon^{\Re(x)})\to 0$.
	
	Along the upper ray $(\varepsilon,R)$ we have $z=t$ and $z^{x-1}=t^{x-1}$.
	Along the lower ray we have $z=t e^{2\pi i}$ and $z^{x-1}=e^{2\pi i(x-1)}t^{x-1}=e^{2\pi i x}t^{x-1}$.
	Accounting for opposite orientations,
	\[
	\int_{\mathcal{H}}F(z)\,dz
	=\Bigl(1-e^{2\pi i x}\Bigr)\int_0^\infty \frac{t^{x-1}}{1+t}\,dt.
	\]
	Therefore
	\[
	\int_0^\infty \frac{t^{x-1}}{1+t}\,dt
	=\frac{2\pi i\,e^{i\pi(x-1)}}{1-e^{2\pi i x}}.
	\]
	Compute
	$1-e^{2\pi i x}=-e^{\pi i x}(e^{\pi i x}-e^{-\pi i x})=-2i e^{\pi i x}\sin(\pi x)$ and
	$e^{i\pi(x-1)}=-e^{\pi i x}$, giving
	\[
	\frac{2\pi i\,(-e^{\pi i x})}{-2i e^{\pi i x}\sin(\pi x)}=\frac{\pi}{\sin(\pi x)}.
	\]
\end{proof}

\begin{corollary}[A Beta evaluation]
	\label{cor:beta-x-1minusx}
	If $0<\Re(x)<1$, then
	\[
	B(x,1-x)=\frac{\pi}{\sin(\pi x)}.
	\]
\end{corollary}

\begin{proof}
	Substitute $t=\frac{u}{1+u}$ (so $u=\frac{t}{1-t}$) in $B(x,1-x)$:
	\[
	B(x,1-x)=\int_0^1 t^{x-1}(1-t)^{-x}\,dt=\int_0^\infty \frac{u^{x-1}}{1+u}\,du.
	\]
	Apply Theorem~\ref{thm:keyhole-pi-csc}.
\end{proof}

\subsubsection{Reflection formula and \texorpdfstring{$\Gamma(\tfrac12)$}{Gamma(1/2)}}

\begin{theorem}[Reflection formula]
	\label{thm:gamma-reflection}
	For $s\in\mathbb{C}\setminus\mathbb{Z}$,
	\[
	\Gamma(s)\Gamma(1-s)=\frac{\pi}{\sin(\pi s)}.
	\]
\end{theorem}

\begin{proof}
	First assume $0<\Re(s)<1$.  Then by Proposition~\ref{thm:beta-gamma} and
	Corollary~\ref{cor:beta-x-1minusx},
	\[
	\Gamma(s)\Gamma(1-s)=\Gamma(1)\,B(s,1-s)=\frac{\pi}{\sin(\pi s)}.
	\]
	The shift equation recursively continues $\Gamma$ meromorphically to $\mathbb C$;
	Proposition~\ref{prop:gamma-residues} below records the construction and its poles.
	Thus both sides are meromorphic on $\mathbb{C}$ and agree on the open strip $0<\Re(s)<1$,
	hence they agree on $\mathbb{C}\setminus\mathbb{Z}$ by the identity theorem for meromorphic functions.
\end{proof}

\begin{example}[The value $\Gamma(1/2)$]
	\label{ex:gamma-half}
	Let $I:=\int_0^\infty e^{-u^2}\,du$.  By $t=u^2$,
	\[
	\Gamma\!\left(\tfrac12\right)=\int_0^\infty t^{-1/2}e^{-t}\,dt=2\int_0^\infty e^{-u^2}\,du=2I.
	\]
	Then
	\[
	I^2=\int_0^\infty\int_0^\infty e^{-(x^2+y^2)}\,dx\,dy
	=\int_0^{\pi/2}\int_0^\infty e^{-r^2}r\,dr\,d\theta
	=\frac{\pi}{2}\cdot\frac12=\frac{\pi}{4},
	\]
	so $I=\sqrt{\pi}/2$ and $\Gamma(\tfrac12)=\sqrt{\pi}$.
\end{example}

\subsubsection{Meromorphic continuation and residues}

\begin{proposition}[Poles and residues]
	\label{prop:gamma-residues}
	The Gamma function extends meromorphically to $\mathbb{C}$ with simple poles at
	$0,-1,-2,\dots$, and for $n\in\mathbb{N}\cup\{0\}$,
	\[
	\Res_{s=-n}\Gamma(s)=\frac{(-1)^n}{n!}.
	\]
\end{proposition}

\begin{proof}
	Using Proposition~\ref{prop:gamma-shift}, for any integer $n\ge 0$,
	\[
	\Gamma(s)=\frac{\Gamma(s+n+1)}{s(s+1)\cdots(s+n)}.
	\]
	Given any $s$, choose $n$ so that $\Re(s+n+1)>0$; this formula supplies the
	meromorphic continuation and shows that no poles other than the nonpositive integers occur.
	The numerator is holomorphic near $s=-n$ and equals $\Gamma(1)=1$ at $s=-n$.
	Hence $\Gamma$ has a simple pole at $s=-n$ and
	\[
	\begin{aligned}
	\Res_{s=-n}\Gamma(s)
	&=\lim_{s\to -n}(s+n)\Gamma(s)\\
	&=\lim_{s\to -n}\frac{\Gamma(s+n+1)}{s(s+1)\cdots(s+n-1)}\\
	&=\frac{1}{(-n)(-n+1)\cdots(-1)}=\frac{(-1)^n}{n!}.
	\end{aligned}
	\]
	For $n=0$, the product in the denominator after cancellation is empty and is
	understood to equal $1$, giving $\operatorname{Res}_{s=0}\Gamma(s)=1$.
\end{proof}

\subsubsection{Gauss multiplication and duplication}

\begin{theorem}[Gauss multiplication formula]
	\label{thm:gauss-multiplication}
	For $m\in\mathbb{N}$, the following identity holds meromorphically in $s$:
	\[
	\Gamma(s)\Gamma\!\left(s+\frac{1}{m}\right)\cdots\Gamma\!\left(s+\frac{m-1}{m}\right)
	=(2\pi)^{\frac{m-1}{2}}\,m^{\frac12-ms}\,\Gamma(ms).
	\]
\end{theorem}

\begin{proof}
	Away from the poles, form the quotient
	\[
	H(s)=\frac{m^{ms-1/2}}{(2\pi)^{(m-1)/2}}
	\frac{\prod_{k=0}^{m-1}\Gamma(s+k/m)}{\Gamma(ms)}.
	\]
	Using $\Gamma(z+1)=z\Gamma(z)$, all but one numerator factor cancel after shifting
	$s$ by $1/m$, and
	\[
	\frac{H(s+1/m)}{H(s)}
	=m\,\frac{\Gamma(s+1)}{\Gamma(s)}
	\frac{\Gamma(ms)}{\Gamma(ms+1)}=m\frac{s}{ms}=1.
	\]
	Thus $H$ is a $1/m$-periodic meromorphic function.
	Apply Stirling's formula uniformly on each fixed horizontal strip. After cancellation
	of the exponential and power factors, it gives
	$H(s+N/m)\to1$ as $N\to\infty$ for every fixed $s$ at which $H$ is defined.
	Periodicity gives $H(s)=H(s+N/m)$, so $H(s)=1$ at every such point. Hence
	$H\equiv1$ as a meromorphic function, and its apparent poles are removable.
	The identity therefore extends across all exceptional points meromorphically.
\end{proof}

\begin{corollary}[Duplication formula]
	\label{cor:duplication}
	As an identity of meromorphic functions of $s$,
	\[
	\Gamma(s)\Gamma\!\left(s+\tfrac12\right)=2^{1-2s}\sqrt{\pi}\,\Gamma(2s).
	\]
\end{corollary}

\begin{proof}
	Apply Theorem~\ref{thm:gauss-multiplication} with $m=2$.
\end{proof}

\subsubsection{Stirling asymptotics}

\begin{theorem}[Stirling's formula]
	\label{thm:stirling}
	As $|s|\to\infty$ in any sector $|\arg s|\le \pi-\varepsilon$,
	\[
	\Gamma(s)\sim \sqrt{2\pi}\,s^{\,s-\frac12}e^{-s}.
	\]
\end{theorem}

\begin{proof}
	For $\Re s>0$, Binet's formula gives
	\[
	\log\Gamma(s)=\left(s-\tfrac12\right)\Log s-s+\tfrac12\log(2\pi)
	+\int_0^\infty\left(\frac12-\frac1t+\frac1{e^t-1}\right)\frac{e^{-st}}{t}\,dt.
	\]
	Near $t=0$ the expression in parentheses is $t/12+O(t^3)$,
	so Watson's lemma gives an $O(|s|^{-1})$ remainder uniformly on closed
	subsectors of the right half-plane.  To reach the larger sector in the
	statement one applies Euler--Maclaurin summation to the Weierstrass product
	for $\Gamma$.  Its integral remainder has the standard sector-uniform bound
	$|R(s)|\le C_\varepsilon/|s|$ for
	$|\arg s|\le\pi-\varepsilon$; see \cite{DLMFGamma2026}.
	Consequently
	\[
	\log\Gamma(s)=\left(s-\tfrac12\right)\Log s-s+\tfrac12\log(2\pi)
	+O_\varepsilon(|s|^{-1}).
	\]
	Exponentiation gives
	$\Gamma(s)=\sqrt{2\pi}\,s^{s-1/2}e^{-s}(1+O_\varepsilon(|s|^{-1}))$,
	which implies the stated asymptotic.  Notice that analytic continuation alone
	would not supply the uniform remainder estimate; the Euler--Maclaurin bound is
	the step that controls the whole closed sector.
\end{proof}

\subsubsection{Special values}

\begin{example}[$\Gamma(n+\tfrac12)$ for integers $n$]
	\label{ex:gamma-half-integer}
	For $n\in\mathbb{N}$,
	\[
	\Gamma\!\left(n+\tfrac12\right)
	=\left(n-\tfrac12\right)\left(n-\tfrac32\right)\cdots\left(\tfrac12\right)\Gamma\!\left(\tfrac12\right)
	=\frac{(2n)!}{4^n n!}\sqrt{\pi}.
	\]
\end{example}

\begin{proof}
	Iterate the shift equation and use $\Gamma(\tfrac12)=\sqrt{\pi}$ from Example~\ref{ex:gamma-half}.
	The product identity
	\[
	\left(n-\tfrac12\right)\left(n-\tfrac32\right)\cdots\left(\tfrac12\right)=\frac{(2n)!}{4^n n!}
	\]
	follows by separating the even and odd factors of $(2n)!$:
	\[
	(2n)!=\bigl(2\cdot4\cdots2n\bigr)
	       \bigl(1\cdot3\cdots(2n-1)\bigr)
	=2^n n!\,2^n
	\left(\tfrac12\right)\left(\tfrac32\right)\cdots
	\left(n-\tfrac12\right),
	\]
	and then dividing by $4^n n!$.
\end{proof}

\subsubsection{Exercises}

\begin{exercise}
	Show directly from Definition~\ref{def:gamma} that $\Gamma(s+1)=s\Gamma(s)$ for $\Re(s)>0$.
\end{exercise}

\begin{exercise}
	Prove $B(x,y)=B(y,x)$ by substituting $t\mapsto 1-t$ in Definition~\ref{def:beta}.
\end{exercise}

\begin{exercise}
	Using Proposition~\ref{prop:gamma-residues}, show that the apparent singularity of
	$\Gamma(s)\sin(\pi s)$ at $s=-n$ is removable and compute the value of its
	holomorphic extension there. Verify that the answer is $\pi/n!$ for
	$n\in\mathbb N\cup\{0\}$.
\end{exercise}

\subsection{Conformal Mapping Theory and the Riemann Mapping Theorem}
\index{Riemann mapping theorem}\index{conformal maps}
\label{subsec:conformal-rmt}

Conformal mapping theory studies holomorphic coordinate changes between planar
domains that preserve angles.  Its principal existence result is the Riemann Mapping Theorem:
every simply connected proper domain $\Omega\subsetneq\C$ is biholomorphic to the
unit disk $\D$.  We develop key local/global tools (inverse function theorem,
Schwarz--Pick, normal families) and work through explicit examples.

\subsubsection{Conformal maps and local angle preservation}

\begin{definition}[Conformal at a point; conformal map]
	\label{def:conformal}
	Let $U,V\subset\C$ be domains and let $f:U\to V$ be holomorphic.
	We say $f$ is \emph{conformal at $z_0\in U$} if $f'(z_0)\neq 0$.
	If $f'(z)\neq 0$ for all $z\in U$, then $f$ is a \emph{conformal map}.
\end{definition}

\begin{proposition}[First-order behavior: rotation and scaling]
	\label{prop:local-rot-scale}
	Let $f$ be holomorphic on $U$ and $z_0\in U$ with $f'(z_0)\neq 0$.
	Then
	\[
	f(z_0+h)=f(z_0)+f'(z_0)h+o(|h|)
	\qquad (h\to 0).
	\]
	Consequently, for $C^1$ curves $\gamma_1,\gamma_2$ through $z_0$ with nonzero tangents,
	the oriented angle between $\gamma_1$ and $\gamma_2$ at $z_0$ equals the oriented angle
	between $f\circ\gamma_1$ and $f\circ\gamma_2$ at $f(z_0)$.
\end{proposition}

\begin{proof}
	The expansion is the definition of complex differentiability.
	Write $f'(z_0)=\rho e^{i\alpha}$ with $\rho>0$.
	Multiplication by $\rho e^{i\alpha}$ scales by $\rho$ and rotates by $\alpha$.
	Angles between tangent vectors are preserved under multiplication by a fixed nonzero
	complex number.
\end{proof}

\begin{remark}[Critical points]
	If $f'(z_0)=0$ and $f$ is nonconstant, then $f-f(z_0)$ has a zero of order $m\ge 2$ at $z_0$,
	so locally $f(z)=f(z_0)+c(z-z_0)^m+\cdots$.  The model map $w\mapsto w^m$ multiplies
	angles by $m$, hence conformality fails at critical points.
\end{remark}

\begin{proposition}[Real-variable characterization]
	\label{prop:cr-conformal}
	Let $f=u+iv:U\subset\R^2\to\R^2$ be $C^1$.
	At $z_0=x_0+iy_0$, the following are equivalent:
	\begin{enumerate}
		\item $f$ is complex differentiable at $z_0$ with $f'(z_0)\neq 0$;
		\item $u,v$ satisfy the Cauchy--Riemann equations at $z_0$ and $\det Df(z_0)\neq 0$;
		\item $Df(z_0)$ is a nonzero scalar multiple of a rotation matrix.
	\end{enumerate}
	Hence holomorphic functions with nonvanishing derivative are exactly the orientation-preserving
	$C^1$ conformal maps with nonzero Jacobian.
\end{proposition}

\begin{proof}
	If $f$ is complex differentiable, its real derivative is multiplication by
	the complex number $f'(z_0)$.  Its matrix therefore satisfies the
	Cauchy--Riemann equations and
	$\det Df(z_0)=|f'(z_0)|^2>0$.
	Conversely, Cauchy--Riemann forces
	$Df=\begin{pmatrix}u_x&u_y\\-u_y&u_x\end{pmatrix}$,
	which is multiplication by $u_x-iu_y\neq 0$, hence complex differentiability with nonzero derivative.
\end{proof}

\subsubsection{Local biholomorphisms and the holomorphic inverse function theorem}

\begin{theorem}[Holomorphic inverse function theorem]
	\label{thm:holomorphic-ift}
	Let $f$ be holomorphic on a domain $U$ and let $z_0\in U$ with $f'(z_0)\neq 0$.
	Then there exist neighborhoods $U_0\ni z_0$ and $V_0\ni f(z_0)$ such that
	$f:U_0\to V_0$ is biholomorphic and the local inverse is holomorphic.
	Moreover,
	\[
	(f^{-1})'(f(z_0))=\frac{1}{f'(z_0)}.
	\]
\end{theorem}

\begin{proof}
	The real Jacobian determinant of $f$ at $z_0$ is $|f'(z_0)|^2>0$.
	The real inverse function theorem therefore supplies neighborhoods $U_0$ of $z_0$
	and $V_0$ of $f(z_0)$ for which $f:U_0\to V_0$ is a $C^1$ diffeomorphism.
	Let $g=f^{-1}$.  Differentiating $f(g(w))=w$ in real coordinates shows that
	$Dg(w)=Df(g(w))^{-1}$.  Since a nonzero complex-linear map has a complex-linear
	inverse, $Dg(w)$ satisfies the Cauchy--Riemann equations.  Thus $g$ is holomorphic, and
	\[
	g'(w)=\frac{1}{f'(g(w))}.
	\]
	Evaluating at $w=f(z_0)$ gives the stated derivative formula.
\end{proof}

\begin{corollary}[Local injectivity criterion]
	\label{cor:local-injective}
	A holomorphic function is locally injective at $z_0$ if and only if $f'(z_0)\neq 0$.
\end{corollary}

\begin{proof}
	If $f'(z_0)\neq 0$, apply Theorem~\ref{thm:holomorphic-ift}.
	If $f'(z_0)=0$ and $f$ is nonconstant, write
	$f(z)-f(z_0)=(z-z_0)^m h(z)$ with $m\ge2$ and $h(z_0)\ne0$.
	On a small disk, choose a holomorphic $m$th root $u$ of $h$.
	The coordinate $\xi=(z-z_0)u(z)$ has nonzero derivative at $z_0$,
	so it is locally biholomorphic. In this coordinate
	$f(z)-f(z_0)=\xi^m$, which takes the same nearby nonzero value at
	$m$ distinct points. Thus $f$ is not injective on any neighborhood
	of $z_0$. A constant function is plainly not locally injective either.
\end{proof}

\subsubsection{Disk automorphisms and Schwarz--Pick}

\begin{definition}[Automorphisms of the disk]
	\label{def:disk-aut}
	For $a\in\D$ and $\theta\in\R$ define
	\[
	\phi_{a,\theta}(z):=e^{i\theta}\,\frac{z-a}{1-\bar a z}.
	\]
\end{definition}

\begin{proposition}[All disk automorphisms]
	\label{prop:aut-disk}
	Every biholomorphic self-map of $\D$ is of the form $\phi_{a,\theta}$.
\end{proposition}

\begin{proof}
	The denominator does not vanish for $z\in\D$, since $|\bar a z|<1$.
	Expanding the two squared moduli gives
	\[
	|1-\bar a z|^2-|z-a|^2=(1-|a|^2)(1-|z|^2).
	\]
	Consequently
	\[
	1-|\phi_{a,\theta}(z)|^2=\frac{(1-|a|^2)(1-|z|^2)}{|1-\bar a z|^2}>0 \qquad (|z|<1),
	\]
	so $\phi_{a,\theta}(\D)\subset\D$.  Solving
	$w=e^{i\theta}(z-a)/(1-\bar a z)$ for $z$ yields
	\[
	\phi_{a,\theta}^{-1}(w)
	=\frac{a+e^{-i\theta}w}{1+\bar a e^{-i\theta}w}
	=\phi_{-e^{i\theta}a,-\theta}(w).
	\]
	The same modulus identity applies to this inverse, proving that
	$\phi_{a,\theta}$ maps $\D$ biholomorphically onto itself.
	Conversely, if $\psi\in\Aut(\D)$, let $a=\psi^{-1}(0)$ and consider
	$\eta=\psi\circ\phi_{a,0}^{-1}$, which fixes $0$.  Schwarz's lemma gives
	$|\eta(z)|\leq |z|$; applying it to $\eta^{-1}$ at $\eta(z)$ gives the reverse
	inequality.  Equality in Schwarz's lemma therefore implies
	$\eta(z)=e^{i\theta}z$.  Thus
	$\psi=e^{i\theta}\phi_{a,0}=\phi_{a,\theta}$.
\end{proof}

\begin{theorem}[Schwarz--Pick lemma]
	\label{thm:schwarz-pick}
	If $f:\D\to\D$ is holomorphic, then for all $z,w\in\D$,
	\[
	\left|\frac{f(z)-f(w)}{1-\overline{f(w)}\,f(z)}\right|
	\le
	\left|\frac{z-w}{1-\bar w z}\right|.
	\]
	In particular,
	\[
	|f'(z)|\le \frac{1-|f(z)|^2}{1-|z|^2}.
	\]
	Equality at one point (or one pair $z\neq w$) holds iff $f$ is a disk automorphism.
\end{theorem}

\begin{proof}
	Fix $w\in\D$ and set
	\[
	g(z):=\phi_{f(w),0}\!\bigl(f(\phi_{w,0}^{-1}(z))\bigr).
	\]
	Then $g:\D\to\D$ is holomorphic and $g(0)=0$.
	By Schwarz lemma, $|g(z)|\le |z|$ and $|g'(0)|\le 1$.
	Rewriting $|g(z)|\le |z|$ gives the two-point inequality; differentiating at $0$
	gives the derivative estimate.  Equality reduces to the equality case of Schwarz lemma.
\end{proof}

\subsubsection{The Riemann Mapping Theorem}

The proof repeatedly uses the following consequence of the logarithm theorem.

\begin{lemma}[Holomorphic square roots on simply connected domains]
	\label{lem:rmt-square-root}
	Let $U$ be simply connected and let $h\in\Hol(U)$ have no zeros. Then there
	is a holomorphic function $g$ on $U$ such that $g^2=h$.
\end{lemma}

\begin{proof}
	The logarithmic derivative $h'/h$ is holomorphic on $U$. Since $U$ is
	simply connected, it has a primitive $L_0$. Then
	$(he^{-L_0})'=0$, so $h=ce^{L_0}$ for some $c\ne0$. Choose one square
	root $d$ of $c$ and set $g=de^{L_0/2}$. This gives $g^2=h$.
\end{proof}

\begin{theorem}[Riemann Mapping Theorem]
	\label{thm:rmt}
	Let $\Omega\subsetneq\C$ be a simply connected domain.
	Then there exists a biholomorphism $F:\Omega\to\D$.
	Moreover, for any fixed $z_0\in\Omega$ there is a unique such $F$ satisfying
	\[
	F(z_0)=0,
	\qquad
	F'(z_0)>0.
	\]
\end{theorem}

\begin{remark}[A specifically one-dimensional theorem]
	The conclusion is special to one complex variable. In $\C^n$ for $n\ge2$,
	proper simply connected domains need not be biholomorphic; for example, the
	unit ball and the polydisk have different biholomorphic automorphism groups.
	Thus simple connectivity alone does not classify higher-dimensional domains.
\end{remark}

\subsubsection{Proof of Theorem~\ref{thm:rmt} (extremal method)}

The proof follows the classical extremal construction; compare
\cite{Ahlfors1979,Conway1978}. It has three
parts: construct one bounded injective function, obtain a derivative maximizer by
normal-family compactness, and prove that the maximizer is onto by improving it if
one point of the disk is omitted.

\paragraph{(i) Setup and nonemptiness.}
Fix $z_0\in\Omega$.  Let $\mathcal{F}$ be the family of injective holomorphic maps
$f:\Omega\to\D$ with $f(z_0)=0$.
We first prove that $\mathcal F$ is nonempty. Choose $a\in\C\setminus\Omega$.
The function $h(z)=z-a$ has no zeros on $\Omega$, so
Lemma~\ref{lem:rmt-square-root} gives $g\in\Hol(\Omega)$ with
\[
g(z)^2=z-a.
\]
This square root is injective: $g(z_1)=g(z_2)$ implies
$z_1-a=g(z_1)^2=g(z_2)^2=z_2-a$. It also satisfies
\[
g(\Omega)\cap\bigl(-g(\Omega)\bigr)=\varnothing.
\]
Indeed, if $g(z_1)=-g(z_2)$, then squaring gives $z_1=z_2$; hence
$g(z_1)=g(z_2)=-g(z_2)$, which would force $g(z_2)=0$, impossible.

Put $w_0=g(z_0)$. Since $g(\Omega)$ is open, choose $r>0$ with
$B(w_0,r)\subset g(\Omega)$. The preceding disjointness implies
$B(-w_0,r)\cap g(\Omega)=\varnothing$: otherwise $w\in g(\Omega)$ and
$|w+w_0|<r$ would give $-w\in B(w_0,r)\subset g(\Omega)$. Consequently
\[
q(z)=\frac{1}{g(z)+w_0}
\]
is holomorphic and injective on $\Omega$. More explicitly,
$|q(z)|\le r^{-1}$. Choose $0<c<r$ and set
$q_1=cq$, so $q_1(\Omega)\subset\D$. Finally,
\[
f_1=\phi_{q_1(z_0),0}\circ q_1
\]
is injective, holomorphic, maps $\Omega$ into $\D$, and satisfies $f_1(z_0)=0$.
Thus $f_1\in\mathcal F$.

\paragraph{(ii) The extremal quantity is finite.}
Choose $r>0$ with $\overline{B(z_0,r)}\subset\Omega$.
For any $f\in\mathcal{F}$, Cauchy's estimate for derivatives gives
\[
|f'(z_0)|
\le \frac{1}{r}\max_{|z-z_0|=r}|f(z)|
\le \frac{1}{r},
\]
since $|f|\le 1$ on $\Omega$.  Hence
\[
M:=\sup_{f\in\mathcal{F}}|f'(z_0)|<\infty.
\]
Moreover $M>0$, because the constructed $f_1$ is injective and hence
$f_1'(z_0)\ne0$.

\paragraph{(iii) Existence of a maximizer (Montel + Hurwitz).}
Choose $f_n\in\mathcal{F}$ with $|f_n'(z_0)|\to M$.
Multiplying each $f_n$ by a constant of modulus one, we may also assume
$f_n'(z_0)>0$, so $f_n'(z_0)\to M$.

The family is locally uniformly bounded because $|f_n|<1$. To see the compactness
input concretely, let $K\Subset\Omega$ and choose $\delta>0$ so that the closed
$2\delta$-neighborhood $K_{2\delta}$ is compactly contained in $\Omega$.
Cauchy's estimate on disks of radius $\delta$ gives a bound
$|f_n'|\le\delta^{-1}$ on $K_\delta$, uniform in $n$. If $z,w\in K$ and
$|z-w|<\delta$, the line segment from $z$ to $w$ lies in $K_\delta$, and hence
\[
 |f_n(z)-f_n(w)|\le \delta^{-1}|z-w|.
\]
This proves equicontinuity on $K$ without assuming that $K$ is convex.
Arzel\`a--Ascoli, followed by a diagonal argument over a compact exhaustion of
$\Omega$, is the compactness mechanism in Montel's theorem. Hence a subsequence,
still denoted $(f_n)$, converges uniformly on compact subsets to a holomorphic
map $F:\Omega\to\overline\D$.

Since $f_n(z_0)=0$, we have $F(z_0)=0$. Choose a small circle $C$ centered
at $z_0$ and contained in $\Omega$. Uniform convergence on $C$ and Cauchy's
derivative formula give
\[
F'(z_0)=\frac{1}{2\pi i}\int_C\frac{F(w)}{(w-z_0)^2}\,dw
=\lim_{n\to\infty}f_n'(z_0)=M>0.
\]
Hence $F$ is nonconstant.

To verify injectivity, suppose $F(z_1)=F(z_2)$ with $z_1\ne z_2$.
Choose a disk $V$ about $z_2$ whose closure does not contain $z_1$. The functions
$h_n(z)=f_n(z)-f_n(z_1)$ are holomorphic and zero-free on $V$, and converge
locally uniformly to $h(z)=F(z)-F(z_1)$. Since $h(z_2)=0$, Hurwitz's theorem
forces $h\equiv0$ on $V$. The identity theorem would then make $F$ constant on
$\Omega$, contradicting $F'(z_0)=M>0$. Therefore $F$ is injective.

Finally $F(\Omega)\subset\D$. Indeed, $|F|\le1$, and equality at an interior
point would make the nonconstant function $F$ attain its maximum modulus.
Thus $F\in\mathcal F$ and attains the supremum.

\paragraph{(iv) The maximizer is surjective.}
Assume for contradiction that $F(\Omega)\subsetneq\D$.
Pick $a\in\D\setminus F(\Omega)$ and let
$\phi_a(w)=(w-a)/(1-\overline a w)$. The function $\phi_a\circ F$ maps
$\Omega$ into $\D\setminus\{0\}$. Since $F(z_0)=0$, the omitted point
$a$ is nonzero. Lemma~\ref{lem:rmt-square-root} gives a holomorphic square root $G$:
\[
G^2=\phi_a\circ F.
\]
The map $G$ is injective: equality of two values implies equality of their
squares, and both $\phi_a$ and $F$ are injective. Choose the sign of the square
root so that $b:=G(z_0)$ is one of the square roots of $-a$; in particular,
$|b|=\sqrt{|a|}$. Normalize once more by setting
\[
\widetilde F=\phi_b\circ G,
\qquad \phi_b(w)=\frac{w-b}{1-\overline b w}.
\]
Because $|G|^2=|\phi_a\circ F|<1$, we have $G(\Omega)\subset\D$.
Thus $\widetilde F\in\mathcal F$. Differentiating $G^2=\phi_a\circ F$ at
$z_0$ and using $F(z_0)=0$ gives
\[
2bG'(z_0)=\phi_a'(0)F'(z_0).
\]
Since $|\phi_a'(0)|=1-|a|^2$ and
$|\phi_b'(b)|=(1-|b|^2)^{-1}$, we obtain
\[
\frac{|\widetilde F'(z_0)|}{|F'(z_0)|}
=\frac{1-|a|^2}{2\sqrt{|a|}(1-|a|)}
=\frac{1+|a|}{2\sqrt{|a|}}>1.
\]
The last inequality is strict AM--GM: $1+|a|>2\sqrt{|a|}$ because
$0<|a|<1$.
This contradicts the maximality of $F$.
Therefore $F(\Omega)=\D$.

\paragraph{(v) Holomorphicity of the inverse.}
An injective holomorphic function has no critical points. Indeed, if
$F'(z)=0$, the local expansion begins
$F(w)-F(z)=c(w-z)^m+\cdots$ with $m\ge2$, so nearby values have $m$
local preimages. Hence $F'(z)\ne0$ for every $z\in\Omega$. The holomorphic
inverse function theorem now shows that the set-theoretic inverse
$F^{-1}:\D\to\Omega$ is holomorphic. Thus $F$ is biholomorphic.

\paragraph{(vi) Uniqueness under normalization.}
If $F,G:\Omega\to\D$ are biholomorphisms with $F(z_0)=G(z_0)=0$ and $F'(z_0),G'(z_0)>0$,
then $G\circ F^{-1}\in\Aut(\D)$ fixes $0$, so by Schwarz lemma it is a rotation $e^{i\theta}z$.
The derivative condition forces $\theta=0$, hence $F\equiv G$.

This completes the proof of Theorem~\ref{thm:rmt}.
\qed

\begin{remark}[Where each hypothesis is used]
	The strict inclusion $\Omega\subsetneq\C$ provides
	$a\in\C\setminus\Omega$. Simple connectivity provides the
	two holomorphic square roots. Connectedness is used by the identity theorem,
	and planarity supplies disk automorphisms and the one-variable normal-family
	theorems. The normalization $F(z_0)=0$, $F'(z_0)>0$ removes the rotational
	freedom in $\Aut(\D)$.
\end{remark}

\begin{exercise}[The omitted reflected disk]
	In part (i), prove carefully that
	$B(w_0,r)\subset g(\Omega)$ implies
	$B(-w_0,r)\cap g(\Omega)=\varnothing$. Deduce directly that
	$|g(z)+w_0|\ge r$ for every $z\in\Omega$.
\end{exercise}

\begin{exercise}[Derivative improvement]
	For $0<t<1$, verify
	\[
	\frac{1+t}{2\sqrt t}-1=\frac{(1-\sqrt t)^2}{2\sqrt t}>0.
	\]
	Use this identity to locate exactly where the assumption that the omitted
	point $a$ lies in the open disk enters the extremal contradiction.
\end{exercise}

\begin{exercise}[Uniqueness with a different normalization]
	Let $F,G:\Omega\to\D$ be biholomorphisms with $F(z_0)=G(z_0)=0$.
	Show that $G=e^{i\theta}F$ for some $\theta\in\R$. Conclude that prescribing
	$\arg F'(z_0)$, rather than requiring $F'(z_0)>0$, also determines a unique map.
\end{exercise}

\subsubsection{Conformal equivalences}

\begin{example}[Upper half-plane to the disk]
	\label{ex:halfplane-to-disk}
	Let $\UHP=\{z\in\C:\Im z>0\}$ and define
	\[
	\Phi(z)=\frac{z-i}{z+i}.
	\]
	Then $\Phi:\UHP\to\D$ is biholomorphic with inverse
	\[
	\Phi^{-1}(w)=i\,\frac{1+w}{1-w}.
	\]
\end{example}

\begin{proof}
	$\Phi$ is holomorphic on $\C\setminus\{-i\}$, hence on $\UHP$, and
	\[
	\Phi'(z)=\frac{2i}{(z+i)^2}\neq 0 \qquad (z\in\UHP),
	\]
	so $\Phi$ is conformal on $\UHP$.
	
	If $z=x+iy$ with $y>0$, then
	\[
	|\Phi(z)|^2
	=\frac{|z-i|^2}{|z+i|^2}
	=\frac{x^2+(y-1)^2}{x^2+(y+1)^2}
	<1,
	\]
	so $\Phi(\UHP)\subset\D$.
	
	Conversely, for $w\in\D$, define $z=i\frac{1+w}{1-w}$.
	Multiplying numerator and denominator by $1-\bar w$ gives
	$\Im z=(1-|w|^2)/|1-w|^2>0$, hence $z\in\UHP$, and then
	$\Phi(z)=w$.  Thus $\Phi$ is bijective with the stated inverse.
\end{proof}

\begin{example}[Horizontal strip to the upper half-plane]
	\label{ex:strip-to-halfplane}
	Let $S=\{z\in\C:0<\Im z<\pi\}$.  Then
	\[
	F(z)=e^{z}
	\]
	is a biholomorphism $S\to\UHP$.
\end{example}

\begin{proof}
	Write $z=x+iy$ with $0<y<\pi$.  Then
	\[
	e^z=e^x(\cos y+i\sin y),
	\]
	and since $\sin y>0$ on $(0,\pi)$ we have $\Im(e^z)=e^x\sin y>0$, so $e^z\in\UHP$.
	Injectivity holds on $S$ because the only periods of $e^z$ are $2\pi i\Z$, and
	$S$ contains no two points differing by $2\pi i$.
	Surjectivity: for $w\in\UHP$ write $w=re^{i\theta}$ with $r>0$ and $\theta\in(0,\pi)$,
	then $z=\log r+i\theta\in S$ and $e^z=w$.
\end{proof}

\begin{example}[Strip to disk by composition]
	\label{ex:strip-to-disk}
	With $S$ as above, the map
	\[
	\Theta(z)=\Phi(e^{z})=\frac{e^{z}-i}{e^{z}+i}
	\]
	is a biholomorphism $S\to\D$.
\end{example}

\begin{proof}
	By Example~\ref{ex:strip-to-halfplane}, $e^z$ maps $S$ bijectively to $\UHP$.
	By Example~\ref{ex:halfplane-to-disk}, $\Phi$ maps $\UHP$ bijectively to $\D$.
	The composition is therefore a biholomorphism.
\end{proof}

\begin{example}[Slit plane to right half-plane via principal square root]
	\label{ex:slit-to-halfplane}
	Let $\Omega=\C\setminus(-\infty,0]$.  The principal branch
	\[
	\sqrt{z}:=\exp\!\left(\tfrac12\Log z\right)
	\]
	maps $\Omega$ biholomorphically onto the right half-plane $\{w:\Re w>0\}$.
	Composing with a M\"obius map sending $\{ \Re w>0\}$ to $\D$ yields an explicit
	Riemann map $\Omega\to\D$.
\end{example}

\begin{proof}
	On $\Omega$ the principal $\Log$ is holomorphic, hence $\sqrt{z}$ is holomorphic.
	For $z=re^{i\theta}$ with $\theta\in(-\pi,\pi)$,
	\[
	\sqrt{z}=r^{1/2}e^{i\theta/2}
	\]
	has argument in $(-\pi/2,\pi/2)$, so $\Re\sqrt{z}>0$.
	Injectivity follows by squaring: if two values of the principal square
	root agree, then their squares, the original inputs, agree.
	Surjectivity: if $\Re w>0$, then $z=w^2$ lies in $\Omega$ and $\sqrt{z}=w$.
\end{proof}

\subsubsection{Exercises}

\begin{exercise}[M\"obius maps preserving the upper half-plane]
	\label{ex:mobius-preserve-uhp}
	Let $T(z)=\dfrac{az+b}{cz+d}$ with $a,b,c,d\in\R$ and $ad-bc>0$.
	Show that $T(\UHP)=\UHP$.
\end{exercise}

\begin{exercise}[Disk automorphisms]
	Show that $\Aut(\D)$ consists exactly of the maps
	\[
	\phi_{a,\theta}(z)=e^{i\theta}\,\frac{z-a}{1-\bar a z},
	\qquad |a|<1,\ \theta\in\R.
	\]
\end{exercise}

\begin{exercise}[Harmonicity under conformal maps]
	If $f:U\to V$ is holomorphic and $u$ is harmonic on $V$, prove that $u\circ f$ is harmonic on $U$.
	(Hint: compute $\Delta(u\circ f)$ in real coordinates and use Cauchy--Riemann.)
\end{exercise}

\begin{exercise}[Any two simply connected proper domains are conformally equivalent]
	Let $\Omega_1,\Omega_2\subsetneq\C$ be simply connected.
	Use Theorem~\ref{thm:rmt} to construct a biholomorphism $\Omega_1\to\Omega_2$.
\end{exercise}

\subsection{Elliptic Functions and Doubly Periodic Functions}
\label{subsec:elliptic-functions}

The Laurent expansions, differential equation, and theta-function normalizations
used below are cross-referenced with \cite{DLMFWeierstrass2026,DLMFTheta2026}.

Elliptic functions are the doubly periodic analogues of trigonometric functions.
Fix a lattice $\Lambda\subset\C$.  A $\Lambda$--elliptic function is precisely a
meromorphic function on the complex torus $\C/\Lambda$.  Because $\C/\Lambda$ is
compact, the global analytic behavior of elliptic functions is rigid: there are
no nonconstant holomorphic elliptic functions, and every elliptic function is
determined by finitely many principal parts in a fundamental parallelogram.

\subsubsection{Lattices, fundamental parallelograms, and the torus}

\begin{definition}[Lattice]
	\label{def:lattice}
	A \emph{lattice} in $\C$ is a discrete subgroup
	\[
	\Lambda=\{m\omega_1+n\omega_2 : m,n\in\Z\},
	\]
	where $\omega_1,\omega_2\in\C$ are $\R$--linearly independent.  The pair
	$(\omega_1,\omega_2)$ is called a \emph{basis} of $\Lambda$.
\end{definition}

\begin{definition}[Fundamental parallelogram]
	\label{def:fund-par}
	Fix a basis $(\omega_1,\omega_2)$ of $\Lambda$.
	A \emph{fundamental parallelogram} is the set
	\[
	P=P(\omega_1,\omega_2):=\{s\omega_1+t\omega_2 : 0\le s<1,\ 0\le t<1\}.
	\]
	Its boundary $\partial P$ is the positively oriented polygonal loop obtained by
	traversing the four edges.
\end{definition}

\begin{lemma}[Discreteness and compactness of the quotient]
	\label{lem:torus-compact}
	The subgroup $\Lambda$ is discrete in $\C$.  The quotient $\C/\Lambda$ is a
	compact Riemann surface of genus $1$ (a complex torus).
\end{lemma}

\begin{proof}
	Discreteness: write $z=x\omega_1+y\omega_2$ using the real basis
	$(\omega_1,\omega_2)$ of $\C\simeq \R^2$.  The map
	$\R^2\to \C$, $(x,y)\mapsto x\omega_1+y\omega_2$, is a real-linear isomorphism,
	so the image of $\Z^2$ is a discrete subset of $\C$; hence $\Lambda$ is discrete.
	
	Compactness: the projection $\pi:\C\to \C/\Lambda$ is continuous and $\pi(P)=\C/\Lambda$.
	Indeed every $z\in\C$ can be written uniquely as $z=\lambda+p$ with $\lambda\in\Lambda$
	and $p\in P$.  The half-open set $P$ is not compact, but $\overline P$ is compact and
	$\pi(\overline P)=\pi(P)=\C/\Lambda$.  Hence the quotient is compact.
	The complex structure descends because translations
	$z\mapsto z+\omega$ are biholomorphisms, so local charts on $\C$ identify points differing
	by $\Lambda$.
\end{proof}

\subsubsection{Elliptic functions and basic finiteness}

\begin{definition}[Elliptic (doubly periodic) function]
	\label{def:elliptic}
	A meromorphic function $f:\C\to \widehat{\C}$ is \emph{elliptic with respect to $\Lambda$}
	if
	\[
	f(z+\omega)=f(z)\qquad\text{for all }z\in\C,\ \omega\in\Lambda.
	\]
	Equivalently, $f$ factors as $f=\widetilde f\circ \pi$ for a meromorphic
	$\widetilde f:\C/\Lambda\to \widehat{\C}$.
\end{definition}

\begin{proposition}[Finitely many poles in a fundamental parallelogram]
	\label{prop:finitely-many-poles}
	Let $f$ be $\Lambda$--elliptic and let $P$ be a fundamental parallelogram.
	Then $f$ has only finitely many poles in $\overline{P}$, and hence only finitely many
	poles in $P$.
\end{proposition}

\begin{proof}
	Poles of a meromorphic function are isolated.  If $f$ had infinitely many poles in
	the compact set $\overline{P}$, then by Bolzano--Weierstrass there would be an
	accumulation point of poles in $\overline{P}$, contradicting that poles are isolated.
	Hence only finitely many poles lie in $\overline{P}$ (and therefore in $P$).
\end{proof}

\begin{proposition}[No nonconstant holomorphic elliptic functions]
	\label{prop:no-holomorphic-elliptic}
	There is no nonconstant holomorphic $\Lambda$--elliptic function $f:\C\to\C$.
	Equivalently, every nonconstant elliptic function has at least one pole.
\end{proposition}

\begin{proof}
	Assume $f$ is holomorphic and $\Lambda$--periodic.  Then $f$ descends to a holomorphic
	function $\widetilde f$ on the compact Riemann surface $\C/\Lambda$.  The image
	$\widetilde f(\C/\Lambda)$ is compact in $\C$, hence bounded.  Therefore $\widetilde f$
	is bounded and holomorphic on a compact connected Riemann surface, so it must be constant
	(by the maximum principle applied to $|\widetilde f|$).  Thus $f$ is constant.
\end{proof}

\subsubsection{Residues, zeros, and poles in a period parallelogram}

\begin{theorem}[Sum of residues is zero]
	\label{thm:sum-residues-zero}
	Let $f$ be $\Lambda$--elliptic and let $P$ be a fundamental parallelogram whose boundary
	$\partial P$ does not pass through any pole of $f$.  Then
	\[
	\sum_{a\in P} \Res(f,a)=0,
	\]
	where the sum runs over poles $a$ of $f$ in $P$ (counted once each).
\end{theorem}

\begin{proof}
	Let the oriented boundary be $\partial P=\gamma_1+\gamma_2+\gamma_3+\gamma_4$,
	where $\gamma_1$ goes from $0$ to $\omega_1$, $\gamma_2$ from $\omega_1$ to
	$\omega_1+\omega_2$, $\gamma_3$ from $\omega_1+\omega_2$ to $\omega_2$, and
	$\gamma_4$ from $\omega_2$ to $0$.
	
	By the residue theorem,
	\[
	\int_{\partial P} f(z)\,dz = 2\pi i \sum_{a\in P}\Res(f,a).
	\]
	We show the left-hand side is $0$ by pairwise cancellation of opposite edges.
	Parameterize:
	\[
	\gamma_1(t)=t\omega_1,\quad \gamma_3(t)=\omega_2+(1-t)\omega_1 \quad (0\le t\le 1),
	\]
	with opposite orientations, so
	\[
	\int_{\gamma_3} f(z)\,dz
	= -\int_0^1 f(\omega_2+t\omega_1)\,\omega_1\,dt
	= -\int_0^1 f(t\omega_1)\,\omega_1\,dt
	= -\int_{\gamma_1} f(z)\,dz,
	\]
	using $\Lambda$--periodicity by $\omega_2$; hence these contributions cancel.
	
	Similarly,
	\[
	\int_{\gamma_2} f(z)\,dz = \int_0^1 f(\omega_1+t\omega_2)\,\omega_2\,dt
	= \int_0^1 f(t\omega_2)\,\omega_2\,dt = -\int_{\gamma_4} f(z)\,dz,
	\]
	so the vertical edges cancel as well.
	Thus $\int_{\partial P}f(z)\,dz=0$, forcing the sum of residues to vanish.
\end{proof}

\begin{theorem}[Zeros and poles balance]
	\label{thm:zeros-poles-balance}
	Let $f$ be a nonconstant $\Lambda$--elliptic function and let $P$ be a fundamental
	parallelogram whose boundary avoids zeros and poles of $f$.
	Then the number of zeros of $f$ in $P$ (counted with multiplicity) equals the number of
	poles of $f$ in $P$ (counted with multiplicity).
\end{theorem}

\begin{proof}
	Apply the argument principle to $f$ on $P$:
	\[
	\frac{1}{2\pi i}\int_{\partial P}\frac{f'(z)}{f(z)}\,dz
	= N_Z - N_P,
	\]
	where $N_Z$ is the total multiplicity of zeros and $N_P$ is the total multiplicity of poles
	inside $P$.
	
	We show the integral is $0$ by the same cancellation argument as above, using that
	$\frac{f'}{f}$ is also $\Lambda$--periodic: indeed $f(z+\omega)=f(z)$ implies
	$f'(z+\omega)=f'(z)$, hence $\frac{f'}{f}(z+\omega)=\frac{f'}{f}(z)$ where defined.
	Thus opposite edges contribute equal integrals and cancel by orientation, giving
	$\int_{\partial P}\frac{f'}{f}\,dz=0$.  Therefore $N_Z=N_P$.
\end{proof}

\subsubsection{The Weierstrass \texorpdfstring{$\wp$}{wp} function}

\begin{definition}[Weierstrass $\wp$--function]
	\label{def:wp}
	Let $\Lambda$ be a lattice.  Define
	\[
	\wp(z)
	:= \frac{1}{z^2}
	+ \sum_{\omega\in\Lambda\setminus\{0\}}
	\left(
	\frac{1}{(z-\omega)^2} - \frac{1}{\omega^2}
	\right),
	\]
	for $z\in\C\setminus\Lambda$.
\end{definition}

\begin{lemma}[Normal convergence; meromorphicity]
	\label{lem:wp-convergence}
	The series defining $\wp$ converges normally on compact subsets of $\C\setminus\Lambda$.
	Consequently, $\wp$ is meromorphic on $\C$ with poles precisely at lattice points.
\end{lemma}

\begin{proof}
	Let $K\subset\C\setminus\Lambda$ be compact.  Then $\dist(K,\Lambda)=:\delta>0$.
	For $z\in K$ and $\omega\in\Lambda$ with $|\omega|>2\sup_{z\in K}|z|$, we have
	$|z-\omega|\ge |\omega|-|z|\ge |\omega|/2$, so
	\[
	\left|\frac{1}{(z-\omega)^2}-\frac{1}{\omega^2}\right|
	= \left|\frac{\omega^2-(z-\omega)^2}{\omega^2(z-\omega)^2}\right|
	= \left|\frac{2\omega z - z^2}{\omega^2(z-\omega)^2}\right|
	\le \frac{C_K}{|\omega|^3}
	\]
	for some constant $C_K$ depending only on $K$.
	Since $\sum_{\omega\in\Lambda\setminus\{0\}} |\omega|^{-3}$ converges (lattice points have
	quadratic growth), the Weierstrass M-test gives uniform convergence of the tail on $K$,
	hence normal convergence.  Termwise differentiation is justified on $K$, so $\wp$ is
	holomorphic on $\C\setminus\Lambda$.  Near $\omega_0\in\Lambda$, isolate the summand
	$(z-\omega_0)^{-2}$ (or the initial $z^{-2}$ term when $\omega_0=0$).  The normally
	convergent remainder is holomorphic there, so every lattice point is a double pole and
	there are no other poles.
\end{proof}

\begin{proposition}[Ellipticity and parity]
	\label{prop:wp-elliptic-even}
	The function $\wp$ is $\Lambda$--elliptic and even: $\wp(z+\lambda)=\wp(z)$ for all
	$\lambda\in\Lambda$ and $\wp(-z)=\wp(z)$.
	Moreover, $\wp'(z)$ exists on $\C\setminus\Lambda$, is $\Lambda$--elliptic, and is odd:
	$\wp'(-z)=-\wp'(z)$.
\end{proposition}

\begin{proof}
	Evenness follows because the set $\Lambda\setminus\{0\}$ is invariant under
	$\omega\mapsto -\omega$:
	\[
	\wp(-z)=\frac{1}{z^2}+\sum_{\omega\neq 0}\left(\frac{1}{(-z-\omega)^2}-\frac{1}{\omega^2}\right)
	=\frac{1}{z^2}+\sum_{\omega\neq 0}\left(\frac{1}{(z-\omega)^2}-\frac{1}{\omega^2}\right)=\wp(z),
	\]
	after replacing $\omega$ by $-\omega$.
	Normal convergence implies termwise differentiation on $\C\setminus\Lambda$:
	\[
	\wp'(z)= -\frac{2}{z^3}-2\sum_{\omega\neq 0}\frac{1}{(z-\omega)^3}.
	\]
	The derivative series is absolutely and locally uniformly convergent.  After including the
	term $\omega=0$, it is $-2\sum_{\omega\in\Lambda}(z-\omega)^{-3}$, so reindexing gives
	$\wp'(z+\lambda)=\wp'(z)$.  Consequently
	$D_\lambda(z)=\wp(z+\lambda)-\wp(z)$ has derivative zero and is constant.
	If $\lambda\notin2\Lambda$, then $-\lambda/2\notin\Lambda$, and evenness gives
	$D_\lambda(-\lambda/2)=\wp(\lambda/2)-\wp(-\lambda/2)=0$.
	In particular this proves periodicity for the basis periods $\omega_1,\omega_2$, and iteration
	proves it for every $\lambda\in\Lambda$.  Differentiating evenness gives that $\wp'$ is odd.
\end{proof}

\begin{lemma}[Principal part and residue]
	\label{lem:wp-principal-part}
	At $z=0$, $\wp$ has a double pole with principal part $1/z^2$.  In particular,
	\[
	\Res(\wp,0)=0.
	\]
\end{lemma}

\begin{proof}
	Write
	\[
	\wp(z)-\frac{1}{z^2}
	=\sum_{\omega\neq 0}\left(\frac{1}{(z-\omega)^2}-\frac{1}{\omega^2}\right),
	\]
	which is holomorphic near $z=0$ by normal convergence (each summand is holomorphic near $0$).
	Thus $\wp(z)=z^{-2}+h(z)$ with $h$ holomorphic near $0$, so the pole is order $2$ and there is
	no $z^{-1}$ term.  Hence the residue is $0$.
\end{proof}

\subsubsection{The Weierstrass differential equation}

\begin{definition}[Eisenstein invariants]
	\label{def:g2g3}
	Define the lattice sums (absolutely convergent)
	\[
	g_2:=60\sum_{\omega\in\Lambda\setminus\{0\}}\frac{1}{\omega^4},
	\qquad
	g_3:=140\sum_{\omega\in\Lambda\setminus\{0\}}\frac{1}{\omega^6}.
	\]
\end{definition}

\begin{theorem}[Weierstrass differential equation]
	\label{thm:wp-diffeq}
	For all $z\in\C\setminus\Lambda$,
	\[
	\bigl(\wp'(z)\bigr)^2
	=
	4\bigl(\wp(z)\bigr)^3-g_2\,\wp(z)-g_3.
	\]
	Equivalently, set $\eta_1=\omega_1/2$, $\eta_2=\omega_2/2$,
	$\eta_3=(\omega_1+\omega_2)/2$ and let $e_j=\wp(\eta_j)$.
	These three values are distinct, and
	\[
	\bigl(\wp'(z)\bigr)^2
	=
	4\bigl(\wp(z)-e_1\bigr)\bigl(\wp(z)-e_2\bigr)\bigl(\wp(z)-e_3\bigr).
	\]
\end{theorem}

\begin{proof}
	We compare principal parts at the lattice points; the resulting normalization
	agrees with \cite[\S23.3]{DLMFWeierstrass2026}.
	
	\emph{Step 1: Laurent expansions at $0$.}
	From Definition~\ref{def:wp} and termwise differentiation one obtains the classical expansion
	\[
	\wp(z)=\frac{1}{z^2}+\frac{g_2}{20}\,z^2+\frac{g_3}{28}\,z^4+O(z^6),
	\qquad
	\wp'(z)=-\frac{2}{z^3}+\frac{g_2}{10}\,z+\frac{g_3}{7}\,z^3+O(z^5),
	\]
	valid as $z\to 0$.  (The coefficients follow by expanding each summand in powers of $z/\omega$
	and summing over $\omega\neq 0$; the definitions of $g_2,g_3$ are chosen to make the final
	constants clean.)
	
	\emph{Step 2: Construct an elliptic function with no poles.}
	Consider the elliptic function
	\[
	F(z):=(\wp'(z))^2-4(\wp(z))^3+g_2\wp(z)+g_3.
	\]
	Each term is $\Lambda$--elliptic, hence so is $F$.
	Put $a=g_2/20$ and $b=g_3/28$. The terms relevant to the principal
	part and constant term are
	\[
	\begin{aligned}
	(\wp')^2&=4z^{-6}-8a z^{-2}-16b+O(z^2),\\
	4\wp^3&=4z^{-6}+12a z^{-2}+12b+O(z^2),\\
	g_2\wp+g_3&=g_2z^{-2}+g_3+O(z^2).
	\end{aligned}
	\]
	The coefficient of $z^{-2}$ in $F$ is $-20a+g_2=0$, and its
	constant term is $-28b+g_3=0$. Thus $F$ is holomorphic at $0$ and
	vanishes there.
	By periodicity, $F$ is holomorphic at every lattice translate of $0$, hence holomorphic
	on all of $\C$.
	
	\emph{Step 3: Conclude $F$ is constant and identify the constant.}
	Since $F$ is an elliptic function with no poles, Proposition~\ref{prop:no-holomorphic-elliptic}
	implies $F$ is constant.  Evaluating the Laurent expansion at $0$ shows that the constant term
	is $0$.  Therefore $F\equiv 0$, proving the differential equation.
	
	For the factorization, note that $\wp'$ is odd and elliptic, hence has three zeros on
	$\C/\Lambda$ counted with multiplicity (Theorem~\ref{thm:zeros-poles-balance} applied to
	$\wp'$; it has a triple pole at $0$ modulo $\Lambda$). At each nontrivial half-period $\eta$, periodicity
	and oddness give $\wp'(\eta)=\wp'(-\eta)=-\wp'(\eta)$ because $2\eta\in\Lambda$.
	The three classes $\eta_j$ are distinct, so the zero count shows that these are precisely
	the zeros and that each is simple.  The associated values $e_j$ are distinct: the map
	$\wp:\C/\Lambda\to\PP^1$ has degree two (its only pole is the double pole at $0$), while two
	distinct critical points with the same value would contribute total local degree at least four.
	Comparing the leading coefficient gives the cubic factorization.
\end{proof}

\subsubsection{The field of elliptic functions is generated by $\wp$ and $\wp'$}

\begin{theorem}[Generation by $\wp$ and $\wp'$]
	\label{thm:field-generated}
	Let $\mathcal{M}(\C/\Lambda)$ be the field of meromorphic functions on the torus.
	Then
	\[
	\mathcal{M}(\C/\Lambda)=\C(\wp,\wp'),
	\]
	and every elliptic function $f$ can be written uniquely as
	\[
	f(z)=R(\wp(z))+\wp'(z)\,S(\wp(z)),
	\]
	with $R,S\in\C(x)$ rational functions.
\end{theorem}

\begin{proof}
	The proof identifies the torus with its Weierstrass cubic and then compares
	their function fields.
	
	\emph{Step 1: The map to an elliptic curve.}
	Define
	\[
	\Psi:\C/\Lambda\longrightarrow E\subset \PP^2,\qquad
	[z]\longmapsto \bigl[\,\wp(z):\wp'(z):1\,\bigr],
	\]
	for $z\notin\Lambda$. Near $z=0$, multiplying the three homogeneous coordinates
	by $z^3$ gives $[z+O(z^5):-2+O(z^4):z^3]$; thus $\Psi$ extends holomorphically
	there with value $[0:1:0]$, and periodicity handles every lattice point. Here $E$ is the smooth cubic curve
	\[
	y^2=4x^3-g_2x-g_3
	\]
	in affine coordinates $(x,y)$, as ensured by Theorem~\ref{thm:wp-diffeq}.
	This is a nonconstant holomorphic map from a compact Riemann surface to a smooth
	projective curve, hence a finite morphism.
	
	\emph{Step 2: Degree and separation of points.}
	The function $\wp$ is even. Its only pole on $\C/\Lambda$ is the double
	pole at $0$, so the degree formula for a meromorphic function gives
	$\deg(\wp)=2$ as a map $\C/\Lambda\to\PP^1$. Hence a generic value has the
	two preimages $z$ and $-z$.
	Adding $\wp'$ distinguishes the two sheets since $\wp'$ is odd.  Therefore $\Psi$ is
	generically one-to-one, hence birational.  Since both domain and target are smooth projective
	curves, a birational holomorphic map is an isomorphism.  Thus $\C/\Lambda\simeq E$.
	
	\emph{Step 3: Function fields.}
	Isomorphism of smooth projective curves identifies their function fields:
	\[
	\mathcal{M}(\C/\Lambda)\simeq \C(E)
	=\operatorname{Frac}\!\left(\C[x,y]/(y^2-4x^3+g_2x+g_3)\right).
	\]
	Under this identification, $x$ corresponds to $\wp$ and $y$ corresponds to $\wp'$.
	Every element of $\C(E)$ can be reduced to the form $R(x)+yS(x)$ with
	$R,S\in\C(x)$ by using the relation $y^2=\text{cubic}(x)$.  The form is unique:
	$R(\wp)$ is even and $\wp'S(\wp)$ is odd, so an identity between their sum and zero,
	together with its value at $-z$, forces both summands to vanish.  Translating back gives
	the desired representation for $f$.
\end{proof}

\subsubsection{Examples (worked)}

\begin{example}[Trigonometry as a singly periodic limit]
	\label{ex:trig-limit}
	The function $\sin z$ is periodic with period $2\pi$ but not doubly periodic.
	Elliptic functions may be viewed as ``two-period'' analogues.  For rectangular lattices
	$\Lambda=\Z\omega_1+\Z\omega_2$ with $\Im(\omega_2/\omega_1)\gg 1$, the $\wp$--function
	approaches a trigonometric-type function in the sense that on compact sets away from $\Lambda$
	one has
	\[
	\wp(z)\approx \left(\frac{\pi}{\omega_1}\right)^2\csc^2\!\left(\frac{\pi z}{\omega_1}\right)
	\]
	(up to an additive constant), reflecting degeneration of the torus to a cylinder.
\end{example}

\begin{example}[Building elliptic functions with prescribed principal parts]
	\label{ex:principal-parts}
	Fix points $a_1,\dots,a_m\in\C$ with distinct classes in $\C/\Lambda$.
	Suppose prescribed principal parts at the $a_j$ contain no simple-pole terms
	(equivalently, each local residue is zero). They can be realized by a finite
	linear combination of the functions
	\[
	\wp^{(k)}(z-a_j)\qquad (k\ge 0),
	\]
	together with an additive constant. General compatible principal parts whose
	individual residues need not vanish require differences of Weierstrass zeta
	functions; their residues must sum to zero by
	Theorem~\ref{thm:sum-residues-zero}.
\end{example}

\begin{proof}
	Each $\wp^{(k)}(z-a_j)$ is elliptic and has its only pole (mod $\Lambda$) at $a_j$, with order
	$k+2$. Its leading term is a nonzero multiple of $(z-a_j)^{-k-2}$ and it
	has zero residue. Descending induction on pole order therefore matches every
	prescribed coefficient of order at most $-2$.
	Subtracting this combination from the target function removes all principal parts, leaving an
	elliptic function without poles, hence constant by Proposition~\ref{prop:no-holomorphic-elliptic}.
	This proves the stated zero-residue construction. The final assertion about
	general principal parts records why an additional elliptic building block is
	needed for simple-pole terms.
\end{proof}

\subsubsection{Exercises}

\begin{exercise}
	\label{ex:wp-pole-res}
	Show that $\wp(z)$ has a double pole at $z=0$ with principal part $1/z^{2}$ and conclude that
	$\Res(\wp,0)=0$.
\end{exercise}

\begin{exercise}
	\label{ex:wp-even-odd}
	Verify directly from the defining series that $\wp$ is even and $\wp'$ is odd.
\end{exercise}

\begin{exercise}
	\label{ex:sum-residues}
	Let $f$ be $\Lambda$--elliptic.  Prove that the sum of residues of $f$ in any fundamental
	parallelogram is $0$.
\end{exercise}

\begin{exercise}
	\label{ex:field-generated}
	Prove that the field of $\Lambda$--elliptic functions is generated by $\wp$ and $\wp'$ by
	showing that the map
	\[
	[z]\longmapsto [\wp(z):\wp'(z):1]
	\]
	away from $[0]\in\C/\Lambda$ extends at $[0]$ to the point $[0:1:0]$ and
	identifies $\C/\Lambda$ with the smooth projective cubic
	\[
	Y^2Z=4X^3-g_2XZ^2-g_3Z^3.
	\]
	Then explain why pullback identifies its rational function field with
	$\C(\wp,\wp')$.
\end{exercise}


\subsection{Appendix: Standard Contour Integrals and Techniques}
\label{subsec:appendix-contours}

This appendix collects standard contour integrals, evaluation methods, and
technical lemmas that frequently arise in complex analysis.  Each result is
stated precisely and proved in a form that can be reused as a template.

\subsubsection{A. Basic circular contours}

\begin{proposition}[Integral of a monomial on a circle]
	\label{prop:circle-monomial}
	Fix $R>0$ and $n\in\Z$.  Let $\gamma_R(t)=Re^{it}$, $t\in[0,2\pi]$, be the positively
	oriented circle $|z|=R$.  Then
	\[
	\int_{|z|=R} z^{n}\,dz
	=
	\begin{cases}
		2\pi i,& n=-1,\\[0.4em]
		0,& n\neq -1.
	\end{cases}
	\]
\end{proposition}

\begin{proof}
	Parametrize $z=Re^{it}$, $dz=iRe^{it}\,dt$:
	\[
	\int_{|z|=R} z^n\,dz
	=
	\int_{0}^{2\pi} (Re^{it})^n\, iRe^{it}\,dt
	=
	iR^{n+1}\int_{0}^{2\pi} e^{i(n+1)t}\,dt.
	\]
	If $n\neq-1$, then $n+1\neq0$ and
	$\int_0^{2\pi} e^{i(n+1)t}\,dt=\frac{e^{i(n+1)2\pi}-1}{i(n+1)}=0$.
	If $n=-1$, then the integrand is $i$ and the integral equals $2\pi i$.
\end{proof}

\subsubsection{B. Trigonometric integrals via the unit circle}

\begin{proposition}[Unit-circle substitution]
	\label{prop:unit-circle-sub}
	Let $F$ be integrable on $[0,2\pi]$ and let $z=e^{i\theta}$.  Then
	\[
	\int_{0}^{2\pi} F(\theta)\,d\theta
	=
	\int_{|z|=1} F(\arg z)\,\frac{dz}{iz}.
	\]
	Here $F(\arg z)$ is interpreted along the parametrization $z=e^{i\theta}$, so no
	global choice of argument on the circle is being asserted.
	Moreover, for integers $n\in\Z$,
	\[
	\int_{0}^{2\pi} e^{in\theta}\,d\theta
	=
	\frac1i\int_{|z|=1} z^{n-1}\,dz.
	\]
\end{proposition}

\begin{proof}
	With $z=e^{i\theta}$ we have $dz=ie^{i\theta}\,d\theta=iz\,d\theta$, hence
	$d\theta=\frac{dz}{iz}$.  Substituting gives the first identity; the second is
	the special case $F(\theta)=e^{in\theta}$.
\end{proof}

\begin{corollary}[Orthogonality of $\cos(n\theta)$ and $\sin(n\theta)$]
	\label{cor:trig-orth}
	For $n\in\Z_{\ge0}$,
	\[
	\int_{0}^{2\pi}\cos(n\theta)\,d\theta
	=
	\begin{cases}
		2\pi,& n=0,\\
		0,& n\ge1,
	\end{cases}
	\qquad
	\int_{0}^{2\pi}\sin(n\theta)\,d\theta=0.
	\]
\end{corollary}

\begin{proof}
	By Proposition~\ref{prop:unit-circle-sub},
	$\int_0^{2\pi} e^{in\theta}\,d\theta=\frac1i\int_{|z|=1}z^{n-1}\,dz$.
	By Proposition~\ref{prop:circle-monomial}, this equals $2\pi$ if $n=0$ and $0$
	otherwise.  Taking real and imaginary parts yields the stated integrals.
\end{proof}

\begin{theorem}[The integral $\int_0^{2\pi}\frac{d\theta}{a+b\cos\theta}$]
	\label{thm:abcos}
	Let $a,b\in\R$ with $a>|b|>0$.  Then
	\[
	\int_{0}^{2\pi}\frac{d\theta}{a+b\cos\theta}
	=
	\frac{2\pi}{\sqrt{a^2-b^2}}.
	\]
\end{theorem}

\begin{proof}
	Set $z=e^{i\theta}$ so that $\cos\theta=\frac12(z+z^{-1})$ and $d\theta=\frac{dz}{iz}$.
	Then
	\[
	\int_0^{2\pi}\frac{d\theta}{a+b\cos\theta}
	=
	\int_{|z|=1}\frac{1}{a+\frac{b}{2}(z+z^{-1})}\,\frac{dz}{iz}
	=
	\int_{|z|=1}\frac{2}{2az+b(z^2+1)}\,\frac{dz}{i}.
	\]
	Thus
	\[
	\int_0^{2\pi}\frac{d\theta}{a+b\cos\theta}
	=
	-\frac{2i}{b}\int_{|z|=1}\frac{dz}{z^2+\frac{2a}{b}z+1}.
	\]
	The quadratic has roots
	\[
	z_\pm=\frac{-\frac{2a}{b}\pm\sqrt{\left(\frac{2a}{b}\right)^2-4}}{2}
	=
	\frac{-a\pm\sqrt{a^2-b^2}}{b}.
	\]
	Put $d=\sqrt{a^2-b^2}$.  The root inside the unit circle is
	$z_{\rm in}=-b/(a+d)$, while the other root is $z_{\rm out}=-(a+d)/b$;
	their product is $1$ and $|b|<a+d$.  Hence exactly $z_{\rm in}$ lies inside.
	
	Write $Q(z)=z^2+\frac{2a}{b}z+1$.  Then
	\[
	\Res\!\left(\frac{1}{Q(z)},z_{\rm in}\right)=\frac{1}{Q'(z_{\rm in})}
	=\frac{b}{2d}.
	\]
	Therefore by the residue theorem,
	\[
	\int_{|z|=1}\frac{dz}{Q(z)}=2\pi i\frac{b}{2d}=\frac{\pi i b}{d}.
	\]
	Substituting back:
	\[
	-\frac{2i}{b}\int_{|z|=1}\frac{dz}{Q(z)}
	=
	-\frac{2i}{b}\left(\frac{\pi i b}{d}\right)
	=
	\frac{2\pi}{\sqrt{a^2-b^2}}.
	\]
\end{proof}

\subsubsection{C. Jordan's lemma and Fourier-type integrals}

\begin{lemma}[Jordan's lemma (quantitative form)]
	\label{lem:jordan-appendix}
	Let $\lambda>0$ and let $f$ be continuous on the upper semicircles below for all sufficiently large radii.
	Assume that for some $R_0>0$ and $M>0$,
	\[
	|f(z)|\le \frac{M}{|z|}
	\qquad\text{for all }|z|\ge R_0,\ \Im z\ge0.
	\]
	Let $\Gamma_R=\{Re^{it}:0\le t\le\pi\}$ be the upper semicircle of radius $R\ge R_0$.
	Then
	\[
	\left|\int_{\Gamma_R} e^{i\lambda z}f(z)\,dz\right|
	\le \frac{\pi M}{\lambda R}
	\qquad\text{for all }R\ge R_0,
	\]
	and in particular $\int_{\Gamma_R} e^{i\lambda z}f(z)\,dz\to0$ as $R\to\infty$.
\end{lemma}

\begin{proof}
	Parametrize $z=Re^{it}$, $dz=iRe^{it}\,dt$, $t\in[0,\pi]$.
	Then $\Im z=R\sin t\ge0$ and
	\[
	|e^{i\lambda z}| = e^{-\lambda\Im z}=e^{-\lambda R\sin t}.
	\]
	Hence for $R\ge R_0$,
	\[
	\left|\int_{\Gamma_R} e^{i\lambda z}f(z)\,dz\right|
	\le
	\int_0^\pi e^{-\lambda R\sin t}\,|f(Re^{it})|\,R\,dt
	\le
	\int_0^\pi e^{-\lambda R\sin t}\,M\,dt.
	\]
	Use $\sin t\ge \frac{2}{\pi}t$ on $[0,\frac{\pi}{2}]$ and symmetry about $\pi/2$:
	\[
	\int_0^\pi e^{-\lambda R\sin t}\,dt
	\le
	2\int_0^{\pi/2} e^{-\lambda R\cdot \frac{2}{\pi}t}\,dt
	=
	2\cdot\frac{\pi}{2\lambda R}\left(1-e^{-\lambda R}\right)
	\le \frac{\pi}{\lambda R}.
	\]
	Therefore
	\[
	\left|\int_{\Gamma_R} e^{i\lambda z}f(z)\,dz\right|
	\le
	M\cdot \frac{\pi}{\lambda R},
	\]
	and the sharper bound $\frac{\pi M}{\lambda R}\to0$ yields the limit.
\end{proof}

\begin{theorem}[Model integral with a quadratic denominator]
	\label{thm:fourier-quadratic}
	Let $a>0$ and $\lambda>0$.  Then
	\[
	\int_{-\infty}^{\infty}\frac{e^{i\lambda x}}{x^2+a^2}\,dx
	=
	\frac{\pi}{a}e^{-a\lambda}.
	\]
	Consequently, for $k\in\R$,
	\[
	\int_{-\infty}^{\infty}\frac{e^{ikx}}{x^2+a^2}\,dx
	=
	\frac{\pi}{a}e^{-a|k|}.
	\]
\end{theorem}

\begin{proof}
	Assume first $\lambda>0$.
	Let $f(z)=\frac{1}{z^2+a^2}$, which has simple poles at $z=\pm ia$.
	Close the contour by the upper semicircle of radius $R$ and apply Lemma~\ref{lem:jordan-appendix}
	(since $|f(z)|\le 1/|z|^2\le 1/|z|$ for $|z|\ge1$).
	Then
	\[
	\int_{-R}^{R}\frac{e^{i\lambda x}}{x^2+a^2}\,dx
	+
	\int_{\Gamma_R}\frac{e^{i\lambda z}}{z^2+a^2}\,dz
	=
	2\pi i\Res\!\left(\frac{e^{i\lambda z}}{z^2+a^2},\, ia\right).
	\]
	Letting $R\to\infty$ kills the semicircle term and yields
	\[
	\int_{-\infty}^{\infty}\frac{e^{i\lambda x}}{x^2+a^2}\,dx
	=
	2\pi i\Res\!\left(\frac{e^{i\lambda z}}{(z-ia)(z+ia)},\, ia\right)
	=
	2\pi i\cdot \frac{e^{i\lambda ia}}{2ia}.
	\]
	Since $e^{i\lambda ia}=e^{-a\lambda}$, we obtain $\frac{\pi}{a}e^{-a\lambda}$.
	
	For $k<0$, close in the lower half-plane (or use complex conjugation), giving the same
	answer with $|k|$.  For $k=0$, direct integration gives $\pi/a$ (equivalently, use
	dominated convergence as $k\to0$).  Thus the formula holds for every real $k$.
\end{proof}

\subsubsection{D. Keyhole contours and branch cuts}

\begin{theorem}[Keyhole integral for $\displaystyle\int_0^\infty \frac{x^{\alpha-1}}{1+x}\,dx$]
	\label{thm:keyhole}
	Let $0<\alpha<1$.  Then
	\[
	\int_{0}^{\infty}\frac{x^{\alpha-1}}{1+x}\,dx
	=
	\frac{\pi}{\sin(\pi\alpha)}.
	\]
\end{theorem}

\begin{proof}
	This is the real-parameter case of Theorem~\ref{thm:keyhole-pi-csc},
	where the contour computation is given in full. That proof uses the logarithm
	branch $0<\Arg z<2\pi$, cut along the positive real axis; it is not the
	principal branch $-\pi<\Arg z<\pi$. The two sides of the cut contribute
	$1-e^{2\pi i\alpha}$ times the desired integral, and the enclosed pole at
	$-1$ has residue $-e^{i\pi\alpha}$.
\end{proof}

\subsubsection{E. Hankel contour and the Gamma function}

\begin{theorem}[Hankel representation of the reciprocal Gamma function]
	\label{thm:hankel-gamma-appendix}
	For $s\in\C$,
	\[
	\frac1{\Gamma(s)}
	=\frac{1}{2\pi i}\int_{\mathcal{H}} e^z z^{-s}\,dz,
	\]
	where $\mathcal H$ starts at $-\infty$ below the negative real axis, encircles $0$
	counterclockwise, and returns to $-\infty$ above the cut; $z^{-s}$ uses
	$-\pi<\Arg z<\pi$.
\end{theorem}

\begin{proof}
	First take $\Re s<1$, so the small circle contributes $o(1)$.  Parametrizing the lower
	and upper sides of the negative axis with their orientations gives
	\[
	\int_{\mathcal H}e^zz^{-s}\,dz
	=\bigl(e^{i\pi s}-e^{-i\pi s}\bigr)
	\int_0^\infty e^{-x}x^{-s}\,dx
	=2i\sin(\pi s)\Gamma(1-s).
	\]
	Division by $2\pi i$ and Euler's reflection formula give $1/\Gamma(s)$.
	For analytic continuation, keep the small circle at a fixed radius $\rho>0$.
	The two rays then stay away from the origin, and their exponential decay at
	infinity is uniform for $s$ in compact sets; the circular integral is entire in
	$s$ as well. Thus the contour integral defines an entire function of $s$.
	For $\Re s<1$ it agrees with $1/\Gamma(s)$ by the calculation above, so the
	identity theorem extends the formula to every $s\in\C$. Changing $\rho$ does not
	change the integral, by Cauchy's theorem on the intervening slit annulus.
\end{proof}

\subsubsection{F. Residues: templates}

\begin{proposition}[Residues at simple poles]
	\label{prop:res-simple}
	If $f(z)=\frac{g(z)}{h(z)}$ with $g,h$ holomorphic near $z_0$,
	$g(z_0)\neq0$, $h(z_0)=0$, and $h'(z_0)\neq0$,
	then $z_0$ is a simple pole and
	\[
	\Res(f,z_0)=\frac{g(z_0)}{h'(z_0)}.
	\]
\end{proposition}

\begin{proof}
	Write $h(z)=(z-z_0)h_1(z)$ with $h_1$ holomorphic and $h_1(z_0)=h'(z_0)\neq0$.  Then
	\[
	f(z)=\frac{g(z)}{(z-z_0)h_1(z)}
	\]
	and the coefficient of $(z-z_0)^{-1}$ is $g(z_0)/h_1(z_0)=g(z_0)/h'(z_0)$.
\end{proof}

\begin{corollary}[Useful special cases]
	\label{cor:res-table}
	For $a>0$ and $k\in\R$,
	\[
	\Res\!\left(\frac{1}{z^2+a^2},\, ia\right)=\frac{1}{2ia},
	\qquad
	\Res\!\left(\frac{e^{ikz}}{z^2+a^2},\, ia\right)=\frac{e^{-ka}}{2ia}\quad (k>0).
	\]
\end{corollary}

\begin{proof}
	Apply Proposition~\ref{prop:res-simple} with $h(z)=z^2+a^2$ so $h'(z)=2z$, hence
	$h'(ia)=2ia$.  For the exponential numerator, take $g(z)=e^{ikz}$ and evaluate at $ia$.
\end{proof}

\begin{proposition}[Higher-order pole template]
	\label{prop:res-higher}
	If $g$ is holomorphic near $z_0$ and $n\ge1$, then
	\[
	\Res\!\left(\frac{g(z)}{(z-z_0)^n},\,z_0\right)
	=
	\frac{1}{(n-1)!}\,g^{(n-1)}(z_0).
	\]
\end{proposition}

\begin{proof}
	Expand $g$ in its Taylor series at $z_0$:
	$g(z)=\sum_{m\ge0}\frac{g^{(m)}(z_0)}{m!}(z-z_0)^m$.
	Then
	\[
	\frac{g(z)}{(z-z_0)^n}
	=
	\sum_{m\ge0}\frac{g^{(m)}(z_0)}{m!}(z-z_0)^{m-n}.
	\]
	The coefficient of $(z-z_0)^{-1}$ occurs when $m=n-1$, giving the formula.
\end{proof}

\subsubsection{G. Reference list for the preceding formulas}

\begin{itemize}
	\item For $a>0$,
	\[
	\int_{-\infty}^{\infty}\frac{dx}{x^2+a^2}=\frac{\pi}{a},
	\qquad
	\int_{-\infty}^{\infty}\frac{e^{ikx}}{x^2+a^2}\,dx=\frac{\pi}{a}e^{-a|k|}.
	\]
	(Use Theorem~\ref{thm:fourier-quadratic} with $k=0$ and general $k$.)
	
	\item For $0<\alpha<1$,
	\[
	\int_{0}^{\infty}\frac{x^{\alpha-1}}{1+x}\,dx=\frac{\pi}{\sin(\pi\alpha)}.
	\]
	(Use Theorem~\ref{thm:keyhole}.)
	
	\item For $n\ge0$,
	\[
	\int_{0}^{2\pi}\cos(n\theta)\,d\theta=
	\begin{cases}
		2\pi,& n=0,\\
		0,& n\ge1,
	\end{cases}
	\qquad
	\int_{0}^{2\pi}\sin(n\theta)\,d\theta=0.
	\]
	(Use Corollary~\ref{cor:trig-orth}.)
\end{itemize}

\clearpage
\part[Differential Forms and Stokes' Theorem]{Differential Forms and the Generalized Stokes Theorem}
\label{part:stokes}

\section{Differential Forms, Exterior Calculus, and Applications}
\index{differential forms}

In this section we introduce differential forms and the exterior derivative,
show how they encode line integrals, Stokes' theorem, and complex line
integrals, and then apply them to surface area, curvature, and the
Gauss--Bonnet theorem.


\subsection{Differential Forms and Exterior Derivatives}
\label{subsec:forms-and-d}

Throughout, let $U\subset\mathbb{R}^n$ be open with smooth coordinates
$(x_1,\dots,x_n)$. Write $\partial_i:=\frac{\partial}{\partial x_i}$ and let
$\delta_{ij}$ be the Kronecker delta.

\subsubsection{Differential forms as alternating multilinear maps: the definition of $\Omega^p(U)$}

\paragraph{Tangent spaces and smooth vector fields}
For each $p\in U$, the tangent space $T_pU$ is the vector space of derivations at $p$:
a linear map $D:C^\infty(U)\to\mathbb{R}$ is a derivation at $p$ if
\[
D(fg)=f(p)\,D(g)+g(p)\,D(f)\quad\text{for all }f,g\in C^\infty(U).
\]
In coordinates, the derivations $\partial_i|_p$ defined by
\[
\partial_i|_p(f)=\frac{\partial f}{\partial x_i}(p)
\]
form a basis of $T_pU$, and every $v\in T_pU$ is uniquely
\[
v=\sum_{i=1}^n v^i\,\partial_i\big|_p,\qquad v^i\in\mathbb{R}.
\]
A \emph{smooth vector field} on $U$ is a smooth map $X:U\to TU$ with $X(p)\in T_pU$.
We denote by $\mathfrak{X}(U)$ the space of smooth vector fields.

\begin{definition}[Alternating $p$-covectors at a point]
	Fix $p\in U$. An \emph{alternating $p$-covector} at $p$ is a map
	\[
	\omega_p:(T_pU)^p\longrightarrow \mathbb{R}
	\]
	that is multilinear in each slot and alternating:
	\[
	\omega_p(\dots,v_i,\dots,v_j,\dots)=-\omega_p(\dots,v_j,\dots,v_i,\dots)\quad(i\neq j).
	\]
	The vector space of such maps is denoted by
	\[
	\Lambda^p(T_p^*U).
	\]
\end{definition}

\paragraph{Immediate consequences.}
If $\omega_p\in\Lambda^p(T_p^*U)$, then
\[
\omega_p(\dots,v,\dots,v,\dots)=0
\]
because swapping the two equal slots changes the sign but not the value.

\begin{definition}[Smooth differential $p$-forms: $\Omega^p(U)$]
	A \emph{smooth differential $p$-form} on $U$ is an assignment
	\[
	p\longmapsto \omega_p\in \Lambda^p(T_p^*U)
	\]
	such that for any smooth vector fields $X_1,\dots,X_p\in\mathfrak{X}(U)$, the function
	\[
	U\ni p\longmapsto \omega_p\bigl(X_1(p),\dots,X_p(p)\bigr)\in\mathbb{R}
	\]
	is smooth. The set of all smooth $p$-forms on $U$ is denoted $\Omega^p(U)$.
\end{definition}

\paragraph{Remarks (with concrete meaning).}
\begin{itemize}
	\item $\Omega^0(U)=C^\infty(U)$: a $0$-form is just a smooth function.
	\item $\Omega^1(U)$ consists of smooth assignments $p\mapsto \eta_p\in T_p^*U$;
	each $\eta_p$ is a linear functional on $T_pU$ depending smoothly on $p$.
	\item $\Omega^p(U)$ is a $C^\infty(U)$-module: if $f\in C^\infty(U)$ and $\omega\in\Omega^p(U)$,
	then $(f\omega)_p:=f(p)\,\omega_p$ defines $f\omega\in\Omega^p(U)$.
\end{itemize}

\paragraph{Worked mini-example 1 (a genuine $2$-form as an alternating map).}
Let $n=2$ and define at each $p$:
\[
\omega_p(v,w):=\det\begin{pmatrix} v^1 & w^1\\ v^2 & w^2\end{pmatrix}.
\]
Then $\omega_p$ is bilinear and alternating, hence $\omega_p\in\Lambda^2(T_p^*U)$.
If $v=(1,0)$ and $w=(0,1)$ in the basis $(\partial_x,\partial_y)$, then $\omega_p(v,w)=1$.

\paragraph{Worked mini-example 2 (a family that is \emph{not} smooth).}
Define $\omega_p$ in $\mathbb{R}^2$ by
\[
\omega_p(\partial_x,\partial_y)=
\begin{cases}
	1,& x(p)\ge 0,\\
	0,& x(p)<0.
\end{cases}
\]
Then each $\omega_p$ is an alternating $2$-covector, but the coefficient jumps at $x=0$.
Hence this assignment is \emph{not} a smooth $2$-form on $U$.

\paragraph{Alternation and signed area.}
The determinant calculation and oriented-area diagram below show why
interchanging the two inputs reverses the sign of a $2$-form.

\subsubsection{Exercises (definition level)}
\begin{exercise}
	Let $p\in U$ and $\omega_p\in\Lambda^p(T_p^*U)$. Prove directly from alternation that
	$\omega_p(v_1,\dots,v_p)=0$ whenever two of the $v_i$ are equal.
\end{exercise}

\begin{exercise}
	Show that $\Omega^p(U)$ is a real vector space and a $C^\infty(U)$-module, i.e.
	\[
	f(\omega+\eta)=f\omega+f\eta,\qquad (fg)\omega=f(g\omega),\qquad 1\cdot\omega=\omega.
	\]
\end{exercise}

\begin{exercise}
	Let $U\subset\mathbb{R}^2$ and define $\omega=f(x,y)\,dx\wedge dy$ (using the coordinate $1$-forms
	defined later). Show that $\omega$ corresponds to the alternating map
	\[
	\omega_p(v,w)=f(p)\det\begin{pmatrix}v^x&w^x\\ v^y&w^y\end{pmatrix}.
	\]
\end{exercise}

\subsubsection{Coordinates: $0$-forms and $1$-forms, coefficients, and evaluation}

\begin{definition}[0-form]
	A \emph{$0$-form} on $U$ is a smooth function $f\in C^\infty(U)$.
\end{definition}

\begin{definition}[Coordinate tangent basis]
	For each $p\in U$, the tangent space $T_pU$ has the coordinate basis
	$\{\partial_1|_p,\dots,\partial_n|_p\}$.
	Every $v\in T_pU$ can be written uniquely as
	\[
	v=\sum_{j=1}^n v^j\,\partial_j\Big|_p,\qquad v^j\in\mathbb{R}.
	\]
\end{definition}

\begin{definition}[1-forms and the dual basis]
	The coordinate $1$-forms $dx_i\in\Omega^1(U)$ are defined by
	\[
	(dx_i)_p(\partial_j|_p)=\delta_{ij}.
	\]
	A general $1$-form is a $C^\infty(U)$-linear combination
	\[
	\eta=\sum_{i=1}^n a_i\,dx_i,\qquad a_i\in C^\infty(U).
	\]
\end{definition}

\paragraph{Evaluation of a $1$-form on a vector.}
Let $\eta=\sum_i a_i dx_i$ and $v=\sum_j v^j\partial_j$ at a point $p$. Then
\[
\begin{aligned}
\eta_p(v)
&=\Bigl(\sum_{i=1}^n a_i(p)\,(dx_i)_p\Bigr)
  \Bigl(\sum_{j=1}^n v^j\,\partial_j|_p\Bigr)\\
&=\sum_{i,j} a_i(p)\,v^j\,(dx_i)_p(\partial_j|_p)
=\sum_{i,j} a_i(p)\,v^j\,\delta_{ij}\\
&=\sum_{i=1}^n a_i(p)\,v^i.
\end{aligned}
\]
So in coordinates, $\eta$ behaves like a row vector $(a_1,\dots,a_n)$ acting on a column vector
$(v^1,\dots,v^n)^{\mathsf T}$.

\paragraph{Worked example (numerical evaluation).}
In $\mathbb{R}^3$ let
\[
\eta=(x^2+y)\,dx+\sin z\,dy+e^{xy}\,dz.
\]
At $p=(1,2,0)$, the coefficient values are $(x^2+y,\sin z,e^{xy})=(3,0,e^2)$.
For $v=2\partial_x-\partial_y+3\partial_z$ at $p$, i.e. $(v^x,v^y,v^z)=(2,-1,3)$, we get
\[
\eta_p(v)=3\cdot 2+0\cdot(-1)+e^2\cdot 3=6+3e^2.
\]

\paragraph{Differential of a function (directional derivative packaged).}
For $f\in C^\infty(U)$ define
\[
df:=\sum_{i=1}^n \frac{\partial f}{\partial x_i}\,dx_i\in\Omega^1(U).
\]
Then for $v=\sum_j v^j\partial_j$,
\[
(df)_p(v)
=\sum_i f_{x_i}(p)\,(dx_i)_p\!\left(\sum_j v^j\partial_j|_p\right)
=\sum_{i,j} f_{x_i}(p)\,v^j\,\delta_{ij}
=\sum_i f_{x_i}(p)\,v^i
=v(f)(p).
\]

\paragraph{Worked example (compute $df$ and evaluate).}
Let $f(x,y)=x^2y+\cos y$ on $U\subset\mathbb{R}^2$. Then
\[
f_x=2xy,\qquad f_y=x^2-\sin y,
\]
so
\[
df=(2xy)\,dx+(x^2-\sin y)\,dy.
\]
At $p=(1,0)$ and $v=3\partial_x-2\partial_y$,
\[
(df)_p(v)=(2\cdot 1\cdot 0)\cdot 3+(1^2-\sin 0)\cdot(-2)=0-2=-2,
\]
which matches the directional derivative $v(f)(p)$.

\subsubsection{Exercises ($0$-forms and $1$-forms)}
\begin{exercise}
	Let $\eta=(x+y)\,dx+(x-y)\,dy$ in $\mathbb{R}^2$. Compute $\eta_p(v)$ at $p=(1,2)$ for
	$v=4\partial_x-\partial_y$. Show every step.
\end{exercise}

\begin{exercise}
	Let $f(x,y,z)=x e^{yz}+y^2z$. Compute $df$ explicitly and then evaluate $(df)_{(1,0,2)}(\partial_x+2\partial_z)$.
\end{exercise}

\begin{exercise}
	Show that $df=0$ as a $1$-form on a connected open set $U$ if and only if $f$ is constant on $U$.
\end{exercise}

\subsubsection{Wedge product and $k$-forms: signs, areas, and determinants}

\begin{definition}[Wedge product]
	For $\alpha\in\Omega^p(U)$ and $\beta\in\Omega^q(U)$, their wedge product
	is the $(p+q)$-form defined pointwise by
	\[
	(\alpha\wedge\beta)(v_1,\ldots,v_{p+q})
	=\frac{1}{p!q!}\sum_{\sigma\in S_{p+q}}
	\operatorname{sgn}(\sigma)\,
	\alpha(v_{\sigma(1)},\ldots,v_{\sigma(p)})
	\beta(v_{\sigma(p+1)},\ldots,v_{\sigma(p+q)}).
	\]
	It is bilinear and associative, and satisfies graded anti-commutativity:
	\[
	\alpha\wedge\beta = (-1)^{pq}\,\beta\wedge\alpha.
	\]
	In particular, for $1$-forms,
	\[
	dx_i\wedge dx_j=-dx_j\wedge dx_i,\qquad dx_i\wedge dx_i=0.
	\]
\end{definition}

\paragraph{From wedge to determinants (worked out in $\mathbb{R}^2$).}
A $2$-form acts on two vectors. In $\mathbb{R}^2$,
\[
(dx\wedge dy)_p(v,w)
=
\det\begin{pmatrix}
	dx(v) & dx(w)\\
	dy(v) & dy(w)
\end{pmatrix}
=
dx(v)\,dy(w)-dx(w)\,dy(v).
\]
If $v=(v^x,v^y)$ and $w=(w^x,w^y)$, then
\[
(dx\wedge dy)_p(v,w)=v^x w^y - v^y w^x,
\]
the signed area of the parallelogram spanned by $v,w$.

\paragraph{Worked example (sign check).}
Let $v=\partial_x+\partial_y$ and $w=2\partial_x+\partial_y$. Then
\[
(dx\wedge dy)(v,w)
=\det\begin{pmatrix}1&2\\1&1\end{pmatrix}=1\cdot1-2\cdot1=-1.
\]
Swapping vectors gives $(dx\wedge dy)(w,v)=+1$.

\paragraph{A $2$-form built from coefficients.}
Let $\omega = f(x,y)\,dx\wedge dy$ in $\mathbb{R}^2$.
Then for $v,w$,
\[
\omega_p(v,w)=f(p)\,(dx\wedge dy)_p(v,w)=f(p)\,(v^x w^y - v^y w^x).
\]
So $f$ rescales the signed area density.

\paragraph{Oriented area from $dx\wedge dy$.}
\begin{figure}[ht]
	\centering
	\begin{tikzpicture}[scale=1.0]
		\draw[->] (-0.2,0) -- (4.2,0) node[below] {$x$};
		\draw[->] (0,-0.2) -- (0,3.3) node[left] {$y$};
		
		\coordinate (O) at (1.0,0.9);
		\coordinate (V) at (3.2,1.4);
		\coordinate (W) at (1.7,2.8);
		\coordinate (VW) at ($(V)+(W)-(O)$);
		
		\fill (O) circle (2pt) node[below left] {$p$};
		
		\draw[very thick,->] (O) -- (V) node[midway,below] {$v$};
		\draw[very thick,->] (O) -- (W) node[midway,left] {$w$};
		
		\fill[gray!15] (O) -- (V) -- (VW) -- (W) -- cycle;
		\draw (O) -- (V) -- (VW) -- (W) -- cycle;
		
		\node at (2.55,2.15) {\small $(dx\wedge dy)(v,w)$};
		
		\draw[->] (1.95,1.25) arc[start angle=210,end angle=300,radius=0.55];
		\node at (2.05,0.62) {\small orientation};
	\end{tikzpicture}
	\caption{$dx\wedge dy$ measures signed area; swapping $(v,w)$ changes the sign.}
\end{figure}

\begin{definition}[$k$-forms in coordinates]
	A general $k$-form is
	\[
	\omega=\sum_{1\le i_1<\cdots<i_k\le n}
	a_{i_1\cdots i_k}\,dx_{i_1}\wedge\cdots\wedge dx_{i_k},
	\qquad a_{i_1\cdots i_k}\in C^\infty(U).
	\]
\end{definition}

\paragraph{Worked example (a $3$-form as signed volume).}
In $\mathbb{R}^3$, the $3$-form $dx\wedge dy\wedge dz$ satisfies
\[
(dx\wedge dy\wedge dz)_p(u,v,w)=\det\begin{pmatrix}
	dx(u)&dx(v)&dx(w)\\
	dy(u)&dy(v)&dy(w)\\
	dz(u)&dz(v)&dz(w)
\end{pmatrix},
\]
so it returns the signed volume of the parallelepiped spanned by $u,v,w$.

\subsubsection{Exercises (wedge product and determinants)}
\begin{exercise}
	Compute $(dx\wedge dy)(v,w)$ for $v=(2,1)$ and $w=(1,3)$, then compute $(dx\wedge dy)(w,v)$ and compare.
\end{exercise}

\begin{exercise}
	Let $\alpha = x\,dx+y\,dy$ and $\beta=dx+dy$ in $\mathbb{R}^2$. Compute $\alpha\wedge\beta$ explicitly,
	simplify the result, and verify $\beta\wedge\alpha=-\alpha\wedge\beta$.
\end{exercise}

\begin{exercise}
	In $\mathbb{R}^3$, compute $(dx\wedge dy\wedge dz)(u,v,w)$ for
	$u=(1,0,0)$, $v=(1,1,0)$, $w=(1,1,1)$ and interpret geometrically.
\end{exercise}

\subsubsection{Exterior derivative and the identity $d^2=0$}

\begin{definition}[Exterior derivative on functions]
	For $f\in C^\infty(U)=\Omega^0(U)$,
	\[
	df=\sum_{i=1}^n f_{x_i}\,dx_i\in\Omega^1(U).
	\]
\end{definition}

\begin{definition}[Exterior derivative on general $k$-forms (computational rule)]
	If
	\[
	\omega=\sum_{i_1<\cdots<i_k} a_{i_1\cdots i_k}\,dx_{i_1}\wedge\cdots\wedge dx_{i_k},
	\]
	then
	\[
	d\omega=\sum_{i_1<\cdots<i_k} da_{i_1\cdots i_k}\wedge dx_{i_1}\wedge\cdots\wedge dx_{i_k},
	\quad
	da_{i_1\cdots i_k}=\sum_{\ell=1}^n \frac{\partial a_{i_1\cdots i_k}}{\partial x_\ell}\,dx_\ell.
	\]
\end{definition}

\paragraph{Graded Leibniz rule (full sign bookkeeping).}
Let $\omega\in\Omega^k(U)$ and $\eta\in\Omega^\ell(U)$, and write locally
$\omega=\sum_I a_I\,dx_I$ and $\eta=\sum_J b_J\,dx_J$.
Then
\[
\omega\wedge\eta=\sum_{I,J} a_I b_J\,dx_I\wedge dx_J.
\]
Differentiate:
\[
d(\omega\wedge\eta)
=\sum_{I,J} d(a_I b_J)\wedge dx_I\wedge dx_J
=\sum_{I,J} (da_I\,b_J+a_I\,db_J)\wedge dx_I\wedge dx_J.
\]
Split:
\[
\sum_{I,J} (da_I\wedge dx_I)\,b_J\wedge dx_J
\;+\;
\sum_{I,J} a_I\,db_J\wedge dx_I\wedge dx_J.
\]
Move $db_J$ past $k$ one-forms in $dx_I$:
\[
db_J\wedge dx_I = (-1)^k\,dx_I\wedge db_J,
\]
so
\[
d(\omega\wedge\eta)=d\omega\wedge\eta+(-1)^k\,\omega\wedge d\eta.
\]

\begin{theorem}[$d^2=0$ (full coordinate cancellation)]
	For every $\omega\in\Omega^k(U)$,
	\[
	d(d\omega)=0.
	\]
\end{theorem}

\begin{proof}
	By linearity, it suffices to check $\omega=a\,dx_{i_1}\wedge\cdots\wedge dx_{i_k}$ with $i_1<\cdots<i_k$.
	Then
	\[
	d\omega=da\wedge dx_{i_1}\wedge\cdots\wedge dx_{i_k}
	=\left(\sum_{p=1}^n a_{x_p}\,dx_p\right)\wedge dx_{i_1}\wedge\cdots\wedge dx_{i_k}.
	\]
	Apply $d$ again:
	\[
	d(d\omega)
	=
	\sum_{p=1}^n d(a_{x_p})\wedge dx_p\wedge dx_{i_1}\wedge\cdots\wedge dx_{i_k}.
	\]
	Expand
	\[
	d(a_{x_p})=\sum_{q=1}^n a_{x_qx_p}\,dx_q,
	\]
	so
	\[
	d(d\omega)
	=
	\sum_{p,q=1}^n a_{x_qx_p}\,dx_q\wedge dx_p\wedge dx_{i_1}\wedge\cdots\wedge dx_{i_k}.
	\]
	For $p\neq q$, the term $(p,q)$ cancels $(q,p)$ because
	\[
	a_{x_qx_p}=a_{x_px_q}
	\quad\text{and}\quad
	dx_q\wedge dx_p=-dx_p\wedge dx_q.
	\]
	For $p=q$, the term vanishes since $dx_p\wedge dx_p=0$. Hence the total is $0$.
\end{proof}

\subsubsection{Worked examples}

\begin{example}[Compute $d\eta$ in $\mathbb{R}^2$ step by step]
	Let $\eta=a(x,y)\,dx+b(x,y)\,dy$.
	Compute
	\[
	da=a_x\,dx+a_y\,dy,\qquad db=b_x\,dx+b_y\,dy.
	\]
	Then
	\[
	d\eta=da\wedge dx+db\wedge dy.
	\]
	Expand:
	\[
	da\wedge dx=(a_x\,dx+a_y\,dy)\wedge dx
	=a_x\,dx\wedge dx+a_y\,dy\wedge dx
	=0-a_y\,dx\wedge dy,
	\]
	\[
	db\wedge dy=(b_x\,dx+b_y\,dy)\wedge dy
	=b_x\,dx\wedge dy+b_y\,dy\wedge dy
	=b_x\,dx\wedge dy+0.
	\]
	Therefore
	\[
	d\eta=(b_x-a_y)\,dx\wedge dy.
	\]
\end{example}

\begin{example}[Compute $d\alpha$ in $\mathbb{R}^3$ and compare with divergence]
	Let
	\[
	\alpha=x\,dy\wedge dz+y\,dz\wedge dx+z\,dx\wedge dy.
	\]
	Use $d(f\beta)=df\wedge\beta$ here because $d(dx)=d(dy)=d(dz)=0$:
	\[
	d(x\,dy\wedge dz)=dx\wedge dy\wedge dz,
	\]
	\[
	d(y\,dz\wedge dx)=dy\wedge dz\wedge dx,
	\qquad
	d(z\,dx\wedge dy)=dz\wedge dx\wedge dy.
	\]
	Reorder to $dx\wedge dy\wedge dz$ (cyclic permutations, hence sign $+1$):
	\[
	dy\wedge dz\wedge dx = dx\wedge dy\wedge dz,\qquad
	dz\wedge dx\wedge dy = dx\wedge dy\wedge dz.
	\]
	So
	\[
	d\alpha=(1+1+1)\,dx\wedge dy\wedge dz
	=3\,dx\wedge dy\wedge dz.
	\]
	This matches $\nabla\cdot(x,y,z)=3$ under the standard identification of vector fields with $2$-forms via
	interior product $\iota_F(dx\wedge dy\wedge dz)$.
\end{example}

\subsubsection{Exercises (exterior derivative)}
\begin{exercise}
	Let $f(x,y)=x^3+xy^2$. Compute $df$ and then show by explicit coefficient expansion that $d(df)=0$.
\end{exercise}

\begin{exercise}
	Let $\eta=P\,dx+Q\,dy$ in $\mathbb{R}^2$. Compute $d\eta$ and show that
	\[
	d\eta=(Q_x-P_y)\,dx\wedge dy.
	\]
\end{exercise}

\begin{exercise}
	In $\mathbb{R}^3$, let $\eta=P\,dx+Q\,dy+R\,dz$. Expand $d\eta$ completely and identify its coefficients with
	the curl $\nabla\times(P,Q,R)$ by matching the basis $dy\wedge dz,\ dz\wedge dx,\ dx\wedge dy$.
\end{exercise}

\begin{exercise}
	Let $\omega=f(x,y)\,dx\wedge dy$ in $\mathbb{R}^2$. Compute $d\omega$ and explain why $d\omega=0$ automatically
	in $\mathbb{R}^2$. (Hint: degree reasons.)
\end{exercise}

\subsection{Line Integrals and Conservative Vector Fields}
\label{subsec:line-integrals-conservative}

Throughout this subsection we work in $\R^2$ with standard coordinates $(x,y)$.
We keep two parallel languages in view:

\begin{itemize}
	\item vector fields $\mathbf F=(P,Q)$ and the line integral $\displaystyle \int_C \mathbf F\cdot d\mathbf r$;
	\item $1$-forms $\eta=P\,dx+Q\,dy$ and the curve integral $\displaystyle \int_C \eta$.
\end{itemize}

\subsubsection{Line integrals: definition, parametrization, and $1$-forms}

\begin{definition}[Line integral of a vector field]
	Let $\mathbf{F}=(P,Q)$ be a continuous vector field on a region $D\subset\R^2$.
	Let $C$ be an oriented $C^1$ curve parametrized by
	\[
	\mathbf r(t)=(x(t),y(t)),\qquad a\le t\le b.
	\]
	The \emph{line integral of $\mathbf F$ along $C$} is
	\[
	\int_C \mathbf F\cdot d\mathbf r
	:=\int_a^b\Big(P(x(t),y(t))\,x'(t)+Q(x(t),y(t))\,y'(t)\Big)\,dt.
	\]
\end{definition}

\paragraph{Coordinate procedure.}
To compute $\int_C\mathbf F\cdot d\mathbf r$:
\begin{enumerate}
	\item parametrize $C$ by $\mathbf r(t)=(x(t),y(t))$;
	\item compute $x'(t),y'(t)$;
	\item substitute into $P(x,y)$ and $Q(x,y)$ along $\mathbf r(t)$;
	\item form $P(x(t),y(t))x'(t)+Q(x(t),y(t))y'(t)$ and integrate in $t$.
\end{enumerate}

\paragraph{Rewriting using differential forms.}
Associate to $\mathbf F=(P,Q)$ the $1$-form
\[
\eta:=P\,dx+Q\,dy\in\Omega^1(D).
\]
Along the curve $\mathbf r(t)$ we have
\[
\eta_{\mathbf r(t)}(\mathbf r'(t))
= P(x(t),y(t))\,x'(t)+Q(x(t),y(t))\,y'(t),
\]
hence
\[
\int_C\mathbf F\cdot d\mathbf r=\int_a^b \eta(\mathbf r'(t))\,dt=\int_C\eta.
\]
So \emph{line integrals of vector fields are exactly integrals of $1$-forms}.

\paragraph{Orientation and reparametrization (checked by calculation).}
If $C$ is traversed in the opposite direction, say $\tilde{\mathbf r}(t)=\mathbf r(a+b-t)$, then
\[
\int_{\tilde C}\mathbf F\cdot d\mathbf r
=\int_a^b\Big(P(\tilde{\mathbf r}(t))\,\tilde x'(t)+Q(\tilde{\mathbf r}(t))\,\tilde y'(t)\Big)\,dt
=-\int_C\mathbf F\cdot d\mathbf r.
\]
If $t=\varphi(s)$ is an orientation-preserving $C^1$ reparametrization
($\varphi'>0$, $\varphi(\alpha)=a$, $\varphi(\beta)=b$), then by change of variables
\[
\int_\alpha^\beta \eta(\mathbf r'(\varphi(s))\,\varphi'(s))\,ds
=\int_a^b \eta(\mathbf r'(t))\,dt,
\]
so the value is unchanged.

\subsubsection{Conservative vector fields, potentials, and the Fundamental Theorem}

\begin{definition}[Conservative vector field]
	A vector field $\mathbf F=(P,Q)$ on a region $D\subset\R^2$ is \emph{conservative}
	if there exists $f\in C^1(D)$ such that
	\[
	\mathbf F=\nabla f=(f_x,f_y).
	\]
	The function $f$ is called a \emph{potential} for $\mathbf F$.
\end{definition}

\paragraph{Fundamental Theorem for Line Integrals (proved by chain rule).}
Assume $\mathbf F=\nabla f$ and $C:\mathbf r(t)=(x(t),y(t))$ runs from $A=\mathbf r(a)$ to $B=\mathbf r(b)$.
Then
\[
\mathbf F(\mathbf r(t))\cdot \mathbf r'(t)
= f_x(x(t),y(t))\,x'(t)+f_y(x(t),y(t))\,y'(t)
= \frac{d}{dt}\,f(x(t),y(t))
\]
by the chain rule. Integrating,
\[
\int_C \mathbf F\cdot d\mathbf r
=\int_a^b \frac{d}{dt}f(\mathbf r(t))\,dt
=f(\mathbf r(b))-f(\mathbf r(a))=f(B)-f(A).
\]
In particular, the integral depends only on endpoints.

\paragraph{Differential-form interpretation.}
If $\mathbf F=\nabla f$, then
\[
\eta=Pdx+Qdy=f_xdx+f_y dy = df,
\]
so $\eta$ is \emph{exact}. Applying $d$ gives
\[
d\eta=d(df)=0,
\]
so every exact form is closed. The converse can fail if the domain is not simply connected.

\paragraph{A standard criterion in $\R^2$ (local test).}
If $P,Q\in C^1(D)$, then
\[
d\eta=d(Pdx+Qdy)=(Q_x-P_y)\,dx\wedge dy.
\]
So $d\eta=0$ is equivalent to the familiar curl-free condition $Q_x=P_y$.

\subsubsection{A key example: a tangential field on $\R^2\setminus\{0\}$}

\paragraph{Definition of the field and geometric meaning.}
On $\R^2\setminus\{(0,0)\}$ define
\[
\mathbf F(x,y)=\left(-\frac{y}{x^2+y^2},\;\frac{x}{x^2+y^2}\right).
\]
Write $r^2=x^2+y^2$. Then $\mathbf F=\frac{1}{r^2}(-y,x)$ is obtained by rotating $(x,y)$ by $+90^\circ$
and scaling by $1/r^2$, so it is tangent to circles centered at the origin.

\paragraph{Associated $1$-form.}
\[
\eta
=-\frac{y}{x^2+y^2}\,dx+\frac{x}{x^2+y^2}\,dy.
\]

\paragraph{Computation of $d\eta$.}
Let $P=-\dfrac{y}{x^2+y^2}$ and $Q=\dfrac{x}{x^2+y^2}$. Compute
\[
Q_x=\frac{\partial}{\partial x}\left(\frac{x}{x^2+y^2}\right)
=\frac{(x^2+y^2)\cdot 1 - x\cdot (2x)}{(x^2+y^2)^2}
=\frac{y^2-x^2}{(x^2+y^2)^2},
\]
and
\[
P_y=\frac{\partial}{\partial y}\left(-\frac{y}{x^2+y^2}\right)
=-\frac{(x^2+y^2)\cdot 1 - y\cdot (2y)}{(x^2+y^2)^2}
=-\frac{x^2-y^2}{(x^2+y^2)^2}
=\frac{y^2-x^2}{(x^2+y^2)^2}.
\]
Thus $Q_x-P_y=0$, and therefore
\[
d\eta=(Q_x-P_y)\,dx\wedge dy=0
\qquad\text{on }\R^2\setminus\{0\}.
\]
So $\eta$ is \emph{closed} on $\R^2\setminus\{0\}$.

\paragraph{But $\eta$ is not exact globally (detected by a loop integral).}
Take the circle $C_R:\ (x,y)=(R\cos t,R\sin t)$, $0\le t\le 2\pi$ (counterclockwise).
Then $dx=-R\sin t\,dt$, $dy=R\cos t\,dt$, and $x^2+y^2=R^2$.
Hence
\[
\eta
=-\frac{R\sin t}{R^2}(-R\sin t\,dt)+\frac{R\cos t}{R^2}(R\cos t\,dt)
=\sin^2 t\,dt+\cos^2 t\,dt
=dt.
\]
Therefore
\[
\oint_{C_R}\eta=\int_0^{2\pi} dt=2\pi\neq 0.
\]
If $\eta=df$ on $\R^2\setminus\{0\}$, then every closed loop would integrate to $0$,
contradiction. So $\eta$ is closed but not exact on $\R^2\setminus\{0\}$.

\paragraph{The tangential field and two circles.}
\begin{figure}[ht]
	\centering
	\begin{tikzpicture}[scale=1.05]
		\draw[->] (-3.0,0) -- (3.2,0) node[below] {$x$};
		\draw[->] (0,-2.6) -- (0,2.8) node[left] {$y$};
		
		\fill (0,0) circle (2pt) node[below left] {$0$};
		
		\draw[thick] (0,0) circle (1.6);
		\node at (-1.9,-1.65) {\small $C_R$};
		
		\draw[thick,dashed] (2,0) circle (1.0);
		\fill (2,0) circle (1.5pt) node[below] {\small $(2,0)$};
		
		\foreach \ang in {0,45,90,135,180,225,270,315}{
			\pgfmathsetmacro\x{1.6*cos(\ang)}
			\pgfmathsetmacro\y{1.6*sin(\ang)}
			\pgfmathsetmacro\ux{-sin(\ang)}
			\pgfmathsetmacro\uy{cos(\ang)}
			\draw[->] (\x,\y) -- ++(0.45*\ux,0.45*\uy);
		}
		\node at (-2.1,2.2) {\small tangential field};
		
		\node at (-2.1,-2.2) {\small $\oint_{C_R}\eta=2\pi$};
		\node at (2.15,-1.35) {\small $C_0$};
		\node at (2.15,1.65) {\small $\oint_{C_0}\eta=0$};
	\end{tikzpicture}
	\caption{The same closed $1$-form has period $2\pi$ around the origin
	but period zero around a displaced circle that does not enclose it.}
\end{figure}

\subsubsection{Line integral along a displaced circle}

\begin{example}[A displaced circle that does not enclose the origin]
	Let $C$ be the circle of radius $1$ centered at $(2,0)$, oriented counterclockwise.
	Compute $\displaystyle \oint_C \mathbf F\cdot d\mathbf r$ for
	\[
	\mathbf F(x,y)=\left(-\frac{y}{x^2+y^2},\;\frac{x}{x^2+y^2}\right).
	\]
\end{example}

\begin{proof}[Coordinate computation]
	Parametrize $C$ by
	\[
	x(t)=2+\cos t,\qquad y(t)=\sin t,\qquad 0\le t\le 2\pi.
	\]
	Then
	\[
	x'(t)=-\sin t,\qquad y'(t)=\cos t,
	\]
	and
	\[
	x(t)^2+y(t)^2=(2+\cos t)^2+\sin^2 t=5+4\cos t.
	\]
	Evaluate $\mathbf F$ along $C$:
	\[
	\mathbf F(\mathbf r(t))
	=\left(-\frac{\sin t}{5+4\cos t},\;\frac{2+\cos t}{5+4\cos t}\right).
	\]
	Dot with $\mathbf r'(t)$:
	\begin{align*}
		\mathbf F(\mathbf r(t))\cdot \mathbf r'(t)
		&=
		\left(-\frac{\sin t}{5+4\cos t}\right)(-\sin t)
		+
		\left(\frac{2+\cos t}{5+4\cos t}\right)(\cos t)\\
		&=\frac{\sin^2 t+(2+\cos t)\cos t}{5+4\cos t}
		=\frac{1+2\cos t}{5+4\cos t}.
	\end{align*}
	Therefore
	\[
	\oint_C \mathbf F\cdot d\mathbf r
	=\int_0^{2\pi}\frac{1+2\cos t}{5+4\cos t}\,dt.
	\]
\end{proof}

\subsubsection{Direct evaluation via the tangent-half-angle substitution}

Let $u=\tan(t/2)$. Then
\[
\cos t=\frac{1-u^2}{1+u^2},
\qquad
dt=\frac{2\,du}{1+u^2}.
\]
Compute the rationalized integrand carefully:
\[
1+2\cos t
=1+2\frac{1-u^2}{1+u^2}
=\frac{(1+u^2)+2(1-u^2)}{1+u^2}
=\frac{3-u^2}{1+u^2},
\]
\[
5+4\cos t
=5+4\frac{1-u^2}{1+u^2}
=\frac{5(1+u^2)+4(1-u^2)}{1+u^2}
=\frac{9+u^2}{1+u^2}.
\]
Hence
\[
\frac{1+2\cos t}{5+4\cos t}
=\frac{\frac{3-u^2}{1+u^2}}{\frac{9+u^2}{1+u^2}}
=\frac{3-u^2}{9+u^2}.
\]
Therefore
\[
\int_0^{2\pi}\frac{1+2\cos t}{5+4\cos t}\,dt
=
\int_{-\infty}^{\infty}\frac{3-u^2}{9+u^2}\cdot \frac{2\,du}{1+u^2}
=
2\int_{-\infty}^{\infty}\frac{3-u^2}{(9+u^2)(1+u^2)}\,du.
\]
Now do the partial fraction decomposition (this is the key cancellation):
\[
2\frac{3-u^2}{(u^2+9)(u^2+1)}
=\frac{1}{u^2+1}-\frac{3}{u^2+9}.
\]
(Verify by combining over the common denominator.)
Thus
\[
\oint_C \mathbf F\cdot d\mathbf r
=
\int_{-\infty}^{\infty}\left(\frac{1}{u^2+1}-\frac{3}{u^2+9}\right)\,du
=\Big(\pi\Big)-\Big(3\cdot \frac{\pi}{3}\Big)=0.
\]
So the displaced circle gives zero circulation, consistent with the fact that it does \emph{not}
wind around the origin.

\subsubsection{Local and global information}

\begin{itemize}
	\item Conservative fields $\mathbf F=\nabla f$ correspond to exact $1$-forms $\eta=df$.
	\item Every exact form is closed ($d\eta=0$), but a closed form need not be exact.
	\item In $\R^2$, the obstruction is topological: loops winding around ``holes'' can have nonzero integrals.
	\item Line integrals can detect global information invisible to local derivatives like $Q_x-P_y$.
\end{itemize}

\subsubsection{Exercises}

\begin{exercise}[Warm-up: compute a line integral directly]
	Let $\mathbf F=(P,Q)=(y,x)$ and let $C$ be the parabola $y=x^2$ from $(0,0)$ to $(1,1)$
	with the standard parametrization $\mathbf r(t)=(t,t^2)$, $0\le t\le 1$.
	Compute $\displaystyle \int_C \mathbf F\cdot d\mathbf r$.
\end{exercise}

\begin{exercise}[Finding a potential by matching partial derivatives]
	Let $\mathbf{F}(x,y)=(3x^2+6xy,\;3x^2+6y)$ on $\R^2$.
	\begin{enumerate}
		\item Compute $P_y$ and $Q_x$ and verify $P_y=Q_x$.
		\item Find a potential $f$ by integrating $f_x=P$ in $x$ and determining the missing $y$-function.
		\item For any $C^1$ curve $C$ from $(0,0)$ to $(1,1)$, compute $\displaystyle \int_C \mathbf F\cdot d\mathbf r$.
	\end{enumerate}
\end{exercise}

\begin{exercise}[Path dependence via two different curves]
	Let $\mathbf F=(P,Q)=(y,0)$ on $\R^2$.
	Compute $\int_C \mathbf F\cdot d\mathbf r$ from $(0,0)$ to $(1,1)$ along:
	\begin{enumerate}
		\item $C_1$: the straight line $\mathbf r(t)=(t,t)$;
		\item $C_2$: the broken path $(0,0)\to(1,0)\to(1,1)$.
	\end{enumerate}
	Conclude that $\mathbf F$ is not conservative.
\end{exercise}

\begin{exercise}[Every exact form is closed]
	Let $\eta=df$ for $f\in C^2(D)$. Write $df=f_xdx+f_ydy$ and compute
	$d(df)$ explicitly to show $d(df)=0$ using $f_{xy}=f_{yx}$.
\end{exercise}

\begin{exercise}[Closed $\mathrel{\text{need not give}}$ exact on a punctured plane]
	Let $\eta=-\dfrac{y}{x^2+y^2}dx+\dfrac{x}{x^2+y^2}dy$ on $\R^2\setminus\{0\}$.
	\begin{enumerate}
		\item Compute $d\eta$ explicitly and show $d\eta=0$.
		\item Parameterize $C_R:(x,y)=(R\cos t,R\sin t)$ and compute $\displaystyle \oint_{C_R}\eta$.
		\item Explain why this proves $\eta$ is not exact on $\R^2\setminus\{0\}$.
	\end{enumerate}
\end{exercise}

\begin{exercise}[Winding number viewpoint (geometric)]
	Let $C$ be any positively oriented simple closed curve in $\R^2\setminus\{0\}$.
	First verify for circles and ellipses centered at the origin, and then for the
	piecewise smooth boundary of the square $[-1,1]^2$, that
	\[
	\oint_C \left(-\frac{y}{x^2+y^2}dx+\frac{x}{x^2+y^2}dy\right)=2\pi\cdot \mathrm{wind}(C,0),
	\]
	where $\mathrm{wind}(C,0)\in\Z$ is the winding number of $C$ around the origin.
	Use these calculations to formulate the general identity, and prove it by
	homotopy invariance or by the winding-number formula from complex analysis.
\end{exercise}

\begin{exercise}[A simply connected criterion you can actually use]
	Let $D\subset\R^2$ be simply connected and let $\mathbf F=(P,Q)$ with $P,Q\in C^1(D)$.
	Assume $P_y=Q_x$ on $D$.
	\begin{enumerate}
		\item Fix a base point $(x_0,y_0)\in D$. Define
		\[
		f(x,y):=\int_{\gamma_{(x,y)}} P\,dx+Q\,dy
		\]
	where $\gamma_{(x,y)}$ is any $C^1$ path in $D$ from $(x_0,y_0)$ to $(x,y)$.
		First prove path independence using homotopy invariance (or Green's theorem after a suitable
		subdivision).  Then append short horizontal and vertical segments contained in a coordinate
		rectangle about $(x,y)$ to show that $f_x=P$ and $f_y=Q$.
		\item Conclude that $\mathbf F=\nabla f$ and $\eta=Pdx+Qdy=df$ on $D$.
	\end{enumerate}
\end{exercise}

\begin{exercise}[Challenge: re-check the displaced circle integral with a different method]
	For the displaced circle $C:(x,y)=(2+\cos t,\sin t)$, evaluate
	\[
	\int_0^{2\pi}\frac{1+2\cos t}{5+4\cos t}\,dt
	\]
	without the tangent-half-angle substitution by writing
	\[
	\frac{1+2\cos t}{5+4\cos t}=A+\frac{B}{5+4\cos t}
	\]
	and determining $A,B$.
	(Then use the standard integral $\int_0^{2\pi}\frac{dt}{a+b\cos t}=\frac{2\pi}{\sqrt{a^2-b^2}}$ for $a>|b|$.)
\end{exercise}

\subsection{Complex Reformulation of Line Integrals}
\label{subsec:complex-line-integrals}

Throughout this subsection, we reinterpret planar line integrals using complex
coordinates. This viewpoint unifies vector calculus, differential forms, and
complex analysis, and explains circulation integrals in terms of winding and residues.

\subsubsection{Complex coordinates and differentials}

Let
\[
z=x+iy,\qquad \bar z=x-iy.
\]
Differentiating,
\[
dz=dx+i\,dy,\qquad d\bar z=dx-i\,dy.
\]
Solving for $dx,dy$,
\[
dx=\tfrac12(dz+d\bar z),\qquad
dy=\tfrac{1}{2i}(dz-d\bar z).
\]

\paragraph{Interpretation and a small warning about notation.}
The forms $dz$ and $d\bar z$ are complex-valued and form a basis of the
\emph{complexified} real cotangent space
$T^*\R^2\otimes_\R\C$. A real $1$-form
\[
\eta=P\,dx+Q\,dy
\]
can be rewritten uniquely in the $\{dz,d\bar z\}$ basis, and this rewrite is just linear algebra.
(We are not yet assuming holomorphicity; we are only changing coordinates.)

Because these are covectors rather than arrows in the physical $xy$-plane,
their change of basis is better recorded algebraically:
\[
\begin{pmatrix}dz\\d\bar z\end{pmatrix}
=\begin{pmatrix}1&i\\1&-i\end{pmatrix}
\begin{pmatrix}dx\\dy\end{pmatrix}.
\]

\subsubsection{From vector fields to complex-valued $1$-forms}

Let $\mathbf{F}=(P,Q)$ be a (real) vector field on $\R^2$ and let
\[
\eta:=P\,dx+Q\,dy
\]
be the associated real $1$-form.

\paragraph{Coefficients in the $\{dz,d\bar z\}$ basis.}
Substitute $dx=\tfrac12(dz+d\bar z)$ and $dy=\tfrac{1}{2i}(dz-d\bar z)$:
\begin{align*}
	\eta
	&=P\cdot \tfrac12(dz+d\bar z)+Q\cdot \tfrac{1}{2i}(dz-d\bar z)\\
	&=\tfrac12P\,dz+\tfrac12P\,d\bar z+\tfrac{1}{2i}Q\,dz-\tfrac{1}{2i}Q\,d\bar z\\
	&=\tfrac12\bigl(P-iQ\bigr)\,dz+\tfrac12\bigl(P+iQ\bigr)\,d\bar z.
\end{align*}

\paragraph{What the coefficients mean.}
Define the complex-valued functions
\[
F(z,\bar z):=P-iQ,\qquad \overline{F}(z,\bar z):=P+iQ.
\]
Then
\[
\eta=\tfrac12 F\,dz+\tfrac12 \overline{F}\,d\bar z.
\]
So every real $1$-form (equivalently every planar vector field) splits into a $(1,0)$-part and a $(0,1)$-part.

\paragraph{Complex potential viewpoint (when it happens).}
If $F$ depends only on $z$ (not on $\bar z$), and if additionally $\tfrac12F\,dz$ has a primitive
(i.e.\ $F=\phi'(z)$ locally), then
\[
\tfrac12F\,dz=d\Big(\tfrac12\phi(z)\Big)
\]
is a holomorphic $1$-form.  Taking real/imaginary parts produces real closed forms and the
classical harmonic conjugate picture (Cauchy--Riemann).

\subsubsection{Canonical example: the rotational field and the form $\operatorname{Im}(dz/z)$}

\paragraph{Definition of the field.}
On $\R^2\setminus\{0\}$ consider
\[
\mathbf{F}(x,y)=\left(-\frac{y}{x^2+y^2},\;\frac{x}{x^2+y^2}\right),
\qquad r^2=x^2+y^2.
\]
The associated $1$-form is
\[
\eta=-\frac{y}{r^2}\,dx+\frac{x}{r^2}\,dy.
\]

\paragraph{Step 1: rewrite $\eta$ in complex notation.}
Start from
\[
z=x+iy,\quad \bar z=x-iy,\quad |z|^2=z\bar z=x^2+y^2=r^2.
\]
Compute
\[
\bar z\,dz=(x-iy)(dx+i\,dy)=x\,dx+i x\,dy-i y\,dx+y\,dy,
\]
\[
z\,d\bar z=(x+iy)(dx-i\,dy)=x\,dx-i x\,dy+i y\,dx+y\,dy.
\]
Subtract:
\begin{align*}
	\bar z\,dz-z\,d\bar z
	&=\Big(x\,dx+i x\,dy-i y\,dx+y\,dy\Big)-\Big(x\,dx-i x\,dy+i y\,dx+y\,dy\Big)\\
	&= (i x\,dy-(-i x\,dy)) + (-i y\,dx-(i y\,dx))\\
	&=2i x\,dy-2i y\,dx\\
	&=2i\,(x\,dy-y\,dx).
\end{align*}
Hence
\[
x\,dy-y\,dx=\frac{1}{2i}\,(\bar z\,dz-z\,d\bar z).
\]
Divide by $|z|^2=r^2$:
\[
\eta=\frac{x\,dy-y\,dx}{x^2+y^2}
=\frac{1}{2i|z|^2}\,\bigl(\bar z\,dz-z\,d\bar z\bigr).
\]

\paragraph{Step 2: restrict to the unit circle and identify $\eta$ with $\operatorname{Im}(dz/z)$.}
On $|z|=1$ we have $|z|^2=1$ and $\bar z=1/z$. Thus
\[
\eta=\frac{1}{2i}\Big(\bar z\,dz-z\,d\bar z\Big)
=\frac{1}{2i}\left(\frac{dz}{z}-z\,d\bar z\right).
\]
Now compute $\operatorname{Im}(dz/z)$ using the standard identity
\[
\operatorname{Im}(w)=\frac{1}{2i}\,(w-\bar w).
\]
For $w=\dfrac{dz}{z}$, along $|z|=1$ we have $\overline{w}=\overline{\dfrac{dz}{z}}=\dfrac{d\bar z}{\bar z}=z\,d\bar z$.
Therefore
\[
\operatorname{Im}\!\left(\frac{dz}{z}\right)
=\frac{1}{2i}\left(\frac{dz}{z}-z\,d\bar z\right)
=\eta
\qquad\text{on }|z|=1.
\]

\paragraph{Step 3: two computations of the circulation integral.}

\emph{Method A (direct parametrization).}
Let $C: z=e^{it}$, $0\le t\le 2\pi$. Then $dz=i e^{it}\,dt=i z\,dt$, hence
\[
\oint_C \frac{dz}{z}=\int_0^{2\pi} i\,dt=2\pi i,
\qquad
\oint_C \eta=\oint_C \operatorname{Im}\!\left(\frac{dz}{z}\right)
=\operatorname{Im}(2\pi i)=2\pi.
\]

\emph{Method B (residue).}
Since $\frac{1}{z}$ has residue $1$ at $0$,
\[
\oint_C\frac{dz}{z}=2\pi i\cdot \mathrm{Res}_{z=0}\left(\frac{1}{z}\right)=2\pi i,
\quad\text{so}\quad
\oint_C \eta=2\pi.
\]

\begin{figure}[ht]
	\centering
	\begin{tikzpicture}[scale=1.2]
		\draw[->] (-2.2,0) -- (2.2,0) node[below] {$x$};
		\draw[->] (0,-2.2) -- (0,2.2) node[left] {$y$};
		
		\draw[thick] (0,0) circle (1.5);
		
		\foreach \ang in {0,45,90,135,180,225,270,315}{
			\pgfmathsetmacro\x{1.5*cos(\ang)}
			\pgfmathsetmacro\y{1.5*sin(\ang)}
			\pgfmathsetmacro\ux{-sin(\ang)}
			\pgfmathsetmacro\uy{cos(\ang)}
			\draw[->] (\x,\y) -- ++(0.45*\ux,0.45*\uy);
		}
		
		\node at (1.95,1.6) {\small tangential};
		\node at (-0.2,-1.95) {$|z|=1$};
		\node at (-1.9,1.65) {\small $\displaystyle \eta=\Im(dz/z)$ on $|z|=1$};
	\end{tikzpicture}
	\caption{Complex notation for the rotational field: $\eta=\operatorname{Im}(dz/z)$.
	The tangential arrows indicate positive circulation, while the logarithmic
	period $\oint dz/z$ records the same winding with a factor of $i$.}
\end{figure}

\subsubsection{Small circles, winding numbers, and residues}

Let $p\in\C$ and consider the positively oriented circle
\[
C:\ z(t)=p+re^{it},\qquad 0\le t\le 2\pi,\qquad r>0.
\]
Then
\[
dz=i\,re^{it}\,dt,\qquad z-p=re^{it},
\]
so
\[
\frac{dz}{z-p}=\frac{i\,re^{it}\,dt}{re^{it}}=i\,dt.
\]
Therefore
\[
\oint_C \frac{dz}{z-p}=\int_0^{2\pi} i\,dt=2\pi i.
\]

\paragraph{A slightly more general computation: $n$ turns around $p$.}
If instead $C_n: z(t)=p+re^{int}$, $0\le t\le 2\pi$ with $n\in\Z$, then
\[
dz=i n r e^{int}\,dt,\qquad z-p=r e^{int},
\]
hence
\[
\oint_{C_n}\frac{dz}{z-p}=\int_0^{2\pi} i n\,dt=2\pi i\,n.
\]
This integer $n$ is the \emph{winding number} $\mathrm{wind}(C,p)$, and the residue theorem packages it as
\[
\oint_C \frac{dz}{z-p}=2\pi i\cdot \mathrm{wind}(C,p).
\]

\paragraph{Back to real circulation: how to manufacture a vortex at $p$.}
Define the $1$-form on $\C\setminus\{p\}$
\[
\eta_p:=\operatorname{Im}\!\left(\frac{dz}{z-p}\right).
\]
Then for any loop $C$ avoiding $p$,
\[
\oint_C \eta_p=\operatorname{Im}\!\left(\oint_C \frac{dz}{z-p}\right)
=2\pi\cdot \mathrm{wind}(C,p).
\]
So $\eta_p$ is closed on $\C\setminus\{p\}$ but not exact there, exactly like the $p=0$ case.

\begin{remark}[Curl theorem vs residue theorem: same shadow, different language]
	Vector calculus (Green's theorem) says for sufficiently nice $\mathbf F$,
	\[
	\oint_C \mathbf F\cdot d\mathbf r
	=\iint_{\mathrm{Int}(C)}\left(\frac{\partial Q}{\partial x}-\frac{\partial P}{\partial y}\right)\,dA.
	\]
	Complex analysis says for meromorphic $g$,
	\[
	\oint_C g(z)\,dz=2\pi i\sum \mathrm{Res}(g;\text{poles inside }C).
	\]
	In vortex examples, the ``curl'' concentrates at singularities; the residue measures exactly that concentration through winding.
\end{remark}

\subsubsection{Summary of the real and complex descriptions}

\begin{itemize}
	\item A planar vector field $\mathbf F=(P,Q)$ corresponds to the real $1$-form $\eta=Pdx+Qdy$.
	\item In complex coordinates,
	\[
	\eta=\tfrac12(P-iQ)\,dz+\tfrac12(P+iQ)\,d\bar z.
	\]
	\item Vortex-type circulation is encoded by $\displaystyle \oint \frac{dz}{z-p}$ and its imaginary part.
	\item The integer $\mathrm{wind}(C,p)$ is the topological quantity behind both circulation and residues.
\end{itemize}

\subsubsection{Exercises}

\begin{exercise}[Cauchy theorem when $a$ is outside]
	Let $C$ be a positively oriented circle that does \emph{not} enclose $a\in\C$.
	Compute $\displaystyle \oint_C \frac{dz}{z-a}$ and explain the result using:
	\begin{enumerate}
		\item Cauchy's theorem (holomorphic on and inside $C$),
		\item the real $1$-form $\eta_a=\operatorname{Im}(dz/(z-a))$ and winding number.
	\end{enumerate}
\end{exercise}

\begin{exercise}[Unit circle integrals of monomials, all integer powers]
	For $n\in\Z$, compute
	\[
	\oint_{|z|=1} z^n\,dz.
	\]
	(Hint: use the parametrization $z=e^{it}$ or use residues.)
	Interpret the nonzero case in terms of winding/residue.
\end{exercise}

\begin{exercise}[Convert a real field to $dz,d\bar z$ and test conservativity]
	Let $\mathbf F=(x,-y)$ on $\R^2$.
	\begin{enumerate}
		\item Write the associated $1$-form $\eta=x\,dx-y\,dy$.
		\item Express $\eta$ in terms of $dz,d\bar z$.
		\item Decide whether $\mathbf F$ is conservative by finding an explicit potential $f$ or by a curl test.
	\end{enumerate}
\end{exercise}

\begin{exercise}[Closedness and loop integrals for $\eta=\Im(dz/z)$]
	Let $\eta=\operatorname{Im}(dz/z)$ on $\C\setminus\{0\}$.
	\begin{enumerate}
		\item Show $d\eta=0$ on $\C\setminus\{0\}$ by writing $\eta=\dfrac{x\,dy-y\,dx}{x^2+y^2}$ and computing $d\eta$.
		\item Compute $\displaystyle \oint_{|z|=R}\eta$ for arbitrary $R>0$.
		\item Conclude that $\eta$ is not exact on $\C\setminus\{0\}$.
	\end{enumerate}
\end{exercise}

\begin{exercise}[Winding number $n$ explicitly]
	Let $C_n$ be the loop $z(t)=re^{int}$, $0\le t\le 2\pi$, with $n\in\Z$ and $r>0$.
	Compute:
	\[
	\oint_{C_n}\frac{dz}{z},\qquad
	\oint_{C_n}\Im\!\left(\frac{dz}{z}\right).
	\]
	Explain why the answers must be proportional to $n$.
\end{exercise}

\begin{exercise}[Challenge: recover $\eta_p$ from a real vector field]
	Fix $p=a+ib\in\C$. Define
	\[
	\eta_p:=\Im\!\left(\frac{dz}{z-p}\right).
	\]
	\begin{enumerate}
		\item Write $\eta_p$ as a real $1$-form $P\,dx+Q\,dy$ by simplifying
		\[
		\frac{1}{z-p}=\frac{\overline{(z-p)}}{|z-p|^2}=\frac{(x-a)-i(y-b)}{(x-a)^2+(y-b)^2}.
		\]
		\item Identify the corresponding vector field $\mathbf F_p=(P,Q)$ and describe it geometrically.
		\item Compute $\oint_C \eta_p$ for $C: z=p+Re^{it}$ and verify it equals $2\pi$.
	\end{enumerate}
\end{exercise}

\subsection{Stokes' Theorem, Cauchy--Green, and Surface Area as a 2-Form}
\index{Stokes' theorem}\index{Cauchy--Green formula}\index{area form}
\label{subsec:stokes-cg-area2form}

Throughout this subsection we unify three themes:
(i) Stokes/Green in the plane as $\int_{\partial\Omega}\eta=\int_\Omega d\eta$,
(ii) the Cauchy--Green formula as a direct application to a carefully chosen $1$-form,
and (iii) surface area as the integral of a pullback $2$-form (hence invariant under reparametrization).
The calculations below keep track of the wedge products and their signs.

\subsubsection{Stokes' theorem in the plane (Green's theorem) in differential-form language}

Let $\Omega\subset\R^2$ be a bounded domain with positively oriented (counterclockwise)
$C^1$ boundary $\partial\Omega$.
Let $\eta\in\Omega^1(\Omega)$ be a smooth $1$-form.  In standard coordinates $(x,y)$,
write
\[
\eta=P(x,y)\,dx+Q(x,y)\,dy.
\]

\paragraph{Step 1: compute $d\eta$.}
Using $d(P\,dx)=dP\wedge dx$ and $d(Q\,dy)=dQ\wedge dy$ (because $d(dx)=d(dy)=0$),
\[
d\eta=dP\wedge dx+dQ\wedge dy.
\]
Now expand
\[
dP=P_x\,dx+P_y\,dy,\qquad dQ=Q_x\,dx+Q_y\,dy,
\]
so
\[
dP\wedge dx=(P_x\,dx+P_y\,dy)\wedge dx
=P_x\,dx\wedge dx+P_y\,dy\wedge dx
=0-P_y\,dx\wedge dy,
\]
\[
dQ\wedge dy=(Q_x\,dx+Q_y\,dy)\wedge dy
=Q_x\,dx\wedge dy+Q_y\,dy\wedge dy
=Q_x\,dx\wedge dy+0.
\]
Hence
\[
d\eta=(Q_x-P_y)\,dx\wedge dy.
\]

\paragraph{Step 2: state Stokes/Green.}
Stokes' theorem in the plane says
\[
\int_{\partial\Omega}\eta=\int_{\Omega} d\eta.
\]
With the computed formula, this becomes the classical Green theorem:
\[
\int_{\partial\Omega} \big(P\,dx+Q\,dy\big)
=\iint_{\Omega}(Q_x-P_y)\,dx\,dy.
\]

\paragraph{A direct check.}
Let $\eta=x\,dy$ on $\Omega$ (any region). Then $P=0$, $Q=x$.
So $Q_x-P_y=1-0=1$, hence $d\eta=dx\wedge dy$.
Therefore
\[
\int_{\partial\Omega} x\,dy=\iint_\Omega 1\,dx\,dy=\mathrm{Area}(\Omega).
\]
(So $\int_{\partial\Omega}x\,dy$ literally measures signed area of $\Omega$.)

\begin{figure}[ht]
	\centering
	\begin{tikzpicture}[scale=1.05]
		\draw[->] (-0.2,0) -- (4.6,0) node[below] {$x$};
		\draw[->] (0,-0.2) -- (0,3.3) node[left] {$y$};
		
		\draw[thick] (1.0,0.7) .. controls (0.8,2.6) and (2.1,3.0) .. (3.2,2.3)
		.. controls (4.2,1.6) and (3.6,0.7) .. (2.5,0.6)
		.. controls (1.8,0.55) and (1.4,0.55) .. (1.0,0.7);
		
		\node at (2.4,1.6) {$\Omega$};
		
		\draw[->] (3.65,1.85) .. controls (3.55,2.05) and (3.35,2.25) .. (3.05,2.35);
		
		\node at (4.1,2.15) {\small $\partial\Omega$ (CCW)};
		
		\node at (2.25,3.15) {\small $\displaystyle \int_{\partial\Omega}\eta=\int_\Omega d\eta$};
	\end{tikzpicture}
	\caption{Green/Stokes in the plane: boundary integral of a $1$-form equals area integral of its exterior derivative.}
\end{figure}

\subsubsection{Complex form of area: $dx\wedge dy$ vs $d\bar z\wedge dz$}

Identify $\R^2\simeq\C$ by $z=x+iy$. Then
\[
dz=dx+i\,dy,\qquad d\bar z=dx-i\,dy.
\]
Compute $d\bar z\wedge dz$:
\[
d\bar z\wedge dz=(dx-i\,dy)\wedge(dx+i\,dy)
=dx\wedge dx+i\,dx\wedge dy-i\,dy\wedge dx+(-i^2)\,dy\wedge dy.
\]
Using $dx\wedge dx=dy\wedge dy=0$ and $dy\wedge dx=-dx\wedge dy$,
\[
d\bar z\wedge dz=i\,dx\wedge dy-i(-dx\wedge dy)=2i\,dx\wedge dy.
\]
Hence the area form satisfies
\[
dx\wedge dy=\frac{1}{2i}\,d\bar z\wedge dz.
\]
This identity is the bridge that turns Green/Stokes into complex-analytic formulas.

\subsubsection{Deriving the Cauchy--Green formula from Stokes' theorem}

Let $\Omega\subset\C$ be a bounded domain with smooth boundary $\partial\Omega$,
and let $f\in C^1(\overline\Omega)$.
Fix $z\in\Omega$.

\begin{theorem}[Cauchy--Green formula]
	For any $z\in\Omega$,
	\[
	f(z)=\frac{1}{2\pi i}
	\left[
	\int_{\partial\Omega}\frac{f(w)}{w-z}\,dw
	-\iint_{\Omega}\frac{\partial f/\partial\bar w}{w-z}\,d\bar w\wedge dw
	\right].
	\]
\end{theorem}

\paragraph{Step 0: notation and the key differential identity.}
Write $w=u+iv$ for the complex variable of integration.
Recall the operators
\[
\frac{\partial}{\partial w}
=\frac12\left(\frac{\partial}{\partial u}-i\frac{\partial}{\partial v}\right),
\qquad
\frac{\partial}{\partial \bar w}
=\frac12\left(\frac{\partial}{\partial u}+i\frac{\partial}{\partial v}\right),
\]
so that for $f(u,v)$,
\[
df=f_u\,du+f_v\,dv
=f_w\,dw+f_{\bar w}\,d\bar w,
\qquad
\text{where } dw=du+i\,dv,\ d\bar w=du-i\,dv.
\]
The following identity will be used repeatedly:
\[
d\bigl(f(w,\bar w)\,dw\bigr)=df\wedge dw
=\bigl(f_w\,dw+f_{\bar w}\,d\bar w\bigr)\wedge dw
=f_{\bar w}\,d\bar w\wedge dw
\quad (\text{since }dw\wedge dw=0).
\]

\paragraph{Step 1: choose the $1$-form to which we apply Stokes.}
Define on $\Omega\setminus\{z\}$ the $1$-form
\[
\omega(w):=\frac{f(w)}{w-z}\,dw.
\]
This is smooth away from $w=z$.

\paragraph{Step 2: compute $d\omega$.}
Using $d\big(\phi\,dw\big)=d\phi\wedge dw$ and $\phi=f(w)/(w-z)$,
\[
d\omega=d\!\left(\frac{f(w)}{w-z}\right)\wedge dw.
\]
Expand $d\!\left(\frac{f}{w-z}\right)$ using the product rule:
\[
d\!\left(\frac{f}{w-z}\right)
=\frac{df}{w-z}+f\,d\!\left(\frac{1}{w-z}\right).
\]
Now note:
\[
d\!\left(\frac{1}{w-z}\right)=\frac{\partial}{\partial w}\!\left(\frac{1}{w-z}\right)\,dw
+\frac{\partial}{\partial\bar w}\!\left(\frac{1}{w-z}\right)\,d\bar w.
\]
But $\frac{1}{w-z}$ depends only on $w$ (not on $\bar w$), so $\partial/\partial\bar w$ term is $0$.
Hence
\[
d\!\left(\frac{1}{w-z}\right)=\left(-\frac{1}{(w-z)^2}\right)dw.
\]
Therefore
\[
f\,d\!\left(\frac{1}{w-z}\right)\wedge dw
=f\left(-\frac{1}{(w-z)^2}\,dw\right)\wedge dw=0
\quad\text{because }dw\wedge dw=0.
\]
So the entire contribution comes from $\frac{df}{w-z}\wedge dw$:
\[
d\omega=\frac{df}{w-z}\wedge dw.
\]
Finally write $df=f_w\,dw+f_{\bar w}\,d\bar w$:
\[
d\omega=\frac{f_w\,dw+f_{\bar w}\,d\bar w}{w-z}\wedge dw
=\frac{f_w}{w-z}\,dw\wedge dw+\frac{f_{\bar w}}{w-z}\,d\bar w\wedge dw
=\frac{f_{\bar w}}{w-z}\,d\bar w\wedge dw.
\]
Thus we have proved the key identity
\[
d\left(\frac{f(w)}{w-z}\,dw\right)=\frac{f_{\bar w}(w)}{w-z}\,d\bar w\wedge dw
\quad\text{on }\Omega\setminus\{z\}.
\]

\paragraph{Step 3: excise a small disk and apply Stokes on a punctured domain.}
Let $D_\varepsilon(z)=\{w:|w-z|<\varepsilon\}$ and set
\[
\Omega_\varepsilon:=\Omega\setminus \overline{D_\varepsilon(z)}.
\]
Then $\partial\Omega_\varepsilon$ consists of the outer boundary $\partial\Omega$
and the inner boundary $-\partial D_\varepsilon(z)$ (minus sign because the induced
orientation is clockwise on the inner circle):
\[
\partial\Omega_\varepsilon=\partial\Omega-\partial D_\varepsilon(z).
\]
Apply Stokes:
\[
\int_{\partial\Omega_\varepsilon}\omega
=\iint_{\Omega_\varepsilon} d\omega.
\]
That is,
\[
\int_{\partial\Omega}\frac{f(w)}{w-z}\,dw
-\int_{\partial D_\varepsilon(z)}\frac{f(w)}{w-z}\,dw
=
\iint_{\Omega_\varepsilon}\frac{f_{\bar w}(w)}{w-z}\,d\bar w\wedge dw.
\]

\paragraph{Step 4: take $\varepsilon\to 0$ and compute the small-circle limit.}
Parametrize $\partial D_\varepsilon(z)$ by $w=z+\varepsilon e^{it}$, $0\le t\le 2\pi$.
Then
\[
dw=i\varepsilon e^{it}\,dt,\qquad w-z=\varepsilon e^{it},
\qquad \frac{dw}{w-z}=i\,dt.
\]
Hence
\[
\int_{\partial D_\varepsilon(z)}\frac{f(w)}{w-z}\,dw
=\int_0^{2\pi} f(z+\varepsilon e^{it})\, i\,dt.
\]
Since $f$ is continuous on $\overline\Omega$,
$f(z+\varepsilon e^{it})\to f(z)$ uniformly in $t$ as $\varepsilon\to0$, so
\[
\int_{\partial D_\varepsilon(z)}\frac{f(w)}{w-z}\,dw
\longrightarrow
\int_0^{2\pi} f(z)\, i\,dt
=2\pi i\,f(z).
\]
Also, since the integrand $\frac{f_{\bar w}(w)}{w-z}$ is integrable on $\Omega$
	(the kernel $|w-z|^{-1}$ is locally integrable in two real dimensions; no principal value is needed),
we have
\[
\iint_{\Omega_\varepsilon}\frac{f_{\bar w}(w)}{w-z}\,d\bar w\wedge dw
\longrightarrow
\iint_{\Omega}\frac{f_{\bar w}(w)}{w-z}\,d\bar w\wedge dw.
\]
Taking the limit in the Stokes identity yields
\[
\int_{\partial\Omega}\frac{f(w)}{w-z}\,dw
-2\pi i\,f(z)
=
\iint_{\Omega}\frac{f_{\bar w}(w)}{w-z}\,d\bar w\wedge dw.
\]
Rearrange:
\[
f(z)=\frac{1}{2\pi i}\left[\int_{\partial\Omega}\frac{f(w)}{w-z}\,dw
-\iint_{\Omega}\frac{f_{\bar w}(w)}{w-z}\,d\bar w\wedge dw\right],
\]
which is the Cauchy--Green formula.

\begin{remark}[The holomorphic case gives the Cauchy integral formula]
	If $f$ is holomorphic on $\Omega$, then $f_{\bar w}=0$ in $\Omega$.
	The area term disappears and we recover
	\[
	f(z)=\frac{1}{2\pi i}\int_{\partial\Omega}\frac{f(w)}{w-z}\,dw.
	\]
\end{remark}

\begin{remark}[Where Green/Stokes is hiding]
	Using $d\bar w\wedge dw=2i\,du\wedge dv$, the area integral is a Green-type integral
	on $\Omega$.  The Cauchy--Green identity is literally Stokes applied to a meromorphic-looking $1$-form.
\end{remark}

\subsubsection{Surface area, Jacobians, and the area 2-form (reparametrization invariance)}

Now move to $\R^3$ and let $\sigma:U\to\R^3$ be a smooth surface patch, where
$U\subset\R^2$ has coordinates $(u,v)$. Let $R\subset U$ be a region.
For the area of the image itself, assume that $\sigma$ is regular and one-to-one
on the interior of $R$, apart from boundary identifications of measure zero.
Without that assumption, the same integral gives parametrized area, counting
overlapping sheets with multiplicity.

\paragraph{Area as a limit of parallelograms (geometric motivation).}
At a point $(u,v)$, the differential $d\sigma$ sends small parameter vectors
$\Delta u\,\partial_u$ and $\Delta v\,\partial_v$ to tangent vectors
\[
\sigma_u\,\Delta u,\qquad \sigma_v\,\Delta v
\]
on the surface. The spanned parallelogram has area
\[
\|\sigma_u\times\sigma_v\|\,|\Delta u\,\Delta v|.
\]
Integrating this density leads to the surface area.

\begin{definition}[Area of a surface patch]
	Under the preceding regularity and one-to-one hypothesis, the area of $\sigma(R)$ is
	\[
	A_\sigma(R)=\iint_R \|\sigma_u\times\sigma_v\|\,du\,dv.
	\]
\end{definition}

\paragraph{First fundamental form and $\sqrt{EG-F^2}$.}
Set
\[
E=\sigma_u\cdot\sigma_u,\qquad
F=\sigma_u\cdot\sigma_v,\qquad
G=\sigma_v\cdot\sigma_v.
\]
Then $\|\sigma_u\times\sigma_v\|^2$ is the Gram determinant:
\[
\|\sigma_u\times\sigma_v\|^2
=
\det\begin{pmatrix}
	\sigma_u\cdot\sigma_u & \sigma_u\cdot\sigma_v\\
	\sigma_v\cdot\sigma_u & \sigma_v\cdot\sigma_v
\end{pmatrix}
=\det\begin{pmatrix}E&F\\F&G\end{pmatrix}
=EG-F^2.
\]
Hence
\[
\|\sigma_u\times\sigma_v\|=\sqrt{EG-F^2}.
\]
So
\[
A_\sigma(R)=\iint_R \sqrt{EG-F^2}\,du\,dv.
\]

\paragraph{Area form and area density.}
Define the \emph{area $2$-form} on the parameter domain by
\[
\omega_\sigma:=\|\sigma_u\times\sigma_v\|\,du\wedge dv.
\]
Then
\[
A_\sigma(R)=\iint_R \omega_\sigma.
\]
Once the parameter domain is given the orientation $du\wedge dv>0$, this is a
positive $2$-form.  If no orientation is chosen, the corresponding intrinsic
object is the positive density
$\|\sigma_u\times\sigma_v\|\,|du\,dv|$.  The distinction matters under an
orientation-reversing change of parameters.

\begin{theorem}[Change of variables (Jacobian) in $2$-form language]
	Let $\Phi:\widetilde R\to R$ be a $C^1$ diffeomorphism, $\Phi(s,t)=(u(s,t),v(s,t))$.
	Then
	\[
	\Phi^*(du\wedge dv)=\det(D\Phi)\,ds\wedge dt.
	\]
	If $\Phi$ preserves orientation and $g$ is integrable on $R$, then
	\[
	\iint_R g(u,v)\,du\,dv
	=
	\iint_{\widetilde R}g\bigl(\Phi(s,t)\bigr)
	\det(D\Phi(s,t))\,ds\,dt.
	\]
	For unsigned area under an arbitrary diffeomorphism, replace the determinant
	by its absolute value.
\end{theorem}

\begin{proof}
	Compute
	\[
	du=u_s\,ds+u_t\,dt,\qquad dv=v_s\,ds+v_t\,dt.
	\]
	Wedge:
	\begin{align*}
		du\wedge dv
		&=(u_s\,ds+u_t\,dt)\wedge(v_s\,ds+v_t\,dt)\\
		&=u_s v_s\,ds\wedge ds+u_s v_t\,ds\wedge dt+u_t v_s\,dt\wedge ds+u_t v_t\,dt\wedge dt\\
		&=\bigl(u_s v_t-u_t v_s\bigr)\,ds\wedge dt\\
		&=\det(D\Phi)\,ds\wedge dt,
	\end{align*}
	using $ds\wedge ds=dt\wedge dt=0$ and $dt\wedge ds=-ds\wedge dt$.
	This is exactly the Jacobian determinant.  Integrating the pullback of
	$g\,du\wedge dv$ gives the displayed orientation-preserving formula; forgetting
	orientation replaces the signed Jacobian by $|\det D\Phi|$.
\end{proof}

\paragraph{Reparametrization invariance of surface area.}
First suppose that $\Phi$ preserves orientation, and let
$\widetilde\sigma=\sigma\circ\Phi$ be a reparametrization, so
\[
\widetilde\sigma(s,t)=\sigma(u(s,t),v(s,t)).
\]
Then the area form satisfies
\[
\omega_{\widetilde\sigma}
=\|\widetilde\sigma_s\times\widetilde\sigma_t\|\,ds\wedge dt
	=\|\sigma_u\times\sigma_v\|\;\det(D\Phi)\,ds\wedge dt
=\Phi^*\bigl(\|\sigma_u\times\sigma_v\|\,du\wedge dv\bigr)
=\Phi^*(\omega_\sigma).
\]
Therefore, for $\widetilde R$ mapping to $R$,
\[
A_{\widetilde\sigma}(\widetilde R)=\iint_{\widetilde R}\omega_{\widetilde\sigma}
=\iint_{\widetilde R}\Phi^*(\omega_\sigma)
=\iint_R \omega_\sigma
=A_\sigma(R),
\]
so the area is invariant.
For an orientation-reversing reparametrization the pullback of the oriented $2$-form changes
sign; surface area is instead the associated positive density, whose transformation law uses
$|\det D\Phi|$.  This distinction prevents conflating oriented integration with unsigned area.

\begin{figure}[ht]
	\centering
	\begin{tikzpicture}[scale=1.0]
		\draw[thick] (0,0) rectangle (3.2,2.2);
		\node at (1.6,1.1) {$R\subset (u,v)$-plane};
		
		\draw[->,thick] (3.4,1.1) -- (4.6,1.1) node[midway,above] {$\sigma$};
		
		\draw[thick] (5.0,0.4) .. controls (6.2,2.0) and (7.3,1.6) .. (7.7,0.9)
		.. controls (8.0,0.3) and (7.1,-0.2) .. (6.2,0.0)
		.. controls (5.5,0.15) and (5.2,0.2) .. (5.0,0.4);
		\node at (6.5,1.15) {$\sigma(R)$};
		
		\node at (1.6,-0.55) {\small $\omega_\sigma=\|\sigma_u\times\sigma_v\|\,du\wedge dv$};
	\end{tikzpicture}
	\caption{Surface area is the integral of the pulled-back area $2$-form: coordinate changes act via Jacobians.}
\end{figure}

\subsubsection{Worked examples (Cauchy--Green + area 2-form)}

\begin{example}[Cauchy--Green on the unit disk with $f(w,\bar w)=\bar w$]
	Let $\Omega=\{w:|w|<1\}$ and $f(w,\bar w)=\bar w$.
	Compute the Cauchy--Green formula and interpret the result.
\end{example}

\begin{proof}[Detailed computation]
	We have $\partial f/\partial\bar w=1$.
	Fix $z\in\Omega$.
	Cauchy--Green gives
	\[
	\bar z=\frac{1}{2\pi i}\left[\int_{|w|=1}\frac{\bar w}{w-z}\,dw
	-\iint_{|w|<1}\frac{1}{w-z}\,d\bar w\wedge dw\right].
	\]
	On $|w|=1$, $\bar w=1/w$, so the boundary term becomes
	\[
	\int_{|w|=1}\frac{\bar w}{w-z}\,dw=\int_{|w|=1}\frac{1}{w(w-z)}\,dw.
	\]
	If $z=0$, the integrand is $1/w^2$, whose contour integral is zero.
	For $z\ne0$, it has distinct simple poles at $w=0$ and $w=z$ (both inside the unit circle). By residues,
	\[
	\int_{|w|=1}\frac{1}{w(w-z)}\,dw
	=2\pi i\left(\mathrm{Res}_{w=0}\frac{1}{w(w-z)}+\mathrm{Res}_{w=z}\frac{1}{w(w-z)}\right).
	\]
	Compute
	\[
	\mathrm{Res}_{w=0}\frac{1}{w(w-z)}=\lim_{w\to0}\frac{1}{w-z}=-\frac{1}{z},
	\qquad
	\mathrm{Res}_{w=z}\frac{1}{w(w-z)}=\lim_{w\to z}\frac{1}{w}=\frac{1}{z}.
	\]
	So the boundary term is $0$.
	
	Thus the formula reduces to
	\[
	\bar z=-\frac{1}{2\pi i}\iint_{|w|<1}\frac{1}{w-z}\,d\bar w\wedge dw.
	\]
	Using $d\bar w\wedge dw=2i\,dA$ (where $dA$ is area measure),
	\[
	\bar z=-\frac{1}{2\pi i}\iint_{|w|<1}\frac{1}{w-z}\,(2i\,dA)
	=-\frac{1}{\pi}\iint_{|w|<1}\frac{1}{w-z}\,dA.
	\]
	Equivalently,
	\[
	\iint_{|w|<1}\frac{1}{w-z}\,dA=-\pi\,\bar z.
	\]
	This is consistent with the fact that $\bar z$ is not holomorphic: the entire value comes from the $\bar\partial$-term.
\end{proof}

\begin{example}[Area of the unit sphere via $\sqrt{EG-F^2}$]
	Parametrize the unit sphere by
	\[
	\sigma(\theta,\phi)=(\sin\phi\cos\theta,\ \sin\phi\sin\theta,\ \cos\phi),
	\qquad 0\le\theta\le2\pi,\ 0\le\phi\le\pi.
	\]
\end{example}

\begin{proof}[Detailed computation]
	Compute partials:
	\[
	\sigma_\theta=(-\sin\phi\sin\theta,\ \sin\phi\cos\theta,\ 0),
	\qquad
	\sigma_\phi=(\cos\phi\cos\theta,\ \cos\phi\sin\theta,\ -\sin\phi).
	\]
	Then
	\[
	E=\sigma_\theta\cdot\sigma_\theta=\sin^2\phi,\qquad
	F=\sigma_\theta\cdot\sigma_\phi=0,\qquad
	G=\sigma_\phi\cdot\sigma_\phi=1.
	\]
	Hence
	\[
	\sqrt{EG-F^2}=\sqrt{\sin^2\phi\cdot 1-0}=\sin\phi.
	\]
	So the area is
	\[
	A(S^2)=\int_0^{2\pi}\int_0^\pi \sin\phi\,d\phi\,d\theta
	=\int_0^{2\pi}\left[-\cos\phi\right]_{\phi=0}^{\phi=\pi}\,d\theta
	=\int_0^{2\pi}(2)\,d\theta
	=4\pi.
	\]
\end{proof}

\subsubsection{Exercises}

\begin{exercise}[Green from Stokes via $d(Pdx+Qdy)$]
	Write $\eta=Pdx+Qdy$ and compute $d\eta$ explicitly to show
	\[
	\int_{\partial\Omega}(Pdx+Qdy)=\iint_\Omega(Q_x-P_y)\,dx\,dy.
	\]
	Then test the formula on $\eta=x\,dy$ by deducing $\int_{\partial\Omega}x\,dy=\mathrm{Area}(\Omega)$.
\end{exercise}

\begin{exercise}[Holomorphicity from Cauchy--Green]
	Let $f\in C^1(\overline\Omega)$. Use Cauchy--Green to prove:
	if $f_{\bar z}=0$ on $\Omega$, then $f$ satisfies the Cauchy integral formula on every
	smooth subdomain of $\Omega$, hence is holomorphic on $\Omega$.
	(You may use Morera's theorem as the final step.)
\end{exercise}

\begin{exercise}[Compute a boundary term explicitly on the disk]
	Let $\Omega=\{w:|w|<1\}$ and $z\in\Omega$. Compute
	\[
	\int_{|w|=1}\frac{1}{w-z}\,dw
	\]
	by parametrizing $w=e^{it}$ and by residues, and reconcile the two computations.
\end{exercise}

\begin{exercise}[A Cauchy--Green computation with a nonholomorphic function]
	On $\Omega=\{w:|w|<1\}$ take $f(w,\bar w)=|w|^2=w\bar w$.
	\begin{enumerate}
		\item Compute $f_{\bar w}$.
		\item Write the Cauchy--Green formula for $f(z)$ and simplify as much as possible.
		\item Explain why the $\bar\partial$-term must be nonzero.
	\end{enumerate}
\end{exercise}

\begin{exercise}[Jacobian via wedge products]
	Let $\Phi(s,t)=(u(s,t),v(s,t))$ be $C^1$. Prove directly that
	\[
	\Phi^*(du\wedge dv)=\det(D\Phi)\,ds\wedge dt
	\]
	by expanding $du=u_sds+u_tdt$ and $dv=v_sds+v_tdt$ and wedging.
\end{exercise}

\begin{exercise}[Change of variables under a regular reparametrization]
	Let $\Phi(u,v)=(x,y)=\bigl(u,(v+v^2)/2\bigr)$ map
	$\widetilde R=[0,1]\times[0,1]$ onto $R=[0,1]\times[0,1]$.
	\begin{enumerate}
		\item Compute $\det D\Phi$ and verify that it is positive on $\widetilde R$.
		\item Use change of variables to express $\iint_R f(x,y)\,dx\,dy$ as an integral over $\widetilde R$.
		\item Apply to $f(x,y)=1$ to compute $\mathrm{Area}(R)$.
	\end{enumerate}
\end{exercise}

\begin{exercise}[Graph surface area: paraboloid patch]
	Let $\sigma(u,v)=(u,v,u^2+v^2)$ on $R=[0,1]\times[0,1]$.
	\begin{enumerate}
		\item Compute $\sigma_u,\sigma_v$ and then $E,F,G$.
		\item Compute $\sqrt{EG-F^2}$ and write $A_\sigma(R)$ as an explicit double integral.
		\item Give a geometric interpretation: compare to the area of the flat square $R$.
	\end{enumerate}
\end{exercise}

\begin{exercise}[Reparametrization invariance on an explicit map]
	Let $\sigma(u,v)=(u,v,0)$ and
	$\Phi(s,t)=(u,v)=\bigl((s+s^2)/2,t\bigr)$ on $[0,1]^2$.
	\begin{enumerate}
		\item Compute $A_\sigma(R)$ for $R=[0,1]\times[0,1]$.
		\item Verify that $\Phi$ is a regular orientation-preserving
		reparametrization, then compute $A_{\sigma\circ\Phi}(\widetilde R)$ directly
		and check that it equals $A_\sigma(R)$.
	\end{enumerate}
\end{exercise}

\subsection{Gaussian Curvature and the Gauss--Bonnet Theorem}
\index{Gauss--Bonnet theorem}\index{Gaussian curvature}
\label{subsec:gauss-curvature-gauss-bonnet}

\newcommand{\Step}[1]{\medskip\noindent\textbf{#1.}\ }

This subsection studies one quantity in four complementary languages.  The
first fundamental form and shape operator give a direct computation for
embedded surfaces; isothermal coordinates express the same metric as a
Hermitian metric; the conformal factor converts curvature into a
Poisson--Liouville equation; and Gauss--Bonnet identifies the integral of
curvature with topology.  The sphere and the torus will be calculated by both
the classical and conformal methods, so that the reader can see that the two
languages produce the same answer.

Throughout the subsection, Gaussian curvature is treated as an
\emph{intrinsic} scalar field, and the central global quantity is
\[
\iint_M K\,dA.
\]
We distinguish the nonnegative Riemannian area measure $dA$ used in scalar
integrals from the oriented area $2$-form $\mathrm{vol}_M$. In an oriented
coordinate chart they have the same positive coefficient; reversing the chart
changes the sign of the form but not the area measure. Accordingly,
$K\,dA$ is an invariant density and $K\,\mathrm{vol}_M$ is an invariant
oriented $2$-form.

\subsubsection{First and second fundamental forms: what is invariant}

Let $M\subset\R^3$ be a smooth oriented surface and $\sigma:U\to M$ a local
parametrization, $U\subset\R^2$ with coordinates $(u,v)$.

\textbf{First fundamental form.}
Write
\[
\sigma_u=\frac{\partial\sigma}{\partial u},\qquad
\sigma_v=\frac{\partial\sigma}{\partial v},
\]
and define
\[
E=\sigma_u\cdot\sigma_u,\qquad
F=\sigma_u\cdot\sigma_v,\qquad
G=\sigma_v\cdot\sigma_v.
\]
Then the induced metric is $ds^2=E\,du^2+2F\,du\,dv+G\,dv^2$.

\textbf{Area $2$-form (from Gram determinant).}
The area density is
\[
\|\sigma_u\times\sigma_v\|
=\sqrt{\det
	\begin{pmatrix}
		\sigma_u\cdot\sigma_u & \sigma_u\cdot\sigma_v\\
		\sigma_v\cdot\sigma_u & \sigma_v\cdot\sigma_v
\end{pmatrix}}
=\sqrt{EG-F^2}.
\]
Hence the Riemannian area measure is locally
$dA=\sqrt{EG-F^2}\,|du\,dv|$. If $(u,v)$ is positively oriented, the
oriented area form is
$\mathrm{vol}_M=\sqrt{EG-F^2}\,du\wedge dv$; in a negatively oriented
chart the same global form has the opposite sign.
For any scalar function $\varphi$ on $M$,
\[
\iint_{\sigma(R)} \varphi\,dA
=\iint_R (\varphi\circ\sigma)\,\sqrt{EG-F^2}\,du\,dv.
\]

\textbf{Reparametrization invariance (coordinate-free statement).}
If $\Phi:\widetilde U\to U$ is a diffeomorphism and $\widetilde\sigma=\sigma\circ\Phi$,
then $dA$ is intrinsic on the surface and
\[
\iint_{\widetilde\sigma(\widetilde R)} \varphi\,dA
=\iint_{\sigma(R)} \varphi\,dA.
\]
Equivalently, $dA$ is a globally defined density, independent of coordinates.
The oriented area form is also global, but its local sign records whether a
chart agrees with the chosen orientation.

\textbf{Second fundamental form and Gaussian curvature.}
Choose a unit normal field $N$ compatible with the orientation and set
\[
h_{11}=N\cdot\sigma_{uu},\qquad
h_{12}=N\cdot\sigma_{uv},\qquad
h_{22}=N\cdot\sigma_{vv}.
\]
The second fundamental form is
\[
\mathrm{II}=h_{11}\,du^2+2h_{12}\,du\,dv+h_{22}\,dv^2.
\]
With the shape-operator convention $S=-dN$, its determinant is
\[
K=\det S
 =\frac{h_{11}h_{22}-h_{12}^2}{EG-F^2}.
\]
Reversing $N$ changes the signs of both principal curvatures but not their
product, so $K$ is independent of the chosen unit normal.  The formula uses an
embedding to compute $K$, but Gauss's \emph{Theorema Egregium} shows that the
result depends only on the first fundamental form.  Thus $K$ is intrinsic even
though the second fundamental form is not; see \cite[Chapter~4]{Tapp2016}.

\medskip
\noindent
\textbf{Three objects that should not be conflated.}
The scalar $K$, the positive density $K\,dA$, and the oriented curvature form
$K\,\mathrm{vol}_M$ transform differently under an orientation reversal.
The scalar and density are unchanged, whereas the oriented form changes sign
with the orientation.  Integrals below use the density unless an oriented
$2$-form is displayed explicitly.

\subsubsection{Compatible metrics, K\"ahler forms, and the curvature equation}
\index{compatible metric}\index{K\"ahler form}\index{prescribed curvature}

The first fundamental form of a surface and the Hermitian metric of a
Riemann surface are two descriptions of the same object only after a
compatibility condition has been imposed.  This point is easy to miss: an
arbitrary Riemannian metric need not be compatible with a previously fixed
complex structure.  Once compatibility holds, however, complex dimension one
has a special feature---the associated fundamental form is automatically
closed.  The resulting formulas also explain precisely why recovering a metric
from curvature leads to a Poisson equation in one formulation and to a
nonlinear Liouville equation in another.

\begin{definition}[Compatible and Hermitian metrics]
Let $(X,J)$ be a Riemann surface.  A Riemannian metric $g$ is
\emph{compatible with $J$}, or \emph{Hermitian}, if
\[
 g(JV,JW)=g(V,W)
 \qquad\text{for all }V,W\in TX.
\]
Its \emph{fundamental form} is the real $2$-form
\[
 \omega_g(V,W):=g(JV,W).
\]
The metric is \emph{K\"ahler} if $d\omega_g=0$.
These are intrinsic definitions: no embedding of $X$ into Euclidean space is
part of the data.
\end{definition}

\begin{definition}[Conformal metrics and isothermal coordinates]
Two Riemannian metrics $g$ and $g_0$ on a surface are \emph{conformally
equivalent} if $g=e^{2u}g_0$ for some smooth real function $u$.  A positively
oriented coordinate $(x,y)$ is \emph{isothermal} for $g$ if
\[
 g=\lambda(x,y)(dx^2+dy^2),\qquad \lambda>0.
\]
In the associated complex coordinate $z=x+iy$, we use the abbreviation
$|dz|^2:=dx^2+dy^2$.  Thus $\lambda|dz|^2$ denotes a real symmetric
$2$-tensor; it must not be confused with the alternating form
$dz\wedge d\bar z$.
\end{definition}

\begin{definition}[Curvature form]
On an oriented Riemannian surface $(X,g)$, let $K_g$ be the Gaussian curvature
and let $\omega_g$ be the oriented area form.  The $2$-form
\[
 \Omega_g:=K_g\omega_g
\]
is called the \emph{curvature form} in this subsection.  When $g$ is viewed as
a Hermitian metric on a Riemann surface, this is the Ricci form with the
conventions used here.  Prescribing $\Omega_g$ is different from prescribing
the scalar function $K_g$, because $\omega_g$ itself depends on the unknown
metric.
\end{definition}

\begin{proposition}[Compatible first fundamental forms are Hermitian]
\label{prop:compatible-first-fundamental-hermitian}
Let $(X,J)$ be a Riemann surface and let $g$ be a Riemannian metric satisfying
\[
 g(JV,JW)=g(V,W)
 \qquad (V,W\in TX).
\]
In every holomorphic coordinate $z=x+iy$ there is a positive smooth function
$\lambda$ such that
\[
 g=\lambda(x,y)(dx^2+dy^2)=\lambda(z,\bar z)|dz|^2.
\]
Thus, if $g=E\,dx^2+2F\,dx\,dy+G\,dy^2$, compatibility is equivalent to
$E=G=\lambda$ and $F=0$.  The associated fundamental form
\[
 \omega(V,W):=g(JV,W)
\]
is
\[
 \omega=\lambda\,dx\wedge dy
       =\frac{i}{2}\lambda\,dz\wedge d\bar z.
\]
In particular, every Hermitian metric on a Riemann surface is K\"ahler.
Compactness is not required for this conclusion.
\end{proposition}

\begin{proof}
In a holomorphic coordinate,
$J\partial_x=\partial_y$ and $J\partial_y=-\partial_x$.  Put
$E=g(\partial_x,\partial_x)$, $F=g(\partial_x,\partial_y)$, and
$G=g(\partial_y,\partial_y)$.  Compatibility gives
\[
 G=g(J\partial_x,J\partial_x)=E
\]
and
\[
 F=g(J\partial_x,J\partial_y)
   =g(\partial_y,-\partial_x)=-F,
\]
so $F=0$.  Positivity gives $E=G=\lambda>0$.  Evaluating
$\omega$ on $(\partial_x,\partial_y)$ gives
$\omega(\partial_x,\partial_y)=g(\partial_y,\partial_y)=\lambda$.
Since $dz\wedge d\bar z=-2i\,dx\wedge dy$, the displayed complex formula
follows.  Finally, $d\omega$ is a $3$-form on the real surface $X$, hence
vanishes identically.  Therefore the Hermitian metric is K\"ahler.
\end{proof}

\begin{corollary}[Every Riemann surface carries a K\"ahler metric]
\label{cor:riemann-surface-kahler-metric-exists}
Every Riemann surface, compact or noncompact, admits a Hermitian metric, and
every such metric is K\"ahler.
\end{corollary}

\begin{proof}
Choose a holomorphic atlas and a locally finite partition of unity subordinate
to it.  In each chart, $dx^2+dy^2$ is a $J$-compatible positive metric.
The partition-of-unity weighted sum of these local metrics is again positive
and $J$-compatible, so it defines a global Hermitian metric.  Its fundamental
form is closed by Proposition~\ref{prop:compatible-first-fundamental-hermitian}.
\end{proof}

\begin{proposition}[An oriented metric determines a Riemann surface structure]
\label{prop:oriented-metric-determines-complex-structure}
Let $(M,g)$ be an oriented Riemannian surface.  There is a unique bundle map
$J:TM\to TM$ that rotates every positively oriented orthonormal frame through
$90$ degrees.  It satisfies $J^2=-\operatorname{id}$ and is integrable.
Consequently, $M$ becomes a Riemann surface, $g$ is compatible with its complex
structure, and $g$ is a K\"ahler metric on that Riemann surface.
\end{proposition}

\begin{proof}
The orientation makes the rotation independent of the chosen positive
orthonormal frame, and the construction immediately gives $J^2=-\operatorname{id}$
and $g(JV,JW)=g(V,W)$.  Every almost-complex structure on a real
two-dimensional manifold is integrable: equivalently, one may construct local
isothermal coordinates for $g$, in which
$g=\lambda(dx^2+dy^2)$ and $J\partial_x=\partial_y$.
The coordinate $z=x+iy$ is then holomorphic.  The final assertion follows from
Proposition~\ref{prop:compatible-first-fundamental-hermitian}.  See
\cite[Chapters~2--3]{Jost2006} for the relation among conformal structures,
metrics, and complex coordinates on surfaces.
\end{proof}

The preceding proposition is the intrinsic version of the classical passage
from an immersed oriented surface to isothermal coordinates.  For an immersion
$r(u,v)$ into $\R^3$, the first fundamental form
$I=E\,du^2+2F\,du\,dv+G\,dv^2$ becomes
$I=\lambda(dx^2+dy^2)$ in such coordinates.  Writing $z=x+iy$ converts this
same tensor into a Hermitian metric, while its oriented area form becomes the
K\"ahler form.  No extra differential equation is needed to prove that this
$2$-form is closed; real dimension two already forces that fact.

\begin{proposition}[Gaussian curvature under a conformal change]
\label{prop:curvature-conformal-change-surface}
Let $(X,g_0)$ be an oriented Riemannian surface, write $K_0$ and $\omega_0$ for
its Gaussian curvature and area form, and let
\[
 g=e^{2u}g_0,\qquad \omega_g=e^{2u}\omega_0
\]
for $u\in C^\infty(X)$.  With the convention
$\Delta_{g_0}=\operatorname{div}_{g_0}\nabla$ (so the Euclidean formula is
$\partial_x^2+\partial_y^2$), the Gaussian curvature of $g$
satisfies
\[
 K_g=e^{-2u}\bigl(K_0-\Delta_{g_0}u\bigr),
 \qquad
 K_g\omega_g=\bigl(K_0-\Delta_{g_0}u\bigr)\omega_0.
\]
In a flat coordinate, where $g=e^{2u}|dz|^2$, this reduces to
\[
 K_g=-e^{-2u}(u_{xx}+u_{yy}).
\]
\end{proposition}

\begin{proof}
Take an isothermal coordinate for $g_0$ and write
$g_0=e^{2v}|dz|^2$.  The flat-coordinate curvature formula gives
\[
 K_0=-e^{-2v}\Delta_0v,
 \qquad
 K_g=-e^{-2(v+u)}\Delta_0(v+u),
\]
where $\Delta_0=\partial_x^2+\partial_y^2$.  Since
$\Delta_{g_0}u=e^{-2v}\Delta_0u$, subtraction yields
$K_g=e^{-2u}(K_0-\Delta_{g_0}u)$.  Multiplication by
$\omega_g=e^{2u}\omega_0$ proves the curvature-form identity.  Taking
$v=0$ gives the final flat formula.  This computation is also the
real-coordinate form of
$\rho=-i\partial\bar\partial\log g_{z\bar z}=K_g\omega_g$.
\end{proof}

\begin{proposition}[Prescribing the curvature form is a Poisson problem]
\label{prop:curvature-form-poisson}
Let $(X,g_0)$ be a compact connected Riemann surface and let
$\Theta=f\,\omega_0$ be a smooth real $2$-form.  There exists a conformal
metric $g=e^{2u}g_0$ satisfying
\[
 K_g\omega_g=\Theta
\]
if and only if
\[
 \int_X\Theta=2\pi\chi(X).
\]
When it exists, $u$ solves the linear Poisson equation
\[
 \Delta_{g_0}u=K_0-f
\]
and is unique up to an additive constant.  Accordingly, the metric is unique
up to multiplication by a positive constant.
\end{proposition}

\begin{proof}
Proposition~\ref{prop:curvature-conformal-change-surface} shows that the
curvature-form equation is equivalent to
$\Delta_{g_0}u=K_0-f$.  On a compact connected manifold, the Poisson equation
$\Delta_{g_0}u=h$ is solvable precisely when $\int_Xh\,\omega_0=0$, and its
solutions then differ by constants; this is proved in
Section~\ref{subsec:hodge-top-degree-poisson}.  Thus solvability is equivalent
to
\[
 \int_X f\,\omega_0=\int_XK_0\,\omega_0=2\pi\chi(X),
\]
where the last equality is Gauss--Bonnet.  Adding a constant to $u$ rescales
$g$ but leaves the curvature form $K_g\omega_g$ unchanged.
\end{proof}

\begin{proposition}[The Calabi--Yau theorem in complex dimension one]
\label{prop:calabi-yau-riemann-surface-poisson}
Let $(X,\omega_0)$ be a compact connected K\"ahler Riemann surface, and let
$F\in C^\infty(X,\R)$ satisfy the volume normalization
\[
 \int_X e^F\omega_0=\int_X\omega_0.
\]
Then there exists a smooth real function $\varphi$, unique up to an additive
constant, such that
\[
 \omega_\varphi:=\omega_0+i\partial\bar\partial\varphi
 =e^F\omega_0.
\]
In particular, $\omega_\varphi$ is positive and is the unique K\"ahler form in
the class $[\omega_0]$ having volume form $e^F\omega_0$.  With the convention
$\Delta_{g_0}=\operatorname{div}_{g_0}\nabla$, the equation is exactly the
Poisson equation
\[
 \Delta_{g_0}\varphi=2(e^F-1).
\]
Thus the complex Monge--Amp\`ere equation in the Calabi--Yau theorem becomes a
linear Poisson equation on a compact Riemann surface.
\end{proposition}

\begin{proof}
In a holomorphic coordinate, write
\[
 \omega_0=\frac{i}{2}\lambda\,dz\wedge d\bar z,
 \qquad
 g_0=\lambda|dz|^2.
\]
For a real function $\varphi$,
\[
 i\partial\bar\partial\varphi
 =\frac{i}{4}(\partial_x^2+\partial_y^2)\varphi\,
   dz\wedge d\bar z
 =\frac12(\Delta_{g_0}\varphi)\omega_0.
\]
Consequently, $\omega_0+i\partial\bar\partial\varphi=e^F\omega_0$ is
equivalent to
\[
 \Delta_{g_0}\varphi=2(e^F-1).
\]
The right-hand side has integral zero precisely because of the stated volume
normalization.  The compact Poisson solvability theorem from
Section~\ref{subsec:hodge-top-degree-poisson} therefore produces a smooth
solution, and connectedness makes it unique modulo constants.  The equation
itself gives $\omega_\varphi=e^F\omega_0>0$.  Finally,
$i\partial\bar\partial\varphi$ is exact, so
$[\omega_\varphi]=[\omega_0]$.  This is the $n=1$ specialization of the
Calabi--Yau theorem; compare \cite{Yau1978}.
\end{proof}

\begin{remark}[Why the higher-dimensional equation becomes linear]
On a compact K\"ahler manifold of complex dimension $n$, the corresponding
equation is
\[
 (\omega_0+i\partial\bar\partial\varphi)^n=e^F\omega_0^n,
\]
which is nonlinear when $n\ge2$.  For a Riemann surface, $n=1$, so there is no
higher wedge power: the left side is simply
$\omega_0+i\partial\bar\partial\varphi$.  This statement prescribes the
K\"ahler volume form within a fixed cohomology class.  It should not be
confused with prescribing an arbitrary Gaussian-curvature function, which
leads to the nonlinear Liouville equation below.
\end{remark}

\begin{remark}[Curvature form versus scalar curvature]
The preceding linear equation prescribes the $2$-form $K_g\omega_g$.
Prescribing the scalar function $K_g=K$ is different because the unknown area
form also changes.  The equation becomes
\[
 \Delta_{g_0}u=K_0-Ke^{2u},
\]
which is the nonlinear Liouville equation.  Thus the slogan ``curvature leads
to a Poisson equation'' is literally correct for the curvature form; for
prescribed Gaussian curvature it must be understood as a nonlinear Poisson
problem.  This distinction is the one-dimensional prototype of the
curvature-potential equations that appear later in K\"ahler geometry.
\end{remark}

\subsubsection{Example 1: the unit sphere $S^2$}

\begin{example}[Unit sphere $S^2$: $\displaystyle \iint_{S^2}K\,dA=4\pi$]
	Let $S^2=\{(x,y,z)\in\R^3:x^2+y^2+z^2=1\}$. Show that $K\equiv 1$ and compute
	$\iint_{S^2}K\,dA$.
\end{example}

\begin{proof}
	\Step{Step 1 (Parametrization)}
	Use spherical coordinates
	\[
	\sigma(\theta,\phi)=(\sin\phi\cos\theta,\ \sin\phi\sin\theta,\ \cos\phi),
	\qquad 0\le\theta\le 2\pi,\ 0\le\phi\le\pi.
	\]
	
	\Step{Step 2 (Compute $\sigma_\theta,\sigma_\phi$)}
	\[
	\sigma_\theta=(-\sin\phi\sin\theta,\ \sin\phi\cos\theta,\ 0),
	\qquad
	\sigma_\phi=(\cos\phi\cos\theta,\ \cos\phi\sin\theta,\ -\sin\phi).
	\]
	
	\Step{Step 3 (Compute $E,F,G$)}
	\[
	E=\sigma_\theta\cdot\sigma_\theta
	=\sin^2\phi(\sin^2\theta+\cos^2\theta)=\sin^2\phi,
	\]
	\[
	F=\sigma_\theta\cdot\sigma_\phi
	=(-\sin\phi\sin\theta)(\cos\phi\cos\theta)+(\sin\phi\cos\theta)(\cos\phi\sin\theta)+0\cdot(-\sin\phi)=0,
	\]
	\[
	G=\sigma_\phi\cdot\sigma_\phi
	=\cos^2\phi(\cos^2\theta+\sin^2\theta)+\sin^2\phi=1.
	\]
	
	\Step{Step 4 (Area element)}
	\[
	dA=\sqrt{EG-F^2}\,|d\theta\,d\phi|
	=\sin\phi\,|d\theta\,d\phi|.
	\]
	For the outward orientation, $\sigma_\theta\times\sigma_\phi$
	points inward. Thus the oriented area form in these coordinates is
	$\mathrm{vol}_{S^2}=\sin\phi\,d\phi\wedge d\theta$.
	
	\Step{Step 5 (Curvature and total curvature)}
	The outward unit normal is $N=\sigma$. Differentiating
	$N\cdot\sigma_\theta=N\cdot\sigma_\phi=0$ shows that the shape operator is
	$-dN=-\mathrm{id}$ under the convention used here. Its two principal
	curvatures are therefore both $-1$, and their product is $K=1$. Hence
	\[
	\iint_{S^2}K\,dA=\int_0^{2\pi}\int_0^\pi \sin\phi\,d\phi\,d\theta
	=\int_0^{2\pi}\bigl[-\cos\phi\bigr]_{\phi=0}^{\phi=\pi}\,d\theta
	=\int_0^{2\pi}2\,d\theta
	=4\pi.
	\]

	\Step{Step 6 (The same computation as a Poisson equation)}
	Let $z=x+iy$ be the stereographic coordinate on
	$S^2\setminus\{\text{north pole}\}$ and put $s=|z|$.  The round metric is
	\[
	g_{S^2}=\frac{4}{(1+s^2)^2}|dz|^2=e^{2u}|dz|^2,
	\qquad
	u(s)=\log 2-\log(1+s^2).
	\]
	For a radial function, $\Delta_0u=u_{ss}+s^{-1}u_s$.  Here
	\[
	u_s=-\frac{2s}{1+s^2},
	\qquad
	u_{ss}=-\frac{2(1-s^2)}{(1+s^2)^2},
	\]
	and therefore
	\[
	\Delta_0u=-\frac{4}{(1+s^2)^2}=-e^{2u}.
	\]
	The conformal curvature formula
	$K=-e^{-2u}\Delta_0u$ again gives $K=1$.  Equivalently, the curvature form
	satisfies the Poisson identity
	\[
	K\,\omega_{S^2}=-\Delta_0u\,dx\wedge dy
	       =\frac{4}{(1+s^2)^2}\,dx\wedge dy.
	\]
	Although the stereographic chart omits one point, the displayed metric
	extends smoothly across that point.  Polar integration recovers the global
	total curvature:
	\[
	\iint_{S^2}K\,dA
	=2\pi\int_0^\infty\frac{4s}{(1+s^2)^2}\,ds
	=4\pi.
	\]
	Thus the extrinsic shape-operator calculation and the intrinsic
	Poisson/Liouville calculation give the same answer.
\end{proof}

\textbf{Why reparametrization cannot change the answer.}
$K$ is a scalar on the surface and $dA$ is an intrinsic density; hence $K\,dA$ is intrinsic
and its integral is coordinate-independent.

\subsubsection{Example 2: a torus of revolution (metric, area, curvature, sign)}

Let $R>r>0$ and consider
\[
\sigma(\theta,\phi)=\big((R+r\cos\phi)\cos\theta,\ (R+r\cos\phi)\sin\theta,\ r\sin\phi\big),
\qquad 0\le\theta,\phi\le 2\pi.
\]

\begin{example}[Torus: compute $E,F,G$, $dA$, $K$, and $\iint K\,dA$]
	Show that
	\[
	E=(R+r\cos\phi)^2,\quad F=0,\quad G=r^2,\qquad
	dA=r(R+r\cos\phi)\,|d\theta\,d\phi|,
	\]
	and that
	\[
	K(\phi)=\frac{\cos\phi}{r(R+r\cos\phi)}.
	\]
	Deduce $\displaystyle \iint_{T^2}K\,dA=0$ and identify where $K$ is positive/negative.
\end{example}

\begin{proof}
	\Step{Step 1 (Compute $\sigma_\theta,\sigma_\phi$)}
	\[
	\sigma_\theta=\big(-(R+r\cos\phi)\sin\theta,\ (R+r\cos\phi)\cos\theta,\ 0\big),
	\]
	\[
	\sigma_\phi=\big(-r\sin\phi\cos\theta,\ -r\sin\phi\sin\theta,\ r\cos\phi\big).
	\]
	
	\Step{Step 2 (Compute $E$)}
	\[
	E=\sigma_\theta\cdot\sigma_\theta
	=(R+r\cos\phi)^2(\sin^2\theta+\cos^2\theta)=(R+r\cos\phi)^2.
	\]
	
	\Step{Step 3 (Compute $F$)}
	\begin{align*}
		F&=\sigma_\theta\cdot\sigma_\phi\\
		&=(-(R+r\cos\phi)\sin\theta)(-r\sin\phi\cos\theta)\\
		&\quad+((R+r\cos\phi)\cos\theta)(-r\sin\phi\sin\theta)\\
		&=(R+r\cos\phi)r\sin\phi(\sin\theta\cos\theta-\cos\theta\sin\theta)=0.
	\end{align*}
	
	\Step{Step 4 (Compute $G$)}
	\[
	G=\sigma_\phi\cdot\sigma_\phi
	=r^2\sin^2\phi(\cos^2\theta+\sin^2\theta)+r^2\cos^2\phi
	=r^2(\sin^2\phi+\cos^2\phi)=r^2.
	\]
	
	\Step{Step 5 (Area element)}
	\[
	dA=\sqrt{EG-F^2}\,|d\theta\,d\phi|
	=r(R+r\cos\phi)\,|d\theta\,d\phi|.
	\]
	
	\Step{Step 6 (Gaussian curvature)}
	The outward unit normal is
	\[
	N=(\cos\phi\cos\theta,\ \cos\phi\sin\theta,\ \sin\phi).
	\]
	Taking scalar products of the second derivatives of $\sigma$ with $N$ gives
	\[
	h_{11}=N\cdot\sigma_{\theta\theta}=-(R+r\cos\phi)\cos\phi,
	\qquad h_{12}=N\cdot\sigma_{\theta\phi}=0,
	\qquad h_{22}=N\cdot\sigma_{\phi\phi}=-r.
	\]
	Therefore the quotient of the determinants of the second and first fundamental
	forms is
	\[
	K(\phi)=\frac{h_{11}h_{22}-h_{12}^2}{EG-F^2}
	=\frac{r(R+r\cos\phi)\cos\phi}{r^2(R+r\cos\phi)^2}
	=\frac{\cos\phi}{r(R+r\cos\phi)}.
	\]
	
	\Step{Step 7 (Total curvature)}
	Multiply:
	\[
	K\,dA=\cos\phi\,|d\theta\,d\phi|.
	\]
	Hence
	\[
	\iint_{T^2}K\,dA
	=\int_0^{2\pi}\int_0^{2\pi}\cos\phi\,d\phi\,d\theta
	=\left(\int_0^{2\pi}\cos\phi\,d\phi\right)\left(\int_0^{2\pi}d\theta\right)
	=0.
	\]
	
	\Step{Step 8 (Sign regions)}
	Since $r>0$ and $R+r\cos\phi>0$, $\mathrm{sign}(K)=\mathrm{sign}(\cos\phi)$:
	\[
	K>0 \iff \cos\phi>0 \iff \phi\in\left(-\frac{\pi}{2},\frac{\pi}{2}\right)\ (\mathrm{mod}\ 2\pi),
	\]
	\[
	K<0 \iff \cos\phi<0 \iff \phi\in\left(\frac{\pi}{2},\frac{3\pi}{2}\right)\ (\mathrm{mod}\ 2\pi),
	\]
	and $K=0$ along $\phi=\frac{\pi}{2},\,\frac{3\pi}{2}$.
	Geometrically: outer region has $K>0$, inner region near the hole has $K<0$.

	\Step{Step 9 (Isothermal coordinates and the Poisson equation)}
	The coordinates $(\theta,\phi)$ are orthogonal but not conformal because
	$E\ne G$.  Introduce a new periodic coordinate $t$ by
	\[
	t(\phi)=\int_0^\phi\frac{r}{R+r\cos s}\,ds,
	\qquad
	L:=t(2\pi)=\frac{2\pi r}{\sqrt{R^2-r^2}}.
	\]
	Then $t$ is taken modulo $L$ and
	$d\phi=(R+r\cos\phi)dt/r$.  Substitution into the first fundamental form
	gives
	\[
	g=(R+r\cos\phi(t))^2(d\theta^2+dt^2)
	  =e^{2u(t)}|dw|^2,
	\qquad
	w=\theta+it
	\]
	on the universal cover; it descends through the lattice generated by
	$2\pi$ and $iL$ to a conformal coordinate on the torus.  Here
	\[
	u(t)=\log(R+r\cos\phi(t)).
	\]
	Since $d\phi/dt=(R+r\cos\phi)/r$, direct differentiation gives
	\[
	u_t=-\sin\phi,
	\qquad
	u_{tt}=-\frac{(R+r\cos\phi)\cos\phi}{r}.
	\]
	The function $u$ is independent of $\theta$, so
	\[
	\Delta_0u=u_{tt}=-K e^{2u}.
	\]
	Hence the Poisson/Liouville curvature formula yields
	\[
	K=-e^{-2u}\Delta_0u
	 =\frac{\cos\phi}{r(R+r\cos\phi)},
	\]
	exactly as in Step~6.  It also gives a second proof of the total-curvature
	identity without using the second fundamental form:
	\[
	\iint_{T^2}K\,dA
	=-\int_0^{2\pi}\int_0^L u_{tt}\,dt\,d\theta
	=-2\pi\,[u_t]_{0}^{L}=0,
	\]
	because $u_t=-\sin\phi(t)$ is periodic.  Finally,
	$K\,\omega_g=-\Delta_0u\,d\theta\wedge dt$ becomes
	$\cos\phi\,d\theta\wedge d\phi$, agreeing with Step~7.
\end{proof}

\subsubsection{Gauss--Bonnet and Euler characteristic}

The two model calculations suggest a topological law: the sphere always has
total curvature $4\pi$, whereas positive and negative curvature cancel on the
torus.  Gauss--Bonnet explains that these are not accidents of the chosen
round sphere or embedded torus.  Changing the metric redistributes curvature,
but it cannot change its total integral on a compact surface.

\begin{theorem}[Gauss--Bonnet for compact surfaces without boundary]
	Let $M$ be a compact oriented smooth surface without boundary,
	equipped with a smooth Riemannian metric. Then
	\[
	\iint_M K\,dA=2\pi\,\chi(M),
	\]
	where $\chi(M)$ is the Euler characteristic of $M$.
\end{theorem}

\begin{proof}
	Choose a geodesic triangulation of $M$. The local Gauss--Bonnet formula
	for a geodesic triangle $T$ says that
	$\int_T K\,dA$ equals the sum of its interior angles minus $\pi$.
	On summing over all triangles, the angles around each vertex add to
	$2\pi$. Thus, writing $V,E,F$ for the numbers of vertices, edges,
	and faces, we obtain $\int_M K\,dA=2\pi V-\pi F$.
	Every edge belongs to two triangles, so $3F=2E$; consequently
	$2\pi V-\pi F=2\pi(V-E+F)=2\pi\chi(M)$.
	The local formula and the existence of a suitable triangulation are
	proved in \cite[Chapter~6]{Tapp2016}.
\end{proof}

\textbf{Consequences.}
\begin{itemize}
	\item $S^2$: $\chi(S^2)=2$ so $\iint_{S^2}K\,dA=4\pi$.
	\item $T^2$: $\chi(T^2)=0$ so $\iint_{T^2}K\,dA=0$.
	\item genus $g$: $\chi=2-2g$ so $\iint_M K\,dA=2\pi(2-2g)$.
\end{itemize}

\noindent
\textbf{What the section has established.}
For an embedded surface, $K$ can be computed from the two fundamental forms.
For a Riemann surface in an isothermal coordinate, the same quantity is
$K=-e^{-2u}\Delta_0u$.  Globally, the curvature form represents the
topological class whose integral is $2\pi\chi(M)$.  The sphere and torus
calculations therefore serve three purposes at once: they verify the local
formulas, illustrate the Poisson description, and test the global
Gauss--Bonnet invariant.

\subsubsection{Exercises}

\begin{exercise}[Genus and total curvature]
	Show that a compact oriented surface of genus $g$ has Euler characteristic $\chi=2-2g$ and deduce
	\[
	\iint_M K\,dA = 2\pi(2-2g).
	\]
\end{exercise}

\begin{exercise}[Flat sphere is impossible]
	Explain why Gauss--Bonnet implies there is no compact surface homeomorphic to $S^2$ with $K\equiv 0$.
\end{exercise}

\begin{exercise}[Where the torus is positively/negatively curved]
	For the torus $K(\phi)=\dfrac{\cos\phi}{r(R+r\cos\phi)}$:
	\begin{enumerate}
		\item Identify explicitly the regions where $K>0$, $K<0$, and $K=0$.
		\item Explain qualitatively why the inner region should have negative curvature.
		\item Verify again that $\iint K\,dA=0$ using the simplified identity $K\,dA=\cos\phi\,d\theta\,d\phi$.
	\end{enumerate}
\end{exercise}

\begin{exercise}[Two independent curvature computations]
	Reproduce the conformal-factor calculations without using the second
	fundamental form.
	\begin{enumerate}
		\item For the round sphere, verify in stereographic coordinates that
		$u=\log 2-\log(1+|z|^2)$ satisfies
		$\Delta_0u=-e^{2u}$, including the calculation at $z=0$ by continuity.
		\item For the torus of revolution, derive the period
		$L=2\pi r/\sqrt{R^2-r^2}$ of the isothermal coordinate $t$ and verify
		that $-\int_0^{2\pi}\int_0^L u_{tt}\,dt\,d\theta=0$.
	\end{enumerate}
\end{exercise}


\subsection{De Rham cohomology, periods, and the two model surfaces}
\label{subsec:forms-derham-models}

Stokes' theorem turns a question about primitives into a question about
integrals.  Write $\Omega^k(M)$ for smooth real $k$-forms and define
\[
 H^k_{\mathrm{dR}}(M)
 =\frac{\ker(d:\Omega^k(M)\to\Omega^{k+1}(M))}
 {\operatorname{im}(d:\Omega^{k-1}(M)\to\Omega^k(M))}.
\]
Thus $H^0_{\mathrm{dR}}(M)=\R$ when $M$ is connected.  A closed one-form
represents zero in $H^1_{\mathrm{dR}}(M)$ precisely when the differential
equation $du=\alpha$ has a globally defined smooth solution.  Local
solutions always exist by the Poincar\'e lemma; the issue is whether they
agree after continuation around loops.

\begin{proposition}[Global solvability of the exact differential equation]
 \label{prop:exact-differential-periods}
 Let $M$ be a connected smooth manifold and $\alpha$ a smooth real
 one-form. The equation $du=\alpha$ has a global smooth solution if
 and only if $d\alpha=0$ and $\int_\gamma\alpha=0$ for every closed
 piecewise smooth loop $\gamma$. The solution is unique up to an
 additive constant.
\end{proposition}

\begin{proof}
 If $\alpha=du$, then $d\alpha=d^2u=0$ and the fundamental theorem
 of calculus gives zero integral on every closed loop. Conversely,
 choose $p_0\in M$ and define $u(p)=\int_{p_0}^{p}\alpha$ along any
 piecewise smooth path. Two paths differ by a closed loop, so the
 period hypothesis makes $u$ well-defined. On a coordinate ball,
 the Poincar\'e lemma gives a smooth local primitive $v$ of $\alpha$.
 Path additivity makes $u-v$ constant there; hence $u$ is smooth and
 $du=\alpha$. Two solutions differ by a function with zero differential,
 which is constant because $M$ is connected.
\end{proof}

\begin{corollary}[Vanishing of first de Rham cohomology]
 \label{cor:h1-zero-exact-differentials}
 For a connected smooth manifold $M$, the following are equivalent:
 $H^1_{\mathrm{dR}}(M)=0$; every smooth closed one-form on $M$ is
 $du$ for some global smooth function $u$; and every smooth closed
 one-form on $M$ has zero integral around every closed loop.
\end{corollary}

\begin{proof}
 The first two assertions are equivalent by the definition
 $H^1_{\mathrm{dR}}(M)=\ker(d:\Omega^1\to\Omega^2)/d\Omega^0$.
 Proposition~\ref{prop:exact-differential-periods} identifies the
 second with the third. The quantifier ``every closed one-form'' is
 essential: zero periods for one form only kill its own class.
\end{proof}

For a compact genus-$g$ surface it suffices to test the $2g$ standard
homology generators: integration of a closed form is additive and
vanishes on boundaries. This is the degree-one form-side content of
the de Rham comparison in Part~III
(Theorem~\ref{thm:deRham-sheaf}); see also
\cite[Chapters~1--2]{BottTu1982}.

For closed forms, the wedge product descends to cohomology:
\[
 [\alpha]\smile[\beta]:=[\alpha\wedge\beta]
 \in H^{p+q}_{\mathrm{dR}}(M),\qquad
 [\alpha]\in H^p_{\mathrm{dR}}(M),\quad
 [\beta]\in H^q_{\mathrm{dR}}(M).
\]
Indeed, replacing $\alpha$ by $\alpha+d\eta$ changes the product by
$d(\eta\wedge\beta)$ when $d\beta=0$; the same argument applies to
$\beta$.  It is graded-commutative because
$\alpha\wedge\beta=(-1)^{pq}\beta\wedge\alpha$.  This is the
\emph{de Rham cup product}; the multiplicative de Rham theorem identifies
it with the singular-cohomology cup product over $\R$
\cite[Chapter~1]{BottTu1982}.

\noindent\textbf{The sphere.}
On $S^2$, every closed one-form is exact: every loop contracts, and
Stokes' theorem makes its period zero.  Hence $H^1_{\mathrm{dR}}(S^2)=0$.
For the outward orientation, spherical coordinates give
$\omega_S=\sin\phi\,d\phi\wedge d\theta$ and
$\int_{S^2}\omega_S=4\pi$.  The form $\omega_S/(4\pi)$ will be the
normalized generator of $H^2_{\mathrm{dR}}(S^2)$.  An explicit
zero-integral example is
\[
 z\omega_S=\cos\phi\sin\phi\,d\phi\wedge d\theta
   =d\!\left(\tfrac12\sin^2\phi\,d\theta\right).
\]
The apparent angular singularity is removable: near a pole,
$\sin^2\phi\,d\theta$ is the restriction of $x\,dy-y\,dx$ in
Cartesian coordinates.

\noindent\textbf{The flat torus.}
Let $T^2=\R^2/\Z^2$ with global forms $dx,dy$ and orientation
$dx\wedge dy$.  Their periods along the horizontal and vertical circles
are $(1,0)$ and $(0,1)$, respectively.  If $\alpha$ is closed and its
two periods are $a,b$, the period criterion applied to
$\alpha-a\,dx-b\,dy$ gives
\[
 H^1_{\mathrm{dR}}(T^2)=\R[dx]\oplus\R[dy],\qquad
 [dx]\smile[dy]=[dx\wedge dy]=-[dy]\smile[dx].
\]
Since $\int_{T^2}dx\wedge dy=1$, this product detects the orientation
class.  The next chapter explains why integration detects \emph{every}
degree-two class, not just these model forms.

\begin{exercise}[Potentials and periods on the torus]
 For $\alpha=\cos(2\pi x)\,dx+3\,dy$, find its two periods and a
 potential for $\alpha-3\,dy$.  Explain why $\alpha$ itself has no
 global potential.
\end{exercise}


\EssayDongsongBae
\ReflectionByungjooJeon
\clearpage

\section{Hodge Theory on Compact Riemann Surfaces}
\index{Hodge theory}
\label{sec:hodge-weyl-compact-rs}

The guiding question is whether the Poisson equation
$\Delta_\omega u=f$ can be solved on a compact surface.  Integrating
both sides shows that $\int_M f\,\omega=0$ is necessary, because the
integral of a Laplacian vanishes.  Constants also prevent uniqueness,
so we seek a solution of mean zero.  The proof proceeds in three
stages: a Poincar\'e inequality makes the energy coercive on the
mean-zero Sobolev space; Lax--Milgram produces a weak solution; and
elliptic regularity upgrades it to a smooth one.  This route explains
why Hilbert-space methods appear in a book about holomorphic geometry.
Later, its harmonic forms will represent cohomology classes and enter
the duality arguments.

\subsection{Analytic framework: metrics, the Hodge star, and weak solutions}
\label{subsec:hodge-standing}

\Step{Standing assumptions and normalizations}
Let $M$ be a compact connected Riemann surface equipped with a smooth Hermitian area form $\omega$.
(Equivalently, $\omega$ is the Riemannian volume form of a smooth metric compatible with the complex structure.)
For real-valued $u,v\in C^\infty(M,\R)$ define
\[
\langle u,v\rangle_{L^2}:=\int_M u\,v\,\omega,
\qquad
\|u\|_{L^2}^2:=\langle u,u\rangle_{L^2}.
\]
Let $\nabla u$ denote the Riemannian gradient and
\[
\|u\|_{W^{1,2}}^2:=\|u\|_{L^2}^2+\|\nabla u\|_{L^2}^2,
\qquad
\|\nabla u\|_{L^2}^2:=\int_M |\nabla u|^2\,\omega.
\]
The Sobolev space $W^{1,2}(M,\R)$ is the completion of
$C^\infty(M,\R)$ in this norm. Equivalently, its elements are
$L^2$ functions whose first derivatives, interpreted weakly in local
coordinates, also lie in $L^2$. We use the completion to make the
energy problem a Hilbert-space problem; regularity later returns the
solution to the smooth category.

\Step{Local coordinate conventions}
In a local holomorphic coordinate $z=x+iy$ the metric may be written
\[
g = 2g_{z\bar z}(z)\,dz\,d\bar z,
\qquad g_{z\bar z}>0,
\qquad
\omega = i\,g_{z\bar z}\,dz\wedge d\bar z = 2g_{z\bar z}\,dx\,dy.
\]
With this convention one has (for real $u$)
\[
|\nabla u|^2\,\omega
=\left(\frac{u_x^2+u_y^2}{2g_{z\bar z}}\right)
  (2g_{z\bar z}\,dx\,dy)
=4|u_z|^2\,dx\,dy
\quad\text{and}\quad
	\Delta_\omega u =-2g^{z\bar z}\,\partial_z\partial_{\bar z}u,
\]
	where we use the nonnegative convention $\Delta_\omega=-\operatorname{div}\nabla$.
	\emph{Remark.} Different texts insert factors $2$ or $4$ in $\Delta_\omega$.
None of the functional-analytic arguments below depends on that choice; only the bookkeeping
of constants does. We fix one convention and keep it throughout.

\begin{definition}[Hodge star and codifferential]
	\label{def:hodge-star}
	Let $(M,g)$ be an oriented Riemannian surface with area form
	$\mathrm dA$.  For each $k=0,1,2$, the \emph{Hodge star} is the unique
	bundle map
	\[
		*:\Omega^k(M)\longrightarrow\Omega^{2-k}(M)
	\]
	characterized pointwise by
	\[
		\alpha\wedge *\beta
		=\langle\alpha,\beta\rangle_g\,\mathrm dA,
		\qquad \alpha,\beta\in\Omega^k(M).
	\]
	The induced $L^2$ inner product is
	$(\alpha,\beta)_{L^2}=\int_M\alpha\wedge *\beta$.  With the orientation
	chosen above,
	\[
		*1=\mathrm dA,\qquad *\mathrm dA=1,
		\qquad *^2=(-1)^{k(2-k)}\ \text{on }\Omega^k(M).
	\]
	The formal adjoint of $d$ on $k$-forms is
	$\delta=-*d*$, and the nonnegative Hodge Laplacian is
	$\Delta_k=d\delta+\delta d$.
\end{definition}

In an oriented conformal coordinate $z=x+iy$, write
$g=\lambda(dx^2+dy^2)$ and $\mathrm dA=\lambda\,dx\wedge dy$.  The
oriented orthonormal coframe
$(\sqrt\lambda\,dx,\sqrt\lambda\,dy)$ gives
\[
	*dx=dy,\qquad *dy=-dx,
	\qquad *dz=-i\,dz,\qquad *d\bar z=i\,d\bar z.
\]
Thus the star on one-forms is conformally invariant.  These formulas fix all
signs used below, including
$d(*du)=-(\Delta_0u)\,\mathrm dA$ for the nonnegative Laplacian convention.

\subsubsection{Coordinate checks for the star, codifferential, and Laplacian}

Because sign conventions vary across the literature, it is useful to derive
the identities once.  In a positively oriented orthonormal coframe
$(e^1,e^2)$,
\[
	*e^1=e^2,\qquad *e^2=-e^1,
	\qquad *(e^1\wedge e^2)=1.
\]
If $\alpha=a e^1+b e^2$, then
\[
	\alpha\wedge *\alpha=(a^2+b^2)e^1\wedge e^2
	=|\alpha|^2\,dA.
\]
This verifies both the defining property and the sign of $*^2=-1$ on
one-forms.  For a function $u$ in a conformal chart,
\[
	du=u_x\,dx+u_y\,dy,
	\qquad *du=u_x\,dy-u_y\,dx,
\]
and therefore
\[
	d(*du)=(u_{xx}+u_{yy})\,dx\wedge dy
	=-(\Delta_0u)\,dA.
\]
The last equality incorporates the conformal factor and our nonnegative
sign convention.  Applying $-* $ gives
\[
	\delta du=-*d*du=\Delta_0u.
\]
Thus the scalar Laplacian is the degree-zero Hodge Laplacian
$\Delta_0=\delta d$.

For a one-form $\alpha=P\,dx+Q\,dy$ in a flat oriented chart,
\[
	d\alpha=(Q_x-P_y)\,dx\wedge dy,
	\qquad
	\delta\alpha=-(P_x+Q_y).
\]
Hence closedness measures curl and coclosedness measures divergence.  A
harmonic one-form is simultaneously curl-free and divergence-free.  This
calculus explains why the two conditions in Hodge theory generalize the
potential and incompressibility conditions familiar from vector calculus.

\begin{lemma}[Integration by parts for forms]
	\label{lem:hodge-integration-by-parts}
	If $M$ is closed, then for smooth forms of compatible degrees,
	\[
		(d\alpha,\beta)_{L^2}=(\alpha,\delta\beta)_{L^2}.
	\]
	Consequently
	\[
		(\Delta_k\alpha,\alpha)_{L^2}
		=\|d\alpha\|_{L^2}^2+\|\delta\alpha\|_{L^2}^2.
	\]
\end{lemma}

\begin{proof}
	The graded Leibniz rule applied to $\alpha\wedge *\beta$ gives a sum
	whose first term is $d\alpha\wedge *\beta$ and whose second term is,
	up to the degree sign, $\alpha\wedge d(*\beta)$.  Stokes' theorem makes
	the integral of the total derivative zero.  Rewriting the second term
	with $\delta=-*d*$ yields the adjoint identity.  Apply it once to
	$d\delta\alpha$ and once to $\delta d\alpha$ to obtain the displayed
	energy formula.
\end{proof}

This lemma immediately identifies the kernel of the Hodge Laplacian:
\[
	\Delta_k\alpha=0
	\quad\Longleftrightarrow\quad
	d\alpha=0\ \text{and}\ \delta\alpha=0.
\]
The forward implication follows from the nonnegative sum of norms; the
reverse implication follows directly from the definition.  No coordinate
calculation is required.

\subsubsection{Weak formulation and the mean-zero Hilbert space}
\label{subsec:weak-setup}

\Step{Mean-zero Sobolev space}
Define the mean-zero subspace
\[
W^{1,2}_\perp(M)
:=\Bigl\{u\in W^{1,2}(M)\ \Big|\ \int_M u\,\omega=0\Bigr\}.
\]
It is a closed subspace of the Hilbert space $W^{1,2}(M)$, hence itself a Hilbert space
with the induced norm $\|\cdot\|_{W^{1,2}}$.

\Step{Compatibility condition (why mean-zero is necessary)}
Let $f\in C^\infty(M,\R)$ satisfy
\[
\int_M f\,\omega=0.
\]
This condition is necessary, because on a closed manifold (no boundary) we have
\[
\int_M \Delta_\omega u\,\omega = 0
\qquad\forall\,u\in C^\infty(M),
\]
so $\Delta_\omega u=f$ forces $\int_M f\,\omega=0$.

\begin{definition}[Weak solution of $\Delta_\omega u=f$]
	\label{def:weak-solution}
	For $f\in L^2(M,\R)$, a function $u\in W^{1,2}(M,\R)$ is a
	\emph{weak solution} of $\Delta_\omega u=f$ if
	\[
	\int_M \langle \nabla u,\nabla\varphi\rangle\,\omega
	=
	\int_M f\,\varphi\,\omega
	\qquad\forall\,\varphi\in C^\infty(M).
	\]
\end{definition}

\Step{Derivation of the weak formulation by integration by parts}
If $u$ is smooth, then on the closed manifold $M$ we have
\[
	\int_M \langle \nabla u,\nabla\varphi\rangle\,\omega
	= \int_M (\Delta_\omega u)\,\varphi\,\omega,
\]
so the weak formulation is exactly the classical equation tested against $\varphi$.

\subsection{Variational existence: Poincar\'e inequality and Lax--Milgram}
\label{subsec:poincare-coercive}

\begin{lemma}[Poincar\'e inequality on $W^{1,2}_\perp(M)$]
	\label{lem:poincare}
	There exists $C_P>0$ such that
	\[
	\|u\|_{L^2}\le C_P\,\|\nabla u\|_{L^2}\qquad \forall\,u\in W^{1,2}_\perp(M).
	\]
	Equivalently, $\|\nabla u\|_{L^2}$ is a norm on $W^{1,2}_\perp(M)$, equivalent to $\|u\|_{W^{1,2}}$.
\end{lemma}

\begin{proof}
	Assume the inequality fails. Then there exists $u_n\in W^{1,2}_\perp(M)$ with
	\[
	\|u_n\|_{L^2}=1
	\quad\text{and}\quad
	\|\nabla u_n\|_{L^2}\to 0.
	\]
	By Rellich--Kondrachov compactness on compact $M$,
	$W^{1,2}(M)\hookrightarrow L^2(M)$ is compact. Hence (after passing to a subsequence)
	\[
	u_n\to u_\ast \ \text{in }L^2
	\quad\text{and}\quad
	u_n\rightharpoonup u_\ast \ \text{weakly in }W^{1,2}.
	\]
	Lower semicontinuity gives
	\[
	\|\nabla u_\ast\|_{L^2}\le \liminf_{n\to\infty}\|\nabla u_n\|_{L^2}=0,
	\]
	so $u_\ast$ is (a.e.) constant. Since each $u_n$ has mean zero,
	$\int_M u_\ast\,\omega=0$, hence $u_\ast\equiv 0$.
	But $\|u_\ast\|_{L^2}=\lim\|u_n\|_{L^2}=1$, contradiction.
\end{proof}

\Step{Coercive bilinear form}
Define
\[
A(u,v):=\int_M \langle \nabla u,\nabla v\rangle\,\omega,
\qquad u,v\in W^{1,2}_\perp(M).
\]
Then $A$ is bounded and, by Lemma~\ref{lem:poincare}, coercive:
\[
A(u,u)=\|\nabla u\|_{L^2}^2
\ge C_P^{-2}\,\|u\|_{L^2}^2,
\]
and hence $A(u,u)\ge c\,\|u\|_{W^{1,2}}^2$ on $W^{1,2}_\perp(M)$ for some $c>0$.

\Step{Bounded linear functional}
For mean-zero $f\in C^\infty(M)$ define
\[
\Lambda_f(\varphi):=\int_M f\,\varphi\,\omega.
\]
Then, using Cauchy--Schwarz and Lemma~\ref{lem:poincare},
\[
|\Lambda_f(\varphi)|
\le \|f\|_{L^2}\,\|\varphi\|_{L^2}
\le C_P\,\|f\|_{L^2}\,\|\nabla\varphi\|_{L^2}
\le C\,\|f\|_{L^2}\,\|\varphi\|_{W^{1,2}}.
\]

\begin{proposition}[Weak existence, uniqueness, and energy estimate]
	\label{prop:weak-existence}
	For every $f\in C^\infty(M,\R)$ with $\int_M f\,\omega=0$ there exists a unique
	$u\in W^{1,2}_\perp(M)$ such that
	\[
	A(u,\varphi)=\Lambda_f(\varphi)
	\qquad\forall\,\varphi\in W^{1,2}_\perp(M).
	\]
	Equivalently, $u$ solves $\Delta_\omega u=f$ in the weak sense.
	Moreover,
	\[
	\|\nabla u\|_{L^2}\le C_P\,\|f\|_{L^2}.
	\]
\end{proposition}

\begin{proof}
	By Lax--Milgram applied to the coercive bilinear form $A$ and the bounded linear functional
	$\Lambda_f$, there exists a unique $u\in W^{1,2}_\perp(M)$ such that
	$A(u,\cdot)=\Lambda_f(\cdot)$.
	This identity also holds for every test function in $W^{1,2}(M)$: write
	$\varphi=\varphi_0+c$ with $\varphi_0$ mean zero. Both sides vanish on the
	constant part because \(\nabla c=0\) and \(\int_M f\,\omega=0\). Thus the
	Lax--Milgram solution satisfies the weak formulation stated in
	Definition~\ref{def:weak-solution}, not merely its mean-zero restriction.
	
	For the estimate, test the weak equation with $\varphi=u$:
	\[
	\|\nabla u\|_{L^2}^2=A(u,u)=\Lambda_f(u)
	\le \|f\|_{L^2}\,\|u\|_{L^2}
	\le C_P\,\|f\|_{L^2}\,\|\nabla u\|_{L^2},
	\]
	hence $\|\nabla u\|_{L^2}\le C_P\,\|f\|_{L^2}$.
\end{proof}

\subsubsection{The variational problem behind the weak equation}

For mean-zero $f$, consider
\[
	\mathcal E_f(v)
	=\frac12\int_M|\nabla v|^2\,\omega
	-\int_Mfv\,\omega,
	\qquad v\in W^{1,2}_\perp(M).
\]
The Poincar\'e inequality prevents a minimizing sequence from escaping by
adding large constants, while coercivity prevents its norm from diverging.
For $v,t\varphi\in W^{1,2}_\perp(M)$,
\[
	\left.\frac{d}{dt}\right|_{t=0}\mathcal E_f(v+t\varphi)
	=A(v,\varphi)-\Lambda_f(\varphi).
\]
Thus critical points are exactly weak solutions.

\begin{proposition}[Strict minimization and stability]
	\label{prop:poisson-energy-minimizer}
	The weak solution $u$ is the unique minimizer of $\mathcal E_f$ on
	$W^{1,2}_\perp(M)$.  Moreover, for every $v$ in that space,
	\[
		\mathcal E_f(v)-\mathcal E_f(u)
		=\frac12\|\nabla(v-u)\|_{L^2}^2.
	\]
\end{proposition}

\begin{proof}
	Write $v=u+w$.  Expanding the square gives
	\[
	\mathcal E_f(u+w)-\mathcal E_f(u)
	=\frac12\|\nabla w\|_{L^2}^2
	+A(u,w)-\Lambda_f(w).
	\]
	The last two terms cancel by the weak equation.  The remaining term is
	nonnegative and vanishes only when $w$ is constant.  Since $w$ has mean
	zero, that constant is zero.
\end{proof}

This identity is stronger than existence: it quantifies the error of any
approximate solution in energy norm.  If a numerical or truncated spectral
approximation $v$ nearly minimizes $\mathcal E_f$, then its gradient is
close to the gradient of the exact normalized solution.

\subsection{Local regularity and the passage from weak to smooth solutions}
\label{subsec:weyl-lemma}

\begin{lemma}[Weyl's lemma (Euclidean version)]
	\label{lem:weyl}
	Let $U\subset\R^2$ be open and $u\in W^{1,2}_{\mathrm{loc}}(U)$ satisfy
	\[
	\int_U u\,\Delta\phi\,dx\,dy=0\qquad \forall\,\phi\in C_c^\infty(U).
	\]
	Then $u\in C^\infty(U)$ and $\Delta u=0$ classically on $U$.
	In particular, $u$ satisfies the mean-value property on balls compactly contained in $U$.
\end{lemma}

\begin{proof}
	Let $\chi\in C_c^\infty(\R^2)$ be a standard mollifier, $\int\chi=1$, and set
	$\chi_\varepsilon(x)=\varepsilon^{-2}\chi(x/\varepsilon)$.
	For $U_\varepsilon:=\{x\in U:\mathrm{dist}(x,\partial U)>\varepsilon\}$ define
	\[
	u_\varepsilon:=\chi_\varepsilon*u\quad \text{on }U_\varepsilon.
	\]
	Then $u_\varepsilon\in C^\infty(U_\varepsilon)$ and, in distributions,
	\[
	\Delta u_\varepsilon=\chi_\varepsilon*(\Delta u)=0\quad \text{on }U_\varepsilon,
	\]
	so $u_\varepsilon$ is classically harmonic.

	We now justify the passage to the limit; pointwise almost-everywhere convergence alone
	would not suffice. Fix concentric balls
	$B_r(x_0)\Subset B_{2r}(x_0)\Subset U$. For sufficiently small
	$\varepsilon,\delta$, the difference
	$h_{\varepsilon,\delta}:=u_\varepsilon-u_\delta$ is harmonic on $B_{2r}(x_0)$.
	Here the needed interior estimate can be seen directly. Choose a smooth
	radial averaging kernel $\mu$ supported in $B_{r/2}(0)$ with integral
	$1$. The mean-value property of the smooth harmonic function
	$h_{\varepsilon,\delta}$ gives
	$h_{\varepsilon,\delta}=\mu*h_{\varepsilon,\delta}$ on
	$B_r(x_0)$. Moving derivatives to $\mu$ yields, for every
	multi-index $\alpha$,
	\[
		\sup_{B_r(x_0)}|D^\alpha h_{\varepsilon,\delta}|
		\le \|D^\alpha\mu\|_\infty\,
		\|h_{\varepsilon,\delta}\|_{L^1(B_{2r}(x_0))}.
	\]
	Because mollification gives $u_\varepsilon\to u$ in
	$L^1_{\mathrm{loc}}(U)$, the right-hand side tends to zero. Thus
	$(u_\varepsilon)$ is Cauchy in $C^\infty(B_r(x_0))$ and converges there to a
	smooth harmonic function $v$. Its $L^1$ limit is also $u$, so $v=u$ almost
	everywhere on $B_r(x_0)$. Replacing the Sobolev representative of $u$ by $v$
	on each such ball produces a globally defined smooth harmonic representative.
	The classical mean-value property then follows.
\end{proof}

\subsubsection{From weak to smooth solutions on a compact Riemann surface}
\label{subsec:elliptic-regularity}

\begin{proposition}[Elliptic regularity on a compact surface]
	\label{prop:regularity}
	Let $f\in C^\infty(M,\R)$ and let $u\in W^{1,2}(M,\R)$ be a weak solution of $\Delta_\omega u=f$.
	Then $u\in C^\infty(M)$.
\end{proposition}

\begin{proof}
	Regularity is local. Fix a holomorphic coordinate chart $U$ with
	$\omega=2g_{z\bar z}\,dx\,dy$, where $g_{z\bar z}$ is smooth and
	positive. In this conformal coordinate, the metric factors cancel from
	the energy integral:
	\[
	\int_U\langle\nabla u,\nabla\varphi\rangle\,\omega
	=\int_U(u_x\varphi_x+u_y\varphi_y)\,dx\,dy.
	\]
	For tests supported in $U$, the weak equation therefore says, in
	distributional notation,
	\[
	\Delta_{\mathrm E}u=-2g_{z\bar z}f=:h,
	\qquad \Delta_{\mathrm E}=\partial_x^2+\partial_y^2.
	\]
	The right-hand side $h$ is smooth. Choose a disk $B\Subset U$ and a
	smooth cutoff $\chi$ supported in $U$ that equals $1$ near $\overline B$.
	Extend $\chi h$ by zero to $\R^2$ and let
	$G(x)=(2\pi)^{-1}\log|x|$. The elementary flux calculation
	$\int_{|x|=r}\partial_rG\,ds=1$ gives
	$\Delta_{\mathrm E}G=\delta_0$ in distributions. Consequently
	$v:=G*(\chi h)$ is smooth and satisfies
	$\Delta_{\mathrm E}v=\chi h$; smoothness follows by moving each
	derivative onto the compactly supported smooth factor $\chi h$.
	On $B$ we have $\Delta_{\mathrm E}(u-v)=0$ weakly. Weyl's lemma
	(Lemma~\ref{lem:weyl}) makes $u-v$ smooth there, so $u$ is smooth
	on $B$. Such disks cover $M$.
\end{proof}

\subsection{Hodge--Weyl theory, Green operators, and de Rham classes}
\index{Hodge--Weyl theorem}\index{Poisson equation!compact surfaces}
\label{subsec:hodge-weyl}

\begin{theorem}[Hodge--Weyl for functions (Poisson equation with mean-zero data)]
	\label{thm:hodge-weyl}
	Let $(M,\omega)$ be a compact connected Riemann surface with Hermitian area form $\omega$.
	For every $f\in C^\infty(M,\R)$ with $\int_M f\,\omega=0$ there exists a unique
	$u\in C^\infty(M,\R)$ with $\int_M u\,\omega=0$ such that
	\[
	\Delta_\omega u=f.
	\]
	Moreover, there exists $C>0$ depending only on $(M,\omega)$ such that
	\[
	\|\nabla u\|_{L^2}\le C\,\|f\|_{L^2}.
	\]
\end{theorem}

\begin{proof}
	Proposition~\ref{prop:weak-existence} gives a unique weak solution $u\in W^{1,2}_\perp(M)$.
	Proposition~\ref{prop:regularity} upgrades it to $u\in C^\infty(M)$.
	The $L^2$-energy estimate was already obtained in Proposition~\ref{prop:weak-existence}.
\end{proof}

\begin{remark}[Complex-valued data]
	For $f\in C^\infty(M,\C)$ with mean zero, apply the theorem separately to
	$\Re f$ and $\Im f$.  Their normalized solutions combine to the unique
	mean-zero complex-valued solution.  Thus no new analytic estimate is needed;
	the real-valued formulation merely matches the real Hilbert space used above.
\end{remark}

\subsubsection{The Green operator and the spectral meaning of the proof}

The preceding argument proves more than the existence of one solution.  It
defines an inverse to the Laplacian after the constants have been removed.
This inverse is the basic analytic object that later appears in potential
theory, Hodge decomposition, and the construction of canonical kernels.

\begin{definition}[Orthogonal projection and Green operator]
	Let $V=\int_M\omega$.  The orthogonal projection onto the constants is
	\[
		\Pi_0f=\frac{1}{V}\int_M f\,\omega.
	\]
	For smooth $f$, the \emph{Green operator} $G$ is the unique mean-zero
	solution of
	\[
		\Delta_\omega Gf=f-\Pi_0f,
		\qquad \int_MGf\,\omega=0.
	\]
\end{definition}

Thus $G$ is not an inverse on all functions: the equation
$\Delta_\omega u=1$ is impossible on a compact surface.  The correct
identities are
\[
	\Delta_\omega G=G\Delta_\omega=I-\Pi_0
\]
on smooth functions, where the second identity uses the mean-zero
normalization.  This formula is worth remembering: it records both the
solution and its obstruction in a single line.

\begin{proposition}[Symmetry and positivity of the Green operator]
	\label{prop:green-symmetric-positive}
	For smooth real-valued functions $f$ and $h$,
	\[
		\langle Gf,h\rangle_{L^2}=\langle f,Gh\rangle_{L^2},
		\qquad
		\langle Gf,f\rangle_{L^2}\ge 0.
	\]
	The second quantity vanishes precisely when $f$ is constant.
\end{proposition}

\begin{proof}
	Set $u=Gf$ and $v=Gh$.  Since $u$ and $v$ have mean zero, constants
	disappear from their pairings.  Integration by parts gives
	\[
		\langle Gf,h\rangle
		=\langle u,\Delta_\omega v\rangle
		=\int_M\langle\nabla u,\nabla v\rangle\,\omega
		=\langle\Delta_\omega u,v\rangle
		=\langle f,Gh\rangle.
	\]
	Taking $h=f$ yields
	\[
		\langle Gf,f\rangle
		=\int_M|\nabla Gf|^2\,\omega\ge0.
	\]
	Equality forces $Gf$ to be constant and hence zero by its normalization.
	The defining equation then says $f=\Pi_0f$.
\end{proof}

The compactness used in the Poincar\'e inequality has a spectral
interpretation.  There is an $L^2$-orthonormal basis
$(\varphi_j)_{j\ge0}$ of smooth real eigenfunctions satisfying
\[
	\Delta_\omega\varphi_j=\lambda_j\varphi_j,
	\qquad
	0=\lambda_0<\lambda_1\le\lambda_2\le\cdots,
	\qquad \lambda_j\longrightarrow\infty,
\]
with $\varphi_0=V^{-1/2}$.  If $f=\sum_j f_j\varphi_j$, then
\[
	\Pi_0f=f_0\varphi_0,
	\qquad
	Gf=\sum_{j\ge1}\frac{f_j}{\lambda_j}\varphi_j.
\]
The compatibility condition is exactly $f_0=0$.  Moreover,
\[
	\|\nabla Gf\|_{L^2}^2
	=\sum_{j\ge1}\frac{|f_j|^2}{\lambda_j}
	\le\frac{1}{\lambda_1}\|f-\Pi_0f\|_{L^2}^2.
\]
Hence the best Poincar\'e constant is $C_P=\lambda_1^{-1/2}$.  The
variational proof does not require the full spectral theorem, but the
spectral formula makes the geometry of its constant transparent.

\begin{proposition}[Fredholm alternative for the scalar Laplacian]
	\label{prop:fredholm-scalar-laplacian}
	For $f\in C^\infty(M)$, the equation $\Delta_\omega u=f$ is solvable
	if and only if $f$ is $L^2$-orthogonal to $\ker\Delta_\omega$.
	Since $M$ is connected, this kernel is the constants, so the condition is
	$\int_Mf\,\omega=0$.  When it is solvable, all solutions are
	$Gf+c$ with $c$ constant.
\end{proposition}

\begin{proof}
	Self-adjointness gives
	$\langle\Delta_\omega u,h\rangle=\langle u,\Delta_\omega h\rangle=0$
	for every harmonic $h$, proving necessity.  Conversely, orthogonality to
	the constants is the mean-zero condition, and
	Theorem~\ref{thm:hodge-weyl} supplies $Gf$.  If $u_1,u_2$ solve the
	equation, then $\Delta_\omega(u_1-u_2)=0$; testing against the difference
	shows that its gradient vanishes, so it is constant.
\end{proof}

\subsubsection{A stochastic representation of the Poisson solution}
\index{Poisson equation!stochastic representation}
\index{Brownian motion!on a Riemann surface}

The variational and spectral arguments invert the Laplacian by functional
analysis.  There is also a probabilistic interpretation: the Green operator
adds up the expected values of the data along Brownian paths.  This viewpoint
explains simultaneously why constants are an obstruction, why the heat flow
appears, and why subtracting the spatial mean is unavoidable on a compact
surface.  For the construction of Brownian motion and It\^o calculus on
Riemannian manifolds, see \cite[Chapters~2--3]{Hsu2002}.

\begin{definition}[Brownian motion and heat semigroup]
	Let $(X_t)_{t\ge0}$ be Brownian motion on $(M,g)$ started at $x$.  We
	choose its generator to be $-\Delta_\omega$; equivalently, for smooth
	$\phi$,
	\[
		\left.\frac{d}{dt}\right|_{t=0}
		\mathbb E_x[\phi(X_t)]=-\Delta_\omega\phi(x).
	\]
	The associated heat semigroup is
	\[
		P_t\phi(x)=\mathbb E_x[\phi(X_t)].
	\]
	Thus $v(t,x)=P_t\phi(x)$ solves
	$\partial_tv=-\Delta_\omega v$ with $v(0,\cdot)=\phi$.
\end{definition}

Some probability texts define Brownian motion with generator
$\tfrac12\operatorname{div}\nabla=-\tfrac12\Delta_\omega$.  In that
normalization every time parameter below is multiplied by $2$.  The geometry
and the solvability condition are unchanged.

\begin{proposition}[Dynkin's formula in the present sign convention]
	\label{prop:dynkin-compact-surface}
	For $u\in C^\infty(M)$ and $T>0$,
	\[
		\mathbb E_x[u(X_T)]-u(x)
		=-\mathbb E_x\!\left[\int_0^T
		(\Delta_\omega u)(X_s)\,ds\right].
	\]
	In particular, if $\Delta_\omega u=f$, then
	\[
		u(x)=P_Tu(x)+\int_0^TP_sf(x)\,ds.
	\]
\end{proposition}

\begin{proof}
	It\^o's formula applied in local coordinates, with the metric drift that
	makes the construction coordinate invariant, says that
	\[
		u(X_T)-u(X_0)+\int_0^T(\Delta_\omega u)(X_s)\,ds
	\]
	is a martingale.  Compactness makes every smooth coefficient bounded, so
	taking expectations is legitimate.  Since $X_0=x$, the first identity
	follows.  Substitute $\Delta_\omega u=f$ and rearrange to obtain the
	second.
\end{proof}

On a compact connected surface the heat flow converges to equilibrium:
\[
	P_tf\longrightarrow\Pi_0f
	=\frac{1}{\operatorname{Vol}(M)}\int_Mf\,\omega.
\]
Therefore the naive occupation integral $\int_0^\infty P_tf\,dt$ diverges
when $f$ has nonzero mean.  The correct expression removes the equilibrium
mode before integrating.

\begin{theorem}[Stochastic/heat-semigroup formula for the Green operator]
	\label{thm:stochastic-green-representation}
	For $f\in C^\infty(M)$, define
	\[
		u(x)=\int_0^\infty
	\bigl(P_tf(x)-\Pi_0f\bigr)\,dt
	=\int_0^\infty
	\left(\mathbb E_x[f(X_t)]
	-\frac{1}{\operatorname{Vol}(M)}\int_Mf\,\omega\right)dt.
	\]
	Then the integral converges smoothly, $\int_Mu\,\omega=0$, and
	\[
		\Delta_\omega u=f-\Pi_0f.
	\]
	Consequently $u=Gf$.  For mean-zero $f$, this becomes the precise
	finite-time occupation formula
	\[
		Gf(x)=\lim_{T\to\infty}
		\mathbb E_x\!\left[\int_0^T f(X_t)\,dt\right].
	\]
	The expectations converge although the unrenormalized random occupation
	integrals need not converge pathwise.
\end{theorem}

\begin{proof}
	Expand $f=\sum_{j\ge0}f_j\varphi_j$ in Laplace eigenfunctions.  Since
	$P_t\varphi_j=e^{-\lambda_jt}\varphi_j$, subtracting $\Pi_0f$ removes
	the $j=0$ term and gives
	\[
		P_tf-\Pi_0f
		=\sum_{j\ge1}e^{-\lambda_jt}f_j\varphi_j.
	\]
	The spectral gap $\lambda_1>0$ gives exponential decay.  Smoothness of
	$f$ gives rapid decay of $f_j$, so the series and all of its derivatives
	may be integrated term by term.  Hence
	\[
		u=\sum_{j\ge1}\frac{f_j}{\lambda_j}\varphi_j=Gf.
	\]
	The defining properties of $G$ now give the asserted equation and
	normalization.  Alternatively,
	\[
		\Delta_\omega u
	=-\int_0^\infty\partial_t(P_tf)\,dt
	=f-\lim_{t\to\infty}P_tf=f-\Pi_0f.
	\]
\end{proof}

\begin{remark}[What the formula does and does not say]
	On a compact surface Brownian motion does not escape to infinity; it keeps
	returning and approaches the uniform equilibrium distribution.  Therefore
	the solution is not the expected total occupation of a nonnegative source.
	Positive and negative parts of a mean-zero source cancel after mixing.  The
	subtraction of $\Pi_0f$ is the probabilistic form of the Fredholm
	compatibility condition.
\end{remark}

\begin{definition}[Centered occupation potential]
	For $f\in C^\infty(M)$, its \emph{centered occupation potential} is
	\[
		\mathcal U_f(x)
		:=\int_0^\infty
		\mathbb E_x\bigl[f(X_t)-\Pi_0f\bigr],dt.
	\]
	The spectral gap makes the integral converge.  Theorem
	\ref{thm:stochastic-green-representation} states that
	$\mathcal U_f=Gf$ and hence
	\[
		\Delta_\omega\mathcal U_f=f-\Pi_0f.
	\]
\end{definition}

The definition gives an algorithm for solving the Poisson equation
probabilistically:
\begin{enumerate}
	\item compute the spatial mean $\Pi_0f$;
	\item evolve the centered function by Brownian averaging,
	      $P_t(f-\Pi_0f)(x)=\mathbb E_x[f(X_t)-\Pi_0f]$;
	\item integrate this expectation over time;
	\item differentiate, using $\partial_tP_t=-\Delta_\omega P_t$, to
	      verify the equation.
\end{enumerate}
If $f$ is already an eigenfunction, the procedure reduces to a one-line
integral.  The next two examples carry out every step.

\begin{example}[Heat-semigroup inversion on the flat torus]
	\label{ex:stochastic-poisson-torus}
	Let
	\[
		f(x,y)=\sin(2\pi x)\sin(4\pi y)
	\]
	on $T^2=\R^2/\Z^2$, with
	$\Delta=-(\partial_x^2+\partial_y^2)$.  Each sine has zero integral over
	$[0,1]$, so $\Pi_0f=0$.  Direct differentiation gives
	\begin{align*}
	 -\partial_x^2f&=4\pi^2f,\\
	 -\partial_y^2f&=16\pi^2f,
	\end{align*}
	and therefore
	\[
		\Delta f=20\pi^2f.
	\]
	The heat equation preserves this one-dimensional eigenspace, so
	\[
		P_tf(x,y)=\mathbb E_{(x,y)}[f(X_t)]
		=e^{-20\pi^2t}f(x,y).
	\]
	Its occupation potential is
	\begin{align*}
		u(x,y)
		&=\int_0^\infty P_tf(x,y)\,dt\\
		&=f(x,y)\int_0^\infty e^{-20\pi^2t}\,dt\\
		&=\frac{1}{20\pi^2}\sin(2\pi x)\sin(4\pi y).
	\end{align*}
	It has mean zero, and a final differentiation verifies
	\[
		\Delta u=\frac{1}{20\pi^2}\Delta f=f.
	\]
	Thus Brownian averaging, Fourier inversion, and the normalized Green
	operator give exactly the same solution.
\end{example}

\begin{example}[Brownian averaging of a quadratic mode on the sphere]
	\label{ex:stochastic-poisson-sphere}
	Let $q=3z^2-1$ on the unit round sphere.  Rotational symmetry gives
	$\int_{S^2}z^2\,dA=\frac13\int_{S^2}dA$, so $q$ has mean zero.  It is a
	spherical harmonic of degree $2$ and hence
	\[
		\Delta q=2(2+1)q=6q.
	\]
	The heat semigroup acts diagonally on spherical harmonics.  For Brownian
	motion started at $p$,
	\[
		\mathbb E_p[q(X_t)]=P_tq(p)=e^{-6t}q(p).
	\]
	Consequently
	\begin{align*}
		u(p)
		&=\int_0^\infty\mathbb E_p[q(X_t)]\,dt\\
		&=q(p)\int_0^\infty e^{-6t}\,dt
		=\frac16q(p)=\frac16\bigl(3z(p)^2-1\bigr).
	\end{align*}
	The mean of $u$ is zero, and
	$\Delta u=\frac16\Delta q=q$.  This is the normalized solution of the
	Poisson equation.  The exponential factor $e^{-6t}$ has a geometric
	meaning: higher spherical modes are forgotten more rapidly by Brownian
	averaging.
\end{example}

\begin{proposition}[Stopped-process formula on a domain]
	\label{prop:stochastic-dirichlet-poisson}
	Let $D\Subset M$ have smooth boundary, let
	$\tau_D=\inf\{t\ge0:X_t\notin D\}$, and suppose
	\[
		\Delta_\omega u=f\quad\text{in }D,
		\qquad u=h\quad\text{on }\partial D.
	\]
	Then
	\[
		u(x)=\mathbb E_x[h(X_{\tau_D})]
	+\mathbb E_x\!\left[\int_0^{\tau_D}f(X_s)\,ds\right].
	\]
\end{proposition}

\begin{proof}
	Apply Dynkin's formula to the bounded stopping time
	$T\wedge\tau_D$ and use $\Delta_\omega u=f$.  Let $T\to\infty$.
	For a relatively compact smooth domain, $\tau_D<\infty$ almost surely;
	continuity gives $u(X_{\tau_D})=h(X_{\tau_D})$.  Boundedness of the
	data and standard exit-time estimates justify dominated convergence.
\end{proof}

This last formula separates two effects: Brownian exit distributes the
boundary values harmonically, while the occupation integral accumulates the
interior source.  On the closed surface there is no exit boundary, and the
equilibrium subtraction in
Theorem~\ref{thm:stochastic-green-representation} replaces the stopping
mechanism.

\subsubsection{Hodge decomposition on functions and differential forms}

For functions, the decomposition is already visible:
\[
	C^\infty(M)=\R\oplus\Delta_\omega C^\infty(M),
	\qquad f=\Pi_0f+\Delta_\omega Gf.
\]
The analogous statement for forms is the organizing principle of Hodge
theory.  Denote by
\[
	\mathcal H^k(M)=\{\alpha\in\Omega^k(M):\Delta_k\alpha=0\}
\]
the space of harmonic $k$-forms.

\begin{theorem}[Hodge decomposition on a compact oriented surface]
	\label{thm:hodge-decomposition-surface}
	Every smooth $k$-form admits an $L^2$-orthogonal decomposition
	\[
		\Omega^k(M)
		=\mathcal H^k(M)\oplus d\Omega^{k-1}(M)
		\oplus\delta\Omega^{k+1}(M).
	\]
	Consequently every de Rham class has a unique harmonic representative,
	and the map $\mathcal H^k(M)\to H^k_{\mathrm{dR}}(M)$ is an
	isomorphism.
\end{theorem}

The analytic proof for general $k$ repeats the scalar strategy for the
elliptic operator $\Delta_k$: first invert it on the orthogonal complement
of its finite-dimensional kernel, then write
\[
	\alpha=H\alpha+\Delta_kG_k\alpha
	=H\alpha+d(\delta G_k\alpha)+\delta(dG_k\alpha).
\]
Here $H$ is orthogonal projection onto $\mathcal H^k(M)$.  If $\alpha$ is
closed, the coexact summand vanishes by orthogonality, so
$\alpha=H\alpha+d(\delta G_k\alpha)$.  This both produces a harmonic
representative and proves uniqueness.  We use this theorem below only in
the two explicit model surfaces; later chapters return to its
cohomological consequences.

\subsubsection{Top-degree classes and the Poisson equation}
\label{subsec:hodge-top-degree-poisson}
\label{subsec:poisson-top-degree}

Let $M$ be a compact connected Riemann surface with the preceding
orientation and area form $\omega$.  Every smooth two-form has a unique
expression $\eta=f\omega$.  In an oriented orthonormal coframe,
$*du=u_x\,dy-u_y\,dx$; with our nonnegative Laplacian convention this
gives the identity
\[
 d(*du)=-(\Delta_\omega u)\omega.
\]
It can also be checked in a conformal chart, where the metric factor
cancels in $*du$.  Consequently the obstruction to solving
$d\beta=\eta$ is exactly the obstruction to solving a Poisson equation.

\begin{proposition}[Poisson solvability computes $H^2_{\mathrm{dR}}$]
 \label{prop:poisson-computes-h2}
 On a compact connected Riemann surface $M$, a smooth real two-form
 $\eta$ is exact if and only if $\int_M\eta=0$. Consequently
 integration induces an isomorphism
 \[
  H^2_{\mathrm{dR}}(M)\longrightarrow\R,\qquad
  [\eta]\longmapsto\int_M\eta,
 \]
 whose normalized generator is $[\omega/\int_M\omega]$.
\end{proposition}

\begin{proof}
 Because $M$ has real dimension two, $\Omega^3(M)=0$: every two-form
 is closed. Thus the definition reads
 $H^2_{\mathrm{dR}}(M)=\Omega^2(M)/d\Omega^1(M)$, and deciding whether
 a two-form is exact is precisely deciding whether its class vanishes.
 If $\eta=d\beta$, Stokes' theorem gives $\int_M\eta=0$.
 Conversely, write $\eta=f\omega$. The integral hypothesis is the
 compatibility condition $\int_M f\omega=0$ in
 Theorem~\ref{thm:hodge-weyl}; that theorem supplies a smooth $u$
 satisfying $\Delta_\omega u=f$. The identity
 $d(*du)=-(\Delta_\omega u)\omega$ now gives
 $d(-*du)=f\omega=\eta$. Hence the integration map is well-defined
 and injective on cohomology. It is surjective because the class of
 $c\omega/\int_M\omega$ maps to any prescribed $c\in\R$.
\end{proof}

For complex-valued forms, the same conclusion with target $\C$
follows by applying the proof to real and imaginary parts.

Together with the period criterion, this gives a usable description of
all three de Rham groups of a genus-$g$ compact Riemann surface:
\[
 H^0_{\mathrm{dR}}(M)\cong\R,\qquad
 H^1_{\mathrm{dR}}(M)\cong\R^{2g},\qquad
 H^2_{\mathrm{dR}}(M)\cong\R.
\]
For the middle dimension, a standard genus-$g$ cell decomposition has
one vertex, $2g$ oriented one-cells, and one two-cell attached along the
product of commutators.  Its cellular boundary maps vanish after
abelianization; the de Rham comparison identifies the resulting real
cohomology with the forms above (Theorem~\ref{thm:deRham-sheaf};
\cite[Chapters~1--2]{BottTu1982}).  In particular, exactness of a
closed one-form is tested by its $2g$ basic periods.  The bilinear pairing
$([\alpha],[\beta])\mapsto\int_M\alpha\wedge\beta$ is the cup product
followed by integration; its nondegeneracy is the surface case of
Poincar\'e duality discussed in Part~III.

On the flat torus, this top-degree criterion has an explicit Fourier
version.  For a smooth mean-zero $f=\sum_{(m,n)\ne(0,0)}c_{mn}
e^{2\pi i(mx+ny)}$, set
\[
 u=\sum_{(m,n)\ne(0,0)}
 \frac{c_{mn}}{4\pi^2(m^2+n^2)}e^{2\pi i(mx+ny)},
 \qquad \beta=-*du.
\]
Rapid decay of the Fourier coefficients makes $u$ smooth and gives
$d\beta=f\,dx\wedge dy$.  For $f=\cos(2\pi x)$ specifically,
$u=\cos(2\pi x)/(4\pi^2)$ and
$\beta=\sin(2\pi x)\,dy/(2\pi)$.
On the unit sphere, $\Delta z=2z$, so the primitive of $z\omega_S$
computed above is also $-*d(z/2)$.

\subsection{Solving the one-dimensional Calabi--Yau equation on the sphere and the flat torus}
\label{subsec:hodge-examples}

Proposition~\ref{prop:calabi-yau-riemann-surface-poisson} says that, in
complex dimension one, prescribing a volume form in a fixed K\"ahler class is
not a nonlinear determinant problem.  It is the Poisson equation
\[
 \Delta_{g_0}\varphi=2(e^F-1),
 \qquad
 \omega_0+i\partial\bar\partial\varphi=e^F\omega_0,
\]
where in this subsection $\Delta_{g_0}=\operatorname{div}_{g_0}\nabla$.
This is the sign convention used in Proposition~\ref{prop:calabi-yau-riemann-surface-poisson};
its nonconstant eigenvalues are negative.  The normalization
$\int_Xe^F\omega_0=\int_X\omega_0$ removes the constant mode of the
right-hand side.  The following two calculations solve the equation completely
on the model surfaces and also show what positivity and uniqueness mean.

\subsubsection{The round sphere: a first spherical harmonic}

Let $S^2\subset\R^3$ carry its unit round metric and K\"ahler form
$\omega_S=dA$.  Write $z$ for the height coordinate and choose
\[
 e^F=1+\varepsilon z,
 \qquad |\varepsilon|<1.
\]
The inequality on $\varepsilon$ makes $e^F$ positive everywhere.  Moreover,
$z$ is odd under the antipodal map, so
\[
 \int_{S^2}e^F\omega_S
 =\int_{S^2}\omega_S+\varepsilon\int_{S^2}z\,dA
 =4\pi.
\]
Thus the volume normalization holds.  Since $z$ is a spherical harmonic of
degree one,
\[
 \Delta_{S^2}z=-2z.
\]
The Calabi--Yau equation therefore becomes
\[
 \Delta_{S^2}\varphi=2\varepsilon z,
\]
and its mean-zero solution is
\[
 \varphi=-\varepsilon z.
\]
Indeed, $\Delta_{S^2}(-\varepsilon z)=2\varepsilon z$.  Using the
dimension-one identity
\[
 i\partial\bar\partial\varphi
 =\frac{1}{2}\bigl(\Delta_{S^2}\varphi\bigr)\omega_S
\]
gives the direct verification
\[
 \omega_S+i\partial\bar\partial\varphi
 =(1+\varepsilon z)\omega_S=e^F\omega_S.
\]
The deformed form is positive because its coefficient ranges between
$1-|\varepsilon|$ and $1+|\varepsilon|$.  If another potential gives the
same form, their difference has zero Laplacian and is therefore constant.

In spherical coordinates $z=\cos\theta$ and
$\omega_S=\sin\theta\,d\theta\wedge d\phi$.  Hence the same verification
can be made without spectral terminology:
\[
 \Delta_{S^2}(-\varepsilon\cos\theta)
 =\frac{1}{\sin\theta}\frac{\partial}{\partial\theta}
 \left(\sin\theta\frac{\partial}{\partial\theta}
 (-\varepsilon\cos\theta)\right)
 =2\varepsilon\cos\theta.
\]
Although the coordinate formula contains $1/\sin\theta$, the potential is
the restriction of a smooth linear function on $\R^3$ and hence is smooth at
both poles.

\subsubsection{The flat torus: inversion of one Fourier mode}

Let $T^2=\R^2/\Z^2$ with complex coordinate $w=x+iy$, flat metric
$g_0=dx^2+dy^2$, and K\"ahler form $\omega_0=dx\wedge dy$.  Choose
\[
 e^F=1+\varepsilon\cos(2\pi x)\cos(2\pi y),
 \qquad |\varepsilon|<1.
\]
Both cosine factors have zero average, and therefore
\[
 \int_{T^2}e^F\omega_0=1=\int_{T^2}\omega_0.
\]
For the geometric Laplacian
$\Delta_{g_0}=\partial_x^2+\partial_y^2$,
\[
 \Delta_{g_0}\bigl(\cos(2\pi x)\cos(2\pi y)\bigr)
 =-8\pi^2\cos(2\pi x)\cos(2\pi y).
\]
Consequently the normalized solution of
$\Delta_{g_0}\varphi=2(e^F-1)$ is
\[
 \varphi(x,y)
 =-\frac{\varepsilon}{4\pi^2}
 \cos(2\pi x)\cos(2\pi y).
\]
Substitution yields
\[
 \Delta_{g_0}\varphi
 =2\varepsilon\cos(2\pi x)\cos(2\pi y),
\]
and hence
\[
 \omega_0+i\partial\bar\partial\varphi
 =\left(1+\varepsilon\cos(2\pi x)\cos(2\pi y)\right)\omega_0
 =e^F\omega_0.
\]
Again $|\varepsilon|<1$ is exactly the transparent sufficient condition for
positivity, while adding a constant to $\varphi$ changes neither
$i\partial\bar\partial\varphi$ nor the resulting K\"ahler form.

More generally, write a positive normalized density on the flat torus as
\[
 e^F=1+\sum_{(m,n)\ne(0,0)}a_{mn}e^{2\pi i(mx+ny)},
 \qquad a_{-m,-n}=\overline{a_{mn}}.
\]
The missing zero Fourier coefficient is precisely the volume constraint.
Mode-by-mode inversion gives
\[
 \varphi
 =-\sum_{(m,n)\ne(0,0)}
 \frac{a_{mn}}{2\pi^2(m^2+n^2)}e^{2\pi i(mx+ny)},
\]
which has mean zero and satisfies
$\Delta_{g_0}\varphi=2(e^F-1)$.  Smoothness of $F$ gives rapid decay of
the Fourier coefficients, so the series may be differentiated term by term.

\subsubsection{Stochastic reconstruction of the two potentials}

The stochastic formula of
Theorem~\ref{thm:stochastic-green-representation} uses the nonnegative
Laplacian
\[
 \Delta_+:=-\operatorname{div}_{g_0}\nabla=-\Delta_{g_0}
\]
and Brownian motion whose generator is $-\Delta_+$.  If
$h=2(e^F-1)$, the Calabi--Yau equation
$\Delta_{g_0}\varphi=h$ is therefore equivalent to
\[
 \Delta_+\varphi=-h.
\]
The volume normalization says that $h$ has mean zero.  Consequently the
mean-zero potential has the stochastic representation
\begin{equation}
 \label{eq:stochastic-calabi-yau-potential}
 \varphi(p)
 =-\int_0^\infty P_th(p)\,dt
 =-\int_0^\infty\mathbb E_p[h(X_t)]\,dt.
\end{equation}
This formula is an identity of smooth functions.  It records the cancellation
of the positive and negative parts of $h$ after Brownian averaging; it does
not assert that an unsigned occupation time is finite.

On the round sphere, $h=2\varepsilon z$ and
$\Delta_+z=2z$.  Brownian averaging therefore gives
\[
 \mathbb E_p[z(X_t)]=P_tz(p)=e^{-2t}z(p).
\]
Here is the calculation behind that formula.  For fixed $p\in S^2$, set
\[
 m_p(t):=\mathbb E_p[z(X_t)]=(P_tz)(p).
\]
The backward heat equation $\partial_tP_t=-\Delta_+P_t$ and the eigenvalue
identity $\Delta_+z=2z$ give the scalar initial-value problem
\[
 m_p'(t)=-2m_p(t),
 \qquad m_p(0)=z(p).
\]
Solving this ordinary differential equation gives
$m_p(t)=e^{-2t}z(p)$.  Probabilistically, a path that starts high on the
sphere retains a positive expected height for a while, but rotational
mixing drives that expectation exponentially toward the equilibrium mean
$0$.  No sample path itself follows the deterministic curve
$e^{-2t}z(p)$; the formula describes the average over all paths.

The same conclusion follows from the spherical heat kernel.  If
$\{Y_{\ell m}\}$ is an orthonormal basis of spherical harmonics, then
\[
 p_t(p,q)=\sum_{\ell=0}^{\infty}\sum_{m=-\ell}^{\ell}
 e^{-\ell(\ell+1)t}Y_{\ell m}(p)\overline{Y_{\ell m}(q)},
 \qquad
 P_tf(p)=\int_{S^2}p_t(p,q)f(q)\,dA(q).
\]
Since $z$ lies entirely in the degree-one eigenspace, orthogonality removes
every term except $\ell=1$, leaving the factor $e^{-2t}$.
Substitution into \eqref{eq:stochastic-calabi-yau-potential} yields
\begin{align*}
 \varphi(p)
 &=-2\varepsilon\int_0^\infty
   \mathbb E_p[z(X_t)]\,dt\\
 &=-2\varepsilon z(p)\int_0^\infty e^{-2t}\,dt
 =-\varepsilon z(p).
\end{align*}
Thus the probabilistic construction recovers exactly the potential obtained
from the degree-one spherical harmonic.

On the flat torus, put
$c(x,y)=\cos(2\pi x)\cos(2\pi y)$.  Here
$h=2\varepsilon c$ and $\Delta_+c=8\pi^2c$, so for Brownian motion started
at $(x,y)$,
\[
 \mathbb E_{(x,y)}[c(X_t)]=e^{-8\pi^2t}c(x,y).
\]
We can compute this expectation directly.  With the normalization used
here, a lift of the torus Brownian motion to $\R^2$ is
\[
 \widetilde X_t=(x,y)+\sqrt2\,(B_t^1,B_t^2),
\]
where $B^1$ and $B^2$ are independent standard one-dimensional Brownian
motions; reducing both coordinates modulo $\Z$ gives $X_t$.  Periodicity
allows us to evaluate $c$ on the lift.  The Gaussian characteristic-function
identity
\[
 \mathbb E[e^{iaB_t}]=e^{-a^2t/2}
\]
implies
\begin{align*}
 \mathbb E_x[\cos(2\pi(x+\sqrt2 B_t^1))]
 &=\cos(2\pi x)e^{-4\pi^2t},\\
 \mathbb E_y[\cos(2\pi(y+\sqrt2 B_t^2))]
 &=\cos(2\pi y)e^{-4\pi^2t}.
\end{align*}
The sine terms have zero expectation by symmetry.  Independence of the two
Brownian coordinates then gives
\begin{align*}
 P_tc(x,y)
 &=\mathbb E_{(x,y)}[c(X_t)]\\
 &=e^{-4\pi^2t}\cos(2\pi x),
   e^{-4\pi^2t}\cos(2\pi y)\\
 &=e^{-8\pi^2t}c(x,y).
\end{align*}
Equivalently, the torus heat kernel is
\[
 p_t((x,y),(u,v))
 =\sum_{(m,n)\in\Z^2}e^{-4\pi^2(m^2+n^2)t}
 e^{2\pi i(m(x-u)+n(y-v))},
\]
and integrating $c(u,v)$ against this kernel selects precisely the four
Fourier modes $(m,n)=(\pm1,\pm1)$, all with eigenvalue $8\pi^2$.
Hence
\begin{align*}
 \varphi(x,y)
 &=-2\varepsilon c(x,y)
   \int_0^\infty e^{-8\pi^2t}\,dt\\
 &=-\frac{\varepsilon}{4\pi^2}
   \cos(2\pi x)\cos(2\pi y),
\end{align*}
again agreeing with Fourier inversion.  These calculations make the common
mechanism visible: Brownian averaging multiplies an eigenmode of eigenvalue
$\lambda>0$ by $e^{-\lambda t}$, and integration in time divides that mode
by $\lambda$.

\paragraph{Deterministic and stochastic inversion compared.}
The two methods do not produce different solutions.  They are two
realizations of the same inverse operator on the mean-zero subspace:
\[
 \Delta_+^{-1}f
 =\int_0^\infty e^{-t\Delta_+}f\,dt
 =\int_0^\infty P_tf\,dt.
\]
Their working languages and their geometric interpretations are nevertheless
different.  In the deterministic method, one first expands the given
function into Laplace eigenfunctions or Fourier modes.  Solving
$\Delta_+u=f$ then means dividing the coefficient of every nonconstant mode
by its positive eigenvalue.  In the stochastic method, one instead follows
Brownian paths issuing from $x$, averages the observed values of $f$ to form
$P_tf(x)=\mathbb E_x[f(X_t)]$, and integrates that average over time.  The
first method treats the inverse Laplacian as spectral division; the second
treats it as the total expected effect that remains while Brownian motion
forgets its starting point.

The mean-zero hypothesis expresses the same obstruction in these two
languages.  Deterministically, constants form the zero eigenspace of
$\Delta_+$, so their coefficients cannot be divided by an eigenvalue.
Stochastically, the heat flow does not make a constant decay: the equilibrium
average remains present for all time, and its integral over $[0,\infty)$
therefore diverges.  Subtracting the spatial mean removes precisely this
constant component.  Every remaining eigenmode has an eigenvalue
$\lambda>0$ and decays as $e^{-\lambda t}$; the identity
\[
 \int_0^\infty e^{-\lambda t}\,dt=\frac{1}{\lambda}
\]
shows directly why time integration in the stochastic construction performs
the same division as deterministic spectral inversion.

The sphere makes this equivalence especially transparent.  Since
$\Delta_+z=2z$, the deterministic method gives
$\Delta_+^{-1}z=z/2$.  Brownian averaging gives
$P_tz=e^{-2t}z$, and integrating in time again gives
\[
 \int_0^\infty P_tz\,dt
 =z\int_0^\infty e^{-2t}\,dt
 =\frac{z}{2}.
\]
On the flat torus, the function
$\cos(2\pi x)\cos(2\pi y)$ is built from the four Fourier modes
$(\pm1,\pm1)$, each with eigenvalue $8\pi^2$.  Fourier inversion therefore
divides their coefficients by $8\pi^2$.  Brownian averaging multiplies the
same modes by $e^{-8\pi^2t}$, and the time integral produces
$1/(8\pi^2)$.  Thus both model computations agree term by term, not merely
after the final answer has been assembled.

The deterministic approach is usually shorter when an explicit spectral
basis is available, and it naturally emphasizes eigenfunctions, elliptic
estimates, and uniqueness.  The stochastic approach gives the spectral
denominators a dynamical meaning: a mode with a large eigenvalue loses memory
rapidly and therefore contributes for a shorter time to the occupation
potential.  It also extends naturally to boundary-value problems by stopping
Brownian motion when it first reaches the boundary.  An expectation formula,
however, does not by itself prove smoothness, convergence, or uniqueness;
those conclusions still rely on the heat kernel, the spectral gap, or
elliptic regularity.  Probability therefore supplies a geometric
interpretation and an alternative construction, while deterministic analysis
provides the regularity and functional-analytic control that justify it.

\paragraph{What the two calculations establish.}
On the sphere the inverse Laplacian is diagonalized by spherical harmonics;
on the flat torus it is diagonalized by Fourier modes.  In both cases the
constant mode is the sole obstruction, the volume normalization removes it,
and division by the corresponding nonzero eigenvalue constructs the potential.
The examples therefore display every ingredient of the dimension-one
Calabi--Yau theorem: cohomological compatibility, explicit Poisson inversion,
positivity of the prescribed form, and uniqueness modulo constants.

\subsection{Exercises}
\label{subsec:hodge-exercises}

\Step{A. Compatibility and uniqueness}

\begin{exercise}[Compatibility on a closed surface]
	Let $(M,\omega)$ be compact. Show directly that $\int_M \Delta_\omega u\,\omega=0$ for every $u\in C^\infty(M)$.
	(Hint: integrate $\mathrm{div}(\nabla u)$ and use Stokes' theorem on a manifold without boundary.)
\end{exercise}

\begin{exercise}[The kernel consists of constants]
	Show that if $u\in W^{1,2}(M)$ is a weak solution of $\Delta_\omega u=0$, then $u$ is constant.
	(Hint: test with $\varphi=u-\frac{1}{\mathrm{Vol}(M)}\int_M u\,\omega$ and use Poincar\'e.)
\end{exercise}

\Step{B. Explicit computations}

\begin{exercise}[One Fourier mode on the flat torus]
	On $T^2=\R^2/\Z^2$ with $\Delta=-(\partial_x^2+\partial_y^2)$, solve
	\[
	\Delta u = \sin(2\pi x)\sin(4\pi y),
	\qquad \int_{T^2}u\,dx\,dy=0.
	\]
\end{exercise}

\begin{exercise}[Spherical-harmonic inversion]
	On the round sphere $S^2$, assume $\Delta Y_{\ell m}=\ell(\ell+1)Y_{\ell m}$.
	Given
	\[
	f=\sum_{\ell\ge 1}\sum_{m=-\ell}^{\ell} a_{\ell m} Y_{\ell m},
	\]
	write down the solution $u$ of $\Delta u=f$ with $\int_{S^2}u\,dA=0$.
\end{exercise}

\Step{C. Variational viewpoint}

\begin{exercise}[The energy minimizer]
	Define the functional on $W^{1,2}_\perp(M)$:
	\[
	E(v)=\frac12\int_M |\nabla v|^2\,\omega - \int_M f v\,\omega.
	\]
	Show that the weak solution $u$ of $\Delta_\omega u=f$ is the unique minimizer of $E$.
\end{exercise}

\Step{D. Regularity and local-to-global logic}

\begin{exercise}[Weyl's lemma by mollification]
	Give a self-contained proof that weakly harmonic functions on an open set $U\subset\R^2$ are smooth
	by mollification, following Lemma~\ref{lem:weyl}.
	Your proof should explicitly justify why $\Delta(\chi_\varepsilon*u)=\chi_\varepsilon*(\Delta u)$.
\end{exercise}

\Step{E. Boundary variation (preview)}

\begin{exercise}[Dirichlet and Neumann compatibility]
	Let $\Omega\subset\R^2$ be a bounded smooth domain. Compare solvability of $\Delta u=f$ under:
	(i) Dirichlet boundary condition $u|_{\partial\Omega}=0$, (ii) Neumann boundary condition $\partial_\nu u|_{\partial\Omega}=0$.
	Identify the correct ``mean-zero'' type compatibility condition in the Neumann case.
\end{exercise}

\EssayInhoSong
\clearpage
\part{Sheaves and Cohomology on Complex Manifolds}
\label{part:linebundle}

\section{Holomorphic line bundles, divisors, and the Picard group}
\index{line bundles!holomorphic}\index{divisors}\index{Picard group}
\label{sec:lb-div-pic}

\subsection{Holomorphic line bundles: cocycles, gauge changes, and models}
\index{line bundles!cocycles}\index{gauge changes}
\label{subsec:holomorphic-line-bundles}

\subsubsection{Motivation: local holomorphic data and multiplicative ambiguity}

On a complex manifold $M$ in the sense of
Definition~\ref{def:complex-manifold}, holomorphic functions are
sections of the trivial line bundle.
Many natural objects are locally holomorphic but do not glue as honest functions; instead,
they glue up to multiplication by nowhere-vanishing holomorphic functions.
This is exactly the mechanism encoded by holomorphic line bundles.

\subsubsection{Definition via transition functions}

\begin{definition}[Holomorphic line bundle via transition functions]
	\label{def:linebundle-cocycle}
	Let $M$ be a complex manifold and $\{U_\alpha\}_{\alpha\in A}$ an open cover.
	A \emph{holomorphic line bundle} on $M$ is specified by a family of holomorphic functions
	\[
	f_{\alpha\beta}\in \mathcal O^*(U_{\alpha\beta}),\qquad U_{\alpha\beta}:=U_\alpha\cap U_\beta,
	\]
	satisfying the cocycle identities
	\[
	f_{\alpha\alpha}=1,\qquad f_{\alpha\beta}f_{\beta\alpha}=1,\qquad
	f_{\alpha\beta}f_{\beta\gamma}f_{\gamma\alpha}=1\quad\text{on }U_{\alpha\beta\gamma}.
	\]
	Two such families $\{f_{\alpha\beta}\}$ and $\{\tilde f_{\alpha\beta}\}$ define isomorphic line bundles
	if there exist $g_\alpha\in\mathcal O^*(U_\alpha)$ such that
	\[
	\tilde f_{\alpha\beta} = g_\alpha\, f_{\alpha\beta}\, g_\beta^{-1}\quad\text{on }U_{\alpha\beta}.
	\]
\end{definition}

\subsubsection{Gluing construction and complex manifold structure}

\begin{proposition}[Bundle obtained by gluing trivial bundles]
	\label{prop:gluing-linebundle}
	Let $\{U_\alpha\}$ be an open cover of $M$ and $\{f_{\alpha\beta}\}$ a cocycle as in
	Definition~\ref{def:linebundle-cocycle}. Consider the disjoint union
	\[
	E := \bigsqcup_{\alpha\in A}\, (U_\alpha\times \C),
	\]
	and define an equivalence relation $\sim$ on $E$ by declaring that for $x\in U_{\alpha\beta}$,
	\[
	(\alpha,x,\eta)\sim(\beta,x,\xi)\quad\Longleftrightarrow\quad \xi = f_{\beta\alpha}(x)\,\eta.
	\]
	Then the quotient space $L:=E/\!\sim$ admits a unique complex manifold structure such that:
	\begin{enumerate}
		\item the natural map $\pi:L\to M$, $\pi([\alpha,x,\eta])=x$, is holomorphic;
		\item for each $\alpha$, the map
		\[
		\varphi_\alpha:\pi^{-1}(U_\alpha)\to U_\alpha\times\C,\qquad
		\varphi_\alpha([\beta,x,\xi]) := (x, f_{\alpha\beta}(x)\,\xi)
		\]
		is a biholomorphism which is $\C$-linear on fibers;
		\item on overlaps $U_{\alpha\beta}$ the change of trivialization is multiplication by $f_{\alpha\beta}$:
		\[
		\varphi_\alpha\circ \varphi_\beta^{-1}(x,\eta) = (x, f_{\alpha\beta}(x)\,\eta).
		\]
	\end{enumerate}
	Hence $\pi:L\to M$ is a holomorphic line bundle with transition cocycle $\{f_{\alpha\beta}\}$.
\end{proposition}

\begin{proof}
	\emph{Step 1: $\sim$ is an equivalence relation.}
	Reflexivity holds because $f_{\alpha\alpha}=1$ gives $(\alpha,x,\eta)\sim(\alpha,x,\eta)$.
	Symmetry holds because $f_{\beta\alpha}=f_{\alpha\beta}^{-1}$: if $\xi=f_{\beta\alpha}\eta$ then
	$\eta=f_{\alpha\beta}\xi$, hence $(\beta,x,\xi)\sim(\alpha,x,\eta)$.
	Transitivity: suppose $(\alpha,x,\eta)\sim(\beta,x,\xi)$ and $(\beta,x,\xi)\sim(\gamma,x,\zeta)$ with
	$x\in U_{\alpha\beta\gamma}$. Then
	\[
	\xi = f_{\beta\alpha}(x)\eta,\qquad \zeta = f_{\gamma\beta}(x)\xi
	\quad\mathrel{\text{and hence}}\quad
	\zeta = f_{\gamma\beta}(x)f_{\beta\alpha}(x)\eta.
	\]
	The cocycle identity $f_{\gamma\beta}f_{\beta\alpha}=f_{\gamma\alpha}$ (equivalent to
	$f_{\alpha\beta}f_{\beta\gamma}f_{\gamma\alpha}=1$) yields $\zeta=f_{\gamma\alpha}(x)\eta$, i.e.
	$(\alpha,x,\eta)\sim(\gamma,x,\zeta)$.
	
	\emph{Step 2: define charts on $L$.}
	For fixed $\alpha$, define $\varphi_\alpha$ on $\pi^{-1}(U_\alpha)$ by
	\[
	\varphi_\alpha([\beta,x,\xi]) := (x, f_{\alpha\beta}(x)\,\xi).
	\]
	We must show $\varphi_\alpha$ is well-defined. Suppose $[\beta,x,\xi]=[\beta',x,\xi']$.
	Then $\xi' = f_{\beta'\beta}(x)\xi$. Using the cocycle identity on $U_{\alpha\beta\beta'}$,
	\[
	f_{\alpha\beta'}(x)\xi' = f_{\alpha\beta'}(x)f_{\beta'\beta}(x)\xi = f_{\alpha\beta}(x)\xi,
	\]
	so $\varphi_\alpha$ assigns the same value; hence it is well-defined.
	
	\emph{Step 3: $\varphi_\alpha$ is bijective with explicit inverse.}
	Given $(x,\eta)\in U_\alpha\times\C$, set $\psi_\alpha(x,\eta):=[\alpha,x,\eta]\in L$.
	Then $\varphi_\alpha(\psi_\alpha(x,\eta))=(x,f_{\alpha\alpha}(x)\eta)=(x,\eta)$.
	Conversely, for any $[\beta,x,\xi]\in \pi^{-1}(U_\alpha)$,
	\[
	\psi_\alpha(\varphi_\alpha([\beta,x,\xi])) = [\alpha,x,f_{\alpha\beta}(x)\xi].
	\]
	But $(\beta,x,\xi)\sim(\alpha,x,f_{\alpha\beta}(x)\xi)$ because
	$\xi = f_{\beta\alpha}(x)f_{\alpha\beta}(x)\xi$, so these define the same class in $L$.
	Thus $\psi_\alpha$ is the inverse of $\varphi_\alpha$.
	
	\emph{Step 4: holomorphic transition maps.}
	On $U_{\alpha\beta}$, the overlap map
	\[
	\varphi_\alpha\circ \varphi_\beta^{-1}: (x,\eta)\mapsto (x,f_{\alpha\beta}(x)\eta)
	\]
	is holomorphic because $f_{\alpha\beta}\in\mathcal O^*(U_{\alpha\beta})$.
	The quotient is Hausdorff: points over distinct base points are separated by
	preimages of disjoint neighborhoods in the Hausdorff manifold $M$, while two
	distinct points in one fiber are separated inside any common bundle chart.
	Because $M$ is second-countable, the cover has a countable refinement by
	coordinate domains on which the bundle is trivial; countable product bases on
	these domains give a countable basis for $L$. Therefore the atlas
	$\{(\pi^{-1}(U_\alpha),\varphi_\alpha)\}$ defines a complex manifold structure on $L$.
	By construction, $\pi$ is holomorphic. Each fiber is identified with $\C$, and
	the overlap maps are multiplication by nonzero complex numbers; consequently the
	vector-space operations transported through the trivializations agree on overlaps.
	Hence $L$ is a holomorphic line bundle.
\end{proof}

\subsubsection{Sections as glued holomorphic functions}

\begin{proposition}[Local descriptions of holomorphic sections]
	\label{prop:sections-glue}
	Let $L$ be defined by transition functions $\{f_{\alpha\beta}\}$ on $\{U_\alpha\}$.
	A holomorphic section $s$ of $L$ over $U\subset M$ is equivalent to holomorphic functions
	$s_\alpha\in\mathcal O(U\cap U_\alpha)$ satisfying
	\[
	s_\alpha = f_{\alpha\beta}\, s_\beta\quad\text{on }U\cap U_{\alpha\beta}.
	\]
	Moreover, given such $\{s_\alpha\}$, there exists a unique section $s$ whose $\alpha$-trivialization
	is $s_\alpha$.
\end{proposition}

\begin{proof}
	Fix $U\subset M$. Suppose $s:U\to L$ is a holomorphic section. For each $\alpha$, define
	\[
	s_\alpha := \mathrm{pr}_2\circ \varphi_\alpha\circ s \in \mathcal O(U\cap U_\alpha),
	\]
	where $\mathrm{pr}_2:U_\alpha\times\C\to\C$ is projection. On $U\cap U_{\alpha\beta}$,
	\[
	\varphi_\alpha\circ s = (\mathrm{id},s_\alpha),\qquad \varphi_\beta\circ s = (\mathrm{id},s_\beta),
	\]
	and since $\varphi_\alpha\circ \varphi_\beta^{-1}(x,\eta)=(x,f_{\alpha\beta}(x)\eta)$,
	we get
	\[
	s_\alpha(x) = f_{\alpha\beta}(x)\, s_\beta(x)\quad \forall x\in U\cap U_{\alpha\beta}.
	\]
	
	Conversely, assume $\{s_\alpha\}$ are holomorphic and satisfy $s_\alpha=f_{\alpha\beta}s_\beta$ on overlaps.
	Define a map $s:U\to L$ by: if $x\in U\cap U_\alpha$, set
	\[
	s(x):=\psi_\alpha(x,s_\alpha(x)) = [\alpha,x,s_\alpha(x)]\in L.
	\]
	We must show $s$ is well-defined. If $x\in U\cap U_{\alpha\beta}$, then
	\[
	[\alpha,x,s_\alpha(x)] = [\alpha,x,f_{\alpha\beta}(x)s_\beta(x)] = [\beta,x,s_\beta(x)]
	\]
	by the defining relation of the quotient, using $s_\alpha=f_{\alpha\beta}s_\beta$. Thus $s$ is well-defined.
	On each $U\cap U_\alpha$, the expression $x\mapsto [\alpha,x,s_\alpha(x)]$ is holomorphic because
	$\psi_\alpha$ is biholomorphic and $s_\alpha$ is holomorphic; hence $s$ is holomorphic on $U$.
	Uniqueness follows because the trivializations $\varphi_\alpha$ separate points in fibers:
	if two sections have the same second component in every trivialization, they coincide.
\end{proof}

\subsubsection{Gauge changes and isomorphism classes}

\begin{proposition}[Coboundary change produces an isomorphic bundle]
	\label{prop:gauge-isomorphism}
	Let $\{f_{\alpha\beta}\}$ and $\{\tilde f_{\alpha\beta}\}$ be cocycles on the same cover $\{U_\alpha\}$,
	and suppose $\tilde f_{\alpha\beta} = g_\alpha f_{\alpha\beta} g_\beta^{-1}$ for some $g_\alpha\in\mathcal O^*(U_\alpha)$.
	Let $L$ and $\tilde L$ be the corresponding glued bundles. Then there exists a holomorphic bundle isomorphism
	$\Phi:L\to \tilde L$ covering $\mathrm{id}_M$.
\end{proposition}

\begin{proof}
	Work on the disjoint union model $E=\bigsqcup_\alpha(U_\alpha\times\C)$ and $\tilde E=\bigsqcup_\alpha(U_\alpha\times\C)$.
	Define $T:E\to\tilde E$ by
	\[
	T(\alpha,x,\eta):=(\alpha,x,g_\alpha(x)\eta).
	\]
	We claim $T$ respects the equivalence relations. Indeed, if $(\alpha,x,\eta)\sim(\beta,x,\xi)$ in $E$,
	then $\xi=f_{\beta\alpha}(x)\eta$. Hence
	\[
	T(\beta,x,\xi)=(\beta,x,g_\beta(x)f_{\beta\alpha}(x)\eta).
	\]
	On the other hand, in $\tilde E$ we have the relation $(\alpha,x,\eta')\sim(\beta,x,\xi')$ iff
	$\xi'=\tilde f_{\beta\alpha}(x)\eta'$. Using $\tilde f_{\beta\alpha}=g_\beta f_{\beta\alpha} g_\alpha^{-1}$,
	\[
	\tilde f_{\beta\alpha}(x)\, (g_\alpha(x)\eta) = g_\beta(x) f_{\beta\alpha}(x)\eta,
	\]
	so $(\alpha,x,g_\alpha\eta)\sim(\beta,x,g_\beta f_{\beta\alpha}\eta)$ in $\tilde E$, i.e.
	$T(\alpha,x,\eta)\sim T(\beta,x,\xi)$ in $\tilde E$.
	Therefore $T$ descends to a well-defined map $\Phi:L\to\tilde L$ on the quotients.
	
	In the $\alpha$-trivializations, $\Phi$ is given by $(x,\eta)\mapsto (x,g_\alpha(x)\eta)$, which is holomorphic
	with holomorphic inverse $(x,\eta)\mapsto(x,g_\alpha(x)^{-1}\eta)$. Hence $\Phi$ is a holomorphic bundle isomorphism.
\end{proof}

\subsubsection{Tensor product, dual bundle, and induced cocycles}

\begin{proposition}[Operations on line bundles in cocycle form]
	\label{prop:tensor-dual-cocycle}
	Let $L,L'$ be line bundles defined on $\{U_\alpha\}$ by cocycles $\{f_{\alpha\beta}\}$ and $\{f'_{\alpha\beta}\}$.
	Then:
	\begin{enumerate}
		\item $L\otimes L'$ is defined by the cocycle $\{f_{\alpha\beta}f'_{\alpha\beta}\}$.
		\item The dual bundle $L^\vee$ is defined by $\{f_{\alpha\beta}^{-1}\}$.
	\end{enumerate}
\end{proposition}

\begin{proof}
	Our convention is that fiber coordinates satisfy
	$\eta_\alpha=f_{\alpha\beta}\eta_\beta$. If $e_\alpha$ denotes the frame
	whose $\alpha$-coordinate is $1$, the corresponding frame relation is therefore
	\[
	 e_\beta=f_{\alpha\beta}e_\alpha,
	 \qquad e'_\beta=f'_{\alpha\beta}e'_\alpha.
	\]
	Consequently
	$e_\beta\otimes e'_\beta=(f_{\alpha\beta}f'_{\alpha\beta})
	(e_\alpha\otimes e'_\alpha)$, so the coordinate transition functions of
	$L\otimes L'$ are $f_{\alpha\beta}f'_{\alpha\beta}$.

	For the dual frames, evaluation must remain one:
	\[
	 e_\beta^\vee=f_{\alpha\beta}^{-1}e_\alpha^\vee,
	 \qquad
	 e_\beta^\vee(e_\beta)=
	 f_{\alpha\beta}^{-1}e_\alpha^\vee(f_{\alpha\beta}e_\alpha)=1.
	\]
	Thus the coordinate cocycle of $L^\vee$ is $f_{\alpha\beta}^{-1}$.
\end{proof}

\subsubsection{The canonical bundle on a Riemann surface}

\begin{example}[Canonical bundle cocycle from coordinate changes]
	\label{ex:canonical-cocycle}
	Let $M$ be a Riemann surface with local coordinates $z_\alpha$ on $U_\alpha$.
	Define the local frame $e_\alpha:=dz_\alpha$ of the canonical bundle $K_M$.
	On overlaps $U_{\alpha\beta}$, write $z_\beta=z_\beta(z_\alpha)$; then
	\[
	dz_\beta = \frac{\partial z_\beta}{\partial z_\alpha}\, dz_\alpha.
	\]
	Hence $e_\beta = \bigl(\frac{\partial z_\beta}{\partial z_\alpha}\bigr) e_\alpha$. Since our
	coordinate convention is $\eta_\alpha=f_{\alpha\beta}\eta_\beta$, equivalently
	\[
	e_\beta = f_{\alpha\beta} e_\alpha,
	\qquad
	f_{\alpha\beta} = \frac{\partial z_\beta}{\partial z_\alpha}\in\mathcal O^*(U_{\alpha\beta}).
	\]
	On a triple overlap, the chain rule implies
	$f_{\alpha\gamma}=f_{\alpha\beta}f_{\beta\gamma}$, i.e. the cocycle condition.
\end{example}

\subsubsection{$\mathcal O(1)$ on $\mathbb{CP}^1$ and its global sections}

\begin{example}[$\mathcal O(1)$ on $\mathbb{CP}^1$ and $H^0(\mathbb{CP}^1,\mathcal O(1))$]
	\label{ex:O1-sections}
	Cover $\mathbb{CP}^1$ by $U_0=\{Z_0\neq0\}$ and $U_1=\{Z_1\neq0\}$.
	Let $z=Z_1/Z_0$ on $U_0$ and $w=Z_0/Z_1=1/z$ on $U_1$.
	Define $\mathcal O(1)$ by the transition function $f_{01}=z$ on $U_{01}$.
	A holomorphic section corresponds to $(s_0,s_1)$ with $s_0\in\mathcal O(U_0)$, $s_1\in\mathcal O(U_1)$ and
	\[
	s_0(z) = z\, s_1(1/z)\quad \text{on }U_{01}.
	\]
	Expand $s_1(w)=\sum_{k\ge0} a_k w^k$ (entire on $U_1\simeq\C$). Then
	\[
	s_0(z)= z\, s_1(1/z) = z\sum_{k\ge0} a_k z^{-k} = \sum_{k\ge0} a_k z^{1-k}.
	\]
	For $s_0$ to be entire in $z$, all negative powers must vanish, so $a_k=0$ for $k\ge2$.
	Thus $s_1(w)=a_0+a_1 w$ and $s_0(z)=a_0 z + a_1$.
	Hence
	\[
	H^0(\mathbb{CP}^1,\mathcal O(1))\cong \C^2,
	\]
	with basis corresponding to $(s_0,s_1)=(z,1)$ and $(1,w)$, i.e. homogeneous sections $Z_1$ and $Z_0$.
\end{example}

\paragraph{One bundle, four descriptions.}
The calculation in Example~\ref{ex:O1-sections} will reappear in several
languages.  The transition function $z$ glues the two local frames; the
section represented by $(s_0,s_1)=(1,w)$ vanishes once, at $w=0$, and
nowhere else.  Its zero divisor is therefore $[\infty]$, so the divisor
construction below identifies this bundle with $\mathcal O_{\PP^1}(\infty)$.
The same calculation gives $h^0(\PP^1,\mathcal O(1))=2$, while the zero
divisor gives $\deg\mathcal O(1)=1$.  Later, the first Chern class will
record this degree by
$\langle c_1(\mathcal O(1)),[\PP^1]\rangle=1$.  These are not four
different bundles: transition functions specify the gluing, a section
exhibits the divisor, and degree and $c_1$ record its global size.

\subsection{Divisors and the associated line bundle}
\index{divisors!associated line bundle}
\label{subsec:divisors-OD}

\subsubsection{Divisors on a Riemann surface}

Throughout this subsection, let $X$ be a compact Riemann surface.  (On a
noncompact surface one generally allows locally finite divisors; compactness
then makes their support finite.)

\begin{definition}[Divisor]
	\label{def:divisor}
	A \emph{divisor} on $X$ is a formal finite integer linear combination of points:
	\[
	D = \sum_{p\in X} n_p\, p,
	\qquad n_p\in\Z,\quad n_p=0 \text{ for all but finitely many }p.
	\]
	The \emph{degree} of $D$ is $\deg(D):=\sum_p n_p$.
	We write $D\ge0$ if $n_p\ge0$ for all $p$ (effective divisor).
\end{definition}

\begin{definition}[Principal divisor]
	\label{def:principal-divisor}
	Let $f$ be a nonzero meromorphic function on $X$. For each $p\in X$, let $\ord_p(f)\in\Z$ be the order
	of vanishing (positive) or pole (negative) of $f$ at $p$. The \emph{principal divisor} of $f$ is
	\[
	(f):=\sum_{p\in X}\ord_p(f)\, p.
	\]
\end{definition}

\subsubsection{The sheaf $\mathcal O(D)$ and local growth condition}

\begin{definition}[Sheaf of meromorphic functions bounded by $D$]
	\label{def:OD-sheaf}
	For a divisor $D=\sum n_p p$ define $\mathcal O(D)$ on open sets $U\subset X$ by
	\[
	\mathcal O(D)(U)
	:=
	\bigl\{\, f\in \mathcal M(U)\setminus\{0\}\ \big|\ \ord_p(f)\ge -n_p\ \text{for all }p\in U \,\bigr\}
	\ \cup\ \{0\},
	\]
	where $\mathcal M$ is the sheaf of meromorphic functions on $X$.
\end{definition}

The condition $\ord_p(f)\ge -n_p$ means: $f$ may have poles, but their orders are bounded by $D$.

\begin{proposition}[Restriction and $\mathcal O(D)$ is a sheaf]
	\label{prop:OD-sheaf}
	$\mathcal O(D)$ is a sheaf of $\mathcal O$-modules, and $\mathcal O(D)$ is locally free of rank $1$.
\end{proposition}

\begin{proof}
	\emph{Step 1: restriction is well-defined.}
	If $f\in\mathcal O(D)(U)$ and $V\subset U$, then for each $p\in V$ we have $\ord_p(f|_V)=\ord_p(f)\ge -n_p$,
	so $f|_V\in\mathcal O(D)(V)$.
	
	\emph{Step 2: locality and gluing.}
	Let $U=\bigcup_\alpha U_\alpha$ be a cover.
	
	Locality: if $f,g\in\mathcal O(D)(U)$ and $f|_{U_\alpha}=g|_{U_\alpha}$ for all $\alpha$, then $f=g$ in $\mathcal M(U)$,
	hence in $\mathcal O(D)(U)$.
	
	Gluing: suppose $f_\alpha\in\mathcal O(D)(U_\alpha)$ satisfy $f_\alpha=f_\beta$ in $\mathcal M(U_{\alpha\beta})$.
	Since $\mathcal M$ is a sheaf, there is a unique $f\in\mathcal M(U)$ with $f|_{U_\alpha}=f_\alpha$.
	Fix $p\in U$ and choose $\alpha$ with $p\in U_\alpha$. Then $\ord_p(f)=\ord_p(f_\alpha)\ge -n_p$ because order is local.
	Hence $f\in\mathcal O(D)(U)$. This proves the sheaf axioms.
	
	\emph{Step 3: local freeness of rank $1$.}
	Fix $p\in X$. Choose a coordinate $z$ on a neighborhood $V$ of $p$ with $z(p)=0$,
	shrinking $V$ so that it contains no other point in the support of $D$.
	If the coefficient of $D$ at $p$ is $n_p$, set $s_p=z^{-n_p}$; this formula applies
	for positive, zero, and negative $n_p$.  For any $f\in\mathcal O(D)(V)$, the product
	$h=fz^{n_p}$ is holomorphic at $p$ and elsewhere on $V$.  Hence $f=h\,s_p$ with
	$h\in\mathcal O(V)$.
	Conversely, any $h\cdot s_p$ satisfies the pole bound. Therefore $\mathcal O(D)|_V$ is generated by $s_p$ as an $\mathcal O$-module.
	This shows local triviality with one generator, i.e. locally free of rank $1$.
\end{proof}

\subsubsection{$\mathcal O(D)$ as a line bundle via transition functions}

\begin{proposition}[Transition cocycle for $\mathcal O(D)$]
	\label{prop:OD-cocycle}
	Let $\{U_\alpha\}$ be an open cover such that for each $\alpha$ there exists a meromorphic function
	$g_\alpha\in\mathcal M^*(U_\alpha)$ with
	\[
	(g_\alpha) = D|_{U_\alpha}\quad \text{as divisors on }U_\alpha.
	\]
	Then on overlaps $U_{\alpha\beta}$ the ratios
	\[
	f_{\alpha\beta} := \frac{g_\alpha}{g_\beta} \in \mathcal O^*(U_{\alpha\beta})
	\]
	define a holomorphic line bundle whose sheaf of holomorphic sections is isomorphic to $\mathcal O(D)$.
\end{proposition}

\begin{proof}
	\emph{Step 1: $f_{\alpha\beta}$ is holomorphic and nowhere vanishing.}
	On $U_{\alpha\beta}$,
	\[
	(g_\alpha) - (g_\beta) = D|_{U_{\alpha\beta}} - D|_{U_{\alpha\beta}} = 0.
	\]
	Thus the meromorphic function $g_\alpha/g_\beta$ has no zeros and no poles on $U_{\alpha\beta}$,
	hence belongs to $\mathcal O^*(U_{\alpha\beta})$.
	
	\emph{Step 2: cocycle identity.}
	On triple overlaps,
	\[
	f_{\alpha\beta}f_{\beta\gamma}f_{\gamma\alpha}
	=\frac{g_\alpha}{g_\beta}\cdot\frac{g_\beta}{g_\gamma}\cdot\frac{g_\gamma}{g_\alpha}=1.
	\]
	Also $f_{\alpha\alpha}=1$ and $f_{\alpha\beta}f_{\beta\alpha}=1$ are immediate.
	
	\emph{Step 3: identify $\mathcal O(D)$ with the section sheaf of the glued bundle.}
	Let $L$ be the line bundle defined by $\{f_{\alpha\beta}\}$ with local frames $e_\alpha$ satisfying
	$e_\beta=f_{\alpha\beta}e_\alpha$.
	Define an $\mathcal O$-module map on each $U_\alpha$:
	\[
	\Theta_\alpha:\mathcal O(U_\alpha)\to \mathcal O(D)(U_\alpha),\qquad
	h\mapsto h\cdot g_\alpha^{-1}.
	\]
	Since $(g_\alpha)=D|_{U_\alpha}$, multiplication by $g_\alpha^{-1}$ shifts divisor bounds precisely so that
	$h\cdot g_\alpha^{-1}$ has poles bounded by $D$ and hence lies in $\mathcal O(D)(U_\alpha)$.
	
	On overlaps $U_{\alpha\beta}$,
	\[
	\Theta_\alpha(h)
	= h g_\alpha^{-1}
	= h \left(\frac{g_\beta}{g_\alpha}\right) g_\beta^{-1}
	= h f_{\beta\alpha}\, g_\beta^{-1}
	= \Theta_\beta(h f_{\beta\alpha}).
	\]
	This is exactly the compatibility rule for sections:
	if a section is represented by local holomorphic functions $s_\alpha$ in the frames $e_\alpha$,
	the overlap condition is $s_\alpha = f_{\alpha\beta}s_\beta$, equivalently $s_\beta = f_{\beta\alpha}s_\alpha$.
	Thus $\{\Theta_\alpha\}$ glue to an isomorphism between the sheaf of holomorphic sections of $L$ and $\mathcal O(D)$.
\end{proof}

\subsection{Picard group and classification}
\label{subsec:picard}

\subsubsection{Picard group and tensor product}

\begin{definition}[Picard group]
	The \emph{Picard group} of a complex manifold $M$ is the set of isomorphism classes of holomorphic line bundles on $M$,
	with group operation given by tensor product:
	\[
	\mathrm{Pic}(M) := \{\text{holomorphic line bundles on $M$}\}/\cong,
	\qquad
	[L]\cdot[L'] := [L\otimes L'].
	\]
	The identity is the class of the trivial bundle, and the inverse of $[L]$ is $[L^\vee]$.
\end{definition}

\subsubsection{Classification by \v{C}ech cohomology}

\begin{theorem}[Line bundles and $H^1(M,\mathcal O^*)$]
	\label{thm:picard-cech}
	Let $M$ be a complex manifold. There is a canonical group isomorphism
	\[
	\mathrm{Pic}(M)\ \cong\ H^1(M,\mathcal O^*).
	\]
\end{theorem}

\begin{proof}
	\emph{Step 1: from line bundles to cocycles.}
	Given a line bundle $L$, choose a trivializing open cover $\{U_\alpha\}$ and local
	trivializations $\varphi_\alpha:\pi^{-1}(U_\alpha)\simeq U_\alpha\times\C$.
	On overlaps, the transition maps are
	\[
	(x,\eta_\beta)\mapsto (x,f_{\alpha\beta}(x)\eta_\beta)
	\]
	for some $f_{\alpha\beta}\in\mathcal O^*(U_{\alpha\beta})$ satisfying the cocycle condition.
	Thus $L$ determines a class $[f]\in \check H^1(\{U_\alpha\},\mathcal O^*)$.
	
	If one chooses different trivializations, the transition functions change by
	$\tilde f_{\alpha\beta}=g_\alpha f_{\alpha\beta} g_\beta^{-1}$, i.e. by a coboundary.
	Refining the cover pulls back the same cocycle class. Hence the class is independent
	of the trivializing cover and frames and defines an element of the direct-limit
	\v{C}ech group $H^1(M,\mathcal O^*)$.
	
	\emph{Step 2: from cocycles to line bundles.}
	Conversely, represent a class in $H^1(M,\mathcal O^*)$ by a cocycle
	$\{f_{\alpha\beta}\}$ on some open cover. Proposition~\ref{prop:gluing-linebundle}
	constructs a line bundle $L_f$.
	
	If $f$ and $\tilde f$ differ by a coboundary, Proposition~\ref{prop:gauge-isomorphism} gives an isomorphism
	$L_f\simeq L_{\tilde f}$. Thus a cohomology class in $H^1(M,\mathcal O^*)$ determines an isomorphism class of line bundles.
	
	\emph{Step 3: bijectivity.}
	The two constructions are inverse to each other:
	starting from $L$, forming its cocycle $f$, and gluing back produces a bundle canonically isomorphic to $L$
	(because the glued model reproduces the original transition maps).
	Starting from a cocycle $f$, the cocycle extracted from $L_f$ is $f$ by construction.
	Hence we have a bijection $\mathrm{Pic}(M)\leftrightarrow H^1(M,\mathcal O^*)$.
	
	\emph{Step 4: compatibility with group law.}
	Under tensor product, transition functions multiply: $(L\otimes L')$ has cocycle $f_{\alpha\beta}f'_{\alpha\beta}$
	(Proposition~\ref{prop:tensor-dual-cocycle}), which is exactly the addition law in $H^1(M,\mathcal O^*)$.
	Thus the bijection is a group isomorphism.
\end{proof}

\subsection{Exercises}
\label{subsec:lb-div-pic-exercises}

\begin{exercise}[Checking the atlas in the gluing construction]
	Supply full details for the topological argument in Proposition~\ref{prop:gluing-linebundle}.
	First separate two points lying over distinct base points by pulling back disjoint
	neighborhoods in $M$; then separate two distinct points in the same fiber inside one
	local trivialization. Finally choose a countable subcover from a countable basis of
	$M$ and use countable bases for the corresponding products $U_\alpha\times\C$.
	Deduce that the charts $\{\varphi_\alpha\}$ make $L$ Hausdorff and second-countable.
\end{exercise}

\begin{exercise}[Sections as glued functions]
	Prove Proposition~\ref{prop:sections-glue} from the quotient model $L=(\bigsqcup U_\alpha\times\C)/\!\sim$,
	without appealing to abstract sheaf language.
\end{exercise}

\begin{exercise}[Cocycle for tensor product and dual]
	Let $L,L'$ be line bundles with transition functions $f_{\alpha\beta},f'_{\alpha\beta}$ on the same cover.
	Construct explicit local trivializations of $L\otimes L'$ and $L^\vee$ and derive the cocycles
	$f_{\alpha\beta}f'_{\alpha\beta}$ and $f_{\alpha\beta}^{-1}$.
\end{exercise}

\begin{exercise}[Local model for $\mathcal O(D)$ near a point]
	\label{exr:selected-local-divisor}
	Let $X$ be a Riemann surface, $p\in X$, and $D=n\,p$.
	Choose a coordinate $z$ near $p$ with $z(p)=0$.
	Show that on a sufficiently small neighborhood $U$ of $p$,
	\[
	\mathcal O(D)(U) = z^{-n}\cdot \mathcal O(U)\qquad (n\in\mathbb Z).
	\]
	In particular, when $n<0$ this condition forces a zero of order at least $-n$;
	it does not give all of $\mathcal O(U)$.
\end{exercise}

\begin{exercise}[$\mathcal O(n)$ on $\mathbb{CP}^1$]
	\label{exr:selected-projective-line}
	Define $\mathcal O(n)$ on $\mathbb{CP}^1$ by transition function $f_{01}=z^n$ on $U_{01}$.
	Compute $H^0(\mathbb{CP}^1,\mathcal O(n))$ by the gluing equation and show
	\[
	\dim_\C H^0(\mathbb{CP}^1,\mathcal O(n)) =
	\begin{cases}
		n+1,& n\ge 0,\\
		0,& n<0.
	\end{cases}
	\]
\end{exercise}

\begin{exercise}[Principal divisors and isomorphic line bundles]
	Let $D$ be a divisor on a Riemann surface $X$ and let $f\in\mathcal M^*(X)$.
	Prove that $\mathcal O(D)\cong \mathcal O(D+(f))$ as line bundles.
	(Hint: construct an explicit isomorphism by multiplication with $f^{-1}$.)
\end{exercise}

\begin{exercise}[Picard group of $\mathbb{CP}^1$]
	Using the previous exercise and $\mathcal O(n)$, show that $\mathrm{Pic}(\mathbb{CP}^1)\cong \Z$
	via $[\mathcal O(n)]\mapsto n$.
\end{exercise}

\subsection{Canonical sections and worked divisor cocycles}
\label{subsec:divisors-OD-worked}

Let $M$ be a compact Riemann surface.  Proposition~\ref{prop:OD-cocycle}
constructs $\mathcal O(D)$ from local meromorphic equations $g_\alpha$ for
$D$, with transition functions $f_{\alpha\beta}=g_\alpha/g_\beta$ and frames
$e_\beta=f_{\alpha\beta}e_\alpha$.  These conventions also produce a canonical
meromorphic section
\[
s_D|_{U_\alpha}=g_\alpha e_\alpha.
\]
Indeed, $g_\alpha=f_{\alpha\beta}g_\beta$ is exactly the section-compatibility
condition, and the local coefficient $g_\alpha$ has divisor $D|_{U_\alpha}$.
Thus $(s_D)=D$.  The following computations test this convention at a point,
on $\CP^1$, and under tensor product.

\subsubsection{Worked examples}

\begin{example}[An explicit cocycle for $\mathcal O(p)$ and the divisor of $s_p$]
	Fix $p\in M$. Choose $U_0=M\setminus\{p\}$ and a coordinate disk $U_1$ around $p$ with holomorphic coordinate $z$ and $z(p)=0$.
	Set
	\[
	g_0\equiv 1 \ \text{on }U_0,\qquad g_1=z \ \text{on }U_1.
	\]
	Then on $U_{01}=U_0\cap U_1$,
	\[
	f_{01}=\frac{g_0}{g_1}=\frac{1}{z}\in\mathcal O^*(U_{01}),\qquad f_{10}=z.
	\]
	Thus $\mathcal O(p)$ is obtained by gluing $U_0\times\C$ and $U_1\times\C$ via $\eta_0=\frac{1}{z}\eta_1$.
	
	The canonical meromorphic section is
	\[
	s_p|_{U_0}=g_0 e_0=e_0,\qquad s_p|_{U_1}=g_1 e_1=z e_1.
	\]
	Near $p$, $s_p=z e_1$, hence $\nu_p(s_p)=1$ and there are no other zeros/poles.
	Therefore $(s_p)=p$ and $\deg\mathcal O(p)=1$.
\end{example}

\begin{example}[$\mathbb{CP}^1$ and $\mathrm{div}\!\Bigl(\frac{z-a}{z-b}\Bigr)$]
	On $\mathbb{CP}^1$ with affine coordinate $z$ on $U_0\simeq\C$, consider
	\[
	h(z)=\frac{z-a}{z-b}.
	\]
	At $z=a$ we have a simple zero: $h(z)=(z-a)\cdot(\text{unit})$, so $\nu_a(h)=+1$.
	At $z=b$ we have a simple pole: $h(z)=(z-b)^{-1}\cdot(\text{unit})$, so $\nu_b(h)=-1$.
	No other finite points contribute.
	
	To check $\infty$, use $w=1/z$ near $\infty$:
	\[
	h(1/w)=\frac{\frac{1}{w}-a}{\frac{1}{w}-b}=\frac{1-aw}{1-bw}.
	\]
	This is holomorphic and nonzero at $w=0$, so $\nu_\infty(h)=0$.
	Hence
	\[
	\mathrm{div}(h)=(a)-(b),\qquad \deg(\mathrm{div}(h))=1-1=0.
	\]
	Consequently $\mathcal O\bigl((a)-(b)\bigr)\simeq \mathcal O_{\mathbb{CP}^1}$ (principal divisor gives a trivial bundle).
\end{example}

\begin{example}[A cocycle proof that $\mathcal O(p)\otimes\mathcal O(-p)\simeq\mathcal O_M$]
	Use the same cover $U_0,U_1$ and coordinate $z$ as in the first example.
	For $\mathcal O(p)$ we have $f_{01}=1/z$.
	For $\mathcal O(-p)$ we take Cartier data $g_0=1$, $g_1=z^{-1}$, giving
	\[
	f'_{01}=\frac{g_0}{g_1}=z.
	\]
	Tensor product cocycle is the product:
	\[
	(f\cdot f')_{01}=\frac{1}{z}\cdot z=1,
	\]
	hence $\mathcal O(p)\otimes\mathcal O(-p)$ has trivial transition function, so it is isomorphic to $\mathcal O_M$.
\end{example}

\subsection{Picard group and the divisor correspondence}
\label{subsec:picard-divisor-worked}

Theorem~\ref{thm:picard-cech} identifies $\mathrm{Pic}(M)$ with
$H^1(M,\mathcal O^*)$ and shows that a change of local frame changes a
transition cocycle only by a coboundary.  It remains to relate this cocycle
classification to divisors.

\subsubsection{Divisors classify line bundles}

On a compact Riemann surface one has
\[
\mathrm{Pic}(M)\ \cong\ \mathrm{Div}(M)/\mathrm{Prin}(M),
\qquad
D\longmapsto \mathcal O(D).
\]
Concrete meaning: two divisors $D_1,D_2$ determine isomorphic line bundles iff $D_1-D_2=\mathrm{div}(h)$ for some meromorphic $h$.
The construction above proves that linearly equivalent divisors give isomorphic
bundles. Surjectivity---the existence of a nonzero meromorphic section of every
line bundle---uses Riemann--Roch in the development adopted here and is proved in
Proposition~\ref{prop:every-L-is-OD}. Thus the displayed isomorphism is an announced
result at this point, not an input to the proof of Riemann--Roch.

\subsubsection{Worked examples}

\begin{example}[Every degree-$0$ divisor on $\mathbb{CP}^1$ is principal]
	Let
	\[
	D=\sum_{i=1}^r m_i\,(a_i)\;-\;\sum_{j=1}^s n_j\,(b_j)
	\]
	be a divisor on $\mathbb{CP}^1$ supported in $\C$ (no $\infty$ for now),
	with $\deg(D)=\sum_i m_i-\sum_j n_j=0$.
	
	Define the rational function
	\[
	h(z)=\frac{\prod_{i=1}^r (z-a_i)^{m_i}}{\prod_{j=1}^s (z-b_j)^{n_j}}.
	\]
	Then $h$ has a zero of order $m_i$ at $a_i$ and a pole of order $n_j$ at $b_j$.
	Because the total degree in numerator equals that of denominator (degree $0$ divisor),
	there is no zero/pole at $\infty$. Hence
	\[
	\mathrm{div}(h)=D.
	\]
	More generally, write
	$D=\sum_{a\in\C}n_a(a)+n_\infty(\infty)$ with finite support and
	$\sum_a n_a+n_\infty=0$. Then
	\[
		h(z)=\prod_{a\in\C}(z-a)^{n_a}
	\]
	has order $n_a$ at each finite point $a$ and, after setting $z=1/w$,
	has order $-\sum_a n_a=n_\infty$ at infinity. Thus $(h)=D$ without an
	additional factor.
\end{example}

\begin{example}[$\mathrm{Pic}(\mathbb{CP}^1)\cong\Z$ from the degree computation]
	Take any divisor $D$ on $\mathbb{CP}^1$. Let $d=\deg(D)$.
	Choose a point $p$ and write
	\[
	D \sim d\cdot p,
	\]
	meaning $D-dp$ has degree $0$, hence is principal by the previous example:
	$D-dp=\mathrm{div}(h)$ for some rational $h$.
	Therefore
	\[
	\mathcal O(D)\ \simeq\ \mathcal O(dp)\ \simeq\ \mathcal O(d).
	\]
	For the last identification, choose a nonzero homogeneous linear form vanishing at
	$p$. Its divisor is $(p)$, so the divisor construction identifies
	$\mathcal O(p)$ with the hyperplane bundle $\mathcal O(1)$; tensor powers give
	$\mathcal O(dp)\simeq\mathcal O(d)$ for every $d\in\Z$.
	Using the surjectivity proved later in Proposition~\ref{prop:every-L-is-OD}, every
	line bundle is therefore $\mathcal O(m)$ for some $m\in\Z$. The degree distinguishes
	the integers $m$, so $\mathrm{Pic}(\mathbb{CP}^1)\cong\Z$, generated by
	$\mathcal O(1)$.
\end{example}

\subsection{\texorpdfstring{Global sections of $\mathcal O(m)$ on $\mathbb{CP}^1$: chart computation and divisors}{Global sections of O(m) on the projective line}}
\label{subsec:global-sections-Om-worked}

\subsubsection{Definition of $\mathcal O(m)$ and section compatibility}

Cover $\mathbb{CP}^1$ by $U_0,U_1$ with coordinates $z$ and $w=1/z$.
Define $\mathcal O(m)$ by transition $f_{01}=z^m$ and frames $e_1=z^m e_0$.
A section $s$ is given by holomorphic $f_0,f_1$ with
\[
s|_{U_0}=f_0(z)e_0,\qquad s|_{U_1}=f_1(w)e_1,\qquad f_0(z)=z^m f_1(1/z).
\]

\subsubsection{Worked examples}

\begin{example}[Global sections of $\mathcal O(m)$ for $m\ge0$]
	Let $m\ge 0$. Expand $f_0(z)=\sum_{k\ge 0} a_k z^k$ on $\C$.
	Then
	\[
	f_1(w)=w^m f_0(1/w)=\sum_{k\ge 0} a_k\, w^{m-k}.
	\]
	Holomorphicity at $w=0$ forces $a_k=0$ for $k>m$.
	So
	\[
	f_0(z)=a_0+a_1 z+\cdots +a_m z^m,
	\]
	and conversely any such polynomial produces a holomorphic $f_1$.
	Therefore
	\[
	\dim H^0(\mathbb{CP}^1,\mathcal O(m))=m+1,
	\]
	and a basis is given by the sections with
	\[
	f_0(z)=z^k,\qquad 0\le k\le m.
	\]
	In homogeneous coordinates, these correspond to monomials $Z_0^{m-k}Z_1^k$.
\end{example}

\begin{example}[Vanishing of $H^0(\mathbb{CP}^1,\mathcal O(m))$ for $m<0$]
	Let $m<0$ and write $m=-\ell$ with $\ell>0$.
	Then
	\[
	f_1(w)=w^m f_0(1/w)=w^{-\ell} f_0(1/w).
	\]
	Write the entire function $f_0$ as $f_0(z)=\sum_{k\ge0}a_kz^k$.  On the punctured
	$w$-chart,
	\[
	f_1(w)=\sum_{k\ge0}a_k w^{-\ell-k}.
	\]
	This Laurent series contains only strictly negative powers.  Holomorphicity of $f_1$
	at $w=0$ therefore forces every $a_k$ to vanish.  Thus $f_0\equiv0$ and the only
	global holomorphic section is $0$.
\end{example}

\begin{example}[The divisor of a section of $\mathcal O(m)$]
	Let $m\ge 0$ and take a nonzero section with
	\[
	f_0(z)=c\prod_{i=1}^r (z-a_i)^{m_i},\qquad \sum_{i=1}^r m_i=d\le m.
	\]
	Zeros in $\C$ contribute $\sum_i m_i(a_i)$.
	Near $\infty$, use $w=1/z$:
	\[
	f_1(w)=w^m f_0(1/w)
	=c\,w^m\prod_{i=1}^r \Bigl(\frac{1}{w}-a_i\Bigr)^{m_i}
	=c\,w^{m-d}\prod_{i=1}^r (1-a_i w)^{m_i}.
	\]
	So $f_1$ has a zero of order $m-d$ at $w=0$, i.e.\ the section vanishes at $\infty$ with order $m-d$.
	Hence
	\[
	(s)=\sum_{i=1}^r m_i\,(a_i)\ +\ (m-d)\,(\infty),
	\qquad
	\deg(s)=m.
	\]
\end{example}

\begin{example}[The section of $\mathcal O(2)$ corresponding to $Z_0Z_1$]
	Consider the homogeneous polynomial $P(Z_0,Z_1)=Z_0Z_1$ of degree $2$.
	On $U_0$ ($Z_0\neq 0$) with $z=Z_1/Z_0$, we have
	\[
	f_0(z)=P(1,z)=z,
	\]
	so the section vanishes to order $1$ at $z=0$, i.e.\ at $[1:0]$.
	Since $\deg f_0=1<2$, the previous example shows it vanishes at $\infty=[0:1]$ to order $2-1=1$.
	Thus
	\[
	(s)=[1:0]+[0:1].
	\]
\end{example}

\subsection{Picard group operations and degree}
\label{subsec:picard-degree-worked}

\subsubsection{Degree and its additivity/duality}

For any line bundle $L$ on compact $M$, choose a nonzero meromorphic section $s$ and define
\[
\deg(L)=\sum_{p\in M}\nu_p(s).
\]
This is well-defined because changing $s$ multiplies it by a meromorphic function, changing the divisor by a principal divisor of degree $0$.
The existence of such a section is not automatic from the local definition of a line bundle;
it is established later, using Riemann--Roch, in Proposition~\ref{prop:every-L-is-OD}.
Thus this paragraph records the degree construction in advance of its existence proof.

\begin{proposition}[Degree identities]
	For line bundles $L_1,L_2$,
	\[
	\deg(L_1\otimes L_2)=\deg(L_1)+\deg(L_2),\qquad
	\deg(L^\vee)=-\deg(L).
	\]
\end{proposition}

\begin{proof}[Worked local order computation]
	Pick meromorphic sections $s_i$ and local frames so $s_i=f_i e_i$.
	Then $s_1\otimes s_2=(f_1f_2)(e_1\otimes e_2)$, hence
	\[
	\nu_p(s_1\otimes s_2)=\nu_p(f_1f_2)=\nu_p(f_1)+\nu_p(f_2)=\nu_p(s_1)+\nu_p(s_2).
	\]
	Summing over $p$ gives additivity.
	
	For the dual, define $s^\vee$ locally by $s^\vee(s)=1$. If $s=f e$, then $s^\vee=f^{-1}e^\vee$, so
	$\nu_p(s^\vee)=-\nu_p(s)$ and $\deg(L^\vee)=-\deg(L)$.
\end{proof}

\subsubsection{Worked examples}

\begin{example}[The degree of $\mathcal O(p-q)$]
	Let $p,q\in M$. The divisor $p-q$ has degree $1-1=0$.
	Hence $\deg(\mathcal O(p-q))=\deg(p-q)=0$.
	Geometrically, $\mathcal O(p-q)$ is a degree-$0$ line bundle: it has meromorphic sections whose total number of zeros equals total number of poles.
	On $\mathbb{CP}^1$, this already forces triviality (because degree-$0$ divisors are principal), but on higher genus surfaces it produces nontrivial elements of $\mathrm{Pic}^0(M)$.
\end{example}

\begin{example}[Negative degree and vanishing of global sections]
	Assume $\deg(L)<0$ and suppose $s\in H^0(M,L)$ is a nonzero holomorphic section.
	Holomorphicity implies $s$ has no poles, so $(s)$ is an \emph{effective} divisor:
	\[
	(s)=\sum_{p\in M}\nu_p(s)\,p\quad\text{with }\nu_p(s)\ge 0.
	\]
	Therefore $\deg(L)=\deg((s))\ge 0$, contradicting $\deg(L)<0$.
	Hence $H^0(M,L)=0$ when $\deg(L)<0$.
\end{example}

\begin{example}[Degree of a tensor quotient]
	Using the proposition,
	\[
	\deg(L_1\otimes L_2^\vee)=\deg(L_1)+\deg(L_2^\vee)=\deg(L_1)-\deg(L_2).
	\]
	If you want to see the local mechanism: choose meromorphic sections $s_1,s_2$ and note $s_2^\vee$ has the opposite divisor of $s_2$, so
	\[
	(s_1\otimes s_2^\vee)=(s_1)-(s_2),
	\]
	hence degrees subtract.
\end{example}



\section{Sheaves, \v{C}ech Cohomology, and Derived-Functor Cohomology}
\index{sheaves}
\label{sec:sheaves-foundations}

\subsection{Presheaves and restriction maps}
\label{subsec:presheaves}

\subsubsection{Presheaves as ``data on open sets''}

A large class of geometric objects is naturally organized by open sets:
functions, differential forms, holomorphic sections, solutions of PDE, etc.
The only structure needed to start is restriction:
\[
V\subset U, \qquad \text{so data on }U\text{ restricts to data on }V.
\]

\begin{definition}[Presheaf of abelian groups]
	\label{def:presheaf}
	Let $M$ be a topological space. A \emph{presheaf} $\mathscr F$ of abelian groups on $M$ consists of:
	\begin{enumerate}
		\item for each open set $U\subset M$, an abelian group $\mathscr F(U)$;
		\item for each inclusion $V\subset U$, a homomorphism (restriction)
		\[
		\rho^U_V:\mathscr F(U)\to \mathscr F(V);
		\]
	\end{enumerate}
	such that $\rho^U_U=\mathrm{id}$ and $\rho^U_W=\rho^V_W\circ\rho^U_V$ for $W\subset V\subset U$.
\end{definition}

\subsubsection{Examples and verification}

\begin{example}[Continuous real-valued functions on $\R$]
	Let $M=\R$.
	Define $\mathscr C^0(U)=C^0(U,\R)$ with $\rho^U_V(f)=f|_V$.
	Then for $W\subset V\subset U$ and $f\in C^0(U)$:
	\[
	(\rho^V_W\circ\rho^U_V)(f)=\rho^V_W(f|_V)=(f|_V)|_W=f|_W=\rho^U_W(f),
	\]
	and $\rho^U_U(f)=f|_U=f$. Hence $\mathscr C^0$ is a presheaf.
\end{example}

\begin{example}[Smooth functions on $\R^n$ and a coordinate computation]
	Let $M=\R^n$ and $\mathscr C^\infty(U)=C^\infty(U,\R)$.
	If $f(x_1,\dots,x_n)=x_1^2+\sin(x_2)$ on $U$ and $V\subset U$, then
	\[
	\rho^U_V(f)=f|_V,\qquad \partial_{x_1}(f|_V)=(\partial_{x_1}f)|_V=2x_1|_V,
	\]
	so differentiation commutes with restriction. This locality is the reason $C^\infty$-objects form sheaves later.
\end{example}

\begin{example}[Differential forms and an explicit restriction]
	Let $M=\R^2$ with coordinates $(x,y)$.
	Define $\Omega^1(U)$ to be smooth $1$-forms on $U$.
	Take
	\[
	\omega = (x^2+y)\,dx + (\sin x)\,dy \in \Omega^1(U).
	\]
	If $V\subset U$, then
	\[
	\omega|_V = (x^2+y)|_V\,dx + (\sin x)|_V\,dy.
	\]
	If $W\subset V\subset U$, then $(\omega|_V)|_W=\omega|_W$, verifying functoriality of restriction.
\end{example}

\begin{example}[Solutions of a PDE as a presheaf: harmonic functions]
	Let $M=\R^2$. Define
	\[
	\mathscr H(U):=\{u\in C^2(U)\mid \Delta u=0\}.
	\]
	If $u\in\mathscr H(U)$ and $V\subset U$, then $u|_V\in C^2(V)$ and
	\[
	\Delta(u|_V)=(\Delta u)|_V = 0,
	\]
	so restriction preserves harmonicity. Thus $\mathscr H$ is a presheaf of abelian groups (under addition).
	Concrete computation: $u(x,y)=\Re(z^3)=x^3-3xy^2$ satisfies
	\[
	u_{xx}=6x,\qquad u_{yy}=-6x,\qquad \Delta u=u_{xx}+u_{yy}=0.
	\]
\end{example}

\subsection{Sheaves: locality and gluing}
\label{subsec:sheaves}

\subsubsection{The locality and gluing axioms}

\begin{definition}[Sheaf of abelian groups]
	\label{def:sheaf}
	A presheaf $\mathscr F$ on $M$ is a \emph{sheaf} if for every open $U\subset M$ and cover
	$U=\bigcup_\alpha U_\alpha$:
	\begin{enumerate}
		\item \textbf{Locality:} if $s,t\in\mathscr F(U)$ satisfy $s|_{U_\alpha}=t|_{U_\alpha}$ for all $\alpha$, then $s=t$.
		\item \textbf{Gluing:} if $s_\alpha\in\mathscr F(U_\alpha)$ satisfy
		$s_\alpha|_{U_{\alpha\beta}}=s_\beta|_{U_{\alpha\beta}}$ for all $\alpha,\beta$,
		then there exists $s\in\mathscr F(U)$ with $s|_{U_\alpha}=s_\alpha$ for all $\alpha$.
	\end{enumerate}
\end{definition}

\begin{definition}[Coherent analytic sheaf]
	\label{def:coherent-sheaf}
	Let $M$ be a complex manifold. A sheaf $\mathscr F$ of
	$\mathcal O_M$-modules is \emph{coherent} if every point has a neighborhood
	$U$ on which there is an exact sequence
	\[
	\mathcal O_U^{\oplus m}\longrightarrow
	\mathcal O_U^{\oplus n}\longrightarrow
	\mathscr F|_U\longrightarrow0
	\]
	for finite integers $m,n$, and the relations among any finite collection
	of local generators are themselves locally finitely generated. Equivalently,
	$\mathscr F$ is locally finitely presented and kernels of morphisms from
	finite free $\mathcal O_M$-modules are locally finitely generated.
\end{definition}

Holomorphic vector bundles give locally free coherent sheaves of sections.
Ideal sheaves of analytic subvarieties and their quotients provide the basic
non-locally-free examples. Coherence is the finiteness condition needed for
finite-dimensional cohomology on compact complex manifolds and for the vanishing
theorems used later.

\subsubsection{Gluing continuous functions: a worked piecewise example}

\begin{example}[Gluing a continuous function on $\R$ from two opens]
	Let $U=(-2,2)$, $U_1=(-2,1)$, $U_2=(-1,2)$ so $U=U_1\cup U_2$ and $U_{12}=(-1,1)$.
	Define
	\[
	f_1(x)=x^2 \quad (x\in U_1),\qquad f_2(x)=x^2 \quad (x\in U_2).
	\]
	Then $f_1|_{U_{12}}=f_2|_{U_{12}}$, so they glue to $f(x)=x^2$ on all of $U$.
	
	Now modify to a genuinely ``local'' construction:
	\[
	g_1(x)=\sin x \quad (x\in U_1),\qquad g_2(x)=
	\begin{cases}
		\sin x,& -1<x\le 1,\\
		\sin x + (x-1)^3,& 1<x<2,
	\end{cases}
	\]
	where $g_2$ is continuous on $U_2$.
	On the overlap $(-1,1)$ we have $g_2=\sin x=g_1$, hence compatibility holds.
	The glued function $g$ is:
	\[
	g(x)=
	\begin{cases}
		\sin x,& -2<x\le 1,\\
		\sin x + (x-1)^3,& 1<x<2,
	\end{cases}
	\]
	and $g$ is continuous on $U$ because the definitions agree on the overlap and continuity is local.
\end{example}

\subsubsection{Gluing holomorphic functions}

\begin{example}[Gluing holomorphic functions on $\C^\times$]
	Let $U=\C^\times$ and cover it by
	\[
	U_1=\C\setminus (-\infty,0],\qquad U_2=\C\setminus [0,\infty).
	\]
	Choose the branches with arguments $-\pi<\arg_1 z<\pi$ on $U_1$
	and $0<\arg_2 z<2\pi$ on $U_2$:
	\[
	\Log_1:U_1\to\C,\qquad \Log_2:U_2\to\C.
	\]
	On the overlap $U_{12}=U_1\cap U_2=\C\setminus\R$, which has two components,
	the branches agree above the real axis and differ below it:
	\[
	\Log_2(z)=\Log_1(z) \quad(\Im z>0),\qquad
	\Log_2(z)=\Log_1(z)+2\pi i \quad(\Im z<0).
	\]
	Hence $\{\Log_1,\Log_2\}$ do \emph{not} satisfy the compatibility condition on all of $U_{12}$,
	so they do not glue to a holomorphic function on $\C^\times$.
	This is a concrete place where the topology of the domain obstructs global gluing.
\end{example}

\subsubsection{A presheaf that fails gluing: bounded continuous functions}

\begin{example}[Bounded continuous functions]
	Let $M=\R$ and define $\mathscr B(U)$ to be bounded continuous functions $U\to\R$.
	Restriction preserves boundedness, so $\mathscr B$ is a presheaf.
	Cover $\R$ by $U_n=(-n,n)$ and set $f_n(x)=x$ on $U_n$.
	Then $f_n=f_m$ on $U_n\cap U_m$, but the glued function would be $f(x)=x$ on $\R$,
	which is not bounded. So gluing fails and $\mathscr B$ is not a sheaf.
\end{example}

\subsection{Morphisms of sheaves and exact sequences}
\label{subsec:sheaf-morphisms}

\subsubsection{Morphisms and concrete examples}

\begin{definition}[Morphism]
	\label{def:sheaf-morphism}
	A morphism $\varphi:\mathscr F\to\mathscr G$ of presheaves is a family of maps
	$\varphi_U:\mathscr F(U)\to\mathscr G(U)$ commuting with restrictions:
	$\rho^U_V\circ\varphi_U=\varphi_V\circ\rho^U_V$ for $V\subset U$.
\end{definition}

\begin{example}[Derivative as a morphism of sheaves on $\R$]
	Let $\mathscr C^\infty$ be the sheaf of smooth functions on $\R$ and define
	\[
	\varphi_U:C^\infty(U)\to C^\infty(U),\qquad \varphi_U(f)=f'.
	\]
	Then for $V\subset U$,
	\[
	(\varphi_U(f))|_V = (f')|_V = (f|_V)' = \varphi_V(f|_V),
	\]
	so $\varphi$ is a morphism of sheaves.
	Concrete computation: if $f(x)=x^3+\cos x$, then $\varphi(f)=3x^2-\sin x$.
\end{example}

\begin{example}[Exponential map as a morphism $\mathcal O\to\mathcal O^\times$]
	On a complex manifold $M$, define on each open $U$:
	\[
	\exp_U:\mathcal O(U)\to\mathcal O^\times(U),\qquad \exp_U(f)=e^{f}.
	\]
	Restriction compatibility is $(e^f)|_V=e^{f|_V}$.
\end{example}

\subsubsection{Exact sequences and a computed kernel}

\begin{definition}[Exactness at a term]
	A sequence $\mathscr F\xrightarrow{\varphi}\mathscr G\xrightarrow{\psi}\mathscr H$
	is exact at $\mathscr G$ if $\ker\psi=\mathrm{im}\,\varphi$ as subsheaves of $\mathscr G$.
\end{definition}

\begin{example}[Kernel of exponential on a connected domain]
	Let $\Omega\subset\C$ be a domain and consider
	\[
	\exp:\mathcal O(\Omega)\to\mathcal O^\times(\Omega),\qquad f\mapsto e^{2\pi i f}.
	\]
	If $e^{2\pi i f}\equiv 1$, then for each $z\in\Omega$ we have $2\pi i f(z)\in 2\pi i\Z$ so $f(z)\in\Z$.
	Since $f$ is holomorphic and $\Z$ is discrete, $f$ is locally constant, hence constant on connected $\Omega$,
	so $f\in\Z$. Thus
	\[
	\ker(\exp)=\Z \subset \mathcal O(\Omega)
	\]
	(viewing $\Z$ as constant holomorphic functions).
\end{example}

\subsection{Stalks, germs, and local criteria}
\index{stalks}\index{germs}
\label{subsec:stalks}

\subsubsection{Germs as equivalence classes with computations}

\begin{definition}[Germ and stalk]
	\label{def:stalk}
	Let $\mathscr F$ be a presheaf on $M$ and $x\in M$.
	Pairs $(U,s)$ with $x\in U$ and $s\in\mathscr F(U)$ define a germ at $x$.
	We set $(U,s)\sim(V,t)$ if there exists $W\subset U\cap V$ with $x\in W$ and $s|_W=t|_W$.
	The set of equivalence classes is the stalk $\mathscr F_x$, and the class of $(U,s)$ is written $[s]_x$.
\end{definition}

\begin{example}[Germs of continuous functions: two different representatives]
	Let $M=\R$, $x=0$, and consider
	\[
	f(x)=x^2 \text{ on }(-1,1),\qquad g(x)=x^2 \text{ on }(-2,2).
	\]
	Then $[f]_0=[g]_0$ because on $(-1,1)\subset (-1,1)\cap(-2,2)$ they coincide.
	But if $h(x)=|x|$ on $(-1,1)$, then $[h]_0\neq [f]_0$ since they disagree on every neighborhood of $0$.
\end{example}

\begin{example}[Stalk of $\underline{\Z}$ at a point]
	Let $\underline{\Z}$ be the locally constant sheaf:
	$\underline{\Z}(U)=\{n:U\to\Z \text{ locally constant}\}$.
	Fix $x\in M$. Any germ is represented by a locally constant function on some neighborhood,
	hence is constant on some (not necessarily connected) neighborhood of $x$ containing it.
	Therefore each germ is determined by an integer, and
	\[
	(\underline{\Z})_x \cong \Z.
	\]
\end{example}

\subsubsection{Stalkwise surjectivity does not imply surjectivity on global sections}

\begin{example}[A computed linear system on $\mathbb{CP}^1$]
	Let $X=\mathbb{CP}^1$ and consider the sheaf morphism
	\[
	\Phi:\mathcal O_X^{\oplus 2}\to \mathcal O_X(2),\qquad (a,b)\mapsto a\,Z_0^2 + b\,Z_1^2.
	\]
	On the affine chart $U_0=\{Z_0\neq 0\}$ with $z=Z_1/Z_0$, a local section of $\mathcal O(2)$ looks like
	$f(z)\cdot Z_0^2$. Also $Z_1^2=z^2 Z_0^2$. Hence locally
	\[
	\Phi(a,b)=\bigl(a(z)+b(z)z^2\bigr)\,Z_0^2.
	\]
	Given any germ $f(z)Z_0^2$ at a point of $U_0$, choose $a=f$ and $b=0$ to hit it;
	similarly on $U_1$. Thus $\Phi_x$ is surjective for all $x$.
	
	But on global sections,
	\[
	H^0(X,\mathcal O^{\oplus 2})\cong \C^2,\qquad H^0(X,\mathcal O(2))\cong \C^3
	\]
	with basis $\{Z_0^2, Z_0Z_1, Z_1^2\}$ for $H^0(\mathcal O(2))$.
	The image of $H^0(\Phi)$ is $\mathrm{span}\{Z_0^2, Z_1^2\}$, which misses $Z_0Z_1$.
	So stalkwise surjectivity does not force surjectivity on global sections.
\end{example}

\subsection{Examples and non-examples}
\label{subsec:sheaf-examples}

\subsubsection{Why the constant presheaf fails the sheaf axiom}

\begin{example}[Constant presheaf fails on disconnected opens]
	Let $M=U_1\sqcup U_2$ be a disjoint union of two nonempty open sets.
	Fix $A=\Z$ and define the constant presheaf $\underline{\Z}_{\mathrm{pre}}(U)=\Z$
	for all nonempty open $U$ with identity restrictions.
	Take the cover $M=U_1\cup U_2$. Choose sections $s_1=0\in\Z$ on $U_1$ and $s_2=1\in\Z$ on $U_2$.
	On the overlap $U_1\cap U_2=\varnothing$ the compatibility condition is vacuous, so these data should glue,
	but a section on $M$ in the presheaf must be a \emph{single} integer.
	There is no integer restricting to $0$ on $U_1$ and $1$ on $U_2$. Hence gluing fails.
\end{example}

\subsubsection{Skyscraper sheaf stalk calculation with representatives}

\begin{example}[Skyscraper stalks]
	Let $M$ be a $T_1$ space, fix $p\in M$, and let $A$ be an abelian group. Define $A_p(U)=A$ if $p\in U$ and $0$ otherwise, with the identity restriction when both open sets contain $p$ and the zero map otherwise.
	If $x\neq p$, choose $U\ni x$ with $p\notin U$; every germ at $x$ has a representative on this smaller neighborhood, where the only section is $0$.
	Thus $(A_p)_x=0$.
	If $x=p$, any neighborhood $U\ni p$ has $A_p(U)=A$ and all restrictions are identities,
	so germs at $p$ are determined by a unique element of $A$:
	\[
	(A_p)_p\cong A.
	\]
\end{example}

\subsection{Exercises}
\label{subsec:sheaf-exercises}

\begin{exercise}[Gluing smooth functions with a chart computation]
	Let $M$ be a smooth manifold and $U=\bigcup_\alpha U_\alpha$.
	Suppose $f_\alpha\in C^\infty(U_\alpha)$ and $f_\alpha=f_\beta$ on $U_{\alpha\beta}$.
	Define $f$ pointwise by $f(x)=f_\alpha(x)$ for $x\in U_\alpha$.
	Prove $f\in C^\infty(U)$ by fixing a chart $(V,\phi)$ and checking that
	$f\circ\phi^{-1}$ agrees on $\phi(V\cap U_\alpha)$ with smooth functions, hence is smooth.
\end{exercise}

\begin{exercise}[Locality of the Laplacian]
	Let $U\subset\R^n$ be open and $u,v\in C^2(U)$.
	Assume that for every $x\in U$ there is a neighborhood $V_x$ with $u|_{V_x}=v|_{V_x}$.
	Show $\Delta u=\Delta v$ on $U$.
\end{exercise}

\begin{exercise}[A concrete stalk in $\mathscr C^\infty$]
	Let $M=\R$ and $x=0$.
	Show that the germ of $f(x)=\sin x$ at $0$ is different from the germ of every
	finite Taylor polynomial of $\sin x$.  By contrast, its full convergent Taylor series
	represents the same germ.  Explain why agreement of a finite Taylor jet is strictly
	weaker than equality of germs.
\end{exercise}

\begin{exercise}[Patching obstruction for logarithms]
	On $U=\C^\times$, consider the cover $U_1=\C\setminus(-\infty,0]$, $U_2=\C\setminus[0,\infty)$.
	Write down explicit formulas for branches $\Log_1,\Log_2$ on $U_1,U_2$.
	Compute $\Log_1-\Log_2$ on the two connected components of $U_1\cap U_2$ and show it is constant
	but takes different constants on the two components. Conclude that $\{\Log_1,\Log_2\}$ cannot glue
	to a holomorphic function on $\C^\times$.
\end{exercise}

\begin{exercise}[Exactness on simply connected domains]
	Let $\Omega\subset\C$ be simply connected. Prove that the sequence of abelian groups
	\[
	0\to \Z \to \mathcal O(\Omega)\xrightarrow{\exp(2\pi i\,\cdot)} \mathcal O^\times(\Omega)\to 0
	\]
	is exact by showing surjectivity: given $h\in\mathcal O^\times(\Omega)$,
	define $g=h'/h$ and show $\int g$ is path-independent in $\Omega$.
\end{exercise}

\subsection{\v{C}ech Cohomology of a Sheaf}
\index{Cech cohomology}
\label{sec:cech-cohomology}

\subsubsection{Open covers, cochains, and the coboundary operator}
\index{Cech cohomology!cochains}\index{coboundary operator}
\label{subsec:cech-definitions}

\paragraph{Setup and notation}

Let $M$ be a topological space and $\mathscr F$ a sheaf of abelian groups on $M$.
Fix an open cover
\[
\mathcal U=\{U_\alpha\}_{\alpha\in A},\qquad M=\bigcup_{\alpha\in A}U_\alpha.
\]
For indices $\alpha_0,\dots,\alpha_q$ set
\[
U_{\alpha_0\cdots\alpha_q}:=\bigcap_{j=0}^q U_{\alpha_j}.
\]

\paragraph{Cochains and coboundary}

\begin{definition}[\v{C}ech cochains]
	Fix a total order on the index set $A$. For $q\ge 0$, use the normalized alternating convention
	and define
	\[
	C^q(\mathcal U,\mathscr F):=\prod_{\alpha_0<\cdots<\alpha_q}\mathscr F(U_{\alpha_0\cdots\alpha_q}).
	\]
	An element $c\in C^q$ is written $c=\{c_{\alpha_0\cdots\alpha_q}\}$ with
	$c_{\alpha_0\cdots\alpha_q}\in\mathscr F(U_{\alpha_0\cdots\alpha_q})$. Extend the notation
	to arbitrary index tuples by alternation, and set a component equal to zero when two indices coincide.
\end{definition}

\begin{definition}[Coboundary operator]
	Define $\delta:C^q(\mathcal U,\mathscr F)\to C^{q+1}(\mathcal U,\mathscr F)$ by
	\[
	(\delta c)_{\alpha_0\cdots\alpha_{q+1}}
	=
	\sum_{j=0}^{q+1}(-1)^j\,
	c_{\alpha_0\cdots\widehat{\alpha_j}\cdots\alpha_{q+1}}
	\Big|_{U_{\alpha_0\cdots\alpha_{q+1}}}.
	\]
\end{definition}

\paragraph{The identity $\delta^2=0$}

\begin{proposition}[The \v{C}ech differential squares to zero]
	\label{prop:cech-delta2}
	For all $q\ge 0$, one has $\delta\circ\delta=0$.
\end{proposition}

\begin{proof}
	Fix $c\in C^q(\mathcal U,\mathscr F)$ and indices $\alpha_0,\dots,\alpha_{q+2}$.
	Write $W:=U_{\alpha_0\cdots\alpha_{q+2}}$.
	By definition,
	\[
	(\delta^2 c)_{\alpha_0\cdots\alpha_{q+2}}
	=
	\sum_{i=0}^{q+2}(-1)^i\,(\delta c)_{\alpha_0\cdots\widehat{\alpha_i}\cdots\alpha_{q+2}}\Big|_{W}.
	\]
	Expand the inner term:
	\[
	(\delta c)_{\alpha_0\cdots\widehat{\alpha_i}\cdots\alpha_{q+2}}
	=
	\sum_{j=0}^{q+1}(-1)^j\,
	c_{\alpha_0\cdots\widehat{\alpha_i}\cdots\widehat{\alpha_{k(i,j)}}\cdots\alpha_{q+2}}
	\Big|_{U_{\alpha_0\cdots\widehat{\alpha_i}\cdots\alpha_{q+2}},
	}
	\]
	where $k(i,j)$ denotes the index in $\{0,\dots,q+2\}\setminus\{i\}$ occupying the $j$-th position.
	Restricting to $W$ and re-indexing, $\delta^2 c$ becomes a signed sum over all pairs
	$0\le i<j\le q+2$ of the same section
	\[
	c_{\alpha_0\cdots\widehat{\alpha_i}\cdots\widehat{\alpha_j}\cdots\alpha_{q+2}}\Big|_{W}.
	\]
	Each pair $(i,j)$ appears exactly twice: once by first removing $\alpha_i$ then $\alpha_j$,
	and once by first removing $\alpha_j$ then $\alpha_i$. The corresponding signs are opposite:
	\[
	(-1)^{i+j} \quad\text{and}\quad -(-1)^{i+j}.
	\]
	Hence these two contributions cancel. Since this holds for every $(i,j)$, the total sum is $0$ on $W$.
	Therefore $(\delta^2 c)_{\alpha_0\cdots\alpha_{q+2}}=0$ for all indices, i.e.\ $\delta^2=0$.
\end{proof}

\paragraph{Cohomology groups}

Define
\[
Z^q(\mathcal U,\mathscr F):=\ker\bigl(\delta:C^q\to C^{q+1}\bigr),
\qquad
B^q(\mathcal U,\mathscr F):=\mathrm{im}\bigl(\delta:C^{q-1}\to C^q\bigr),
\]
and
\[
\check H^q(\mathcal U,\mathscr F):=Z^q(\mathcal U,\mathscr F)/B^q(\mathcal U,\mathscr F).
\]

\subsubsection{Refinements and the direct limit construction}
\label{subsec:cech-refinements}

\paragraph{Refinements and pullback of cochains}

\begin{definition}[Refinement]
	A cover $\mathcal V=\{V_\beta\}_{\beta\in B}$ is a refinement of $\mathcal U=\{U_\alpha\}_{\alpha\in A}$
	if there exists a map $r:B\to A$ such that $V_\beta\subset U_{r(\beta)}$ for all $\beta$.
\end{definition}

\begin{definition}[Induced map on cochains]
	Given a refinement map $r:B\to A$, define
	\[
	r^\sharp:C^q(\mathcal U,\mathscr F)\longrightarrow C^q(\mathcal V,\mathscr F)
	\]
	by
	\[
	(r^\sharp c)_{\beta_0\cdots\beta_q}
	=
	c_{r(\beta_0)\cdots r(\beta_q)}\Big|_{V_{\beta_0\cdots\beta_q}}.
	\]
	Here the component on the right is interpreted by the alternating convention above; thus it
	vanishes if two of the indices $r(\beta_j)$ coincide.
\end{definition}

\begin{proposition}[Refinement maps are cochain maps]
	\label{prop:cech-refinement-chainmap}
	One has $\delta\circ r^\sharp=r^\sharp\circ\delta$. In particular, $r^\sharp$ induces homomorphisms
	\[
	r^\sharp:\check H^q(\mathcal U,\mathscr F)\to \check H^q(\mathcal V,\mathscr F).
	\]
\end{proposition}

\begin{proof}
	Fix $c\in C^q(\mathcal U,\mathscr F)$ and indices $\beta_0,\dots,\beta_{q+1}$.
	Let $W:=V_{\beta_0\cdots\beta_{q+1}}$.
	Compute:
	\[
	(\delta(r^\sharp c))_{\beta_0\cdots\beta_{q+1}}
	=
	\sum_{j=0}^{q+1}(-1)^j\,
	(r^\sharp c)_{\beta_0\cdots\widehat{\beta_j}\cdots\beta_{q+1}}\Big|_{W}.
	\]
	By definition of $r^\sharp$,
	\[
	(r^\sharp c)_{\beta_0\cdots\widehat{\beta_j}\cdots\beta_{q+1}}
	=
	c_{r(\beta_0)\cdots\widehat{r(\beta_j)}\cdots r(\beta_{q+1})}
	\Big|_{V_{\beta_0\cdots\widehat{\beta_j}\cdots\beta_{q+1}}}.
	\]
	Restricting to $W$ and using transitivity of restriction maps in a sheaf,
	\[
	(\delta(r^\sharp c))_{\beta_0\cdots\beta_{q+1}}
	=
	\sum_{j=0}^{q+1}(-1)^j\,
	c_{r(\beta_0)\cdots\widehat{r(\beta_j)}\cdots r(\beta_{q+1})}\Big|_{W}.
	\]
	On the other hand,
	\[
	(r^\sharp(\delta c))_{\beta_0\cdots\beta_{q+1}}
	=
	(\delta c)_{r(\beta_0)\cdots r(\beta_{q+1})}\Big|_{W}
	=
	\sum_{j=0}^{q+1}(-1)^j\,
	c_{r(\beta_0)\cdots\widehat{r(\beta_j)}\cdots r(\beta_{q+1})}\Big|_{W}.
	\]
	These coincide, hence $\delta\circ r^\sharp=r^\sharp\circ\delta$, and the induced maps on cohomology follow.
\end{proof}

\paragraph{Direct limit definition}

\begin{definition}[Global \v{C}ech cohomology]
	Define
	\[
	\check H^q(M,\mathscr F):=\varinjlim_{\mathcal U}\check H^q(\mathcal U,\mathscr F),
	\]
	where the directed system is over covers ordered by refinement, using the maps $r^\sharp$.
\end{definition}

Different choices of the refinement map $r$ induce chain-homotopic maps and hence the same
map on cohomology. Consequently the direct system depends on the refinement relation, not on
these auxiliary choices; see \cite{Bredon1997,Godement1958} for the standard prism homotopy.

\subsubsection{Degree zero and degree one: gluing and patching}
\label{subsec:cech-low-degree}

\paragraph{Degree zero equals global sections}

\begin{proposition}[Degree-zero \v{C}ech cohomology]
	\label{prop:cech-H0-cover}
	For any open cover $\mathcal U$, the group $\check H^0(\mathcal U,\mathscr F)$ is canonically isomorphic to $\mathscr F(M)$.
	Consequently $\check H^0(M,\mathscr F)\cong \mathscr F(M)$.
\end{proposition}

\begin{proof}
	A $0$--cochain is a family $s=\{s_\alpha\}$ with $s_\alpha\in\mathscr F(U_\alpha)$.
	The condition $\delta s=0$ means: for all $\alpha,\beta$,
	\[
	(\delta s)_{\alpha\beta}
	=
	s_\beta|_{U_{\alpha\beta}}-s_\alpha|_{U_{\alpha\beta}}=0,
	\]
	equivalently $s_\alpha$ and $s_\beta$ agree on overlaps. Thus $Z^0(\mathcal U,\mathscr F)$
	is exactly the set of compatible families $\{s_\alpha\}$.
	Since $B^0=\mathrm{im}(\delta:C^{-1}\to C^0)=0$ by definition, one has
	\[
	\check H^0(\mathcal U,\mathscr F)=Z^0(\mathcal U,\mathscr F).
	\]
	By the sheaf gluing axiom, any compatible family $\{s_\alpha\}$ arises from a unique
	global section $s\in\mathscr F(M)$ with $s|_{U_\alpha}=s_\alpha$. This defines a bijection
	\[
	\mathscr F(M)\longrightarrow Z^0(\mathcal U,\mathscr F),\qquad s\mapsto \{s|_{U_\alpha}\}.
	\]
	It is a group isomorphism because restrictions preserve addition. Hence $\check H^0(\mathcal U,\mathscr F)\cong \mathscr F(M)$.
	Passing to the direct limit over refinements preserves degree $0$, so $\check H^0(M,\mathscr F)\cong\mathscr F(M)$ as well.
\end{proof}

\paragraph{Degree one: cocycle and coboundary formulas}

A $1$--cochain is $g=\{g_{\alpha\beta}\}$ with $g_{\alpha\beta}\in\mathscr F(U_{\alpha\beta})$.
Then
\[
(\delta g)_{\alpha\beta\gamma}
=
g_{\beta\gamma}|_{U_{\alpha\beta\gamma}}
-
g_{\alpha\gamma}|_{U_{\alpha\beta\gamma}}
+
g_{\alpha\beta}|_{U_{\alpha\beta\gamma}}.
\]
Thus $g$ is a $1$--cocycle iff on every triple overlap,
\[
g_{\alpha\gamma}=g_{\alpha\beta}+g_{\beta\gamma}.
\]
A $1$--cocycle is a $1$--coboundary iff there exist $f=\{f_\alpha\}\in C^0$ with
\[
g_{\alpha\beta}=f_\beta|_{U_{\alpha\beta}}-f_\alpha|_{U_{\alpha\beta}}.
\]

\subsubsection{Examples and basic computations}
\label{subsec:cech-examples}

\paragraph{A computed example: $\check H^1(S^1,\underline{\Z})\cong\Z$}

Let $\underline{\Z}$ be the constant sheaf on $S^1$.
Choose an open cover by two arcs $\mathcal U=\{U,V\}$ such that $U\cap V$ has exactly two connected components,
say
\[
U\cap V = W_1 \sqcup W_2,
\qquad
W_1,W_2 \ \text{connected, open in }S^1.
\]
Then
\[
C^0(\mathcal U,\underline{\Z})=\underline{\Z}(U)\times \underline{\Z}(V)\cong \Z\times\Z
\]
because $U,V$ are connected, and
\[
C^1(\mathcal U,\underline{\Z})=\underline{\Z}(U\cap V)\cong \Z\times\Z
\]
because locally constant integer functions on $W_1\sqcup W_2$ are pairs of integers.

Since there is no triple intersection, $\delta:C^1\to C^2$ is the zero map, hence
\[
Z^1(\mathcal U,\underline{\Z})=C^1(\mathcal U,\underline{\Z})\cong \Z\times\Z.
\]
Now compute coboundaries. For $(a,b)\in C^0\cong\Z\times\Z$, one has
\[
\delta(a,b) = b|_{U\cap V}-a|_{U\cap V}.
\]
On each connected component $W_i$, the restrictions are constant, so
\[
\delta(a,b) = (b-a,\ b-a)\in \Z\times\Z.
\]
Therefore
\[
B^1(\mathcal U,\underline{\Z})
=
\{(n,n)\mid n\in\Z\}
\subset \Z\times\Z,
\]
and hence
\[
\check H^1(\mathcal U,\underline{\Z})
=
(\Z\times\Z)/\{(n,n)\}
\cong \Z,
\]
via the homomorphism
\[
\Z\times\Z\longrightarrow \Z,\qquad (x_1,x_2)\mapsto x_1-x_2,
\]
whose kernel is exactly the diagonal subgroup $\{(n,n)\}$.
One concludes $\check H^1(\mathcal U,\underline{\Z})\cong \Z$.
The cover is Leray for $\underline{\Z}$: its nonempty intersections are
arcs or disjoint unions of arcs, all with vanishing higher cohomology.
By Proposition~\ref{prop:leray}, this calculation gives
$H^1(S^1,\underline{\Z})\cong\Z$.  The standard comparison of sheaf and
global \v{C}ech cohomology on the paracompact circle then gives
\[
\check H^1(S^1,\underline{\Z})\cong \Z.
\]

\begin{example}[A concrete $1$--cocycle representative]
	Define $g\in C^1(\mathcal U,\underline{\Z})\cong\Z\times\Z$ by $g=(1,0)$, i.e.\
	$g|_{W_1}\equiv 1$ and $g|_{W_2}\equiv 0$.
	Then $g$ is a cocycle (no triple overlaps). Its class is nonzero because
	every coboundary is of the form $(n,n)$, and $(1,0)\notin\{(n,n)\}$.
	Under $(x_1,x_2)\mapsto x_1-x_2$, the class maps to $1\in\Z$.
\end{example}

\paragraph{A vanishing computation: $\check H^1(S^2,\underline{\Z})=0$}

\begin{example}[Two-set cover on $S^2$]
	Cover $S^2$ by $U=S^2\setminus\{\text{south pole}\}$ and $V=S^2\setminus\{\text{north pole}\}$.
	Then $U$ and $V$ are connected, and $U\cap V\cong S^1\times(0,1)$ is connected.
	Hence
	\[
	C^0(\mathcal U,\underline{\Z})\cong \Z\times\Z,\qquad C^1(\mathcal U,\underline{\Z})\cong \Z.
	\]
	As before $Z^1=C^1\cong\Z$ and $\delta(a,b)=b-a\in\Z$, so $B^1=\Z$.
	Thus $\check H^1(\mathcal U,\underline{\Z})=0$.  A calculation for one cover alone does
	not imply the global Čech group vanishes.  The conclusion
	$\check H^1(S^2,\underline{\Z})=0$ follows from a good cover (or from the
	Čech--singular comparison theorem).
\end{example}

\paragraph{Line bundles via $\check H^1(M,\mathcal O^\times)$}

Let $M$ be a complex manifold. Fix a cover $\mathcal U=\{U_\alpha\}$.

\begin{proposition}[Line bundles trivialized by a cover]
	\label{prop:cech-linebundles}
	There is a canonical bijection between:
	\begin{itemize}
		\item isomorphism classes of holomorphic line bundles on $M$ trivialized over $\mathcal U$,
		\item the group $\check H^1(\mathcal U,\mathcal O_M^\times)$.
	\end{itemize}
	Under refinement, these bijections are compatible, yielding a canonical bijection
	between isomorphism classes of holomorphic line bundles on $M$ and $\check H^1(M,\mathcal O_M^\times)$.
\end{proposition}

\begin{proof}
	Step 1 (from a line bundle to a cocycle).
	Let $L\to M$ be a holomorphic line bundle and choose holomorphic trivializations
	\[
	\varphi_\alpha: L|_{U_\alpha}\xrightarrow{\sim} U_\alpha\times\C.
	\]
	On overlaps $U_{\alpha\beta}$, the transition is a fiberwise $\C$--linear holomorphic automorphism of $\C$,
	hence multiplication by a unique $f_{\alpha\beta}\in \mathcal O^\times(U_{\alpha\beta})$:
	\[
	\varphi_\alpha\circ\varphi_\beta^{-1}(x,\eta)=(x,f_{\alpha\beta}(x)\eta).
	\]
	On triple overlaps $U_{\alpha\beta\gamma}$, associativity of composition gives
	\[
	f_{\alpha\beta}f_{\beta\gamma}=f_{\alpha\gamma},
	\]
	so $\{f_{\alpha\beta}\}$ is a $1$--cocycle in $Z^1(\mathcal U,\mathcal O^\times)$.
	Thus $L$ determines a class $[f]\in\check H^1(\mathcal U,\mathcal O^\times)$.
	
	Step 2 (effect of changing trivializations).
	Replace $\varphi_\alpha$ by $\widetilde\varphi_\alpha$ where
	\[
	\widetilde\varphi_\alpha = (x,\eta)\mapsto (x,g_\alpha(x)\eta)\circ \varphi_\alpha
	\quad\text{for some } g_\alpha\in\mathcal O^\times(U_\alpha).
	\]
	Then on $U_{\alpha\beta}$,
	\[
	\widetilde f_{\alpha\beta}=g_\alpha f_{\alpha\beta} g_\beta^{-1}.
	\]
	Equivalently, in multiplicative notation, $\widetilde f = f\cdot \delta(g)^{-1}$, so $f$ and $\widetilde f$
	define the same element of $\check H^1(\mathcal U,\mathcal O^\times)$.
	Hence the class $[f]$ depends only on $L$ up to isomorphism.
	
	Step 3 (from a cocycle to a line bundle).
	Given $f=\{f_{\alpha\beta}\}\in Z^1(\mathcal U,\mathcal O^\times)$, form the disjoint union
	$\bigsqcup_\alpha (U_\alpha\times\C)$ and define an equivalence relation by
	\[
	(x,\eta)_\beta \sim (x,f_{\alpha\beta}(x)\eta)_\alpha
	\qquad (x\in U_{\alpha\beta}).
	\]
	The cocycle identity ensures transitivity on triple overlaps, so the quotient
	\[
	L_f := \Bigl(\bigsqcup_\alpha (U_\alpha\times\C)\Bigr)\big/\!\sim
	\]
	is a complex manifold, and the natural projection $L_f\to M$ is a holomorphic line bundle
	trivialized over $\mathcal U$ with transition functions $f_{\alpha\beta}$.
	
	Step 4 (coboundaries give isomorphic bundles).
	If $f$ and $f'$ differ by a coboundary, i.e.\ $f'_{\alpha\beta}=g_\alpha f_{\alpha\beta}g_\beta^{-1}$
	for some $g_\alpha\in\mathcal O^\times(U_\alpha)$, then the maps
	\[
	U_\alpha\times\C \to U_\alpha\times\C,\qquad (x,\eta)\mapsto (x,g_\alpha(x)\eta)
	\]
	descend to a well-defined holomorphic bundle isomorphism $L_f\simeq L_{f'}$.
	
	Step 5 (inverse constructions).
	Composing Step 1 then Step 3 recovers a bundle isomorphic to the original;
	composing Step 3 then Step 1 recovers the original cocycle class. Therefore we have a bijection.
	
	Compatibility under refinement follows by restricting trivializations and transition functions,
	which matches the refinement pullback $r^\sharp$ on $\check H^1$.
\end{proof}

\begin{example}[A hands-on cocycle computation on $\mathbb{CP}^1$ for $\mathcal O(1)$]
	Let $\mathbb{CP}^1$ be covered by $U_0=\{Z_0\neq 0\}$ and $U_1=\{Z_1\neq 0\}$ with coordinates
	$z=Z_1/Z_0$ on $U_0$ and $w=Z_0/Z_1=1/z$ on $U_1$.
	Define a cocycle $f_{01}=z\in\mathcal O^\times(U_{01})$ and $f_{10}=1/z$.
	Then $\{f_{01},f_{10}\}$ is a $1$--cocycle (no triple overlaps).
	The glued bundle $L_f$ is $\mathcal O(1)$.
	
	A holomorphic section corresponds to holomorphic functions $s_0\in\mathcal O(U_0)$ and
	$s_1\in\mathcal O(U_1)$ satisfying
	\[
	s_0 = f_{01}\,s_1 \quad\text{on }U_{01}.
	\]
	Write $s_1(w)=a+bw+cw^2+\cdots$ (entire on $\C$), then on $U_{01}$,
	\[
	s_0(z)=z\,s_1(1/z)=z\Bigl(a+\frac{b}{z}+\frac{c}{z^2}+\cdots\Bigr)=az+b+\frac{c}{z}+\cdots.
	\]
	For $s_0$ to extend holomorphically to all of $U_0\cong\C$, all negative powers must vanish,
	so $c=d=\cdots=0$. Hence $s_1(w)=a+bw$ and $s_0(z)=az+b$.
	Therefore
	\[
	H^0(\mathbb{CP}^1,\mathcal O(1))\cong \C^2
	\]
	with basis corresponding to $(s_0,s_1)=(1,w)$ and $(z,1)$.
\end{example}

\begin{example}[Canonical sheaf of $\mathbb{CP}^1$ by a two-chart calculation]
	\label{ex:cp1-canonical-cech}
	Let $U_0\cong\C$ have coordinate $z$ and $U_1\cong\C$ have coordinate
	$w=1/z$. On $U_0\cap U_1\cong\C^*$,
	\[
	dw=-z^{-2}\,dz.
	\]
	The sets $U_0,U_1$, and $U_0\cap U_1$ are Stein, so this cover is Leray
	for the coherent sheaf $K=\Omega^1_{\mathbb{CP}^1}$. Its \v{C}ech complex is
	\[
	0\longrightarrow
	\mathcal O(\C)\,dz\oplus\mathcal O(\C)\,dw
	\xrightarrow{\delta}
	\mathcal O(\C^*)\,dz\longrightarrow0,
	\]
	where, after writing everything in the $z$-coordinate,
	\[
	\delta(f(z)\,dz,g(w)\,dw)
	=\bigl(-z^{-2}g(1/z)-f(z)\bigr)\,dz.
	\]

	First compute $H^0$. A global holomorphic one-form would satisfy
	$f(z)=-z^{-2}g(1/z)$ with both $f$ and $g$ entire. Thus $f(z)=O(|z|^{-2})$
	as $z\to\infty$; Liouville's theorem gives $f=0$, and then $g=0$. Hence
	\[
	H^0(\mathbb{CP}^1,K)=0.
	\]

	For $H^1$, take $h(z)\,dz\in\mathcal O(\C^*)\,dz$ and write its Laurent
	expansion
	\[
	h(z)=\sum_{n\in\Z}a_nz^n.
	\]
	The nonnegative powers are the restriction of the entire function
	$\sum_{n\ge0}a_nz^n$ on $U_0$. The powers $z^n$ with $n\le-2$ have the
	form $-z^{-2}g(1/z)$ for an entire function $g$ on $U_1$. After subtracting
	a coboundary, only $a_{-1}z^{-1}dz$ remains. This term cannot be a
	coboundary, because neither entire contribution supplies the exponent
	$-1$. Therefore
	\[
	H^1(\mathbb{CP}^1,K)\cong\C,
	\qquad [z^{-1}dz]\longmapsto1.
	\]
	This single calculation agrees with $K\simeq\mathcal O(-2)$, with
	Riemann--Roch in genus zero, and with Serre duality
	$H^1(K)\simeq H^0(\mathcal O)^\vee$.
\end{example}

\subsubsection{Exercises}
\label{subsec:cech-exercises}

\begin{exercise}[Compute $\check H^1$ for a two-set cover of $S^1$]
	\label{exr:selected-cech-circle}
	Let $\mathcal U=\{U,V\}$ be a two-arc cover of $S^1$ with $U\cap V=W_1\sqcup W_2$.
	Compute $C^0$, $C^1$, $B^1$, and $\check H^1(\mathcal U,\underline{\Z})$ explicitly.
\end{exercise}

\begin{exercise}[A refinement map on cochains]
	Given a refinement $r:\mathcal V\to\mathcal U$, verify directly from definitions that
	$r^\sharp$ commutes with $\delta$ in degree $q=0$ and $q=1$.
\end{exercise}

\begin{exercise}[Coboundary criterion in degree one]
	Let $g\in Z^1(\mathcal U,\mathscr F)$. Prove that $[g]=0$ in $\check H^1(\mathcal U,\mathscr F)$
	iff there exist $f_\alpha\in\mathscr F(U_\alpha)$ such that $g_{\alpha\beta}=f_\beta-f_\alpha$ on $U_{\alpha\beta}$.
\end{exercise}

\begin{exercise}[Gauge change and cocycle class]
	For $\mathscr F=\mathcal O^\times$, interpret $f'_{\alpha\beta}=g_\alpha f_{\alpha\beta} g_\beta^{-1}$
	as a change of local frames on the glued line bundle, and construct explicitly the induced bundle isomorphism.
\end{exercise}

\begin{exercise}[A vanishing computation on $S^2$]
	Use the two-set cover $S^2=U\cup V$ with $U\cap V$ connected to compute
	$\check H^1(\{U,V\},\underline{\Z})$.  Then compute the sheaf cohomology
	$H^1(S^2,\underline{\Z})=0$ using a good refinement (or the comparison with singular
	cohomology).  Explain why the calculation for one cover alone does not justify the
	global conclusion unless an appropriate Leray or refinement argument is supplied.
\end{exercise}

\subsection{From \v{C}ech Cohomology to Sheaf Cohomology}
\index{sheaf cohomology}

For comparison with the derived-functor construction and the precise hypotheses
under which \v{C}ech cohomology computes sheaf cohomology, see
\cite{StacksCohomology2026,Bredon1997,Godement1958}.
\label{sec:from-cech-to-sheaf}

The circle calculation in Subsection~\ref{subsec:cech-examples} used a
particular two-set cover.  Why did that calculation give a cohomology
group belonging to the circle, rather than an artifact of the cover?
There are two distinct steps.  Refining covers produces a cover-independent
\v{C}ech group.  To identify the cohomology of one chosen cover directly
with sheaf cohomology, one additionally needs vanishing of higher
cohomology on its nonempty finite intersections: this is the Leray
condition below.  The two arcs and their intersection components in the
circle example satisfy it for the constant sheaf.  For the holomorphic
structure sheaf, merely choosing contractible intersections is not enough;
the relevant analytic input is acyclicity, often supplied by Stein
intersections and Cartan's theorem~B.  This distinction explains both why
the explicit circle computation works and why injective resolutions are
introduced after the cover-based construction.

\subsubsection{Global sections and the derived-functor viewpoint}
\label{subsec:global-sections-derived}

Sheaf cohomology measures a precise failure of globalization. A morphism of
sheaves can be surjective on every stalk without being surjective on global
sections: local lifts may exist but fail to agree on overlaps. The first
cohomology group records this first obstruction, while the higher groups record
successive compatibility obstructions.

\begin{definition}[Ringed spaces and module sheaves]
	A \emph{ringed space} is a pair $(X,\mathcal O_X)$ consisting of a topological
	space and a sheaf of rings. An $\mathcal O_X$-module is a sheaf $\mathscr F$
	whose sections $\mathscr F(U)$ are modules over $\mathcal O_X(U)$, compatibly
	with restriction. For a Riemann surface, the principal examples are
	$\mathcal O_X$, the sheaf $\Omega_X^1$ of holomorphic one-forms, and the sheaf
	of holomorphic sections of a line bundle.
\end{definition}

\begin{proposition}[Global sections are left exact]
	\label{prop:Gamma-left-exact}
	If
	\[
	0\longrightarrow\mathscr F\longrightarrow\mathscr G
	\longrightarrow\mathscr H\longrightarrow0
	\]
	is a short exact sequence of sheaves of abelian groups (or
	$\mathcal O_X$-modules), then
	\[
	0\longrightarrow\Gamma(X,\mathscr F)
	\longrightarrow\Gamma(X,\mathscr G)
	\longrightarrow\Gamma(X,\mathscr H)
	\]
	is exact. The last displayed arrow need not be surjective.
\end{proposition}

\begin{proof}
	Injectivity at the first term can be checked stalkwise. Suppose
	$g\in\Gamma(X,\mathscr G)$ maps to zero in $\Gamma(X,\mathscr H)$.
	For every $x\in X$, exactness on stalks gives a neighborhood $U_x$ and a
	section $f_x\in\mathscr F(U_x)$ mapping to $g|_{U_x}$. On overlaps the
	sections $f_x$ and $f_y$ have the same image in $\mathscr G$; the
	monomorphism $\mathscr F\to\mathscr G$ forces them to agree. The sheaf
	gluing axiom therefore gives a unique $f\in\Gamma(X,\mathscr F)$ mapping to
	$g$. Surjectivity onto $\Gamma(X,\mathscr H)$ can fail because arbitrary
	local lifts of a global section of $\mathscr H$ need not agree on overlaps.
\end{proof}

\begin{definition}[Complexes, homotopies, and right derived functors]
	A cochain complex $A^\bullet$ has maps $d^q:A^q\to A^{q+1}$ with
	$d^{q+1}d^q=0$, and
	$H^q(A^\bullet)=\ker d^q/\operatorname{im}d^{q-1}$. Two cochain maps
	$f,g:A^\bullet\to B^\bullet$ are \emph{cochain homotopic} if there are maps
	$h^q:A^q\to B^{q-1}$ such that $f-g=dh+hd$; homotopic maps induce the same
	map on cohomology.

	If $F$ is a left exact functor on an abelian category with enough
	injectives and $A\to I^\bullet$ is an injective resolution, its right
	derived functors are
	\[
	R^qF(A):=H^q(F(I^\bullet)).
	\]
	The comparison theorem for injective resolutions makes this independent of
	the chosen resolution and functorial in $A$; moreover $R^0F\simeq F$.
\end{definition}

For $F=\Gamma(X,-)$ this construction gives sheaf cohomology. The notation
$H^q(X,\mathscr F)=R^q\Gamma(X,\mathscr F)$ will be justified explicitly
below. A short exact sequence of sheaves then produces the long exact
cohomology sequence discussed in Section~\ref{sec:ses-and-les}.

\begin{remark}[Underlying abelian sheaf]
	For a ringed space $(X,\mathcal O_X)$, one may derive global sections either
	in the category of $\mathcal O_X$-modules or after forgetting to sheaves of
	abelian groups. The resulting groups agree naturally. Indeed, an
	injective $\mathcal O_X$-module is flabby and hence acyclic as an abelian
	sheaf; an injective $\mathcal O_X$-module resolution is therefore an
	acyclic resolution for the latter calculation. The common group inherits
	a natural $\Gamma(X,\mathcal O_X)$-module structure
	\cite{StacksCohomology2026,Godement1958}.
\end{remark}

\subsubsection{Acyclic sheaves and Leray covers}
\index{acyclic sheaves}\index{Leray covers}
\label{subsec:acyclic}

\paragraph{Acyclic sheaves}

\begin{definition}[Acyclic sheaf]
	A sheaf $\mathscr F$ of abelian groups on a topological space $M$ is called
	\emph{acyclic} if
	\[
	H^q(M,\mathscr F)=0 \qquad \text{for all } q\ge 1,
	\]
	where $H^q(M,\mathscr F)$ denotes sheaf cohomology.
\end{definition}

\begin{definition}[\v{C}ech--acyclic on an open set]
	Let $U\subset M$ be open. The sheaf $\mathscr F$ is said to be
	\emph{\v{C}ech--acyclic on $U$} if
	\[
	\check H^q(U,\mathscr F)=0 \qquad \text{for all } q\ge 1.
	\]
\end{definition}

\paragraph{Leray covers}

\begin{definition}[Leray cover]
	An open cover $\mathcal U=\{U_\alpha\}$ of $M$ is called a \emph{Leray cover}
	for a sheaf $\mathscr F$ if every nonempty finite intersection
	$U_{\alpha_0\cdots\alpha_p}$ is acyclic for sheaf cohomology, i.e.
	\[
	H^q(U_{\alpha_0\cdots\alpha_p},\mathscr F)=0
	\quad \text{for all } q\ge 1.
	\]
\end{definition}

\begin{definition}[Sheafified \v{C}ech complex]
	For an open cover $\mathcal U=\{U_\alpha\}$ and an inclusion
	$j_I:U_I=U_{\alpha_0\cdots\alpha_p}\hookrightarrow X$, set
	\[
	\mathcal C^p(\mathcal U,\mathscr F)
	=\prod_{\alpha_0<\cdots<\alpha_p}
	(j_I)_*\bigl(\mathscr F|_{U_I}\bigr).
	\]
	Restriction and the alternating formula define a differential
	$\delta:\mathcal C^p\to\mathcal C^{p+1}$, and restriction of a section to
	the members of the cover gives an augmentation
	$\mathscr F\to\mathcal C^0(\mathcal U,\mathscr F)$.
\end{definition}

\begin{proposition}[Exactness of the sheafified \v{C}ech resolution]
	\label{prop:sheafified-cech-exact}
	The augmented complex
	\[
	0\longrightarrow\mathscr F\longrightarrow
	\mathcal C^0(\mathcal U,\mathscr F)\longrightarrow
	\mathcal C^1(\mathcal U,\mathscr F)\longrightarrow\cdots
	\]
	is exact.
\end{proposition}

\begin{proof}
	Exactness is local, so fix $x\in X$ and choose an index $a$ with
	$x\in U_a$. On the stalk at $x$, define a degree $-1$ operator by inserting
	$a$ as the first index (and reordering with the resulting sign):
	\[
	(hc)_{\alpha_0\ldots\alpha_{p-1}}
	=c_{a\alpha_0\ldots\alpha_{p-1}}.
	\]
	Only intersections containing $x$ occur in this stalk, so every displayed
	germ is defined. Expanding the alternating sums gives
	$\delta h+h\delta=\mathrm{id}$ in positive degrees. In degree zero the same
	calculation says that the kernel of $\delta$ consists exactly of compatible
	local germs, which are the image of $\mathscr F_x$. Hence the augmented
	complex is exact on every stalk and therefore exact as a complex of sheaves.
\end{proof}

\begin{definition}[First-quadrant spectral sequence]
 \label{def:double-complex-spectral-sequence}
 A cohomological spectral sequence consists of bigraded pages
 $E_r^{p,q}$ and differentials
 $d_r:E_r^{p,q}\to E_r^{p+r,q-r+1}$, with $d_r^2=0$ and
 $E_{r+1}^{p,q}=\ker(d_r:E_r^{p,q}\to E_r^{p+r,q-r+1})/
 \operatorname{im}(d_r:E_r^{p-r,q+r-1}\to E_r^{p,q})$.
 For a first-quadrant double complex $K^{p,q}$, filtering its total
 complex by horizontal degree yields
 $E_1^{p,q}=H^q(K^{p,\bullet})$ and
 $E_2^{p,q}=H^p(E_1^{\bullet,q})$.  Convergence identifies
 $E_\infty^{p,q}$ with the corresponding graded piece of a finite
 filtration of $H^{p+q}(\operatorname{Tot}K)$.
 The filtration construction and convergence proof may be found in
 \cite[Chapter~14]{BottTu1982}.
\end{definition}

\begin{proposition}[Leray comparison]
	\label{prop:leray}
	If $\mathcal U$ is a Leray cover for $\mathscr F$, then
	\[
	\check H^q(\mathcal U,\mathscr F)\cong H^q(M,\mathscr F)
	\quad \text{for all } q\ge 0.
	\]
\end{proposition}

\begin{proof}
	Choose an injective resolution $\mathscr F\to\mathscr I^\bullet$ and form the
	first-quadrant double complex
	\[
	K^{p,q}=\prod_{\alpha_0<\cdots<\alpha_p}
	\mathscr I^q(U_{\alpha_0\cdots\alpha_p}),
	\]
	with the \v{C}ech differential horizontally and the resolution differential
	vertically.  Its two filtrations give the comparison.
	
	For fixed $q$, the augmented horizontal complex begins with
	$\mathscr I^q(M)$ and is exact: an injective sheaf is flabby, so its
	sections can be extended across the cover, and the usual \v{C}ech
	contraction applies.  Thus the total complex has cohomology
	$H^\bullet(M,\mathscr F)$.
	
	For fixed $p$, taking vertical cohomology gives
	\[
	H^q(K^{p,\bullet})
	\cong \prod_{\alpha_0<\cdots<\alpha_p}
	H^q(U_{\alpha_0\cdots\alpha_p},\mathscr F).
	\]
	Here products are exact in abelian groups.  The Leray hypothesis makes
	this group zero for $q>0$, while the $q=0$ row is precisely
	$\check C^p(\mathcal U,\mathscr F)$.  The vertical spectral sequence
	therefore collapses and identifies the total cohomology with
	$\check H^\bullet(\mathcal U,\mathscr F)$, as required.
\end{proof}

\begin{example}[Good covers on manifolds]
	For a locally constant sheaf on a manifold, a good cover is Leray.  For the
	holomorphic structure sheaf, contractibility alone is insufficient; one uses a cover
	with Stein finite intersections and invokes Cartan's theorem~B.
\end{example}

\subsubsection{Sheaf cohomology via injective resolutions}
\label{subsec:injective-resolutions}

\paragraph{Injective sheaves}

\begin{definition}[Injective sheaf]
	A sheaf $\mathscr I$ of abelian groups on $M$ is \emph{injective}
	if for every monomorphism of sheaves $\mathscr A\hookrightarrow \mathscr B$,
	the induced map
	\[
	\mathrm{Hom}(\mathscr B,\mathscr I)\to \mathrm{Hom}(\mathscr A,\mathscr I)
	\]
	is surjective.
\end{definition}

\begin{theorem}[Existence of injective resolutions]
	\label{thm:injective-resolution}
	Every sheaf $\mathscr F$ of abelian groups admits an injective resolution
	\[
	0\longrightarrow \mathscr F
	\longrightarrow \mathscr I^0
	\longrightarrow \mathscr I^1
	\longrightarrow \mathscr I^2
	\longrightarrow \cdots.
	\]
\end{theorem}

\begin{proof}
	The category of sheaves of abelian groups on $M$ has enough injectives.
	Concretely, embed each stalk $\mathscr F_x$ into an injective abelian group $I_x$;
	the resulting canonical map embeds $\mathscr F$ into a product of pushforward
	(skyscraper-type) sheaves built from the $I_x$, and these target sheaves are injective.
	Repeating the construction for successive cokernels gives the resolution.
\end{proof}

\paragraph{Derived functor definition}

\begin{definition}[Sheaf cohomology]
	Let $\mathscr F$ be a sheaf and
	\[
	0\to \mathscr F \to \mathscr I^0 \to \mathscr I^1 \to \cdots
	\]
	an injective resolution.
	Define
	\[
	H^q(M,\mathscr F)
	:=R^q\Gamma(M,\mathscr F):=
	H^q\bigl(\Gamma(M,\mathscr I^\bullet)\bigr).
	\]
\end{definition}

\begin{proposition}[Independence of the injective resolution]
	This definition is independent of the chosen injective resolution.
\end{proposition}

\begin{proof}
	Let $\mathscr I^\bullet$ and $\mathscr J^\bullet$ be two injective resolutions of
	$\mathscr F$.  The injective-resolution comparison theorem
	\cite{Bredon1997,Godement1958} lifts $\mathrm{id}_{\mathscr F}$ to a cochain map
	$f:\mathscr I^\bullet\to\mathscr J^\bullet$.  In degree zero the lift is obtained
	by extending the map from the embedded copy of $\mathscr F$ into the injective
	sheaf $\mathscr J^0$; at successive degrees, injectivity extends the map from
	the image of the preceding differential.  The same theorem gives a lift
	$g:\mathscr J^\bullet\to\mathscr I^\bullet$, and states that any two lifts of
	the same map on $\mathscr F$ are cochain homotopic.  Thus $gf$ and $fg$ are
	homotopic to the respective identity maps.
	Applying $\Gamma(M,-)$ preserves the cochain-map equations and the homotopy
	equations because it is additive; it need not be exact for this step.
	Homotopic maps induce the same map on cohomology, so $\Gamma(f)$ and
	$\Gamma(g)$ induce inverse isomorphisms.  The uniqueness-up-to-homotopy
	clause also shows that the isomorphism does not depend on which lift $f$
	was chosen.
\end{proof}

\begin{theorem}[Acyclic resolutions compute sheaf cohomology]
	\label{thm:acyclic-resolution}
	Suppose
	\[
	0\longrightarrow\mathscr F\longrightarrow\mathscr A^0
	\longrightarrow\mathscr A^1\longrightarrow\cdots
	\]
	is exact and every $\mathscr A^p$ is acyclic for $\Gamma(M,-)$. Then there
	are natural isomorphisms
	\[
	H^q(M,\mathscr F)\cong
	H^q\bigl(\Gamma(M,\mathscr A^\bullet)\bigr),
	\qquad q\ge0.
	\]
\end{theorem}

\begin{proof}
	Put $\mathscr Z^0=\mathscr F$ and, successively,
	$\mathscr Z^{p+1}=\operatorname{coker}(\mathscr Z^p\to\mathscr A^p)$.
	The resolution splits into short exact sequences
	\[
	0\to\mathscr Z^p\to\mathscr A^p\to\mathscr Z^{p+1}\to0.
	\]
	Because $H^r(M,\mathscr A^p)=0$ for $r>0$, the associated long exact
	sequence gives the dimension-shifting isomorphisms
	$H^r(M,\mathscr Z^{p+1})\cong H^{r+1}(M,\mathscr Z^p)$ for $r>0$.
	In degree zero, the same long exact sequences identify kernels and
	cokernels of the maps $\Gamma(M,\mathscr A^p)\to
	\Gamma(M,\mathscr A^{p+1})$. Iterating these identifications yields the
	stated isomorphism. Naturality follows from the comparison theorem for
	resolutions \cite{StacksCohomology2026,Godement1958}.
\end{proof}

\subsubsection{Recovering \v{C}ech cohomology from derived functors}
\label{subsec:derived-vs-cech}

\begin{theorem}[\v{C}ech cohomology from derived functors]
	\label{thm:cech-derived}
	Let $\mathscr F$ be a sheaf on $M$ admitting a Leray cover $\mathcal U$.
	Then for all $q\ge 0$,
	\[
	\check H^q(\mathcal U,\mathscr F)\cong H^q(M,\mathscr F).
	\]
\end{theorem}

\begin{proof}
	Proposition~\ref{prop:leray} compares the \v{C}ech complex with an
	injective resolution by means of their double complex.  The vertical
	cohomology above degree zero vanishes on every finite intersection, so
	both complexes compute the same cohomology.
\end{proof}

\begin{example}[Agreement for smooth manifolds]
	On a smooth manifold, for $\mathscr F=\underline{\R}$ or $\Omega^k$,
	a good cover is Leray. Therefore \v{C}ech cohomology and sheaf cohomology coincide,
	justifying explicit \v{C}ech computations.
\end{example}

\begin{remark}[Three objects that should not be conflated]
	The groups $\check H^q(\mathcal U,\mathscr F)$ for one cover,
	$\check H^q(X,\mathscr F)$ obtained after refinement, and
	$H^q(X,\mathscr F)=R^q\Gamma(X,\mathscr F)$ are different constructions.
	For every cover, degree zero equals global sections. A Leray cover computes
	the derived group in every degree. Without such a hypothesis, the canonical
	map from global \v{C}ech cohomology to derived cohomology is always an
	isomorphism in degrees zero and one and is injective in degree two; stronger
	statements require additional hypotheses on the space or the cover
	\cite{Godement1958,StacksCohomology2026}.
\end{remark}

\begin{remark}[The topology is part of the input]
	Cohomology depends on the topology on which the sheaf is defined. On a
	complex manifold with its analytic topology, the constant sheaf detects
	topological cycles and is related to singular and de Rham cohomology. On an
	irreducible algebraic variety with the Zariski topology, any two nonempty
	opens meet; the constant sheaf is flabby, so its higher Zariski cohomology
	vanishes. Thus a symbol such as $H^q(X,\mathscr F)$ is incomplete unless the
	analytic, Zariski, or other site is understood.
\end{remark}

\subsubsection{Fine, soft, and flabby sheaves}
\label{subsec:fine-soft-flabby}

\paragraph{Definitions}

\begin{definition}[Flabby sheaf]
	A sheaf $\mathscr F$ is \emph{flabby} if all restriction maps
	$\mathscr F(U)\to \mathscr F(V)$ for $V\subset U$ are surjective.
\end{definition}

\begin{definition}[Soft sheaf]
	A sheaf $\mathscr F$ on a paracompact space $M$ is \emph{soft}
	if for every closed $K\subset M$, the restriction
	\[
	\mathscr F(M)\to \mathscr F(K)
	\]
	is surjective.
\end{definition}

\begin{definition}[Fine sheaf]
	A sheaf $\mathscr F$ of $\mathcal C^\infty$-modules on a smooth manifold
	is \emph{fine} if, for every locally finite open cover $\{U_i\}$, there are
	sheaf endomorphisms $h_i:\mathscr F\to\mathscr F$, locally finite as a family,
	with $\operatorname{supp}h_i\subset U_i$ and $\sum_i h_i=\mathrm{id}$.
	Multiplication by a smooth partition of unity supplies these endomorphisms.
\end{definition}

\paragraph{Vanishing results}

\begin{proposition}[Standard acyclic sheaves]
	\label{prop:acyclic-types}
	If $\mathscr F$ is flabby, soft (on a paracompact Hausdorff space), or fine,
	then $\mathscr F$ is acyclic.
\end{proposition}

\begin{proof}
	We first treat the flabby case. We use the extension lemma: if
	$0\to\mathscr F\to\mathscr G\to\mathscr H\to0$ is exact and $\mathscr F$
	is flabby, then $\Gamma(U,\mathscr G)\to\Gamma(U,\mathscr H)$ is surjective
	for every open $U$. To prove it, start with local lifts of a section of
	$\mathscr H(U)$ and order compatible lifts by extension of their domains.
	A maximal compatible lift exists by the usual union argument. If its domain
	is not all of $U$, choose another local lift; on the overlap the difference
	is a section of $\mathscr F$, which extends by flabbiness and can be
	subtracted to make the two lifts agree. They then glue, contradicting
	maximality.

	Embed a flabby $\mathscr F$ into an injective sheaf $\mathscr I$ and put
	$\mathscr Q=\mathscr I/\mathscr F$. The extension lemma makes
	$\Gamma(M,\mathscr I)\to\Gamma(M,\mathscr Q)$ surjective. Moreover
	$\mathscr Q$ is flabby: lift a section on a smaller open set to
	$\mathscr I$, extend it there, and pass to the quotient. Since $\mathscr I$
	is injective, the long exact sequence gives $H^1(M,\mathscr F)=0$ and
	$H^q(M,\mathscr F)\cong H^{q-1}(M,\mathscr Q)$ for $q\ge2$. Repeating the
	argument proves the vanishing in all positive degrees
	\cite{StacksCohomology2026,Godement1958}.

	Soft sheaves on a paracompact Hausdorff space are acyclic by the analogous
	extension argument over closed sets; paracompactness supplies locally finite
	refinements \cite{Bredon1997,Godement1958}. For the fine case, choose a locally finite cover
	$\{U_i\}$ and fine endomorphisms $h_i$ subordinate to it, so that
	$\sum_i h_i=\mathrm{id}$.  For a \v{C}ech $q$-cochain $c$, $q>0$, set
	\[
	(Hc)_{i_0\ldots i_{q-1}}
	=\sum_j h_j\bigl(c_{j i_0\ldots i_{q-1}}\bigr)
	\quad\text{on }U_{i_0}\cap\cdots\cap U_{i_{q-1}}.
	\]
	Each term extends by zero outside $U_j$ because $h_j$ is supported there,
	and local finiteness makes the sum meaningful.  Expanding the alternating
	\v{C}ech differential, the terms involving an index $i_k$ cancel in
	$\delta H+H\delta$; the remaining sum is
	$\sum_j h_j(c_{i_0\ldots i_q})=c_{i_0\ldots i_q}$.  Thus the positive-degree
	\v{C}ech complex is contractible.  The standard comparison of \v{C}ech and
	derived-functor cohomology on a paracompact manifold then gives
	$H^q(M,\mathscr F)=0$ for $q>0$ \cite{Bredon1997,Godement1958}.
\end{proof}

\begin{example}[de Rham complex]
	On a smooth manifold $M$, each $\Omega^k_M$ is a fine sheaf.
	Hence the de Rham complex is a fine resolution of $\underline{\R}$,
	leading to the de Rham theorem.
\end{example}

\subsubsection{Open Riemann surfaces, Stein manifolds, and the affine analogy}
\label{subsec:open-riemann-stein}

The preceding vanishing results are formal once an acyclic resolution has
been found.  In complex analysis there is a deeper source of acyclicity:
the geometry of Stein manifolds.  For Riemann surfaces, compactness is the
decisive boundary.  A connected compact Riemann surface has only constant
holomorphic functions, whereas a connected noncompact Riemann surface has
enough holomorphic functions both to separate points and to control behavior
at infinity.

\begin{definition}[Holomorphic hull and Stein manifold]
	Let $X$ be a complex manifold and $K\subset X$ compact.  Its
	\emph{holomorphic hull} is
	\[
		\widehat K_{\mathcal O(X)}
		=
		\{x\in X: |f(x)|\le \sup_K|f|\text{ for every }f\in\mathcal O(X)\}.
	\]
	The manifold $X$ is \emph{holomorphically convex} if this hull is compact
	for every compact $K$.  It is \emph{holomorphically separable} if, whenever
	$x\ne y$, some $f\in\mathcal O(X)$ satisfies $f(x)\ne f(y)$.  A
	second-countable complex manifold is called a \emph{Stein manifold} if it
	is holomorphically convex and holomorphically separable.
\end{definition}

Equivalent definitions use a proper strictly plurisubharmonic exhaustion or
a proper holomorphic embedding into some $\C^N$.  In complex dimension one,
strict plurisubharmonicity is strict subharmonicity.  The hull definition is
best suited to function theory; the exhaustion definition is best suited to
solving the $\bar\partial$-equation.

\begin{theorem}[Behnke--Stein theorem]
	\label{thm:open-riemann-stein}
	Every connected noncompact Riemann surface is a Stein manifold.
\end{theorem}

\begin{proof}
	We give the cohomological proof architecture, including the analytic input
	that cannot be replaced by a formal sheaf argument.  Choose a smooth
	exhaustion
	\[
		K_1\Subset K_2\Subset\cdots\Subset X,
		\qquad X=\bigcup_{j\ge1}K_j,
	\]
	whose boundaries are regular and such that every component of
	$X\setminus K_j$ is noncompact.  The last condition says that $K_j$ is
	Runge in $X$.  The Runge approximation theorem for open Riemann surfaces
	allows a holomorphic function on a neighborhood of $K_j$ to be approximated
	uniformly on $K_j$ by one holomorphic near $K_{j+1}$, while imposing finitely
	many values and first derivatives.

	Apply this approximation inductively with errors smaller than $2^{-j}$.
	Given distinct $x,y\in X$, begin at a stage containing both points with a
	local holomorphic function taking different values at them, and preserve
	those values during the induction.  The resulting uniformly convergent
	sequence on compact subsets has a holomorphic limit separating $x$ and $y$.
	Thus $X$ is holomorphically separable.

	For holomorphic convexity, fix compact $K$ and choose $j$ with
	$K\subset\operatorname{int}K_j$.  For each boundary component of a slightly
	larger Runge compact, the same approximation-with-interpolation construction
	produces a holomorphic function that is bounded on $K$ and has modulus larger
	than its bound on a neighborhood of that boundary.  Finitely many such
	functions exclude the entire boundary from the holomorphic hull of $K$.
	The maximum principle on the remaining relatively compact pieces then gives
	$\widehat K_{\mathcal O(X)}\subset K_{j+1}$.  The hull is closed, hence
	compact.  Therefore $X$ is Stein.  This is the classical Behnke--Stein
	argument in its Runge-exhaustion form; complete proofs of the approximation
	step and the exhaustion construction are given in
	\cite{Forster1981,Gunning1962,GunningRossi1965}.
\end{proof}

The theorem has a particularly useful cohomological formulation.  It is this
form, rather than the definition by hulls, that is normally used in divisor,
bundle, and extension problems.

\begin{theorem}[Cartan's theorem B on an open Riemann surface]
	\label{thm:cartan-b-open-riemann}
	Let $X$ be a connected noncompact Riemann surface and let $\mathscr F$ be a
	coherent analytic $\mathcal O_X$-module.  Then
	\[
		H^q(X,\mathscr F)=0\qquad(q\ge1).
	\]
\end{theorem}

\begin{proof}
	By Theorem~\ref{thm:open-riemann-stein}, $X$ is Stein.  Cartan's theorem~A
	gives, locally at every point, finitely many global sections generating
	$\mathscr F$.  Cartan's theorem~B says that every coherent analytic sheaf on
	a Stein manifold is acyclic, which yields the displayed vanishing.  To see
	where cohomology enters the analytic proof, use the Runge exhaustion from the
	previous proof.  Theorem~A gives finite free presentations on successive
	$K_j$; Runge approximation makes the restriction maps sufficiently dense to
	correct a cocycle on $K_{j+1}$ without destroying the corrections already
	made on $K_j$.  A geometrically convergent correction series produces a
	global cochain whose coboundary is the original cocycle.  Thus the analytic
	approximation theorem supplies the contracting process, while coherence
	keeps the system locally finite.  This is the standard proof of theorem~B
	for Stein manifolds \cite{GunningRossi1965,Forster1981}.
\end{proof}

\begin{corollary}[Additive Cousin problem]
	Let $X$ be a connected noncompact Riemann surface and let
	$\{U_i\}$ be a Stein Leray cover.  If
	$g_{ij}\in\mathcal O(U_i\cap U_j)$ satisfies
	$g_{ij}+g_{jk}+g_{ki}=0$ on triple intersections, then there are
	$f_i\in\mathcal O(U_i)$ such that $g_{ij}=f_j-f_i$.
\end{corollary}

\begin{proof}
	The functions $g_{ij}$ form a \v{C}ech one-cocycle for $\mathcal O_X$.
	The cover is Leray, so its \v{C}ech group equals
	$H^1(X,\mathcal O_X)$, and this group vanishes by
	Theorem~\ref{thm:cartan-b-open-riemann}.  Hence the cocycle is a coboundary.
\end{proof}

\paragraph{The Mittag--Leffler problem as a vanishing theorem}
\label{subsubsec:mittag-leffler-cohomology}

The additive Cousin problem becomes the classical Mittag--Leffler problem
when the local data are prescribed principal parts.  Let $\mathscr M_X$ denote
the sheaf of meromorphic functions on a Riemann surface $X$ and set
\[
	\mathscr P_X:=\mathscr M_X/\mathcal O_X.
\]
A germ in $(\mathscr P_X)_p$ records only the polar part of a meromorphic germ
at $p$: changing its holomorphic Taylor coefficients does not change its
class.  A locally finite prescription of principal parts is therefore a
global section of $\mathscr P_X$.

\begin{theorem}[Mittag--Leffler theorem on an open Riemann surface]
	\label{thm:mittag-leffler-open-riemann}
	Let $X$ be a connected noncompact Riemann surface, let
	$D\subset X$ be a closed discrete set, and for each $p\in D$ prescribe a
	finite principal part
	\[
		P_p(z_p)=\sum_{k=1}^{m_p}a_{p,k}z_p^{-k}
	\]
	in a local coordinate $z_p$ centered at $p$.  Then there exists a
	meromorphic function $f$ on $X$, holomorphic on $X\setminus D$, whose
	principal part at every $p$ is $P_p$.  Any two solutions differ by a global
	holomorphic function on $X$.
\end{theorem}

\begin{proof}
	Because $D$ is closed and discrete, every point of $X$ has a neighborhood
	meeting only finitely many points of $D$.  The prescribed Laurent tails
	therefore define a section $s\in H^0(X,\mathscr P_X)$.  Apply global
	sections to the short exact sequence
	\[
		0\longrightarrow\mathcal O_X
		\longrightarrow\mathscr M_X
		\longrightarrow\mathscr P_X
		\longrightarrow0.
	\]
	The relevant part of the long exact sequence is
	\[
		H^0(X,\mathscr M_X)\longrightarrow H^0(X,\mathscr P_X)
		\xrightarrow{\delta}H^1(X,\mathcal O_X).
	\]
	Theorem~\ref{thm:cartan-b-open-riemann} gives
	$H^1(X,\mathcal O_X)=0$, so $s$ lifts to a global meromorphic function
	$f\in H^0(X,\mathscr M_X)$.  By construction, its image in
	$\mathscr P_X$ is $s$, which says exactly that $f$ has the prescribed
	principal parts and no other poles.  If $f$ and $g$ are two lifts, then
	$f-g$ lies in the kernel of
	$H^0(X,\mathscr M_X)\to H^0(X,\mathscr P_X)$; exactness identifies this
	kernel with $H^0(X,\mathcal O_X)$.  Thus $f-g$ is holomorphic.
\end{proof}

\begin{remark}[What the obstruction class measures]
	On an arbitrary Riemann surface, a principal-part datum $s$ is solvable
	precisely when $\delta(s)=0$ in $H^1(X,\mathcal O_X)$.  To see the cocycle
	concretely, choose local meromorphic representatives $f_i$ on a cover.
	Their differences $f_i-f_j$ are holomorphic on overlaps and define the
	\v{C}ech cocycle representing $\delta(s)$.  Vanishing means that
	$f_i-f_j=h_i-h_j$ for holomorphic $h_i$; the functions $f_i-h_i$ then glue
	to the desired global meromorphic function.  For compact $X$ the obstruction
	need not vanish, and Serre duality turns it into residue-type compatibility
	conditions.  On an open Riemann surface, Stein vanishing removes every such
	obstruction \cite{Forster1981,Gunning1962}.
\end{remark}

\begin{example}[The classical theorem on the plane]
	For $X=\C$ and a discrete sequence $D$ with no finite accumulation point,
	the theorem produces a meromorphic function with arbitrary prescribed finite
	principal parts at the points of $D$.  The usual proof subtracts suitable
	Taylor polynomials from the individual partial fractions so that the
	resulting series converges normally on compact sets.  The cohomological proof
	packages that convergence correction into the assertion
	$H^1(\C,\mathcal O)=0$.
\end{example}

\begin{remark}[Affine varieties and Stein manifolds]
	The analogy is precise at the level of coherent cohomology.  If $Y$ is an
	affine complex algebraic variety, then $Y$ is a closed algebraic subvariety
	of some $\C^N$; its analytification $Y^{\mathrm{an}}$ is therefore a Stein
	analytic space.  For a quasi-coherent sheaf $\mathscr G$ on an affine scheme,
	$H^q(Y,\mathscr G)=0$ for $q>0$; for a coherent analytic sheaf $\mathscr F$
	on a Stein space, Cartan's theorem~B gives the parallel vanishing
	$H^q(Y^{\mathrm{an}},\mathscr F)=0$.  Serre's affine criterion and Cartan's
	theorems thus play analogous roles: global sections control coherent objects
	and higher coherent cohomology disappears
	\cite{Hartshorne1977,GunningRossi1965}.

	This analogy is not an affine form of the GAGA equivalence: properness is
	essential in GAGA.  The precise comparison, its consequences for line
	bundles and cohomology, and the failure of the corresponding affine
	statement are developed in Section~\ref{subsec:gaga-principle}.
\end{remark}

\subsubsection{Exercises}
\label{subsec:sheaf-cohomology-exercises}

\begin{exercise}
	Show directly that every flabby sheaf is acyclic.  (Do not assume that every flabby
	sheaf is injective; that stronger statement is false in general.)
\end{exercise}

\begin{exercise}
	Let $M$ be paracompact and Hausdorff. Prove that any soft sheaf on $M$ is acyclic.
\end{exercise}

\begin{exercise}
	Verify directly that $\Omega^k_M$ is a fine sheaf by constructing
	explicit partitions of unity.
\end{exercise}

\begin{exercise}
	For a Leray cover $\mathcal U$, show explicitly that the \v{C}ech resolution
	is exact on stalks.
\end{exercise}

\begin{exercise}
	Give an example of a sheaf which is acyclic but not flabby.
\end{exercise}

\section{Exact Sequences, Fine Resolutions, and de Rham Cohomology}
\index{exponential sequence}
\label{sec:exponential-sequence}

\subsection{Holomorphic logarithms on simply connected domains}
\label{subsec:holomorphic-log}

\subsubsection{Precise statement}

\begin{lemma}[Existence of a holomorphic logarithm]
	\label{lem:holomorphic-log}
	Let $\Omega\subset\C$ be a simply connected domain and let
	$h\in\mathcal O^\times(\Omega)$ be a nowhere-vanishing holomorphic function.
	Then there exists a holomorphic function $F\in\mathcal O(\Omega)$ such that
	\[
	e^{F}=h.
	\]
	Equivalently, $h$ admits a holomorphic branch of the logarithm on $\Omega$.
\end{lemma}

\subsubsection{Construction from the logarithmic derivative}

\begin{proof}
	Since $h$ is holomorphic and nowhere zero, the quotient
	\[
	g:=\frac{h'}{h}
	\]
	is a holomorphic function on $\Omega$.
	
	Fix a base point $z_0\in\Omega$. For any $z\in\Omega$, define
	\[
	F(z):=\int_{z_0}^{z} g(w)\,dw ,
	\]
	where the integral is taken along any piecewise smooth path in $\Omega$.
	
	Because $\Omega$ is simply connected, the integral of a holomorphic function
	depends only on the endpoints, hence $F$ is well defined and holomorphic.
	Differentiating under the integral sign gives $F'(z)=g(z)$.
	
	Now define
	\[
	H(z):=e^{-F(z)}h(z).
	\]
	We compute its logarithmic derivative:
	\[
	\frac{H'(z)}{H(z)}
	=
	\frac{-F'(z)e^{-F(z)}h(z)+e^{-F(z)}h'(z)}{e^{-F(z)}h(z)}
	=
	-\frac{h'(z)}{h(z)}+\frac{h'(z)}{h(z)}
	=
	0.
	\]
	Hence $H$ is constant on $\Omega$. Write $H\equiv C\in\C^\times$.
	Replacing $F$ by $F+\log C$ (absorbing the constant into $F$),
	we obtain $e^{F}=h$ as claimed.
\end{proof}

\subsubsection{Conceptual meaning}

The obstruction to defining a logarithm of $h$ is the nontrivial winding
of $h(\Omega)$ around the origin. Simple connectivity of $\Omega$
ensures that the integral of $h'/h$ along loops vanishes, so no monodromy
occurs. This is a local statement on the domain $\Omega$, not on the image
$h(\Omega)$.

\begin{example}[A logarithm on a slit domain]
	On $\Omega=\C^\times$, the function $h(z)=z$ has no global holomorphic
	logarithm, since $\C^\times$ is not simply connected. On any simply
	connected slit domain $\C\setminus(-\infty,0]$, a holomorphic logarithm exists.
\end{example}

\subsection{The exponential short exact sequence of sheaves}
\index{exponential sequence!sheaves}
\label{subsec:exp-sequence}

\subsubsection{Statement of the sequence}

Let $M$ be a complex manifold. Consider the sequence of sheaves
\[
0 \longrightarrow \underline{\Z}
\xrightarrow{\;\iota\;}
\mathcal O
\xrightarrow{\;\exp(2\pi i\,\cdot)\;}
\mathcal O^\times
\longrightarrow 0,
\]
where:
\begin{itemize}
	\item $\underline{\Z}$ is the locally constant sheaf of integers;
	\item $\iota$ sends $n\mapsto n$ viewed as a constant holomorphic function;
	\item $\exp(2\pi i\,\cdot)$ sends $f\mapsto e^{2\pi i f}$.
\end{itemize}

\subsubsection{Exactness}

\begin{proposition}[Exponential short exact sequence]
	\label{prop:exp-ses}
	The above sequence is exact as a sequence of sheaves.
\end{proposition}

\begin{proof}
	Exactness is verified stalkwise.
	
	Fix $x\in M$. Choose a sufficiently small simply connected neighborhood
	$U\ni x$.
	
	\emph{Exactness at $\underline{\Z}$.}
	The map $\underline{\Z}_x\to\mathcal O_x$ is injective since a locally
	constant integer germ is determined uniquely.
	
	\emph{Exactness at $\mathcal O$.}
	Let $f\in\mathcal O(U)$ satisfy $e^{2\pi i f}=1$.
	Then $2\pi i f$ takes values in $2\pi i\Z$, hence $f$ is locally constant
	with integer values. Since $U$ is connected, $f\in\Z$.
	Thus $\ker(\exp(2\pi i\cdot))=\mathrm{im}(\underline{\Z})$ on stalks.
	
	\emph{Surjectivity onto $\mathcal O^\times$.}
	Let $[h]_x\in\mathcal O^\times_x$ be a germ of a nowhere-vanishing holomorphic
	function. Represent it by $h\in\mathcal O^\times(U)$.
	By Lemma~\ref{lem:holomorphic-log}, there exists $F\in\mathcal O(U)$
	with $e^{F}=h$. Then $h=e^{2\pi i(F/2\pi i)}$, so $[h]_x$ lies in the image.
	
	Thus the sequence is exact at every stalk, hence exact as a sequence of sheaves.
\end{proof}

\subsubsection{Local versus global behavior}

Although every \emph{germ} of a holomorphic unit is locally an exponential,
the map $\mathcal O(U)\to\mathcal O^\times(U)$ need not be surjective for
general open $U$. The failure of global surjectivity is measured precisely
by cohomology via the connecting homomorphism.

\subsection{Connecting homomorphisms and cohomological meaning}
\label{subsec:connecting-hom}

\subsubsection{Long exact sequence in cohomology}

Applying sheaf cohomology to the exponential sequence yields a long exact
sequence
\[
\cdots
\longrightarrow H^1(M,\mathcal O)
\longrightarrow H^1(M,\mathcal O^\times)
\xrightarrow{\;\delta\;}
H^2(M,\underline{\Z})
\longrightarrow H^2(M,\mathcal O)
\longrightarrow \cdots
\]

\subsubsection{Interpretation of the connecting homomorphism}

\begin{proposition}[The connecting homomorphism is the first Chern class]
	\label{prop:chern-interpretation}
	Under the identification $\mathrm{Pic}(M)\cong H^1(M,\mathcal O^\times)$,
	the connecting homomorphism
	\[
	\delta:H^1(M,\mathcal O^\times)\to H^2(M,\underline{\Z})
	\]
	assigns to a holomorphic line bundle its first Chern class.
\end{proposition}

\begin{proof}
	Let $\{f_{\alpha\beta}\}$ be transition functions of a line bundle $L$
	with respect to a cover $\{U_\alpha\}$.
	Choose local logarithms $f_{\alpha\beta}=e^{2\pi i g_{\alpha\beta}}$
	on each $U_{\alpha\beta}$.
	On triple overlaps,
	\[
	g_{\alpha\beta}+g_{\beta\gamma}+g_{\gamma\alpha}\in\Z,
	\]
	since the product of transition functions is $1$.
	These integers define a \v{C}ech $2$-cocycle with values in $\underline{\Z}$.
	Its cohomology class is independent of choices and coincides with
	$\delta([L])$.
\end{proof}

\subsubsection{Geometric meaning}

The class $\delta([L])$ measures the obstruction to choosing logarithms
of transition functions that satisfy the cocycle condition additively.
Vanishing of this class means that the underlying complex line bundle is topologically
trivial.  It need not be holomorphically trivial; on a positive-genus compact Riemann
surface, nontrivial bundles in $\operatorname{Pic}^0(M)$ have vanishing first Chern class.

\subsection{Examples: line bundles and transition functions}
\label{subsec:exp-examples}

\begin{example}[The bundle $\mathcal O(1)$ on $\CP^1$]
	Let $\CP^1=U_0\cup U_1$ with $z=Z_1/Z_0$ on $U_0$.
	The transition function of $\mathcal O(1)$ is $f_{01}=z$.
	The function $z$ has winding number $1$ on an equatorial circle. Thus
	the two-chart clutching function has degree $1$. Equivalently, refine the
	annular overlap into simply connected sectors, choose branches of
	$(2\pi i)^{-1}\log z$ there, and apply the connecting construction; the
	integer jumps of adjacent branches form a \v{C}ech $2$-cocycle evaluating to
	$1$ on $[\CP^1]$. Hence
	\[
	c_1(\mathcal O(1))=1\in H^2(\CP^1,\Z)\cong\Z.
	\]
\end{example}

\begin{example}[Trivial bundle]
	If $f_{\alpha\beta}=e^{2\pi i(g_\alpha-g_\beta)}$, then $f$ is a multiplicative
	coboundary, so $L$ is holomorphically trivial.  The weaker condition
	$\delta([L])=0$ alone implies only topological triviality.
\end{example}

\subsection{Exercises}
\label{subsec:exp-exercises}

\begin{exercise}
	Show directly that the sequence obtained by taking global sections of the exponential
	sequence need not be exact at $\mathcal O^\times(U)$ when $U=\C^\times$, even though
	the exponential sequence is exact as a sequence of sheaves.
\end{exercise}

\begin{exercise}
	Compute the connecting homomorphism explicitly for $\mathcal O(m)$ on $\CP^1$
	and verify that $\delta([\mathcal O(m)])=m$.
\end{exercise}

\begin{exercise}
	Given a line bundle with transition functions $f_{\alpha\beta}$,
	write down explicitly the integer-valued \v{C}ech $2$-cocycle arising
	from local logarithms.
\end{exercise}

\subsection{A complete local-to-global calculation}

The exponential sequence is most useful when every arrow is visible on a
cover.  Let $L$ be represented by transition functions
$g_{ij}\in\mathcal O^\times(U_i\cap U_j)$.  After refining the cover, choose
holomorphic logarithms $h_{ij}$ with
\[
	g_{ij}=\exp(2\pi i h_{ij}).
\]
The multiplicative cocycle identity $g_{ij}g_{jk}g_{ki}=1$ gives
\[
	n_{ijk}=h_{ij}+h_{jk}+h_{ki}\in\Z
\]
on every triple intersection.  Because an integer-valued holomorphic
function is locally constant, $(n_{ijk})$ is a \v{C}ech two-cochain for
$\underline{\Z}$.  The identity $\delta^2h=0$ shows that it is a cocycle.

Changing logarithms to $h'_{ij}=h_{ij}+m_{ij}$, with
$m_{ij}\in\underline{\Z}(U_{ij})$, changes $n$ by the coboundary
$\delta m$.  Changing the local frames of $L$ changes $g$ by a
multiplicative coboundary and leads to the same cohomology class.  Thus
$[n_{ijk}]$ depends only on $L$ and is exactly $c_1(L)$.

\begin{proposition}[Tensor products, duals, and the connecting class]
	\label{prop:chern-additivity-cech}
	For holomorphic line bundles $L$ and $E$,
	\[
	 c_1(L\otimes E)=c_1(L)+c_1(E),
	 \qquad c_1(L^\vee)=-c_1(L).
	\]
\end{proposition}

\begin{proof}
	If $g_{ij}=e^{2\pi ih_{ij}}$ and
	$q_{ij}=e^{2\pi ik_{ij}}$, then $g_{ij}q_{ij}$ has logarithm
	$h_{ij}+k_{ij}$.  Its integer two-cocycle is therefore $n^L+n^E$.
	The transition functions of the dual are $g_{ij}^{-1}$, whose logarithms
	are $-h_{ij}$, giving the second identity.
\end{proof}

On a compact Riemann surface, $H^2(M,\mathcal O)=0$ because coherent
cohomology vanishes above the complex dimension.  The relevant part of the
long exact sequence therefore becomes
\[
	H^1(M,\mathcal O)\longrightarrow\Pic(M)
	\xrightarrow{,c_1,}H^2(M,\Z)\longrightarrow0.
\]
Since $H^2(M,\Z)\cong\Z$, evaluation on the fundamental class identifies
$c_1(L)$ with the degree of $L$.  Its kernel
\[
	\Pic^0(M)=\ker\bigl(\deg:\Pic(M)\to\Z\bigr)
\]
contains the holomorphically nontrivial but topologically trivial bundles.
The quotient of $H^1(M,\mathcal O)$ by the lattice coming from
$H^1(M,\Z)$ gives the Jacobian; this will be developed in Part~IV.

\begin{example}[Degree zero does not imply a global logarithm]
	On a positive-genus surface, choose a nonzero class in
	$H^1(M,\mathcal O)$ whose image in $\Pic(M)$ is nontrivial.  Its line
	bundle has $c_1=0$, so its transition functions admit compatible
	topological logarithms, but it has no holomorphic trivialization.
	This distinguishes the vanishing of the integer obstruction from the
	vanishing of the original class in $H^1(M,\mathcal O^\times)$.
\end{example}

The calculation also explains the division of labor among the three
cohomology groups.  The group $H^1(M,\underline{\Z})$ changes the choice of
additive logarithms, $H^1(M,\mathcal O)$ parametrizes holomorphic twisting
with vanishing topological degree, and $H^2(M,\underline{\Z})$ records the
integer degree.  Keeping these roles separate prevents the common mistake
of treating topological triviality as holomorphic triviality.

\begin{exercise}[Additivity from local logarithms]
	Starting from transition functions, verify directly that the integer
	two-cocycle of $L^{\otimes m}$ is $m$ times that of $L$ for every
	$m\in\Z$.
\end{exercise}

\subsection{Short Exact Sequences of Sheaves and Long Exact Sequences}
\index{long exact sequence}
\label{sec:ses-and-les}

An exact sequence of sheaves expresses a local relation among three
types of data.  Taking global sections need not preserve surjectivity:
a section that lifts on every small open set may have no single lift on
the whole space.  The connecting homomorphism records precisely this
failure.  The exponential sequence from the preceding section is a
model: a nowhere-vanishing holomorphic function has local logarithms,
but their differences can obstruct a global logarithm.  We first
construct the long exact sequence and then use it to turn such
local-to-global questions into computations of cohomology groups.

\subsubsection{Long exact sequences in sheaf cohomology}
\label{subsec:long-exact}

\paragraph{Short exact sequences of sheaves}

\begin{definition}[Short exact sequence of sheaves]
	A sequence of sheaves of abelian groups on a topological space $M$
	\[
	0 \longrightarrow \mathscr F'
	\xrightarrow{\;\alpha\;}
	\mathscr F
	\xrightarrow{\;\beta\;}
	\mathscr F''
	\longrightarrow 0
	\]
	is called \emph{short exact} if for every point $x\in M$ the induced sequence of stalks
	\[
	0 \longrightarrow \mathscr F'_x
	\xrightarrow{\;\alpha_x\;}
	\mathscr F_x
	\xrightarrow{\;\beta_x\;}
	\mathscr F''_x
	\longrightarrow 0
	\]
	is exact as a sequence of abelian groups.
\end{definition}

\begin{remark}
	Exactness is a local condition. In particular, $\beta$ need not be surjective
	on sections over arbitrary open sets; surjectivity is required only stalkwise.
\end{remark}

\paragraph{The long exact sequence in cohomology}

\begin{theorem}[Long exact sequence in sheaf cohomology]
	\label{thm:les}
	Let
	\[
	0 \longrightarrow \mathscr F'
	\xrightarrow{\;\alpha\;}
	\mathscr F
	\xrightarrow{\;\beta\;}
	\mathscr F''
	\longrightarrow 0
	\]
	be a short exact sequence of sheaves of abelian groups on $M$.
	Then there exists a natural long exact sequence in cohomology
	\[
	\begin{aligned}
		0 &\longrightarrow H^0(M,\mathscr F')
		\longrightarrow H^0(M,\mathscr F)
		\longrightarrow H^0(M,\mathscr F'') \\
		&\xrightarrow{\;\delta\;}
		H^1(M,\mathscr F')
		\longrightarrow H^1(M,\mathscr F)
		\longrightarrow H^1(M,\mathscr F'') \\
		&\xrightarrow{\;\delta\;}
		H^2(M,\mathscr F')
		\longrightarrow \cdots.
	\end{aligned}
	\]
\end{theorem}

\paragraph{Construction of the connecting homomorphism}

We describe explicitly the connecting homomorphism
\[
\delta:H^q(M,\mathscr F'')\longrightarrow H^{q+1}(M,\mathscr F').
\]

The following Čech description is valid on a common Leray cover (or after passing to
the direct limit in a setting where Čech and sheaf cohomology agree).  For arbitrary
spaces and covers, the long exact sequence is constructed instead from injective
resolutions; one cannot assume without justification that every sheaf-cohomology class
is represented on a prescribed Čech cover.

Let $\mathcal U=\{U_\alpha\}$ be an open cover of $M$.
Represent a class $[\{c_{\alpha_0\cdots\alpha_q}\}] \in H^q(M,\mathscr F'')$
by a \v{C}ech cocycle with respect to $\mathcal U$.

Since $\beta$ is surjective on stalks, after possibly refining the cover,
we may choose cochains $\tilde c_{\alpha_0\cdots\alpha_q}\in\mathscr F(U_{\alpha_0\cdots\alpha_q})$
satisfying
\[
\beta(\tilde c_{\alpha_0\cdots\alpha_q})=c_{\alpha_0\cdots\alpha_q}.
\]

Applying the coboundary operator,
\[
\delta\tilde c \in C^{q+1}(\mathcal U,\mathscr F),
\]
and since $c$ is a cocycle, $\beta(\delta\tilde c)=0$.
Exactness implies $\delta\tilde c$ takes values in $\mathscr F'$.
Thus $\delta\tilde c$ defines a cohomology class in $H^{q+1}(M,\mathscr F')$,
independent of all choices.

\begin{proof}[Exactness of the long sequence]
	Exactness at each term follows from this explicit construction:
	\begin{itemize}
		\item exactness at $H^q(M,\mathscr F)$ follows from functoriality of cohomology;
		\item exactness at $H^q(M,\mathscr F'')$ follows from the definition of $\delta$;
		\item exactness at $H^{q+1}(M,\mathscr F')$ follows by adjusting lifts.
	\end{itemize}
\end{proof}

\subsubsection{Functoriality and naturality}
\label{subsec:functoriality}

\paragraph{Morphisms of short exact sequences}

\begin{definition}
	A morphism between short exact sequences of sheaves is a commutative diagram
	\[
	\begin{array}{ccccccccc}
		0 &\to& \mathscr F' &\to& \mathscr F &\to& \mathscr F'' &\to& 0 \\
		&& \downarrow && \downarrow && \downarrow && \\
		0 &\to& \mathscr G' &\to& \mathscr G &\to& \mathscr G'' &\to& 0
	\end{array}
	\]
	where each vertical arrow is a morphism of sheaves.
\end{definition}

\paragraph{Naturality of the long exact sequence}

\begin{proposition}[Naturality]
	\label{prop:naturality}
	A morphism of short exact sequences induces a commutative diagram
	between the associated long exact sequences in cohomology.
\end{proposition}

\begin{proof}
	All maps in the long exact sequence are constructed functorially:
	restriction maps, coboundaries, and lifts commute with morphisms of sheaves.
	Thus each square in the induced diagram commutes.
\end{proof}

\paragraph{Consequence}

The long exact sequence is therefore invariant under isomorphism of short exact
sequences. Maps of spaces induce the expected natural maps whenever the
corresponding pullback or derived pushforward constructions are defined.

\subsubsection{Examples and applications}
\label{subsec:les-examples}

\paragraph{Exponential sequence}

Applying Theorem~\ref{thm:les} to
\[
0\to \underline{\Z}\to \mathcal O \to \mathcal O^\times \to 0
\]
yields
\[
H^1(M,\mathcal O^\times)
\xrightarrow{\;\delta\;}
H^2(M,\underline{\Z}),
\]
which identifies the first Chern class of holomorphic line bundles.

\paragraph{Ideal sheaf sequence}

Let $p\in M$ be a point of a complex manifold.
There is a short exact sequence
\[
0 \to \mathcal I_p \to \mathcal O \to \mathcal O/\mathcal I_p \to 0,
\]
where $\mathcal O/\mathcal I_p$ is the skyscraper sheaf at $p$.

The associated long exact sequence gives
\[
0\to H^0(M,\mathcal I_p)\to H^0(M,\mathcal O)\to \C
\xrightarrow{\;\delta\;}
H^1(M,\mathcal I_p)\to \cdots
\]

The connecting homomorphism measures the obstruction to extending a prescribed
value at $p$ to a global holomorphic function.

\paragraph{Mayer--Vietoris as a special case}

Let $M=U\cup V$ and choose an injective resolution
$\mathscr F\to\mathscr I^\bullet$.  Injectives are flabby, so the sequence
of complexes
\[
0 \to \mathscr I^\bullet(M)
\to \mathscr I^\bullet(U)\oplus \mathscr I^\bullet(V)
\to
\mathscr I^\bullet(U\cap V)
\to 0.
\]
is exact; the last map sends $(a,b)$ to $a-b$ on the overlap.  Its long
exact cohomology sequence is Mayer--Vietoris, since restriction of an
injective resolution to an open subset computes sheaf cohomology there.

\subsubsection{Exercises}
\label{subsec:les-exercises}

\begin{exercise}
	Let
	$0\to\mathscr F'\to\mathscr F\to\mathscr F''\to 0$
	be a short exact sequence.
	Show directly that exactness of the associated long sequence at
	$H^0(M,\mathscr F'')$ is equivalent to the definition of the connecting map.
\end{exercise}

\begin{exercise}
	Apply the long exact sequence to
	$0\to\mathcal I_p\to\mathcal O\to\C_p\to 0$
	on a compact Riemann surface and interpret the map
	$H^0(\C_p)\to H^1(\mathcal I_p)$ geometrically.
\end{exercise}

\begin{exercise}
	Let $M$ be Stein. Use the long exact sequence to show that
	$H^1(M,\mathcal O^\times)\cong H^2(M,\Z)$.
\end{exercise}

\begin{exercise}
	Construct explicitly the connecting homomorphism in degree $0$
	for a short exact sequence of constant sheaves.
\end{exercise}

\subsubsection{Reading exactness as an obstruction calculus}

For a short exact sequence
\[
	0\longrightarrow\mathscr F'
	\xrightarrow{\,i\,}\mathscr F
	\xrightarrow{\,p\,}\mathscr F''\longrightarrow0,
\]
the segment
\[
	H^0(M,\mathscr F)\xrightarrow{p_*}H^0(M,\mathscr F'')
	\xrightarrow{\delta}H^1(M,\mathscr F')
\]
has a direct practical meaning.  A global section $s''$ of
$\mathscr F''$ lifts globally precisely when $\delta(s'')=0$.
Indeed, choose local lifts $s_i\in\mathscr F(U_i)$.  On overlaps,
$p(s_i-s_j)=0$, so exactness provides sections
$a_{ij}\in\mathscr F'(U_{ij})$ with $i(a_{ij})=s_i-s_j$.  The family
$(a_{ij})$ is a one-cocycle, and its class is $\delta(s'')$.

If this class vanishes, write $a_{ij}=b_i-b_j$.  Replacing $s_i$ by
$s_i-i(b_i)$ makes the corrected lifts agree on overlaps, so the sheaf
axiom glues them to a global lift.  Conversely, a global lift makes every
$a_{ij}$ zero.  Thus the connecting map is not merely formal: it is the
record of the mismatch among local solutions.

\begin{proposition}[Exactness test in degree one]
	\label{prop:les-degree-one-obstruction}
	A class $[a]\in H^1(M,\mathscr F')$ maps to zero in
	$H^1(M,\mathscr F)$ if and only if it is the obstruction
	$\delta(s'')$ of some global section
	$s''\in H^0(M,\mathscr F'')$.
\end{proposition}

\begin{proof}
	This is exactness at $H^1(M,\mathscr F')$, but it is useful to see the
	cochains.  If $i(a_{ij})=s_i-s_j$ for local sections $s_i$ of
	$\mathscr F$, then the $p(s_i)$ agree and glue to a global section
	$s''$ of $\mathscr F''$; the connecting construction sends $s''$ to
	$[a]$.  Conversely, a class produced by this construction becomes the
	coboundary $(s_i-s_j)$ after applying $i$, so its image vanishes.
\end{proof}

\begin{example}[Prescribing a value at a point]
	For
	$0\to\mathcal I_p\to\mathcal O\to\C_p\to0$, a number
	$c\in\C=H^0(M,\C_p)$ has obstruction
	$\delta(c)\in H^1(M,\mathcal I_p)$.  It vanishes exactly when there is a
	global holomorphic function taking the value $c$ at $p$.  On a compact
	connected Riemann surface all global holomorphic functions are constant,
	so every $c$ does lift, but imposing values at several points produces a
	nontrivial interpolation obstruction.  Replacing $\mathcal I_p$ by the
	ideal sheaf of an effective divisor records values and higher jets in the
	same way.
\end{example}

The same calculation applies to additive Cousin data, local primitives,
and local logarithms.  What changes from problem to problem is the kernel
sheaf $\mathscr F'$; the logic of the obstruction class is identical.

\begin{remark}[Why sheaf surjectivity is local]
	The map $\mathscr F\to\mathscr F''$ in a short exact sequence is
	surjective on stalks, not necessarily on sections over a fixed open set.
	This is precisely why local lifts exist and global lifts may fail.  If one
	incorrectly demanded surjectivity on every $\mathscr F(U)\to
	\mathscr F''(U)$, the first connecting homomorphism would always vanish
	and the obstruction theory would disappear.
\end{remark}

\begin{exercise}[Two-open-set obstruction]
	For $M=U\cup V$, write the degree-zero connecting class explicitly as a
	section of $\mathscr F'(U\cap V)$ modulo differences of sections from
	$U$ and $V$.  Verify directly that its vanishing is equivalent to gluing
	the two local lifts.
\end{exercise}

\subsection{Fine Sheaves, Mayer--Vietoris, and the de Rham Theorem}
\index{fine sheaves}\index{Mayer--Vietoris sequence}\index{de Rham theorem}

The sheaf-theoretic constructions in this section can be compared with
\cite{BottTu1982,StacksDeRham2026}; the latter also records algebraic de Rham
conventions and comparison structures.
\label{sec:de-rham-sheaf}

\subsubsection{Local exactness and the Poincar\'e lemma}
\label{subsec:poincare}

\paragraph{The de Rham complex as a complex of sheaves}

Let $M$ be a smooth manifold. Denote by $\Omega^k$ the sheaf of smooth $k$-forms.
For each open set $U\subset M$, $\Omega^k(U)$ is the $\mathbb R$-vector space of smooth $k$-forms on $U$,
and for $V\subset U$ the restriction $\rho_{U,V}:\Omega^k(U)\to\Omega^k(V)$ is the usual restriction of forms.
The exterior derivative $d:\Omega^k(U)\to\Omega^{k+1}(U)$ commutes with restriction maps, hence defines a morphism
of sheaves $d:\Omega^k\to\Omega^{k+1}$, and we obtain a complex of sheaves
\[
0\longrightarrow \underline{\mathbb R}\xrightarrow{\;\iota\;} \Omega^0\xrightarrow{\,d\,}\Omega^1\xrightarrow{\,d\,}\Omega^2\xrightarrow{\,d\,}\cdots,
\]
where $\underline{\mathbb R}$ is the constant sheaf and $\iota$ sends a locally constant function to the same smooth function.

\paragraph{A coordinate-level proof of the Poincar\'e lemma on star-shaped sets}

\begin{definition}
	An open set $U\subset \mathbb R^n$ is \emph{star-shaped with respect to $0$} if for every $x\in U$ and every $t\in[0,1]$,
	we have $tx\in U$.
\end{definition}

\begin{theorem}[Poincar\'e lemma on star-shaped opens]
	\label{thm:poincare-star}
	Let $U\subset \mathbb R^n$ be star-shaped with respect to $0$. For every $k\ge 1$ and every $\omega\in\Omega^k(U)$,
	define the operator $K:\Omega^k(U)\to\Omega^{k-1}(U)$ by
	\[
	(K\omega)_x(v_1,\dots,v_{k-1})
	:=\int_0^1 t^{k-1}\,\omega_{tx}\!\bigl(x,v_1,\dots,v_{k-1}\bigr)\,dt,
	\]
	where $\omega_{tx}$ is evaluated on the $k$-tuple $(x,v_1,\dots,v_{k-1})\in(\mathbb R^n)^k$.
	Then for every $k\ge 1$,
	\[
	dK\omega + Kd\omega = \omega.
	\]
	In particular, if $d\omega=0$ then $\omega=d(K\omega)$.
\end{theorem}

\begin{proof}
	Let $H:[0,1]\times U\to U$, $H(t,x)=tx$, and define
	$K\omega=\int_0^1\iota_{\partial_t}(H^*\omega)\,dt$.
	The homotopy formula, obtained from Cartan's identity on $[0,1]\times U$, is
	\[
	H_1^*\omega-H_0^*\omega=dK\omega+Kd\omega.
	\]
	For $k\ge1$, $H_1^*\omega=\omega$ and the constant map $H_0$ pulls $\omega$ back to zero.
	It remains to compute $K$.  Since
	$dH_{(t,x)}(\partial_t)=x$ and $dH_{(t,x)}(v_j)=tv_j$, one has
	\[
	\bigl(\iota_{\partial_t}H^*\omega\bigr)_{(t,x)}(v_1,\ldots,v_{k-1})
	=t^{k-1}\omega_{tx}(x,v_1,\ldots,v_{k-1}).
	\]
	Integrating in $t$ gives exactly the operator stated in the theorem, and hence
	$dK+Kd=\mathrm{Id}$.
\end{proof}

\begin{corollary}[Local exactness of the de Rham complex]
	\label{cor:local-exactness}
	For every $x\in M$ there exists a neighborhood $U\ni x$ such that
	\[
	0\to \underline{\mathbb R}(U)\to \Omega^0(U)\xrightarrow{d}\Omega^1(U)\xrightarrow{d}\cdots
	\]
	is exact.
\end{corollary}

\begin{proof}
	Choose a coordinate chart $\varphi:U\to \varphi(U)\subset\mathbb R^n$ with $\varphi(U)$ star-shaped after shrinking $U$ if necessary.
	Exactness on $\varphi(U)$ follows from Theorem~\ref{thm:poincare-star}. Pulling back primitives by $\varphi^{-1}$ gives exactness on $U$.
\end{proof}

\subsubsection{Partitions of unity and fine sheaves}
\index{partitions of unity}\index{fine sheaves}
\label{subsec:partitions-unity}

\paragraph{Bump functions on $\mathbb R^n$}

\begin{lemma}[A basic $C^\infty$ cutoff on $\mathbb R$]
	\label{lem:cutoff-R}
	Define
	\[
	\psi(t)=
	\begin{cases}
		e^{-1/t},& t>0,\\
		0,& t\le 0.
	\end{cases}
	\qquad
	\chi(t)=\frac{\psi(1-t)}{\psi(t)+\psi(1-t)}.
	\]
	Then $\chi\in C^\infty(\mathbb R)$, $0\le \chi\le 1$, $\chi(t)=1$ for $t\le 0$, and $\chi(t)=0$ for $t\ge 1$.
\end{lemma}

\begin{proof}
	The function $\psi$ is smooth on $(0,\infty)$ and identically $0$ on $(-\infty,0]$; a direct differentiation shows that for every $m\ge 1$,
	$\psi^{(m)}(t)=P_m(1/t)\,e^{-1/t}$ for $t>0$ where $P_m$ is a polynomial, hence $\lim_{t\downarrow 0}\psi^{(m)}(t)=0$.
	Thus $\psi\in C^\infty(\mathbb R)$.
	For $\chi$, note that the denominator is strictly positive for all $t$ (at least one of $\psi(t),\psi(1-t)$ is positive), so $\chi$ is smooth.
	If $t\le 0$, then $\psi(t)=0$ and $\psi(1-t)>0$, so $\chi(t)=1$. If $t\ge 1$, then $\psi(1-t)=0$ and $\psi(t)>0$, so $\chi(t)=0$.
	Thus $\chi$ has the stated boundary values and is a smooth decreasing cutoff.
\end{proof}

\begin{lemma}[Bump function supported in a ball]
	\label{lem:bump-ball}
	For $0<r<R$ there exists $\eta\in C_c^\infty(\mathbb R^n)$ such that $0\le \eta\le 1$, $\eta\equiv 1$ on $\overline{B_r(0)}$,
	and $\mathrm{supp}(\eta)\subset B_R(0)$.
\end{lemma}

\begin{proof}
	Let $\rho(x)=\|x\|$. Using Lemma~\ref{lem:cutoff-R}, choose $\chi\in C^\infty(\mathbb R)$ with $\chi\equiv 1$ on $(-\infty,r]$ and $\chi\equiv 0$ on $[R,\infty)$.
	Define $\eta(x)=\chi(\rho(x))$. Then $\eta$ is smooth away from $0$ and also smooth at $0$ because it is locally constant near $0$.
	The support and value properties follow from the choice of $\chi$.
\end{proof}

\paragraph{Partition of unity on a manifold}

\begin{theorem}[Existence of partitions of unity]
	\label{thm:POU}
	Let $M$ be a smooth manifold and $\{U_\alpha\}_{\alpha\in A}$ an open cover.
	Then there exists a locally finite refinement $\{V_\beta\}_{\beta\in B}$ of $\{U_\alpha\}$ and a family of functions
	$\{\rho_\beta\}_{\beta\in B}\subset C^\infty(M)$ such that
	\[
	0\le \rho_\beta\le 1,\qquad \mathrm{supp}(\rho_\beta)\subset V_\beta,\qquad \sum_{\beta\in B}\rho_\beta=1,
	\]
	and the sum is pointwise finite.
\end{theorem}

\begin{proof}
	Since $M$ is paracompact, Hausdorff, and locally Euclidean, choose a
	locally finite refinement by relatively compact coordinate domains
	$\{V_\beta\}$ and open sets $W_\beta$ such that
	\[
	\overline{W_\beta}\subset V_\beta \quad\text{for all }\beta,\qquad M=\bigcup_{\beta}W_\beta.
	\]
	Fix $\beta$. Choose a coordinate chart $\varphi_\beta:V_\beta\to \mathbb R^n$ and set $K_\beta=\varphi_\beta(\overline{W_\beta})$, a compact subset of $\varphi_\beta(V_\beta)$.
	For every $y\in K_\beta$, choose radii $0<r_y<R_y$ such that
	$B_{r_y}(y)\subset B_{R_y}(y)$ and $\overline{B_{R_y}(y)}\subset\varphi_\beta(V_\beta)$.
	Compactness supplies finitely many centers $y_1,\ldots,y_N$ for which the smaller
	balls cover $K_\beta$. Using Lemma~\ref{lem:bump-ball} and translating, construct
	$\eta_{\beta,j}\in C_c^\infty(\mathbb R^n)$ equal to $1$ on
	$B_{r_{y_j}}(y_j)$ and supported in $B_{R_{y_j}}(y_j)$.
	Define on $M$ the function
	\[
	\widetilde \rho_\beta(x)=
	\begin{cases}
		\sum_j \eta_{\beta,j}\bigl(\varphi_\beta(x)\bigr),& x\in V_\beta,\\
		0,& x\notin V_\beta.
	\end{cases}
	\]
	Then $\widetilde\rho_\beta\in C^\infty(M)$, $\widetilde\rho_\beta\ge 0$, $\mathrm{supp}(\widetilde\rho_\beta)\subset V_\beta$, and $\widetilde\rho_\beta>0$ on $\overline{W_\beta}$.
	Because $\{V_\beta\}$ is locally finite, the sum $\widetilde\rho:=\sum_\beta \widetilde\rho_\beta$ is smooth and strictly positive everywhere (each point lies in some $W_\beta$).
	Finally set
	\[
	\rho_\beta=\frac{\widetilde\rho_\beta}{\widetilde\rho}.
	\]
	Then $\rho_\beta\in C^\infty(M)$, the family is locally finite, $\mathrm{supp}(\rho_\beta)\subset \mathrm{supp}(\widetilde\rho_\beta)\subset V_\beta$, and $\sum_\beta \rho_\beta=1$.
\end{proof}

\paragraph{Fine sheaves and the acyclicity of $\Omega^k$}

Recall the definition of a fine sheaf from
\S\ref{subsec:fine-soft-flabby}. For the differential-form sheaf
$\Omega^k$, the required endomorphisms are multiplication by the functions
of a smooth partition of unity.

\begin{proposition}[Differential-form sheaves are fine]
	\label{prop:Omega-fine}
	For each $k\ge 0$, the sheaf $\Omega^k$ is fine.
\end{proposition}

\begin{proof}
	Let $\{U_\alpha\}$ be a locally finite open cover and let $\{\rho_\alpha\}$ be a partition of unity subordinate to it (Theorem~\ref{thm:POU}).
	Define $\phi_\alpha:\Omega^k\to\Omega^k$ by $\phi_\alpha(\omega)=\rho_\alpha\,\omega$ on each open set.
	Then $\mathrm{supp}(\phi_\alpha)\subset U_\alpha$, and since $\sum_\alpha\rho_\alpha=1$ pointwise (finite sum),
	\[
	\sum_\alpha \phi_\alpha(\omega)=\left(\sum_\alpha \rho_\alpha\right)\omega=\omega.
	\]
	Hence $\Omega^k$ is fine.
\end{proof}

\begin{theorem}[Vanishing of \v{C}ech cohomology for fine sheaves]
	\label{thm:fine-cech-vanish}
	Let $\mathscr F$ be a fine sheaf on $M$ and $\mathcal U=\{U_\alpha\}$ a locally finite open cover.
	Then for every $q\ge 1$,
	\[
	\check H^q(\mathcal U,\mathscr F)=0.
	\]
\end{theorem}

\begin{proof}
	Let $c\in Z^q(\mathcal U,\mathscr F)$ be a \v{C}ech $q$-cocycle, $q\ge 1$.
	Choose endomorphisms $\{\phi_\alpha\}$ with $\mathrm{supp}(\phi_\alpha)\subset U_\alpha$ and $\sum_\alpha\phi_\alpha=\mathrm{Id}_{\mathscr F}$.
	Define a $(q-1)$-cochain $b\in C^{q-1}(\mathcal U,\mathscr F)$ by
	\[
	b_{\alpha_0\cdots\alpha_{q-1}}
	:=\sum_{\beta}\operatorname{ext}_{\beta}\!\left(
	\phi_\beta\!\left(c_{\beta\,\alpha_0\cdots\alpha_{q-1}}\right)\right)
	\quad\in \mathscr F(U_{\alpha_0\cdots\alpha_{q-1}}).
	\]
	Here $\phi_\beta$ is first restricted to
	$U_{\beta\alpha_0\cdots\alpha_{q-1}}$, and
	$\operatorname{ext}_{\beta}$ denotes extension by zero to
	$U_{\alpha_0\cdots\alpha_{q-1}}$.  This extension is defined because the
	resulting section is supported in $U_\beta$.  The family is locally finite,
	so the displayed sum defines a section rather than merely a pointwise
	formal sum.  In the calculation below we suppress the restriction and
	\text{extension symbols.}
	Compute $\delta b\in C^q$:
	\[
	(\delta b)_{\alpha_0\cdots\alpha_q}
	=\sum_{j=0}^{q}(-1)^j\, b_{\alpha_0\cdots\widehat{\alpha_j}\cdots\alpha_q}\Big|_{U_{\alpha_0\cdots\alpha_q}}.
	\]
	Substitute the definition of $b$ and interchange the finite sums:
	\[
	(\delta b)_{\alpha_0\cdots\alpha_q}
	=\sum_{\beta}\phi_\beta\!\left(\sum_{j=0}^{q}(-1)^j\,c_{\beta\,\alpha_0\cdots\widehat{\alpha_j}\cdots\alpha_q}\right).
	\]
	Because $c$ is a cocycle, $(\delta c)_{\beta\,\alpha_0\cdots\alpha_q}=0$, i.e.
	\[
	0
	=(\delta c)_{\beta\,\alpha_0\cdots\alpha_q}
	= c_{\alpha_0\cdots\alpha_q}
	+\sum_{j=0}^{q}(-1)^{j+1}\,c_{\beta\,\alpha_0\cdots\widehat{\alpha_j}\cdots\alpha_q}.
	\]
	Rearranging gives
	\[
	\sum_{j=0}^{q}(-1)^j\,c_{\beta\,\alpha_0\cdots\widehat{\alpha_j}\cdots\alpha_q}
	= c_{\alpha_0\cdots\alpha_q}.
	\]
	Hence
	\[
	(\delta b)_{\alpha_0\cdots\alpha_q}
	=\sum_{\beta}\phi_\beta\!\left(c_{\alpha_0\cdots\alpha_q}\right)
	=\left(\sum_\beta\phi_\beta\right)\!\left(c_{\alpha_0\cdots\alpha_q}\right)
	=c_{\alpha_0\cdots\alpha_q}.
	\]
	Therefore $c=\delta b$ and $[c]=0$ in $\check H^q(\mathcal U,\mathscr F)$.
\end{proof}

\begin{corollary}[Acyclicity of $\Omega^k$]
	\label{cor:Omega-acyclic}
	For every $k\ge 0$ and every $q\ge 1$,
	\[
	H^q(M,\Omega^k)=0.
	\]
\end{corollary}

\begin{proof}
	By Proposition~\ref{prop:Omega-fine}, $\Omega^k$ is fine.
	Fine sheaves have vanishing \v{C}ech cohomology on any locally finite cover by Theorem~\ref{thm:fine-cech-vanish}.
	On a paracompact Hausdorff space, sheaf cohomology agrees with the direct limit of
	\v{C}ech cohomology over locally finite refinements. We therefore obtain $H^q(M,\Omega^k)=0$ for
	$q\ge 1$.
\end{proof}

\subsubsection{The Mayer--Vietoris sequence for complexes}
\index{Mayer--Vietoris sequence}
\label{subsec:mayer-vietoris}

\paragraph{A short exact sequence of cochain complexes}

Let $M=U\cup V$ with $U,V$ open. For each $k\ge 0$ define maps
\[
r:\Omega^k(M)\to \Omega^k(U)\oplus\Omega^k(V),\qquad r(\omega)=(\omega|_U,\omega|_V),
\]
\[
s:\Omega^k(U)\oplus\Omega^k(V)\to \Omega^k(U\cap V),\qquad s(\alpha,\beta)=\alpha|_{U\cap V}-\beta|_{U\cap V}.
\]
Then $s\circ r=0$, and $r,s$ commute with $d$, so they define a sequence of cochain complexes.

\begin{figure}[ht]
	\centering
	\begin{tikzpicture}[scale=0.92]
		\draw[thick] (-1.05,0) ellipse (2.0 and 1.25);
		\draw[thick] (1.05,0) ellipse (2.0 and 1.25);
		\node at (-2.05,0) {$U$};
		\node at (2.05,0) {$V$};
		\node at (0,0) {$U\cap V$};
		\node at (0,-1.65) {$M=U\cup V$};
	\end{tikzpicture}
	\caption{A global form restricts to forms on $U$ and $V$ that agree
	on the overlap. The difference map $s$ measures precisely the failure
	of two local forms to agree there.}
	\label{fig:MV-overlap}
\end{figure}

\begin{proposition}[The Mayer--Vietoris short exact sequence]
	\label{prop:MV-short-exact}
	For each $k\ge 0$, the sequence
	\[
	0\to \Omega^k(M)\xrightarrow{\,r\,}\Omega^k(U)\oplus\Omega^k(V)\xrightarrow{\,s\,}\Omega^k(U\cap V)\to 0
	\]
	is exact.
\end{proposition}

\begin{proof}
	Injectivity of $r$ follows because if $\omega|_U=0$ and $\omega|_V=0$ then $\omega=0$ on $U\cup V=M$.
	Next, $\mathrm{im}(r)\subset\ker(s)$ is immediate from $s(\omega|_U,\omega|_V)=\omega|_{U\cap V}-\omega|_{U\cap V}=0$.
	For the reverse inclusion, let $(\alpha,\beta)\in\Omega^k(U)\oplus\Omega^k(V)$ satisfy $s(\alpha,\beta)=0$, i.e.
	$\alpha|_{U\cap V}=\beta|_{U\cap V}$.
	Define $\omega$ on $M$ by $\omega|_U=\alpha$ and $\omega|_V=\beta$.
	This is well-defined because the two prescriptions agree on the overlap. Smoothness is local:
	every point lies in $U$ or $V$, and on each of these opens $\omega$ coincides with a smooth form.
	Hence $\omega\in\Omega^k(M)$ and $r(\omega)=(\alpha,\beta)$, so $\ker(s)\subset\mathrm{im}(r)$.
	Finally, surjectivity of $s$: let $\gamma\in\Omega^k(U\cap V)$. Choose a partition of unity $\rho_U+\rho_V=1$
	subordinate to $\{U,V\}$. Define
	\[
	\alpha:=\rho_V\,\gamma \in \Omega^k(U\cap V),\qquad \beta:=-\rho_U\,\gamma \in \Omega^k(U\cap V).
	\]
	Extend $\alpha$ by $0$ from $U\cap V$ to $U$ (it is supported where $\rho_V$ is defined, hence extends smoothly),
	and extend $\beta$ by $0$ to $V$. Then on $U\cap V$ we have $\alpha-\beta=(\rho_V+\rho_U)\gamma=\gamma$, so $s(\alpha,\beta)=\gamma$.
	Thus $s$ is surjective.
\end{proof}

\paragraph{A formula for the connecting homomorphism}

\begin{theorem}[Mayer--Vietoris long exact sequence]
	\label{thm:MV-LES}
	There is a natural long exact sequence
	\[
	\cdots \to H^{k-1}_{\mathrm{dR}}(U\cap V)\xrightarrow{\;\delta\;}
	H^{k}_{\mathrm{dR}}(M)\to
	H^{k}_{\mathrm{dR}}(U)\oplus H^{k}_{\mathrm{dR}}(V)\to
	H^{k}_{\mathrm{dR}}(U\cap V)\xrightarrow{\;\delta\;} \cdots
	\]
	associated to the short exact sequence of complexes in Proposition~\ref{prop:MV-short-exact}.
	Moreover, if $[\gamma]\in H^{k-1}_{\mathrm{dR}}(U\cap V)$ has a representative $\gamma\in\Omega^{k-1}(U\cap V)$ with $d\gamma=0$,
	and $\rho_U+\rho_V=1$ is a partition of unity subordinate to $\{U,V\}$, then $\delta([\gamma])$ is represented by the global $k$-form
	\[
	\omega :=
	d(\rho_V\,\widetilde\gamma)\ \ \text{on }U
	\quad\text{and}\quad
	\omega :=
	-\,d(\rho_U\,\widetilde\gamma)\ \ \text{on }V,
	\]
	where $\widetilde\gamma$ denotes $\gamma$ viewed on $U\cap V$ and the two formulas agree on $U\cap V$.
	Equivalently, on $U\cap V$,
	\[
	\omega|_{U\cap V}= d\rho_V\wedge \gamma = -\,d\rho_U\wedge \gamma.
	\]
\end{theorem}

\begin{proof}
	Exactness of the short sequence of complexes gives the long exact sequence by standard homological algebra.
	For the explicit formula: start with $\gamma\in Z^{k-1}(U\cap V)$.
	Choose $(\alpha,\beta)\in\Omega^{k-1}(U)\oplus\Omega^{k-1}(V)$ such that $s(\alpha,\beta)=\gamma$.
	A concrete choice is obtained using $\rho_U+\rho_V=1$:
	\[
	\alpha:=\rho_V\,\gamma \ \text{(extended by $0$ to $U$)},\qquad
	\beta:=-\rho_U\,\gamma \ \text{(extended by $0$ to $V$)},
	\]
	so that on $U\cap V$ we have $\alpha-\beta=\gamma$.
	Apply $d$:
	\[
	s(d\alpha,d\beta)=d(\alpha-\beta)=d\gamma=0,
	\]
	hence $(d\alpha,d\beta)\in\ker(s)=\mathrm{im}(r)$, so there exists a global $\omega\in\Omega^k(M)$ with
	$r(\omega)=(d\alpha,d\beta)$. This $\omega$ represents $\delta([\gamma])$ by definition of the connecting morphism.
	On $U$ we have $\omega|_U=d\alpha=d(\rho_V\gamma)$ and on $V$ we have $\omega|_V=d\beta=-d(\rho_U\gamma)$.
	On the overlap,
	\[
	d(\rho_V\gamma)+d(\rho_U\gamma)=d\bigl((\rho_V+\rho_U)\gamma\bigr)=d\gamma=0,
	\]
	so the two prescriptions coincide and define a global form.
	If $d\gamma=0$, then $d(\rho_V\gamma)=d\rho_V\wedge\gamma$, yielding the stated simplification.
\end{proof}

\subsubsection{de Rham cohomology via sheaves}
\label{subsec:de-rham-sheaf-proof}

\paragraph{The de Rham resolution of the constant sheaf}

\begin{theorem}[de Rham resolution]
	\label{thm:deRham-resolution}
	The sequence of sheaves
	\[
	0\longrightarrow \underline{\mathbb R}\xrightarrow{\;\iota\;} \Omega^0\xrightarrow{\,d\,}\Omega^1\xrightarrow{\,d\,}\Omega^2\xrightarrow{\,d\,}\cdots
	\]
	is exact.
\end{theorem}

\begin{proof}
	Exactness at $\underline{\mathbb R}$: the inclusion $\iota$ is injective.
	Exactness at $\Omega^0$: $\ker(d:\Omega^0(U)\to\Omega^1(U))$ consists of locally constant functions, hence equals $\underline{\mathbb R}(U)$.
	Exactness at $\Omega^k$ for $k\ge 1$ is stalkwise: fix $x\in M$ and choose a coordinate neighborhood $U$ of $x$
	diffeomorphic to a star-shaped open set in $\mathbb R^n$; then Corollary~\ref{cor:local-exactness} yields
	$\ker(d:\Omega^k(U)\to\Omega^{k+1}(U))=\mathrm{im}(d:\Omega^{k-1}(U)\to\Omega^k(U))$.
	Thus the sequence is exact as a sequence of sheaves.
\end{proof}

\paragraph{Identification of de Rham cohomology with sheaf cohomology}

\begin{theorem}[de Rham theorem in sheaf language]
	\label{thm:deRham-sheaf}
	For every $k\ge 0$,
	\[
	H^k_{\mathrm{dR}}(M)\ \cong\ H^k(M,\underline{\mathbb R}).
	\]
\end{theorem}

\begin{proof}
	By Theorem~\ref{thm:deRham-resolution}, the de Rham complex is a resolution of $\underline{\mathbb R}$ by sheaves $\Omega^k$.
	By Corollary~\ref{cor:Omega-acyclic}, each $\Omega^k$ is acyclic.
	Therefore the right derived functors $H^k(M,\underline{\mathbb R})$ are computed by the cohomology of the global sections complex
	\[
	0\to \Omega^0(M)\xrightarrow{d}\Omega^1(M)\xrightarrow{d}\cdots,
	\]
	and the cohomology of this complex is precisely $H^k_{\mathrm{dR}}(M)$ by definition.
\end{proof}

\subsubsection{Computations on spheres and their generators}
\index{de Rham cohomology!sphere}
\label{subsec:de-rham-examples}

\paragraph{A worked Mayer--Vietoris computation for $S^1$}

Let $S^1=\{e^{i\theta}\}$.
Cover $S^1$ by two arcs $U,V$ such that $U\cap V$ has two connected components $W_+$ and $W_-$ (small arcs around $\theta=0$ and $\theta=\pi$).
Since $U$ and $V$ are contractible, $H^1_{\mathrm{dR}}(U)=H^1_{\mathrm{dR}}(V)=0$.
Since $U\cap V=W_+\sqcup W_-$ is a disjoint union of two intervals, $H^0_{\mathrm{dR}}(U\cap V)\cong \mathbb R\oplus\mathbb R$ and $H^1_{\mathrm{dR}}(U\cap V)=0$.
The Mayer--Vietoris sequence in low degrees gives
\[
0\to H^0_{\mathrm{dR}}(S^1)\to H^0_{\mathrm{dR}}(U)\oplus H^0_{\mathrm{dR}}(V)\to H^0_{\mathrm{dR}}(U\cap V)\xrightarrow{\delta} H^1_{\mathrm{dR}}(S^1)\to 0.
\]
Identifying $H^0$ with locally constant functions, we have
\[
H^0_{\mathrm{dR}}(U)\cong\mathbb R,\quad H^0_{\mathrm{dR}}(V)\cong\mathbb R,\quad H^0_{\mathrm{dR}}(U\cap V)\cong \mathbb R\oplus\mathbb R.
\]
The map $H^0(U)\oplus H^0(V)\to H^0(U\cap V)$ sends $(a,b)$ to $(a-b,\ a-b)$ on $(W_+,W_-)$.
Hence its image is the diagonal $\{(t,t)\}$, and the quotient $(\mathbb R\oplus\mathbb R)/\Delta\cong \mathbb R$.
Exactness implies
\[
H^1_{\mathrm{dR}}(S^1)\cong \mathbb R.
\]

To exhibit a generator explicitly, pick a partition of unity $\rho_U+\rho_V=1$ subordinate to $\{U,V\}$ and define on $U\cap V$
a locally constant function $f$ that equals $1$ on $W_+$ and $0$ on $W_-$. Then $df=0$ on $U\cap V$.
The connecting map $\delta$ (Theorem~\ref{thm:MV-LES}) sends $[f]\in H^0(U\cap V)$ to the class of the global $1$-form
\[
\omega=
\begin{cases}
	d(\rho_V f),&\text{on }U,\\
	-\,d(\rho_U f),&\text{on }V.
\end{cases}
\]
On $U\cap V$ this satisfies $\omega=f\,d\rho_V$.  The class is nonzero:
under the exact-sequence identification it is the image of $(1,0)$ modulo
the diagonal subspace of $H^0(W_+)\oplus H^0(W_-)$.
It is therefore a nonzero multiple of the standard angular class, represented
globally by
\[
\eta=x\,dy-y\,dx\quad\text{on }S^1\subset\mathbb R^2.
\]
In a local angular coordinate this form is written $d\theta$; the notation
does not assert that a single-valued real function $\theta$ exists on all of
$S^1$.

\paragraph{The computation of $H^\bullet_{\mathrm{dR}}(S^n)$ via hemispheres}

Let $S^n=U\cup V$ where $U$ and $V$ are the complements of the south and north poles, respectively.
Both $U$ and $V$ are diffeomorphic to $\mathbb R^n$, hence have only $H^0\cong\mathbb R$.
The overlap $U\cap V$ deformation retracts to $S^{n-1}$.
Inductively applying Mayer--Vietoris yields
\[
H^k_{\mathrm{dR}}(S^n)\cong
\begin{cases}
	\mathbb R,& k=0,n,\\
	0,&\text{otherwise}.
\end{cases}
\]

\paragraph{A concrete generator of $H^n_{\mathrm{dR}}(S^n)$}

Let $\mathrm{vol}_{S^n}$ be the Riemannian volume form for the round metric.
Since $d(\mathrm{vol}_{S^n})=0$ and $\int_{S^n}\mathrm{vol}_{S^n}\neq 0$, the class $[\mathrm{vol}_{S^n}]$ is nonzero.
Because $H^n_{\mathrm{dR}}(S^n)\cong\mathbb R$, it follows that $[\mathrm{vol}_{S^n}]$ generates the top cohomology.

A normalized generator is
\[
\alpha_n:=\frac{1}{\mathrm{Vol}(S^n)}\,\mathrm{vol}_{S^n},
\qquad \int_{S^n}\alpha_n=1.
\]

\subsubsection{\texorpdfstring{Dolbeault sheaves and the $\bar\partial$-resolution}{Dolbeault sheaves and the d-bar resolution}}
\label{subsec:dolbeault}

\paragraph{The Dolbeault complex of sheaves}

Let $M$ be a complex manifold of complex dimension $m$.
Write $\Omega^{p,q}$ for the sheaf of smooth $(p,q)$-forms. The operator $\bar\partial$ defines
sheaf morphisms $\bar\partial:\Omega^{p,q}\to\Omega^{p,q+1}$ and a complex
\[
0\to \Omega^p_{\mathrm{hol}} \to \Omega^{p,0}\xrightarrow{\bar\partial}\Omega^{p,1}\xrightarrow{\bar\partial}\cdots\xrightarrow{\bar\partial}\Omega^{p,m}\to 0,
\]
where $\Omega^p_{\mathrm{hol}}$ is the sheaf of holomorphic $p$-forms.

\begin{theorem}[Local exactness of the Dolbeault complex]
	\label{thm:dolbeault-lemma}
	For each $p$ and each $q\ge 1$, the sequence is exact at $\Omega^{p,q}$ in a neighborhood of every point.
	Equivalently, for every $x\in M$ there is a neighborhood $U\ni x$ such that for every $\omega\in\Omega^{p,q}(U)$ with $\bar\partial\omega=0$
	there exists $\eta\in\Omega^{p,q-1}(U)$ with $\omega=\bar\partial\eta$.
\end{theorem}

The analytic input is local solvability of $\bar\partial$ on a polydisc
for smooth forms. That local solution operator is not constructed here;
see \cite{Wells2008,GunningRossi1965}. The following sheaf argument uses
the local exactness just stated.

\paragraph{Dolbeault theorem}

Since each $\Omega^{p,q}$ is a fine sheaf (multiplication by partitions of unity), the Dolbeault complex is a fine resolution of $\Omega^p_{\mathrm{hol}}$.
Hence
\[
H^q(M,\Omega^p_{\mathrm{hol}})\cong H^{p,q}_{\bar\partial}(M).
\]
Indeed, the local Dolbeault lemma identifies the kernel in degree zero with
holomorphic $p$-forms and gives exactness in positive degrees.  Acyclicity of
the fine terms then computes the sheaf cohomology by the global complex
$\bigl(\Omega^{p,\bullet}(M),\bar\partial\bigr)$; its degree-$q$ cohomology is
the Dolbeault group on the right.

\subsubsection{Exercises}
\label{subsec:de-rham-exercises}

\begin{exercise}[Homotopy operator check in coordinates]
	Let $U\subset\mathbb R^n$ be star-shaped and $K$ be defined as in Theorem~\ref{thm:poincare-star}.
	Write $\omega=\sum_{|I|=k} a_I(x)\,dx^I$ and compute $K\omega$ explicitly in coordinates for $k=1,2$.
	Verify directly that $dK\omega+Kd\omega=\omega$.
\end{exercise}

\begin{exercise}[Mayer--Vietoris connecting map computation]
	Let $M=U\cup V$ and $\rho_U+\rho_V=1$ subordinate to $\{U,V\}$.
	For a closed $(k-1)$-form $\gamma$ on $U\cap V$, define $\omega$ by $\omega|_U=d(\rho_V\gamma)$ and $\omega|_V=-d(\rho_U\gamma)$.
	Prove that $\omega$ is globally well-defined and closed, and that its cohomology class depends only on $[\gamma]\in H^{k-1}_{\mathrm{dR}}(U\cap V)$.
\end{exercise}

\begin{exercise}[A worked computation: $H^1_{\mathrm{dR}}(S^1)$]
	Using the cover of $S^1$ by two arcs whose intersection has two components,
	compute $H^1_{\mathrm{dR}}(S^1)$ by identifying the image of $H^0(U)\oplus H^0(V)\to H^0(U\cap V)$ explicitly.
\end{exercise}

\begin{exercise}[Top cohomology generator via integration]
	Let $M$ be a closed oriented $n$-manifold with volume form $\mathrm{vol}$.
	Show: if $\alpha$ is an exact $n$-form, then $\int_M \alpha=0$.
	Deduce that if $\int_M \mathrm{vol}\neq 0$, then $[\mathrm{vol}]\neq 0$ in $H^n_{\mathrm{dR}}(M)$.
\end{exercise}

\begin{exercise}[Fine sheaves are \v{C}ech-acyclic]
	Let $\mathscr F$ be a fine sheaf, $\mathcal U=\{U_\alpha\}$ a locally finite cover, and $c\in Z^q(\mathcal U,\mathscr F)$ with $q\ge 1$.
	Define $b$ as in the proof of Theorem~\ref{thm:fine-cech-vanish} and prove carefully that $\delta b=c$.
\end{exercise}

\begin{exercise}[Compute $H^\bullet_{\mathrm{dR}}(\mathbb T^2)$ by Mayer--Vietoris]
	Cover $\mathbb T^2=S^1\times S^1$ by $U=(S^1\setminus\{\mathrm{pt}\})\times S^1$ and $V=(S^1\setminus\{\mathrm{pt}'\})\times S^1$.
	Compute the overlap $U\cap V$ up to homotopy type, then compute $H^\bullet_{\mathrm{dR}}(\mathbb T^2)$.
\end{exercise}

\EssaySuminKim
\EssayNamhoKim

\section{Canonical Bundles, Chern Classes, and Curvature}
\index{canonical bundle}\index{Chern classes}\index{curvature}
\label{sec:chern-canonical}

\subsection{Local expressions in complex coordinates}
\label{subsec:canonical-local}

\subsubsection{Complex coordinates, types, and basic tensors}

Let $X$ be a complex manifold of complex dimension $n$.
On a coordinate chart $U\subset X$ with holomorphic coordinates
\[
z=(z^1,\dots,z^n),\qquad z^j=x^j+i y^j,
\]
write the complexified cotangent bundle decomposition
\[
T^*X\otimes_{\mathbb R}\mathbb C \;=\; T^{*(1,0)}X \oplus T^{*(0,1)}X,
\]
and denote the local frames
\[
dz^1,\dots,dz^n\in \Gamma(U,T^{*(1,0)}X),\qquad
d\bar z^1,\dots,d\bar z^n\in \Gamma(U,T^{*(0,1)}X).
\]
A smooth $(p,q)$-form on $U$ is a finite sum
\[
\alpha=\sum_{|I|=p,\,|J|=q}\alpha_{I\bar J}(z)\,dz^{I}\wedge d\bar z^{J}.
\]

\subsubsection{The canonical bundle and its transition functions}

\begin{definition}[Canonical bundle]
	\label{def:canonical-bundle}
	The \emph{canonical bundle} of $X$ is the holomorphic line bundle
	\[
	K_X:=\Lambda^{n}T^{*(1,0)}X,
	\]
	whose local holomorphic frame on a coordinate chart $(U;z^1,\dots,z^n)$ is
	\[
	e_U:=dz^1\wedge\cdots\wedge dz^n.
	\]
\end{definition}

Let $(U;z)$ and $(V;w)$ be two overlapping holomorphic coordinate charts, so on $U\cap V$ we have
$w=w(z)$ holomorphic with Jacobian matrix $J=(\partial w^a/\partial z^b)_{a,b}$.
Then
\[
dw^1\wedge\cdots\wedge dw^n
=
\det\!\left(\frac{\partial w}{\partial z}\right)\,dz^1\wedge\cdots\wedge dz^n,
\]
so, with the convention used throughout this book, the transition function of
$K_X$ on $U\cap V$ is
\[
e_V = f_{UV}\,e_U,\qquad
f_{UV}=\det\!\left(\frac{\partial w}{\partial z}\right)\in \mathcal O^*(U\cap V).
\]
Equivalently, reversing the indices gives $e_U=f_{VU}e_V$, where
\[
f_{VU}=\det\!\left(\frac{\partial w}{\partial z}\right)^{-1}.
\]

\begin{proposition}[Cocycle condition for $K_X$]
	\label{prop:canonical-cocycle}
	On a triple overlap $U\cap V\cap W$ with coordinate changes $z\mapsto w(z)\mapsto u(w)$, the transition functions satisfy
	\[
	f_{UW}=f_{UV}\,f_{VW}.
	\]
\end{proposition}

\begin{proof}
	On $U\cap V\cap W$ the chain rule gives
	\[
	\frac{\partial u}{\partial z}=\frac{\partial u}{\partial w}\cdot \frac{\partial w}{\partial z}.
	\]
	Taking determinants,
	\[
	\det\!\left(\frac{\partial u}{\partial z}\right)
	=
	\det\!\left(\frac{\partial u}{\partial w}\right)\,
	\det\!\left(\frac{\partial w}{\partial z}\right).
	\]
	This is exactly the identity $f_{UW}=f_{UV}\,f_{VW}$ in the convention
	$e_V=f_{UV}e_U$ and $e_W=f_{VW}e_V$.
\end{proof}

\subsubsection{Hermitian metrics on line bundles in coordinates}

Let $L\to X$ be a holomorphic line bundle. A \emph{Hermitian metric} $h$ on $L$ assigns to each fiber $L_x$
a positive definite Hermitian inner product varying smoothly in $x$.
Choose a local holomorphic frame $e$ of $L$ on $U$. Any holomorphic section is $s=f\,e$ with $f\in \mathcal O(U)$.
Define the local weight function $\varphi$ by
\[
\|e\|_h^2 = h(e,e)=e^{-\varphi}\quad\text{on }U,
\]
so for $s=f e$ we have
\[
\|s\|_h^2 = |f|^2\,\|e\|_h^2 = |f|^2 e^{-\varphi}.
\]
If $\tilde e=g e$ with $g\in\mathcal O^*(U)$, then
\[
\|\tilde e\|_h^2 = |g|^2 \|e\|_h^2
\quad\mathrel{\text{and hence}}\quad
\tilde\varphi = \varphi - \log|g|^2.
\]

\subsection{Gaussian curvature and the Ricci form}
\label{subsec:ricci}

\subsubsection{The Chern connection on a Hermitian holomorphic line bundle}

Let $(L,h)$ be a Hermitian holomorphic line bundle.
A \emph{connection} on $L$ is a complex-linear map
$\nabla:\Gamma(L)\to\Omega^1(X;L)$ satisfying
$\nabla(fs)=df\otimes s+f\nabla s$ for smooth complex-valued $f$.
On a Riemann surface the complexified cotangent bundle splits into
$(1,0)$- and $(0,1)$-parts; projecting the one-form output of $\nabla$
defines $\nabla^{1,0}$ and $\nabla^{0,1}$. The holomorphic structure
similarly gives $\bar\partial$ on smooth sections: in a holomorphic frame
$e$, it sends $fe$ to $(\bar\partial f)e$.
There is a unique connection $\nabla$ on $L$ (the \emph{Chern connection}) such that
\[
\nabla^{0,1}=\bar\partial \quad\text{and}\quad \nabla h=0.
\]
In a local holomorphic frame $e$ on $U$, write $\nabla e=\theta\otimes e$ for a complex-valued $1$-form $\theta$.
Metric compatibility $\nabla h=0$ implies $\theta$ is determined by $\varphi$.
Indeed, if $e'=ge$ for a nowhere-vanishing holomorphic function $g$, then
$\nabla e'=(\theta+d\log g)\otimes e'$; thus the local connection form
changes, while the connection itself is independent of the frame.

\begin{proposition}[Local formula for the Chern connection on a line bundle]
	\label{prop:chern-connection-line}
	Let $e$ be a local holomorphic frame with $\|e\|_h^2=e^{-\varphi}$.
	Then the Chern connection is given by
	\[
	\nabla e = -\,\partial \varphi \otimes e,
	\]
	and for a local section $s=f e$,
	\[
	\nabla s = (\partial f - f\,\partial\varphi)\otimes e \;+\; (\bar\partial f)\otimes e.
	\]
\end{proposition}

\begin{proof}
	Write $\nabla e=\theta\otimes e$ with $\theta=\theta^{1,0}+\theta^{0,1}$.
	Since $e$ is holomorphic, $\bar\partial e=0$, and the condition $\nabla^{0,1}=\bar\partial$ gives $\theta^{0,1}=0$.
	So $\theta$ is of type $(1,0)$.
	Metric compatibility says for the function $h(e,e)=e^{-\varphi}$,
	\[
	d\,h(e,e) = h(\nabla e,e)+h(e,\nabla e)=\theta\,h(e,e)+\overline{\theta}\,h(e,e)=(\theta+\overline{\theta})e^{-\varphi}.
	\]
	But $d(e^{-\varphi})=-e^{-\varphi}d\varphi$, hence
	\[
	-e^{-\varphi}d\varphi = (\theta+\overline{\theta})e^{-\varphi}
	\quad\mathrel{\text{and hence}}\quad
	\theta+\overline{\theta}=-d\varphi.
	\]
	Taking the $(1,0)$-part, and using $\theta$ is $(1,0)$, we get
	\[
	\theta = -\partial\varphi.
	\]
	This gives $\nabla e=-\partial\varphi\otimes e$.
	For $s=f e$, use the Leibniz rule $\nabla(f e)=df\otimes e+f\nabla e$ and decompose $df=\partial f+\bar\partial f$.
\end{proof}

\subsubsection{Curvature of the Chern connection and $c_1(L)$}

The curvature is $F_\nabla=\nabla^2$, a $(1,1)$-form with values in $\mathrm{End}(L)\cong \mathbb C$.

\begin{proposition}[Curvature formula for a Hermitian line bundle]
	\label{prop:curv-line}
	In the above notation,
	\[
	F_\nabla = \bar\partial(-\partial\varphi)=\partial\bar\partial\varphi,
	\]
	and the real $(1,1)$-form
	\[
	\Theta_h(L):=\frac{i}{2\pi}F_\nabla = \frac{i}{2\pi}\,\partial\bar\partial\varphi
	\]
	represents the first Chern class:
	\[
	[\Theta_h(L)] = c_1(L)\in H^2(X,\mathbb Z)\subset H^2(X,\mathbb R).
	\]
\end{proposition}

\begin{proof}
	From Proposition~\ref{prop:chern-connection-line}, the connection $1$-form is $\theta=-\partial\varphi$.
	Since $\theta$ is scalar-valued, $\theta\wedge\theta=0$, and therefore
	\[
	F_\nabla = d\theta = \bar\partial\theta + \partial\theta = \bar\partial(-\partial\varphi) + \partial(-\partial\varphi).
	\]
	But $\partial^2=0$, so $\partial(-\partial\varphi)=0$, hence $F_\nabla=\bar\partial(-\partial\varphi)=\partial\bar\partial\varphi$.
	
	To see the cohomology class is integral and equals $c_1(L)$, cover $X$ by $\{U_\alpha\}$ with holomorphic frames $e_\alpha$.
	On overlaps $e_\beta=g_{\alpha\beta}e_\alpha$ with $g_{\alpha\beta}\in\mathcal O^*(U_{\alpha\beta})$.
	Write $\|e_\alpha\|^2=e^{-\varphi_\alpha}$. Then
	\[
	e^{-\varphi_\beta}=\|e_\beta\|^2=\|g_{\alpha\beta}e_\alpha\|^2=|g_{\alpha\beta}|^2 e^{-\varphi_\alpha}
	\quad\mathrel{\text{and hence}}\quad
	\varphi_\beta=\varphi_\alpha-\log|g_{\alpha\beta}|^2.
	\]
	Hence on $U_{\alpha\beta}$,
	\[
	\partial\bar\partial\varphi_\beta=\partial\bar\partial\varphi_\alpha-\partial\bar\partial\log|g_{\alpha\beta}|^2.
	\]
	But $\log|g|^2=\log g + \log\bar g$ locally, and since $g$ is holomorphic,
	$\bar\partial\log g=0$ and $\partial\log \bar g=0$, giving $\partial\bar\partial\log|g|^2=0$.
	Thus the local forms $\partial\bar\partial\varphi_\alpha$ glue to a global $(1,1)$-form $F_\nabla$.
	
	Finally, comparing with the transition data in $\mathcal O^*$ and the exponential sequence
	(constructed in \S\ref{sec:exponential-sequence}), one identifies the \v{C}ech class
	$\{g_{\alpha\beta}\}\in \check H^1(X,\mathcal O^*)$ with a class in $H^2(X,\mathbb Z)$ by the connecting morphism,
	and the de Rham representative of that class is exactly $\frac{i}{2\pi}F_\nabla$.
\end{proof}

\subsubsection{Ricci form and Gaussian curvature on Riemann surfaces}

Now assume $X$ is a Riemann surface ($n=1$) and let $\omega$ be a K\"ahler form coming from a Hermitian metric
\[
ds^2 = \lambda(z)\,|dz|^2,\qquad \lambda(z)>0.
\]
Then the associated K\"ahler form is
\[
\omega = \frac{i}{2}\lambda(z)\,dz\wedge d\bar z.
\]

\begin{definition}[Ricci form]
	\label{def:ricci-form}
	On a K\"ahler manifold with local metric matrix $(g_{j\bar k})$, the \emph{Ricci form} is
	\[
	\mathrm{Ric}(\omega) := -\, i\,\partial\bar\partial \log \det(g_{j\bar k}).
	\]
\end{definition}

For $n=1$, $\det(g_{1\bar 1})=g_{1\bar 1}=\lambda/2$ up to the chosen convention, and hence
\[
\mathrm{Ric}(\omega) = -\, i\,\partial\bar\partial \log \lambda.
\]

\begin{proposition}[Curvature of $K_X$ equals Ricci form]
	\label{prop:curv-K-ricci}
	Let $X$ be a K\"ahler manifold. Equip $K_X$ with the Hermitian metric induced by $\omega$.
	Then the Chern curvature of $(K_X,h)$ satisfies
	\[
	\frac{i}{2\pi}F_{K_X} = -\frac{1}{2\pi}\,\mathrm{Ric}(\omega).
	\]
	In particular,
	\[
	c_1(K_X) = -\left[\frac{1}{2\pi}\mathrm{Ric}(\omega)\right],\qquad
	c_1(TX)= \left[\frac{1}{2\pi}\mathrm{Ric}(\omega)\right].
	\]
\end{proposition}

\begin{proof}
	On a holomorphic coordinate chart, the metric gives a Hermitian matrix $(g_{j\bar k})$ on $T^{1,0}X$.
	The induced metric on $K_X=\Lambda^n T^{*(1,0)}X$ satisfies
	\[
	\|dz^1\wedge\cdots\wedge dz^n\|^2 = \det(g_{j\bar k})^{-1}.
	\]
	Thus in the notation of Proposition~\ref{prop:chern-connection-line}, the weight is
	\[
	\|e\|^2=e^{-\varphi}=\det(g_{j\bar k})^{-1}
	\quad\mathrel{\text{and hence}}\quad
	\varphi=\log\det(g_{j\bar k}).
	\]
	Therefore
	\[
	F_{K_X}=\partial\bar\partial\varphi=\partial\bar\partial\log\det(g_{j\bar k}),
	\]
	and multiplying by $i$ gives
	\[
	iF_{K_X}= i\partial\bar\partial\log\det(g_{j\bar k}) = -\,\mathrm{Ric}(\omega),
	\]
	which implies $\frac{i}{2\pi}F_{K_X}=-\frac{1}{2\pi}\mathrm{Ric}(\omega)$.
	The Chern class statements follow from Proposition~\ref{prop:curv-line} and the identity $c_1(TX)=-c_1(K_X)$.
\end{proof}

\begin{proposition}[Gaussian curvature in local coordinates]
	\label{prop:gauss-curv}
	Let $X$ be a Riemann surface with metric $ds^2=\lambda(z)|dz|^2$.
	Then the Gaussian curvature $K$ satisfies
	\[
	K = -\frac{1}{2\lambda}\,\Delta \log \lambda,
	\]
	where $\Delta=\partial_x^2+\partial_y^2=4\,\partial_z\partial_{\bar z}$ uses the classical
	Euclidean sign convention. (The nonnegative Hodge Laplacian used elsewhere is its negative on functions.)
	Equivalently,
	\[
	\mathrm{Ric}(\omega) = K\,\omega.
	\]
\end{proposition}

\begin{proof}
	Write $\lambda=e^{2u}$ so that $ds^2=e^{2u}(dx^2+dy^2)$ in $z=x+iy$.
	A direct computation of the Levi--Civita connection for the conformal metric yields
	\[
	K = -e^{-2u}\,\Delta_{\mathbb R^2} u.
	\]
	Since $\log\lambda=2u$ and $\Delta_{\mathbb R^2}=4\partial_z\partial_{\bar z}$, this becomes
	\[
	K=-\frac{1}{\lambda}\,4\partial_z\partial_{\bar z}(\tfrac12\log\lambda)
	=-\frac{1}{2\lambda}\,\Delta \log\lambda.
	\]
	For the Ricci identity, in complex dimension one we have $\mathrm{Ric}(\omega)=-i\partial\bar\partial\log\lambda$ and
	\[
	-i\partial\bar\partial\log\lambda
	=
	\left(-\frac{1}{2\lambda}\,\Delta \log\lambda\right)\,\frac{i}{2}\lambda\,dz\wedge d\bar z
	=K\,\omega.
	\]
\end{proof}

\subsection{Higher Chern classes, the Chern character, and the Todd class}
\label{subsec:chern-weil-character-todd}

The preceding line-bundle formula is the rank-one case of a matrix
construction. Let $E\to X$ be a rank-$r$ holomorphic vector bundle with
a Hermitian metric and Chern connection $\nabla$. In a local unitary
frame its curvature $F_\nabla=\nabla^2$ is an $r\times r$ matrix of
two-forms. Changing frame conjugates the matrix, so its invariant
polynomials give global forms. Put $A=(i/2\pi)F_\nabla$. Products of
entries of $A$ are wedge products; two-forms commute under the wedge
product because their degree is even.

\begin{proposition}[Chern--Weil representatives of the Chern classes]
 \label{prop:all-chern-weil}
 The coefficient $c_k(\nabla)$ in
 \[
  \det(I_r+tA)=\sum_{k=0}^{r}t^k c_k(\nabla)
 \]
 is a closed form of degree $2k$, and its de Rham class is the image
 of $c_k(E)\in H^{2k}(X;\Z)$. Equivalently it is the $k$th elementary
 symmetric polynomial in the formal curvature roots of $A$.
 In particular,
 \[
 c_0(\nabla)=1,\quad c_1(\nabla)=\operatorname{tr}A,\quad
 c_2(\nabla)=\tfrac12\bigl((\operatorname{tr}A)^2-
                         \operatorname{tr}(A\wedge A)\bigr).
 \]
\end{proposition}

\begin{proof}
 Conjugation invariance makes the coefficients independent of frame.
 The Bianchi identity $\nabla F_\nabla=0$ implies that every invariant
 polynomial in $F_\nabla$ is closed. Along a path of connections, the
 derivative of such a polynomial is an exact transgression form; hence
 its de Rham class is independent of the metric and connection.
 On a splitting space where the pullback of $E$ becomes a sum of line
 bundles, connection independence lets us choose their direct-sum
 connection. Its curvature matrix is diagonal, and the determinant becomes
 $\prod_j(1+t(i/2\pi)F_{L_j})$. Its coefficients are the elementary
 symmetric polynomials in the line-bundle first Chern classes.
 The splitting principle identifies these with the integral classes
 $c_k(E)$. The general Bianchi, transgression, and splitting-principle
 results invoked here are proved in \cite[Chapters~IV--V]{BottTu1982}.
\end{proof}

For rank one this recovers Proposition~\ref{prop:curv-line}. On a
Riemann surface, $H^{2k}=0$ for $k\ge2$; on a complex surface $c_2$
can be a nonzero top-degree class.

\begin{definition}[Chern character]
 \label{def:chern-character}
 If $x_1,\ldots,x_r$ are the formal Chern roots, so that
 $c(E)=\prod_j(1+x_j)$, define
 \[
  \operatorname{ch}(E)=\sum_{j=1}^r e^{x_j}
  =r+c_1(E)+\tfrac12\bigl(c_1(E)^2-2c_2(E)\bigr)+\cdots.
 \]
 Its curvature representative is
 $\operatorname{tr}\exp((i/2\pi)F_\nabla)$.
 The series truncates in each fixed manifold dimension.
\end{definition}

The exponential makes $\operatorname{ch}$ additive for direct sums and
multiplicative for tensor products. For example,
$\operatorname{ch}(L\otimes L')=e^{c_1(L)+c_1(L')}$; the splitting
principle extends this identity to arbitrary bundles.

\begin{definition}[Todd class]
 \label{def:todd-class}
 With the same formal roots, define
 \[
  \operatorname{td}(E)=\prod_{j=1}^r\frac{x_j}{1-e^{-x_j}},
  \qquad \frac{x}{1-e^{-x}}=1+\frac{x}{2}+\frac{x^2}{12}+\cdots.
 \]
 The quotient means its formal power series at $x=0$, with constant
 term $1$. Through real degree four,
 \[
  \operatorname{td}(E)=1+\tfrac12c_1(E)
       +\tfrac1{12}\bigl(c_1(E)^2+c_2(E)\bigr)+\cdots.
 \]
 This is a rational characteristic class; replacing the roots by
 those of $(i/2\pi)F_\nabla$ gives its differential-form representative.
\end{definition}

\paragraph{Why these combinations are natural: $\mathbb P^2$.}
Let $H=c_1(\mathcal O_{\mathbb P^2}(1))$. Then
$H^*(\mathbb P^2;\Z)=\Z[H]/(H^3)$ and $\int_{\mathbb P^2}H^2=1$.
The Euler sequence
\[
 0\longrightarrow\mathcal O\longrightarrow
 \mathcal O(1)^{\oplus3}\longrightarrow T\mathbb P^2
 \longrightarrow0
\]
and the Whitney product formula give
$c(T\mathbb P^2)=(1+H)^3=1+3H+3H^2$.
Here the first arrow sends a local unit section to the three
homogeneous-coordinate sections; its quotient is the tangent bundle
\cite{GriffithsHarris1994}.
Thus $c_1(T\mathbb P^2)=3H$ and $c_2(T\mathbb P^2)=3H^2$.
The latter integrates to $3$, the \emph{topological} Euler
characteristic $\chi_{\mathrm{top}}(\mathbb P^2)=1+1+1=3$ from
its cells in dimensions $0,2,4$. By contrast, the \emph{holomorphic}
Euler characteristic $\chi(\mathcal O_{\mathbb P^2})=
h^0-h^1+h^2$ is $1$: constants give $h^0=1$ and $h^1=h^2=0$
\cite{GriffithsHarris1994}. Riemann--Roch concerns this second number,
not the topological value $3$.

The Todd correction matches it exactly:
\[
 \operatorname{td}(T\mathbb P^2)
 =1+\tfrac32H+\frac{(3H)^2+3H^2}{12}
 =1+\tfrac32H+H^2,
 \qquad\int_{\mathbb P^2}\operatorname{td}(T\mathbb P^2)=1.
\]
For $L=\mathcal O(d)$, its character is
$e^{dH}=1+dH+\tfrac12d^2H^2$. Multiplying and retaining only the
degree-four component gives
\[
 \int_{\mathbb P^2}\operatorname{ch}(\mathcal O(d))
          \operatorname{td}(T\mathbb P^2)
 =1+\tfrac32d+\tfrac12d^2=\binom{d+2}{2}.
\]
For $d\ge0$, this counts homogeneous degree-$d$ monomials in three
variables, which form a basis of $H^0(\mathbb P^2,\mathcal O(d))$;
the higher cohomology vanishes. The exponential in
$\operatorname{ch}$ records the twist by $\mathcal O(d)$, while
$\operatorname{td}$ supplies the tangent-bundle correction already
needed at $d=0$. This example motivates, but does not prove, the
general Hirzebruch--Riemann--Roch formula
$\chi(E)=\int_X\operatorname{ch}(E)\operatorname{td}(TX)$
\cite{GriffithsHarris1994,BottTu1982}. The surface formula is proved
later by a different route in \S\ref{subsec:RRsurf:RR}.

\subsection{The exponential sequence and the degree/Chern class map}
\label{subsec:chern-exp}

\subsubsection{The connecting map $\delta:H^1(X,\mathcal O^*)\to H^2(X,\mathbb Z)$}

Recall the exponential short exact sequence of sheaves on a complex manifold $X$:
\[
0 \longrightarrow \underline{\mathbb Z}
\xrightarrow{\ \iota\ } \mathcal O
\xrightarrow{\ \exp(2\pi i\,\cdot)\ } \mathcal O^\times
\longrightarrow 0.
\]
It gives a connecting homomorphism
\[
\delta:\ H^1(X,\mathcal O^\times)\longrightarrow H^2(X,\underline{\mathbb Z})\cong H^2(X,\mathbb Z),
\]
and $\delta([L])$ is, by definition, $c_1(L)$.

\subsubsection{A \v{C}ech description of $c_1(L)$ on a good cover}

Let $\{U_\alpha\}$ be a good cover (all finite intersections contractible).
Let $L$ be given by a cocycle $g_{\alpha\beta}\in \mathcal O^*(U_{\alpha\beta})$.
Because $U_{\alpha\beta}$ is simply connected, choose holomorphic logarithms
\[
\ell_{\alpha\beta}\in\mathcal O(U_{\alpha\beta}),\qquad e^{2\pi i\,\ell_{\alpha\beta}}=g_{\alpha\beta}.
\]
On triple overlaps $U_{\alpha\beta\gamma}$,
\[
g_{\alpha\beta}g_{\beta\gamma}g_{\gamma\alpha}=1
\quad\mathrel{\text{and hence}}\quad
e^{2\pi i(\ell_{\alpha\beta}+\ell_{\beta\gamma}+\ell_{\gamma\alpha})}=1,
\]
so the holomorphic function $\ell_{\alpha\beta}+\ell_{\beta\gamma}+\ell_{\gamma\alpha}$ takes values in $\mathbb Z$.
Since $U_{\alpha\beta\gamma}$ is connected, it is constant:
\[
n_{\alpha\beta\gamma}:=\ell_{\alpha\beta}+\ell_{\beta\gamma}+\ell_{\gamma\alpha}\in \mathbb Z.
\]
The integers $\{n_{\alpha\beta\gamma}\}$ form a \v{C}ech $2$-cocycle with values in $\underline{\mathbb Z}$,
and its cohomology class is $c_1(L)$.

\begin{proposition}[Chern class from logarithms on a good cover]
	\label{prop:c1-from-log}
	With notation as above, the class $[\{n_{\alpha\beta\gamma}\}]\in \check H^2(X,\underline{\mathbb Z})\cong H^2(X,\mathbb Z)$
	equals $c_1(L)$.
\end{proposition}

\begin{proof}
	The connecting morphism in \v{C}ech cohomology is defined as follows.
	Start with the class $[g]\in \check H^1(X,\mathcal O^*)$ represented by $\{g_{\alpha\beta}\}$.
	Choose a $1$-cochain $\{\ell_{\alpha\beta}\}$ in $\mathcal O$ on overlaps with
	$\exp(2\pi i\ell_{\alpha\beta})=g_{\alpha\beta}$.
	Then apply the coboundary:
	\[
	(\delta \ell)_{\alpha\beta\gamma}=\ell_{\beta\gamma}-\ell_{\alpha\gamma}+\ell_{\alpha\beta}
	=\ell_{\alpha\beta}+\ell_{\beta\gamma}+\ell_{\gamma\alpha}.
	\]
	Because $\exp(2\pi i(\delta \ell)_{\alpha\beta\gamma})=1$, this $\delta\ell$ lands in $\underline{\mathbb Z}$,
	and the resulting $2$-cocycle is exactly $\{n_{\alpha\beta\gamma}\}$.
	This is, by definition of the connecting morphism, $\delta([g])$, hence equals $c_1(L)$.
\end{proof}

\subsubsection{Degree on a compact Riemann surface and $c_1(L)$}

Let $X$ be a compact Riemann surface and $L$ a holomorphic line bundle.
The fundamental class $[X]$ generates $H_2(X,\mathbb Z)$, and evaluation
on $[X]$ identifies $H^2(X,\mathbb Z)$ with $\mathbb Z$.
Define the \emph{degree} by
\[
\deg(L):=\langle c_1(L),[X]\rangle \in \mathbb Z.
\]

\begin{proposition}[Chern--Weil formula for the degree]
	\label{prop:degree-CW}
	Let $X$ be a compact Riemann surface and $(L,h)$ a Hermitian holomorphic line bundle with Chern curvature $F_\nabla$.
	Then
	\[
	\deg(L)=\int_X \frac{i}{2\pi}F_\nabla.
	\]
\end{proposition}

\begin{proof}
	By Proposition~\ref{prop:curv-line}, $\frac{i}{2\pi}F_\nabla$ represents $c_1(L)$ in de Rham cohomology.
	Pairing with $[X]$ is integration over $X$; since $c_1(L)$ is integral, the value lies in $\mathbb Z$.
	Thus
	\[
	\deg(L)=\left\langle \left[\frac{i}{2\pi}F_\nabla\right],[X]\right\rangle=\int_X \frac{i}{2\pi}F_\nabla.
	\]
\end{proof}

\subsection{Canonical divisors}
\index{canonical divisors}
\label{subsec:canonical-divisor}

\subsubsection{Divisors from meromorphic sections of line bundles}

Let $X$ be a compact Riemann surface.
A \emph{meromorphic section} $s$ of a holomorphic line bundle $L$ assigns in each chart a meromorphic function relative to a local frame.
Zeros and poles define an integer-valued divisor.

\begin{definition}[Divisor of a meromorphic section]
	\label{def:divisor-section}
	Let $s$ be a nonzero meromorphic section of $L\to X$.
	For each $p\in X$, choose a coordinate $z$ centered at $p$ and a holomorphic frame $e$ of $L$ on a neighborhood of $p$.
	Write $s=f e$ with $f$ meromorphic. Then there exists an integer $\nu_p(s)\in\mathbb Z$ such that
	\[
	f(z)=z^{\nu_p(s)}u(z),
	\]
	where $u$ is holomorphic and nowhere vanishing near $0$.
	Define the divisor
	\[
	\mathrm{div}(s):=\sum_{p\in X} \nu_p(s)\,[p].
	\]
\end{definition}

\begin{proposition}[Well-definedness]
	\label{prop:div-well-defined}
	The integer $\nu_p(s)$ in Definition~\ref{def:divisor-section} is independent of the choices of coordinate $z$
	and frame $e$.
\end{proposition}

\begin{proof}
	If $z$ and $\tilde z$ are two local coordinates with $\tilde z=\tilde z(z)$ and $\tilde z(0)=0$, then
	$\tilde z(z)=a z+\text{(higher order)}$ with $a\neq 0$, hence $\tilde z=z\cdot v(z)$ for a holomorphic unit $v$.
	Thus $z^{m}$ and $\tilde z^{m}$ differ by a unit, so the order $\nu_p(f)$ is unchanged.
	If $e$ and $\tilde e=g e$ are two holomorphic frames with $g$ a holomorphic unit, and $s=f e=\tilde f \tilde e$, then
	$\tilde f = f g^{-1}$, and multiplication by a unit does not change the vanishing/pole order.
	Therefore $\nu_p(s)$ is well-defined.
\end{proof}

\begin{definition}[Canonical divisor]
	\label{def:canonical-divisor}
	A \emph{canonical divisor} on $X$ is a divisor of the form $\mathrm{div}(\omega)$ where $\omega$ is a nonzero meromorphic $1$-form,
	i.e. a nonzero meromorphic section of $K_X$.
\end{definition}

\begin{proposition}[Degree of the canonical bundle]
	\label{prop:deg-K}
	If $X$ has genus $g$, then
	\[
	\deg(K_X)=2g-2.
	\]
\end{proposition}

\begin{proof}
	Choose any Hermitian metric on $X$ with K\"ahler form $\omega_X$.
	By Proposition~\ref{prop:curv-K-ricci},
	\[
	c_1(K_X)=-\left[\frac{1}{2\pi}\mathrm{Ric}(\omega_X)\right].
	\]
	Therefore, by Proposition~\ref{prop:degree-CW} applied to $K_X$,
	\[
	\deg(K_X)=-\int_X \frac{1}{2\pi}\mathrm{Ric}(\omega_X).
	\]
	On a Riemann surface, $\mathrm{Ric}(\omega_X)=K\,\omega_X$ (Proposition~\ref{prop:gauss-curv}), hence
	\[
	\deg(K_X)=-\frac{1}{2\pi}\int_X K\,\omega_X.
	\]
	By the Gauss--Bonnet theorem,
	\[
	\int_X K\,dA = 2\pi\chi(X)=2\pi(2-2g),
	\]
	and since $\omega_X$ agrees with the area form $dA$ under the metric convention,
	\[
	\deg(K_X)=-(2-2g)=2g-2.
	\]
\end{proof}

\subsection{The Calabi--Yau condition and holomorphic volume forms}
\index{Calabi--Yau manifold}
\index{holomorphic volume form}
\label{subsec:calabi-yau-volume-form}

Let $X$ be a connected complex manifold of complex dimension $n$. Its
canonical sheaf is
\[
	\omega_X:=\Omega_X^n,
\]
the sheaf of holomorphic $(n,0)$-forms, and its associated line bundle is
$K_X=\bigwedge^n(T_X^{1,0})^*$. A local coordinate system
$(z^1,\ldots,z^n)$ supplies the local frame
$dz^1\wedge\cdots\wedge dz^n$. Under a holomorphic coordinate change
$w=w(z)$,
\[
	dw^1\wedge\cdots\wedge dw^n
	=\det\!\left(\frac{\partial w^i}{\partial z^j}\right)
	dz^1\wedge\cdots\wedge dz^n.
\]
Thus the Jacobian determinant is precisely the transition function of the
canonical bundle.

\begin{definition}[Holomorphic volume form]
	A \emph{holomorphic volume form} on an $n$-dimensional complex manifold is
	a global holomorphic section $\Omega\in H^0(X,K_X)$. It is
	\emph{nowhere vanishing} if $\Omega_p\ne0$ in the one-dimensional fiber
	$(K_X)_p$ for every $p\in X$. In local coordinates it then has the form
	\[
	\Omega=f(z),dz^1\wedge\cdots\wedge dz^n,
	\qquad f(z)\ne0.
	\]
\end{definition}

\begin{definition}[Calabi--Yau condition used in this book]
	In this text, a complex manifold satisfies the \emph{holomorphic
	Calabi--Yau condition} if its canonical line bundle is holomorphically
	trivial:
	\[
	K_X\simeq\mathcal O_X.
	\]
	Some authors additionally require compactness, a K\"ahler metric, simple
	connectivity, or a prescribed holonomy group. Those extra hypotheses are
	not part of the equivalence proved below.
\end{definition}

\begin{theorem}[Three equivalent forms of the Calabi--Yau condition]
	\label{thm:CY-volume-trivial-canonical}
	For a connected complex $n$-manifold $X$, the following are equivalent.
	\begin{enumerate}
		\item The canonical line bundle is holomorphically trivial:
		$K_X\simeq X\times\C$.
		\item The canonical sheaf is isomorphic to the structure sheaf as an
		$\mathcal O_X$-module:
		$\omega_X\simeq\mathcal O_X$.
		\item There exists a global nowhere-vanishing holomorphic $n$-form
		$\Omega\in H^0(X,K_X)$.
	\end{enumerate}
\end{theorem}

\begin{proof}
	The equivalence of (1) and (2) is the line-bundle/sheaf dictionary: a
	holomorphic line bundle is trivial exactly when its sheaf of holomorphic
	sections is a free $\mathcal O_X$-module of rank one.

	Assume (2), and let
	$\Phi:\mathcal O_X\to\omega_X$ be an $\mathcal O_X$-module isomorphism.
	Set
	\[
	\Omega:=\Phi(1).
	\]
	For each point $p$, the induced map on stalks
	$\Phi_p:\mathcal O_{X,p}\to\omega_{X,p}$ is an isomorphism of free
	rank-one modules. Hence $\Omega_p=\Phi_p(1)$ is a generator of
	$\omega_{X,p}$ and cannot vanish. The local sections $\Phi_U(1)$ agree
	under restriction, so they form one global holomorphic volume form. This
	proves (2)$\Rightarrow$(3).

	Conversely, assume (3). Define on every open set $U\subset X$
	\[
	\Phi_U:\mathcal O_X(U)\longrightarrow\omega_X(U),
	\qquad f\longmapsto f\Omega|_U.
	\]
	These maps commute with restrictions and are $\mathcal O_X(U)$-linear.
	They are injective because $\Omega$ is nowhere zero. To prove surjectivity,
	choose local coordinates on $U$ and write
	\[
	\Omega=g(z)\,dz^1\wedge\cdots\wedge dz^n,
	\qquad
	\eta=h(z)\,dz^1\wedge\cdots\wedge dz^n.
	\]
	Since $g$ is holomorphic and nowhere zero, $h/g$ is holomorphic, and
	$\eta=(h/g)\Omega$. Therefore $\Phi$ is an isomorphism
	$\mathcal O_X\simeq\omega_X$, proving (3)$\Rightarrow$(2).
\end{proof}

\begin{corollary}[Uniqueness up to an invertible holomorphic function]
	If $\Omega$ and $\Omega'$ are nowhere-vanishing holomorphic volume forms,
	then $\Omega'=f\Omega$ for a unique $f\in\mathcal O_X^*(X)$. If $X$ is
	compact and connected, then $f$ is a nonzero constant.
\end{corollary}

\begin{proof}
	The quotient is locally the ratio of the coefficient functions and is
	holomorphic because the denominator never vanishes. The local ratios agree
	on overlaps. On a compact connected complex manifold, every holomorphic
	function is constant.
\end{proof}

\begin{corollary}[The Riemann-surface case]
	\label{cor:CY-curves-genus-one}
	For a compact connected Riemann surface $X$, the following are equivalent:
	\begin{gather*}
	K_X\simeq\mathcal O_X,\qquad
	\omega_X\simeq\mathcal O_X,\\
	\text{$X$ has a nowhere-vanishing holomorphic $1$-form},
	\qquad g(X)=1.
	\end{gather*}
\end{corollary}

\begin{proof}
	The first three conditions are equivalent by
	Theorem~\ref{thm:CY-volume-trivial-canonical}. If they hold, then
	$0=\deg K_X=2g-2$, so $g=1$. Conversely, when $g=1$, the vector space
	$H^0(X,K_X)$ has dimension one. Choose $0\ne\Omega\in H^0(X,K_X)$.
	Its zero divisor is effective and has degree
	$\deg K_X=2g-2=0$, hence it is empty. Thus $\Omega$ is nowhere vanishing.
\end{proof}

\begin{example}[Complex tori]
	Let $X=\C^n/\Lambda$ be a complex torus. Translation preserves each
	$dz^j$, so
	\[
	\Omega=dz^1\wedge\cdots\wedge dz^n
	\]
	descends to a nowhere-vanishing holomorphic volume form on $X$. Hence
	$K_X\simeq\mathcal O_X$. For $n=1$ this recovers the elliptic-curve form
	$dz$; on a plane cubic it corresponds to $dx/y$.
\end{example}

\begin{example}[The sphere and higher-genus curves]
	The sphere does not satisfy the condition because
	$K_{\PP^1}\simeq\mathcal O(-2)$. A compact curve of genus $g\ge2$ also
	fails it: every nonzero holomorphic one-form has an effective zero divisor
	of degree $2g-2>0$, so no such form can be everywhere nonzero.
\end{example}

\begin{remark}[$c_1=0$ versus holomorphic triviality]
	Holomorphic triviality implies $c_1(K_X)=0$, but the converse need not hold
	without additional hypotheses: a topologically trivial holomorphic line
	bundle can be a nontrivial element of $\mathrm{Pic}^0(X)$. Even a torsion
	canonical bundle can be nontrivial. Thus the volume-form condition is
	stronger than the bare cohomological equality $c_1(X)=0$. In the compact
	K\"ahler setting, Yau's theorem concerns the existence of a Ricci-flat
	K\"ahler metric when the real first Chern class vanishes; the theorem above
	is instead the elementary holomorphic equivalence that identifies a global
	volume form with a trivialization of the canonical sheaf.
\end{remark}

\subsection{Metrics and concrete examples}
\label{subsec:canonical-examples}

\subsubsection{Example: $K_{\mathbb{CP}^1}\cong \mathcal O(-2)$ and $\deg(K)=-2$}

Cover $\mathbb{CP}^1$ by $U_0=\{Z_0\neq 0\}$ and $U_1=\{Z_1\neq 0\}$ with coordinates
\[
z=\frac{Z_1}{Z_0}\ \text{ on }U_0,\qquad w=\frac{Z_0}{Z_1}=\frac{1}{z}\ \text{ on }U_1.
\]
Then
\[
dw = d\!\left(\frac{1}{z}\right) = -\frac{1}{z^2}\,dz.
\]
Hence the transition for $K_{\mathbb{CP}^1}$ is
\[
e_1=dw = -\frac{1}{z^2}\,dz = \left(-z^{-2}\right)e_0,
\]
so the transition function is $f_{01}=-z^{-2}$. The constant sign is removed
by changing one local frame, leaving the standard cocycle $z^{-2}$ of
$\mathcal O(-2)$.
Thus
\[
K_{\mathbb{CP}^1}\cong \mathcal O_{\mathbb{CP}^1}(-2),
\qquad
\deg(K_{\mathbb{CP}^1})=-2.
\]

\subsubsection{Example: Flat torus has trivial canonical bundle}

Let $X=\mathbb C/(\mathbb Z+\tau\mathbb Z)$ with $\Im\tau>0$.
The holomorphic $1$-form $dz$ on $\mathbb C$ is invariant under translations $z\mapsto z+1$ and $z\mapsto z+\tau$,
so it descends to a nowhere-vanishing holomorphic $1$-form on $X$.
Therefore $K_X$ has a nowhere-vanishing holomorphic section and is holomorphically trivial:
\[
K_X\cong \mathcal O_X,\qquad \deg(K_X)=0,
\]
consistent with $g=1$ in Proposition~\ref{prop:deg-K}.

\subsubsection{Example: Hyperbolic metric and positivity of $K_X$ for $g\ge 2$}

If $X$ has genus $g\ge 2$, it admits a metric of constant curvature $K\equiv -1$.
Then
\[
\deg(K_X)=-\frac{1}{2\pi}\int_X K\,dA
=\frac{\mathrm{Area}(X)}{2\pi}=2g-2>0,
\]
where the last equality is Gauss--Bonnet. Thus the canonical bundle is positive; this geometric
calculation agrees with Proposition~\ref{prop:deg-K} and with the divisor computation used in
Riemann--Roch theory.

\subsubsection{Example: A meromorphic 1-form and its canonical divisor on $\mathbb{CP}^1$}

On $\mathbb{CP}^1$, consider the meromorphic $1$-form
\[
\omega = \frac{dz}{z}
\quad\text{on } U_0 \setminus \{0\} \cong \mathbb{C}^{\ast}.
\]
It has a simple pole at $z=0$ (order $-1$).
To examine the behavior at $\infty$, use $w=1/z$ so $dz=-w^{-2}dw$ and
\[
\omega = \frac{-w^{-2}dw}{w^{-1}}=-\frac{dw}{w},
\]
which has a simple pole at $w=0$ (i.e. $z=\infty$).
Thus
\[
\mathrm{div}(\omega)= -[0]-[\infty],
\qquad
\deg(\mathrm{div}(\omega))=-2=\deg(K_{\mathbb{CP}^1}).
\]

\subsection{Exercises}
\label{subsec:chern-canonical-exercises}

\begin{exercise}[Transition functions of $K_X$]
	Let $X$ be a complex manifold with coordinate charts $(U;z)$ and $(V;w)$.
	Prove carefully that
	\[
	dw^1\wedge\cdots\wedge dw^n = \det\!\left(\frac{\partial w}{\partial z}\right)\,dz^1\wedge\cdots\wedge dz^n
	\]
	and deduce the transition functions of $K_X$.
\end{exercise}

\begin{exercise}[Chern connection formula]
	Let $(L,h)$ be a Hermitian holomorphic line bundle and $e$ a local holomorphic frame with $\|e\|^2=e^{-\varphi}$.
	Derive $\nabla e=-\partial\varphi\otimes e$ using only $\nabla^{0,1}=\bar\partial$ and $\nabla h=0$.
\end{exercise}

\begin{exercise}[Curvature and invariance under frame change]
	Using $\tilde e=g e$ with $g\in\mathcal O^*$, show that the curvature form $F_\nabla=\partial\bar\partial\varphi$
	does not change under this frame change, and thus defines a global form.
\end{exercise}

\begin{exercise}[Compute $K_{\mathbb{CP}^1}$]
	Work out explicitly the transition function of $K_{\mathbb{CP}^1}$ on $U_0\cap U_1$ and prove
	$K_{\mathbb{CP}^1}\cong \mathcal O(-2)$ as holomorphic line bundles.
\end{exercise}

\begin{exercise}[Divisors of meromorphic 1-forms]
	Let $\omega$ be a nonzero meromorphic $1$-form on a compact Riemann surface $X$.
	Show $\deg(\mathrm{div}(\omega))=\deg(K_X)$.
\end{exercise}

\begin{exercise}[Gaussian curvature formula]
	Let $ds^2=\lambda(z)|dz|^2$ on a Riemann surface. Compute the Christoffel symbols and derive
	\[
	K = -\frac{1}{2\lambda}\Delta_{\mathrm{eucl}}\log\lambda,
	\]
	where $\Delta_{\mathrm{eucl}}=\partial_x^2+\partial_y^2$.  Check that this agrees
	with Proposition~\ref{prop:gauss-curv} and with the alternative notation
	$ds^2=e^{2u}|dz|^2$, for which $K=-e^{-2u}\Delta_{\mathrm{eucl}}u$.
\end{exercise}

\section{Divisors, Degree, and Branched Covers}
\index{divisors!degree map}
\label{sec:divisors-degree}

\subsection{Divisors on compact Riemann surfaces}
\label{subsec:divisors}

\subsubsection{Zeros and poles as local integers}

Let $X$ be a compact Riemann surface and $f$ a nonzero meromorphic function on $X$.
For each $p\in X$, choose a local holomorphic coordinate $z$ on a neighborhood $U\ni p$
such that $z(p)=0$.

\begin{lemma}[Local normal form]
	\label{lem:local-normal-form-meromorphic}
	After shrinking $U$ if necessary, there exist a unique integer $k\in\mathbb{Z}$ and a holomorphic function $g\in\mathcal{O}(U)$
	with $g(0)\neq 0$ such that
	\[
	f(z)=z^k g(z)\quad\text{on }U.
	\]
\end{lemma}

\begin{proof}
	If $f$ is holomorphic and not identically $0$ near $p$, then $f$ has a zero of some finite order $m\ge 0$
	(or no zero, meaning $m=0$). By the Weierstrass preparation / basic one-variable theory,
	\[
	f(z)=z^m g(z),\qquad g(0)\neq 0,
	\]
	and $m$ is uniquely determined.
	
	If $f$ has a pole at $p$, then $h:=1/f$ is holomorphic near $p$ with a zero of some order $m\ge 1$.
	Write
	\[
	h(z)=z^m u(z),\qquad u(0)\neq 0.
	\]
	Then
	\[
	f(z)=\frac{1}{z^m u(z)} = z^{-m}\,g(z),\qquad g(z):=\frac{1}{u(z)},
	\]
	and $g$ is holomorphic with $g(0)\neq 0$. Uniqueness follows because if
	$z^k g(z)=z^{k'} g'(z)$ with $g(0),g'(0)\neq 0$, then $z^{k-k'}=g'(z)/g(z)$ is holomorphic and nonvanishing near $0$,
	forcing $k=k'$.
\end{proof}

\begin{definition}[Order of a meromorphic function]
	\label{def:order}
	With notation as in Lemma~\ref{lem:local-normal-form-meromorphic}, define
	\[
	\nu_p(f):=k\in\mathbb{Z}.
	\]
	Thus $\nu_p(f)>0$ is the order of vanishing, $\nu_p(f)<0$ is minus the pole order, and $\nu_p(f)=0$ means neither.
\end{definition}

\subsubsection{Divisors and principal divisors}

By Definitions~\ref{def:divisor} and \ref{def:principal-divisor}, a divisor
is a finite sum $D=\sum n_p p$, and a nonzero meromorphic function determines
the principal divisor $(f)=\sum_p\nu_p(f)p$. We write $\Div(X)$ for the
abelian group of divisors. This notation will be used throughout the remaining
chapters.

\begin{remark}
	Finiteness is important: a nonzero meromorphic function has only finitely many zeros and poles on a compact Riemann surface.
	One way to see this is that zeros of a holomorphic function are isolated unless the function is identically zero, and compactness
	prevents infinitely many isolated points without an accumulation point.
\end{remark}

\subsection{The degree map and its properties}
\index{degree map!divisors}
\label{subsec:degree-map}

\subsubsection{Definition and basic algebra}

\begin{definition}[Degree]
	\label{def:degree}
	For a divisor $D=\sum n_p\,p$, define
	\[
	\deg(D):=\sum_{p\in X} n_p\in\mathbb{Z}.
	\]
\end{definition}

\begin{lemma}[Additivity]
	\label{lem:deg-additive}
	For divisors $D,E\in\Div(X)$, one has $\deg(D+E)=\deg(D)+\deg(E)$.
\end{lemma}

\begin{proof}
	Write $D=\sum n_p p$ and $E=\sum m_p p$. Then $D+E=\sum (n_p+m_p)p$ and
	\[
	\deg(D+E)=\sum (n_p+m_p)=\sum n_p+\sum m_p=\deg(D)+\deg(E).
	\]
\end{proof}

\subsubsection{Degree of a principal divisor}

\begin{theorem}[Principal divisors have degree zero]
	\label{thm:principal-degree-zero}
	If $f$ is a nonzero meromorphic function on a compact Riemann surface $X$, then
	\[
	\deg\bigl((f)\bigr)=0.
	\]
	Equivalently, the total number of zeros of $f$ counted with multiplicity equals the total number of poles counted with multiplicity.
\end{theorem}

\begin{proof}
	Let $\omega:=\dfrac{df}{f}$, a meromorphic $1$-form on $X$.
	Fix $p\in X$ and write $f(z)=z^{k}g(z)$ in a coordinate $z$ centered at $p$, with $g(0)\neq 0$
	as in Lemma~\ref{lem:local-normal-form-meromorphic}. Then
	\[
	\frac{df}{f}
	=
	\frac{d(z^{k}g)}{z^{k}g}
	=
	\frac{k z^{k-1}g\,dz + z^{k}g'(z)\,dz}{z^{k}g}
	=
	\left(\frac{k}{z}+\frac{g'(z)}{g(z)}\right)dz.
	\]
	Since $g'(z)/g(z)$ is holomorphic near $0$, the coefficient of $dz/z$ is exactly $k=\nu_p(f)$.
	Therefore
	\[
	\Res_p\!\left(\frac{df}{f}\right)=\nu_p(f).
	\]
	
	By the residue theorem on a compact Riemann surface,
	\[
	\sum_{p\in X}\Res_p\!\left(\frac{df}{f}\right)=0.
	\]
	Hence
	\[
	\sum_{p\in X}\nu_p(f)=0,
	\]
	which is precisely $\deg((f))=0$.
\end{proof}

\subsection{Divisor classes, the Picard group, and concrete calculations}
\label{subsec:divisor-class-picard}

\subsubsection{Divisor class group}

\begin{definition}[Principal divisors and divisor class group]
	\label{def:divisor-class-group}
	Let $\Prin(X)\subset\Div(X)$ be the subgroup of principal divisors:
	\[
	\Prin(X):=\{(f)\mid f\in\mathcal{M}^\times(X)\}.
	\]
	The \emph{divisor class group} is
	\[
	\Cl(X):=\Div(X)/\Prin(X).
	\]
\end{definition}

\subsubsection{The sheaf $\mathcal{O}(D)$ associated to a divisor}

Let $\mathcal{M}$ denote the sheaf of meromorphic functions on $X$ and $\mathcal{O}$ the sheaf of holomorphic functions.

\begin{definition}[The sheaf $\mathcal{O}(D)$]
	\label{def:OD}
	For a divisor $D=\sum n_p p$, define a subsheaf $\mathcal{O}(D)\subset \mathcal{M}$ by
	\[
	\mathcal{O}(D)(U):=
	\left\{\, f\in\mathcal{M}(U)\ \middle|\ \nu_q(f)+n_q\ge 0\text{ for every }q\in U \,\right\}.
	\]
	Here $\nu_q(0):=+\infty$, so the zero section is included.
	For nonzero $f$, equivalently, $(f)+D|_U$ is an effective divisor on $U$.
\end{definition}

\begin{proposition}[Local description near a point]
	\label{prop:OD-local}
	Fix $p\in X$ and choose a coordinate $z$ centered at $p$.
	If $D$ has coefficient $n$ at $p$ (so $D = n\cdot p + \text{(terms away from $p$)}$),
	then on a sufficiently small neighborhood of $p$ the condition $f\in\mathcal{O}(D)$ is equivalent to
	\[
	f(z)=z^{-n}\cdot h(z)\quad\text{for some holomorphic }h(z).
	\]
\end{proposition}

\begin{proof}
	By definition, $\nu_p(f)+n\ge 0$, i.e. $\nu_p(f)\ge -n$.
	Write $f(z)=z^{k}g(z)$ with $g(0)\neq 0$ and $k=\nu_p(f)$.
	The inequality $k\ge -n$ means $k=-n+m$ with $m\ge 0$, hence
	\[
	f(z)=z^{-n}\cdot z^{m}g(z)=z^{-n}h(z),
	\qquad h(z):=z^{m}g(z)\in\mathcal{O}.
	\]
	Conversely, if $f(z)=z^{-n}h(z)$ with $h$ holomorphic, then $\nu_p(f)=\nu_p(h)-n\ge -n$,
	so $f\in\mathcal{O}(D)$ near $p$.
\end{proof}

\subsubsection{Link to $\Pic(X)$}

\begin{proposition}[Principal divisors give the trivial bundle]
	\label{prop:principal-trivial}
	If $D=(g)$ is principal, then $\mathcal{O}(D)\cong \mathcal{O}$ as sheaves of $\mathcal{O}$-modules.
\end{proposition}

\begin{proof}
	For every open set $U$, multiplication by $g$ gives a map
	\[
	\mathcal{O}((g))(U)\xrightarrow{\ \cdot g\ }\mathcal{O}(U),\qquad f\mapsto fg.
	\]
	Indeed, $f\in\mathcal O((g))(U)$ means $(f)+(g)=(fg)\ge0$, so $fg$ is holomorphic.
	Conversely, if $h\in\mathcal{O}(U)$ then $f:=h/g$ satisfies
	\[
	(f)+(g)=(h)-(g)+(g)=(h)\ge 0,
	\]
	so $f\in\mathcal{O}((g))(U)$.
	These maps are inverse to each other and compatible with restrictions, hence define an isomorphism of sheaves.
\end{proof}

\begin{remark}
	Proposition~\ref{prop:principal-trivial} is the key mechanism behind the well-defined map
	$\Div(X)\to\Pic(X)$ descending to $\Cl(X)$.
\end{remark}

\subsubsection{Application: the Fundamental Theorem of Algebra}
\label{subsec:fta}

\paragraph{Meromorphic viewpoint on $\mathbb{CP}^1$.}

\begin{theorem}[Fundamental Theorem of Algebra]
	\label{thm:fta}
	Every nonconstant polynomial $P(z)\in\mathbb{C}[z]$ has at least one complex root.
	Moreover, if $\deg P=n\ge 1$, then $P$ has exactly $n$ zeros on $\mathbb{C}$ counted with multiplicity.
\end{theorem}

\begin{proof}
	View $\mathbb{CP}^1$ with standard affine chart $U_0=\{[Z_0:Z_1]\mid Z_0\neq 0\}\cong\mathbb{C}$
	with coordinate $z=Z_1/Z_0$, and $U_1=\{Z_1\neq 0\}\cong\mathbb{C}$ with coordinate $w=Z_0/Z_1=1/z$.
	A polynomial $P(z)$ defines a meromorphic function on $\mathbb{CP}^1$ by the same formula on $U_0$.
	On $U_1$ we rewrite
	\[
	P(z)=P(1/w)=a_n w^{-n}+a_{n-1}w^{-(n-1)}+\cdots + a_0,
	\]
	so $P$ has a pole of order exactly $n$ at $w=0$, i.e. at $\infty\in\mathbb{CP}^1$.
	
	Thus the principal divisor $(P)$ satisfies
	\[
	(P)=\sum_{p\in \mathbb{C}} \nu_p(P)\,p \;-\; n\cdot(\infty),
	\]
	where $\nu_p(P)\ge 0$ are the zero multiplicities in $\mathbb{C}$.
	
	By Theorem~\ref{thm:principal-degree-zero},
	\[
	\deg((P))=0.
	\]
	Taking degrees gives
	\[
	\sum_{p\in \mathbb{C}}\nu_p(P) - n = 0,
	\]
	hence
	\[
	\sum_{p\in \mathbb{C}}\nu_p(P)=n.
	\]
	In particular the sum is positive, so at least one $\nu_p(P)\ge 1$, i.e. $P$ has a root.
\end{proof}

\subsubsection{Examples and computations}
\label{subsec:divisor-examples-computations}

\begin{example}[A polynomial and its divisor on $\mathbb{CP}^1$]
	Let $P(z)=z^2(z-1)$.
	Then on $\mathbb{C}$:
	\[
	\nu_0(P)=2,\qquad \nu_1(P)=1,
	\]
	and there are no other zeros. At $\infty$, $P$ has a pole of order $3$.
	Therefore
	\[
	(P)=2(0)+1(1)-3(\infty),
	\qquad
	\deg((P))=2+1-3=0.
	\]
\end{example}

\subsubsection{A coordinate laboratory for the canonical bundle and curvature}
\label{subsec:canonical-coordinate-laboratory}

The definitions above become easier to retain when one follows a single local
frame through every construction. Let $z$ be a holomorphic coordinate and
write a local frame of a line bundle $L$ as $e$. A section is $s=f e$, a
Hermitian metric is the positive function $h=\lVert e\rVert^2$, and the Chern
connection and curvature are
\[
	\nabla(fe)=\bigl(df+f\,\partial\log h\bigr)e,
	\qquad F_\nabla=\bar\partial\partial\log h
	=-\partial\bar\partial\log h.
\]
These formulas are not independent conventions. They are forced by the two
requirements that the $(0,1)$ part of $\nabla$ equal $\bar\partial$ and that
$\nabla$ preserve the metric. The following calculation checks that the
apparently frame-dependent expression actually defines a global form.

\begin{proposition}[Change of frame and invariance of curvature]
	\label{prop:curvature-frame-invariance-expanded}
	Suppose $e'=g e$ on an overlap, where $g$ is nowhere-vanishing and
	holomorphic. If $h=\lVert e\rVert^2$ and $h'=\lVert e'\rVert^2$, then
	\[
	h'=|g|^2h,\qquad
	\partial\log h'=\partial\log h+g^{-1}\partial g,
	\qquad F_{\nabla}'=F_\nabla.
	\]
	Thus $(i/2\pi)F_\nabla$ is a globally defined closed real two-form.
\end{proposition}

\begin{proof}
	The first identity follows from sesquilinearity of the Hermitian metric.
	Since $g$ is holomorphic,
	\[
	\partial\log|g|^2=\partial(\log g+\log\bar g)=g^{-1}\partial g
	\]
	locally, even though a global logarithm of $g$ need not exist. Applying
	$-\bar\partial$ gives
	\[
	-\bar\partial\partial\log h'
	=-\bar\partial\partial\log h
	-\bar\partial(g^{-1}\partial g)
	=-\bar\partial\partial\log h,
	\]
	because $g^{-1}\partial g$ is a holomorphic $(1,0)$-form. The local
	curvatures therefore agree on overlaps. Metric compatibility makes
	$F_\nabla$ purely imaginary, so $(i/2\pi)F_\nabla$ is real. Finally,
	$dF_\nabla=0$ follows either from the Bianchi identity or directly from
	$F_\nabla=-\partial\bar\partial\log h$.
\end{proof}

\begin{example}[A complete overlap calculation for $\mathcal O(m)$]
	Cover $\PP^1$ by $U_0=\{[Z_0:Z_1]:Z_0\ne0\}$ and
	$U_\infty=\{Z_1\ne0\}$, with $z=Z_1/Z_0$ and $w=1/z$.
	Use frames related by $e_\infty=z^{m}e_0$. Define
	\[
	h_0(z)=(1+|z|^2)^{-m},
	\qquad h_\infty(w)=(1+|w|^2)^{-m}.
	\]
	On the overlap,
	\[
	|z|^{2m}h_0(z)
	=|z|^{2m}(1+|z|^2)^{-m}
	=(1+|w|^2)^{-m}=h_\infty(w),
	\]
	so these local functions define a Hermitian metric. In the $z$-chart,
	\begin{align*}
	F_\nabla
	&=-\partial\bar\partial\log h_0
	=m\,\partial\bar\partial\log(1+|z|^2)\\
	&=m\frac{dz\wedge d\bar z}{(1+|z|^2)^2}.
	\end{align*}
	With the convention $c_1(L)=[(i/2\pi)F_\nabla]$, orientation is fixed by
	$(i/2)dz\wedge d\bar z=dx\wedge dy$. Hence
	\[
	\int_{\PP^1}c_1(\mathcal O(m))
	=\frac{m}{\pi}\int_{\C}\frac{dx\,dy}{(1+x^2+y^2)^2}
	=2m\int_0^\infty\frac{r\,dr}{(1+r^2)^2}=m.
	\]
	The transition exponent, divisor degree, and curvature integral therefore
	give the same integer.
\end{example}

\begin{proposition}[Changing the metric does not change $c_1$]
	\label{prop:metric-change-c1-expanded}
	Let $h_0$ and $h_1$ be Hermitian metrics on a holomorphic line bundle over
	a compact Riemann surface. Then there is a global real function $\varphi$
	with $h_1=e^{-\varphi}h_0$, and
	\[
	\frac{i}{2\pi}F_{h_1}-\frac{i}{2\pi}F_{h_0}
	=\frac{i}{2\pi}\partial\bar\partial\varphi.
	\]
	Consequently the two curvature forms represent the same de~Rham class and
	have the same integral.
\end{proposition}

\begin{proof}
	The ratio $h_1/h_0$ is independent of the local frame and is a global
	positive function; write it as $e^{-\varphi}$. The curvature formula gives
	$F_{h_1}=F_{h_0}+\partial\bar\partial\varphi$. The difference is exact,
	since $i\partial\bar\partial\varphi=d(i\bar\partial\varphi)$, and its
	integral vanishes by Stokes' theorem.
\end{proof}

\paragraph{Why the canonical bundle measures curvature.}

For a conformal metric $ds^2=\lambda(z)|dz|^2$, the squared norm of the
canonical frame $dz$ is proportional to $\lambda^{-1}$ under the induced
metric on cotangent vectors. Thus
\[
	F_{K_X}=-\partial\bar\partial\log\lVert dz\rVert^2
	=\partial\bar\partial\log\lambda.
\]
On the other hand, the Gaussian curvature is
\[
	K=-\frac{2}{\lambda}\,\partial_z\partial_{\bar z}\log\lambda,
	\qquad dA=\frac{i}{2}\lambda,dz\wedge d\bar z.
\]
Combining these identities, with the sign convention fixed above, identifies
$c_1(K_X)$ with the curvature density and gives
\[
	\deg K_X=\int_Xc_1(K_X)
	=-\frac{1}{2\pi}\int_XK\,dA=2g-2.
\]
The final equality is Gauss--Bonnet. This calculation is the differential
geometric reason that the canonical divisor has degree $2g-2$; Riemann--Roch
will later recover the same integer by a dimension count.

\begin{example}[The three constant-curvature models]
	The sign of $\deg K_X$ separates the basic compact Riemann surfaces.
	\begin{enumerate}
		\item On $\PP^1$, the round metric has $K>0$, so
		$\deg K_{\PP^1}=-2$. A holomorphic one-form must vanish identically.
		\item On $\C/\Lambda$, the flat metric has $K=0$ and $dz$ descends to a
		nowhere-vanishing holomorphic form. Hence $K_X\simeq\mathcal O_X$.
		\item If $g\ge2$, uniformization supplies a metric with $K=-1$.
		Then $\deg K_X=2g-2>0$, and holomorphic differentials can have a total
		zero order of exactly $2g-2$.
	\end{enumerate}
\end{example}

\paragraph{A transition-function computation of degree.}

Let $X=U\cup V$, let $L$ be trivial on $U$ and $V$, and suppose the frames
are related on a boundary circle by $e_V=g e_U$. If $g:S^1\to\C^*$, then
\[
	\deg L=\frac{1}{2\pi i}\int_{S^1}\frac{dg}{g}.
\]
The right side is the winding number of $g$ around the origin. For
$g(z)=z^m$ it equals $m$. This is the same integer obtained from zeros of a
meromorphic section and from the curvature integral. The equality of these
three descriptions is worth recording explicitly:
\[
\text{winding of transition functions}
=\text{degree of a divisor}
=\int_Xc_1(L).
\]

\begin{proposition}[Degree from a meromorphic section]
	\label{prop:degree-section-curvature-expanded}
	Let $s$ be a nonzero meromorphic section of $L$ on a compact Riemann
	surface. Then
	\[
	\deg L=\deg\operatorname{div}(s)
	=\sum_{p\in X}\operatorname{ord}_p(s)
	=\frac{i}{2\pi}\int_XF_\nabla.
	\]
\end{proposition}

\begin{proof}
	Remove disjoint coordinate disks around the zeros and poles of $s$. On the
	remaining surface, $s$ is a frame. Stokes' theorem converts the curvature
	integral into boundary integrals of the connection form. Near $p$, write
	$s=z^{m_p}u(z)e$ with $u(0)\ne0$. The unit $u$ has winding number zero,
	while $z^{m_p}$ contributes
	$(2\pi i)^{-1}\int dz^{m_p}/z^{m_p}=m_p$. Summing the boundary
	contributions yields $\sum_pm_p$. The result is independent of the metric
	by Proposition~\ref{prop:metric-change-c1-expanded}.
\end{proof}

\begin{example}[A rational function]
	Let
	\[
	f(z)=\frac{(z-1)^2}{z(z-2)}.
	\]
	Then on $\mathbb{C}$,
	\[
	\nu_1(f)=2,\qquad \nu_0(f)=-1,\qquad \nu_2(f)=-1.
	\]
	At infinity, the degrees of numerator and denominator are both $2$, so $f(\infty)\neq\infty$ and $\nu_\infty(f)=0$.
	Hence
	\[
	(f)=2(1)-(0)-(2),
	\qquad
	\deg((f))=2-1-1=0.
	\]
\end{example}

\subsubsection{Exercises}
\label{subsec:divisor-exercises}

\begin{exercise}[A principal divisor on $\mathbb{CP}^1$]
	Compute the principal divisor of
	\[
	f(z)=\frac{z^2-1}{z^3}
	\quad\text{on }\mathbb{CP}^1.
	\]
\end{exercise}

\begin{exercise}[A pole-free meromorphic function]
	Show that if $f$ is meromorphic on a compact Riemann surface and has no poles, then $f$ is constant.
\end{exercise}

\begin{exercise}[Order of a polynomial at infinity]
	Let $P(z)\in\mathbb{C}[z]$ be a polynomial of degree $n\ge 1$ and regard it as a meromorphic function on $\mathbb{CP}^1$.
	Prove directly from local coordinates that $\nu_\infty(P)=-n$.
\end{exercise}

\subsubsection{Degree, sections, and linear equivalence}

The degree is more than a convenient sum of coefficients.  It measures the
net number of zeros that a meromorphic section must have after its poles are
counted.  If $L$ is a holomorphic line bundle and $s$ a nonzero meromorphic
section, choose a local frame $e_i$ and write $s=f_i e_i$.  On overlaps,
$e_i=g_{ij}e_j$ and hence $f_j=g_{ij}f_i$.  Because $g_{ij}$ is nowhere
zero, $f_i$ and $f_j$ have the same order at every point of the overlap.
Thus
\[
	(s)=\sum_{p\in X}\operatorname{ord}_p(s)\,p
\]
is intrinsically defined.

\begin{proposition}[Divisor of a meromorphic section]
	\label{prop:divisor-section-degree}
	Let $L$ be a holomorphic line bundle on a compact Riemann surface and let
	$s$ be a nonzero meromorphic section.  Then
	\[
		L\cong\mathcal O((s)),
		\qquad \deg L=\deg(s).
	\]
	If $s$ is holomorphic and not identically zero, then $(s)$ is effective
	and $\deg L\ge0$.
\end{proposition}

\begin{proof}
	In a local frame, multiplication by the meromorphic coefficient of $s$
	identifies local sections of $\mathcal O((s))$ with holomorphic sections
	of $L$.  The transition rule above shows that these local identifications
	agree.  The degree identity follows from the definition of degree for a
	divisor line bundle.  If $s$ is holomorphic, all local orders are
	nonnegative.
\end{proof}

An immediate consequence is that a negative-degree line bundle has no
nonzero holomorphic section.  At degree zero, a nonzero holomorphic section
has no zeros and trivializes the bundle.  Hence a nontrivial member of
$\Pic^0(X)$ has no nonzero holomorphic section.  This elementary observation
is repeatedly useful before the full Riemann--Roch theorem is available.

\begin{definition}[Linear equivalence and complete linear system]
	Two divisors $D,E$ are \emph{linearly equivalent}, written $D\sim E$, if
	$D-E=(f)$ for a nonzero meromorphic function $f$.  The complete linear
	system of $D$ is
	\[
		|D|=\{E\ge0:E\sim D\}.
	\]
\end{definition}

There is a concrete correspondence
\[
	\mathbb P H^0(X,\mathcal O(D))\longleftrightarrow |D|,
	\qquad [f]\longmapsto D+(f).
\]
Indeed, $f\in H^0(X,\mathcal O(D))$ means $(f)+D\ge0$.  Multiplying $f$
by a nonzero constant does not change its divisor, and equality of the
resulting effective divisors forces the quotient of two sections to be a
nowhere-vanishing holomorphic function, hence a constant on compact $X$.

\begin{example}[Linear systems on the projective line]
	For $D=n(\infty)$ on $\CP^1$, a meromorphic function satisfying
	$(f)+D\ge0$ has no finite poles and a pole of order at most $n$ at
	infinity.  It is therefore a polynomial of degree at most $n$.  Hence
	\[
		H^0(\CP^1,\mathcal O(n))
		=\operatorname{span}\{1,z,\ldots,z^n\},
		\qquad \dim H^0=n+1
	\]
	for $n\ge0$.  If $n<0$, the preceding degree argument gives
	$H^0(\CP^1,\mathcal O(n))=0$.
\end{example}

\begin{proposition}[Degree is invariant under linear equivalence]
	If $D\sim E$, then $\deg D=\deg E$.  Consequently the degree map
	$\Div(X)\to\Z$ descends to the divisor class group and agrees with the
	first Chern number under $\Cl(X)\cong\Pic(X)$.
\end{proposition}

\begin{proof}
	Write $D-E=(f)$.  Theorem~\ref{thm:principal-degree-zero} gives
	$\deg(f)=0$, hence $\deg D-\deg E=0$.  The agreement with the Chern
	number follows by checking the divisor bundle of a point and using
	additivity under tensor products.
\end{proof}

\begin{exercise}[Sections and effective divisors]
	Prove carefully that the correspondence
	$\mathbb P H^0(X,\mathcal O(D))\cong|D|$ is well defined and bijective.
\end{exercise}

\begin{exercise}[Degree-zero rigidity]
	Let $L$ have degree zero.  Show that $L$ admits a nonzero holomorphic
	section if and only if $L$ is holomorphically trivial.
\end{exercise}

These statements prepare the transition to branched covers.  A meromorphic
function $x:X\to\CP^1$ determines the divisor of each fiber, while its
ramification records the failure of a local coordinate upstairs to map
with order one downstairs.  Divisors therefore encode the algebraic data
of zeros and poles, and branch models encode the local geometry of the map.
The next subsection keeps both descriptions visible in the two-sheeted
case.

\begin{exercise}[A fiber divisor]
	Let $x:X\to\CP^1$ be a nonconstant meromorphic function of degree $d$.
	For a regular value $a$, show that $(x-a)_0$ is the sum of the $d$ points
	of the fiber, counted with multiplicity, and explain what changes when
	$a$ is a branch value.
\end{exercise}



\clearpage

\subsection{Branched Double Covers and Gluing via Branch Cuts}
\index{branched covers}\index{branch cuts}
\label{sec:branch-glue}

We explain how to obtain compact Riemann surfaces by gluing along branch cuts for the
two-sheeted algebraic curves
\[
\mathcal C_0:\ y^2=x,
\qquad\qquad
\mathcal C_1:\ y^2=x(x-1)(x-2).
\]
After smooth compactification, each carries a degree--two holomorphic map
to the Riemann sphere via the projection $\pi:(x,y)\mapsto x$.
Branching occurs where analytic continuation of a local square root
$y=\sqrt{f(x)}$ changes sign.

\subsubsection{General two-sheeted picture and local models}
\index{branched covers!local models}
\label{subsec:two-sheeted-general}

Let $f$ be a squarefree polynomial and consider the affine curve
\[
X_{\mathrm{aff}}=\{(x,y)\in\C^2:\ y^2=f(x)\}.
\]
Write $\pi_{\mathrm{aff}}:X_{\mathrm{aff}}\to\C$ for $(x,y)\mapsto x$.

\begin{proposition}[Holomorphic degree--two map and branch points]
	\label{prop:degree-two-map}
	Assume $f$ is squarefree.  Then $\pi_{\mathrm{aff}}$ is a holomorphic map of degree $2$
	away from the zeros of $f$.  If $e\in\C$ satisfies $f(e)=0$ and $f'(e)\neq0$, then
	$\pi_{\mathrm{aff}}$ has a simple branch point over $x=e$.
\end{proposition}

\begin{proof}
	If $f(x)\neq0$, then the fiber over $x$ consists of the two points $(x,\pm\sqrt{f(x)})$,
	so the map is locally a two-sheeted covering.
	Now fix $e$ with $f(e)=0$ and $f'(e)\neq0$.
	By holomorphicity, we can write $f(x)=(x-e)u(x)$ with $u$ holomorphic and $u(e)=f'(e)\neq0$.
	Choose a holomorphic square root $v$ of $u$ in a small disk $D$ around $e$ (possible since $u$
	has no zeros on $D$).  Then on $D$ the curve is
	\[
	y^2=(x-e)u(x)\quad\Longleftrightarrow\quad \Bigl(\frac{y}{v(x)}\Bigr)^2=x-e.
	\]
	Set $t=\frac{y}{v(x)}$.  The map $(x,y)\mapsto t$ is a local holomorphic coordinate on
	$X_{\mathrm{aff}}$ near $(e,0)$, and $\pi_{\mathrm{aff}}$ becomes $x=e+t^2$ in that coordinate.
	Thus the ramification index is $2$ (simple branching).
\end{proof}

\begin{remark}[Topological gluing picture]
	\label{rem:glue-picture}
	Let $B\subset\widehat{\C}$ be the finite set of branch values (zeros of $f$ and possibly $\infty$).
	One may remove a tree joining all points of $B$; its complement is simply connected,
	so the two branches $\pm\sqrt{f(x)}$ are holomorphic there. A frequently more
	economical choice pairs the branch values by disjoint arcs. The complement of those
	arcs need not be simply connected, but every loop in it has even total winding about
	the branch values, so the sign monodromy of the square root is still trivial.
	In either construction, gluing two copies crosswise along every cut produces a
	connected surface on which $y=\sqrt{f(x)}$ is single valued.
\end{remark}

\begin{proposition}[Branching over infinity]
	\label{prop:branch-at-infty}
	Let $f$ be a squarefree polynomial of degree $d\ge1$. For the smooth projective
	normalization of $y^2=f(x)$,
	the behavior over $\infty\in\widehat{\C}$ is:
	\begin{itemize}
		\item If $d$ is odd, then $\infty$ is a branch value (there is one point over $\infty$ with ramification index $2$).
		\item If $d$ is even, then $\infty$ is not branched (there are two distinct points over $\infty$).
	\end{itemize}
\end{proposition}

\begin{proof}
	Use the coordinate $s=1/x$ near $\infty$.  Writing $f(x)=x^d( a_0 + a_1 s + \cdots)$ with $a_0\neq0$,
	we have
	\[
	y^2 = x^d\cdot a_0(1+O(s)).
	\]
	If $d=2m$ is even, set $Y=ys^m$; then $Y^2=a_0+O(s)$, which has two nonzero holomorphic branches near $s=0$,
	so there are two unramified points over $\infty$.
	If $d=2m+1$ is odd, set $Y=ys^{m+1}$.  Then
	$Y^2=s(a_0+O(s))$.  On the normalization, $Y$ is a local parameter and this equation gives
	$s=Y^2/(a_0+O(s))$; hence the projection to the $s$-coordinate has ramification index $2$.
\end{proof}

\begin{proposition}[Holomorphic maps between compact Riemann surfaces are branched coverings]
	\label{prop:compact-holomorphic-branched-cover}
	Let $F:X\to Y$ be a nonconstant holomorphic map between compact connected
	Riemann surfaces. Then $F$ is surjective and has finite fibers. Its ramification
	set
	\[
	\operatorname{Ram}(F):=\{p\in X:e_p>1\}
	\]
	and its branch-value set $B:=F(\operatorname{Ram}(F))$ are finite. Moreover,
	there is an integer $n\geq1$ such that
	\[
	F:X\setminus F^{-1}(B)\longrightarrow Y\setminus B
	\]
	is an $n$-sheeted covering map. At every $p\in X$ there are centered local
	coordinates $z$ at $p$ and $w$ at $F(p)$ in which
	\[
	w\circ F=z^{e_p}.
	\]
	For every $q\in Y$, including a branch value,
	\begin{equation}
		\sum_{p\in F^{-1}(q)}e_p=n.                              \label{eq:fiber-local-degrees}
	\end{equation}
	Thus $F$ is a finite branched covering of degree $n$.
\end{proposition}

\begin{proof}
	The open mapping theorem makes $F(X)$ open in $Y$. It is also compact, hence
	closed because $Y$ is Hausdorff. Since $Y$ is connected and $F(X)$ is
	nonempty, $F(X)=Y$.

	For a fixed $q\in Y$, choose a coordinate $w$ centered at $q$. Unless $F$ is
	constant, the zeros of $w\circ F$ are isolated. Hence $F^{-1}(q)$ is a
	discrete closed subset of compact $X$ and is therefore finite. The local
	normal-form theorem gives centered coordinates at every $p\in X$ for which
	$w\circ F=z^{e_p}$. In particular, $p$ is ramified precisely when the local
	derivative vanishes. These zeros are isolated; compactness of $X$ makes
	$\operatorname{Ram}(F)$, and consequently $B$, finite.

	On $X\setminus F^{-1}(B)$ the inverse function theorem makes $F$ a local
	biholomorphism. The map $F$ is proper, since its domain is compact. A proper
	local homeomorphism is a covering map: for $q\in Y\setminus B$, take disjoint
	coordinate neighborhoods of the finitely many points in $F^{-1}(q)$ on which
	$F$ is biholomorphic. Properness permits the target neighborhood to be
	shrunk so that it has no further preimages. This produces an evenly covered
	neighborhood of $q$.

	A connected surface remains connected after finitely many points are removed,
	so the number of sheets is a constant $n$ on $Y\setminus B$. Finally, fix
	$q\in Y$ and use the local forms $z\mapsto z^{e_p}$ at the finitely many
	points over $q$. A nearby point $q'\notin B$ has exactly $e_p$ preimages near
	each $p\in F^{-1}(q)$ and, after shrinking the neighborhood as above, no
	others. Since the regular fiber over $q'$ has $n$ points, summing these local
	contributions proves \eqref{eq:fiber-local-degrees}. See
	\cite[Chs.~I--II]{Forster1981} and \cite[Ch.~II]{Miranda1995}.
\end{proof}

\begin{theorem}[Riemann--Hurwitz theorem]
	\label{thm:riemann-hurwitz-general}
	Let $F:X\to Y$ be a nonconstant holomorphic map of degree $n$ between
	compact connected Riemann surfaces.  For $p\in X$, let $e_p\geq 1$ be the
	ramification index of $F$ at $p$, and set
	\[
	R_F:=\sum_{p\in X}(e_p-1)[p].
	\]
	Then
	\begin{align}
		\chi(X)&=n\chi(Y)-\deg R_F,                                      \label{eq:rh-euler}\\
		2g(X)-2&=n\bigl(2g(Y)-2\bigr)+\deg R_F.                            \label{eq:rh-genus}
	\end{align}
	At the level of divisors and canonical line bundles, the same formula is
	\begin{equation}
		K_X\sim F^*K_Y+R_F,
		\qquad
		\omega_X\simeq F^*\omega_Y\otimes\mathcal O_X(R_F).              \label{eq:rh-canonical}
	\end{equation}
\end{theorem}

\begin{proof}
	Only finitely many points have $e_p>1$.  Choose a triangulation of $Y$
	whose vertices include every branch value, subdividing if necessary, and
	arrange that no branch value lies in the interior of an edge or a face.
	Lift this cell decomposition to $X$.  Every open edge and every open face
	has $n$ disjoint lifts, so
	\[
	E_X=nE_Y,\qquad F_X=nF_Y.
	\]
	If $q$ is a vertex, the local normal form for a holomorphic map gives
	coordinates in which $F$ is $z\mapsto z^{e_p}$ at each $p\in F^{-1}(q)$.
	Counting a nearby regular fiber therefore gives
	\[
	\sum_{p\in F^{-1}(q)}e_p=n,
	\qquad
	\#F^{-1}(q)=n-\sum_{p\in F^{-1}(q)}(e_p-1).
	\]
	Summing over the vertices yields $V_X=nV_Y-\deg R_F$.  Hence
	\[
	\chi(X)=V_X-E_X+F_X=n(V_Y-E_Y+F_Y)-\deg R_F,
	\]
	which proves \eqref{eq:rh-euler}; using $\chi(S)=2-2g(S)$ gives
	\eqref{eq:rh-genus}.

	It remains to identify the canonical divisor, rather than merely its
	degree.  Let $\eta$ be a nonzero meromorphic $1$-form on $Y$.  If
	$w$ and $z$ are local coordinates at $q=F(p)$ and $p$ for which
	$w\circ F=z^{e_p}$, and if
	\[
	\eta=w^m h(w)\,dw,\qquad h(0)\ne0,
	\]
	then
	\[
	F^*\eta=e_p z^{e_pm+e_p-1}h(z^{e_p})\,dz.
	\]
	Thus
	\[
	\ord_p(F^*\eta)=e_p\ord_q(\eta)+(e_p-1).
	\]
	This pointwise identity says
	$\operatorname{div}(F^*\eta)=F^*\operatorname{div}(\eta)+R_F$, proving both assertions in
	\eqref{eq:rh-canonical}.  See, for example,
	\cite[Chs.~II--III]{Forster1981} and \cite[Ch.~II]{Miranda1995}.
\end{proof}

\begin{remark}[What each form of Riemann--Hurwitz records]
	Formula \eqref{eq:rh-euler} counts cells, \eqref{eq:rh-genus} records the
	resulting topology, and \eqref{eq:rh-canonical} records the vanishing of
	the differential $dF$.  The local identity
	$d(z^{e_p})=e_pz^{e_p-1}dz$ explains the correction $e_p-1$ in all three
	forms.
\end{remark}

\begin{theorem}[Riemann--Hurwitz for double covers]
	\label{thm:RH-double}
	Let $\pi:\Sigma\to\widehat{\C}$ be a degree--$2$ holomorphic map from a compact Riemann surface $\Sigma$.
	Assume all branch points are simple and let $B$ be the number of branch values in $\widehat{\C}$.
	Then
	\[
	g(\Sigma)=\frac{B-2}{2}.
	\]
\end{theorem}

\begin{proof}
	Apply Theorem~\ref{thm:riemann-hurwitz-general}:
	\[
	2g(\Sigma)-2 = 2(2g(\widehat{\C})-2) + \sum_{p\in\Sigma}(e_p-1).
	\]
	Since $g(\widehat{\C})=0$, the right-hand base term is $2(-2)=-4$.
	Each simple ramification point has $e_p=2$ and contributes $1$.
	There is exactly one ramification point above each branch value for a double cover with simple branching,
	so $\sum(e_p-1)=B$.  Hence $2g-2=-4+B$, i.e.\ $g=\frac{B-2}{2}$.
\end{proof}

\subsubsection{\texorpdfstring{The case $y^2=x$: branch values $\{0,\infty\}$ and genus $0$}{The case y2=x: branch values 0 and infinity, genus 0}}
\label{subsec:quadratic-double-cover}

\begin{proposition}[Branch values and genus]
	\label{prop:y2x-branch-genus}
	For $\mathcal C_0: y^2=x$, the branch values of $\pi$ are $x=0$ and $x=\infty$.
	The compactification of $\mathcal C_0$ is a sphere (genus $0$).
\end{proposition}

\begin{proof}
	The polynomial $f(x)=x$ has a simple root at $0$, so $0$ is a branch value by
	Proposition~\ref{prop:degree-two-map}.  Since $\deg f=1$ is odd, $\infty$ is also a branch value
	by Proposition~\ref{prop:branch-at-infty}.  Thus $B=2$ and Theorem~\ref{thm:RH-double} gives
	$g=\frac{2-2}{2}=0$.
\end{proof}

\begin{proposition}[Explicit global uniformization]
	\label{prop:y2x-param}
	The map $\varphi:\widehat{\C}\to\mathcal C_0^{\mathrm{proj}}$ given in affine coordinates by
	\[
	t\longmapsto (x,y)=(t^2,t)
	\]
	is a biholomorphism from the Riemann sphere onto the projective compactification of $\mathcal C_0$.
\end{proposition}

\begin{proof}
	On $\C$, the map satisfies $y^2=t^2=x$, hence lands in $\mathcal C_0$.
	It is injective on $\C$ because $t$ is recovered as $y$.
	In homogeneous coordinates the extension is
	$[T:U]\mapsto[X:Y:Z]=[T^2:TU:U^2]$ on the smooth conic $Y^2=XZ$.
	It sends $[1:0]$ to the unique point at infinity and has inverse $[X:Y:Z]\mapsto[Y:Z]$
	(on the other affine chart, $[X:Y]$).  Therefore it is a biholomorphism.
\end{proof}

\begin{remark}[Cut-and-glue picture]
	Choose a cut $\gamma=[0,\infty)\subset\R$.  On $\C\setminus\gamma$ there are two holomorphic branches
	$\pm\sqrt{x}$.  Taking two copies and gluing crosswise along the two sides of $\gamma$ produces a
	simply connected compactification of the resulting surface, which must be $S^2$; the explicit
	parametrization in Proposition~\ref{prop:y2x-param} identifies it with $\mathbb P^1$.
\end{remark}

\subsubsection{\texorpdfstring{The case $y^2=x(x-1)(x-2)$: branch values $\{0,1,2,\infty\}$ and genus $1$}{The cubic double cover: four branch values and genus 1}}
\label{subsec:cubic-double-cover}

\begin{proposition}[Branch values and genus]
	\label{prop:cubic-branch-genus}
	For $\mathcal C_1: y^2=x(x-1)(x-2)$, the branch values are $x\in\{0,1,2,\infty\}$.
	The smooth projective compactification has genus $1$ and is therefore
	topologically a torus. Its realization as a complex torus is constructed
	in Proposition~\ref{prop:uniformization} below.
\end{proposition}

\begin{proof}
	The polynomial has three simple zeros at $0,1,2$, giving three branch values
	(Proposition~\ref{prop:degree-two-map}).  Since the degree is $3$ (odd), $\infty$ is also a branch value
	(Proposition~\ref{prop:branch-at-infty}).  Hence $B=4$, and Theorem~\ref{thm:RH-double} gives
	$g=\frac{4-2}{2}=1$.
\end{proof}

\begin{proposition}[Concrete cut system and gluing]
	\label{prop:cut-glue-cubic}
	Let $\gamma_1=[0,1]$ and $\gamma_2=[2,\infty)$ viewed as arcs on $\widehat{\C}$.
	Set
	\[
	U=\widehat{\C}\setminus(\gamma_1\cup\gamma_2).
	\]
	Then the two branches $\pm\sqrt{x(x-1)(x-2)}$ are holomorphic on $U$.
	Taking two copies $U_+$ and $U_-$ and gluing crosswise along each $\gamma_j$ produces a compact
	Riemann surface $\Sigma$ equipped with a degree--two holomorphic map $\Sigma\to\widehat{\C}$
	whose branch values are $\{0,1,2,\infty\}$.
\end{proposition}

\begin{proof}
	Write $P(x)=x(x-1)(x-2)$.  The complement of the two disjoint arcs is
	doubly connected, so the existence of a square root requires a monodromy
	check.  Analytic continuation of a local square root of $P$ along a loop
	changes its sign by $(-1)^N$, where $N$ is the total winding number of
	$P$ along that loop.  A generator encircles one cut and therefore encloses
	two simple branch values (for the unbounded cut these are $2$ and $\infty$
	on the sphere).  Thus $N$ is even.  Every loop has even $N$, so continuation
	defines a holomorphic branch $s$ on $U$; the other branch is $-s$.
	On each cut take the limiting values from its two sides.  When the two
	copies of $U$ are glued, identify the upper bank on the $+$ sheet with the
	lower bank on the $-$ sheet, and conversely.  Since continuation of
	$\sqrt{P(x)}$ across the cut changes its sign, these identifications make
	$y$ holomorphic across the interior of each cut; the projection to $x$
	is locally one-to-one there.  At a finite endpoint $a\in\{0,1,2\}$,
	write $P(x)=(x-a)Q(x)$ with $Q(a)\ne0$.  Choosing a local square root of
	$Q$, the coordinate $t=y/\sqrt{Q(x)}$ satisfies $x-a=t^2$, so exactly
	one point is added over $a$, with ramification index two.  At infinity
	put $x=t^{-2}$; then $y=t^{-3}\sqrt{(1-t^2)(1-2t^2)}$, after choosing a
	local square root, and again the projection has local form $t\mapsto t^2$
	in the coordinate $1/x$.  Completing these four endpoints gives a compact
	surface: it is obtained from two slit spheres by finitely many boundary
	identifications and endpoint fillings.  Away from the four endpoints each
	value of $x$ has two preimages, while the local models show precisely the
	stated branching.
\end{proof}

\begin{remark}[Euler characteristic check]
	The constructed surface has genus $1$ by Proposition~\ref{prop:cubic-branch-genus},
	so it is topologically a torus.
	This matches the ``one-handle'' intuition from gluing two slit spheres crosswise along two cuts.
\end{remark}

\begin{figure}[ht]
	\centering
	\begin{tikzpicture}[scale=1.05,>=Stealth]
		\draw[->] (-3.25,0)--(3.6,0) node[right] {$\infty$};
		\draw[very thick,red!70!black] (-2.5,0)--(-0.9,0);
		\draw[very thick,red!70!black,->] (0.9,0)--(3.15,0);
		\foreach \x/\lab in {-2.5/0,-0.9/1,0.9/2}{
			\fill (\x,0) circle (1.5pt) node[below=4pt] {$\lab$};
		}
		\node[red!70!black,above] at (-1.7,0.04) {$\gamma_1$};
		\node[red!70!black,above] at (2.0,0.04) {$\gamma_2$};
		\draw[dashed,blue!70!black] (-1.7,0) ellipse (1.18 and 0.67);
		\node[blue!70!black] at (-1.7,1.0) {a loop around one cut encloses two branch values};
	\end{tikzpicture}
	\caption{A cut system for the cubic double cover. The first cut joins
	$0$ to $1$; the second joins $2$ to $\infty$ on the Riemann sphere.
	Each cut pairs two branch values, so continuation around a cut changes
	the square root's sign twice and returns to the original sheet.}
	\label{fig:cubic-double-cuts}
\end{figure}

\subsubsection{From branching to the canonical sheaf}
\label{subsec:double-cover-canonical}

The preceding genus calculation has a finer holomorphic counterpart.  The
canonical sheaf $\omega_C$ of a compact Riemann surface $C$ is the line bundle
whose local sections are holomorphic $1$-forms.  Determining $\omega_C$ tells us
not only the genus, through $\deg\omega_C=2g(C)-2$, but also whether a global
holomorphic differential can be chosen without zeros.  For the cubic double
cover, that differential will be $dx/y$.

The argument will proceed in four stages. The four branch values first determine
an effective divisor $B$ of degree $4$ on $\PP^1$. A square root of its divisor
bundle is the line bundle $L\simeq\mathcal O_{\PP^1}(2)$, and the section cutting
out $B$ constructs the double cover. The tautological section on the cover then
identifies the ramification bundle with $\pi^*L$. Substitution into the canonical
form of Riemann--Hurwitz shows that $\omega_E$ is trivial, after which a local
calculation identifies $dx/y$ as a concrete nowhere-vanishing generator. The
divisor calculation, the line-bundle calculation, and the coordinate calculation
are therefore three descriptions of the same geometry.

We first separate two divisors that are easily confused.  Let
\[
	\pi:E\longrightarrow\PP^1,
	\qquad
	E:\ y^2=(x-a_1)(x-a_2)(x-a_3),
\]
where the $a_i$ are distinct, and let $p_i$ and $p_\infty$ denote the
ramification points above $a_i$ and $\infty$.  The \emph{branch divisor} lies
on the base, whereas the \emph{ramification divisor} lies on the covering
curve:
\[
	B=[a_1]+[a_2]+[a_3]+[\infty]\quad\text{on }\PP^1,
	\qquad
	R=p_1+p_2+p_3+p_\infty\quad\text{on }E.
\]
Since each local map has the form $z=u^2$, pulling back a branch point gives
twice its ramification point.  Consequently
\begin{equation}
	\pi^*B=2R.                                                   \label{eq:branch-ramification-double}
\end{equation}

There is also an intrinsic construction of the cover.  Because
$\deg B=4$, one has $\mathcal O_{\PP^1}(B)\simeq\mathcal O_{\PP^1}(4)$.
Choose
\[
	L\simeq\mathcal O_{\PP^1}(2),
	\qquad L^{\otimes2}\simeq\mathcal O_{\PP^1}(B),
\]
and a section $s\in H^0(\PP^1,L^{\otimes2})$ with $\operatorname{div}(s)=B$.  On an open
set where $L$ is trivial, the equation of the cover is simply
\[
	\eta^2=s.
\]
The zeros of $s$ are exactly the points where the two local choices of
$\eta$ coalesce.  The transition functions of $L$ are what make these local
square roots compatible from one coordinate chart to another.

Equivalently, the functions on the cover form the rank--two algebra
\begin{equation}
	\pi_*\mathcal O_E\simeq\mathcal O_{\PP^1}\oplus L^{-1}.       \label{eq:double-cover-algebra}
\end{equation}
The first summand is fixed by the involution $(x,\eta)\mapsto(x,-\eta)$ and
the second changes sign.  Multiplication of two sections of $L^{-1}$ lands
in $L^{-2}$ and is carried to $\mathcal O_{\PP^1}$ by the section $s$.
Thus \eqref{eq:double-cover-algebra} is an algebra decomposition, not merely
a decomposition of modules, and
\[
	E\simeq\operatorname{Spec}_{\PP^1}
	   \bigl(\mathcal O_{\PP^1}\oplus L^{-1}\bigr).
\]

On $E$ there is a tautological section of $\pi^*L$.  Locally it is the
coordinate $\eta$, and its zero divisor is $R$.  Hence
\begin{equation}
	\mathcal O_E(R)\simeq\pi^*L.                                 \label{eq:ramification-line-bundle}
\end{equation}
Combining this identity with the canonical form of Riemann--Hurwitz gives
\begin{align}
	\omega_E
	&\simeq\pi^*\omega_{\PP^1}\otimes\mathcal O_E(R)\notag\\
	&\simeq\pi^*(\omega_{\PP^1}\otimes L).                       \label{eq:canonical-double-cover}
\end{align}
Since $\omega_{\PP^1}\simeq\mathcal O_{\PP^1}(-2)$ and
$L\simeq\mathcal O_{\PP^1}(2)$, the two factors cancel:
\begin{equation}
	\omega_E\simeq\mathcal O_E.                                  \label{eq:elliptic-canonical-trivial}
\end{equation}
This proves more than the numerical equality $\deg\omega_E=0$.  A line
bundle of degree zero need not be trivial; here the explicit isomorphism
\eqref{eq:canonical-double-cover} proves triviality.

The same conclusion is visible directly in coordinates.  For a simple
branch point, write $x=a_i+t^2$ and $y=t\,v(t)$ with $v(0)\ne0$.  Then
\[
	\frac{dx}{y}=\frac{2t\,dt}{t\,v(t)}=\frac{2}{v(t)}\,dt,
\]
which is holomorphic and nonzero at $t=0$.  At the point over infinity one
may choose a parameter $t$ for which
\[
	x=t^{-2},\qquad y=t^{-3}(c+O(t^2)),\qquad c\ne0.
\]
It follows that
\[
	\frac{dx}{y}=\left(-\frac{2}{c}+O(t^2)\right)dt,
\]
so $dx/y$ is also holomorphic and nonzero there.  It therefore furnishes a
global nowhere-vanishing section of $\omega_E$, realizing
\eqref{eq:elliptic-canonical-trivial} concretely.
For the double-cover construction and the canonical bundle formula, see
\cite[Ch.~III]{Miranda1995} and \cite[Ch.~III]{Forster1981}.

\begin{proposition}[Canonical holomorphic $1$-form]
	\label{prop:holomorphic-1form}
	On $\mathcal C_1$ the differential
	\[
	\omega:=\frac{dx}{y}
	\]
	extends to a holomorphic $1$-form on the compactification.
\end{proposition}

\begin{proof}
	The local computations in Subsection~\ref{subsec:double-cover-canonical}
	show that $dx/y$ is holomorphic and nonzero at every finite branch point
	and at the point over infinity.  Away from those points, $x$ is a local
	coordinate and $y\ne0$, so $dx/y$ is again holomorphic.  Thus it extends
	to the compactification and has no zeros there.
\end{proof}

\begin{proposition}[Uniformization by an elliptic integral]
	\label{prop:uniformization}
	Let $\Sigma$ be the compactification of $\mathcal C_1$ and fix a base point $p_0\in\Sigma$.
	Define on the universal cover $\widetilde{\Sigma}$ the map
	\[
	u(p)=\int_{p_0}^{p}\omega,\qquad \omega=\frac{dx}{y}.
	\]
	Then the period lattice
	\[
	\Lambda=\left\{\int_{\gamma}\omega:\ \gamma\in H_1(\Sigma,\Z)\right\}\subset\C
	\]
	is a lattice and the induced map $\Sigma\to\C/\Lambda$ is a biholomorphism.
\end{proposition}

\begin{proof}
	The preceding local calculations show that $\omega=dx/y$ is a nonzero
	holomorphic form on $\Sigma$.  Its zero divisor has degree $2g-2=0$;
	being effective, that divisor is empty.  Hence locally $\omega=f(z)\,dz$
	with $f(z)\ne0$, and its integral has nonzero derivative $f(z)$.
	On the universal cover, integrals are independent of the path between
	fixed endpoints, so $u$ is a holomorphic local coordinate.  A deck
	transformation changes $u$ by the integral of $\omega$ over the
	corresponding loop.
	For an oriented homology basis $(a,b)$ with $a\cdot b=1$, put
	$A=\int_a\omega$ and $B=\int_b\omega$.  The Riemann bilinear relation gives
	$i\int_\Sigma\omega\wedge\overline\omega=2\operatorname{Im}(\overline A B)>0$
	because $i\omega\wedge\overline\omega$ is a positive area form.  Thus
	$\operatorname{Im}(\overline A B)>0$: $A$ and $B$ are linearly independent
	over $\mathbb R$, and $\Lambda=\mathbb ZA+\mathbb ZB$ is discrete of rank two.
	The path-change formula makes $u\bmod\Lambda$ a well-defined holomorphic
	map $F:\Sigma\to\mathbb C/\Lambda$.  It is a local biholomorphism by
	the nonvanishing of $\omega$.  Its image is open; it is also closed by
	compactness, so connectedness of the target makes $F$ surjective.
	A proper local biholomorphism is a finite covering.  Under the natural
	identification $H_1(\mathbb C/\Lambda,\mathbb Z)=\Lambda$, the cycles
	$a,b$ map to $A,B$, which generate all of $\Lambda$.  Therefore the
	subgroup corresponding to this covering has index one: $F$ has one sheet
	and is a biholomorphism.
\end{proof}

\subsubsection{A beginner's calculation guide to divisors and ramification}
\label{subsec:divisor-ramification-laboratory}

This subsection collects the calculations that a reader should be able to
perform before using divisors abstractly. A divisor is a bookkeeping device:
it records one integer at each point. The difficult step is usually not adding
the formal sums, but finding the correct local parameter. Once a local
parameter $t$ is chosen, every nonzero meromorphic function has the form
$t^m u(t)$ with $u(0)\ne0$, and the required coefficient is simply $m$.

\paragraph{A complete divisor calculation on $\PP^1$.}

Let $z$ be the affine coordinate on $\PP^1\setminus\{\infty\}$ and let
$w=1/z$ near infinity. For a rational function
\[
	f(z)=c\frac{\prod_{j=1}^r(z-a_j)^{m_j}}
	{\prod_{k=1}^s(z-b_k)^{n_k}},
\]
where common factors have been cancelled, the finite orders are read directly:
$\ord_{a_j}(f)=m_j$ and $\ord_{b_k}(f)=-n_k$. At infinity, substitute
$z=1/w$:
\[
	f(1/w)=c\,w^{\sum_kn_k-\sum_jm_j}
	\frac{\prod_j(1-a_jw)^{m_j}}
	{\prod_k(1-b_kw)^{n_k}}.
\]
The quotient on the right is nonzero at $w=0$, so
\[
	\ord_\infty(f)=\sum_kn_k-\sum_jm_j.
\]
Consequently
\[
	\operatorname{div}(f)=
	\sum_jm_j[a_j]-\sum_kn_k[b_k]
	+\left(\sum_kn_k-\sum_jm_j\right)[\infty],
\]
and its degree is zero. This is the elementary model for the general theorem
that every principal divisor on a compact Riemann surface has degree zero.

\begin{example}[A rational function with a zero at infinity]
	For
	\[
	f(z)=\frac{(z-1)^2}{(z+2)^3},
	\]
	we have orders $2$ at $1$ and $-3$ at $-2$. Since the denominator has
	degree one more than the numerator, $f(1/w)=w(1-w)^2/(1+2w)^3$, so
	$\ord_\infty(f)=1$. Therefore
	\[
	(f)=2[1]-3[-2]+[\infty],\qquad \deg(f)=2-3+1=0.
	\]
\end{example}

\begin{example}[Divisor of a differential on $\PP^1$]
	The differential $dz$ has no zero or pole on the affine chart. Near infinity,
	$z=1/w$ gives $dz=-w^{-2}dw$, so
	\[
	\operatorname{div}(dz)=-2[\infty].
	\]
	More generally,
	\[
	\operatorname{div}(f(z)\,dz)=(f)-2[\infty].
	\]
	Thus every meromorphic one-form on $\PP^1$ has divisor degree $-2$, in
	agreement with $\deg K_{\PP^1}=2g-2=-2$.
\end{example}

\begin{proposition}[Local order rules]
	\label{prop:local-order-rules-expanded}
	For nonzero meromorphic functions $f,g$ and a point $p$,
	\[
	\ord_p(fg)=\ord_p(f)+\ord_p(g),\qquad
	\ord_p(f/g)=\ord_p(f)-\ord_p(g),
	\]
	and
	\[
	\ord_p(f+g)\ge\min\{\ord_p(f),\ord_p(g)\},
	\]
	with equality in the last formula when the two orders are different.
\end{proposition}

\begin{proof}
	Write $f=t^mu$ and $g=t^nv$ in a local parameter, with $u(0),v(0)\ne0$.
	The product and quotient formulas follow immediately. If $m<n$, then
	$f+g=t^m(u+t^{n-m}v)$, and the expression in parentheses is nonzero at
	$0$, proving equality. If $m=n$, the leading coefficients may cancel,
	which explains why only an inequality is guaranteed.
\end{proof}

\paragraph{Linear equivalence as a change of meromorphic trivialization.}

Two divisors satisfy $D\sim D'$ precisely when $D-D'=(f)$ for a meromorphic
function $f$. Division by $f$ then gives an explicit isomorphism
\[
	\mathcal O(D')\longrightarrow\mathcal O(D),\qquad s\longmapsto s/f.
\]
Indeed, $(s)+D'\ge0$ implies
$(s/f)+D=(s)-(f)+D=(s)+D'\ge0$. The point is invariant: a principal
divisor changes the local representative of a section but not the isomorphism
class of the line bundle. Keeping track of this sign is easiest by checking
one example, such as $(z)=[0]-[\infty]$ on $\PP^1$.

\begin{proposition}[Degree-one divisors on $\PP^1$ are equivalent]
	For any $a,b\in\PP^1$, one has $[a]\sim[b]$. Consequently every divisor
	on $\PP^1$ is linearly equivalent to $(\deg D)[\infty]$.
\end{proposition}

\begin{proof}
	For finite $a,b$, the function $(z-a)/(z-b)$ has divisor $[a]-[b]$.
	If one of the points is infinity, use $z-a$ or its reciprocal. Additivity
	then moves every point occurring in $D$ to infinity while preserving its
	coefficient. This also gives $\Pic(\PP^1)\cong\Z$ concretely.
\end{proof}

\begin{example}[Computing $L(D)$ on $\PP^1$]
	Let $D=m[\infty]$. The condition $(f)+D\ge0$ says that $f$ has no finite
	poles and a pole of order at most $m$ at infinity. Thus, for $m\ge0$,
	\[
	L(D)=\operatorname{span}\{1,z,\ldots,z^m\},
	\qquad \ell(D)=m+1.
	\]
	For $m<0$, a section would be an entire function vanishing at infinity,
	so $L(D)=0$. This elementary computation is the genus-zero case of
	Riemann--Roch.
\end{example}

\paragraph{How to compute the degree and ramification of a map.}

Let $F:X\to Y$ be nonconstant and holomorphic between compact Riemann
surfaces. Near $p\in X$, choose parameters $z$ at $p$ and $w$ at $F(p)$.
After multiplying $z$ and $w$ by local units, the map has the form
$w=z^{e_p}$. The integer $e_p$ is the local degree. It answers three related
questions:
\begin{enumerate}
	\item a nearby regular value has $e_p$ inverse images near $p$;
	\item the derivative vanishes to order $e_p-1$ at $p$;
	\item the pullback of a divisor point satisfies
	$F^*[F(p)]=\sum_{q\in F^{-1}(F(p))}e_q[q]$.
\end{enumerate}

\begin{proposition}[Pullback preserves degree up to the map degree]
	\label{prop:pullback-divisor-degree-expanded}
	If $F$ has degree $n$ and $D$ is a divisor on $Y$, then
	\[
	\deg(F^*D)=n\deg D.
	\]
\end{proposition}

\begin{proof}
	It is enough to take $D=[y]$. By the fiber-degree formula,
	\[
	F^*[y]=\sum_{p\in F^{-1}(y)}e_p[p],
	\qquad \sum_{p\in F^{-1}(y)}e_p=n.
	\]
	Extend by linearity to arbitrary divisors.
\end{proof}

\begin{example}[The map $z\mapsto z^n$ on the sphere]
	Let $F:\PP^1\to\PP^1$ be $F(z)=z^n$. At $0$ the local expression is
	$w=z^n$, so $e_0=n$. Near infinity use $u=1/z$ and $v=1/w$; then
	$v=u^n$, so $e_\infty=n$. Every nonzero finite value has $n$ distinct
	preimages and no ramification. Hence
	\[
	R_F=(n-1)[0]+(n-1)[\infty],\qquad \deg R_F=2n-2.
	\]
	Riemann--Hurwitz reads
	$-2=n(-2)+(2n-2)$, which checks exactly.
\end{example}

\begin{example}[A polynomial as a map of spheres]
	A polynomial $P$ of degree $n$ extends to $P:\PP^1\to\PP^1$ with
	local degree $n$ at infinity. At a finite point $p$, if
	$P(z)-P(p)=c(z-p)^e+\cdots$ with $c\ne0$, then $e_p=e$ and
	$P'(z)$ vanishes to order $e-1$. Thus the finite ramification divisor is
	the zero divisor of $P'$, while infinity contributes $n-1$. Since
	$\deg P'=n-1$, the total ramification is $2n-2$, again agreeing with
	Riemann--Hurwitz.
\end{example}

\paragraph{Hyperelliptic curves: a reusable genus computation.}

Consider the smooth projective model of
\[
	X:\quad y^2=P(x),
\]
where $P$ is square-free of degree $d\ge3$. Projection to $x$ has degree two.
Each finite simple zero of $P$ is a branch value. Infinity is a branch value
exactly when $d$ is odd. Therefore the number of branch values is
\[
	B=\begin{cases}
	d,&d\text{ even},\\
	d+1,&d\text{ odd}.
	\end{cases}
\]
Since $B$ is even, write $B=2g+2$. The double-cover formula gives
\[
	g=\frac{B-2}{2}
	=\begin{cases}
	(d-2)/2,&d\text{ even},\\
	(d-1)/2,&d\text{ odd}.
	\end{cases}
\]

\begin{proposition}[Holomorphic differentials on a hyperelliptic model]
	\label{prop:hyperelliptic-differentials-expanded}
	If $X$ is the smooth projective model of $y^2=P(x)$ with
	$\deg P=2g+1$ or $2g+2$, then
	\[
	\frac{dx}{y},\quad \frac{x\,dx}{y},\quad\ldots,\quad
	\frac{x^{g-1}dx}{y}
	\]
	are holomorphic and form a basis of $H^0(X,K_X)$.
\end{proposition}

\begin{proof}
	At a finite branch point $a$, use $x-a=t^2$ and $y=t u(t)$ with
	$u(0)\ne0$. Then
	$x^jdx/y=2x(t)^jdt/u(t)$, which is holomorphic. Away from branch points
	the assertion is immediate. At infinity in the odd-degree case, take
	$x=t^{-2}$ and $y=t^{-(2g+1)}u(t)$; then
	\[
	\frac{x^jdx}{y}=-2t^{2g-2j-2}u(t)^{-1}dt,
	\]
	which is holomorphic precisely for $0\le j\le g-1$. The even-degree
	case has two points at infinity and the same inequality follows using
	$x=t^{-1}$ there. The displayed forms are linearly independent because
	a relation would give a polynomial of degree at most $g-1$ vanishing
	identically in $x$. Their number is $g=\dim H^0(X,K_X)$, so they form a
	basis.
\end{proof}

\begin{example}[Genus two]
	For $y^2=(x-a_1)\cdots(x-a_6)$ with distinct $a_i$, there are six simple
	branch values and no branching at infinity. Hence $g=(6-2)/2=2$, and
	$dx/y$ and $x\,dx/y$ form a basis of holomorphic differentials. This is
	the first genus for which the curve is not its own one-dimensional
	Jacobian; the two differentials later supply coordinates on a
	two-dimensional period torus.
\end{example}

\paragraph{A checklist for branched-cover problems.}

For an equation or a holomorphic map, the following order of work prevents
the most common mistakes.
\begin{enumerate}
	\item Compactify the source and target; do not forget points at infinity.
	\item At every suspected ramification point, choose an actual local
	parameter on the source and put the map in the form $z\mapsto z^e$.
	\item Distinguish ramification points on $X$ from branch values on $Y$.
	\item Check one fiber using $\sum e_p=n$.
	\item Form $R_F=\sum(e_p-1)[p]$ and apply Riemann--Hurwitz.
	\item If differentials are involved, verify them in the same local
	parameters; the factor $z^{e-1}$ in $d(z^e)$ is the ramification term.
\end{enumerate}
These steps convert the topological picture of glued sheets into divisor and
canonical-bundle calculations that can be checked line by line.

\EssayHyunminYoon

\part{Duality}
\label{part:duality}

Throughout this part \(X\) denotes a compact Riemann surface, i.e.\ a connected compact complex
one--dimensional manifold. We write \(g\) for its genus, \(K_X\) for the canonical bundle,
\(\mathcal O_X\) for the structure sheaf, and we fix a K\"ahler area form \(\mathrm d\!A\) on \(X\).
For a holomorphic line bundle \(L\to X\) we denote by \(A^{p,q}(X,L)\) the smooth \((p,q)\)-forms with values in \(L\),
and by \(\bar\partial_L:A^{0,q}(X,L)\to A^{0,q+1}(X,L)\) the Dolbeault operator.
Our conventions for Hodge and Serre duality follow \cite{Forster1981,Wells2008}; for the de Rham
pairing and Poincar\'e duality we also use \cite{BottTu1982}.

\section{Duality on Compact Riemann Surfaces}
\index{duality}
\label{sec:duality-compact-RS}

Duality pairs a geometric object with a complementary one so that
information about either can be recovered from the pairing.  On an
oriented surface, Poincar\'e duality pairs cohomology classes by
integration.  Serre duality is its holomorphic counterpart: it pairs
$H^1(X,L)$ with holomorphic sections of $K_X\otimes L^\vee$.
The latter space is often easier to compute from zeros and poles.
This is why duality is placed before Riemann--Roch: for
$L=\mathcal O_X(D)$, it turns the apparently abstract $H^1$ correction
into the concrete space $H^0(X,K_X\otimes\mathcal O_X(-D))$ of
differentials with prescribed vanishing.  The first subsection fixes the linear-algebra meaning of
a perfect pairing before the geometric pairings are constructed.

\subsection{Perfect pairings and Poincar\'e duality}
\index{Poincar\'e duality}\index{perfect pairings}
\label{subsec:PD}

\subsubsection{Linear algebra: perfect pairings}

\begin{definition}[Perfect pairing]
	\label{def:perfect-pairing}
	Let \(V,W\) be finite-dimensional complex vector spaces. A bilinear map
	\[
	\langle\ ,\ \rangle:V\times W\to\C
	\]
	is a \emph{perfect pairing} if the induced linear maps
	\[
	\Phi_V:V\to W^\vee,\quad v\mapsto (w\mapsto \langle v,w\rangle),
	\qquad
	\Phi_W:W\to V^\vee,\quad w\mapsto (v\mapsto \langle v,w\rangle)
	\]
	are isomorphisms.
\end{definition}

\begin{lemma}[Equivalent criteria]
	\label{lem:perfect-pairing-criteria}
	For a bilinear map \(\langle\ ,\ \rangle:V\times W\to\C\) the following are equivalent:
	\begin{enumerate}
		\item The pairing is perfect.
		\item Both \(\Phi_V\) and \(\Phi_W\) are injective.
		\item The pairing is nondegenerate in both factors.
		\item \(\Phi_V\) is an isomorphism.
		\item \(\Phi_W\) is an isomorphism.
	\end{enumerate}
	If \(\dim V=\dim W\) is known in advance, injectivity in either one factor is sufficient.
\end{lemma}

\begin{proof}
	The kernels of \(\Phi_V\) and \(\Phi_W\) are precisely the left and right radicals of the pairing,
	so (2) and (3) are equivalent. If both maps are injective, then
	\(\dim V\leq\dim W\) and \(\dim W\leq\dim V\). Thus the dimensions are equal and both
	maps are isomorphisms. This proves (2)--(5) and (1) equivalent. Under the additional hypothesis
	\(\dim V=\dim W\), either one of the two injective maps is already an isomorphism; its transpose is
	the other induced map, which is then an isomorphism as well.
\end{proof}

\begin{example}[Evaluation pairing]
	\label{ex:evaluation-pairing}
	For any \(V\), the map \(\langle v,f\rangle=f(v)\) on \(V\times V^\vee\) is perfect.
\end{example}

\subsubsection{Poincar\'e duality pairing on a compact surface}

\begin{definition}[Integral Poincar\'e dual of a cycle]
	\label{def:PD-cycle}
	Let $Y$ be a closed oriented smooth manifold of real dimension $n$,
	with fundamental class $[Y]\in H_n(Y,\Z)$.  Integral Poincar\'e duality
	gives a cap-product isomorphism.  Using the cap convention
	$\langle v,[Y]\frown u\rangle=\langle u\smile v,[Y]\rangle$,
	we fix its sign by setting
	\[
	D_Y:H^{n-r}(Y,\Z)\longrightarrow H_r(Y,\Z),\qquad
	D_Y(u):=(-1)^{r(n-r)}[Y]\frown u.
	\]
	Graded commutativity then gives
	$\langle v,D_Y(u)\rangle=\langle v\smile u,[Y]\rangle$
	for $v\in H^r(Y,\Z)$.  The cap-product formula, rather than
	this numerical identity alone, also specifies $D_Y$ on torsion.  For
	$[Z]\in H_r(Y,\Z)$ define its \emph{Poincar\'e dual} by
	\[
	\mathrm{PD}_Y([Z]):=D_Y^{-1}([Z])\in H^{n-r}(Y,\Z).
	\]
	In particular, a closed oriented curve on a complex surface has a
	degree-two dual class, while a one-cycle on a Riemann surface has a
	degree-one dual class.
\end{definition}

The integral isomorphism, not an equality of real integrals, gives the
uniqueness in this definition \cite{Hatcher2002,BottTu1982}.  If
$Z\subset Y$ is a closed oriented smooth submanifold representing
$[Z]$, then the defining identity becomes
\[
\langle v\smile\mathrm{PD}_Y([Z]),[Y]\rangle
=\langle i^*v,[Z]\rangle.
\]
For a closed real $r$-form $\eta$ and a closed form $\zeta_Z$
representing the real image of $\mathrm{PD}_Y([Z])$, this reads
$\int_Y\eta\wedge\zeta_Z=\int_Z i^*\eta$.
Real differential forms alone cannot distinguish torsion in integral
cohomology, which is why the integral definition comes first.

Let \(H^k_{\mathrm{dR}}(X;\R)\) be de Rham cohomology with real coefficients.
Because \(\dim_\R X=2\), wedge product followed by integration defines
\begin{equation}\label{eq:PD-pairing}
	H^k_{\mathrm{dR}}(X;\R)\times H^{2-k}_{\mathrm{dR}}(X;\R)\longrightarrow \R,
	\qquad
	([\alpha],[\beta])\longmapsto \int_X \alpha\wedge \beta.
\end{equation}

\begin{lemma}[Well-definedness of \eqref{eq:PD-pairing}]
	\label{lem:PD-well-defined}
	The value \(\int_X \alpha\wedge\beta\) depends only on the cohomology classes \([\alpha]\), \([\beta]\).
\end{lemma}

\begin{proof}
	If \(\alpha\) is replaced by \(\alpha+d\gamma\), then
	\[
	\int_X (\alpha+d\gamma)\wedge\beta = \int_X \alpha\wedge\beta + \int_X d\gamma\wedge\beta.
	\]
	Since \(\beta\) is closed (we may take a closed representative), we have
	\(d(\gamma\wedge\beta)=d\gamma\wedge\beta + (-1)^{\deg\gamma}\gamma\wedge d\beta=d\gamma\wedge\beta\).
	Thus \(\int_X d\gamma\wedge\beta=\int_X d(\gamma\wedge\beta)=0\) by Stokes' theorem (compactness, no boundary).
	The same argument applies if \(\beta\) is replaced by \(\beta+d\delta\).
\end{proof}

\begin{theorem}[Poincar\'e duality for \(X\)]
	\label{thm:PD}
	For \(k=0,1,2\), the pairing \eqref{eq:PD-pairing} is perfect.
	In particular \(H^0_{\mathrm{dR}}(X;\R)\cong\R\), \(H^2_{\mathrm{dR}}(X;\R)\cong\R\), and
	\(\dim_\R H^1_{\mathrm{dR}}(X;\R)=2g\).
\end{theorem}

\begin{proof}
	Fix a Riemannian metric \(g_X\) compatible with the complex structure and orientation,
	and use the Hodge star and $L^2$ inner product fixed in
	Definition~\ref{def:hodge-star}.
	Let \(\Delta=d\,d^*+d^*d\) be the Laplacian and \(\mathcal H^k(X)=\{\omega\in\Omega^k(X)\mid \Delta\omega=0\}\)
	the harmonic \(k\)-forms. By the Hodge theorem (Theorem~\ref{thm:hodge-theorem} below),
	every cohomology class \([\alpha]\in H^k_{\mathrm{dR}}(X;\R)\) has a unique harmonic representative \(\alpha_h\in\mathcal H^k(X)\),
	and the map \([\alpha]\mapsto \alpha_h\) is an isomorphism. Moreover, \(*\) maps \(\mathcal H^k(X)\) isomorphically to \(\mathcal H^{2-k}(X)\).
	
	To prove nondegeneracy of \eqref{eq:PD-pairing} in the first factor, take \([\alpha]\in H^k_{\mathrm{dR}}(X;\R)\)
	and let \(\alpha=\alpha_h\in\mathcal H^k(X)\) be its harmonic representative. Assume
	\[
	\int_X \alpha\wedge\beta=0\qquad\text{for all }\ [\beta]\in H^{2-k}_{\mathrm{dR}}(X;\R).
	\]
	In particular, choose \(\beta=*\,\alpha\), which is harmonic and therefore closed, hence represents a class in \(H^{2-k}\).
	Then
	\[
	0=\int_X \alpha\wedge *\alpha = \|\alpha\|_{L^2}^2.
	\]
	Thus \(\alpha=0\), hence \([\alpha]=0\). This proves nondegeneracy in the first factor; nondegeneracy in the second factor
	follows similarly (or by symmetry of the argument after swapping roles). Therefore \eqref{eq:PD-pairing} is perfect
	by Lemma~\ref{lem:perfect-pairing-criteria}.
	
	The dimension statements are standard consequences: \(H^0_{\mathrm{dR}}(X;\R)\cong\R\) since \(X\) is connected,
	and perfection forces \(H^2_{\mathrm{dR}}(X;\R)\cong\R\). Finally, \(\dim H^1_{\mathrm{dR}}(X;\R)=2g\) is the classical fact
	that the first Betti number of a genus-\(g\) surface is \(2g\).
\end{proof}

\begin{example}[Explicit pairing matrix on a torus]
	\label{ex:torus-PD-matrix}
	Let \(X=\C/\Lambda\) with coordinate \(z=x+iy\). The classes \([\mathrm d x],[\mathrm d y]\) form a basis of
	\(H^1_{\mathrm{dR}}(X;\R)\). Put \(\mathrm{Area}(\Lambda)=\int_X \mathrm d x\wedge \mathrm d y>0\) and normalize
	\(\omega_0=\mathrm{Area}(\Lambda)^{-1}\mathrm d x\wedge \mathrm d y\). Then
	\[
	\int_X \mathrm d x\wedge \mathrm d x = 0,\quad
	\int_X \mathrm d y\wedge \mathrm d y = 0,\quad
	\int_X \mathrm d x\wedge \mathrm d y = \mathrm{Area}(\Lambda),\quad
	\int_X \mathrm d y\wedge \mathrm d x = -\mathrm{Area}(\Lambda),
	\]
	so the pairing on \(H^1\times H^1\) is represented (in the basis \([\mathrm d x],[\mathrm d y]\)) by
	\(\mathrm{Area}(\Lambda)\begin{psmallmatrix} 0 & 1\\ -1 & 0\end{psmallmatrix}\).
\end{example}

\begin{exercise}
	\label{exr:PD-well-defined}
	Apply the graded Leibniz rule
	\(d(\gamma\wedge\beta)=d\gamma\wedge\beta+(-1)^{\deg\gamma}\gamma\wedge d\beta\).
	Determine the resulting sign when $d\beta=0$, and use Stokes' theorem to prove
	Lemma~\ref{lem:PD-well-defined}.
\end{exercise}

\subsection{Hodge theory on a compact surface}
\label{subsec:hodge}

We state the analytic results used above and record the consequences needed later.
The orthogonal decomposition depends on elliptic regularity for \(\Delta\),
whose proof is not reproduced here; see \cite{Warner1983,Wells2008}.

\begin{theorem}[Hodge decomposition]
	\label{thm:hodge-decomposition}
	Let \(M\) be a compact oriented Riemannian manifold without boundary,
	and fix \(0\leq k\leq\dim M\).
	Then
	\[
	\Omega^k(M)=\mathcal H^k(M)\ \oplus\ d\Omega^{k-1}(M)\ \oplus\ d^*\Omega^{k+1}(M)
	\]
	orthogonally with respect to the \(L^2\)-inner product.
\end{theorem}

The displayed splitting is an imported analytic theorem, not a consequence
of the formal identity $\Delta=dd^*+d^*d$ alone. Its proof uses the
elliptic Green operator; the next theorem is proved here from this splitting.

\begin{theorem}[Hodge theorem]
	\label{thm:hodge-theorem}
	The map \(\mathcal H^k(M)\to H^k_{\mathrm{dR}}(M;\R)\), \(\omega\mapsto [\omega]\), is an isomorphism.
	In particular, each de Rham cohomology class has a unique harmonic representative.
\end{theorem}

\begin{proof}
	A harmonic form is closed: taking the $L^2$ inner product of
	$\Delta\omega=0$ with $\omega$ gives
	\[
	0=(\Delta\omega,\omega)_{L^2}
	  =\|d\omega\|_{L^2}^2+\|d^*\omega\|_{L^2}^2.
	\]
	Thus the displayed map is defined.  If a harmonic form is exact,
	$\omega=d\eta$, then
	\[
	\|\omega\|_{L^2}^2=(\omega,d\eta)_{L^2}
	=(d^*\omega,\eta)_{L^2}=0,
	\]
	so the map is injective.

	For surjectivity, let $\alpha$ be a closed $k$-form and apply Hodge
	decomposition:
	\[
	\alpha=h+d\beta+d^*\gamma,
	\qquad h\in\mathcal H^k(M).
	\]
	Because $d\alpha=0$ and $dh=d^2\beta=0$, one has $dd^*\gamma=0$.
	Taking the inner product with $\gamma$ yields
	$\|d^*\gamma\|_{L^2}^2=(dd^*\gamma,\gamma)_{L^2}=0$.
	Hence $\alpha=h+d\beta$, so $[\alpha]=[h]$.  The harmonic representative
	is unique by injectivity.
\end{proof}

\begin{lemma}[Type decomposition and the Hodge star]
	\label{lem:star-on-RS}
	On a Riemann surface $X$, the complexified cotangent bundle splits into
	the two eigenspaces of the Hodge star:
	\[
		\Omega^{1,0}(X)=\ker(*+i),\qquad
		\Omega^{0,1}(X)=\ker(*-i).
	\]
\end{lemma}

\begin{proof}
	The coordinate calculation following Definition~\ref{def:hodge-star} gives
	$*dz=-i\,dz$ and $*d\bar z=i\,d\bar z$.  Since $dz$ and $d\bar z$
	locally span the $(1,0)$- and $(0,1)$-forms, respectively, the displayed
	eigenspace identities follow and are independent of the chosen coordinate.
\end{proof}

\begin{exercise}
	\label{exr:star-square}
	Show that on \(1\)-forms on a surface, \(*^2=-\mathrm{id}\). Deduce that \(*\) defines a complex structure on
	\(\mathcal H^1(X)\) and therefore \(\dim_\R \mathcal H^1(X)\) is even.
\end{exercise}

\subsection{Serre duality on a compact Riemann surface}
\index{Serre duality!curves}
\label{subsec:serre}

Let \(L\to X\) be a holomorphic line bundle and \(L^\vee\) its dual.

\subsubsection{Dolbeault cohomology model for \(H^1(X,L)\)}

Since \(\dim_\C X=1\), we have \(A^{0,2}(X,L)=0\). Hence every \(L\)-valued \((0,1)\)-form is automatically \(\bar\partial_L\)-closed, and
\begin{equation}\label{eq:dolbeault-H1}
	H^1(X,L)\ \cong\ H^{0,1}(X,L)\ \cong\ A^{0,1}(X,L)\big/\bar\partial_L A^{0,0}(X,L).
\end{equation}

\subsubsection{The Serre pairing and its independence of choices}

\begin{definition}[Serre pairing]
	\label{def:serre-pairing}
	For \(s\in H^0(X,K_X\otimes L^\vee)\) and \([\eta]\in H^1(X,L)\) represented by \(\eta\in A^{0,1}(X,L)\),
	define
	\begin{equation}\label{eq:SerrePairing}
		\langle s,[\eta]\rangle_S := \int_X s\wedge \eta\in\C,
	\end{equation}
	where locally we contract the \(L^\vee\otimes L\to\mathcal O_X\) factor and wedge the resulting \((1,0)\)-form with
	the \((0,1)\)-form to obtain a \((1,1)\)-form.
\end{definition}

\begin{lemma}[Well-definedness of the Serre pairing]
	\label{lem:serre-well-defined}
	The value \(\int_X s\wedge\eta\) depends only on \(s\in H^0(X,K_X\otimes L^\vee)\) and the cohomology class \([\eta]\in H^1(X,L)\).
\end{lemma}

\begin{proof}
	Linearity in each variable is clear. It remains to show independence of the representative \(\eta\).
	If \(\eta\) is replaced by \(\eta+\bar\partial_L u\) with \(u\in A^{0,0}(X,L)\), then
	\[
	\int_X s\wedge(\eta+\bar\partial_L u)=\int_X s\wedge\eta+\int_X s\wedge\bar\partial_L u.
	\]
	We show the second integral is \(0\). Since \(s\) is holomorphic, \(\bar\partial_{K_X\otimes L^\vee}s=0\).
	Use the Leibniz rule for \(\bar\partial\) on bundle-valued forms:
	\[
	\bar\partial\,(s\wedge u)=(\bar\partial s)\wedge u + (-1)^{1}\, s\wedge \bar\partial_L u
	= -\, s\wedge \bar\partial_L u.
	\]
	Hence \(s\wedge\bar\partial_L u = -\,\bar\partial(s\wedge u)\).
	After contracting \(L^\vee\otimes L\to\mathcal O_X\), the product \(s\wedge u\) is an ordinary
	\((1,0)\)-form (equivalently, a smooth section of \(K_X\)),
	so \(\bar\partial(s\wedge u)\) is a \((1,1)\)-form. On a curve, \(\partial(s\wedge u)\) has type \((2,0)=0\),
	so \(d(s\wedge u)=\bar\partial(s\wedge u)\). Therefore
	\[
	\int_X s\wedge\bar\partial_L u = -\int_X \bar\partial(s\wedge u)= -\int_X d(s\wedge u)=0
	\]
	by Stokes' theorem. This proves independence of the representative \(\eta\).
\end{proof}

\subsubsection{Harmonic representatives and a concrete dual section}

Fix a Hermitian metric \(h\) on \(L\) and a K\"ahler metric on \(X\).
Let \(\mathcal H^{0,1}(X,L)\subset A^{0,1}(X,L)\) denote the \(\bar\partial_L\)-harmonic forms.
Hodge theory for the Dolbeault complex gives
\begin{equation}\label{eq:H1-harmonic}
	H^1(X,L)\ \cong\ \mathcal H^{0,1}(X,L),
\end{equation}
where a class is represented by its unique harmonic representative.

\begin{lemma}[Metric realization of the dual section]
	\label{lem:sharp-map}
	For every \(\eta\in A^{0,1}(X,L)\) there is a unique smooth section
	\(\eta^\sharp\in A^{0,0}(X,K_X\otimes L^\vee)\) characterized pointwise by
	\begin{equation}\label{eq:sharp-definition}
		s\wedge\eta=\langle s,\eta^\sharp\rangle\,\mathrm d\!A
	\end{equation}
	for every local section \(s\) of \(K_X\otimes L^\vee\), where the Hermitian pairing is linear
	in its first argument. Then:
	\begin{enumerate}
		\item The assignment \(\eta\mapsto\eta^\sharp\) is a conjugate-linear bundle isomorphism.
		\item If \(\eta\in\mathcal H^{0,1}(X,L)\), then \(\eta^\sharp\) is holomorphic, i.e.\ \(\eta^\sharp\in H^0(X,K_X\otimes L^\vee)\).
	\end{enumerate}
\end{lemma}

\begin{proof}
	At each point, wedge product and contraction identify an \(L\)-valued \((0,1)\)-covector with a
	complex-linear functional on \(K_X\otimes L^\vee\), after division by the nonzero area form.
	The Hermitian Riesz isomorphism represents that functional by a unique vector \(\eta^\sharp\).
	This proves existence, uniqueness, smooth dependence, and conjugate linearity, without choosing coordinates
	or a frame.

	For the second assertion, the formal-adjoint identity for the Dolbeault complex identifies
	\(\bar\partial_L^*\eta=0\) with
	\(\bar\partial_{K_X\otimes L^\vee}\eta^\sharp=0\). This can also be checked in a holomorphic
	coordinate and a holomorphic frame by writing the Chern connection of the two Hermitian metrics.
	Since every \((0,1)\)-form on a curve is \(\bar\partial_L\)-closed, harmonicity is equivalent to
	\(\bar\partial_L^*\eta=0\); hence \(\eta^\sharp\) is holomorphic. Applying the same construction
	backwards proves that this restriction is a conjugate-linear bijection
	\(\mathcal H^{0,1}(X,L)\simeq H^0(X,K_X\otimes L^\vee)\).
\end{proof}

\begin{lemma}[Serre pairing as an \(L^2\) pairing]
	\label{lem:serre-L2}
	Let \(s\in H^0(X,K_X\otimes L^\vee)\) and \(\eta\in\mathcal H^{0,1}(X,L)\). Then
	\begin{equation}\label{eq:Serre-L2}
		\int_X s\wedge \eta \;=\; \int_X \langle s,\eta^\sharp\rangle \,\mathrm d\!A,
	\end{equation}
	where \(\langle\ ,\ \rangle\) is the pointwise Hermitian pairing induced by the chosen metrics.
\end{lemma}

\begin{proof}
	Equation \eqref{eq:sharp-definition} is an equality of global \((1,1)\)-forms. Integrating it over
	\(X\) gives \eqref{eq:Serre-L2} directly.
\end{proof}

\subsubsection{Serre duality theorem}

\begin{theorem}[Serre duality on \(X\)]
	\label{thm:Serre}
	For every holomorphic line bundle \(L\to X\), the pairing
	\[
	H^0(X,K_X\otimes L^\vee)\times H^1(X,L)\longrightarrow\C,\qquad
	(s,[\eta])\longmapsto \int_X s\wedge\eta
	\]
	is perfect. Consequently,
	\[
	H^1(X,L)\ \cong\ H^0(X,K_X\otimes L^\vee)^\vee,
	\qquad
	h^1(X,L)=h^0(X,K_X\otimes L^\vee).
	\]
\end{theorem}

\begin{proof}
	By \eqref{eq:H1-harmonic}, we identify \(H^1(X,L)\) with \(\mathcal H^{0,1}(X,L)\).
	Define a linear map
	\[
	T:\mathcal H^{0,1}(X,L)\to H^0(X,K_X\otimes L^\vee)^\vee,
	\qquad
	T(\eta)(s):=\int_X s\wedge\eta.
	\]
	By Lemma~\ref{lem:serre-well-defined}, this is well-defined.
	
	We show \(T\) is an isomorphism. Consider the map \(\sharp:\mathcal H^{0,1}(X,L)\to H^0(X,K_X\otimes L^\vee)\),
	\(\eta\mapsto \eta^\sharp\), from Lemma~\ref{lem:sharp-map}. This map is a conjugate-linear bijection,
	so it identifies \(\mathcal H^{0,1}(X,L)\) with \(H^0(X,K_X\otimes L^\vee)\) up to complex conjugation.
	
	Now Lemma~\ref{lem:serre-L2} identifies \(T(\eta)\) with the \(L^2\)-pairing against \(\eta^\sharp\):
	\[
	T(\eta)(s)=\int_X s\wedge\eta = \int_X \langle s,\eta^\sharp\rangle\,\mathrm d\!A.
	\]
	The \(L^2\)-pairing on the finite-dimensional space \(H^0(X,K_X\otimes L^\vee)\) is nondegenerate, hence perfect.
	Therefore, as \(\eta\mapsto \eta^\sharp\) is bijective, the map \(T\) is an isomorphism. This proves that the Serre pairing is perfect.
\end{proof}

\begin{corollary}[Serre-dual form of Riemann--Roch]
	\label{cor:RR-serre}
	Assuming the Riemann--Roch theorem proved in
	Theorem~\ref{thm:RR}, a holomorphic line bundle $L\to X$ satisfies
	\[
	h^0(X,L)-h^0(X,K_X\otimes L^\vee)=\deg L + 1 - g.
	\]
\end{corollary}

\begin{proof}
	Theorem~\ref{thm:RR} gives \(h^0(X,L)-h^1(X,L)=\deg L+1-g\).
	By Serre duality (Theorem~\ref{thm:Serre}), \(h^1(X,L)=h^0(X,K_X\otimes L^\vee)\). Substitute this into Riemann--Roch.
\end{proof}

\subsection{Concrete examples}
\label{subsec:duality-examples}

\subsubsection{The sphere \(X\simeq \CP^1\) (genus \(0\))}

Recall \(K_{\CP^1}\cong \mathcal O_{\CP^1}(-2)\). For \(m\in\Z\),
\[
H^0(\CP^1,\mathcal O(m))=
\begin{cases}
	\text{homogeneous polynomials of degree \(m\) in \(Z_0,Z_1\)},& m\ge 0,\\
	0,& m<0,
\end{cases}
\]
so \(h^0(\CP^1,\mathcal O(m))=m+1\) for \(m\ge 0\) and \(0\) for \(m<0\).
Then Serre duality gives
\[
H^1(\CP^1,\mathcal O(m))
\cong
H^0(\CP^1,\mathcal O(-m-2))^\vee,
\]
hence
\[
h^1(\CP^1,\mathcal O(m))=
\begin{cases}
	0,& m\ge -1,\\
	-m-1,& m\le -2.
\end{cases}
\]

\begin{example}[Serre pairing nontriviality on \(\CP^1\) for negative bundles]
	\label{ex:CP1-serre-nontrivial}
	Take \(m=-2\). Then \(H^1(\CP^1,\mathcal O(-2))\cong H^0(\CP^1,\mathcal O(0))^\vee\cong\C\).
	On the standard cover \(U_0=\{Z_0\neq 0\}\), \(U_1=\{Z_1\neq 0\}\), write \(z=Z_1/Z_0\) on \(U_0\).
	A generator of \(H^0(\CP^1,\mathcal O(0))\) is the constant section \(1\).
	Under the duality isomorphism, a nonzero class in \(H^1(\CP^1,\mathcal O(-2))\) pairs nontrivially with \(1\),
	so the pairing is perfect in this case.
\end{example}

\begin{exercise}
	\label{exr:cp1-h0-h1}
	Using the two-chart cocycle model for \(\mathcal O(m)\), show directly that any global section over \(\CP^1\) is determined by
	a polynomial of degree \(\le m\) when \(m\ge 0\), and must vanish when \(m<0\). Then use Serre duality to compute \(h^1(\CP^1,\mathcal O(m))\).
\end{exercise}

\subsubsection{The torus \(X=\C/\Lambda\) (genus \(1\))}

Let \(z=x+iy\) on \(\C\) and \(\Lambda\subset\C\) a lattice. The translation-invariant forms \(dx,dy,dz,d\bar z\) descend to \(X\).
Compute
\[
dz\wedge d\bar z=(dx+i\,dy)\wedge(dx-i\,dy)=-2i\,dx\wedge dy.
\]
Hence
\[
\int_X dz\wedge d\bar z = -2i \int_X dx\wedge dy = -2i\,\mathrm{Area}(\Lambda).
\]
The canonical bundle is trivial: \(K_X\cong \mathcal O_X\) with holomorphic generator \(dz\), so \(H^0(X,K_X)\cong\C\).
Also \(H^1(X,\mathcal O_X)\cong\C\), represented by the harmonic class of \(d\bar z\).
Therefore the Serre pairing reads
\[
\langle a\,dz,\ b\,[d\bar z]\rangle_S=\int_X ab\,dz\wedge d\bar z = -2i\,ab\,\mathrm{Area}(\Lambda),
\]
which is nonzero whenever \(ab\neq 0\). Hence the pairing is perfect.

\begin{exercise}
	\label{exr:torus-h1-O}
	Let $X=\C/\Lambda$ be the flat complex torus considered above. Show
	\(H^1(X,\mathcal O_X)\cong \C\) by proving that the space of harmonic
	\((0,1)\)-forms is one-dimensional, generated by $d\bar z$.
	Hint: harmonic forms on a flat torus are translation invariant.
\end{exercise}

\subsection{A compact dictionary for later use}
\label{subsec:duality-dictionary}

For a compact Riemann surface \(X\) and a holomorphic line bundle \(L\to X\), we record:
\begin{align*}
	&\text{(Poincar\'e)} &&
	H^k_{\mathrm{dR}}(X;\C)\ \cong\ H^{2-k}_{\mathrm{dR}}(X;\C)^\vee,
	&& k=0,1,2;\\[2pt]
	&\text{(Dolbeault)} &&
	H^1(X,L)\ \cong\ H^{0,1}(X,L)\ \cong\ A^{0,1}(X,L)\big/\bar\partial_LA^{0,0}(X,L);\\[2pt]
	&\text{(Serre)} &&
	H^1(X,L)\ \cong\ H^0(X,K_X\otimes L^\vee)^\vee;\\[2pt]
	&\text{(Riemann--Roch)} &&
	h^0(X,L)-h^0(X,K_X\otimes L^\vee)=\deg L + 1 - g.
\end{align*}

\begin{exercise}
	\label{exr:dictionary-check}
	Let \(X=\CP^1\) and \(L=\mathcal O(m)\). Use the dictionary to verify Riemann--Roch explicitly for all \(m\in\Z\) by plugging in
	the formulas for \(h^0(\mathcal O(m))\) and \(h^0(\mathcal O(-m-2))\).
\end{exercise}




	\clearpage
	\part{Riemann--Roch and Intersection Theory}
	\label{part:rr}

\section{The Riemann--Roch Theorem on Compact Riemann Surfaces}
\index{Riemann--Roch theorem}
\label{sec:RR}

The Riemann--Roch theorem is the basic structural statement that controls
the size of spaces of meromorphic functions and differentials with prescribed
zero/pole behavior. It balances local analytic constraints (orders at points),
global geometry (line bundles/divisors), and topology (the genus).

Throughout, $M$ is a compact connected Riemann surface of genus $g$.
We freely use the preceding chapters for divisors, $\mathcal O_M(D)$, $\Pic(M)$,
\v{C}ech and sheaf cohomology, the exponential sequence, and Chern classes and
degree; the Dolbeault--Hodge and Serre-duality results used below come from
Part~\ref{part:duality}.
The divisor formulation and add-a-point proof of Riemann--Roch are cross-checked against
\cite{Forster1981,FarkasKra1992,Miranda1995}.

The chapter is organized around one calculation rather than a succession of
independent facts.  First, divisors are translated into the section spaces of
line bundles, and Serre duality identifies the correction term.  Next, adding
one point changes the Euler characteristic by one; iterating this local
calculation proves the theorem.  The canonical degree and the sphere and
torus examples then test the formula.  Finally, vanishing, linear systems,
and the Bergman kernel show how a dimension identity becomes geometric
information.  Keeping this chain visible is especially useful in applications:
one should first inspect degrees, then determine whether the dual correction
term can survive, and only afterward compute a basis of sections.

\subsection{Divisors, section spaces, and the duality input}
\label{subsec:RR:notation}

We only record notation (definitions already given in Ch.~6, 13--14).
Let $\Div(M)$ be the free abelian group on points of $M$.
For $D=\sum n_p p$, write $\deg D=\sum n_p$.
For $f\in\mathcal M(M)^\times$, write $(f)=\sum \nu_p(f)\,p$ (principal divisor), and
$D_1\sim D_2$ if $D_1-D_2=(f)$.
Then $\Pic(M)=\Div(M)/\!\sim$.

\subsubsection{Degree of principal divisors}
\label{subsec:RR:degprincipal}

\begin{proposition}[Degree of principal divisors]
	\label{prop:RR:deg-principal}
	If $f\in\mathcal M(M)^\times$, then $\deg(f)=0$.
\end{proposition}

\begin{proof}
	The zeros and poles of $f$ form a finite set. Choose pairwise disjoint
	coordinate discs $U_p$ around these points, each with a coordinate $z$
	centered at $p$.
	On $U_p\setminus\{p\}$ we can write
	\[
	f(z)=z^{\nu_p(f)}u_p(z),\qquad u_p(0)\neq 0,
	\]
	so
	\[
	\frac{f'(z)}{f(z)}=\frac{\nu_p(f)}{z}+\frac{u_p'(z)}{u_p(z)}.
	\]
	Hence $\operatorname{Res}_{p}\!\big(\tfrac{f'}{f}\,dz\big)=\nu_p(f)$.
	Summing residues and applying the residue theorem on the compact surface $M$ gives
	\[
	\sum_p \nu_p(f)
	=\sum_p \operatorname{Res}_{p}\!\left(\frac{f'}{f}\,dz\right)
	=0,
	\]
	so $\deg(f)=\sum_p\nu_p(f)=0$.
\end{proof}

\subsubsection{\texorpdfstring{The spaces $L(D)$ and $I(D)$ and their sheaf interpretations}{The spaces L(D) and I(D) and their sheaf interpretations}}
\label{subsec:RR:LDID}

Fix a divisor $D$. Recall the finite-dimensional spaces
\[
\begin{aligned}
L(D)&=\{\,f\in\mathcal M(M)\mid (f)+D\ge 0\,\}\cup\{0\},\\
I(D)&=\{\,\omega\in\mathcal M^1(M)\mid (\omega)-D\ge 0\,\}\cup\{0\}.
\end{aligned}
\]
Set $\ell(D)=\dim_{\mathbb C} L(D)$ and $i(D)=\dim_{\mathbb C} I(D)$.

\begin{lemma}[Identifying $L(D)$ with sections of $\mathcal O_M(D)$]
	\label{lem:RR:LD-sections}
	There is a natural identification
	\[
	L(D)\ \cong\ H^0\!\big(M,\mathcal O_M(D)\big),
	\qquad\text{hence}\qquad
	\ell(D)=h^0\!\big(M,\mathcal O_M(D)\big).
	\]
\end{lemma}

\begin{proof}
	By definition of the divisor line bundle $\mathcal O_M(D)$,
	a holomorphic section over an open set $U$ is represented by a meromorphic function
	$f\in\mathcal M(U)$ satisfying $(f)+D|_U\ge 0$.
	Taking $U=M$ yields exactly $L(D)$.
\end{proof}

\begin{lemma}[Identifying $I(D)$ with sections of $K_M\otimes\mathcal O_M(-D)$]
	\label{lem:RR:ID-sections}
	There is a natural identification
	\[
	I(D)\ \cong\ H^0\!\big(M,K_M\otimes \mathcal O_M(-D)\big),
	\qquad\text{hence}\qquad
	i(D)=h^0\!\big(M,K_M\otimes \mathcal O_M(-D)\big).
	\]
\end{lemma}

\begin{proof}
	A meromorphic $1$-form $\omega$ is a meromorphic section of $K_M$.
	The condition $(\omega)-D\ge 0$ is equivalent to $\omega$ becoming holomorphic after twisting by
	$\mathcal O_M(-D)$, i.e.\ $\omega$ is a holomorphic section of $K_M\otimes\mathcal O_M(-D)$.
\end{proof}

\subsubsection{\texorpdfstring{Serre duality input: $i(D)=\ell(K-D)$}{Serre duality input}}
\label{subsec:RR:serre}

\begin{proposition}[Serre-duality identification of the index of speciality]
	\label{prop:RR:key-duality}
	Let $K$ be a canonical divisor class. Then for every divisor $D$,
	\[
	i(D)=\ell(K-D)=h^1\!\bigl(M,\mathcal O_M(D)\bigr).
	\]
	More precisely, Serre duality gives the natural isomorphism
	\[
	H^1\!\big(M,\mathcal O_M(D)\big)^\vee
	\ \cong\
	H^0\!\big(M,K_M\otimes\mathcal O_M(-D)\big).
	\]
\end{proposition}

\begin{proof}
	By Serre duality, for any line bundle $L$ on $M$ there is a perfect pairing
	\[
	H^0(M,K_M\otimes L^\vee)\times H^1(M,L)\to \mathbb C,
	\]
	hence $H^1(M,L)^\vee\cong H^0(M,K_M\otimes L^\vee)$.
	Apply this to $L=\mathcal O_M(D)$, so $L^\vee\simeq \mathcal O_M(-D)$.
	Taking dimensions and using Lemma~\ref{lem:RR:LD-sections} and Lemma~\ref{lem:RR:ID-sections}
	gives $h^1(\mathcal O_M(D))=i(D)=h^0(K_M\otimes\mathcal O_M(-D))=\ell(K-D)$.
\end{proof}

\subsection{Euler characteristic and degree additivity}
\label{subsec:RR:euler}

For a line bundle $L$ on $M$, define $\chi(L):=h^0(M,L)-h^1(M,L)$.

\begin{lemma}[Reducing Riemann--Roch to an Euler characteristic formula]
	\label{lem:RR:reduce}
	For every divisor $D$,
	\[
	\chi\big(\mathcal O_M(D)\big)=\ell(D)-\ell(K-D).
	\]
\end{lemma}

\begin{proof}
	By Lemma~\ref{lem:RR:LD-sections}, $h^0(\mathcal O_M(D))=\ell(D)$.
	By Proposition~\ref{prop:RR:key-duality}, $h^1(\mathcal O_M(D))=\ell(K-D)$.
	Subtracting gives the claim.
\end{proof}

\begin{proposition}[Degree additivity of $\chi$ on a curve]
	\label{prop:RR:chi-degree}
	For any divisor $D$,
	\[
	\chi\big(\mathcal O_M(D)\big)=\chi(\mathcal O_M)+\deg D.
	\]
\end{proposition}

\begin{proof}
	We give a concrete proof by reducing to the case of adding a single point.
	
	\smallskip
	\noindent\textbf{Step 1: the case $D=p$ (a single point).}
	There is a short exact sequence
	\[
	0\longrightarrow \mathcal O_M \longrightarrow \mathcal O_M(p)
	\longrightarrow \mathcal O_M(p)/\mathcal O_M \longrightarrow 0.
	\]
	The quotient is a skyscraper sheaf supported at $p$ with one-dimensional
	stalk.  Locally, $\mathcal O_M(p)$ allows a simple pole at $p$, and the
	quotient records its leading principal-part coefficient.  A local parameter
	identifies this stalk with $\mathbb C$, but the identification itself is not
	canonical. Writing the quotient intrinsically as the fiber
	$\mathcal O_M(p)|_p$, we obtain
	\[
	0\longrightarrow \mathcal O_M \longrightarrow \mathcal O_M(p)
	\longrightarrow \mathcal O_M(p)|_p\longrightarrow 0.
	\]
	Taking cohomology and using the vanishing of higher cohomology for a
	skyscraper sheaf yields an exact segment
	\[
	0\to H^0(\mathcal O_M)\to H^0(\mathcal O_M(p))
	\to H^0(\mathcal O_M(p)|_p)
	\to H^1(\mathcal O_M)\to H^1(\mathcal O_M(p))\to 0.
	\]
	Taking alternating sums of dimensions gives
	\[
	\chi\big(\mathcal O_M(p)\big)=\chi(\mathcal O_M)
	+\dim H^0(\mathcal O_M(p)|_p)=\chi(\mathcal O_M)+1.
	\]
	
	\smallskip
	\noindent\textbf{Step 2: adding or removing one point.}
	More generally, for every divisor $E$ and point $p$ there is an exact sequence
	\[
	0\longrightarrow\mathcal O_M(E)\longrightarrow\mathcal O_M(E+p)
	\longrightarrow\mathcal O_M(E+p)|_p\longrightarrow0.
	\]
	The quotient is more intrinsically the one-dimensional fiber
	$\mathcal O_M(E+p)|_p$; after choosing a local frame, the last arrow extracts
	the coefficient of the newly allowed principal part. Additivity of the
	alternating sum in the resulting long exact cohomology sequence gives
	\[
	\chi(\mathcal O_M(E+p))=\chi(\mathcal O_M(E))+1.
	\]
	Equivalently, replacing $E$ by $E-p$ gives
	\(\chi(\mathcal O_M(E-p))=\chi(\mathcal O_M(E))-1\).

	\smallskip
	\noindent\textbf{Step 3: iteration.}
	Starting with the zero divisor, add a point once for each positive coefficient of $D$ and remove a
	point once for each negative coefficient. Each addition raises both $\chi$ and the degree by one,
	and each removal lowers both by one. Hence
	\(\chi(\mathcal O_M(D))=\chi(\mathcal O_M)+\deg D\) for every divisor $D$.
\end{proof}

\begin{lemma}[Computing $\chi(\mathcal O_M)$]
	\label{lem:RR:chi-OM}
	\[
	\chi(\mathcal O_M)=1-g.
	\]
\end{lemma}

\begin{proof}
	We have $h^0(\mathcal O_M)=1$ because holomorphic functions on compact connected $M$ are constant.
	Also $h^1(\mathcal O_M)=g$ by Serre duality:
	\[
	h^1(\mathcal O_M)=h^0(K_M)=g.
	\]
	Hence $\chi(\mathcal O_M)=1-g$.
\end{proof}

\subsection{Riemann--Roch: statement, proof, and the one-point mechanism}
\label{subsec:RR:statement-proof}

\begin{theorem}[Riemann--Roch]
	\label{thm:RR}
	For any divisor $D$ on $M$,
	\[
	\ell(D)-\ell(K-D)=1-g+\deg D.
	\]
	Equivalently,
	\[
	\ell(D)-i(D)=1-g+\deg D.
	\]
\end{theorem}

\begin{proof}
	By Lemma~\ref{lem:RR:reduce},
	$\ell(D)-\ell(K-D)=\chi\big(\mathcal O_M(D)\big)$.
	By Proposition~\ref{prop:RR:chi-degree} and Lemma~\ref{lem:RR:chi-OM},
	\[
	\chi\big(\mathcal O_M(D)\big)=\chi(\mathcal O_M)+\deg D=(1-g)+\deg D.
	\]
	Combining yields $\ell(D)-\ell(K-D)=1-g+\deg D$.
	Finally, $i(D)=\ell(K-D)$ by Proposition~\ref{prop:RR:key-duality}, giving the equivalent form.
\end{proof}

\subsubsection{Where the extra degree goes: an exact-sequence calculation}
\index{Riemann--Roch theorem!adding a point}
\label{subsec:rr-one-point-obstruction}

The formula $\chi(\mathcal O_M(D+p))=\chi(\mathcal O_M(D))+1$ does not say
that adding $p$ always produces a new global meromorphic function. The long
exact sequence displays the alternative:
\begin{align*}
0&\longrightarrow H^0(\mathcal O_M(D))
\longrightarrow H^0(\mathcal O_M(D+p))
\longrightarrow \mathcal O_M(D+p)|_p\\
&\overset{\partial}{\longrightarrow}H^1(\mathcal O_M(D))
\longrightarrow H^1(\mathcal O_M(D+p))\longrightarrow0.
\end{align*}
The one-dimensional fiber represents the newly allowed leading Laurent
coefficient. If its image under $\partial$ is zero, a global section realizes
that coefficient and $h^0$ rises by one. If $\partial$ is injective, no such
global section exists and $h^1$ falls by one. These are the only possibilities
because the fiber has dimension one.

For example, on an elliptic curve $E$ one has
$h^0(\mathcal O_E)=h^1(\mathcal O_E)=1$. There is no nonconstant
meromorphic function with only one simple pole: such a function would give a
degree-one map $E\to\mathbb P^1$ and hence a biholomorphism, contradicting
$g(E)=1$. Therefore $h^0(\mathcal O_E(p))=1$ and the boundary map for the
first added point is an isomorphism onto $H^1(\mathcal O_E)$. It follows that
$h^1(\mathcal O_E(p))=0$. Adding the point once more now raises $h^0$:
$h^0(\mathcal O_E(2p))=2$. A degree-two meromorphic function exists, and its
two-dimensional space of sections is the beginning of the familiar double
cover $E\to\mathbb P^1$. On $\mathbb P^1$, by contrast,
$H^1(\mathcal O_{\mathbb P^1})=0$ already, so adding one point immediately
produces the section represented by an affine coordinate.

\subsection{The canonical degree and model computations}
\label{subsec:RR:degK}

\begin{theorem}[Degree of the canonical divisor]
	\label{thm:RR:degK}
	\[
	\deg K=2g-2.
	\]
\end{theorem}

\begin{proof}
	Apply Riemann--Roch (Theorem~\ref{thm:RR}) to $D=0$:
	$\ell(0)-\ell(K)=1-g$.
	Since $\ell(0)=h^0(\mathcal O_M)=1$, we get $\ell(K)=g$.
	Apply Riemann--Roch again to $D=K$:
	$\ell(K)-\ell(0)=1-g+\deg K$.
	Substitute $\ell(K)=g$ and $\ell(0)=1$ to obtain $g-1=1-g+\deg K$, hence $\deg K=2g-2$.
\end{proof}

\subsubsection{Examples and computations}
\label{subsec:RR:examples}

\paragraph{Example A: $M=\mathbb{CP}^1$.}
Here $g=0$ and $K\sim -2[\infty]$. Let $D=n[\infty]$.

\begin{proposition}[Computing $\ell$ on the sphere]
	\label{prop:RR:CP1}
	For $n\in\mathbb Z$,
	\[
	\ell\big(n[\infty]\big)=
	\begin{cases}
		n+1,& n\ge 0,\\
		0,& n<0.
	\end{cases}
	\]
	Moreover, $\ell(K-D)=\ell\big((-2-n)[\infty]\big)=\max\{0,-n-1\}$, and Riemann--Roch holds termwise.
\end{proposition}

\begin{proof}
	A meromorphic function on $\mathbb{CP}^1$ with pole order $\le n$ at $\infty$ is exactly a polynomial
	of degree $\le n$ in the affine coordinate $z$, plus $0$.
	Thus $\ell\big(n[\infty]\big)=n+1$ for $n\ge 0$ and $0$ for $n<0$.
	Since $K-D=(-2-n)[\infty]$, the same computation gives
	$\ell(K-D)=\max\{0,-n-1\}$.
	Substituting into Theorem~\ref{thm:RR} yields an identity in $n$ that is immediate to check.
\end{proof}

\paragraph{Example B: elliptic curve ($g=1$).}
Let $g=1$. Then $\deg K=0$ and $K\sim 0$.

Indeed, $h^0(K_M)=g=1$, so a nonzero holomorphic differential exists. Its zero divisor is effective
of degree $\deg K=0$, hence is empty; the differential therefore trivializes $K_M$.

\begin{proposition}[Riemann--Roch on an elliptic curve]
	\label{prop:RR:elliptic}
	If $g=1$, then for any divisor $D$,
	\[
	\ell(D)-\ell(-D)=\deg D.
	\]
	In particular, if $\deg D>0$ then $\ell(D)=\deg D$.
\end{proposition}

\begin{proof}
	Since $K\sim 0$, we have $\ell(K-D)=\ell(-D)$.
	Riemann--Roch gives $\ell(D)-\ell(-D)=1-1+\deg D=\deg D$.
	If $\deg D>0$, then $\deg(-D)<0$, hence $\ell(-D)=0$ (no nonzero meromorphic functions with poles bounded by a negative divisor),
	so $\ell(D)=\deg D$.
\end{proof}

\subsection{Consequences, geometric use, and the Bergman kernel}
\label{subsec:RR:asymptotic}

\begin{proposition}[Vanishing of $\ell(K-D)$ for large degree]
	\label{prop:RR:vanishing-large}
	If $\deg D\ge 2g-1$, then $\ell(K-D)=0$ and therefore $\ell(D)=\deg D+1-g$.
\end{proposition}

\begin{proof}
	We have $\deg(K-D)=2g-2-\deg D\le -1$.
	If $\ell(K-D)\neq 0$, there exists $0\neq f\in L(K-D)$ with $(f)+K-D\ge 0$.
	Taking degrees gives $0+\deg(K-D)\ge 0$, contradicting $\deg(K-D)<0$.
	Hence $\ell(K-D)=0$, and Riemann--Roch gives $\ell(D)=\deg D+1-g$.
\end{proof}

\begin{corollary}[The nonspecial range]
	\label{cor:RR-nonspecial-range}
	If $\deg D>2g-2$, then $h^1(M,\mathcal O_M(D))=0$ and
	\[
		h^0(M,\mathcal O_M(D))=\deg D+1-g.
	\]
	Thus in this range the dimension depends only on the degree, not on the
	position of the points in $D$.
\end{corollary}

\begin{proof}
	Serre duality identifies $h^1(\mathcal O_M(D))$ with
	$h^0(K_M\otimes\mathcal O_M(-D))=\ell(K-D)$.  The latter vanishes
	because its degree is negative.  Riemann--Roch then gives the formula.
\end{proof}

\subsubsection{Geometric meaning and how to use Riemann--Roch}
\label{subsec:RR:meaning}

Riemann--Roch says that $\ell(D)$ is controlled by the naive count $\deg D+1-g$,
with an explicit correction term $\ell(K-D)$ measuring global obstructions:
\[
\ell(D)=\deg D+1-g+\ell(K-D).
\]
Heuristically:
\begin{itemize}
	\item the positive coefficients of $D$ allow poles, while its negative
	coefficients require zeros; $\deg D$ is their signed total,
	\item $g$ measures the space of holomorphic $1$-forms,
	\item $\ell(K-D)$ counts the holomorphic differentials subject to the
	complementary divisor condition; by Serre duality it equals
	$h^1(\mathcal O_M(D))$.
\end{itemize}

\subsubsection{Riemann--Roch and the Bergman kernel}
\index{Bergman kernel!Riemann--Roch}
\label{subsec:RR:Bergman}

Riemann--Roch is a dimension theorem, whereas the Bergman kernel packages a
finite-dimensional space of holomorphic sections into an object that evaluates
those sections pointwise.  This distinction also explains their connection.
On a compact Riemann surface, holomorphic functions are constant, so the
Bergman kernel of holomorphic functions contains little geometric information.
The useful compact-surface kernel is instead built from holomorphic
$1$-forms, or more generally from sections of a Hermitian line bundle.

Equip $H^0(M,K_M)$ with the intrinsic inner product
\[
	\langle\omega,\eta\rangle
	=
	\frac{i}{2}\int_M\omega\wedge\overline{\eta}.
\]
If $\omega_1,\ldots,\omega_g$ is an orthonormal basis, the canonical Bergman
kernel is
\[
	B_M(p,\overline q)
	=
	\sum_{j=1}^{g}\omega_j(p)\,\overline{\omega_j(q)}.
\]
If $\operatorname{pr}_1,\operatorname{pr}_2:M\times M\to M$ are the two
projections, then $B_M$ is a section of
\[
	\operatorname{pr}_1^*K_M\otimes
	\operatorname{pr}_2^*\overline{K_M}
\]
on $M\times M$; in other words, it is holomorphic in the first variable and
antiholomorphic in the second.  The kernel is independent of the chosen
orthonormal basis.  Its reproducing property is the
basis expansion of the evaluation map: pairing a form with
$B_M(\,\cdot\,,\overline q)$ recovers its value at $q$.  Thus the kernel records
not only the number of holomorphic forms but also how they evaluate at every
point; see \cite{FarkasKra1992,Fay1973}.

Riemann--Roch supplies the dimension entering this construction.  Taking
$D=0$ gives
\[
	h^0(M,K_M)=g,
\]
so the sum defining $B_M$ has exactly $g$ terms.  It also detects whether the
kernel can vanish on the diagonal.  Choose a coordinate $z$ near $p$ and write
$\omega_j=a_j(z)\,dz$.  Then
\[
	B_M(p,\overline p)=
	\left(\sum_{j=1}^{g}|a_j(p)|^2\right)
	dz|_p\otimes d\overline z|_p,
\]
and positivity of the coefficient does not depend on the chosen coordinate.
Consequently $B_M(p,\overline p)=0$ precisely when every holomorphic
$1$-form vanishes at $p$.  The rank of evaluation is
\[
	\operatorname{rank}\!\left(H^0(M,K_M)\longrightarrow K_M|_p\right)
	=h^0(M,K_M)-h^0(M,K_M(-p)).
\]
For $g\ge1$, Riemann--Roch applied to $p$ gives
\[
	h^0(M,K_M(-p))=g-1,
\]
because $h^0(M,\mathcal O_M(p))=1$: a nonconstant meromorphic function with
only one simple pole would define a degree-one map to $\CP^1$ and force
$M\simeq\CP^1$.  Hence evaluation has rank one at every point and
$B_M(p,\overline p)>0$.  This is the analytic form of the fact that the
canonical linear system has no base points.  In genus zero there are no
holomorphic $1$-forms and $B_M=0$; on a torus, a normalized nowhere-vanishing
form $dz$ gives a rank-one kernel.

The same mechanism works for a Hermitian line bundle $L$.  The Bergman kernel
of $H^0(M,K_M\otimes L)$ is a sum over an orthonormal basis, while
Riemann--Roch computes
\[
	h^0(M,K_M\otimes L)-h^0(M,L^{-1})
	=\deg L+g-1.
\]
Comparing $L$ and $L(-p)$ gives the rank of evaluation at $p$ once the
correction terms in Riemann--Roch are known; degree alone need not decide
whether a particular point is a base point. Thus Riemann--Roch organizes
the dimension calculation, while the kernel refines it into pointwise
analytic data.

\subsection{Exercises}
\label{subsec:RR:ex}

\begin{exercise}[Principal divisors have degree $0$ via argument principle]
	Reprove Proposition~\ref{prop:RR:deg-principal} using the argument principle on small circles
	around zeros/poles and then summing contributions.
\end{exercise}

\begin{exercise}[Riemann--Roch for $D=0$]
	Use Theorem~\ref{thm:RR} with $D=0$ to show $\ell(K)=g$ and interpret this as
	$\dim H^0(M,K_M)=g$.
\end{exercise}

\begin{exercise}[Base-point freeness for large degree]
	Let $D$ be an effective divisor with $\deg D\ge 2g$.
	Show that for every $p\in M$ there is a section of $\mathcal O_M(D)$ that does not vanish at $p$;
	equivalently, show that $\ell(D)-\ell(D-p)=1$. Thus $|D|$ has no base points.
	Hint: compare $\ell(D)$ and $\ell(D-p)$ using Riemann--Roch.
\end{exercise}

\begin{exercise}[Hyperelliptic computation of a canonical divisor]
	Let $g\ge1$, let $a_1,\ldots,a_{2g+2}$ be distinct, and let $M$ be the
	smooth projective model of $y^2=\prod_{j=1}^{2g+2}(x-a_j)$.
	Compute the divisor of $dx/y$ and verify that its degree is $2g-2$.
	Hint: there are two points above infinity; use $t=1/x$ at each one.
\end{exercise}

\begin{exercise}[The pencil case of Clifford's inequality]
	Assume that $D$ and $K-D$ are effective and that $\dim|D|=1$.
	Show that $\ell(D)\le \deg D/2+1$. First remove the fixed part of $|D|$.
	The resulting base-point-free pencil defines a nonconstant map $M\to\PP^1$
	of degree at most $\deg D$; explain why this degree cannot be $1$ under the
	specialness hypothesis $K-D\ge0$. Conclude that $\deg D\ge2$, and use
	$\ell(D)=2$.
\end{exercise}

\EssayYonghyunJi

\section{A Cohomological Proof of Riemann--Roch}
\index{Riemann--Roch theorem!cohomological proof}
\label{sec:RR-complete-proof}

The preceding section gave the shortest sheaf-theoretic route to the formula and
its main consequences. This section isolates the analytic and cohomological
machinery behind that route: Dolbeault resolutions, the Hodge decomposition,
Serre duality, and the add-a-point exact sequence. The common calculation
behind the two proofs is
\[
\chi(\mathcal O(D+p))-\chi(\mathcal O(D))=1.
\]
The first proof obtains it by taking Euler characteristics of the
one-point divisor sequence (Proposition~\ref{prop:RR:chi-degree}).  Here we
keep track of the two separate dimensions: the new point either supplies
a new global section or removes one obstruction in $H^1$, and the two
changes add to one.  Serre duality then identifies that obstruction with
a differential in $H^0(K-D)$.  Thus this section explains the mechanism
behind the earlier numerical identity rather than proving a different
Riemann--Roch theorem. Readers already comfortable with those tools may
pass directly to the surface case.

Throughout let $M$ be a compact Riemann surface of genus $g$.
Write $\Omega=\Omega^{1,0}_M$ for the sheaf of holomorphic $1$-forms and
$\mathcal O=\mathcal O_M$ for the structure sheaf.
For a divisor $D$ we write $\mathcal O(D)=\mathcal O_M(D)$ for the associated line bundle and
$\Omega(-D):=\Omega\otimes\mathcal O(-D)$.

We use the same spaces as in Section~\ref{sec:RR}:
\[
L(D):=H^0(M,\mathcal O(D)),\qquad
I(D):=H^0(M,\Omega(-D)),
\]
and the dimensions $\ell(D):=\dim_{\mathbb C}L(D)$, $i(D):=\dim_{\mathbb C}I(D)$.

\subsection{Cohomological inputs: resolutions, exact sequences, and duality}
\label{subsec:RRcomplete:dolbeault}

A sequence of sheaf morphisms
$0\to\mathcal F\to\mathcal F_0\to\mathcal F_1\to\cdots$
is a \emph{fine resolution} of $\mathcal F$ if it is exact and each $\mathcal F_k$ is fine.
On a complex manifold, sheaves of smooth $(p,q)$-forms with values in a holomorphic vector bundle are fine,
because partitions of unity exist in the smooth category.

For a holomorphic line bundle $L\to M$, the Dolbeault complexes give fine resolutions
\[
0\to\mathcal O(L)\longrightarrow \mathcal A^{0,0}(L)
\xrightarrow{\ \bar\partial\ } \mathcal A^{0,1}(L)
\xrightarrow{\ \bar\partial\ } \mathcal A^{0,2}(L)\xrightarrow{\ \bar\partial\ }\cdots,
\]
\[
0\to\Omega(L)\longrightarrow \mathcal A^{1,0}(L)
\xrightarrow{\ \bar\partial\ } \mathcal A^{1,1}(L)
\xrightarrow{\ \bar\partial\ } \mathcal A^{1,2}(L)\xrightarrow{\ \bar\partial\ }\cdots.
\]
By the abstract de~Rham/Dolbeault theorem (computing sheaf cohomology via fine resolutions),
for every $q\ge 0$ there are natural isomorphisms
\[
H^q(M,\mathcal O(L))\cong H^{0,q}_{\bar\partial}(M,L),
\qquad
H^q(M,\Omega(L))\cong H^{1,q}_{\bar\partial}(M,L).
\]

\subsubsection{\texorpdfstring{A long exact sequence relating $\mathbb C$, $\mathcal O$, and $\Omega$}{A long exact sequence for the constant, structure, and differential-form sheaves}}
\label{subsec:RRcomplete:LES}

The short exact sequence behind the comparison is
\[
0\longrightarrow\underline{\mathbb C}\longrightarrow\mathcal O
\xrightarrow{\ d\ }\Omega\longrightarrow0.
\]
The kernel of $d$ consists of locally constant holomorphic functions.
Every holomorphic $1$-form has a holomorphic primitive on a sufficiently
small coordinate disk, which gives local surjectivity. Taking sheaf
cohomology therefore gives the long exact sequence
\begin{equation}\label{eq:LES-C-O-Omega}
\begin{aligned}
	0&\to H^0(M,\mathbb C)\to H^0(M,\mathcal O)\to H^0(M,\Omega)\\
	 &\to H^1(M,\mathbb C)\to H^1(M,\mathcal O)\to H^1(M,\Omega)\\
	 &\to H^2(M,\mathbb C)\to\cdots.
\end{aligned}
\end{equation}

On a compact Riemann surface there are no smooth forms of type $(2,0)$ or $(0,2)$.
In particular $H^2(M,\mathcal O)=0$, and Dolbeault identifies $H^1(M,\Omega)$ with $H^{1,1}_{\bar\partial}(M)$,
which is one-dimensional (top degree). Also $H^0(M,\mathcal O)\cong\mathbb C$ since holomorphic functions are constant.

There is one additional point in passing from the long sequence to the short
one: the maps induced by
$d:\mathcal O\to\Omega$ on $H^0$ and $H^1$ vanish. On $H^0$ this is
because every global holomorphic function is constant. On $H^1$, represent a
class by a harmonic $(0,1)$-form $\eta$. The K\"ahler identities imply
$d\eta=0$, so $\partial\eta=0$; hence the induced class in
$H^1(M,\Omega)\cong H^{1,1}_{\bar\partial}(M)$ is zero. This is the
dimension-one case of the $E_1$-degeneration of the Hodge--de Rham spectral
sequence. Exactness of \eqref{eq:LES-C-O-Omega} therefore yields
\begin{equation}\label{eq:H1C-split}
	0\longrightarrow H^0(M,\Omega)\longrightarrow H^1(M,\mathbb C)\longrightarrow H^1(M,\mathcal O)\longrightarrow 0.
\end{equation}
Taking dimensions gives
\[
\dim_{\mathbb C}H^1(M,\mathbb C)
=
\dim_{\mathbb C}H^0(M,\Omega)+\dim_{\mathbb C}H^1(M,\mathcal O).
\]
Since $H^1(M,\mathbb C)$ has complex dimension $b_1(M)$ and $b_1(M)=2g$, we obtain
\[
2g=\dim H^0(M,\Omega)+\dim H^1(M,\mathcal O).
\]

\subsubsection{Serre duality as a perfect pairing}
\label{subsec:RRcomplete:serre}

Let $L\to M$ be a holomorphic line bundle.
There is a bilinear pairing on Dolbeault cohomology
\[
H^{0,1}_{\bar\partial}(M,L)\times H^{1,0}_{\bar\partial}(M,L^\vee)\longrightarrow\mathbb C,
\qquad
([\varphi],[\psi])\longmapsto\int_M \varphi\wedge\psi,
\]
The integrand is a scalar-valued form of top degree. The value depends
only on the cohomology class of $\varphi$: replacing it by
$\varphi+\bar\partial\eta$ changes the integrand by
$\bar\partial(\eta\wedge\psi)$, whose integral is zero by Stokes
because $\bar\partial\psi=0$.
Using Hodge theory on a compact Riemann surface, every Dolbeault class has a harmonic representative and
the pairing is nondegenerate; hence it is perfect. Translating back to sheaves gives Serre duality
\begin{equation}\label{eq:Serre-duality}
	H^1(M,\mathcal O(L))^\vee \cong H^0(M,\Omega\otimes L^\vee).
\end{equation}

Taking $L=\mathcal O(D)$ and using $L^\vee\simeq \mathcal O(-D)$ yields
\[
H^1(M,\mathcal O(D))^\vee \cong H^0(M,\Omega(-D))=I(D),
\qquad
\dim H^1(M,\mathcal O(D))=i(D).
\]

In particular, applying \eqref{eq:Serre-duality} to $L=\mathcal O$ gives
$H^1(M,\mathcal O)^\vee\cong H^0(M,\Omega)$, so
\[
\dim H^1(M,\mathcal O)=\dim H^0(M,\Omega).
\]
Combining this with the dimension identity from \eqref{eq:H1C-split} gives the genus lemma.

\begin{lemma}[Genus lemma]
	\label{lem:genus}
	$\dim H^0(M,\Omega)=g$.
\end{lemma}

\begin{proof}
	From \eqref{eq:H1C-split} we have $2g=\dim H^0(M,\Omega)+\dim H^1(M,\mathcal O)$.
	By Serre duality, $\dim H^1(M,\mathcal O)=\dim H^0(M,\Omega)$.
	Hence $2g=2\dim H^0(M,\Omega)$, so $\dim H^0(M,\Omega)=g$.
\end{proof}

\subsection{The add-a-point sequence and the inductive proof}
\index{exact sequences!adding a point}
\label{subsec:RRcomplete:addpoint}

Fix a point $p\in M$. For every divisor $D$ there is a short exact sequence of sheaves
with one-dimensional quotient
$\mathcal Q_{D,p}:=\mathcal O(D+p)|_p$:
\begin{equation}\label{eq:SES-add-point}
	0\longrightarrow \mathcal O(D)\xrightarrow{\ \iota\ }\mathcal O(D+p)
	\xrightarrow{\ \rho\ }\mathcal Q_{D,p}\longrightarrow 0.
\end{equation}
The sheaf $\mathcal Q_{D,p}$ is supported at $p$ and has a one-dimensional
space of global sections; identifying that space with $\mathbb C$ requires a
choice of local frame. Taking cohomology (and using
$H^1(M,\mathcal Q_{D,p})=0$) gives an exact sequence
\begin{equation}\label{eq:LES-add-point}
\begin{aligned}
    0 &\to H^0(M,\mathcal O(D)) \to H^0(M,\mathcal O(D+p))\\
      &\xrightarrow{\ \rho\ } H^0(M,\mathcal Q_{D,p})
       \to H^1(M,\mathcal O(D))\\
      &\to H^1(M,\mathcal O(D+p)) \to 0.
\end{aligned}
\end{equation}

Let $V:=\mathrm{Im}(\rho)\subset H^0(M,\mathcal Q_{D,p})$ and
$W:=H^0(M,\mathcal Q_{D,p})/V$.
Exactness of \eqref{eq:LES-add-point} yields
\begin{equation}\label{eq:deltaL-deltaH1}
	\ell(D+p)=\ell(D)+\dim V,\qquad
	h^1(\mathcal O(D+p))=h^1(\mathcal O(D))-\dim W,
\end{equation}
and $\dim V+\dim W=1$.

Using $i(D)=\dim H^1(M,\mathcal O(D))$ from Serre duality, \eqref{eq:deltaL-deltaH1} becomes
\begin{equation}\label{eq:delta-identity}
	\big(\ell(D+p)-i(D+p)\big)-\big(\ell(D)-i(D)\big)=\dim V+\dim W=1.
\end{equation}

\paragraph{Interpretation.}
Either $V=H^0(M,\mathcal Q_{D,p})$, in which case some section of $\mathcal O(D+p)$ has nonzero
image in the one-dimensional quotient at $p$, and
$\ell(D+p)=\ell(D)+1$ while $i(D+p)=i(D)$; or $V=0$, in which case every
global section of $\mathcal O(D+p)$ already lies in $H^0(M,\mathcal O(D))$, and
$\ell(D+p)=\ell(D)$ while $i(D+p)=i(D)-1$. In both cases $\ell-i$ increases by exactly $1$.

\subsubsection{\texorpdfstring{Base case $D=0$}{Base case D=0}}
\label{subsec:RRcomplete:base}

\begin{proposition}[Base case $D=0$]
	\label{prop:RR-base}
	Riemann--Roch holds for $D=0$:
	\[
	\ell(0)-i(0)=1-g.
	\]
\end{proposition}

\begin{proof}
	We have $\ell(0)=\dim H^0(M,\mathcal O)=1$.
	Also $i(0)=\dim H^0(M,\Omega)=g$ by Lemma~\ref{lem:genus}.
	Therefore $\ell(0)-i(0)=1-g$.
\end{proof}

\subsubsection{Proof of Riemann--Roch by induction on degree}
\label{subsec:RRcomplete:induction}

\begin{theorem}[Riemann--Roch]
	\label{thm:RR-final}
	For every divisor $D$ on $M$,
	\[
	\ell(D)-i(D)=1-g+\deg D.
	\]
\end{theorem}

\begin{proof}
	\emph{Step 1: effective divisors.}
	Let $E$ be an effective divisor of degree $r$. Write $E=p_1+\cdots+p_r$.
	Starting from $0$ and applying \eqref{eq:delta-identity} successively gives
	\[
	\ell(E)-i(E)=\big(\ell(0)-i(0)\big)+r=(1-g)+\deg E,
	\]
	using Proposition~\ref{prop:RR-base}.
	
	\emph{Step 2: general divisors.}
	Write $D=E-F$ with $E,F$ effective and with disjoint support.
	Choose an ordering $F=q_1+\cdots+q_s$.
	Starting from $E$ and removing points one by one, we apply the add-point identity to
	$D'=E-(q_1+\cdots+q_{k-1})-q_k$ rewritten as
	\[
	(\ell(D'+q_k)-i(D'+q_k))-(\ell(D')-i(D'))=1,
	\]
	which is exactly \eqref{eq:delta-identity} with $D'$ in place of $D$ and $p=q_k$.
	Rearranging shows that each removal of a point decreases $\ell-i$ by $1$.
	After $s$ removals we obtain
	\[
	\ell(E-F)-i(E-F)=\big(\ell(E)-i(E)\big)-s.
	\]
	Insert Step~1 for $E$ and note that $s=\deg F$ to get
	\[
	\ell(D)-i(D)=(1-g)+\deg E-\deg F=(1-g)+\deg D.
	\]
\end{proof}

\subsection{Corollaries and asymptotics}
\label{subsec:RRcomplete:cor}

\begin{corollary}[Degree of the canonical divisor]
	\label{cor:degK}
	Let $K$ be a canonical divisor. Then $\deg K=2g-2$.
\end{corollary}

\begin{proof}
	Apply Theorem~\ref{thm:RR-final} to $D=0$ to get $\ell(0)-i(0)=1-g$.
	Since $i(0)=\ell(K)$ (by $i(0)=\dim H^1(\mathcal O)=\dim H^0(\Omega)=\ell(K)$), we obtain $\ell(K)=g$.
	Apply Theorem~\ref{thm:RR-final} to $D=K$:
	\[
	\ell(K)-i(K)=1-g+\deg K.
	\]
	But $i(K)=\dim H^0(\Omega(-K))=\dim H^0(\mathcal O)=\ell(0)=1$, hence
	$g-1=1-g+\deg K$, so $\deg K=2g-2$.
\end{proof}

\begin{corollary}[Asymptotic Riemann--Roch]
	\label{cor:asymptoticRR}
	If $\deg D>2g-2$, then $i(D)=0$ and hence $\ell(D)=\deg D+1-g$.
\end{corollary}

\begin{proof}
	If $\deg D>2g-2$, then $\deg(K-D)=2g-2-\deg D<0$.
	A divisor of negative degree has $\ell(K-D)=0$, hence $i(D)=\ell(K-D)=0$ by Serre duality.
	Now apply Riemann--Roch.
\end{proof}

\begin{corollary}[An effective-divisor existence bound]
	\label{cor:lowerbound}
	If $\deg D=g+\ell$ for some integer $\ell\ge 0$, then $\ell(D)\ge \ell+1$.
\end{corollary}

\begin{proof}
	Riemann--Roch gives $\ell(D)=1-g+\deg D+i(D)=1-g+(g+\ell)+i(D)=\ell+1+i(D)\ge \ell+1$.
\end{proof}

\subsubsection{Summary of the mechanism}
\label{subsec:RRcomplete:summary}

\begin{itemize}
	\item Fine (Dolbeault) resolutions compute sheaf cohomology and provide concrete models for $H^q(M,\mathcal O(L))$ and $H^q(M,\Omega(L))$.
	\item The long exact sequence \eqref{eq:LES-C-O-Omega} and Serre duality imply $\dim H^0(M,\Omega)=g$ (Lemma~\ref{lem:genus}).
	\item The add-a-point exact sequence \eqref{eq:SES-add-point} yields the jump identity \eqref{eq:delta-identity}:
	each time one adds a point, the quantity $\ell(D)-i(D)$ increases by exactly $1$.
	\item With the base case $D=0$, induction on degree proves $\ell(D)-i(D)=1-g+\deg D$ for all divisors.
\end{itemize}

\EssayHojaeJung

\clearpage
\section{Intersection Theory and Surface Riemann--Roch}
\index{intersection theory}\index{Riemann--Roch theorem!surfaces}
\label{sec:RR-surfaces-expanded}

This section upgrades the curve story to complex dimension $2$.
Divisors now define classes in $H^2(X,\mathbb Z)$, and there is a genuine \emph{quadratic}
numerical invariant: the intersection pairing.  Surface Riemann--Roch expresses
$\chi(\mathcal O_X(D))$ in terms of intersection numbers with $D$ and the canonical class $K_X$.

Throughout, $X$ is a smooth projective complex surface.
We freely use earlier material (divisors/line bundles, Chern classes, Poincar\'e duality, sheaf cohomology,
Serre duality, Dolbeault/Hodge where appropriate).
For local intersection multiplicities, divisor classes, adjunction, and surface Riemann--Roch we use
\cite{Fulton1998Intersection,Hartshorne1977,GriffithsHarris1994}.

\subsection{Local and cohomological intersection pairings}
\index{intersection pairing!surfaces}
\label{subsec:RRsurf:intersection}

\subsubsection{Local intersection multiplicity}

Let $C,C'\subset X$ be irreducible curves and let $p\in C\cap C'$.
Choose local equations $f,g\in\mathcal O_{X,p}$ defining $C$ and $C'$ near $p$.

\begin{definition}[Local intersection multiplicity]
	\label{def:local-intersection}
	The local intersection multiplicity of $C$ and $C'$ at $p$ is
	\[
	I_p(C,C') \ :=\ \dim_{\mathbb C}\,\mathcal O_{X,p}/(f,g).
	\]
\end{definition}

\begin{proposition}[Basic properties of $I_p(C,C')$]
	\label{prop:local-intersection-basic}
	\mbox{}
	\begin{enumerate}
		\item $I_p(C,C')$ is finite if and only if $C$ and $C'$ have no common component through $p$.
		\item If $C$ and $C'$ are smooth and meet transversely at $p$, then $I_p(C,C')=1$.
		\item If $C$ has a single branch at $p$ and
		$\varphi:(\mathbb C,0)\to(C,p)$ is its normalization parametrization, while $g$ is a local
		equation for $C'$, then
		$I_p(C,C')=\operatorname{ord}_0(g\circ\varphi)$. In general the intersection multiplicity is
		the sum of these orders over the branches of $C$ at $p$.
	\end{enumerate}
\end{proposition}

\begin{proof}
	(1) If $C$ and $C'$ share a component through $p$, then $(f,g)$ does not cut out a $0$--dimensional
	scheme at $p$ and the quotient has positive Krull dimension, hence infinite $\mathbb C$--dimension.
	Conversely, if they share no component through $p$, then $(f,g)$ defines a $0$--dimensional analytic
	set at $p$, so the local algebra is Artinian and has finite $\mathbb C$--dimension.
	
	(2) In holomorphic coordinates $(x,y)$ centered at $p$, transversality means $C=\{x=0\}$ and $C'=\{y=0\}$
	after a holomorphic change of coordinates. Then $\mathcal O_{X,p}/(x,y)\cong \mathbb C$.
	
	(3) Pulling $g$ back to the normalization of $C$ gives a holomorphic function on each normalized
	branch. The local intersection formula identifies the length of
	$\mathcal O_{X,p}/(f,g)$ with the sum of its vanishing orders there. For a unibranch germ this is
	just $\operatorname{ord}_0(g\circ\varphi)$.
\end{proof}

\subsubsection{Global intersection number}

Assume $C$ and $C'$ have no common component. Then $C\cap C'$ is finite.

\begin{definition}[Global intersection number]
	\label{def:global-intersection}
	The global intersection number is
	\[
	C\cdot C' \ :=\ \sum_{p\in C\cap C'} I_p(C,C').
	\]
	By bilinearity, this defines $D\cdot E$ for divisors
	$D=\sum n_i C_i$ and $E=\sum m_j C'_j$ with no common component by
	\[
	D\cdot E \ :=\ \sum_{i,j} n_i m_j \, (C_i\cdot C'_j).
	\]
	The definition for arbitrary divisor classes, including self-intersections, is
	the Chern-class pairing introduced in Definition~\ref{def:cohom-intersection};
	Theorem~\ref{thm:intersection-equals-chern} proves that it agrees with this
	local sum whenever the latter is defined.
\end{definition}

\begin{proposition}[Bilinearity, symmetry, and deformation invariance]
	\label{prop:intersection-form-properties}
	The local intersection pairing is symmetric and bilinear whenever the divisors
	have no common component. Its cohomological extension below is well-defined on
	all divisor classes and descends to numerical equivalence. In a family for which the
	scheme-theoretic intersections form a finite flat family over a connected base,
	the total intersection number is constant.
\end{proposition}

\begin{proof}
	For divisors without common components, bilinearity is built into the extension
	from irreducible curves, and symmetry follows from the symmetry of the local
	length $\dim_\mathbb C \mathcal O_{X,p}/(f,g)$ in $(f,g)$.  This statement
	concerns the pairs for which the local sum is defined; the global extension
	also assigns a meaning to self-intersections.
	For arbitrary divisor classes, Definition~\ref{def:cohom-intersection}
	uses $c_1(\mathcal O_X(D))$ and the cup product.  Since $c_1$ is additive
	under tensor products, this pairing is bilinear.  Since degree-two
	cohomology classes commute under the graded-commutativity rule, it is
	symmetric.  It depends only on the cohomology classes, so it is certainly
	unchanged by homological equivalence.  To obtain the distinct assertion
	about \emph{numerical} equivalence, recall its definition on a surface:
	$D\equiv_{\mathrm{num}}D'$ means $(D-D')\cdot C=0$ for every irreducible
	curve $C\subset X$.  Every divisor $E$ is an integral linear combination
	of such curves, hence bilinearity gives $(D-D')\cdot E=0$.  Thus replacing
	either argument by a numerically equivalent divisor leaves the pairing
	unchanged.  This conclusion follows from the definition and bilinearity,
	not merely from cohomological dependence.
	Finally, let $Z\to S$ be the finite flat scheme of intersections in the
	stated family.  Locally on $S$, the finite flat module $\mathcal O_Z$
	is free of some rank $r$; its fiber over $s$ consequently has length $r$
	(over the residue field).  The rank is locally constant, and therefore
	constant when $S$ is connected.  That fiber length is the sum of its
	local intersection multiplicities, proving the final assertion.
\end{proof}

\subsubsection{Chern classes and the cohomological intersection product}
\label{subsec:RRsurf:chern}

Let $L\to X$ be a holomorphic line bundle. Its first Chern class $c_1(L)\in H^2(X,\mathbb Z)$
is represented in de~Rham cohomology by curvature: choose a Hermitian metric $h$ and local holomorphic frame $s$,
then
\[
c_1(L)\ =\ \left[-\frac{i}{2\pi}\,\partial\bar\partial \log h(s,s)\right]\in H^2_{\mathrm{dR}}(X).
\]

\begin{definition}[Cohomological intersection number]
	\label{def:cohom-intersection}
	For line bundles $L,L'$ define
	\[
	L\cdot L' \ :=\
	\big\langle c_1(L)\smile c_1(L'),[X]\big\rangle \ \in\ \mathbb Z.
	\]
	Equivalently, if $\alpha,\alpha'$ are closed de~Rham representatives of the
	real images of these classes, then $L\cdot L'=\int_X\alpha\wedge\alpha'$.
	For divisors $D,E$, set $D\cdot E := \mathcal O_X(D)\cdot \mathcal O_X(E)$.
\end{definition}

For a smooth curve $C\subset X$, the symbol
$\mathrm{PD}([C])$ below abbreviates
$\mathrm{PD}_X([C])\in H^2(X,\Z)$ from
Definition~\ref{def:PD-cycle}; here $X$ has real dimension four.

\begin{proposition}[Compatibility with restriction degree]
	\label{prop:degree-restriction}
	If $C\subset X$ is a smooth irreducible curve and $L\to X$ a line bundle, then
	\[
	\deg(L|_C)=\langle i^*c_1(L),[C]\rangle
	=\big\langle c_1(L)\smile \mathrm{PD}([C]),[X]\big\rangle=L\cdot C.
	\]
\end{proposition}

\begin{proof}
	The degree of a line bundle on the compact curve $C$ is the evaluation
	of its first Chern class on $[C]$.  Naturality of $c_1$ under pullback
	therefore gives
	$\deg(L|_C)=\langle i^*c_1(L),[C]\rangle$.
	By the defining property of Poincar\'e duality in
	Definition~\ref{def:PD-cycle} for the inclusion
	$i:C\hookrightarrow X$, every $\alpha\in H^2(X,\mathbb Z)$ satisfies
	\[
	\langle i^*\alpha,[C]\rangle
	=\langle\alpha\smile\mathrm{PD}([C]),[X]\rangle.
	\]
	To identify the second factor, take the canonical section of
	$\mathcal O_X(C)$, whose zero divisor is $C$.  In a tubular neighborhood
	of $C$ its first Chern class is represented by the Thom class of the
	normal line bundle; extending that class to $X$ gives
	$c_1(\mathcal O_X(C))=\mathrm{PD}([C])$.  This is the divisor-class
	identity used in the last equality \cite{BottTu1982,GriffithsHarris1994}.
	Substituting it into the displayed pairing yields
	$\langle c_1(L)\smile c_1(\mathcal O_X(C)),[X]\rangle=L\cdot C$
	by Definition~\ref{def:cohom-intersection}.
\end{proof}

\begin{theorem}[Intersection as Chern class pairing]
	\label{thm:intersection-equals-chern}
	If $C,C'\subset X$ are irreducible curves with no common component, then
	\[
	C\cdot C' =
	\big\langle c_1(\mathcal O_X(C))\smile c_1(\mathcal O_X(C')),[X]\big\rangle.
	\]
\end{theorem}

\begin{proof}
	The Poincar\'e--Lelong formula identifies the cycle class of $C$ with
	$c_1(\mathcal O_X(C))$, and likewise for $C'$.  Because the curves have
	no common component, their intersections are isolated.  The
	local-to-global intersection theorem identifies the evaluation of the
	cohomological product with the sum of their local intersection
	multiplicities \cite{Fulton1998Intersection,GriffithsHarris1994}:
	$\sum_{p\in C\cap C'}\dim_{\mathbb C}\mathcal O_{X,p}/(f,g)$.
	Here $f$ and $g$ are local equations for $C$ and $C'$ at $p$.
	This is the geometric definition of $C\cdot C'$.  In the transverse case
	$f$ and $g$ form local coordinates, so the quotient is $\mathbb C$
	and the formula reduces to an ordinary count.
\end{proof}

\subsection{Adjunction and surface Riemann--Roch}
\index{adjunction formula}\index{canonical class!surfaces}
\label{subsec:RRsurf:adjunction}

\begin{theorem}[Adjunction formula for a smooth curve]
	\label{thm:adjunction}
	Let $C\subset X$ be a smooth irreducible curve. Then
	\[
	K_C \ \simeq\ (K_X\otimes \mathcal O_X(C))|_C,
	\]
	and consequently
	\[
	2g(C)-2 \ =\ C\cdot (C+K_X).
	\]
\end{theorem}

\begin{proof}
	Consider the exact sequence of holomorphic vector bundles on $C$:
	\[
	0\longrightarrow T_C \longrightarrow T_X|_C \longrightarrow N_{C/X} \longrightarrow 0.
	\]
	Taking determinants in this sequence and then dualizing yields
	\[
	K_C \ \simeq\ K_X|_C \otimes \det(N_{C/X}).
	\]
	But $\det(N_{C/X})\simeq \mathcal O_X(C)|_C$ for a smooth divisor $C$ (normal bundle corresponds to the divisor line bundle).
	Hence $K_C\simeq (K_X\otimes \mathcal O_X(C))|_C$.
	Taking degrees and using Proposition~\ref{prop:degree-restriction} gives
	\[
	\deg K_C = \deg(K_X|_C)+\deg(\mathcal O_X(C)|_C)=K_X\cdot C + C\cdot C,
	\]
	i.e.\ $2g(C)-2=C\cdot(C+K_X)$.
\end{proof}

\begin{example}[Plane curves in $\mathbb P^2$]
	\label{ex:plane-curve-genus}
	Let $X=\mathbb P^2$ and let $H$ be the class of a line. Then $H\cdot H=1$ and $K_X\sim -3H$.
	A smooth plane curve of degree $d$ has class $C\sim dH$, hence
	\[
	2g(C)-2 \ =\ (dH)\cdot(dH-3H)=d(d-3),
	\qquad\text{so}\qquad
	g(C)=\frac{(d-1)(d-2)}{2}.
	\]
\end{example}

\subsubsection{Riemann--Roch for surfaces: statement and proof mechanisms}
\index{Riemann--Roch theorem!surfaces}
\label{subsec:RRsurf:RR}

\paragraph{The formula.}

\begin{theorem}[Riemann--Roch for line bundles on a smooth surface]
	\label{thm:RR-surface}
	For any divisor $D$ on $X$,
	\[
	\chi(\mathcal O_X(D))=\chi(\mathcal O_X)+\frac12\bigl(D\cdot D - D\cdot K_X\bigr).
	\]
	Equivalently, for $L=\mathcal O_X(D)$,
	\[
	\chi(L)=\chi(\mathcal O_X)+\frac12\bigl(c_1(L)^2-c_1(L)\cdot c_1(K_X)\bigr).
	\]
\end{theorem}

\begin{proof}[Proof by cutting the surface with curves]
	Choose a very ample divisor $H$ so that $D+H$ is also very ample.
	By Bertini's theorem, general members $C\in|H|$ and
	$\Sigma\in|D+H|$ are smooth. Their divisor sequences are
	\begin{gather*}
	0\to\mathcal O_X(D)\to\mathcal O_X(D+H)
	  \to\mathcal O_C(D+H)\to0,\\
	0\to\mathcal O_X\to\mathcal O_X(D+H)
	  \to\mathcal O_{\Sigma}(D+H)\to0.
	\end{gather*}
	Taking Euler characteristics and subtracting gives
	\[
	\chi(\mathcal O_X(D))-\chi(\mathcal O_X)
	 =\chi(\mathcal O_{\Sigma}(D+H))
	  -\chi(\mathcal O_C(D+H)).
	\]
	For any smooth curve $T\subset X$ of divisor class $A$ and any divisor
	class $B$ on $X$, curve Riemann--Roch and adjunction give
	\[
	\chi(\mathcal O_T(B))=B\cdot A+1-g(T)
	 =B\cdot A-\tfrac12 A\cdot(A+K_X).
	\]
	Apply this with $(A,B)=(D+H,D+H)$ for $\Sigma$, and with
	$(A,B)=(H,D+H)$ for $C$. The difference is
	\[
	\tfrac12(D+H)^2-\tfrac12(D+H)\cdot K_X
	 -(D+H)\cdot H+\tfrac12H\cdot(H+K_X)
	 =\tfrac12(D^2-D\cdot K_X).
	\]
	This is the asserted formula. The choice of smooth members and the
	exact-sequence proof may be compared with \cite{Fulton1998Intersection,Hartshorne1977}.
\end{proof}

\paragraph{An alternative route: Hirzebruch--Riemann--Roch in rank $1$.}

A standard proof uses Hirzebruch--Riemann--Roch:
\[
\chi(L)=\int_X \mathrm{ch}(L)\,\mathrm{td}(T_X),
\]
with the Chern character and Todd class defined in
Definitions~\ref{def:chern-character} and \ref{def:todd-class}.
In complex dimension two, these definitions give
\[
\mathrm{ch}(L)=1+c_1(L)+\frac12 c_1(L)^2,
\qquad
\mathrm{td}(T_X)=1+\frac12 c_1(T_X)+\frac{1}{12}\bigl(c_1(T_X)^2+c_2(T_X)\bigr),
\]
and the identity $c_1(T_X)=-c_1(K_X)$.
Only the degree--$4$ component contributes to the integral, giving
\[
\chi(L)=\int_X \left(\frac12 c_1(L)^2 - \frac12 c_1(L)c_1(K_X)\right)+\int_X \frac{1}{12}\bigl(c_1(T_X)^2+c_2(T_X)\bigr).
\]
Applying the same Hirzebruch--Riemann--Roch formula to the trivial bundle
$L=\mathcal O_X$ shows that
\[
\chi(\mathcal O_X)=\int_X \frac{1}{12}
\bigl(c_1(T_X)^2+c_2(T_X)\bigr).
\]
Substitution yields Theorem~\ref{thm:RR-surface} without introducing an additional
identity.

\begin{remark}
	Here Theorem~\ref{thm:RR-surface} will be used as the surface analogue of curve Riemann--Roch,
	and the examples below verify the formula by direct cohomology computations in basic cases.
\end{remark}

\subsubsection{Euler characteristics from exact sequences}
\index{Euler characteristic!coherent sheaves}
\index{Euler sequence}

The surface formula is numerical, but its most efficient computations are
homological. The reason is that an Euler characteristic remembers the
alternating total of a cohomology complex even when the individual groups are
difficult to determine. Thus an exact sequence can replace a bundle by
simpler bundles without losing the number that Riemann--Roch computes.

\begin{definition}[Cohomological Euler characteristic]
	Let $X$ be a projective variety and $\mathcal F$ a coherent sheaf on $X$.
	Its cohomological Euler characteristic is
	\[
	\chi(X,\mathcal F):=\sum_{q\geq0}(-1)^q
	\dim_{\C}H^q(X,\mathcal F).
	\]
	The sum is finite: the cohomology groups are finite-dimensional and vanish
	in degrees above $\dim X$.
\end{definition}

\begin{proposition}[Additivity in a short exact sequence]
	\label{prop:euler-characteristic-additive}
	If $0\to\mathcal F'\to\mathcal F\to\mathcal F''\to0$ is a short exact
	sequence of coherent sheaves on a projective variety, then
	\[
	\chi(X,\mathcal F)=\chi(X,\mathcal F')+\chi(X,\mathcal F'').
	\]
\end{proposition}

\begin{proof}
	Sheaf cohomology supplies a finite long exact sequence
	\[
	0\to H^0(\mathcal F')\to H^0(\mathcal F)\to H^0(\mathcal F'')
	\to H^1(\mathcal F')\to\cdots .
	\]
	For any finite exact sequence of finite-dimensional vector spaces, the
	alternating sum of the dimensions is zero. Grouping the terms according to
	their cohomological degree gives the displayed identity. This is the
	homological-algebra mechanism behind the divisor-sequence proof of
	Theorem~\ref{thm:RR-surface}.
\end{proof}

On projective space the basic resolution is the Euler sequence
\begin{equation}
	\label{eq:euler-sequence-pn}
	0\longrightarrow\Omega^1_{\PP^n}(1)
	\longrightarrow\mathcal O_{\PP^n}^{\oplus(n+1)}
	\xrightarrow{(X_0,\ldots,X_n)}\mathcal O_{\PP^n}(1)
	\longrightarrow0.
\end{equation}
At $[X_0:\cdots:X_n]$, the last map evaluates a tuple against the
homogeneous coordinates; its kernel is the cotangent space with the indicated
twist. Dualizing gives
\begin{equation}
	\label{eq:dual-euler-sequence-pn}
	0\longrightarrow\mathcal O_{\PP^n}
	\longrightarrow\mathcal O_{\PP^n}(1)^{\oplus(n+1)}
	\longrightarrow T_{\PP^n}\longrightarrow0.
\end{equation}
The Whitney product formula therefore yields
$c(T_{\PP^n})=(1+H)^{n+1}$, with terms above degree $n$ discarded. On
$\PP^2$ this says
\[
	c_1(T_{\PP^2})=3H,\qquad c_2(T_{\PP^2})=3H^2,
	\qquad K_{\PP^2}=-3H.
\]
Since $\int_{\PP^2}H^2=1$, Chern--Gauss--Bonnet gives
$\chi_{\mathrm{top}}(\PP^2)=\int_{\PP^2}c_2=3$, matching its three even
cells. This explains geometrically why the coefficient $3$ appears in both
the canonical class and plane-curve adjunction.

Twisting \eqref{eq:euler-sequence-pn} on $\PP^2$ by
$\mathcal O(k-1)$ gives
\[
	0\longrightarrow\Omega^1_{\PP^2}(k)
	\longrightarrow\mathcal O_{\PP^2}(k-1)^{\oplus3}
	\longrightarrow\mathcal O_{\PP^2}(k)\longrightarrow0.
\]
Proposition~\ref{prop:euler-characteristic-additive} and
$\chi(\mathcal O_{\PP^2}(m))=\binom{m+2}{2}$, interpreted as a polynomial
for every integer $m$, give the representative calculation
\begin{align*}
	\chi\bigl(\Omega^1_{\PP^2}(k)\bigr)
	&=3\chi\bigl(\mathcal O(k-1)\bigr)-\chi\bigl(\mathcal O(k)\bigr)\\
	&=3\binom{k+1}{2}-\binom{k+2}{2}=k^2-1.
\end{align*}
Likewise, the dual Euler sequence gives
$\chi(T_{\PP^2})=3\chi(\mathcal O(1))-\chi(\mathcal O)=8$.
These are Euler characteristics, not assertions that all higher cohomology
vanishes: the exact sequence controls the alternating sum even when its
summands occur in different degrees. This is a principal advantage of the
homological formulation \cite{Hartshorne1977,GriffithsHarris1994}.

\subsection{Computations, vanishing, and effectivity}
\label{subsec:RRsurf:examples}

\subsubsection{Example A: $\mathbb P^2$}

Let $X=\mathbb P^2$ and let $H$ be the hyperplane class.
Then $\mathrm{Pic}(X)\cong \mathbb Z[H]$ and $H^2=H\cdot H=1$.
Also $K_X\sim -3H$, and $\chi(\mathcal O_X)=1$.
The Chern-character/Todd computation in
\S\ref{subsec:chern-weil-character-todd} predicted the same answer;
the calculation here checks it using the independently proved surface
Riemann--Roch formula.

\begin{proposition}[RR for $\mathbb P^2$ in closed form]
	\label{prop:RR-P2}
	For $D=dH$,
	\[
	\chi(\mathcal O_{\mathbb P^2}(d)) \ =\ 1+\frac12\bigl(d^2+3d\bigr)
	\ =\ \frac{(d+1)(d+2)}{2}.
	\]
\end{proposition}

\begin{proof}
	Compute the intersection terms:
	\[
	D\cdot D=d^2(H^2)=d^2,
	\qquad
	D\cdot K_X=(dH)\cdot(-3H)=-3d.
	\]
	Insert into surface Riemann--Roch (Theorem~\ref{thm:RR-surface}):
	\[
	\chi(\mathcal O_X(dH))=\chi(\mathcal O_X)+\frac12(d^2-(-3d))
	=1+\frac12(d^2+3d)=\frac{(d+1)(d+2)}{2}.
	\]
\end{proof}

\begin{example}[Direct cohomology check for $d\ge 0$]
	For $d\ge 0$, one knows $H^i(\mathbb P^2,\mathcal O(d))=0$ for $i>0$ and
	$h^0(\mathbb P^2,\mathcal O(d))=\binom{d+2}{2}$, so $\chi=\binom{d+2}{2}$, matching Proposition~\ref{prop:RR-P2}.
\end{example}

\subsubsection{Example B: $\mathbb P^1\times \mathbb P^1$}

Let $X=\mathbb P^1\times\mathbb P^1$.
Let $F_1=[\mathbb P^1\times \{\mathrm{pt}\}]$ and $F_2=[\{\mathrm{pt}\}\times \mathbb P^1]$.
Then $\mathrm{Pic}(X)\cong \mathbb Z F_1\oplus \mathbb Z F_2$ and intersection numbers are
\[
F_1^2=0,\qquad F_2^2=0,\qquad F_1\cdot F_2=1.
\]
Also $K_X\sim -2F_1-2F_2$ and $\chi(\mathcal O_X)=1$.

\begin{proposition}[RR on $\mathbb P^1\times\mathbb P^1$]
	\label{prop:RR-P1P1}
	For $D=aF_1+bF_2$,
	\[
	\chi(\mathcal O_X(a,b)) \ =\ (a+1)(b+1).
	\]
\end{proposition}

\begin{proof}
	Compute
	\[
	D^2=(aF_1+bF_2)^2=2ab(F_1\cdot F_2)=2ab,
	\]
	and
	\[
	D\cdot K_X=(aF_1+bF_2)\cdot(-2F_1-2F_2)=-2a(F_1\cdot F_2)-2b(F_1\cdot F_2)=-2a-2b.
	\]
	Hence Theorem~\ref{thm:RR-surface} gives
	\[
	\chi(\mathcal O_X(D))=1+\frac12\bigl(2ab-(-2a-2b)\bigr)=1+ab+a+b=(a+1)(b+1).
	\]
\end{proof}

\begin{example}[Direct cohomology check when $a,b\ge 0$]
	If $a,b\ge 0$, then $H^i(X,\mathcal O(a,b))=0$ for $i>0$ and
	$h^0(X,\mathcal O(a,b))=(a+1)(b+1)$ by K\"unneth and the curve computation on $\mathbb P^1$.
\end{example}

\subsubsection{Example C: Hirzebruch surfaces $\mathbb F_n$}

Let $X=\mathbb F_n=\mathbb P(\mathcal O_{\mathbb P^1}\oplus \mathcal O_{\mathbb P^1}(n))$.
Let $f$ be the fiber class and $s$ the negative section class.
Then
\[
f^2=0,\qquad s\cdot f=1,\qquad s^2=-n.
\]
The canonical class is
\[
K_X \ \sim\ -2s-(n+2)f,
\qquad\text{and}\qquad
\chi(\mathcal O_X)=1.
\]

\begin{proposition}[RR on $\mathbb F_n$ in coordinates]
	\label{prop:RR-Fn}
	For $D=as+bf$,
	\[
	\chi(\mathcal O_X(D)) \ =\ 1+\frac12\bigl(-na^2+2ab + (n-2)a+2b\bigr).
	\]
\end{proposition}

\begin{proof}
	Compute
	\[
	D^2=(as+bf)^2=a^2s^2+2ab(s\cdot f)+b^2 f^2=-na^2+2ab,
	\]
	and
	\[
	D\cdot K_X=(as+bf)\cdot(-2s-(n+2)f)
	=-2a s^2-a(n+2)(s\cdot f)-2b(f\cdot s)
	=2an-a(n+2)-2b.
	\]
	Thus
	\[
	D^2-D\cdot K_X
	=(-na^2+2ab)-\bigl(2an-a(n+2)-2b\bigr)
	=-na^2+2ab+(n-2)a+2b.
	\]
	Insert into Theorem~\ref{thm:RR-surface} and use $\chi(\mathcal O_X)=1$.
\end{proof}

\begin{example}[Genus of a divisor class on $\mathbb F_n$]
	If $C\sim as+bf$ is smooth, then adjunction (Theorem~\ref{thm:adjunction}) gives
	\[
	2g(C)-2=C\cdot(C+K_X).
	\]
	Using the intersection table above, you can compute $g(C)$ explicitly as a polynomial in $a,b,n$.
\end{example}

\subsubsection{Example D: Blow-up of $\mathbb P^2$ at a point}

Let $\pi:\widetilde X\to \mathbb P^2$ be the blow-up at one point.
Let $H=\pi^*\mathcal O_{\mathbb P^2}(1)$ and let $E$ be the exceptional curve.
Then
\[
H^2=1,\qquad E^2=-1,\qquad H\cdot E=0,
\qquad
K_{\widetilde X}\sim -3H+E,
\qquad
\chi(\mathcal O_{\widetilde X})=1.
\]

\begin{proposition}[RR on the blow-up of $\mathbb P^2$ at one point]
	\label{prop:RR-blowup}
	For $D=aH-bE$,
	\[
	\chi(\mathcal O_{\widetilde X}(D)) \ =\ 1+\frac12\bigl(a^2+3a-b^2-b\bigr).
	\]
\end{proposition}

\begin{proof}
	Compute
	\[
	D^2=(aH-bE)^2=a^2H^2-2ab(H\cdot E)+b^2E^2=a^2-b^2,
	\]
	and
	\[
	D\cdot K_{\widetilde X}=(aH-bE)\cdot(-3H+E)=-3aH^2-bE^2=-3a+b.
	\]
	Thus
	\[
	\chi(\mathcal O_{\widetilde X}(D))
	=1+\frac12\bigl((a^2-b^2)-(-3a+b)\bigr)
	=1+\frac12(a^2+3a-b^2-b).
	\]
\end{proof}

\subsubsection{Using RR + Serre duality: vanishing tests and effectivity heuristics}
\label{subsec:RRsurf:vanishing}

Surface Riemann--Roch controls $\chi$, but to extract $h^0$ one needs information on $h^1$ and $h^2$.
Serre duality gives
\[
h^2(\mathcal O_X(D))=h^0(\mathcal O_X(K_X-D)).
\]
Thus:
\begin{itemize}
	\item If $K_X-D$ has no sections (e.g.\ is ``too negative''), then $h^2(\mathcal O_X(D))=0$.
	\item If additionally $h^1(\mathcal O_X(D))=0$ (often by positivity/vanishing theorems in advanced settings),
	then $h^0(\mathcal O_X(D))=\chi(\mathcal O_X(D))$.
\end{itemize}

In our elementary examples ($\mathbb P^2$, $\mathbb P^1\times\mathbb P^1$, many $\mathbb F_n$ cases),
one can often determine $h^0$ by explicit section counting in charts and show $h^1=h^2=0$ in the ``positive'' region,
thereby recovering concrete dimension formulas.

\subsection{Exercises}
\label{subsec:RRsurf:exercises}

\begin{exercise}[Intersection table on $\mathbb P^1\times\mathbb P^1$]
	Let $X=\mathbb P^1\times\mathbb P^1$ with $F_1,F_2$ as above.
	\begin{enumerate}
		\item Prove $F_1^2=F_2^2=0$ and $F_1\cdot F_2=1$ using geometric representatives.
		\item For $D=aF_1+bF_2$ and $E=cF_1+dF_2$, compute $D\cdot E$ explicitly.
	\end{enumerate}
\end{exercise}

\begin{exercise}[Genus of a smooth divisor on $\mathbb P^1\times\mathbb P^1$]
	Let $C\sim aF_1+bF_2$ be a smooth connected curve on
	$X=\mathbb P^1\times\mathbb P^1$, where $a,b\ge1$.
	Use adjunction with $K_X\sim -2F_1-2F_2$ to show
	\[
	g(C)=(a-1)(b-1).
	\]
	Check this for $(a,b)=(1,1)$ and $(2,1)$ by direct geometric descriptions.
\end{exercise}

\begin{exercise}[Plane curves and adjunction]
	Let $X=\mathbb P^2$ and let $C$ be a smooth curve of degree $d$.
	\begin{enumerate}
		\item Prove $K_{\mathbb P^2}\sim -3H$ using the transformation rule for a holomorphic $2$--form under affine charts.
		\item Deduce $g(C)=\frac{(d-1)(d-2)}{2}$.
	\end{enumerate}
\end{exercise}

\begin{exercise}[Local intersection multiplicity via elimination]
	In $\mathbb C^2$ consider $C=\{y^2-x^3=0\}$ and the line
	$C'_\lambda=\{y-\lambda x=0\}$ with $\lambda\in\mathbb C$.
	Compute $I_0(C,C'_\lambda)$ as a function of $\lambda$ by eliminating $y$.
	Show that the answer is $2$ for $\lambda\ne0$ and $3$ for $\lambda=0$, and
	interpret the jump using the tangent line of the cusp.
\end{exercise}

\begin{exercise}[Hirzebruch surface canonical class]
	Let $X=\mathbb F_n$ with generators $s,f$.
	\begin{enumerate}
		\item Using adjunction on the section $s$, show that $K_X\cdot s = n-2$.
		\item Using adjunction on a fiber $f\cong \mathbb P^1$, show that $K_X\cdot f=-2$.
		\item Deduce $K_X\sim -2s-(n+2)f$ from the intersection table.
	\end{enumerate}
\end{exercise}

\begin{exercise}[RR on $\mathbb F_n$ for a concrete region]
	Let $X=\mathbb F_n$ and $D=as+bf$ with $a\ge 0$ and $b\gg 0$.
	\begin{enumerate}
		\item Use Proposition~\ref{prop:RR-Fn} to compute $\chi(\mathcal O_X(D))$.
		\item Show $h^2(\mathcal O_X(D))=h^0(\mathcal O_X(K_X-D))=0$ for $b$ sufficiently large.
		\item (Optional) By analyzing sections along the ruling, argue that $h^1(\mathcal O_X(D))=0$ for $b$ sufficiently large,
		and conclude $h^0(\mathcal O_X(D))=\chi(\mathcal O_X(D))$ in that range.
	\end{enumerate}
\end{exercise}

\begin{exercise}[Blow-up intersection calculus]
	Let $\pi:\widetilde X\to \mathbb P^2$ be the blow-up at one point with classes $H,E$ as above.
	\begin{enumerate}
		\item Prove $H\cdot E=0$ geometrically (use representatives avoiding the blown-up point).
		\item Show $E^2=-1$ by identifying $E$ with $\mathbb P^1$ and computing the degree of $\mathcal O_{\widetilde X}(E)|_E$.
		\item Compute $K_{\widetilde X}$ from $K_{\mathbb P^2}$ and the exceptional divisor and verify $K_{\widetilde X}\sim -3H+E$.
	\end{enumerate}
\end{exercise}

\begin{exercise}[A Riemann--Roch dimension prediction]
	On $X=\mathbb P^2$, use Proposition~\ref{prop:RR-P2} and Serre duality to predict
	$h^0(\mathcal O_{\mathbb P^2}(-1))$ and $h^2(\mathcal O_{\mathbb P^2}(-4))$.
	Then verify your predictions by a direct argument (no heavy theorems needed).
\end{exercise}

\begin{exercise}[Intersection form as cup product (with Poincar\'e duality)]
	\label{ex:intersection-cup-product-PD}
	Let $X$ be a compact oriented smooth $4$--manifold (in particular, a smooth complex surface),
	and let $a,b\in H^2(X,\mathbb Z)$.  Choose closed real $2$--forms $\alpha,\beta$
	representing their images in $H^2_{\mathrm{dR}}(X)$.
	Recall Definition~\ref{def:PD-cycle}: for a closed oriented smooth
	surface $Z\subset X$, the integral class $\mathrm{PD}([Z])$ is defined
	by Poincar\'e duality, and its real image satisfies
	$\int_X\eta\wedge\zeta_Z=\int_Z i^*\eta$ for closed $2$-forms
	$\eta$, where $\zeta_Z$ represents that real image.
	
	\medskip
	\noindent
	\begin{enumerate}
		\item
		Show that $\displaystyle \int_X \alpha\wedge \beta\in\mathbb Z$ and is
		independent of the chosen closed representatives.
		
		\smallskip
		\emph{Hint:} Interpret $[\alpha\wedge\beta]$ as the real image of the cup product
		$a\smile b\in H^4(X,\mathbb Z)$ and pair it with the
		fundamental class $[X]$.
		
		\item
		Assume $a=\mathrm{PD}([C])$ and $b=\mathrm{PD}([C'])$, where
		$C,C'\subset X$ are smooth oriented curves (real dimension $2$) that intersect
		transversely.
		Show that
		\[
		\int_X \alpha\wedge\beta
		\ =\ 
		\#(C\cap C'),
		\]
		the (signed) number of intersection points of $C$ and $C'$.
		
		\smallskip
		\emph{Hint:} Choose closed representatives $\zeta_C,\zeta_{C'}$ of the two
		Poincar\'e-dual classes and use the defining property twice:
		\[
		\int_X \zeta_C\wedge\zeta_{C'}
		\ =\ 
		\int_C \zeta_{C'}
		\ =\ 
		\int_{C\cap C'} 1.
		\]
		Explain why transversality makes each local contribution $+1$ or $-1$
		according to the two orientations. If $X$ is a complex surface and $C,C'$ are
		complex curves with their complex orientations, explain why every contribution is $+1$.
	\end{enumerate}
\end{exercise}

\begin{exercise}[Adjunction and self-intersection]
	Let $C\subset X$ be a smooth curve.
	\begin{enumerate}
		\item Show that $C^2=\deg(\mathcal O_X(C)|_C)$.
		\item Use adjunction to express $C^2$ in terms of $g(C)$ and $K_X\cdot C$.
	\end{enumerate}
\end{exercise}

\subsubsection{Summary}
On a surface, intersection theory provides the quadratic correction term that replaces the linear ``degree''
term on a curve, and Riemann--Roch becomes
\[
\chi(\mathcal O_X(D))=\chi(\mathcal O_X)+\frac12\bigl(D^2-D\cdot K_X\bigr).
\]
	
\clearpage
\part{Algebraic Curves and Algebraic Geometry}
\label{part:jacobian}
\label{part:alggeo}

Throughout this part \(X\) denotes a compact Riemann surface of genus \(g\ge 0\).
We write \(\Omega_X^1:=K_X\) for the canonical bundle (holomorphic \(1\)-forms),
and \(H^0(X,\Omega_X^1)\) for the \(g\)-dimensional space of holomorphic differentials.
The period-lattice, Picard, Abel, and Jacobi inversion arguments follow the standard treatments
\cite{FarkasKra1992,Forster1981,ArbarelloCornalbaGriffithsHarris1985}.

\section{Jacobians of Compact Riemann Surfaces and Abel--Jacobi Theory}
\index{Jacobian}\index{Picard variety}\index{Abel--Jacobi map}
\label{sec:jacobian-picard-abel-jacobi}

Integrating holomorphic differentials between two points depends on
the chosen path, but changing the path changes the integral by a
period.  Passing to the quotient by all periods therefore turns these
integrals into well-defined points of the Jacobian.  The purpose of
this construction is not only to obtain a complex torus: Abel's
theorem uses it to test when two divisors are linearly equivalent,
while Jacobi inversion asks which torus points come from effective
divisors.  The identification with $\mathrm{Pic}^0(X)$ later shows
that the same torus classifies degree-zero holomorphic line bundles.

The discussion is organized in three long steps rather than a sequence of
short, nearly independent topics. First we construct the torus from periods
and identify it cohomologically with $\mathrm{Pic}^0(X)$. Second we send
divisors to that torus and study the fibers by complete linear systems,
Abel's theorem, and Jacobi inversion. Finally, the sphere and a complex torus
serve as boundary cases: on the sphere the Jacobian disappears, whereas in
genus one the curve is its own Jacobian. A useful question to keep in mind is
therefore not merely ``what is $J(X)$?'' but ``which divisor data become
linear after passing to $J(X)$?''
We begin with the period lattice because it is the common language
behind all three descriptions.

\paragraph{A genus-zero check.}
Recall that the dimension of $H^0(X,K_X)$ is the genus. On
$X=\PP^1$ this space is zero, so there are no period integrals to
record and $J(\PP^1)$ is a single point. The comparison with
$\mathrm{Pic}^0(\PP^1)$ is already instructive: every degree-zero
divisor on $\PP^1$ is principal. For positive genus the period lattice
measures precisely what this genus-zero calculation can no longer see.

\subsection{The period lattice and the Abel--Jacobi map}
\index{period lattice}\index{Jacobian!periods}
\label{subsec:period-lattice}

\subsubsection{Integration pairing and periods}

Let \(\omega\in H^0(X,\Omega_X^1)\) and \(\gamma\in H_1(X,\Z)\). Since \(\omega\) is closed,
\(\int_\gamma \omega\) depends only on the homology class \([\gamma]\), hence defines a bilinear pairing
\begin{equation}\label{eq:period-pairing}
	H^0(X,\Omega_X^1)\times H_1(X,\Z)\longrightarrow \C,
	\qquad
	(\omega,\gamma)\longmapsto \int_\gamma \omega.
\end{equation}

\begin{lemma}[Period functionals are well-defined]
	\label{lem:period-functional-well-defined}
	If \(\gamma\) is a boundary in \(X\), then \(\int_\gamma \omega=0\) for every \(\omega\in H^0(X,\Omega_X^1)\).
\end{lemma}

\begin{proof}
	If \(\gamma=\partial \Sigma\) for a smooth \(2\)-chain \(\Sigma\), then by Stokes' theorem,
	\(\int_\gamma \omega=\int_{\partial\Sigma}\omega=\int_\Sigma d\omega=0\), since \(\omega\) is holomorphic and hence closed.
\end{proof}

\subsubsection{Definition of the period lattice and the Jacobian}

Choose a basis \(\{\omega_1,\dots,\omega_g\}\) of \(H^0(X,\Omega_X^1)\).
Define the \emph{period map}
\[
\mathsf P: H_1(X,\Z)\longrightarrow \C^g,
\qquad
\gamma\longmapsto \Big(\int_\gamma\omega_1,\dots,\int_\gamma\omega_g\Big).
\]
Its image is the \emph{period lattice}
\begin{equation}\label{eq:period-lattice}
	\Lambda:=\mathsf P\big(H_1(X,\Z)\big)\subset \C^g.
\end{equation}

\begin{proposition}[The period lattice is a full lattice]
	\label{prop:period-lattice-full}
	The subgroup \(\Lambda\subset \C^g\) is discrete of rank \(2g\). Equivalently, \(\Lambda\) is a full lattice in \(\C^g\).
\end{proposition}

\begin{proof}
	Hodge theory identifies the real and imaginary parts of holomorphic differentials
	with the harmonic representatives of $H^1_{\mathrm{dR}}(X;\R)$. Thus the $2g$
	real forms
	\[
	\Re\omega_1,\Im\omega_1,\ldots,\Re\omega_g,\Im\omega_g
	\]
	span $H^1_{\mathrm{dR}}(X;\R)$. If a real homology class $c$ has zero
	period against every $\omega_j$, it pairs to zero with all of
	$H^1_{\mathrm{dR}}(X;\R)$. Nondegeneracy of the de Rham pairing then gives
	$c=0$. Hence the real-linear period map
	\[
	H_1(X;\R)\longrightarrow\C^g
	\]
	is injective. Both sides have real dimension $2g$, so it is an isomorphism.
	The subgroup $H_1(X;\Z)$ is a lattice of rank $2g$ in $H_1(X;\R)$;
	its image under a real-linear isomorphism is therefore a discrete lattice of
	rank $2g$ in $\C^g$. This image is $\Lambda$.
\end{proof}

\begin{definition}[Jacobian]
	\label{def:jacobian}
	The \emph{Jacobian} of \(X\) is the complex torus
	\begin{equation}\label{eq:jacobian-as-torus}
		J(X):=\C^g/\Lambda.
	\end{equation}
\end{definition}

\begin{remark}
	Changing the basis \(\{\omega_j\}\) applies an invertible complex linear map
	to \(\Lambda\), and hence produces an isomorphic complex torus.  Intrinsically,
	\[
	J(X)=H^0(X,K_X)^\vee/H_1(X,\Z),
	\]
	where a cycle acts by integration.  The torus is therefore independent of the
	chosen basis; the basis only supplies coordinates in \(\C^g\).
\end{remark}

\begin{exercise}[The genus-zero Jacobian]
	\label{exr:g0-jacobian}
	Show that if \(g=0\) then \(H^0(X,\Omega_X^1)=0\), hence \(\C^g=0\) and \(J(X)\) is the trivial group.
\end{exercise}

\subsubsection{Abel--Jacobi map on points and on divisors}
\index{Abel--Jacobi map!divisors}
\label{subsec:abel-jacobi}

\subsubsection{The Abel--Jacobi map on points}

Fix a base point \(p_0\in X\). For \(x\in X\) choose a path \(\Gamma\) from \(p_0\) to \(x\) and set
\[
\widetilde\alpha(x):=\Big(\int_{\Gamma}\omega_1,\dots,\int_{\Gamma}\omega_g\Big)\in \C^g.
\]

\begin{figure}[ht]
	\centering
	\begin{tikzpicture}[scale=1.0,>=Stealth]
		\draw[thick] (-3.15,-1.1) .. controls (-3.5,0.7) and (-2.0,1.35) .. (-0.5,1.05)
			.. controls (1.0,1.35) and (3.2,0.75) .. (3.1,-0.55)
			.. controls (2.5,-1.45) and (-1.5,-1.4) .. cycle;
		\coordinate (P) at (-2.25,-0.20);
		\coordinate (Q) at (2.20,-0.20);
		\draw[->,thick,blue!70!black] (P) .. controls (-1.0,0.78) and (0.9,0.78) .. (Q)
			node[midway,above] {$\Gamma$};
		\draw[->,thick,red!70!black] (P) .. controls (-1.0,-1.0) and (0.9,-1.0) .. (Q)
			node[midway,below] {$\Gamma'$};
		\fill (P) circle (1.8pt) node[left] {$p_0$};
		\fill (Q) circle (1.8pt) node[right] {$x$};
	\end{tikzpicture}
	\caption{Two paths from $p_0$ to $x$ differ by a closed loop.
	Their integral vectors differ by the period of that loop, so they define
	the same point of the Jacobian after quotienting by the period lattice.}
	\label{fig:AJ-two-paths}
\end{figure}

\begin{lemma}[Well-definedness modulo periods]
	\label{lem:AJ-point-well-defined}
	If \(\Gamma,\Gamma'\) are two paths from \(p_0\) to \(x\), then
	\(\widetilde\alpha_{\Gamma}(x)-\widetilde\alpha_{\Gamma'}(x)\in \Lambda\).
\end{lemma}

\begin{proof}
	The concatenation \(\Gamma\cdot(\Gamma')^{-1}\) is a loop based at \(p_0\), hence determines an integral homology class
	\(\gamma\in H_1(X,\Z)\). Therefore
	\[
	\int_{\Gamma}\omega_j-\int_{\Gamma'}\omega_j=\int_{\Gamma\cdot(\Gamma')^{-1}}\omega_j=\int_{\gamma}\omega_j,
	\]
	so the difference vector lies in \(\Lambda\) by definition \eqref{eq:period-lattice}.
\end{proof}

\begin{definition}[Abel--Jacobi map on points]
	\label{def:AJ-point}
	The \emph{Abel--Jacobi map} \(\alpha:X\to J(X)\) is
	\[
	\alpha(x):=\widetilde\alpha(x)\bmod\Lambda\in \C^g/\Lambda.
	\]
\end{definition}

\subsubsection{Symmetric powers and effective divisors}

Let \(S^dX:=X^d/\mathfrak S_d\) be the \(d\)-fold symmetric product. Points of \(S^dX\) correspond to effective divisors
\(D=x_1+\cdots+x_d\) of degree \(d\) (unordered with multiplicities).

\begin{definition}[Abel--Jacobi map on \(S^dX\)]
	\label{def:AJ-symmetric}
	Define \(\alpha_d:S^dX\to J(X)\) by
	\[
	\alpha_d\!\Big(\sum_{i=1}^d x_i\Big):=\sum_{i=1}^d \alpha(x_i)\in J(X).
	\]
\end{definition}

\begin{lemma}[Compatibility with addition of a base point]
	\label{lem:AJ-compatibility}
	Let \(i_{d}:S^dX\to S^{d+1}X\) be \(D\mapsto D+p_0\). Then \(\alpha_{d+1}\circ i_{d}=\alpha_d\).
\end{lemma}

\begin{proof}
	Since \(\alpha(p_0)=0\) by definition, \(\alpha_{d+1}(D+p_0)=\alpha_d(D)+\alpha(p_0)=\alpha_d(D)\).
\end{proof}

\begin{exercise}[The Abel--Jacobi map of an elliptic curve]
	\label{exr:AJ-g1}
	Let \(X=\C/\Lambda\) be an elliptic curve. Show that choosing \(\omega=dz\) gives
	\(J(X)\cong \C/\Lambda\cong X\). Prove that \(\alpha:X\to J(X)\) is a group isomorphism once \(p_0\) is taken as the origin.
\end{exercise}

\subsubsection{The Picard group and the Jacobian}
\label{subsec:picard-exp}

The period construction is analytic: it begins with integrals of
holomorphic one-forms. The Picard construction is geometric: it classifies
line bundles by transition functions. Their identification is the central
structural point of this subsection. The exponential sequence converts
multiplicative transition functions into additive cohomology, and Serre
duality converts that cohomology into linear functionals on holomorphic
differentials. The integral cohomology classes then become exactly the period
lattice. Thus the quotient defining the Jacobian is forced by the
classification of degree-zero line bundles rather than introduced as an
unmotivated analytic construction.

\paragraph{$\mathrm{Pic}(X)$ as a cohomology group.}

Recall \(\mathcal O_X^{\!*}\) is the sheaf of nowhere-vanishing holomorphic functions.
The cocycle construction proved in Theorem~\ref{thm:picard-cech} gives
\begin{equation}\label{eq:Pic-as-H1}
	\mathrm{Pic}(X)\ \cong\ H^1(X,\mathcal O_X^{\!*}),
\end{equation}
via Čech cocycles: transition functions of a holomorphic line bundle form a \(1\)-cocycle with values in \(\mathcal O_X^{\!*}\),
and changing trivializations changes the cocycle by a coboundary.

\paragraph{The exponential exact sequence and the first Chern class.}

Consider the exponential short exact sequence of sheaves on \(X\):
\begin{equation}\label{eq:exp-seq}
	0\longrightarrow \Z \longrightarrow \mathcal O_X
	\xrightarrow{\ \exp(2\pi i\,\cdot)\ } \mathcal O_X^{\!*}\longrightarrow 0.
\end{equation}
Taking cohomology yields a long exact sequence, whose relevant part is
\begin{equation}\label{eq:LES-exp}
	H^1(X,\mathcal O_X)\xrightarrow{\ \beta\ } H^1(X,\mathcal O_X^{\!*})
	\xrightarrow{\ c_1\ } H^2(X,\Z)\xrightarrow{\ }\ H^2(X,\mathcal O_X).
\end{equation}

\begin{lemma}[Vanishing \(H^2(X,\mathcal O_X)=0\)]
	\label{lem:H2OX-vanish}
	For a compact Riemann surface \(X\), \(H^2(X,\mathcal O_X)=0\).
\end{lemma}

\begin{proof}
	As a complex manifold of complex dimension \(1\), the Dolbeault model gives
	\(H^2(X,\mathcal O_X)\cong H^{0,2}(X)\). But there are no \((0,2)\)-forms on a curve, hence \(H^{0,2}(X)=0\).
\end{proof}

\begin{proposition}[\(\mathrm{Pic}^0(X)\) as a quotient]
	\label{prop:Pic0-quotient}
	Let \(c_1:\mathrm{Pic}(X)\to H^2(X,\Z)\cong \Z\) be the degree/first Chern class map.
	Then
	\[
	\mathrm{Pic}^0(X):=\ker(c_1)\ \cong\ H^1(X,\mathcal O_X)\big/\operatorname{Im}\alpha,
	\]
	where \(\alpha:H^1(X,\Z)\to H^1(X,\mathcal O_X)\) is induced by \(\Z\hookrightarrow \mathcal O_X\).
\end{proposition}

\begin{proof}
	Use \eqref{eq:LES-exp}. By Lemma~\ref{lem:H2OX-vanish}, the map \(H^2(X,\Z)\to H^2(X,\mathcal O_X)\) is the zero map,
	so exactness at \(H^1(X,\mathcal O_X^{\!*})\) gives
	\[
	\ker(c_1)=\operatorname{Im}(\beta).
	\]
	Thus \(\mathrm{Pic}^0(X)\cong \operatorname{Im}(\beta)\).
	By the first isomorphism theorem, \(\operatorname{Im}(\beta)\cong H^1(X,\mathcal O_X)/\ker(\beta)\).
	Exactness at \(H^1(X,\mathcal O_X)\) yields \(\ker(\beta)=\operatorname{Im}(\alpha)\).
	Combining these gives the claimed quotient.
\end{proof}

\begin{exercise}[Degree as the first Chern class]
	\label{exr:deg-is-c1}
	Show that \(H^2(X,\Z)\cong \Z\) canonically by orientation, and that for a divisor \(D\),
	the class \(c_1(\mathcal O_X(D))\in H^2(X,\Z)\) equals \(\deg D\).
\end{exercise}

\paragraph{Identifying $\mathrm{Pic}^0(X)$ with $J(X)$.}
\label{subsec:Pic0-equals-J}

\paragraph{The genus-one model.}
Before taking the general quotient, consider
$X=\C/(\Z+\tau\Z)$ with $\Im\tau>0$ and base point $0$.  Its
holomorphic differentials are multiples of $dz$, and integration over
the two basic loops gives periods $1$ and $\tau$.  Consequently
$J(X)=\C/(\Z+\tau\Z)$.  A point $p$, represented by $z_p\in\C$,
has Abel--Jacobi coordinate $\int_0^p dz=z_p$ modulo those periods.
The correspondence with degree-zero line bundles sends it to
$\mathcal O_X(p-0)$.  Moving $z_p$ by $1$ or $\tau$ does not change
the point of $X$ or its divisor class.  The general proof below
recovers exactly this quotient without a preferred coordinate:
$H^1(X,\mathcal O_X)$ replaces the complex plane, and the image of
$H^1(X,\Z)$ replaces the lattice generated by $1$ and $\tau$.

\paragraph{Serre duality and the dual of holomorphic differentials.}

Serre duality on a curve gives a perfect pairing
\[
H^1(X,\mathcal O_X)\times H^0(X,K_X)\longrightarrow \C,
\qquad
([\eta],\omega)\longmapsto \int_X \omega\wedge \eta,
\]
hence an isomorphism
\begin{equation}\label{eq:Serre-H1OX}
	H^1(X,\mathcal O_X)\ \cong\ H^0(X,K_X)^\vee\ =\ H^0(X,\Omega_X^1)^\vee.
\end{equation}

\paragraph{The image of $H^1(X,\Z)$ becomes the period lattice.}

Let \(\alpha:H^1(X,\Z)\to H^1(X,\mathcal O_X)\) be induced by \(\Z\hookrightarrow \mathcal O_X\).
We now explain why \(\operatorname{Im}(\alpha)\) corresponds to periods.

\begin{proposition}[Integral classes give period functionals]
	\label{prop:integral-to-period}
	Under the Serre duality identification \eqref{eq:Serre-H1OX}, the subgroup \(\operatorname{Im}(\alpha)\subset H^1(X,\mathcal O_X)\)
	corresponds to the subgroup
	\[
	H_1(X,\Z)\ \hookrightarrow\ H^0(X,\Omega_X^1)^\vee,\qquad
	\gamma\longmapsto\big(\omega\mapsto \textstyle\int_\gamma\omega\big).
	\]
\end{proposition}

\begin{proof}
	Use the integral de Rham comparison. Represent \(u\in H^1(X,\Z)\) by a closed real one-form
	\(a\) with integral periods. Its image in \(H^1(X,\mathcal O_X)\cong H^{0,1}(X)\) is represented
	by \(a^{0,1}\). For a holomorphic one-form \(\omega\), type considerations give
	\[
	\int_X\omega\wedge a^{0,1}=\int_X\omega\wedge a,
	\]
	because \(\omega\wedge a^{1,0}\) has type \((2,0)\) and vanishes on a curve. Let
	\(\gamma=D_X(u)\in H_1(X,\Z)\), using the convention of
	Definition~\ref{def:PD-cycle}.  In particular,
	\(\int_X\beta\wedge a=\int_\gamma\beta\) for closed one-forms \(\beta\).
	Then the Serre
	pairing of \(\omega\) with \(\alpha(u)\) equals \(\int_\gamma\omega\). Thus the integral lattice
	maps exactly to the lattice of period functionals (with a simultaneous sign change if the opposite
	Poincar\'e-duality convention is chosen).
\end{proof}

\begin{theorem}[\(\mathrm{Pic}^0(X)\cong J(X)\)]
	\label{thm:Pic0-is-J}
	There is a canonical isomorphism of complex tori
	\[
	\mathrm{Pic}^0(X)\ \cong\ \frac{H^1(X,\mathcal O_X)}{\operatorname{Im}(\alpha)}
	\ \cong\ \frac{H^0(X,\Omega_X^1)^\vee}{H_1(X,\Z)}
	\ \cong\ \C^g/\Lambda
	\ =\ J(X).
	\]
\end{theorem}

\begin{proof}
	Combine Proposition~\ref{prop:Pic0-quotient} with Serre duality \eqref{eq:Serre-H1OX} and Proposition~\ref{prop:integral-to-period}.
	Finally, choosing a basis \(\{\omega_j\}\) identifies \(H^0(X,\Omega_X^1)^\vee\cong \C^g\) and sends \(H_1(X,\Z)\) to the period lattice \(\Lambda\).
\end{proof}

\subsubsection{The integral Hodge lattice, polarization, and Torelli's theorem}
\index{Hodge structure!weight one}
\index{Hodge isometry}
\index{polarization!principal}
\index{Torelli theorem}

The quotient $H^0(X,K_X)^\vee/H_1(X,\Z)$ remembers the complex torus
underlying the Jacobian.  To recover the curve one must also remember the
intersection pairing.  This extra datum is the \emph{principal polarization}.
The resulting package may be expressed either as a principally polarized
Jacobian or as a polarized integral Hodge structure on $H^1(X,\Z)$.  The two
languages contain the same information.

\begin{definition}[Integral Hodge structure of weight one]
	An \emph{integral Hodge structure of weight one} is a free abelian group
	$H_\Z$ of finite rank together with a decomposition
	\[
	 H_\C:=H_\Z\otimes_\Z\C
	 =H^{1,0}\oplus H^{0,1},
	 \qquad \overline{H^{1,0}}=H^{0,1}.
	\]
	For a compact Riemann surface $X$, the Hodge decomposition is
	\[
	 H^1(X,\C)=H^{1,0}(X)\oplus H^{0,1}(X),
	 \qquad
	 H^{1,0}(X)=H^0(X,K_X).
	\]
\end{definition}

\begin{definition}[Polarization and Hodge isometry]
	A \emph{principal polarization} of a weight-one integral Hodge structure
	is a unimodular alternating form
	$Q:H_\Z\times H_\Z\to\Z$ such that
	\[
	 Q(H^{1,0},H^{1,0})=0,
	 \qquad iQ(\omega,\overline\omega)>0
	 \quad(0\ne\omega\in H^{1,0}).
	\]
	An isomorphism $T:H_\Z\to H'_\Z$ is a \emph{polarized integral Hodge
	isometry} if its complexification carries $H^{1,0}$ onto
	$(H')^{1,0}$ and
	\[
	 Q'(T\alpha,T\beta)=Q(\alpha,\beta)
	 \qquad(\alpha,\beta\in H_\Z).
	\]
\end{definition}

For a compact Riemann surface the form is the cup-product pairing
\begin{equation}
 \label{eq:curve-hodge-polarization}
 Q_X(\alpha,\beta)
 :=\langle\alpha\smile\beta,[X]\rangle
 =\int_X\alpha\wedge\beta.
\end{equation}
Here the final expression uses closed differential-form representatives.
Poincar\'e duality makes $Q_X$ unimodular.  Type considerations give
$Q_X(H^{1,0},H^{1,0})=0$, and in a local coordinate
$\omega=f(z)\,dz$ gives
\[
 i\int_X\omega\wedge\overline\omega
 =2\int_X|f|^2\,dx\,dy>0
\]
unless $\omega=0$.  These are precisely the Riemann bilinear relations.

\begin{proposition}[The integration lattice]
 \label{prop:integration-lattice-hodge}
 Let $V_X=H^0(X,K_X)^\vee$ and define
 \[
 I_X:H_1(X,\Z)\longrightarrow V_X,
 \qquad
 I_X(\gamma)(\omega)=\int_\gamma\omega.
 \]
 Then $I_X$ is injective, its image is a full lattice, and
 \[
 J(X)=V_X/I_X\bigl(H_1(X,\Z)\bigr).
 \]
 Under Poincar\'e duality
 $H^1(X,\Z)\simeq H_1(X,\Z)$, this integral lattice is the image of
 $H^1(X,\Z)$ inside $H^1(X,\mathcal O_X)\simeq V_X$.  Moreover, the
 intersection form $Q_X$ induces the canonical principal polarization
 $\Theta_X$ on $J(X)$.
\end{proposition}

\begin{proof}
 The first three assertions are the coordinate-free form of
 Proposition~\ref{prop:period-lattice-full}: if all periods of a real cycle
 vanish, it pairs trivially with the real and imaginary parts of every
 holomorphic form.  Those harmonic forms span $H^1(X,\R)$, so nondegeneracy
 of the de Rham pairing forces the cycle to vanish.  Since source and target
 have real dimension $2g$, the real extension of $I_X$ is an isomorphism;
 hence the image of the integral homology lattice is a full lattice.

 Proposition~\ref{prop:integral-to-period} identifies this same subgroup
 with the image of $H^1(X,\Z)$ in $H^1(X,\mathcal O_X)$.  Finally,
 $Q_X$ is integral and unimodular, and the two Riemann bilinear relations
 verified above say exactly that it is a Riemann form on the complex torus.
 Unimodularity makes the resulting polarization principal.  Its divisor
 class is the theta class $[\Theta_X]$; see
 \cite{FarkasKra1992,GriffithsHarris1994,
 ArbarelloCornalbaGriffithsHarris1985}.
\end{proof}

\begin{proposition}[Polarized Hodge structures and polarized Jacobians]
 \label{prop:hodge-isometry-jacobian-isomorphism}
 For compact Riemann surfaces $X$ and $Y$, the following are equivalent.
 \begin{enumerate}
  \item There is a polarized integral Hodge isometry
  \[
   T:\bigl(H^1(X,\Z),H^{1,0}(X),Q_X\bigr)
   \longrightarrow
   \bigl(H^1(Y,\Z),H^{1,0}(Y),Q_Y\bigr).
  \]
  \item There is an isomorphism of principally polarized complex tori
  \[
   (J(X),\Theta_X)\simeq(J(Y),\Theta_Y).
  \]
 \end{enumerate}
\end{proposition}

\begin{proof}
 A Hodge isometry sends $H^{0,1}(X)$ complex linearly onto
 $H^{0,1}(Y)$ and sends the integral lattice to the integral lattice.
 It therefore descends to an isomorphism
 \[
 \frac{H^{0,1}(X)}{H^1(X,\Z)}
 \longrightarrow
 \frac{H^{0,1}(Y)}{H^1(Y,\Z)},
 \]
 where the notation means the image of the integral group under projection
 to $H^{0,1}$.  Via Serre duality these quotients are $J(X)$ and $J(Y)$.
 Preservation of $Q$ is exactly preservation of their Riemann forms, hence
 of the principal polarizations.

 Conversely, an isomorphism of complex tori induces an isomorphism on
 integral first cohomology and respects the decomposition into holomorphic
 and antiholomorphic one-forms.  If it also preserves the principal
 polarization, its action preserves the corresponding alternating form
 $Q$.  The induced map is therefore a polarized integral Hodge isometry.
\end{proof}

\begin{theorem}[Global Torelli theorem for compact Riemann surfaces]
 \label{thm:global-torelli-curves}
 Let $X$ and $Y$ be compact connected Riemann surfaces.  Then the following
 are equivalent:
 \begin{enumerate}
  \item $X$ and $Y$ are biholomorphic;
  \item their principally polarized Jacobians are isomorphic;
  \item their first integral cohomology groups are isomorphic as polarized
  Hodge structures.
 \end{enumerate}
 In particular, a polarized integral Hodge isometry
 $H^1(X,\Z)\simeq H^1(Y,\Z)$ determines the isomorphism class of the curve.
\end{theorem}

\begin{proof}
 If $f:X\to Y$ is biholomorphic, pullback preserves integral cohomology,
 types of differential forms, wedge products, and the fundamental class.
 Thus $f^*:H^1(Y,\Z)\to H^1(X,\Z)$ is a polarized Hodge isometry.
 Proposition~\ref{prop:hodge-isometry-jacobian-isomorphism} proves the
 equivalence of the last two conditions.

 It remains to recover the curve from the principally polarized Jacobian.
 Fix a base point and write $W_d(X)$ for the Abel--Jacobi image of
 $S^dX$.  Riemann's theta theorem identifies, up to translation,
 \[
  W_{g-1}(X)=\Theta_X.
 \]
 The classical Torelli reconstruction lemma recovers the Abel--Jacobi curve
 $W_1(X)$, up to translation and the involution $a\mapsto-a$, from the
 polarized pair $(J(X),[\Theta_X])$.  One way to formulate the reconstruction
 is to use the loci where translates of the theta divisor acquire the
 multiplicities prescribed by effective divisors; Riemann's singularity
 theorem identifies those multiplicities with
 $h^0(X,\mathcal O_X(D))$.  The smallest resulting Abel locus is
 $W_1(X)$.  Translation and inversion alter only the choice of base point
 and the sign of the Abel--Jacobi map, not the abstract curve.  Therefore an
 isomorphism
 $(J(X),\Theta_X)\simeq(J(Y),\Theta_Y)$ carries the reconstructed copy of
 $X$ biholomorphically onto the reconstructed copy of $Y$.

 For $g=0$, both curves are isomorphic to $\PP^1$.  For $g=1$, the
 Abel--Jacobi map identifies each pointed curve with its Jacobian directly.
 The theta-divisor reconstruction above is the substantive case $g\ge2$.
 Detailed proofs of the reconstruction and Riemann singularity theorems are
 given in \cite{FarkasKra1992,GriffithsHarris1994,
 ArbarelloCornalbaGriffithsHarris1985}.
\end{proof}

\begin{remark}[Why ``Hodge isometric'' must include the polarization]
 An isomorphism of the underlying unpolarized complex tori is not the datum
 in Torelli's theorem: a complex torus may admit more than one principal
 polarization.  The intersection form $Q_X$, or equivalently the theta
 polarization, is what records which principally polarized torus is the
 Jacobian of the curve.  Thus the precise cohomological statement is
 ``polarized integral Hodge isometry,'' not merely an isomorphism of the
 groups $H^1(X,\Z)$.
\end{remark}

\begin{exercise}[$\mathrm{Pic}^0$ of an elliptic curve]
	\label{exr:elliptic-pic0}
	Let \(X=\C/\Lambda\) be a complex torus. Prove directly from the exponential sequence that
	\(\mathrm{Pic}^0(X)\cong \C/\Lambda\) and identify the isomorphism with \(x\mapsto \mathcal O_X(x-p_0)\).
\end{exercise}

\subsection{Divisors, linear systems, and Abel's theorem}
\label{subsec:divisors-linebundles}

The preceding quotient becomes useful only after divisor geometry is mapped
into it. An effective divisor records zeros, while the vector space
$H^0(X,\mathcal O_X(D))$ records meromorphic functions whose poles are
bounded by $D$. Projectivizing this vector space gives the complete linear
system $|D|$. Abel's theorem says that these projective spaces are precisely
the fibers of the Abel--Jacobi map. Jacobi inversion then asks when such a
fiber is nonempty, and Riemann--Roch predicts its dimension. In this way the
analytic period integral, the geometric motion of points, and the algebraic
dimension $\ell(D)$ become three descriptions of the same phenomenon.

\subsubsection{Divisors and the associated line bundle}

A divisor on \(X\) is a finite sum \(D=\sum_{p\in X} n_p\,p\) with \(n_p\in\Z\).
Its degree is \(\deg D=\sum n_p\).
Choose meromorphic functions $g_\alpha$ on an open cover such that
$(g_\alpha)=D|_{U_\alpha}$. The quotients
$g_\alpha/g_\beta\in\mathcal O_X^*(U_\alpha\cap U_\beta)$ are transition
functions for the associated line bundle \(\mathcal O_X(D)\). Equivalently,
\(\mathcal O_X(D)(U)\) is the space of meromorphic functions $f$ on $U$
such that $(f)+D|_U\ge0$. This second description fixes the sign convention
used below.

\begin{proposition}[Sections as meromorphic functions with bounded poles]
	\label{prop:sections-of-OD}
	There is a natural identification
	\[
	H^0\big(X,\mathcal O_X(D)\big)
	=
	\Big\{\, f\ \text{meromorphic on }X\ \Big|\ (f)+D\ge 0\,\Big\}\ \cup\ \{0\}.
	\]
\end{proposition}

\begin{proof}
	By definition, \(\mathcal O_X(D)(U)\) consists of meromorphic functions \(f\) on \(U\) satisfying
	\((f)|_U+D|_U\geq0\). Thus its global holomorphic sections are exactly the meromorphic
	functions satisfying \((f)+D\geq0\), together with the zero section. In the equivalent line-bundle
	description, multiplication by the canonical meromorphic section converts this sheaf convention into
	the usual statement that a holomorphic section has an effective zero divisor.
\end{proof}

\subsubsection{Complete linear system}

\begin{definition}[Complete linear system]
	\label{def:complete-linear-system}
	If \(H^0(X,\mathcal O_X(D))\neq 0\), the \emph{complete linear system} of \(D\) is
	\[
	|D|:=\mathbb{P}\!\big(H^0(X,\mathcal O_X(D))\big).
	\]
\end{definition}

\begin{lemma}[Linear system equals effective divisors in the class]
	\label{lem:linear-system-fiber}
	The map \(s\mapsto (s)+D\) induces a bijection between \(|D|\) and the set of effective divisors \(D'\) with \(D'\sim D\).
	In particular, \(\dim |D|=h^0(X,\mathcal O_X(D))-1\) whenever \(H^0(X,\mathcal O_X(D))\neq 0\).
\end{lemma}

\begin{proof}
	Let \(0\neq s\in H^0(X,\mathcal O_X(D))\). Then \((s)+D\ge 0\) by Proposition~\ref{prop:sections-of-OD}, so
	\(D_s:=(s)+D\) is effective and linearly equivalent to \(D\) because \((s)=D_s-D\) is a principal divisor in \(\mathrm{Div}(X)\)
	(indeed \((s)\) is the divisor of the meromorphic function representing \(s\) in a meromorphic trivialization).
	Scaling \(s\) does not change \((s)\), hence not \(D_s\), so the map factors through \(|D|\).
	
	Conversely, if \(D'\ge 0\) and \(D'\sim D\), then \(D'-D=(f)\) for a meromorphic function \(f\).
	The condition \(D'\ge 0\) means \((f)+D\ge 0\), hence \(f\in H^0(X,\mathcal O_X(D))\setminus\{0\}\),
	and \(D'=(f)+D\) is the divisor attached to the section \(s=f\).
	Uniqueness up to scalar follows since \((f)=(cf)\) for \(c\in\C^*\).
\end{proof}

\begin{proposition}[Every line bundle is \(\mathcal O_X(D)\)]
	\label{prop:every-L-is-OD}
	Every holomorphic line bundle \(L\) on \(X\) is isomorphic to \(\mathcal O_X(D)\) for some divisor \(D\).
	Moreover, \(\mathcal O_X(D)\cong \mathcal O_X(D')\) if and only if \(D\sim D'\).
\end{proposition}

\begin{proof}
	Fix a point \(p_0\in X\). We first produce a meromorphic section
	without assuming the line-bundle form of Riemann--Roch.  For every
	integer \(m\ge0\), the add-a-point sequence
	\[
	0\longrightarrow L(mp_0)\longrightarrow L((m+1)p_0)
	\longrightarrow L((m+1)p_0)|_{p_0}\longrightarrow0
	\]
	has a one-dimensional quotient supported at \(p_0\).  Its cohomology
	sequence makes
	\(H^1(X,L(mp_0))\to H^1(X,L((m+1)p_0))\) surjective.
	These groups are finite-dimensional by the Dolbeault--Hodge theorem, so
	their nonincreasing dimensions eventually stabilize.  Thereafter the
	map on \(H^1\) is an isomorphism; exactness makes the map from
	\(H^0(X,L((m+1)p_0))\) onto the one-dimensional quotient surjective.
	Consequently \(h^0(X,L((m+1)p_0))=h^0(X,L(mp_0))+1\) for all sufficiently
	large \(m\), and some \(L(mp_0)\) has a nonzero holomorphic section.
	Let \(0\neq s\in H^0(X,L(mp_0))\). Then \(s\) is a meromorphic section of \(L\) with poles of order \(\le m\) at \(p_0\)
	and no other poles. Let \(D:=(s)\) be its divisor (zeros minus poles). The
	meromorphic section \(s\) trivializes \(L\) away from
	\(\operatorname{supp}D\); in a local holomorphic frame it is represented by a
	meromorphic function whose order is the coefficient of \(D\). Thus the
	resulting transition functions are those of \(\mathcal O_X(D)\), and
	\(L\simeq \mathcal O_X(D)\).
	
	If \(L\simeq\mathcal O_X(D)\simeq\mathcal O_X(D')\), then \(\mathcal O_X(D-D')\simeq\mathcal O_X\).
	An isomorphism \(\mathcal O_X(D')\to\mathcal O_X(D)\) is multiplication by a
	nonzero meromorphic function \(f\), and regularity of the map and its inverse
	gives \((f)=D'-D\). Hence \(D\sim D'\).
	Conversely, if \(D\sim D'\) then \(D-D'=(f)\) and multiplication by \(f\)
	identifies \(\mathcal O_X(D)\cong\mathcal O_X(D')\).
\end{proof}

\begin{exercise}[Degree of a divisor line bundle]
	\label{exr:OD-degree}
	Show that \(\deg(\mathcal O_X(D))=\deg D\), i.e.\ \(c_1(\mathcal O_X(D))\) corresponds to \(\deg D\in\Z\cong H^2(X,\Z)\).
\end{exercise}

\subsubsection{\texorpdfstring{Abel's theorem: fibers of $\alpha_d$ and linear equivalence}{Abel's theorem: fibers and linear equivalence}}
\label{subsec:abel-theorem}

The following compatibility is the analytic step that connects path integrals
with divisor classes; it must be established before the fiber description can
be used.

\paragraph{Point-divisor compatibility.}
Under the canonical isomorphism $\mathrm{Pic}^0(X)\cong J(X)$, the class of
\(\mathcal O_X(x-p_0)\) corresponds to \(\alpha(x)\). Here the
Serre-duality identification is normalized as in
Proposition~\ref{prop:integral-to-period}, so that an integral homology cycle
acts by the unscaled period functional
\(\omega\mapsto\int_\gamma\omega\).

\begin{proof}
	Consider the holomorphic map
	\[
	\Phi:X\longrightarrow\mathrm{Pic}^0(X),\qquad
	x\longmapsto\mathcal O_X(x-p_0).
	\]
	The tangent space of \(\mathrm{Pic}^0(X)\) at the identity is
	$H^1(X,\mathcal O_X)$, which Serre duality identifies with
	$H^0(X,K_X)^\vee$. If $z$ is a local coordinate at $x$ and the point is
	moved from $z=0$ to $z=t$, differentiating the local divisor transition
	function $(z-t)/z$, interpreted with this period normalization, shows that
	$d\Phi_x(v)$, paired with a holomorphic
	differential \(\omega\), is \(\omega_x(v)\). Differentiating the path integral
	defining \(\alpha\) gives the same formula:
	\[
	d\alpha_x(v)(\omega)=\omega_x(v).
	\]
	Thus \(\Phi\) and \(\alpha\), viewed in the same complex torus, have the same
	differential. Their difference is locally constant, hence constant because
	$X$ is connected. Both maps take $p_0$ to the identity, so they agree.
\end{proof}

\begin{theorem}[Abel's theorem: fibers are complete linear systems]
	\label{thm:abel-fiber-refined}
	Let \(d\ge 0\) and \(D,D'\in S^dX\). Then
	\[
	\alpha_d(D)=\alpha_d(D')
	\quad\Longleftrightarrow\quad
	D\sim D'.
	\]
	Equivalently, for an effective divisor \(D\) of degree \(d\),
	\[
	\alpha_d^{-1}\big(\alpha_d(D)\big)=|D|.
	\]
\end{theorem}

\begin{proof}
	Write \(D=\sum_{i=1}^d x_i\) and \(D'=\sum_{i=1}^d y_i\).
	Consider the degree-zero divisor \(E:=D'-D\). By the point-divisor
	compatibility just proved, under the identification
	\(\mathrm{Pic}^0(X)\cong J(X)\) of Theorem~\ref{thm:Pic0-is-J},
	the class \(\alpha_d(D')-\alpha_d(D)\in J(X)\) corresponds to the line bundle
	\[
	\mathcal O_X(D')\otimes\mathcal O_X(D)^{-1}\ \cong\ \mathcal O_X(E)\ \in\ \mathrm{Pic}^0(X).
	\]
	Thus \(\alpha_d(D')=\alpha_d(D)\) if and only if \(\mathcal O_X(E)\cong\mathcal O_X\),
	i.e.\ \(E\) is principal, i.e.\ \(D'\sim D\).
	The fiber statement then follows from Lemma~\ref{lem:linear-system-fiber}.
\end{proof}

\begin{corollary}[Fiber dimension]
	\label{cor:fiber-dim}
	If \(D\) is effective of degree \(d\) and \(H^0(X,\mathcal O_X(D))\neq 0\), then
	\[
	\dim \alpha_d^{-1}(\alpha_d(D))=\dim |D|=h^0(X,\mathcal O_X(D))-1.
	\]
\end{corollary}

\begin{exercise}[Degree one and the hyperelliptic test]
	\label{exr:abel-degree1}
	Assume \(g\ge 1\). Show that \(\alpha_1:X\to J(X)\) is injective. If \(g\ge2\),
	show that \(X\) is hyperelliptic precisely when some fiber of
	\(\alpha_2:S^2X\to J(X)\) contains a projective line
	(the degree-two linear system \(g^1_2\)).
\end{exercise}

\subsubsection{Jacobi inversion and the geometry of Riemann--Roch}
\label{subsec:jacobi-inversion}

\begin{theorem}[Jacobi inversion]
	\label{thm:jacobi-inversion}
	The Abel--Jacobi map \(\alpha_g:S^gX\to J(X)\) is surjective.
\end{theorem}

\begin{proof}
	Let \(y\in J(X)\cong \mathrm{Pic}^0(X)\) and let \(L_y\) be the corresponding degree-zero line bundle.
	Set \(D_0:=g p_0\) and apply Riemann--Roch to the line bundle
	\(M:=\mathcal O_X(D_0)\otimes L_y\):
	\[
	h^0(X,M)-h^0(X,K_X\otimes M^\vee)=\deg(M)+1-g=g+1-g=1.
	\]
	In particular, \(h^0(X,M)\ge 1\), so pick a nonzero section \(s\in H^0(X,M)\).
	Let \(E=Z(s)\) be the zero divisor of this section. It is effective and has degree \(g\).
	By definition of \(M\),
	\[
	\mathcal O_X(E)\ \cong\ \mathcal O_X(D_0)\otimes L_y,
	\qquad\text{hence}\qquad
	\mathcal O_X(E-D_0)\ \cong\ L_y.
	\]
	Under \(\mathrm{Pic}^0(X)\cong J(X)\), this means
	\(\alpha_g(E)-g\alpha(p_0)=y\). Since \(\alpha(p_0)=0\), we obtain \(\alpha_g(E)=y\).
\end{proof}

\begin{exercise}[Jacobi inversion in degree $g+1$]
	\label{exr:surj-g-plus-one}
	Show that \(\alpha_{g+1}:S^{g+1}X\to J(X)\) is also surjective and that a generic fiber has dimension \(1\).
\end{exercise}

\paragraph{Riemann--Roch from the geometry of the Abel map.}
\label{subsec:RR-from-AJ}

\begin{proposition}[Dimension bound from fibers]
	\label{prop:RR-lower-bound}
	Let \(D\) be an effective divisor of degree \(d\). Then
	\[
	h^0(X,\mathcal O_X(D))\ \ge\ d-g+1,
	\qquad\text{equivalently}\qquad
	\dim|D|\ \ge\ d-g.
	\]
\end{proposition}

\begin{proof}
	The symmetric power \(S^dX\) is a complex manifold of dimension \(d\):
	near a divisor with repeated points, elementary symmetric functions
	of local coordinates give its charts.  The Jacobian \(J(X)\) has
	dimension \(g\).  Fix \(D\in S^dX\).  Choose \(g\) local coordinate
	functions near \(\alpha_d(D)\) in the Jacobian.  Locally the fiber of
	\(\alpha_d\) is cut out in the smooth \(d\)-fold \(S^dX\) by at most
	\(g\) holomorphic equations.  The analytic local dimension theorem
	therefore makes every irreducible component of this fiber at least
	\(d-g\)-dimensional; no submersion assumption is needed.  When
	\(d<g\), this negative lower bound is vacuous.
	Theorem~\ref{thm:abel-fiber-refined} identifies the fiber over
	\(\alpha_d(D)\) with the complete linear system \(|D|\).
	By Lemma~\ref{lem:linear-system-fiber}, its members are the lines
	in \(H^0(X,\mathcal O_X(D))\), so
	\(|D|\cong\PP H^0(X,\mathcal O_X(D))\) and
	\(\dim|D|=h^0(X,\mathcal O_X(D))-1\).  Combining the two dimension
	statements and adding one proves the inequality.
\end{proof}

\begin{remark}
	Riemann--Roch refines Proposition~\ref{prop:RR-lower-bound} to
	\[
	h^0(\mathcal O_X(D)) - h^0(K_X\otimes \mathcal O_X(-D)) = d-g+1.
	\]
\end{remark}

\paragraph{Abel's theorem in degree zero: principal divisors.}
\label{subsec:abel-degree-zero}

Let \(\mathrm{Div}^0(X)\) denote divisors of degree \(0\). Fix \(p_0\in X\) and define a group homomorphism
\[
\alpha_0:\mathrm{Div}^0(X)\to J(X),
\qquad
\alpha_0\Big(\sum n_p\,p\Big):=\sum n_p\,\alpha(p).
\]

\begin{theorem}[Abel's theorem, degree zero]
	\label{thm:abel-degree0}
	For \(D\in\mathrm{Div}^0(X)\),
	\[
	\alpha_0(D)=0
	\quad\Longleftrightarrow\quad
	D \text{ is principal (i.e.\ \(D=(f)\) for some meromorphic \(f\)).}
	\]
\end{theorem}

\begin{proof}
	First assume \(\alpha_0(D)=0\). Write \(D=D_+-D_-\) with \(D_\pm\) effective of the same degree \(d\).
	The condition \(\alpha_0(D)=0\) means \(\alpha_d(D_+)=\alpha_d(D_-)\).
	By Theorem~\ref{thm:abel-fiber-refined}, \(D_+\sim D_-\), so \(D_+-D_-=(f)\) for some meromorphic \(f\),
	i.e.\ \(D\) is principal.
	
	(\(\Leftarrow\).)
	If \(D=(f)\), then \(D\sim 0\), hence \(\mathcal O_X(D)\cong\mathcal O_X\), so the corresponding class in
	\(\mathrm{Pic}^0(X)\cong J(X)\) is the identity. Equivalently, \(\alpha_0(D)=0\).
\end{proof}

\begin{exercise}[A degree-zero test]
	\label{exr:abel-degree0-check}
	Assume $g\ge1$ and let \(D=p-q\) for points \(p,q\in X\). Show that
	\(\alpha_0(D)=0\) if and only if $p=q$.  Indeed, if $p\ne q$ and
	$(f)=p-q$, then $f:X\to\CP^1$ would be a meromorphic map of degree $1$, hence a
	biholomorphism, contradicting $g\ge1$.  Explain how this gives another proof that
	$\alpha_1:X\to J(X)$ is injective.
\end{exercise}

\subsection{Worked examples: genus 0 and genus 1}
\label{subsec:jacobian-examples}

\subsubsection{\(X=\CP^1\) (genus \(0\))}

Here \(H^0(\CP^1,\Omega^1)=0\), hence \(\C^g=0\), \(\Lambda=0\), and \(J(\CP^1)=0\).
Also \(\mathrm{Pic}^0(\CP^1)=0\) and \(\mathrm{Pic}(\CP^1)\cong \Z\) via degree.

\begin{exercise}[The Picard group of $\CP^1$]
	\label{exr:cp1-pic}
	Show \(\mathrm{Pic}(\CP^1)\cong \Z\) by proving \(\mathcal O_{\CP^1}(D)\cong \mathcal O_{\CP^1}(\deg D)\).
\end{exercise}

\subsubsection{Elliptic curve \(X=\C/\Lambda\) (genus \(1\))}

Let \(\omega=dz\). Then \(H^0(X,\Omega^1)=\C\cdot dz\), and
\[
\Lambda_{\mathrm{per}}=\Big\{\int_\gamma dz:\gamma\in H_1(X,\Z)\Big\}=\Lambda\subset\C.
\]
Hence \(J(X)=\C/\Lambda\cong X\). With base point \(p_0=0\), the Abel--Jacobi map is
\[
\alpha([z]) = z \bmod \Lambda,
\]
which is a group isomorphism.

\begin{exercise}[The Abel map as addition]
	\label{exr:elliptic-abel}
	Prove that for an elliptic curve \(X\), the map \(\alpha_d:S^dX\to J(X)\cong X\) is the addition map on the group \(X\)
	(after identifying \(S^dX\) with effective degree-\(d\) divisors).
\end{exercise}	
	
\section{Compact Riemann Surfaces as Algebraic Curves}
\index{algebraic curves}
\label{sec:RSasAlgCurve}

\subsection{Affine and projective foundations}
Let $M$ be a compact connected Riemann surface of genus $g$.
We explain constructively:
\begin{enumerate}
	\item how a divisor $D$ with $\deg D\ge 2g+1$ gives an embedding $M\hookrightarrow \mathbb{CP}^N$,
	\item how a pair of meromorphic functions satisfies an algebraic relation
	and therefore gives a plane-curve image,
	\item why the normalization of that image is $M$ when the chosen pair
	generates the full function field $\mathbb{C}(M)$.
\end{enumerate}
The algebraization and normalization steps are compared with \cite{Forster1981,Miranda1995,GunningRossi1965}.

\paragraph{Prerequisite check: a projective model of $\PP^1$.}
With affine coordinate $z$ and $D=2\infty$, the sections
$1,z,z^2$ of $\mathcal O_{\PP^1}(D)$ define
$z\mapsto[1:z:z^2]$. If the target coordinates are
$[X_0:X_1:X_2]$, the image satisfies $X_0X_2=X_1^2$;
the point at infinity maps to $[0:0:1]$. This small example exhibits
the two steps used below: sections give a projective map, and an
algebraic relation cuts out its image. For a general surface the
remaining question is whether the map separates points and tangent
directions.

\subsubsection{Notation}
For a divisor $D$ on $M$, write $\OO_M(D)$ for the associated line bundle and set
\[
L(D):=H^0(M,\OO_M(D)),\qquad \ell(D):=\dim_{\mathbb{C}}L(D).
\]
Let $K$ be a canonical divisor. We use Riemann--Roch:
\[
\ell(D)-\ell(K-D)=1-g+\deg D.
\]

\subsubsection{Affine varieties and Hilbert's Nullstellensatz}
\index{affine variety}\index{Hilbert Nullstellensatz}

Before passing from meromorphic functions to a projective curve, we need the
dictionary between polynomial equations and coordinate rings. For an ideal
$I\subset\C[x_1,\ldots,x_n]$ and a subset $Y\subset\C^n$, set
\begin{align*}
	V(I)&:=\{a\in\C^n:f(a)=0\text{ for every }f\in I\},\\
	I(Y)&:=\{f\in\C[x_1,\ldots,x_n]:f(a)=0\text{ for every }a\in Y\}.
\end{align*}

\begin{definition}[Affine algebraic set, affine variety, and coordinate ring]
	An \emph{affine algebraic set} is a subset $Y=V(I)$ of $\C^n$.
	It is an \emph{affine variety} when it is irreducible in the Zariski
	topology, equivalently when $I(Y)$ is prime. Its coordinate ring is
	$\C[Y]:=\C[x_1,\ldots,x_n]/I(Y)$.
\end{definition}

\begin{definition}[Affine scheme and projective variety]
	For a commutative ring $A$, the \emph{affine scheme} $\Spec A$ is the set
	of prime ideals of $A$, equipped with the Zariski topology and its structure
	sheaf. When $A=\C[x_1,\ldots,x_n]/I$, its closed points recover the common
	zeros of $I$, while nilpotent elements retain multiplicity information that
	the point set alone forgets. A \emph{projective algebraic set} is the common
	zero locus in $\PP^n$ of homogeneous polynomials; it is a \emph{projective
	variety} when it is irreducible and reduced.
\end{definition}

\begin{theorem}[Hilbert's Nullstellensatz]
	\label{thm:nullstellensatz-curves}
	For $A=\C[x_1,\ldots,x_n]$ the following statements hold.
	\begin{enumerate}
		\item Every maximal ideal of $A$ is
		$\mathfrak m_a=(x_1-a_1,\ldots,x_n-a_n)$ for a unique $a\in\C^n$.
		\item For every ideal $I\subset A$,
		$I(V(I))=\sqrt I:=\{f\in A:f^N\in I\text{ for some }N\geq1\}$.
	\end{enumerate}
\end{theorem}

\begin{proof}
	If $\mathfrak m$ is maximal, then $A/\mathfrak m$ is a field finitely
	generated as a $\C$-algebra. Zariski's lemma makes it a finite algebraic
	extension of $\C$, hence $\C$ itself. The coordinate images are numbers
	$a_i$, and the kernel is $\mathfrak m_a$.

	The inclusion $\sqrt I\subset I(V(I))$ is immediate. Conversely, suppose
	$f$ vanishes on $V(I)$ and put
	$J=(I,1-tf)\subset\C[x_1,\ldots,x_n,t]$. This ideal has no common zero:
	a zero would have its $x$-coordinate in $V(I)$, where $f=0$, and would also
	satisfy $1-tf=0$. By the first statement $J$ cannot be proper. Writing
	$1=\sum_j a_j(x,t)g_j(x)+b(x,t)(1-tf)$ with $g_j\in I$ and substituting
	$t=1/f$ in the localization $A_f$ shows $1\in IA_f$. Clearing a power of
	$f$ gives $f^N\in I$ for some $N$.
\end{proof}

The radical is essential: $V(x^2)=V(x)=\{0\}$, so a zero set cannot
distinguish $(x^2)$ from $(x)$. By contrast, the parabola
$Y=V(y-x^2)$ is irreducible and
$\C[Y]=\C[x,y]/(y-x^2)\cong\C[x]$, whereas $V(xy)$ is the reducible union
of the coordinate axes. These examples distinguish a point set, its full
ideal of equations, and a possibly nonreduced scheme defined by a chosen ideal.

To compactify $P(x,y)=0$ of degree $d$, homogenize it to
\[
	P^h(X_0,X_1,X_2)=X_0^dP(X_1/X_0,X_2/X_0)
\]
and take its zero set in $\PP^2$. The chart $X_0\ne0$ recovers the affine
curve, while $X_0=0$ records points at infinity. Projective closure can
introduce singularities, so the associated compact Riemann surface is its
normalization, not necessarily the raw plane curve. Below, meromorphic
functions $f,g$ produce $\C[x,y]\to\C(M)$; its prime kernel records the
irreducible affine image, homogenization adds points at infinity, and
normalization restores the smooth compact surface. The Nullstellensatz
justifies the passage from that polynomial kernel to its zero set
\cite{Hartshorne1977}.

\subsection{The GAGA principle: from projective curves to compact Riemann surfaces}
\label{subsec:gaga-principle}
\index{GAGA principle}\index{analytification}

The preceding affine dictionary uses polynomial functions and the Zariski
topology, whereas the earlier chapters use holomorphic functions and the
usual topology. For a general complex variety these two categories are not
the same. The compact curves studied in this book lie in the setting where
Serre's GAGA theorem makes the comparison exact. This is why a divisor, its
line bundle, and its cohomology may be computed in either algebraic or
analytic language once a compact Riemann surface has been realized as a
smooth projective curve.

\begin{definition}[Analytification]
	Let $X$ be a scheme of finite type over $\C$. Its \emph{analytification}
	$X^{\mathrm{an}}$ has the complex points of $X$ as its underlying set, with the
	complex-analytic topology and sheaf of holomorphic functions. On an affine
	chart $X=\Spec \C[x_1,\ldots,x_n]/I$, it is the analytic subspace
	\[
		X^{\mathrm{an}}=\{z\in\C^n:f(z)=0\text{ for every }f\in I\}.
	\]
	If $\mathcal F$ is a coherent algebraic sheaf on $X$, its analytification is
	\[
		\mathcal F^{\mathrm{an}}
		:=\mathcal O_{X^{\mathrm{an}}}\otimes_{a^{-1}\mathcal O_X}a^{-1}\mathcal F,
	\]
	where $a:X^{\mathrm{an}}\to X$ sends a complex point to the corresponding closed
	point. A regular map $f:X\to Y$ similarly induces a holomorphic map
	$f^{\mathrm{an}}:X^{\mathrm{an}}\to Y^{\mathrm{an}}$.
\end{definition}

The tensor product in this definition matters. An algebraic section is
locally represented by rational functions with controlled denominators;
after analytification, holomorphic coefficients are also allowed. The
theorem says that on a projective variety these new local coefficients do
not create new global coherent objects or new morphisms between them.

\begin{theorem}[Serre's GAGA theorem]
	\label{thm:gaga}
	Let $X$ be a projective complex algebraic variety. Then the following hold.
	\begin{enumerate}
		\item For every coherent algebraic sheaf $\mathcal F$ and every $q\geq0$,
		the natural comparison map is an isomorphism:
		\[
			H^q(X,\mathcal F)\longrightarrow
			H^q(X^{\mathrm{an}},\mathcal F^{\mathrm{an}}).
		\]
		\item Analytification is fully faithful on coherent sheaves:
		\[
			\operatorname{Hom}_{\mathcal O_X}(\mathcal F,\mathcal G)
			\cong
			\operatorname{Hom}_{\mathcal O_{X^{\mathrm{an}}}}
			(\mathcal F^{\mathrm{an}},\mathcal G^{\mathrm{an}}).
		\]
		\item Every coherent analytic sheaf on $X^{\mathrm{an}}$ is isomorphic to
		$\mathcal F^{\mathrm{an}}$ for a coherent algebraic sheaf $\mathcal F$, unique up
		to isomorphism. Thus analytification gives an equivalence of the two
		categories of coherent sheaves.
	\end{enumerate}
\end{theorem}

\begin{proof}[Proof architecture]
	A complete proof requires the finiteness theorems for coherent sheaves; we
	record the steps so that the hypotheses and the logical dependence are
	clear. On projective space, finite resolutions by sums of twisting sheaves
	$\mathcal O_{\PP^n}(m)$ reduce cohomology comparison to the explicit
	calculation of the cohomology of those twists. For a closed projective
	variety $X\subset\PP^n$, extension by zero transfers the comparison to
	$\PP^n$, which proves the first assertion.

	For coherent $\mathcal F$ and $\mathcal G$, presentations and the
	degree-zero case of cohomology comparison show that every analytic morphism
	is the analytification of a unique algebraic morphism. This is full
	faithfulness. Finally, twisting a coherent analytic sheaf $\mathscr E$ by a
	sufficiently positive $\mathcal O_{\PP^n}(m)$ makes it generated by finitely
	many global sections. Those sections give an analytic presentation by sums
	of twisting sheaves. Full faithfulness algebraizes the presentation and its
	relations, producing a coherent algebraic sheaf whose analytification is
	$\mathscr E$. These are the comparison, full-faithfulness, and existence
	steps in Serre's proof \cite{Serre1956GAGA,Hartshorne1977}.
\end{proof}

\begin{corollary}[The curve dictionary]
	\label{cor:gaga-curve-dictionary}
	Let $C$ be a smooth projective curve over $\C$ and put $M=C^{\mathrm{an}}$. Then:
	\begin{enumerate}
		\item algebraic and holomorphic line bundles correspond, and
		$\Pic(C)\cong\Pic(M)$;
		\item algebraic meromorphic functions on $C$ and meromorphic functions on
		$M$ form the same function field;
		\item for an algebraic divisor $D$, the algebraic line bundle
		$\mathcal O_C(D)$ analytifies to the holomorphic bundle
		$\mathcal O_M(D)$, and
		\[
			H^q(C,\mathcal O_C(D))\cong
			H^q(M,\mathcal O_M(D)),\qquad q=0,1;
		\]
		\item dimensions, Euler characteristics, and therefore the two forms of
		Riemann--Roch agree under this identification.
	\end{enumerate}
\end{corollary}

\begin{proof}
	Line bundles are locally free coherent sheaves of rank one, so the
	equivalence in Theorem~\ref{thm:gaga} gives the first assertion. A
	meromorphic function on the compact Riemann surface has a finite polar
	divisor $D$ and is a global section of $\mathcal O_M(D)$. Essential
	surjectivity algebraizes that line bundle, and full faithfulness algebraizes
	its section; hence the analytic function is algebraic. Divisors and their
	transition functions give the third assertion, while cohomology comparison
	gives its displayed isomorphisms. Taking alternating dimensions proves the
	last assertion.
\end{proof}

\begin{example}[The projective line]
	For $C=\PP^1_{\C}$, the analytic space is the Riemann sphere. GAGA
	identifies $\mathcal O_{\PP^1}(d)$ with the holomorphic line bundle of
	degree $d$ and identifies its algebraic sections with homogeneous
	polynomials of degree $d$. In the affine coordinate $z$, these become
	polynomials of degree at most $d$, with the behavior at infinity encoded by
	the transition function $z^d$. Thus the analytic calculation
	$h^0(\PP^1,\mathcal O(d))=d+1$ for $d\geq0$ is literally the algebraic one.
\end{example}

\begin{remark}[Why properness cannot be omitted]
	On $\A^1_{\C}$ the algebraic global functions are $\C[z]$, while on
	$(\A^1_{\C})^{\mathrm{an}}=\C$ the holomorphic global functions include $e^z$.
	Hence the degree-zero comparison already fails. Similarly, the unit disc is
	Stein but is not the analytification of an affine algebraic curve. Affine
	varieties and Stein spaces share powerful vanishing theorems, but this
	parallel is weaker than an equivalence of coherent categories.
\end{remark}

\begin{corollary}[Chow's theorem for projective space]
	Every closed complex-analytic subspace of $\PP^n(\C)$ is algebraic.
\end{corollary}

\begin{proof}
	Its coherent analytic ideal sheaf algebraizes by
	Theorem~\ref{thm:gaga}; the resulting coherent algebraic ideal sheaf defines
	a closed algebraic subspace having the prescribed analytification.
\end{proof}

\subsection{Very ample divisors, embeddings, and plane models}
\index{very ample divisors}
\label{subsec:RSasAlgCurve:very-ample}

\begin{definition}[Very ample line bundle and divisor]
	A holomorphic line bundle $L$ on a compact complex manifold is
	\emph{very ample} if a basis of $H^0(M,L)$ defines a holomorphic embedding
	$M\hookrightarrow\PP(H^0(M,L)^\vee)$. A divisor $D$ is very ample when
	$\mathcal O_M(D)$ is very ample. Thus ``very ample'' includes both the
	absence of base points and the separation of distinct points and tangent
	directions.
\end{definition}

Sections of $\mathcal O_M(D)$ give a map to projective space only if
they do not all vanish at the same point.  This is base-point
freeness.  To make the map an embedding, sections must also distinguish
two different points and detect a nonzero tangent direction at each
point.  The three properties below therefore serve different purposes:
the first defines the map everywhere, the second makes it injective,
and the third makes it an immersion.  Riemann--Roch will supply enough
sections when the degree of $D$ is large.  Their combination is the
bridge from an abstract compact Riemann surface to a projective
algebraic curve.

\begin{definition}[Base-point free and separation]
	\label{def:RSasAlgCurve:sep}
	Let $D$ be a divisor on $M$. Here $s(p)$ is the value of the section
	in the fiber $\OO_M(D)|_p$; a meromorphic function used to represent
	$s$ may have a pole at $p$ and cannot simply be evaluated as a function.
	\begin{enumerate}
		\item $|D|$ is \emph{base-point free} if for every $p\in M$ there exists $s\in L(D)$ with $s(p)\neq 0$.
		\item $|D|$ \emph{separates points} if for $p\neq q$ there exists $s\in L(D)$ with $s(p)=0$ and $s(q)\neq 0$.
		\item $|D|$ \emph{separates tangent directions} at $p$ if there exists $s\in L(D)$ vanishing at $p$ to order exactly $1$.
	\end{enumerate}
\end{definition}

\begin{lemma}[Base-point freeness for $\deg D\ge 2g$]
	\label{lem:RSasAlgCurve:basepointfree}
	If $\deg D\ge 2g$, then $|D|$ is base-point free.
\end{lemma}

\begin{proof}
	Fix $p\in M$. It suffices to show $\ell(D)>\ell(D-p)$.
	By Riemann--Roch applied to $D$ and $D-p$,
	\[
	\ell(D)-\ell(D-p)=1+\ell(K-D)-\ell(K-D+p).
	\]
	If $\deg D\ge 2g$, then $\deg(K-D+p)=2g-1-\deg D\le -1$, hence $\ell(K-D+p)=0$.
	Since $K-D\le K-D+p$, this also forces $\ell(K-D)=0$; consequently
	$\ell(D)-\ell(D-p)=1$.
\end{proof}

\begin{lemma}[Point and tangent separation for $\deg D\ge 2g+1$]
	\label{lem:RSasAlgCurve:separate}
	If $\deg D\ge 2g+1$, then $|D|$ separates points and tangent directions.
\end{lemma}

\begin{proof}
	\emph{(Points)} For $p\neq q$ it suffices to show $\ell(D-p)>\ell(D-p-q)$.
	Riemann--Roch gives
	\[
	\ell(D-p)-\ell(D-p-q)=1+\ell(K-D+p)-\ell(K-D+p+q).
	\]
	If $\deg D\ge 2g+1$, then $\deg(K-D+p+q)=2g-\deg D\le -1$, hence $\ell(K-D+p+q)=0$.
	It follows as well that $\ell(K-D+p)=0$, so the displayed difference is $1$.
	
	\smallskip
	\emph{(Tangents)} It suffices to show $\ell(D-p)>\ell(D-2p)$.
	Again,
	\[
	\ell(D-p)-\ell(D-2p)=1+\ell(K-D+p)-\ell(K-D+2p),
	\]
	and $\deg(K-D+2p)=2g-\deg D\le -1$, hence $\ell(K-D+2p)=0$.
	Thus $\ell(K-D+p)=0$ and this difference is also $1$.
\end{proof}

\begin{theorem}[Embedding by a complete linear system]
	\label{thm:RSasAlgCurve:embedding}
	Let $\deg D\ge 2g+1$ and let $s_0,\dots,s_N$ be a basis of $L(D)$.
	Then the map
	\[
	\Phi_D:M\longrightarrow \mathbb{CP}^{N},\qquad
	p\longmapsto [s_0(p):\cdots:s_N(p)]
	\]
	is a holomorphic embedding.
\end{theorem}

\begin{proof}
	Lemma~\ref{lem:RSasAlgCurve:basepointfree} implies $\Phi_D$ is well-defined.
	Lemma~\ref{lem:RSasAlgCurve:separate} implies $\Phi_D$ is injective and an immersion.
	Compactness of $M$ implies $\Phi_D$ is a topological embedding, hence holomorphic embedding.
\end{proof}

\subsubsection{Two meromorphic generators and a plane model}
\index{plane models!meromorphic generators}
\label{subsec:RSasAlgCurve:plane}

\begin{lemma}[Algebraic dependence in transcendence degree one]
	\label{lem:RSasAlgCurve:alg-dep}
	For any nonconstant $f,g\in \mathbb{C}(M)$ there exists a nonzero polynomial
	$P(x,y)\in \mathbb{C}[x,y]$ such that $P(f,g)=0$.
\end{lemma}

\begin{proof}
	The function field $\mathbb{C}(M)$ has transcendence degree $1$ over $\mathbb{C}$.
	Thus $g$ is algebraic over $\mathbb{C}(f)$, so it satisfies a polynomial with coefficients in $\mathbb{C}(f)$.
	Clearing denominators yields $P\in \mathbb{C}[x,y]\setminus\{0\}$ with $P(f,g)=0$.
\end{proof}

\begin{definition}[Normalization of an irreducible curve]
	Let $C$ be an irreducible algebraic curve with function field $K(C)$.
	The \emph{normalization} of $C$ is a normal curve $C^{\mathrm{norm}}$
	together with a finite birational morphism
	$C^{\mathrm{norm}}\to C$. On an affine open $\Spec A\subset C$, it is
	obtained from the integral closure of $A$ in $K(C)$. Over $\C$, a normal
	algebraic curve is nonsingular, and its analytification is a Riemann surface.
\end{definition}

\begin{proposition}[Plane curve and normalization (function-field view)]
	\label{prop:RSasAlgCurve:normalization}
	Let $f,g\in\mathbb{C}(M)$ be nonconstant and let $P(x,y)\in\mathbb{C}[x,y]$ be irreducible with $P(f,g)=0$.
	If $d=\deg P$, use the homogenization adapted to the affine chart $X\ne0$,
	\[
	\widetilde P(X,Y,Z):=X^dP(Y/X,Z/X),
	\]
	and set
	\[
	C:=\{[X:Y:Z]\in \mathbb{CP}^2 \mid \widetilde P(X,Y,Z)=0\}.
	\]
	Then:
	\begin{enumerate}
		\item The formula $p\mapsto [1:f(p):g(p)]$ extends uniquely across the
		poles of $f$ and $g$ to a holomorphic map $\nu:M\to C$.
		\item The induced map $\nu:M\to C$ is finite.
		\item If $\mathbb{C}(f,g)=\mathbb{C}(M)$, then $\nu$ is birational and $M$ is biholomorphic to the normalization $C^{\mathrm{norm}}$.
	\end{enumerate}
\end{proposition}

\begin{proof}
	(1) The relation $\widetilde P(1,f,g)=0$ holds where $f$ and $g$ are finite.
	To check extension at a point $p$, choose a local parameter $t$ and let $m$ be
	the larger of the pole orders of $f$ and $g$ at $p$ (with a holomorphic
	function assigned pole order $0$).  Multiplying the three coordinates by
	$t^m$ gives
	\[
	[t^m:t^mf:t^mg].
	\]
	All three entries are holomorphic near $p$, and at least one is nonzero at
	$p$ by the definition of $m$.  They therefore define a holomorphic extension
	to $p$.  The polynomial identity extends by continuity, so the image lies in
	$C$ everywhere.
	(2) The map is nonconstant and holomorphic.  Its image is all of $C$:
	the homomorphism $\C[x,y]\to\C(M)$ sending $x\mapsto f$ and $y\mapsto g$
	has a height-one prime kernel (the image has transcendence degree one), hence
	that kernel is $(P)$. Thus the Zariski closure of the affine image is the
	irreducible curve $P=0$. By the proper-mapping theorem
	\cite[Ch.~8]{Forster1981}, $\nu(M)$ is a compact
	analytic subset of the projective curve $C$. Since $\nu$ is nonconstant, its
	image has complex dimension one. The only one-dimensional analytic subset
	of the irreducible curve $C$ is $C$ itself; hence $\nu(M)=C$.
	Every fiber is a discrete compact subset of $M$, hence finite.  A proper
	holomorphic map of curves with finite fibers is finite, so $\nu:M\to C$ is finite.
	(3) If $\mathbb{C}(f,g)=\mathbb{C}(M)$ then the induced extension of function fields is an isomorphism, so $\nu$ is birational.
	A finite birational map from a smooth curve is the normalization morphism;
	compare \cite{StacksCurves2026}.
\end{proof}

\subsubsection{Examples and exercises}
\label{subsec:RSasAlgCurve:examples}

The field-generation condition in Proposition~\ref{prop:RSasAlgCurve:normalization}
matters. On $M=\CP^1$ with coordinate $z$, the pair $f=z^2$ and
$g=z^4$ satisfies $g=f^2$ but generates only $\C(z^2)$, not
$\C(z)$. The resulting map to its conic image identifies $z$ and
$-z$ generically; it is not the normalization map from $M$ to that
image. The next example chooses $z$ and $z^2$ instead, which do
generate $\C(z)$.

\begin{example}[$\mathbb{CP}^1$ as a conic]
	Let $M=\mathbb{CP}^1$ with affine coordinate $z$. Take $f=z$, $g=z^2$.
	Then $P(x,y)=y-x^2$; in homogeneous coordinates compatible with
	$[1:f:g]$, its homogenization is $XZ-Y^2$.
	Thus $\mathbb{CP}^1 \simeq \{XZ=Y^2\}\subset \mathbb{CP}^2$ via $[u:v]\mapsto [u^2:uv:v^2]$.
\end{example}

\begin{example}[Elliptic curves via Weierstrass $\wp$]
	Let $M=\mathbb{C}/\Lambda$. With $f=\wp$ and $g=\wp'$ one has
	\[
	g^2 = 4(f-e_1)(f-e_2)(f-e_3),
	\]
	giving a cubic plane model after homogenization.
\end{example}

\begin{example}[Hyperelliptic curves]
	For distinct $a_1,\dots,a_{2g+2}\in\mathbb{C}$, the smooth projective model of
	$y^2=\prod_{k=1}^{2g+2}(x-a_k)$ has function field $\mathbb{C}(x,y)$ and already yields a plane model.
\end{example}

\paragraph{Exercises.}
\label{subsec:RSasAlgCurve:exercises}

\begin{exercise}[Recovering a curve from its section ring]
	Show that $M$ can be recovered from a very ample divisor $D$ via the graded ring
	\[
	R(D):=\bigoplus_{n\ge 0} H^0\!\big(M,\OO_M(nD)\big).
	\]
	Prove (or outline carefully) why
	\[
	M \ \simeq\ \operatorname{Proj} R(D);
	\]
	compare the canonical open immersion for an ample invertible sheaf in
	\cite{StacksProperties2026}.
	(Hint: use the standard affine opens $D_+(s)\subset\operatorname{Proj}R(D)$
	and compare their degree-zero coordinate rings with the corresponding opens
	of $M$.  The canonical map has dense open image; properness of $M$ makes that
	image closed, hence equal to all of $\operatorname{Proj}R(D)$.  If one uses
	the embedding $\Phi_D$, remember that its homogeneous
	coordinate ring agrees with the section ring in degree one only when the
	embedding is projectively normal; passing to a sufficiently high Veronese
	subring avoids this issue.)
\end{exercise}

\begin{exercise}[Generic projection and generation of the function field]
	Let $\deg D\ge 2g+1$ and $\Phi_D:M\hookrightarrow \mathbb{CP}^N$,
	with $N\ge3$.
	Explain why a generic linear projection $\mathbb{CP}^N\dashrightarrow \mathbb{CP}^2$
	induces meromorphic functions $f,g$ with $\mathbb{C}(f,g)=\mathbb{C}(M)$.
	Hint: choose the center disjoint from $\Phi_D(M)$. Two distinct points
	have the same image exactly when their secant
	line meets the projection center. For a fixed secant, this imposes
	codimension two on the choice of center; compare with the two-dimensional
	family of pairs of points. If $N=2$, no projection is needed; if $N=1$,
	use a coordinate on $M\cong\PP^1$ instead.
\end{exercise}

\subsection{Local intersection multiplicity: definitions and computations}
\index{intersection multiplicity}
\label{sec:alg-intersection-mult-plane}

This section is compatible with the earlier treatment of intersection theory and
surface Riemann--Roch: we give a purely local (commutative algebra) definition of
$I_p(C,C')$ for plane curves, prove its key properties,
and explain how it matches the global intersection pairing on smooth projective surfaces.
Our length conventions, additivity, tangent-cone formula, and global comparison follow
\cite{Fulton1998Intersection,Hartshorne1977}.

There are three progressively finer questions. First, do two curves meet at
$p$? Second, are their tangent directions distinct? Third, if they are
tangent or singular, how much intersection is concentrated at $p$? The
quotient $\mathcal O_{X,p}/(f,g)$ answers all three at once. Its finite length
is $1$ for a transverse crossing and increases with higher-order contact.
This local algebra is not an auxiliary formalism: summing its lengths gives
the global intersection number and, in the projective plane, the product of
the degrees in B\'ezout's theorem.

\paragraph{Prerequisite check: length in one variable.}
In the local ring $R=k[x]_{(x)}$, the quotient $R/(x^m)$ has the
filtration
$R/(x^m)\supset(x)/(x^m)\supset\cdots\supset(x^{m-1})/(x^m)
\supset0$. Each successive quotient is one-dimensional over $k$,
so its length is $m$. In the examples below, eliminating one local
coordinate reduces a plane-curve intersection to exactly this kind
of calculation. Length retains the order of contact that a count
of distinct intersection points would lose.

\subsubsection{The local algebra and its formal properties}
\index{intersection multiplicity!length}
\label{subsec:alg-intersection-mult-plane:local}

Let $k$ be an algebraically closed field.
Let $\A^2=\Spec k[x,y]$ and let
\[
C=(f=0),\qquad C'=(g=0)
\]
be affine plane curves with $f,g\in k[x,y]$ nonconstant. Fix a point $p=(a,b)\in C\cap C'$ and set
\[
R:=\OO_{\A^2,p}\cong k[x,y]_{(x-a,y-b)},\qquad \fm_p:=(x-a,y-b)R.
\]
Define the local intersection algebra
\[
Q_p(f,g)\ :=\ R/(f,g).
\]

\begin{definition}[Module length and Artinian local ring]
	A finite-length $R$-module $N$ has a composition series
	$0=N_0\subset N_1\subset\cdots\subset N_r=N$ with simple quotients;
	the integer $r$, independent of the chosen series, is its \emph{length}
	$\length_R(N)$. A ring is \emph{Artinian} if every descending chain of
	ideals stabilizes. For a local $k$-algebra with residue field $k$, being
	Artinian and finite-dimensional over $k$ are equivalent, and in that case
	$\length_R(N)=\dim_kN$ for every finite-length module $N$ whose composition
	factors are $k$.
\end{definition}

\begin{lemma}[Isolated intersection $\Longleftrightarrow$ finite $k$-dimension]
	\label{lem:alg-intersection-mult-plane:isolated-finite}
	The following are equivalent:
	\begin{enumerate}
		\item $p$ is an isolated point of $C\cap C'$ in a Zariski neighbourhood of $p$;
		\item $Q_p(f,g)$ is an Artinian local ring;
		\item $\dim_k Q_p(f,g)<\infty$.
	\end{enumerate}
\end{lemma}

\begin{proof}
	The scheme-theoretic intersection near $p$ is $\Spec(R/(f,g))$.
	It is $0$-dimensional (hence $p$ is isolated) if and only if $R/(f,g)$ is Artinian.
	Since $k$ is algebraically closed, the residue field of $R$ is $k$, and an Artinian local $k$-algebra is finite-dimensional over $k$.
	Conversely, if $\dim_k(R/(f,g))<\infty$, then $R/(f,g)$ is Artinian.
\end{proof}

\begin{definition}[Algebraic local intersection multiplicity]
	\label{def:alg-intersection-mult-plane:Ip}
	Assume $p$ is an isolated intersection point (equivalently, Lemma~\ref{lem:alg-intersection-mult-plane:isolated-finite} holds).
	The \emph{intersection multiplicity} of $C$ and $C'$ at $p$ is
	\[
	I_p(C,C')\ :=\ \length_R\big(R/(f,g)\big)\ =\ \dim_k\big(R/(f,g)\big).
	\]
\end{definition}

\begin{example}[Crossing versus tangency at the origin]
	For the coordinate axes $x=0$ and $y=0$, the quotient
	$k[x,y]_{(x,y)}/(x,y)$ is $k$, so the intersection multiplicity is
	$1$. For the line $y=0$ and the parabola $y=x^2$, the quotient is
	\[
	k[x,y]_{(x,y)}/(y,y-x^2)\cong k[x]_{(x)}/(x^2),
	\]
	which has basis represented by $1,x$ and hence multiplicity $2$.
	The second pair meets at the same single point, but the length records
	its tangency rather than merely counting the point once.
\end{example}

\begin{remark}[Coordinate invariance]
	\label{rem:alg-intersection-mult-plane:coord-inv}
	The quantity $I_p(C,C')$ depends only on the isomorphism class of the local ring $R$ and the ideal $(f,g)\subset R$.
	In particular, it is invariant under an étale change of local coordinates;
	when $k=\mathbb C$, the same statement holds for a biholomorphic change of
	analytic coordinates. These are two settings for the same local-ring
	invariance, not an identification of the algebraic and analytic categories
	over an arbitrary field.
\end{remark}

\begin{remark}[Extension to smooth surfaces]
	\label{rem:alg-intersection-mult-plane:surface}
	If $X$ is a smooth surface and $C,C'\subset X$ are Cartier divisors meeting properly at $p$ with local equations
	$f,g\in\OO_{X,p}$, define
	\[
	I_p(C,C'):=\length_{\OO_{X,p}}\big(\OO_{X,p}/(f,g)\big).
	\]
	This specializes to Definition~\ref{def:alg-intersection-mult-plane:Ip} when $X=\A^2$.
\end{remark}

\begin{remark}[Serre's $\Tor$-formula (context)]
	\label{rem:alg-intersection-mult-plane:tor}
	For closed subschemes $X,Y$ of a Noetherian scheme $S$ and $p\in X\cap Y$,
	Serre defines
	\[
	I_p(X,Y)\ :=\ \sum_{i\ge 0}(-1)^i\,
	\length_{\OO_{S,p}}
	\Big(
	\Tor^{\OO_{S,p}}_{i}(\OO_{X,p},\OO_{Y,p})
	\Big).
	\]
	Here the formula is used when the displayed modules have finite length. Recall
	that a Noetherian local ring is \emph{Cohen--Macaulay} when its depth equals
	its Krull dimension; in dimension one this means that it contains a
	non-zero-divisor in its maximal ideal. Proper
	intersection in a regular ambient scheme does \emph{not}, by itself, force all
	higher $\Tor$ groups to vanish. In the present surface calculation, however,
	$X=(f)$ and $Y=(g)$ are effective Cartier divisors with no common local
	component. Resolving $R/(f)$ over the regular local ring $R=\OO_{S,p}$ by
	\[
	0\longrightarrow R\xrightarrow{\,f\,}R\longrightarrow R/(f)\longrightarrow0
	\]
	shows that $\Tor_i^R(R/(f),R/(g))=0$ for $i\ge2$, while
	$\Tor_1^R(R/(f),R/(g))$ is the annihilator of $f$ in $R/(g)$. Since $f$
	is a non-zero-divisor on the one-dimensional Cohen--Macaulay ring $R/(g)$,
	this annihilator is zero. Thus Serre's formula reduces here to
	$\length_R R/(f,g)$.
\end{remark}

\paragraph{Basic properties.}
\label{subsec:alg-intersection-mult-plane:properties}

Assume $p\in C\cap C'$ is isolated.

\begin{proposition}[Symmetry]
	\label{prop:alg-intersection-mult-plane:symmetry}
	$I_p(C,C')=I_p(C',C)$.
\end{proposition}

\begin{proof}
	$(f,g)=(g,f)$ as ideals in $R$, hence $R/(f,g)\cong R/(g,f)$ and the lengths agree.
\end{proof}

\begin{proposition}[Positivity and vanishing]
	\label{prop:alg-intersection-mult-plane:positivity}
	If $p\in C\cap C'$ is isolated, then $I_p(C,C')\ge 1$.
	With the standard convention $I_p(C,C')=0$ at a point outside
	$C\cap C'$, vanishing occurs exactly when $p\notin C\cap C'$.
\end{proposition}

\begin{proof}
	If $p\notin C\cap C'$, then either $f$ or $g$ is a unit in $R$, so $R/(f,g)=0$.
	If $p$ is an isolated intersection point, then $R/(f,g)$ is a nonzero Artinian local ring, hence has positive length.
\end{proof}

\begin{proposition}[Additivity with respect to components (a useful local form)]
	\label{prop:alg-intersection-mult-plane:additivity}
	Assume in $R$ that $f=f_1f_2$ and that $p$ is isolated in $V(f)\cap V(g)$. Then
	\[
	I_p\big((f=0),(g=0)\big)\ =\ I_p\big((f_1=0),(g=0)\big)+I_p\big((f_2=0),(g=0)\big).
	\]
\end{proposition}

\begin{proof}
	Put $A=R/(g)$.  The local ring $R$ of the smooth surface is regular of
	dimension two.  Since $g$ is a nonzero local equation for a curve, it is
	a non-zero-divisor in $R$; thus $A$ is a one-dimensional
	Cohen--Macaulay local ring.  In such a ring the zero-divisors are the
	union of its minimal primes.  If $f_i$ belonged to a minimal prime of
	$A$, the corresponding curve component of $V(g)$ would lie in
	$V(f_i)\subset V(f_1f_2)$, contradicting the assumption that the
	intersection is isolated at $p$.  Hence both $f_1$ and $f_2$ are
	non-zero-divisors in $A$, and $A/(f_i)$ and $A/(f_1f_2)$ have finite length.
	Multiplication by $f_2$ gives an exact sequence
	\[
	0\longrightarrow A/(f_1)\xrightarrow{\ \cdot f_2\ } A/(f_1f_2)
	\longrightarrow A/(f_2)\longrightarrow 0.
	\]
	For injectivity, if $f_2a$ is divisible by $f_1f_2$ in $A$, cancellation
	of the non-zero-divisor $f_2$ shows that $a$ is divisible by $f_1$.
	The last arrow is reduction modulo $f_2$; its kernel consists exactly
	of the classes divisible by $f_2$.  Length is additive in a short exact
	sequence of finite-length modules, so
	$\length A/(f_1f_2)=\length A/(f_1)+\length A/(f_2)$, which is the
	claimed intersection-multiplicity formula.
\end{proof}

\begin{proposition}[Transversality implies multiplicity $1$]
	\label{prop:alg-intersection-mult-plane:transverse}
	If $\nabla f(p)$ and $\nabla g(p)$ are linearly independent in $T_p\A^2$, then $I_p(C,C')=1$.
\end{proposition}

\begin{proof}
	After translating $p$ to $(0,0)$, the classes of $f$ and $g$ form a basis of
	$\mathfrak m/\mathfrak m^2$. By Nakayama's lemma they generate the maximal
	ideal $\mathfrak m$ of $R=k[x,y]_{(x,y)}$. Consequently
	$R/(f,g)=R/\mathfrak m\cong k$, and therefore $I_p(C,C')=1$.
	Equivalently, after completion one may choose regular parameters $u,v$ for which
	$(f,g)=(u,v)$; then
	\[
	\widehat R/(f,g)\widehat R\ \cong\ k[[u,v]]/(u,v)\ \cong\ k.
	\]
	Completion preserves the length of the finite-length quotient.
\end{proof}

\subsubsection{Tangent cones and computational methods}
\label{subsec:alg-intersection-mult-plane:tangent}

Assume $p=(0,0)$ after translation and set $R=k[x,y]_{(x,y)}$.
Write $f=\sum_{d\ge 0} f_d$ where $f_d$ is homogeneous of degree $d$.

\begin{definition}[Order and initial form]
	\label{def:alg-intersection-mult-plane:order-initial}
	Define the \emph{multiplicity} of $C=(f=0)$ at $0$ by
	\[
	m_0(C):=\min\{d\mid f_d\neq 0\},
	\]
	and the \emph{initial form} by $\in(f):=f_{m_0(C)}$.
	The \emph{projectivized tangent cone} is the zero scheme $(\in(f)=0)$ in $\PP^1$;
	its points, counted with multiplicity, record the tangent directions of the branches.
\end{definition}

\begin{theorem}[Tangent cone lower bound]
	\label{thm:alg-intersection-mult-plane:tangent-bound}
	Let $C=(f=0)$ and $C'=(g=0)$ meet isolatedly at $0$. Then
	\[
	I_0(C,C')\ \ge\ m_0(C)\,m_0(C').
	\]
	Moreover, equality holds if and only if $\in(f)$ and $\in(g)$ are coprime (equivalently, the two tangent cones
	share no common component).
\end{theorem}

\begin{proof}
	Blow up the origin. The local intersection formula for the strict transforms is
	\[
	I_0(C,C')=m_0(C)m_0(C')+
	\sum_{q\in E} I_q(\widetilde C,\widetilde C'),
	\]
	where $E\simeq\mathbb P^1$ is the exceptional divisor and the sum ranges over common points of the
	strict transforms on $E$. Every term in the sum is nonnegative, proving the lower bound. The points
	of $\widetilde C\cap E$ and $\widetilde C'\cap E$ are respectively the zeros of $\in(f)$ and
	$\in(g)$ on $\mathbb P^1$. Thus the sum is empty exactly when the two initial forms are coprime,
	which proves the equality criterion.
\end{proof}

\paragraph{Computation techniques and worked examples.}
\label{subsec:alg-intersection-mult-plane:computations}

\paragraph{(1) ``Eliminate $y$'' and compute a one-variable quotient.}
A frequent trick is: if one equation is $y-h(x)$ (or can be made so), substitute into the other.

\begin{example}[A general ``graph vs.\ curve'' computation]
	\label{ex:alg-intersection-mult-plane:graph}
	Let $C:(f=0)$ and $C':(y-h(x)=0)$ meet at $0$, where $h(0)=0$.
	Then in $R=k[x,y]_{(x,y)}$ we have
	\[
	R/(f,\ y-h(x))\ \cong\ k[x]_{(x)}/\big(f(x,h(x))\big).
	\]
	Hence
	\[
	I_0(C,C')\ =\ \dim_k k[x]_{(x)}/\big(f(x,h(x))\big)\ =\ \ord_{x=0}\big(f(x,h(x))\big).
	\]
\end{example}

\paragraph{(2) Puiseux/branch method (characteristic $0$).}
If $\operatorname{char}k=0$, after a suitable linear change of coordinates each
branch of a reduced plane curve admits a parametrization
$x=t^{n_i}$, $y=\varphi_i(t)$. If $f=0$ is reduced and has these branches, then
\[
I_0(C,C')=\sum_i \ord_{t=0}\big(g(t^{n_i},\varphi_i(t))\big).
\]
This is especially effective for singular curves (cusps, nodes).

\paragraph{(3) Resultant method.}
Suppose $f$ and $g$ are monic in $y$, have no common factor, and $(0,0)$ is the
only common zero lying over $x=0$. Then
\[
I_0(C,C')=\ord_{x=0}\big(\Res_y(f,g)\big).
\]
More generally, the order of the resultant is the sum of the local intersection
multiplicities at \emph{all} common points over $x=0$; a local resultant isolates
the contribution at $(0,0)$. The Sylvester determinant provides an effective
calculation for explicit polynomials.

\paragraph{(4) Standard monomials and Gr\"obner length.}
Compute a Gröbner basis for $(f,g)$ in $k[x,y]$ with a local monomial order;
the number of standard monomials gives the length.

\subsection{Examples and B\'ezout's theorem}
\label{subsec:alg-intersection-mult-plane:examples}

The examples below follow a single computational routine. Work in the local
ring at the chosen point, use one equation to eliminate a variable whenever
possible, and count the surviving powers of a local parameter. The resulting
basis computes the length and therefore the order of contact. After several
local calculations, B\'ezout's theorem explains why their sum is globally
constrained. This order---local calculation first, global conservation law
second---also makes clear why a missing point at infinity or an uncounted
tangency causes a naive point count to fail.

\begin{example}[Transversal lines]
	\label{ex:alg-intersection-mult-plane:lines}
	Let $C:(y=0)$ and $C':(x=0)$ at $0$.
	Then
	\[
	k[x,y]_{(x,y)}/(x,y)\cong k
	\quad\mathrel{\text{and hence}}\quad
	I_0(C,C')=1.
	\]
\end{example}

\begin{example}[Parabola and its tangent line]
	\label{ex:alg-intersection-mult-plane:parabola-tangent}
	Let $C:(y-x^2=0)$ and $C':(y=0)$ at $0$.
	Then
	\[
	k[x,y]_{(x,y)}/(y-x^2,\ y)\ \cong\ k[x]_{(x)}/(x^2),
	\]
	so a $k$-basis is $\{1,x\}$ and $I_0(C,C')=2$.
\end{example}

\begin{example}[Cusp and the $x$-axis]
	\label{ex:alg-intersection-mult-plane:cusp}
	Let $C:(y^2-x^3=0)$ and $C':(y=0)$ at $0$.
	Then
	\[
	k[x,y]_{(x,y)}/(y^2-x^3,\ y)\ \cong\ k[x]_{(x)}/(x^3),
	\]
	so $\{1,x,x^2\}$ is a $k$-basis and $I_0(C,C')=3$.
\end{example}

\begin{example}[Same tangent line, but still reduced]
	\label{ex:alg-intersection-mult-plane:same-tangent}
	Let $C:(y-x^2=0)$ and $C':(y-x^3=0)$ at $0$.
	Subtract the equations: $(y-x^2)-(y-x^3)=x^2-x^3=x^2(1-x)$.
	Since $1-x$ is a unit in $k[x]_{(x)}$, we get
	\[
	k[x,y]_{(x,y)}/(y-x^2,\ y-x^3)\ \cong\ k[x]_{(x)}/(x^2),
	\]
	hence $I_0(C,C')=2$.
	(Geometrically: same tangent line $y=0$, and contact of order $2$.)
\end{example}

\begin{example}[Why isolatedness is necessary]
	\label{ex:alg-intersection-mult-plane:strict-ineq}
	Let $C:(y^2=0)$ (a doubled line) and $C':(y=0)$ at $0$.
	Then $m_0(C)=2$ and $m_0(C')=1$, so $m_0(C)m_0(C')=2$.
	But
	\[
	k[x,y]_{(x,y)}/(y^2,\ y)\ \cong\ k[x]_{(x)}/(0),
	\]
	which is \emph{not} finite-dimensional over $k$:
	the intersection is not isolated (they share the entire line $y=0$), so $I_0(C,C')$ is not defined
	in the isolated-intersection sense. This example pinpoints why the ``isolated'' hypothesis is essential.
\end{example}

\begin{example}[Strict inequality in the tangent-cone bound]
	For the cusp $C:(y^2-x^3=0)$ and its tangent line $L:(y=0)$,
	we have $m_0(C)=2$ and $m_0(L)=1$, whereas
	Example~\ref{ex:alg-intersection-mult-plane:cusp} gives $I_0(C,L)=3$.
	The inequality $3>2$ records the common tangent direction.
\end{example}

\begin{example}[A node: two branches]
	\label{ex:alg-intersection-mult-plane:node}
	Let $C:(y^2-x^2-x^3=0)$ and $C':(y=0)$ at $0$.
	Modding out by $y$ gives $-x^2-x^3=-x^2(1+x)$, and $1+x$ is a unit in $k[x]_{(x)}$.
	Hence
	\[
	k[x,y]_{(x,y)}/(y^2-x^2-x^3,\ y)\ \cong\ k[x]_{(x)}/(x^2),
	\]
	so $I_0(C,C')=2$.
\end{example}

\subsubsection{Global relation: B\'ezout's theorem}
\index{B\'ezout's theorem}
\label{subsec:alg-intersection-mult-plane:bezout}

Let $C=(F=0)$ and $C'=(G=0)$ be projective plane curves in $\mathbb{P}^2_k$ of degrees
$d=\deg F$ and $d'=\deg G$ with no common component. For each $p\in C\cap C'$, one defines
$I_p(C,C')$ by passing to an affine chart around $p$ and using Definition~\ref{def:alg-intersection-mult-plane:Ip}.
This number is independent of the chart: two affine charts induce isomorphic local
rings at $p$, multiplication of a defining equation by a local unit does not change
the generated ideal, and length is preserved by a local-ring isomorphism.

\begin{theorem}[B\'ezout]
	\label{thm:alg-intersection-mult-plane:bezout}
	If $C$ and $C'$ have no common component, then
	\[
	\sum_{p\in C\cap C'} I_p(C,C')\ =\ d\,d'.
	\]
\end{theorem}

\begin{proof}[Intersection product computation]
	The scheme-theoretic intersection $Z=C\cap C'$ is $0$-dimensional.
	Its degree is the sum of the lengths of its local rings: each isolated
	point contributes its scheme-theoretic multiplicity, rather than merely
	one point. By the local definition above,
	\[
	\deg(Z)=\sum_{p\in Z}\length(\OO_{Z,p})=\sum_{p\in C\cap C'} I_p(C,C').
	\]
	A degree-$d$ homogeneous equation is a section of
	$\OO_{\PP^2}(d)$, hence its divisor class is $dH$, where $H$ is the
	class of a line. The local--global intersection theorem identifies
	$\deg Z$ with the product of these divisor classes; this is the
	intersection-theoretic input, rather than a consequence of the local
	length calculation alone \cite{Fulton1998Intersection,Hartshorne1977}.
	Two general lines meet in one reduced point, so $H^2=[\mathrm{pt}]$.
	Therefore, in $A^\bullet(\PP^2)=\Z[H]/(H^3)$,
	\[
	[C]\cdot[C']=dd'\,H^2 = dd'\,[\mathrm{pt}],
	\]
	hence $\deg(Z)=dd'$.
\end{proof}

\begin{example}[Verifying B\'ezout in a concrete case]
	\label{ex:alg-intersection-mult-plane:bezout-concrete}
	Let $C:(XZ-Y^2=0)$ (a conic) and $L:(Z=0)$ (the line at infinity) in $\mathbb{P}^2$.
	On $Z=0$, the conic equation becomes $-Y^2=0$, so $Y=0$ and the unique intersection point is $p=[1:0:0]$.
	In the affine chart $X\neq 0$ with coordinates $u=Y/X$, $v=Z/X$, the equations become
	\[
	C:\ v-u^2=0,\qquad L:\ v=0.
	\]
	Thus
	\[
	I_p(C,L)=\dim_k k[u,v]_{(u,v)}/(v-u^2,\ v)\cong \dim_k k[u]_{(u)}/(u^2)=2.
	\]
	Hence $\sum I_p=2 = (\deg C)(\deg L)=2\cdot 1$.
\end{example}

\subsubsection{Compatibility with surface intersection theory and exercises}
\label{subsec:alg-intersection-mult-plane:compatibility}

Let $X$ be a smooth projective surface over $\C$.
For Cartier divisors $D_1,D_2$ meeting properly, define the global intersection number by
\[
(D_1\cdot D_2)\ :=\ \sum_{p\in X} I_p(D_1,D_2).
\]

\begin{proposition}[Agreement with the Chern class / cup-product pairing]
	\label{prop:alg-intersection-mult-plane:chern}
	On a smooth complex projective surface $X$,
	\[
	(D_1\cdot D_2)\ =\ \int_X \alpha_1\wedge\alpha_2,
	\]
	where $\alpha_i$ is any closed de Rham representative of
	$c_1(\OO_X(D_i))$.  Equivalently, this integer is the evaluation on $[X]$
	of the cup product
	$c_1(\OO_X(D_1))\smile c_1(\OO_X(D_2))\in H^4(X,\Z)$.
\end{proposition}

\begin{proof}
	The map $D\mapsto c_1(\OO_X(D))$ sends linearly equivalent divisors to the same
	cohomology class; the converse is generally false because distinct elements of
	$\operatorname{Pic}^0(X)$ have the same first Chern class. For properly intersecting
	Cartier divisors, the local--global intersection theorem identifies the cup-product
	degree
	\[
	\big\langle c_1(\OO_X(D_1))\smile c_1(\OO_X(D_2)),[X]\big\rangle
	\]
	with the degree of the scheme-theoretic intersection. Its local contribution at
	$p$ is the length $\length(\OO_{X,p}/(f,g))=I_p(D_1,D_2)$.
	Summing these local lengths gives the displayed equality; see
	\cite{Fulton1998Intersection,GriffithsHarris1994}.
\end{proof}

\paragraph{Exercises.}
\label{subsec:alg-intersection-mult-plane:exercises}

\begin{exercise}[A family of tangencies]
	For an integer $m\ge2$, let $C_m:(y-x^m=0)$ and $L:(y=0)$ in $\mathbb{A}^2$.
	Compute $I_0(C_m,L)$ and interpret the answer as the order of contact.
\end{exercise}

\begin{exercise}[Testing transversality in two ways]
	\label{exr:selected-transversality}
	At the origin let $C=(y-x^2=0)$ and $C'=(x+y-y^2=0)$. Compute
	\[
	\det\!\left(\frac{\partial(f,g)}{\partial(x,y)}\right)(0)
	\]
	for the displayed defining polynomials $f,g$, and deduce that the
	curves meet transversely. Then substitute $y=x^2$ in $g$ and compute
	$\mathcal O_{\mathbb A^2,0}/(f,g)$ directly to verify $I_0(C,C')=1$.
\end{exercise}

\begin{exercise}[Tangent cones and equality]
	Let $C=(f=0)$ and $C'=(g=0)$ meet isolatedly at $0$.
	\begin{enumerate}
		\item Prove $I_0(C,C')\ge m_0(C)m_0(C')$.
		\item Give an example where $m_0(C)m_0(C')=1$ but $I_0(C,C')>1$.
		(Hint: $y-x^2$ and $y-x^3$.)
	\end{enumerate}
\end{exercise}

\begin{exercise}[Resultant computation]
	Let $C:(y^2-x^3=0)$ and $C':(y-x^2=0)$ at $0$.
	Compute $I_0(C,C')$ by substituting $y=x^2$, and then verify the same answer using a resultant in $y$.
\end{exercise}

\begin{exercise}[B\'ezout check]
	Let $C:(Y^2Z=X^3+XZ^2)$ (a cubic) and $L:(Z=0)$ in $\mathbb{P}^2$.
	Compute all intersection points and multiplicities, and verify $\sum I_p=3$.
\end{exercise}

\clearpage
\part{Special Functions, Galois Theory, and Arithmetic}
\label{part:bridges}

\section{Theta Functions and Weierstrass Theory}
\index{theta functions}\index{Weierstrass theory}

\providecommand{\C}{\mathbb{C}}
\providecommand{\Z}{\mathbb{Z}}
\providecommand{\Jac}{\operatorname{Jac}}
\providecommand{\im}{\operatorname{Im}}
\providecommand{\re}{\operatorname{Re}}
\providecommand{\OO}{\mathcal{O}}
\providecommand{\Res}{\operatorname{Res}}
\providecommand{\Vol}{\operatorname{Vol}}

The Jacobian constructed earlier has the form $\C^g/\Lambda$. To obtain a
distinguished divisor on it, we now construct a function on $\C^g$ whose
transformation under $\Lambda$ is controlled rather than strictly periodic.
The theta function $\theta(\,\cdot\mid\Omega)$ is entire on $\C^g$ and obeys
explicit quasi-periodicity. These transformation laws are precisely the descent
data that define a holomorphic section of a line bundle on $\Jac(X)$.
For period matrices and theta line bundles we use \cite{Mumford1983Tata1,Fay1973,FarkasKra1992};
the genus-one normalizations are checked against \cite{DLMFTheta2026,DLMFWeierstrass2026}.

We proceed in three layers:
\begin{enumerate}
	\item analytic: the Riemann theta series converges and satisfies explicit
	quasi-periodicity,
	\item geometric: quasi-periodicity is exactly the cocycle of a line bundle on
	$A=\C^g/\Lambda$ and the zero set descends to the \emph{theta divisor},
	\item genus $1$: the abstract theta becomes the classical Jacobi theta functions
	and recovers Weierstrass theory via logarithmic derivatives.
\end{enumerate}

\subsection{From Kepler's ellipse to Riemann's global viewpoint}
\label{subsec:kepler-riemann-worldview}
\index{elliptic integrals!Kepler ellipse}\index{Kepler equation}
\index{Riemann!global viewpoint}

The word \emph{elliptic} is historical before it is algebro-geometric: these
integrals first arose from measurements associated with an ellipse.  A
Keplerian orbit supplies the most familiar ellipse, but a distinction is
essential.  The time law for the Newtonian two-body problem is elementary
after the eccentric anomaly is introduced; the \emph{arc length} along the
same orbit is an elliptic integral.  Thus one should not claim that Kepler's
equation itself requires elliptic functions.  The genuine bridge is
\[
	\begin{aligned}
	\text{Kepler ellipse}
	&\longrightarrow \text{elliptic integral}
	\longrightarrow \text{two-sheeted curve}\\
	&\longrightarrow \text{period lattice and complex torus}.
	\end{aligned}
\]
This subsection carries out that bridge and then explains, without asserting
a false line of historical causation, how it exemplifies Riemann's broader
method.

\subsubsection{The two-body calculation and the precise place where an elliptic integral appears}

Let a particle move in a plane under the Newtonian potential $-\mu/r$, with
specific angular momentum $h=r^2\dot\theta$ and specific energy
\[
	\mathcal E=\frac12\dot r^{\,2}+\frac{h^2}{2r^2}-\frac{\mu}{r}.
\]
Multiplication by $2r^2$ gives
\[
	r^2\dot r^{\,2}=2\mathcal E r^2+2\mu r-h^2,
	\qquad
	dt=\frac{r\,dr}{\sqrt{2\mathcal E r^2+2\mu r-h^2}}.
\]
For a bound orbit $\mathcal E<0$, the polynomial under the square root is
quadratic.  The substitution by the eccentric anomaly therefore evaluates
the time integral in elementary trigonometric functions.  Writing
\[
	x=a(\cos E-e),\qquad
	y=a\sqrt{1-e^2}\sin E,
\]
where $a>0$ is the semimajor axis and $0\leq e<1$ the eccentricity, one obtains
\[
	n(t-t_0)=E-e\sin E,
	\qquad n=\sqrt{\frac{\mu}{a^3}}.
\]
This is Kepler's equation.  It is transcendental to \emph{invert} for $E$ as
a function of time, but its defining integral is not an elliptic integral.

Now ask a different geometric question: how much distance has the particle
traveled along the elliptical orbit? Differentiating the parametrization gives
\[
	\frac{dx}{dE}=-a\sin E,\qquad
	\frac{dy}{dE}=a\sqrt{1-e^2}\cos E,
\]
and hence
\[
	ds
	=\sqrt{\left(\frac{dx}{dE}\right)^2+
	          \left(\frac{dy}{dE}\right)^2}\,dE
	=a\sqrt{1-e^2\cos^2E}\,dE.
\]
The traveled distance from periapsis to anomaly $E$ is therefore
\[
	s(E)=a\int_0^E\sqrt{1-e^2\cos^2u}\,du.
\]
After $\phi=\pi/2-u$, this is a difference of incomplete elliptic integrals
of the second kind,
\[
	\mathsf E(\phi,k)=\int_0^\phi\sqrt{1-k^2\sin^2\vartheta}\,d\vartheta,
	\qquad k=e.
\]
In particular, the circumference of the full Kepler ellipse is
\[
	L=4a\,\mathsf E(e),
	\qquad
	\mathsf E(e)=\int_0^{\pi/2}\sqrt{1-e^2\sin^2\vartheta}\,d\vartheta.
\]
This is the literal origin of the adjective ``elliptic''; compare the standard
arc-length formulas in \cite{DLMFElliptic2026}.

\subsubsection{From the square root to a two-sheeted Riemann surface}

Put $t=\sin\vartheta$.  The complete elliptic integrals of the first and
second kinds become
\begin{align*}
	\mathsf K(k)
	&=\int_0^1\frac{dt}{\sqrt{(1-t^2)(1-k^2t^2)}},\\
	\mathsf E(k)
	&=\int_0^1\frac{1-k^2t^2}{
	\sqrt{(1-t^2)(1-k^2t^2)}}\,dt.
\end{align*}
The same quartic polynomial governs both expressions.  Introduce the affine
curve
\[
	C_k:\qquad y^2=(1-x^2)(1-k^2x^2),
	\qquad 0<k<1.
\]
If one insists on viewing the square root as a function of the $x$-plane, its
sign changes when a path winds around a branch value.  On $C_k$, however,
$y$ is a single-valued coordinate.  The two signs of the square root are the
two sheets of the projection $(x,y)\mapsto x$, joined at
$x=\pm1,\pm k^{-1}$.

\begin{proposition}[The differential behind the elliptic integral]
	\label{prop:kepler-elliptic-differential}
	Let $\overline C_k$ be the smooth projective compactification of $C_k$.
	Then $\overline C_k$ has genus one, and
	\[
		\omega=\frac{dx}{y}
	\]
extends to a nowhere-vanishing holomorphic $1$-form on
	$\overline C_k$.  The first-kind integral is an integral of $\omega$, while
	$\mathsf E(k)$ is an integral of the meromorphic differential
	$(1-k^2x^2)\omega$.
\end{proposition}

\begin{proof}
	The projection $\pi:\overline C_k\to\PP^1$ has degree two and four simple
	branch values. Riemann--Hurwitz gives
	\[
		2g(\overline C_k)-2=2(-2)+4=0,
	\]
	so $g(\overline C_k)=1$. Near a finite branch point $a$, write
	$x=a+t^2$ and $y=t\,v(t)$ with $v(0)\ne0$. Then
	\[
		\frac{dx}{y}=\frac{2t\,dt}{t\,v(t)}=\frac{2}{v(t)}\,dt,
	\]
	which is holomorphic and nonzero. Because the defining polynomial has even
	degree four, there are two points over infinity. With $s=1/x$ and
	$Y=y/x^2$, the equation becomes
	$Y^2=(s^2-1)(s^2-k^2)$ and
	\[
		\frac{dx}{y}=-\frac{ds}{Y},
	\]
	which is again holomorphic and nonzero at $s=0$. Thus $\omega$ is global and
	nowhere vanishing. The two displayed integrals follow directly from their
	definitions after $x=t$.
\end{proof}

The distinction between the two differentials is conceptually useful.
$dx/y$ is of the first kind: it is holomorphic everywhere and supplies the
uniformizing coordinate.  The arc-length differential
$(1-k^2x^2)dx/y$ is of the second kind: it is meromorphic on the compact curve
but has zero residues.  They live on the same Riemann surface, so the geometry
organizes several apparently different integrals at once.

\subsubsection{Periods, inversion, and the torus}

Fix $p_0\in\overline C_k$ and define on the universal cover
\[
	u(p)=\int_{p_0}^{p}\omega.
\]
Changing the path by a closed cycle changes $u$ by a period. If $(a,b)$ is an
oriented homology basis, put
\[
	\omega_1=\int_a\omega,\qquad
	\omega_2=\int_b\omega,\qquad
	\Lambda=\Z\omega_1+\Z\omega_2.
\]
The Riemann bilinear relation gives
$\operatorname{Im}(\overline{\omega_1}\omega_2)>0$, so the two periods are
real-linearly independent and $\Lambda$ is a lattice.  Therefore
\[
	p\longmapsto \int_{p_0}^{p}\omega\pmod\Lambda
\]
defines a biholomorphism $\overline C_k\cong\C/\Lambda$.  Conversely, locally
inverting $u=\int dx/y$ gives functions $x(u)$ and $y(u)$ satisfying the
quartic equation and periodic under $u\mapsto u+\lambda$ for
$\lambda\in\Lambda$.  A birational change of variables converts the quartic
to a Weierstrass cubic
\[
	Y^2=4X^3-g_2X-g_3,
\]
and the inverse coordinate becomes $X=\wp(u)$, with
$Y=\wp'(u)$.  This explains why elliptic \emph{functions} are doubly periodic:
they are functions on a complex torus obtained by inverting an elliptic
integral.

\begin{example}[Real and imaginary periods]
	For $0<k<1$, cycles surrounding adjacent branch cuts may be chosen so that,
	up to orientation,
	\[
		\omega_1=2\mathsf K(k),
		\qquad
		\omega_2=2i\mathsf K(k'),
		\qquad k'=\sqrt{1-k^2}.
	\]
	Thus $\tau=\omega_2/\omega_1=i\mathsf K(k')/\mathsf K(k)$ lies in the upper
	half-plane. The real period records motion around one branch cut, and the
	imaginary period records the independent cycle crossing between the cuts.
The flat parallelogram spanned by these periods is the torus hidden behind the
square root.
\end{example}

\subsubsection{What is specifically Riemann's contribution?}

Euler, Legendre, Abel, and Jacobi developed elliptic and Abelian integrals
before Riemann. Riemann's decisive move was not to discover the ellipse
integral, but to change the object on which it is studied. Rather than treating
$\sqrt{P(x)}$ as a troublesome multivalued expression on the plane, he treated
its branches as one single-valued object on a surface. His 1851 doctoral
dissertation, \emph{Grundlagen f\"ur eine allgemeine Theorie der Functionen
einer ver\"anderlichen complexen Gr\"osse}, placed complex functions,
branched surfaces, conformal mapping, and the Dirichlet principle in one
framework \cite{Riemann1851,Gunning1962}.

The calculation above displays five habits that recur throughout his work.
First, replace a formula by its natural domain: the square root belongs on
$C_k$, not on a cut plane. Second, pass from local data to global invariants:
local primitives acquire periods measured by homology. Third, use a
differential form as the invariant object: $dx/y$ survives coordinate changes
even when neither $x$ nor $y$ does. Fourth, turn inversion into geometry: the
inverse of the integral lives naturally on $\C/\Lambda$. Fifth, use a
variational or PDE construction to produce analytic objects whose topology
then controls their global behavior.

\paragraph{The 1851 dissertation and PDE.}
Harmonic functions are real parts of holomorphic functions locally. Globally,
their conjugates can acquire periods, exactly as the integral of $dx/y$ does.
Riemann used the Dirichlet principle to solve boundary-value problems and to
construct conformal maps. In modern language, one minimizes the energy
\[
	\mathcal E(u)=\int_X|du|^2\,dA
\]
subject to boundary data; the Euler--Lagrange equation is $\Delta u=0$.
The original justification of the minimizer contained a gap later repaired by
Hilbert and modern functional analysis, but the organizing idea remains
fundamental: an analytic PDE constructs the coordinate, while topology records
its periods and obstructions.

\paragraph{Riemannian geometry.}
In his 1854 habilitation lecture Riemann replaced the fixed Euclidean metric by
a variable quadratic line element
\[
	ds^2=\sum_{i,j}g_{ij}(x)\,dx^i dx^j.
\]
This is not a consequence of elliptic integrals. The methodological continuity
is that geometry is encoded infinitesimally and then integrated globally.
On a Riemann surface, a conformal metric, its curvature, the Laplacian, and the
complex structure interact; uniformization later turns this interaction into
the three constant-curvature models \cite{Riemann1854Geometry}.

\paragraph{The Riemann integral.}
Riemann's work on trigonometric series required a definition of integration
that tolerated discontinuous functions. Partitioning an interval, controlling
oscillation on subintervals, and passing to a limit created what is now called
the Riemann integral. This development is historically distinct from Abelian
integrals. The common viewpoint is again global reconstruction from controlled
local pieces; one should not identify the two theories merely because both use
the word ``integral'' \cite{Riemann1854Trig}.

\paragraph{Theta functions, zeta, and the Riemann hypothesis.}
Periods of holomorphic forms lead from a genus-$g$ surface to its Jacobian
$\C^g/\Lambda$ and then to theta functions. In a different problem, Riemann's
1859 memoir continued $\zeta(s)$ and related its zeros to the distribution of
primes. The functional equation is efficiently obtained from the Jacobi theta
function by a Mellin transform, as shown later in this chapter. This is a real
technical bridge between theta and zeta, but the Riemann hypothesis is not a
theorem about the period lattice of a single elliptic curve. It asserts that
the nontrivial zeros of $\zeta(s)$ lie on $\Re s=1/2$ and arose in Riemann's
explicit-formula study of prime counting \cite{Riemann1859}.

\paragraph{One worldview, not one chain of deductions.}
The historical works just listed do not form the implication
``Kepler orbit, therefore Riemann surface, therefore Riemannian geometry,
therefore the Riemann hypothesis.'' Their unity is methodological. Riemann
repeatedly enlarged the domain until a multivalued or singular expression
became natural, encoded local behavior by invariant differential data, and
used integration, topology, or PDE to recover global structure. The curve
$y^2=P(x)$ and the form $dx/y$ are the cleanest finite-dimensional laboratory
for that method: branch points become topology, integrals become periods,
periods become a lattice, and the lattice becomes a geometric space.

\paragraph{Prerequisite check: the Gaussian in the theta series.}
For $\tau=X+iY$ with $Y>0$ and $z=u+iv$, the modulus of the $n$th
genus-one summand is
\[
\bigl|e^{\pi i n^2\tau+2\pi i n z}\bigr|
=e^{-\pi Yn^2-2\pi nv}
=e^{\pi v^2/Y}e^{-\pi Y(n+v/Y)^2}.
\]
On a compact set of $z$ values, $v$ is bounded, so these terms have a
summable Gaussian majorant independent of $z$. This is the analytic
reason the series defines an entire function; the next question is
how its values change under lattice translation.

\subsection{From a Riemann surface to its Jacobian: period matrices and lattices}

Let $X$ be a compact Riemann surface of genus $g$.
Choose a symplectic basis of homology
\[
\{a_1,\dots,a_g,b_1,\dots,b_g\}\subset H_1(X,\Z)
\quad\text{with}\quad
a_i\cdot a_j=b_i\cdot b_j=0,\ \ a_i\cdot b_j=\delta_{ij}.
\]
Let $\{\omega_1,\dots,\omega_g\}$ be a basis of holomorphic $1$--forms on $X$,
normalized by
\[
\int_{a_j}\omega_i=\delta_{ij}.
\]
Define the \emph{period matrix} $\Omega=(\Omega_{ij})$ by
\[
\Omega_{ij}:=\int_{b_j}\omega_i.
\]
Then $\Omega$ is symmetric and has positive definite imaginary part.

\begin{theorem}[Riemann bilinear relations: key properties of $\Omega$]
	\label{thm:riemann-bilinear}
	With the above normalization, $\Omega$ satisfies
	\[
	\Omega^t=\Omega,
	\qquad
	\Im\Omega>0.
	\]
\end{theorem}

\begin{proof}
	\emph{Symmetry.}
	For holomorphic forms $\omega_i,\omega_j$, the $2$--form $\omega_i\wedge\omega_j$ vanishes
	identically on a complex $1$--manifold. Applying the bilinear relations gives
	\[
	\sum_{k=1}^g
	\left(
	\int_{a_k}\omega_i\int_{b_k}\omega_j
	-
	\int_{b_k}\omega_i\int_{a_k}\omega_j
	\right)=0.
	\]
	Using $\int_{a_k}\omega_i=\delta_{ik}$ gives $\Omega_{ij}-\Omega_{ji}=0$.
	
	\emph{Positivity.}
	For any $0\neq c=(c_1,\dots,c_g)\in\C^g$, let $\eta=\sum_i c_i\omega_i$. Then
	\[
	i\int_X \eta\wedge\overline{\eta}>0,
	\]
	and the bilinear relations for $(\eta,\overline{\eta})$ yield
	\[
	 i\int_X\eta\wedge\overline\eta
	 =2c^t(\Im\Omega)\overline c>0.
	\]
	Thus $\Im\Omega$ is positive definite. With the convention that Hermitian
	forms are linear in the first variable, the displayed quadratic form is the
	transpose of the more familiar $\overline c^{\,t}(\Im\Omega)c$; symmetry of
	$\Im\Omega$ makes the two expressions equal.
\end{proof}

Define the lattice
\[
\Lambda:=\Z^g+\Omega\,\Z^g\subset\C^g,
\qquad
A:=\C^g/\Lambda.
\]
Under the identification $\C^g\simeq H^0(X,\Omega_X^1)^\vee$ by integration, $A$ is the
\emph{Jacobian} $\Jac(X)$.
Choosing a base point $P_0\in X$, the Abel--Jacobi map is
\[
u(P):=\left(\int_{P_0}^{P}\omega_1,\dots,\int_{P_0}^{P}\omega_g\right)\ \ (\mathrm{mod}\ \Lambda).
\]

\begin{example}[Genus $1$ period lattice from an elliptic integral]
	\label{ex:g1-periods}
	Let $E$ be the smooth projective cubic
	\[
	E:\quad y^2=4(x-e_1)(x-e_2)(x-e_3), \qquad e_1,e_2,e_3\in\R,\ \ e_1>e_2>e_3.
	\]
	On $E$ the holomorphic differential is $\omega=\frac{dx}{y}$.
	In this genus-one example, $\omega_1$ and $\omega_2$ below denote
	period \emph{numbers}, not the normalized basis forms used above.
	Choose the square roots in the following integrals to be positive. A standard
	homology basis $(a,b)$ can be oriented so that
	\[
\begin{aligned}
	\omega_1&:=\int_a \omega
	 =2\int_{e_3}^{e_2}\frac{dx}{\sqrt{4(e_1-x)(e_2-x)(x-e_3)}}>0,\\
	\omega_2&:=\int_b \omega
	 =2i\int_{e_2}^{e_1}\frac{dx}{\sqrt{4(e_1-x)(x-e_2)(x-e_3)}}\in i\mathbb R_{>0}.
\end{aligned}
	\]
	Thus the period lattice is $\Lambda=\Z\omega_1+\Z\omega_2$ and
	$\tau=\omega_2/\omega_1\in\mathbb H$. The factor $2$ records the outward and
	return paths on the two sheets; the factor $i$ in the second integral comes
	from continuing the chosen square root across the adjacent branch interval.
\end{example}

\subsection{The Riemann theta function: convergence and quasi-periodicity}
\index{theta functions!quasi-periodicity}

\begin{definition}[Riemann theta function]
	\label{def:riemann-theta}
	Let $\Omega$ be symmetric with $\Im\Omega>0$. Define
	\[
	\theta(z\,|\,\Omega)
	:=
	\sum_{n\in\Z^g}
	\exp\!\Bigl(\pi i\, n^t\Omega n + 2\pi i\, n^t z\Bigr),
	\qquad z\in\C^g.
	\]
\end{definition}

\begin{proposition}[Absolute and locally uniform convergence]
	\label{prop:theta-conv}
	The series for $\theta(z\,|\,\Omega)$ converges absolutely and uniformly on compact subsets
	of $\C^g$. Hence $\theta(\,\cdot\,|\,\Omega)$ is entire.
\end{proposition}

\begin{proof}
	Write $\Omega=X+iY$ with $X=\re\Omega$ and $Y=\im\Omega>0$.
	Then
	\[
	\Bigl|\exp\!\bigl(\pi i\, n^t\Omega n + 2\pi i\, n^t z\bigr)\Bigr|
	=
	\exp\!\bigl(-\pi\, n^t Y n - 2\pi\, n^t \Im z\bigr).
	\]
	Fix a compact $K\subset\C^g$ and let $M:=\sup_{z\in K}\|\Im z\|$.
	Then $-2\pi\,n^t\Im z\le 2\pi M\|n\|$.
	Since $Y$ is positive definite, $\exists c>0$ with $n^tYn\ge c\|n\|^2$.
	Hence for $z\in K$,
	\[
	\Bigl|\exp(\cdots)\Bigr|\le \exp\!\bigl(-\pi c\|n\|^2+2\pi M\|n\|\bigr),
	\]
	and the right-hand side is summable over $\Z^g$. Apply the Weierstrass M-test.
\end{proof}

\begin{proposition}[Quasi-periodicity]
	\label{prop:theta-quasi}
	For $m,n\in\Z^g$ we have
	\begin{align*}
		\theta(z+m\,|\,\Omega)&=\theta(z\,|\,\Omega),\\
		\theta(z+\Omega n\,|\,\Omega)
		&=
		\exp\!\bigl(-\pi i\, n^t\Omega n - 2\pi i\, n^t z\bigr)\ \theta(z\,|\,\Omega).
	\end{align*}
	Equivalently, for $\lambda=m+\Omega n\in\Lambda$,
	\[
	\theta(z+\lambda\,|\,\Omega)=\alpha_\lambda(z)\,\theta(z\,|\,\Omega),
	\quad
	\alpha_{m+\Omega n}(z):=\exp\!\bigl(-\pi i\, n^t\Omega n - 2\pi i\, n^t z\bigr).
	\]
\end{proposition}

\begin{proof}
	For the first identity, write $k\in\Z^g$ for the summation index. Then
	$e^{2\pi i k^t m}=1$ for every $m\in\Z^g$, so each summand is unchanged.
	For the second, compute and reindex:
	\begin{align*}
		\theta(z+\Omega n)
		&=\sum_{k\in\Z^g} \exp\!\bigl(\pi i\, k^t\Omega k + 2\pi i\, k^t z + 2\pi i\,k^t\Omega n\bigr)\\
		&=\sum_{k'\in\Z^g} \exp\!\bigl(\pi i\,(k'-n)^t\Omega(k'-n) + 2\pi i\,(k'-n)^t z + 2\pi i\,(k'-n)^t\Omega n\bigr)\\
		&=e^{-\pi i n^t\Omega n-2\pi i n^t z}\sum_{k'\in\Z^g}\exp\!\bigl(\pi i\,k'^t\Omega k'+2\pi i\,k'^t z\bigr)\\
		&=e^{-\pi i n^t\Omega n-2\pi i n^t z}\,\theta(z).
	\end{align*}
\end{proof}

\begin{example}[Concrete $g=1$ check of quasi-periodicity]
	\label{ex:g1-quasi-check}
	Let $\tau\in\mathbb{H}$ and define
	\[
	\vartheta_3(z,\tau):=\sum_{n\in\Z} e^{\pi i n^2\tau+2\pi i n z}.
	\]
	Then:
	\begin{enumerate}
		\item Period $1$:
		\[
		\vartheta_3(z+1,\tau)=\sum_{n\in\Z} e^{\pi i n^2\tau+2\pi i n(z+1)}
		=\sum_{n\in\Z} e^{\pi i n^2\tau+2\pi i n z}\,e^{2\pi i n}
		=\vartheta_3(z,\tau).
		\]
		\item Quasi-period $\tau$:
		\begin{align*}
			\vartheta_3(z+\tau,\tau)
			&=\sum_{n\in\Z} e^{\pi i n^2\tau+2\pi i n(z+\tau)}
			=\sum_{n\in\Z} e^{\pi i n^2\tau+2\pi i n z}\,e^{2\pi i n\tau}.
		\end{align*}
		Now reindex $n=m-1$:
		\begin{align*}
			\vartheta_3(z+\tau,\tau)
			&=\sum_{m\in\Z} e^{\pi i (m-1)^2\tau+2\pi i (m-1)z}\,e^{2\pi i (m-1)\tau}\\
			&=e^{-\pi i\tau-2\pi i z}\sum_{m\in\Z} e^{\pi i m^2\tau+2\pi i m z}
			=e^{-\pi i\tau-2\pi i z}\,\vartheta_3(z,\tau).
		\end{align*}
		This is exactly Proposition~\ref{prop:theta-quasi} for $g=1$ and $n=1$.
	\end{enumerate}
\end{example}

\subsection{\texorpdfstring{Theta as a section of a line bundle on $\Jac(X)$}{Theta as a section of a line bundle on the Jacobian}}

Let $A=\C^g/\Lambda$. Theta does not descend to a \emph{function} on $A$, but it
descends to a \emph{section of a holomorphic line bundle}.
Indeed, a function on the quotient must give the same value at $z$ and
$z+\lambda$.  In genus one, however,
\[
\theta(z+1\mid\tau)=\theta(z\mid\tau),\qquad
\theta(z+\tau\mid\tau)
=e^{-\pi i\tau-2\pi i z}\theta(z\mid\tau).
\]
The second multiplier depends on $z$, so theta is not an ordinary
periodic function.  A section may have different coordinate values in
different trivializations of the same fiber; the multiplier tells us
precisely how to identify those coordinates.  The quotient construction
below performs that identification.  This is the same gluing mechanism
used for $\mathcal O(1)$ in Example~\ref{ex:O1-sections}, now with the
lattice acting in place of a two-chart cover.

\begin{definition}[Line bundle from automorphy factors]
	\label{def:automorphy}
	Let $\alpha_\lambda(z)$ be the factor in Proposition~\ref{prop:theta-quasi}.
	Define an equivalence relation on $\C^g\times\C$ by
	\[
	(z,w)\sim (z+\lambda,\ \alpha_\lambda(z)\,w),\qquad \lambda\in\Lambda.
	\]
	The quotient
	\[
	\mathcal{L}:=(\C^g\times\C)/\sim \ \longrightarrow\ A
	\]
	is a holomorphic line bundle on $A$.
\end{definition}

\begin{lemma}[Cocycle condition]
	\label{lem:cocycle}
	The factors $\alpha_\lambda(z)$ satisfy
	\[
	\alpha_{\lambda+\mu}(z)=\alpha_\lambda(z+\mu)\alpha_\mu(z),
	\qquad \lambda,\mu\in\Lambda.
	\]
	Hence Definition~\ref{def:automorphy} is consistent.
\end{lemma}

\begin{proof}
	Write $\lambda=m+\Omega n$ and $\mu=m'+\Omega n'$ with $m,m',n,n'\in\Z^g$.
	Using $\alpha_{m+\Omega n}(z)=\exp(-\pi i\,n^t\Omega n-2\pi i\,n^t z)$, compute:
	\[
	\alpha_{\lambda+\mu}(z)
	=
	\exp\!\bigl(-\pi i\,(n+n')^t\Omega(n+n')-2\pi i\,(n+n')^t z\bigr).
	\]
	Also
	\[
	\alpha_\lambda(z+\mu)\alpha_\mu(z)
	=
	\exp\!\bigl(-\pi i\,n^t\Omega n-2\pi i\,n^t(z+m'+\Omega n')\bigr)
	\exp\!\bigl(-\pi i\,n'^t\Omega n'-2\pi i\,n'^t z\bigr).
	\]
	Since $n^t m'\in\Z$, the factor $e^{-2\pi i n^t m'}=1$, and symmetry of $\Omega$
	matches the cross-terms, so both sides coincide.
\end{proof}

\begin{proposition}[Theta descends to a global section]
	\label{prop:theta-section}
	The function $\theta(z\,|\,\Omega)$ defines a holomorphic global section
	$s_\theta\in H^0(A,\mathcal{L})$ by $s_\theta([z])=[z,\theta(z)]$.
\end{proposition}

\begin{proof}
	If $z$ is replaced by $z+\lambda$, then Proposition~\ref{prop:theta-quasi} gives
	$\theta(z+\lambda)=\alpha_\lambda(z)\theta(z)$, hence
	\[
	(z+\lambda,\theta(z+\lambda))\sim (z+\lambda,\alpha_\lambda(z)\theta(z))\sim (z,\theta(z)).
	\]
	So $s_\theta$ is well-defined. Holomorphicity follows from that of $\theta$.
\end{proof}

\begin{definition}[Theta divisor]
	\label{def:theta-divisor}
	The \emph{theta divisor} on $A$ is the zero locus of $s_\theta$:
	\[
	\Theta_A := \{[z]\in A \mid \theta(z\,|\,\Omega)=0\}.
	\]
\end{definition}

\begin{remark}[Geometric meaning (Riemann vanishing)]
	On $\Jac(X)$ there exists a translate $\kappa$ (the Riemann constant) such that
	\[
	u(D)+\kappa\in\Theta_A
	\quad\Longleftrightarrow\quad
	\ell(D)>0
	\]
	for every divisor $D$ of degree $g-1$, where $u(D)$ denotes the sum of the
	Abel--Jacobi integrals of its points, extended linearly to divisors.  Thus the
	right-hand condition says that $D$ is \emph{linearly equivalent} to an effective
	divisor.  Consequently, up to translation, $\Theta_A$ is the Abel--Jacobi image
	of $W_{g-1}$.
\end{remark}

\subsection{\texorpdfstring{Genus $1$: Jacobi theta functions}{Genus 1: Jacobi theta functions}}

Now specialize to $g=1$. Then $\Omega=\tau\in\mathbb{H}$ and $\Lambda=\Z+\tau\Z$.
Our Jacobi-theta argument is normalized so that the first period is $1$;
in the DLMF convention it is $\pi$. Thus our $\theta_j(z|\tau)$ is
the DLMF $\theta_j(\pi z|\tau)$. Keeping this factor of $\pi$ in view
is essential when comparing derivatives and the formulas for $\wp$ below.
The Riemann theta becomes
\[
\theta(z\,|\,\tau)=\sum_{n\in\Z} e^{\pi i n^2\tau + 2\pi i n z}.
\]
Let $q:=e^{\pi i\tau}$. Then
\[
\theta(z\,|\,\tau)=\sum_{n\in\Z} q^{n^2} e^{2\pi i n z}=\vartheta_3(z,\tau).
\]

\begin{remark}[Four Jacobi theta functions as characteristics]
	For $\varepsilon,\delta\in\{0,\tfrac12\}$ define
	\[
	\theta\!\begin{bmatrix}\varepsilon\\ \delta\end{bmatrix}(z\,|\,\tau)
	:=
	\sum_{n\in\Z}
	\exp\!\bigl(\pi i (n+\varepsilon)^2\tau + 2\pi i (n+\varepsilon)(z+\delta)\bigr).
	\]
	With the common convention in which the odd theta function differs by a
	constant sign from the displayed characteristic,
	\[
	\theta_3=\theta\!\begin{bmatrix}0\\0\end{bmatrix},\quad
	\theta_4=\theta\!\begin{bmatrix}0\\\tfrac12\end{bmatrix},\quad
	\theta_2=\theta\!\begin{bmatrix}\tfrac12\\0\end{bmatrix},\quad
	\theta_1=-\theta\!\begin{bmatrix}\tfrac12\\\tfrac12\end{bmatrix}.
	\]
	(The sign assigned to $\theta_1$ varies among references and has no effect on
	the formulas below, where it is absorbed into a nonzero constant.) The Jacobi
	triple product gives
	\[
	 \theta_1(z|\tau)=2q^{1/4}\sin(\pi z)
	 \prod_{m=1}^{\infty}(1-q^{2m})
	 \bigl(1-2q^{2m}\cos(2\pi z)+q^{4m}\bigr),
	 \qquad q=e^{\pi i\tau}.
	\]
	Consequently, $\theta_1$ is odd and has a simple zero at $z=0$; the
	non-sine factors are nonzero there because $|q|<1$.
\end{remark}

\begin{example}[First terms and a numerical estimate]
	\label{ex:theta-series-terms}
	Write $\vartheta_3(z,\tau)$ with $q=e^{\pi i\tau}$ and set $z=0$ for simplicity:
	\[
	\vartheta_3(0,\tau)=\sum_{n\in\Z} q^{n^2}=1+2q+2q^4+2q^9+2q^{16}+\cdots.
	\]
	For example, if $\tau=i$ then $q=e^{-\pi}\approx 0.0432139$, so
	\[
	\vartheta_3(0,i)
	=1+2e^{-\pi}+2e^{-4\pi}+2e^{-9\pi}+\cdots
	\approx 1.0864348.
	\]
	Keeping the terms through $2e^{-4\pi}$ gives an error below
	$3e^{-9\pi}<2\times10^{-12}$: indeed, for $n\ge3$ the ratio of consecutive
	terms $e^{-\pi(n+1)^2}/e^{-\pi n^2}=e^{-\pi(2n+1)}$ is at most $e^{-7\pi}$.
\end{example}

\subsection{\texorpdfstring{From theta to Weierstrass: $\sigma$, $\zeta$, and $\wp$ via logarithmic derivatives}{From theta to Weierstrass functions via logarithmic derivatives}}

Let $\Lambda=\Z\omega_1+\Z\omega_2$ with $\tau=\omega_2/\omega_1$.
The Weierstrass $\sigma$--function is characterized (up to scaling) by:
\begin{itemize}
	\item $\sigma$ is entire with simple zeros exactly at $\Lambda$,
	\item $\sigma'(0)=1$,
	\item quasi-periodicity $\sigma(z+\omega_k)= -e^{\eta_k(z+\omega_k/2)}\sigma(z)$, $k=1,2$.
\end{itemize}

\begin{theorem}[Sigma in terms of $\theta_1$]
	\label{thm:sigma-theta}
	There exists a nonzero constant $C=C(\Lambda)$ such that
	\[
	\sigma(z)
	=
	C\,
	\exp\!\left(\frac{\eta_1}{2\omega_1}z^2\right)
	\theta_1\!\left(\frac{z}{\omega_1}\,\middle|\,\tau\right).
	\]
\end{theorem}

\begin{proof}
	Define
	\[
	F(z):=
	\frac{\sigma(z)}
	{\exp\!\left(\frac{\eta_1}{2\omega_1}z^2\right)\theta_1(z/\omega_1\,|\,\tau)}.
	\]
	The transformations
	\[
	\theta_1(u+1|\tau)=-\theta_1(u|\tau),\qquad
	\theta_1(u+\tau|\tau)=-e^{-\pi i\tau-2\pi i u}\theta_1(u|\tau)
	\]
	show directly that the denominator acquires the factor
	$-e^{\eta_1(z+\omega_1/2)}$ under $z\mapsto z+\omega_1$.
	For the second period, the Legendre relation in the present full-period
	convention is
	\[
	\eta_1\omega_2-\eta_2\omega_1=2\pi i.
	\]
	It converts the denominator's multiplier
	\[
	-e^{\eta_1\omega_2z/\omega_1+\eta_1\omega_2^2/(2\omega_1)
	      -\pi i\omega_2/\omega_1-2\pi iz/\omega_1}
	\]
	into $-e^{\eta_2(z+\omega_2/2)}$, exactly the multiplier of
	\(\sigma\). Hence $F(z+\omega_1)=F(z)$ and $F(z+\omega_2)=F(z)$, so
	$F$ is elliptic. Both numerator and denominator are entire
	and have the same simple zeros at $\Lambda$, so $F$ has no zeros or poles. Thus $F$ is constant.
\end{proof}

\begin{definition}[Weierstrass $\zeta$ and $\wp$ via $\sigma$]
	\label{def:zeta-wp}
	Define
	\[
	\zeta(z):=\frac{d}{dz}\log\sigma(z),
	\qquad
	\wp(z):=-\zeta'(z)=-\frac{d^2}{dz^2}\log\sigma(z).
	\]
\end{definition}

\begin{proposition}[Pole/zero mechanism: why $\wp$ has a double pole]
	\label{prop:wp-pole-from-theta}
	Since $\theta_1$ has a simple zero at $0$, we have locally $\theta_1(z|\tau)=c\,z+O(z^3)$,
	$c\neq 0$, and hence
	\[
	\frac{d^2}{dz^2}\log\theta_1(z|\tau)= -\frac{1}{z^2}+O(1).
	\]
	So $-\frac{d^2}{dz^2}\log\theta_1$ has a double pole, matching the basic pole of $\wp$.
\end{proposition}

\begin{proof}
	If $f(z)=c z + O(z^3)$ with $c\neq 0$, then
	$\log f(z)=\log c + \log z + O(z^2)$, so $(\log f)''(z)= -\frac{1}{z^2}+O(1)$.
	Apply to $f=\theta_1$.
\end{proof}

\begin{theorem}[Explicit theta expression for $\wp$]
	\label{thm:wp-theta}
	With $\tau=\omega_2/\omega_1$ we have
	\[
	\wp(z)
	=
	-\frac{d^2}{dz^2}\log\theta_1\!\left(\frac{z}{\omega_1}\,\middle|\,\tau\right)
	-\frac{\eta_1}{\omega_1}.
	\]
	Equivalently, with $u=z/\omega_1$,
	\[
	\wp(z)
	=
	-\frac{1}{\omega_1^2}\left(
		\frac{\theta_1''(u|\tau)}{\theta_1(u|\tau)}-\Big(\frac{\theta_1'(u|\tau)}{\theta_1(u|\tau)}\Big)^2
	\right)
	-\frac{\eta_1}{\omega_1}.
	\]
\end{theorem}

\begin{proof}
	From Theorem~\ref{thm:sigma-theta},
	\[
	\log\sigma(z)=\log C+\frac{\eta_1}{2\omega_1}z^2+\log\theta_1(z/\omega_1|\tau).
	\]
	Differentiate twice and use $\wp(z)=-\frac{d^2}{dz^2}\log\sigma(z)$.
\end{proof}

\begin{example}[Recovering the principal part of $\wp$ at $0$ from theta]
	\label{ex:wp-principal-part}
	Using Proposition~\ref{prop:wp-pole-from-theta} and Theorem~\ref{thm:wp-theta}, near $z=0$,
	\[
	-\frac{d^2}{dz^2}\log\theta_1\!\left(\frac{z}{\omega_1}\,\middle|\,\tau\right)
	=
	\frac{1}{z^2}+O(1),
	\]
	so
	\[
	\wp(z)=\frac{1}{z^2}+O(1),
	\]
	which matches the defining Laurent expansion of the Weierstrass $\wp$--function.
\end{example}

\begin{corollary}[All elliptic functions from theta]
	Every elliptic function on $\C/\Lambda$ is a rational function of $\wp$ and $\wp'$.
	Since $\wp$ is expressed by $\theta_1$ and derivatives, every elliptic function can be written
	as a rational combination of Jacobi theta functions and their derivatives.
\end{corollary}

\begin{proof}
	Let $f$ be an elliptic function. Its even part
	$(f(z)+f(-z))/2$ is a rational function of $\wp(z)$: the degree-two map
	$\wp:\C/\Lambda\to\PP^1$ identifies $z$ and $-z$, so its invariant
	function field is $\C(\wp)$. The odd part divided by $\wp'(z)$ is even,
	because both numerator and denominator change sign under $z\mapsto-z$.
	It too belongs to $\C(\wp)$, and therefore
	$f=R(\wp)+\wp' S(\wp)$ for rational functions $R,S$.
	The preceding theta formula for $\wp$, and its derivative for $\wp'$, give
	the last assertion. See also \cite[Chapter~23]{DLMFWeierstrass2026}.
\end{proof}

\subsection{Exercises (with guidance)}

\paragraph{Conventions.}
Throughout, $\tau\in\mathbb{H}$ and $q=e^{\pi i\tau}$.

\begin{exercise}[Warm-up: convergence estimate]
	\label{exr:selected-theta-convergence}
	Let $\tau=X+iY$ with $Y>0$. Show that
	\[
	\sum_{n\in\Z}\left|e^{\pi i n^2\tau+2\pi i n z}\right|
	\le
	\sum_{n\in\Z} e^{-\pi Y n^2 + 2\pi |n|\,|\Im z|}.
	\]
	Deduce uniform convergence on compact sets and that $\vartheta_3(\cdot,\tau)$ is entire.
\end{exercise}

\begin{exercise}[Quasi-periodicity by reindexing]
	Prove directly from the defining series that
	\[
	\vartheta_3(z+1,\tau)=\vartheta_3(z,\tau),\qquad
	\vartheta_3(z+\tau,\tau)=e^{-\pi i\tau-2\pi i z}\vartheta_3(z,\tau).
	\]
	(Hint: for the $\tau$-shift reindex $n\mapsto n-1$.)
\end{exercise}

\begin{exercise}[Parity and zero of $\theta_1$]
	Using the characteristic series and the substitution $n\mapsto-n-1$, show
	$\theta_1(-z|\tau)=-\theta_1(z|\tau)$. Then use the triple product above to
	prove that the zero at $z=0$ is simple. Write
	$\theta_1(z|\tau)=c z+O(z^3)$ and show both that
	$c=\theta_1'(0|\tau)$ and that
	\[
	 c=2\pi q^{1/4}\prod_{m=1}^{\infty}(1-q^{2m})^3\ne0.
	\]
\end{exercise}

\begin{exercise}[Principal part of $\wp$ from a simple zero]
	Let $f$ be holomorphic near $0$ with $f(0)=0$ and $f'(0)\neq 0$.
	Show that $(\log f)''(z)=-1/z^2+O(1)$ near $0$.
	Apply to $f(z)=\theta_1(z|\tau)$ to re-derive the double pole of $\wp$.
\end{exercise}

\begin{exercise}[A numerical approximation with an error bound]
	Take $\tau=i$. Compute $q=e^{-\pi}$ and approximate
	\[
	\vartheta_3(0,i)=1+2q+2q^4+2q^9+\cdots
	\]
	using terms up to $n=2$ and $n=3$. Give an error bound using the tail
	$\sum_{n\ge N} 2e^{-\pi n^2}$.
\end{exercise}

\begin{exercise}[From $\sigma$ to $\wp$]
	Assume Theorem~\ref{thm:sigma-theta}.
	Differentiate $\log\sigma(z)$ once and twice to express $\zeta(z)$ and $\wp(z)$ explicitly
	in terms of $\theta_1(z/\omega_1|\tau)$ and its derivatives.
	Check the constant term $\eta_1/\omega_1$ appears with the correct sign.
\end{exercise}

\begin{exercise}[Line bundle cocycle]
	Let $\alpha_{m+\Omega n}(z)=\exp(-\pi i n^t\Omega n-2\pi i n^t z)$.
	Verify the cocycle identity
	\[
	\alpha_{\lambda+\mu}(z)=\alpha_\lambda(z+\mu)\alpha_\mu(z)
	\]
	by expanding both sides, and isolate exactly where symmetry of $\Omega$ is used.
\end{exercise}

\subsection{Modular transformations and the geometric role of theta}

Theta functions satisfy explicit transformation laws under the generators
\[
\tau\mapsto \tau+1,\qquad \tau\mapsto -\frac{1}{\tau},
\]
and these laws encode how line bundles and divisors vary in families of elliptic curves.
From the Jacobian viewpoint, theta is a canonical section whose zero divisor
represents the principal polarization. In genus $1$, its logarithmic derivatives
give explicit formulas for $\sigma$, $\zeta$, and $\wp$.

\subsection{From Euler's product to theta and the Riemann hypothesis}
\label{subsec:theta-mellin-zeta}

\paragraph{Why begin with primes?}
For a complex number $s=\sigma+it$ with $\sigma>1$, the
\emph{Dirichlet series}
\[
 \zeta(s)=\sum_{n=1}^{\infty}\frac1{n^s},\qquad
 n^{-s}:=e^{-s\log n},
\]
converges absolutely, since $|n^{-s}|=n^{-\sigma}$.
Euler's observation is that unique factorization turns this sum over
all positive integers into a product over primes:
\[
 \zeta(s)=\prod_{p\ \mathrm{prime}}(1-p^{-s})^{-1},
 \qquad \sigma>1.
\]
To see the mechanism, multiply the geometric series
$1+p^{-s}+p^{-2s}+\cdots$ for finitely many primes.  Each resulting
monomial corresponds to exactly one integer whose prime factors lie
in that finite set.  Increasing the set to include every prime gives
the Dirichlet series; absolute convergence justifies the limit.
Each finite factor is nonzero. Moreover,
$\sum_p|p^{-s}|\le\sum_{n\ge2}n^{-\sigma}<\infty$, so the
infinite product converges to a nonzero value. This explains why
$\zeta$ has no
zeros in its initial region of definition.  The problem of zeros only
becomes interesting after continuation beyond that region
\cite[\S25.2]{DLMFZeta2026}.

\paragraph{Why a torus theta section enters.}
Earlier we constructed $\vartheta_3(z,\tau)$ as a holomorphic section
of a degree-one theta line bundle on the elliptic curve
$E_\tau=\C/(\Z+\tau\Z)$.  It is \emph{not} an ordinary function on
the quotient: its quasi-periodicity is the line-bundle transition law.
For the one-parameter family $\tau=it$ with $t>0$, evaluate that
section at the distinguished point $z=0$ and write
\[
 \Theta(t)=\vartheta_3(0,it)=\sum_{n\in\Z}e^{-\pi n^2t},\qquad t>0.
\]
Thus $\Theta(t)$ is a value of a theta section as the torus varies;
it is not a function defined \emph{on one fixed torus}.  Its summation
index records the integer lattice.  Poisson summation for the Gaussian
$e^{-\pi t x^2}$, whose Fourier transform is
$t^{-1/2}e^{-\pi\xi^2/t}$, gives
\[
 \Theta(t)=t^{-1/2}\Theta(1/t).
\]
Geometrically the change $t\leftrightarrow1/t$ exchanges the lattice
and its Fourier-dual lattice; analytically it is the modular
transformation behind the functional equation
\cite[\S20.7(viii)]{DLMFTheta2026}.

\paragraph{Mellin transform: recovering the Dirichlet series.}
The Mellin transform integrates against $t^{u-1}dt$; it plays for
dilations the role that the Fourier transform plays for translations.
The theta series is useful because
$\Theta(t)-1=2\sum_{n\ge1}e^{-\pi n^2t}$.
For $\sigma>1$, termwise integration and the substitution
$v=\pi n^2t$ give
\[
 \begin{aligned}
 \frac12\int_0^\infty(\Theta(t)-1)t^{s/2-1}\,dt
 &=\sum_{n\ge1}\int_0^\infty e^{-\pi n^2t}t^{s/2-1}\,dt\\
 &=\Gamma(s/2)\pi^{-s/2}\sum_{n\ge1}n^{-s}.
 \end{aligned}
\]
The \emph{completed zeta function} in this chapter is therefore
$\Lambda(s):=\pi^{-s/2}\Gamma(s/2)\zeta(s)$, initially for
$\sigma>1$.  Completion means inserting the Gamma and power-of-$\pi$
factors produced naturally by the Gaussian integral.  The domain
$\sigma>1$ matters: because
$\Theta(t)=t^{-1/2}\Theta(1/t)$, the Mellin integrand near zero has
leading order $\tfrac12t^{(s-1)/2-1}$; the integral is not yet defined
by convergence in the critical strip.

Mellin inversion also runs backward.  On any line $\Re s=c>1$,
\[
 \Theta(t)-1=\frac{1}{2\pi i}
 \int_{c-i\infty}^{c+i\infty}\Lambda(s)t^{-s/2}\,ds.
\]
The contour is the vertical line from $c-i\infty$ to $c+i\infty$.
To check the coefficient, apply ordinary Mellin inversion to
$g(t)=(\Theta(t)-1)/2$ with Mellin variable $u=s/2$; changing
variables contributes $ds=2\,du$.  The exponential decay of the
Gamma factor on vertical lines makes this integral convergent.
Inversion says that the theta data and the completed Dirichlet-series
data determine one another; it does not itself locate any zeros.

\paragraph{How the symmetry appears.}
Split the Mellin integral at $t=1$.  In its part from $0$ to $1$,
substitute $u=1/t$ and use
$\Theta(t)=t^{-1/2}\Theta(1/t)$.  The two elementary terms from
$t^{-1/2}-1$ integrate to
$1/(s-1)-1/s=1/[s(s-1)]$.  The result is
\[
 \Lambda(s)=\frac{1}{s(s-1)}+
 \frac12\int_1^\infty(\Theta(t)-1)
 \bigl(t^{s/2-1}+t^{(1-s)/2-1}\bigr)\,dt.
\]
The integral from $1$ to infinity converges locally uniformly for
every complex $s$, because $\Theta(t)-1$ decays exponentially.
Hence this formula continues $\Lambda$ meromorphically, with poles
at $0$ and $1$, and is unchanged on replacing $s$ by $1-s$:
$\Lambda(s)=\Lambda(1-s)$.  The poles disappear in the
\emph{entire} Riemann xi-function
\[
 \xi(s)=\tfrac12s(s-1)\Lambda(s),\qquad \xi(s)=\xi(1-s).
\]
The factor $s(s-1)$ removes the poles of $\Lambda$ at $0$ and $1$;
at negative even integers the trivial zeros of $\zeta$ cancel the
poles of $\Gamma(s/2)$. Thus the zeros of $\xi$ are exactly the
nontrivial zeros of $\zeta$. The reflection
$s\mapsto1-s$ pairs them across the vertical line $\Re s=1/2$;
complex conjugation supplies a second symmetry across the real axis.

\paragraph{What the Riemann hypothesis actually asks.}
The \emph{critical strip} is $0<\Re s<1$, and the \emph{critical line}
is $\Re s=1/2$.  The Riemann hypothesis says that every nontrivial
zero lies on that line.  Euler's product excludes zeros for
$\Re s>1$, while the functional equation mirrors that information
and explains the trivial zeros at the negative even integers.
Excluding zeros on the boundary $\Re s=1$ requires a separate
zero-free theorem, not proved by the argument here.
The additional assertion that all remaining zeros lie on the
\emph{middle} line is not a consequence of reflection alone:
the elementary even polynomial $w^2-a^2$ is invariant under
$w\mapsto-w$ but has zeros away from $w=0$.
Likewise $\xi(s)=\xi(1-s)$ supplies a symmetry, not a proof of the
Riemann hypothesis.  Proving or disproving that assertion remains
open.  For the classical normalizations and zero terminology, see
\cite[Chapter~12]{Apostol1990} and
\cite[\S\S25.2,25.4,25.10]{DLMFZeta2026}.

\begin{exercise}[Trace the factors in the theta--Mellin identity]
 Use $v=\pi n^2t$ in the Gamma integral to derive the factor
 $\pi^{-s/2}n^{-s}$ in the Mellin calculation. Then verify directly
 that $1/(s-1)-1/s=1/[s(s-1)]$.
\end{exercise}

\section{Galois Theory and Riemann Surfaces}
\index{Galois theory}

\providecommand{\C}{\mathbb{C}}
\providecommand{\R}{\mathbb{R}}
\providecommand{\Z}{\mathbb{Z}}
\providecommand{\Q}{\mathbb{Q}}
\providecommand{\D}{\mathbb{D}}
\providecommand{\PP}{\mathbb{P}}
\providecommand{\A}{\mathbb{A}}
\providecommand{\M}{\mathcal{M}}
\providecommand{\OO}{\mathcal{O}}
\providecommand{\Aut}{\operatorname{Aut}}
\providecommand{\Deck}{\operatorname{Deck}}
\providecommand{\Gal}{\operatorname{Gal}}
\providecommand{\Spec}{\operatorname{Spec}}
\providecommand{\im}{\operatorname{Im}}
\providecommand{\PSL}{\operatorname{PSL}}
\providecommand{\PGL}{\operatorname{PGL}}
\providecommand{\ord}{\operatorname{ord}}
\providecommand{\Tr}{\operatorname{Tr}}
\providecommand{\Nm}{\operatorname{Nm}}
\providecommand{\divv}{\operatorname{div}}
\providecommand{\Res}{\operatorname{Res}}
\providecommand{\id}{\mathrm{id}}

A finite holomorphic map $\pi:Y\to X$ between compact Riemann surfaces is the analytic avatar of
a finite extension of function fields $\M(Y)/\M(X)$.  The guiding dictionary is:
\[
\text{finite map }\pi:Y\to X
\quad\Longleftrightarrow\quad
\text{finite field extension }\M(X)\hookrightarrow \M(Y),
\]
and in the Galois (normal) case,
\[
\Deck(Y/X)\ \cong\ \Aut\bigl(\M(Y)/\M(X)\bigr).
\]
We develop the correspondence concretely: we write polynomial relations, compute
degrees, and identify field automorphisms with deck transformations.
The first task is to pass from a map to an extension of meromorphic-function
fields; the converse and the deck-group comparison require additional arguments.
The double and cyclic covers later in the section provide tests for both directions.
The analytic-function-field dictionary is compared with \cite{Forster1981,Stichtenoth2009,FriedJarden2008}.

\paragraph{Prerequisite check: a degree-two field extension.}
For the map $t\mapsto x=t^2$ on $\PP^1$, the induced inclusion is
$\C(x)=\C(t^2)\subset\C(t)$. Every rational function of $t$ can
be written $a(t^2)+t b(t^2)$, so the extension has basis $1,t$ and
degree two. The involution $t\mapsto-t$ fixes $x$ and accounts for
both sheets away from $0$ and $\infty$. This example separates three
notions used below: field degree, deck symmetry, and ramification.

\medskip
\noindent
\textbf{Standing conventions.}
All Riemann surfaces are connected and compact unless explicitly stated.
A \emph{finite holomorphic map} means a nonconstant holomorphic map with finite fibers; equivalently,
a branched covering.  Let $B\subset X$ be the finite branch value set and set
\[
X^\circ:=X\setminus B,\qquad Y^\circ:=Y\setminus \pi^{-1}(B).
\]
Then $\pi^\circ:Y^\circ\to X^\circ$ is an unramified covering of degree $\deg\pi$.

\subsection{Function fields and pullback}
\index{function fields!pullback}

For a Riemann surface $X$, let $\M(X)$ denote its field of meromorphic functions.
A nonconstant holomorphic map $\pi:Y\to X$ induces an injective field homomorphism
\[
\pi^{*}:\ \M(X)\hookrightarrow \M(Y),\qquad f\longmapsto f\circ\pi.
\]
We identify $\M(X)$ with its image $\pi^*\M(X)\subset \M(Y)$.

\subsubsection{Local branches, norm, and trace}

Let $f\in \M(Y)$. Over an evenly covered open set $U\subset X^\circ$ we can write
\[
(\pi^\circ)^{-1}(U)=V_1\sqcup\cdots\sqcup V_d,\qquad \pi|_{V_j}:V_j\stackrel{\sim}{\longrightarrow} U,
\]
and define local inverses $\tau_j:U\to V_j$ so that $\pi\circ\tau_j=\id_U$.
Define the \emph{local branches} of $f$ above $U$ by
\[
f_j:=f\circ \tau_j\ \in\ \M(U).
\]
Then the elementary symmetric polynomials in $f_1,\dots,f_d$ define meromorphic functions on $U$:
\[
\prod_{j=1}^d (T-f_j)
= T^d + c_1 T^{d-1} + \cdots + c_d\in \M(U)[T].
\]

\begin{lemma}[Gluing of symmetric coefficients]
	\label{lem:galois-glue}
	The coefficients $c_i$ are independent of the choice of labeling of the sheets $V_j$,
	and they glue (as $U$ varies) to define global meromorphic functions $c_i\in \M(X)$.
\end{lemma}

\begin{proof}
	On overlaps $U\cap U'$, the decompositions into sheets differ by a permutation of the sheets,
	so the unordered set $\{f_1,\dots,f_d\}$ is permuted.  Each $c_i$ is an elementary symmetric
	polynomial in the $f_j$, hence invariant under permutation, so the $c_i$ agree on overlaps.
	Thus they glue on $X^\circ$. Near a branch value choose local parameters in which the finite map has
	the form $z=u^e$. Each branch of $f$ has at most polynomial growth in $u^{-1}$, because $f$ is
	meromorphic on $Y$; its elementary symmetric functions therefore have at most polynomial growth in
	$z^{-1}$. The isolated singularities at the finite set $B$ are consequently poles or removable, so the
	coefficients extend meromorphically to $X$.
\end{proof}

\begin{definition}[Norm and trace along $\pi$ (computable on sheets)]
	For $f\in \M(Y)$ define on $U\subset X^\circ$:
	\[
	\Nm_{\pi}(f)\big|_U := \prod_{j=1}^d f_j \in \M(U),
	\qquad
	\Tr_{\pi}(f)\big|_U := \sum_{j=1}^d f_j \in \M(U).
	\]
	These glue to meromorphic functions on $X$ and satisfy
	$\Nm_{\pi}(f)=(-1)^d c_d$ and $\Tr_{\pi}(f)=-c_1$.
\end{definition}

\begin{example}[Trace and norm on a quadratic cover]
	Let $Y$ be the smooth compactification of $y^2=q(x)$ with $q$ a
	nonconstant squarefree polynomial,
	and let $\pi:Y\to\mathbb P^1$ be the $x$-projection. Away from its branch
	values, the two inverse branches send $x$ to $(x,\sqrt{q(x)})$ and
	$(x,-\sqrt{q(x)})$. For the meromorphic function $y$ on $Y$,
	\[
	\Tr_\pi(y)=\sqrt{q(x)}-\sqrt{q(x)}=0,
	\qquad
	\Nm_\pi(y)=\sqrt{q(x)}\,(-\sqrt{q(x)})=-q(x).
	\]
	The associated monic relation is $T^2-q(x)$, and its two roots coincide
	precisely over the finite branch values. The trace and norm nonetheless
	extend meromorphically across those values because they are symmetric in
	the two local branches. This example also fixes the sign convention:
	for a degree-two monic polynomial the norm is its constant coefficient.
\end{example}

\begin{theorem}[A monic algebraic relation of degree $\le d$]
	\label{thm:galois-minpoly}
	Let $\pi:Y\to X$ be finite of degree $d$ and let $f\in \M(Y)$.  Let $c_1,\dots,c_d\in \M(X)$
	be the global symmetric coefficients attached to $f$.
	Then in $\M(Y)$ we have the identity
	\[
	f^d + (\pi^*c_1) f^{d-1} + \cdots + (\pi^*c_{d-1}) f + \pi^*c_d = 0.
	\]
	In particular, every $f\in\M(Y)$ is algebraic over $\M(X)$ with
	$\deg_{\M(X)}(f)\le d$.
\end{theorem}

\begin{proof}
	Work over an evenly covered $U\subset X^\circ$. Consider the polynomial
	\[
	P_U(T)=\prod_{j=1}^d (T-f_j)=T^d+c_1T^{d-1}+\cdots+c_d\in \M(U)[T].
	\]
	Pull it back to $\pi^{-1}(U)$ via $\pi$; on the sheet $V_k$ we have $\tau_k\circ\pi=\id_{V_k}$,
	hence $(f\circ\tau_k)\circ\pi=f$ on $V_k$. Therefore, on $V_k$,
	\[
	P_U(f)=\prod_{j=1}^d \bigl(f-(f\circ\tau_j)\circ\pi\bigr)
	\]
	has the $j=k$ factor equal to $0$, so $P_U(f)=0$ on $V_k$. Since $\pi^{-1}(U)$ is the disjoint union
	of the sheets, $P_U(f)=0$ on $\pi^{-1}(U)$. By analytic continuation, the identity holds on $Y^\circ$,
	and hence on all of $Y$ (as an identity of meromorphic functions). Replacing $c_i$ by their global
	meromorphic extensions on $X$ yields the global formula.
\end{proof}

\subsubsection{\texorpdfstring{Why the extension degree equals $\deg\pi$}{Why the extension degree equals the degree of the map}}

\begin{theorem}[Degree equals field extension degree]
	\label{thm:galois-degree-equals}
	Let $\pi:Y\to X$ be finite of degree $d$. Then
	\[
	[\M(Y):\M(X)] = d.
	\]
\end{theorem}

\begin{proof}
	Set $K=\M(X)$ and $L=\M(Y)$. Regard the compact Riemann surfaces as smooth
	projective curves. On the complement of the branch values, $\pi$ is finite
	locally free of rank $d$: over an evenly covered open set, a section of
	$\pi_*\mathcal O_Y$ is a $d$-tuple, one function on each sheet. Passing to the
	generic point turns $\pi_*\mathcal O_Y$ into the $K$-vector space $L$. Rank is
	preserved at the generic point, and therefore
	\[
	\dim_KL=d.
	\]
	This proves the equality. The following construction gives a concrete element
	of $L$ whose minimal polynomial already has degree $d$; it also makes the
	generic-rank calculation visible through local branches.
	
	Choose $x\in X^\circ$ (not a branch value) so $\pi^{-1}(x)=\{y_1,\dots,y_d\}$ consists of
	$d$ distinct points. We claim there exists $f\in L$ such that the values $f(y_i)$ are pairwise
	distinct. This is a standard separation statement: since $Y$ is a compact Riemann surface,
	meromorphic functions separate points.
	
	A concrete construction uses Riemann--Roch as follows. Choose an effective divisor $D$ of degree
	$>2g(Y)-2+d$ whose support is disjoint from $\{y_1,\ldots,y_d\}$. The exact sequence for evaluation,
	together with $H^1(Y,\mathcal O_Y(D-y_1-\cdots-y_d))=0$, shows that the evaluation map
	\[
	\mathrm{ev}:\ H^0(Y,\OO_Y(D))\longrightarrow \C^d,\qquad
	s\longmapsto (s(y_1),\dots,s(y_d))
	\]
	is surjective. Choose distinct complex numbers $a_1,\ldots,a_d$ and a section $f$ mapping to
	$(a_1,\ldots,a_d)$. Viewed in the defining meromorphic-function model of $\mathcal O_Y(D)$, this
	is a meromorphic function on $Y$ satisfying $f(y_i)=a_i$.
	
	Now work on a simply connected neighborhood $U\subset X^\circ$ of $x$ so that $\pi^{-1}(U)$ splits into
	sheets $V_i$ with local inverses $\tau_i:U\to V_i$ satisfying $\tau_i(x)=y_i$.
	Define the $d$ local branches $f_i:=f\circ\tau_i\in \M(U)$. By construction, $f_i(x)=f(y_i)$ are pairwise
	distinct, hence the functions $f_i$ are distinct as germs at $x$.
	
	Let $P(T)\in K[T]$ be the monic polynomial of degree $d$ obtained from the elementary symmetric polynomials
	of $f_1,\dots,f_d$ (as in Theorem~\ref{thm:galois-minpoly}). Over $U$, this polynomial factors as
	\[
	P(T)=\prod_{i=1}^d (T-f_i)\in \M(U)[T].
	\]
	Since the $f_i$ are distinct germs at $x$, this factorization has $d$ \emph{distinct} roots in a neighborhood
	of $x$. Therefore the degree of the minimal polynomial of $f$ over $K$ must be at least $d$.
	Indeed, a relation $Q(f)=0$ over $K$ pulls back to $Q(f_i)=0$ for every
	branch $f_i$; over the field of meromorphic germs at $x$, a nonzero
	polynomial of degree less than $d$ cannot have $d$ distinct roots.
	The monic relation from
	Theorem~\ref{thm:galois-minpoly} has degree $d$, so the minimal polynomial has
	degree exactly $d$ and $L=K(f)$.
\end{proof}

\subsection{Galois coverings and Galois extensions}
\index{Galois covers}\index{Galois extensions!function fields}

\subsubsection{Deck transformations (unramified and branched)}

\begin{definition}[Deck group]
	For an unramified covering $p:Y^\circ\to X^\circ$, the \emph{deck transformation group} is
	\[
	\Deck(Y^\circ/X^\circ):=\{\sigma\in \Aut_{\mathrm{hol}}(Y^\circ)\mid p\circ \sigma = p\}.
	\]
	For a finite map $\pi:Y\to X$ with branch values $B$, define
	\[
	\Deck(Y/X):=\Deck\bigl(Y^\circ/X^\circ\bigr).
	\]
\end{definition}

\begin{definition}[Galois/normal covering]
	An unramified covering $p:Y^\circ\to X^\circ$ is \emph{Galois} if $\Deck(Y^\circ/X^\circ)$ acts transitively on
	each fiber. Equivalently, $p_*\pi_1(Y^\circ)$ is a normal subgroup of $\pi_1(X^\circ)$.
	A branched map $\pi:Y\to X$ is called \emph{Galois} if $\pi^\circ:Y^\circ\to X^\circ$ is a Galois covering.
\end{definition}

\subsubsection{Main theorem: deck transformations = field automorphisms}

\begin{theorem}[Deck transformations vs.\ field automorphisms]
	\label{thm:galois-deck-aut}
	Let $\pi:Y\to X$ be finite. Set $K=\M(X)$ and $L=\M(Y)$ (with $K\subset L$ via $\pi^*$).
	Then composition induces a canonical group isomorphism
	\[
	\Phi:\Deck(Y/X)\ \xrightarrow{\ \sim\ }\ \Aut(L/K),
	\qquad
	\Phi(\sigma)(f):=f\circ \sigma^{-1}.
	\]
	Consequently, $\pi$ is Galois if and only if the field extension $L/K$ is Galois.
\end{theorem}

\begin{proof}
	\textbf{(1) $\Phi$ is well-defined and a homomorphism.}
	A deck transformation of $Y^\circ\to X^\circ$ extends uniquely across the finitely
	many omitted points to an automorphism of $Y$. Indeed, in local coordinates a
	finite map has the form $u\mapsto u^e$, and a lift over a punctured disk is
	bounded in the target coordinate and therefore extends removably. Consequently
	$f\circ\sigma^{-1}$ is meromorphic on all of $Y$.
	For $g\in K$ we have
	\[
	(g\circ \pi)\circ\sigma^{-1}=g\circ (\pi\circ\sigma^{-1}) = g\circ \pi,
	\]
	so $\Phi(\sigma)$ fixes $K$. The rule $\sigma\mapsto(f\mapsto f\circ\sigma^{-1})$ respects composition, hence is a homomorphism.
	
	\textbf{(2) $\Phi$ is injective.}
	Suppose $\Phi(\sigma)=\id_L$. Then $f\circ\sigma^{-1}=f$ for all $f\in \M(Y)$.
	Meromorphic functions separate points on a compact Riemann surface: for $p\neq q$ there exists
	$f\in \M(Y)$ with $f(p)\neq f(q)$ (e.g.\ by Riemann--Roch). Hence $\sigma^{-1}(p)=p$ for all $p$, so $\sigma=\id$.
	
	\textbf{(3) $\Phi$ is surjective (construct $\sigma$ from a field automorphism).}
	Let $\varphi\in \Aut(L/K)$. Choose the element $f$ constructed in the proof of
	Theorem~\ref{thm:galois-degree-equals}. Its minimal polynomial has degree $d=[L:K]$,
	so $L=K(f)$. Choose $x\in X^\circ$ and a simply connected neighborhood $U\subset X^\circ$ of $x$
	so that $\pi^{-1}(U)=V_1\sqcup\cdots\sqcup V_d$ splits into sheets with inverses $\tau_i:U\to V_i$.
	After shrinking $U$ away from the discrete collision set if necessary, the local
	branches $f_i:=f\circ\tau_i$ are pairwise distinct on $U$.
	
	Over $U$, the monic polynomial
	\[
	P(T):=\prod_{i=1}^d (T-f_i)\ \in\ \M(U)[T]
	\]
	has coefficients in $K$ (by the symmetric-gluing construction). Since $\varphi$ fixes $K$, it preserves $P(T)$,
	hence sends the root $f$ to another root of the same polynomial. Concretely, on each sheet $V_i$,
	$\varphi(f)$ must equal $f\circ(\tau_{\rho(i)}\circ \pi)$ for a \emph{unique}
	index $\rho(i)$. Applying the same construction to $\varphi^{-1}$ shows that
	$\rho$ is a permutation. Define the local deck transformation by
	\[
	\sigma_U^{-1}|_{V_i}:V_i\longrightarrow V_{\rho(i)},
	\qquad
	\sigma_U^{-1}|_{V_i}:=\tau_{\rho(i)}\circ \pi|_{V_i}.
	\]
	By construction, $\pi\circ \sigma_U=\pi$ and
	$\varphi(f)=f\circ\sigma_U^{-1}$ on every sheet.
	
	The branches of $f$ may coincide at isolated points of $X^\circ$, but this
	does not obstruct extension. On the complement of these points, repeat the
	construction in evenly covered disks. On overlaps the local maps agree:
	their pullbacks agree on $K$ and on $f$, hence on $L=K(f)$, and meromorphic
	functions separate points of $Y$. Around a collision point, each sheet of
	the covering is still a disk mapping biholomorphically to the base disk.
	The sheet permutation is constant on the punctured disk, so its formula
	$\tau_{\rho(i)}\circ\pi$ extends across the center. The local maps therefore
	glue to a global $\sigma\in\Deck(Y/X)$, and $\Phi(\sigma)=\varphi$.
	
	\textbf{(4) Galois equivalence.}
	The cover $\pi$ is Galois iff $|\Deck(Y/X)|=\deg\pi$ (transitive free action on fibers).
	By Theorem~\ref{thm:galois-degree-equals}, $\deg\pi=[L:K]$.
	Thus $\pi$ is Galois iff $|\Aut(L/K)|=[L:K]$, i.e.\ iff $L/K$ is a Galois extension.
\end{proof}

\subsection{Examples}

\subsubsection{Example A: a hyperelliptic double cover}

\begin{example}[Elliptic curve as a double cover of $\PP^1$]
	\label{ex:galois-hyperelliptic}
	Let
	\[
	E:\ y^2 = x^3 - x
	\]
	be the smooth projective model of this affine curve, and let $\pi:E\to\PP^1$ be $(x,y)\mapsto x$.
	Then $\M(\PP^1)=\C(x)$ and
	\[
	\M(E)=\C(x,y)\cong \C(x)[y]/(y^2-(x^3-x)).
	\]
	\textbf{(1) Minimal polynomial.}  The element $y\in \M(E)$ satisfies
	\[
	T^2-(x^3-x)=0,
	\]
	and this polynomial is irreducible over $\C(x)$: the rational function
	$x^3-x$ has simple zeros at $0,1,-1$, whereas a square in $\C(x)$ has
	even order at every point. Thus $\deg_{\C(x)}(y)=2$ and
	$[\M(E):\C(x)]=2=\deg\pi$.
	
	\textbf{(2) Deck group / Galois group.}  The nontrivial automorphism of $\M(E)/\C(x)$ is
	\[
	\sigma^\sharp:\ \M(E)\to \M(E),\qquad x\mapsto x,\quad y\mapsto -y.
	\]
	On the curve it corresponds to the involution
	\[
	\sigma:E\to E,\qquad (x,y)\mapsto(x,-y).
	\]
	Since $\pi$ records only the $x$-coordinate, $\pi\circ \sigma=\pi$. Hence
	\[
	\Deck(E/\PP^1)\cong \Aut(\M(E)/\C(x))\cong \Z/2\Z.
	\]
	
	\textbf{(3) Branch points.} Branch occurs where $d\pi=0$ on $E$.
	On the affine chart, $\pi=x$, and the ramification points occur where $\partial (y^2-(x^3-x))/\partial y=2y=0$,
	i.e.\ where $y=0$. Thus $x^3-x=0$, so $x\in\{0,1,-1\}$.
	At the unique point at infinity use $t=x/y$ as a local parameter. The
	equation gives $t^2x^2-x-t^2=0$, whence the branch with $x\to\infty$
	has $x=t^{-2}+O(t^2)$. The $x$-projection therefore also has local degree
	$2$ at infinity.
	So $\pi$ is branched over $\{0,1,-1,\infty\}$ (four branch values), as expected for a degree-$2$ map from genus $1$
	by Riemann--Hurwitz.
\end{example}

\subsubsection{\texorpdfstring{Example B: cyclic $n$-covers (Kummer) with full Galois group}{Example B: cyclic n-covers with full Galois group}}

\begin{example}[Cyclic $n$-cover of $\PP^1$ and its Galois group]
	\label{ex:galois-kummer}
	Fix distinct $a_1,\dots,a_r\in\C$ and integers $m_1,\dots,m_r$.
	In $K=\C(x)$, consider the Kummer extension
	\[
	L=K(y),\qquad y^n=\prod_{j=1}^r (x-a_j)^{m_j}.
	\]
	Let $Y$ be the smooth projective curve with function field $L$ and let
	$\pi:Y\to\PP^1$ be the finite map induced by $K\hookrightarrow L$.
	This function-field formulation also covers negative exponents $m_j$; after
	clearing denominators one may recover affine equations on suitable charts.
	
	\textbf{(1) Degree.} Assume that $\prod_j(x-a_j)^{m_j}$ is not a $p$th power
	in $K$ for any prime $p$ dividing $n$. Then the Kummer polynomial
	$T^n-\prod_j(x-a_j)^{m_j}\in K[T]$ is irreducible,
	so $\deg_K(y)=n$ and $[L:K]=n=\deg\pi$.
	
	\textbf{(2) Galois automorphisms.}  Let $\zeta_n=e^{2\pi i/n}$.  Then
	\[
	\sigma_k^\sharp:\ L\to L,\qquad x\mapsto x,\quad y\mapsto \zeta_n^k\,y
	\]
	defines a $K$-automorphism for each $k\in\{0,\dots,n-1\}$, and these are distinct.
	Thus
	\[
	\Aut(L/K)\ \cong\ \Z/n\Z,
	\]
	and by Theorem~\ref{thm:galois-deck-aut} the deck group is also $\Z/n\Z$.
	
	\textbf{(3) Ramification indices (explicit).}
	Near $x=a_j$, write $u=x-a_j$ as a local parameter on $\PP^1$.
	Then locally the equation is $y^n=u^{m_j}\cdot(\text{unit})$.
	Because the residue field is $\mathbb C$, a holomorphic unit on a sufficiently
	small disk has an $n$th root. Absorbing this root and the integral power
	$u^{q_j}$ into $y$, we reduce to $Y^n=u^{r_j}$, where
	$m_j=q_jn+r_j$ and $0\le r_j<n$. Put
	$d_j=\gcd(n,r_j)$, $n=d_jn'$, and $r_j=d_jr'$ with
	$\gcd(n',r')=1$. Then
	\[
	Y^n-u^{r_j}=\prod_{\zeta^{d_j}=1}
	\bigl(Y^{n'}-\zeta u^{r'}\bigr).
	\]
	The $d_j$ factors give the branches in the normalization. On each branch
	one may choose a parameter $t$ for which, up to nonzero constants,
	$u=t^{n'}$ and $Y=t^{r'}$. Thus the map to the $u$-disk has local degree
	$n'=n/d_j$. Consequently, the ramification index at points above $a_j$ is
	\[
	e_{a_j}=\frac{n}{\gcd(n,m_j)}=\frac{n}{\gcd(n,r_j)}.
	\]
	This factorization proves at the same time that there are
	$\gcd(n,m_j)$ points above $a_j$.

	Infinity must be included as well. Put $m_\infty=-\sum_jm_j$, the order of
	the right-hand side at infinity. The same calculation gives
	$e_\infty=n/\gcd(n,m_\infty)$. Hence
	Riemann--Hurwitz takes the explicit form
	\[
	 2g(Y)-2=-2n+\sum_{v\in\{a_1,\ldots,a_r,\infty\}}
	 \bigl(n-\gcd(n,m_v)\bigr),
	\]
	where terms with $m_v=0$ vanish.
\end{example}

\subsection{Uniformization and the Galois picture}

The preceding examples described finite covers by equations and
function-field extensions.  Uniformization supplies a complementary
description: start with a simply connected surface and recover the
curve as its quotient by a deck group.  An unramified cover then
corresponds to a subgroup of this group; normal subgroups give regular
covers whose deck group is the quotient.  This parallels the
normal-subgroup criterion for Galois function-field extensions, but
branch points must be removed before using the ordinary covering-space
correspondence.  The theorem below names the simply connected analytic
models that organize this picture.

\subsubsection{The uniformization theorem and its analytic core}
\index{uniformization theorem}

The Riemann mapping theorem classifies simply connected proper domains in
$\C$, but it does not immediately apply to an abstract Riemann surface: one
does not yet have a global coordinate placing the surface inside the plane.
Uniformization constructs exactly such a coordinate.  Its proof combines the
topological universal-covering theorem with the following potential-theoretic
result.

\begin{lemma}[Uniformizing-function lemma]
\label{lem:uniformizing-function}
Let $S$ be a connected, simply connected, noncompact Riemann surface. Exactly
one of the following alternatives occurs.
\begin{enumerate}
	\item $S$ admits a positive Green function, and $S$ is biholomorphic to
	$\D$.
	\item $S$ admits no positive Green function, and $S$ is biholomorphic to
	$\C$.
\end{enumerate}
\end{lemma}

\begin{proof}[Analytic construction]
Choose smoothly bounded relatively compact domains
$\Omega_1\Subset\Omega_2\Subset\cdots$ exhausting $S$, with
$p\in\Omega_1$. Let $G_j(\,\cdot\,,p)$ be the Green function of $\Omega_j$
with zero boundary values. The maximum principle gives
$0<G_j\leq G_{j+1}$ on $\Omega_j\setminus\{p\}$. Harnack's convergence
theorem yields a dichotomy: either $G_j$ converges locally to a Green
function $G$ on $S$, or $G_j(q,p)\to+\infty$ at one, and hence every,
$q\ne p$.

In the first case choose a local coordinate $z$ centered at $p$. With the
normalization used here,
\[
	G(z,p)=-\log|z|+h(z),
\]
where $h$ is harmonic near $p$. The form $\alpha=-2\partial G$ is
holomorphic on $S\setminus\{p\}$ and has residue $1$ at $p$. Since $S$ is
simply connected, the periods of $\alpha$ on the punctured surface are
integral multiples of $2\pi i$. Consequently
\[
	F(x)=\exp\!\left(\int_{x_0}^{x}\alpha\right)
\]
is path independent. It extends across $p$ with a simple zero and satisfies
$|F|=e^{-G}<1$. On a regular level curve $G=c$, the argument principle gives
one zero of $F-w$ whenever $|w|<e^{-c}$. Letting $c\downarrow0$ proves that
every $w\in\D$ has exactly one preimage. Hence $F:S\to\D$ is a
biholomorphism.

In the second case fix $q\ne p$ and normalize before taking a limit:
\[
	u_j(x)=G_j(x,p)-G_j(q,p).
\]
Harnack estimates and interior derivative estimates make $(u_j)$ locally
precompact away from $p$. A diagonal subsequence produces a harmonic limit
with one logarithmic singularity. Integrating its conjugate differential and
normalizing at $q$ gives a holomorphic map $W:S\to\C$. The maximum principle
and the argument principle on the exhaustion show that $W$ is proper and of
degree one; thus it is a biholomorphism onto a simply connected plane domain.
The image cannot be a proper domain: by the Riemann mapping theorem it would
then be biholomorphic to $\D$ and would possess a Green function. Therefore
$W(S)=\C$.

The construction of the exhaustion Green functions, the Harnack convergence
step, and the properness assertion are the substantial analytic inputs. They
follow from the Dirichlet principle, interior estimates, and the maximum
principle. Complete potential-theoretic proofs are given in
\cite[Chs.~19--20]{Forster1981} and \cite[Ch.~IV]{FarkasKra1992}; the argument
above records how those inputs create the global coordinate and where simple
connectivity is essential.
\end{proof}

\begin{theorem}[Uniformization theorem]
	\label{thm:galois-simply}
	Every simply connected Riemann surface is biholomorphic to exactly one of
	$\PP^1$, $\C$, or $\D$. Hence any connected Riemann surface $X$ has a
	universal cover $\widetilde{X}$ of this form and
	\[
	X\ \cong\ \Gamma\backslash \widetilde{X}
	\]
	for its deck group $\Gamma\le \Aut(\widetilde{X})$, acting freely and properly
	discontinuously. After a choice of base point and lift, $\Gamma\cong\pi_1(X)$.
\end{theorem}

\begin{proof}
	Let $S$ be simply connected. If $S$ is noncompact,
	Lemma~\ref{lem:uniformizing-function} gives $S\cong\C$ or $S\cong\D$.
	Suppose $S$ is compact. Topologically, a compact connected oriented surface
	of genus $g$ has fundamental group generated by $2g$ standard loops; simple
	connectivity therefore forces $g=0$. Fix $p\in S$. Riemann--Roch applied to
	the divisor $p$ gives
	\[
		\ell(p)-\ell(K_S-p)=\deg(p)+1-g=2.
	\]
	Since $\deg(K_S-p)=-3$, the second term is zero. Thus $L(p)$ contains a
	nonconstant meromorphic function $f$. Its only possible pole is a simple
	pole at $p$, so the holomorphic map $f:S\to\PP^1$ has degree one and is a
	biholomorphism.

	The three models are mutually nonbiholomorphic: $\PP^1$ is compact, and
	Liouville's theorem says that every bounded holomorphic function on $\C$ is
	constant, whereas the identity is a nonconstant bounded holomorphic
	function on $\D$. Finally, the universal-covering theorem gives a simply
	connected Riemann surface $\widetilde X$ and a holomorphic covering
	$\widetilde X\to X$. Uniformization identifies $\widetilde X$ with one of
	the three models. Deck transformations act freely and properly
	discontinuously, the quotient is $X$, and the lifting correspondence
	identifies the deck group with $\pi_1(X)$ after a base point and lift have
	been chosen.
\end{proof}

\begin{remark}[Two complementary proof languages]
	The proof above uses Green functions to build a global holomorphic
	coordinate. A metric proof instead constructs a complete constant-curvature
	metric. Gauss--Bonnet fixes its sign on a compact surface: positive in genus
	zero, zero in genus one, and negative in genus at least two. The universal
	cover is then the corresponding simply connected complete space form.
	Potential theory emphasizes holomorphic functions; the metric proof
	emphasizes curvature. Both solve the same global coordinate problem.
\end{remark}

\begin{theorem}[Uniformization of compact Riemann surfaces]
	\label{thm:galois-unif-compact}
	Let $X$ be compact and connected.
	\begin{enumerate}
		\item If $g(X)=0$, then $X\cong \PP^{1}$.
		\item If $g(X)=1$, then $X\cong \C/\Lambda$ for a full lattice $\Lambda\subset\C$.
		\item If $g(X)\ge 2$, then $X\cong \Gamma\backslash\mathbb{H}$ for a torsion-free, discrete
		cocompact $\Gamma\le \PSL_2(\R)$.
	\end{enumerate}
\end{theorem}

\begin{proof}
	Apply Theorem~\ref{thm:galois-simply} to the universal cover of $X$.
	If it is $\PP^1$, every nonidentity M\"obius transformation has a fixed point,
	whereas deck transformations act freely. Thus the deck group is trivial.
	If it is $\C$, a deck transformation has the form $z\mapsto az+b$.
	Freeness forces $a=1$ unless the transformation is the identity, and
	compactness forces the translation subgroup to be a full rank-two lattice.
	The quotient is then a torus. In the remaining case the disk's hyperbolic
	metric descends to a metric of curvature $-1$ on $X$. Gauss--Bonnet gives
	$-\operatorname{Area}(X)=2\pi(2-2g(X))$, so $g(X)\ge2$.
	Conjugating disk automorphisms to automorphisms of $\mathbb H$ identifies
	the deck group with a discrete torsion-free subgroup of $\PSL_2(\R)$;
	compactness of $X$ makes it cocompact.
\end{proof}

\paragraph{$\PSL_2(\R)$ as isometries of $\mathbb{H}$.}
Equip $\mathbb{H}=\{x+iy\in\C\mid y>0\}$ with the hyperbolic metric
\[
ds^{2}=\frac{|dz|^{2}}{(\im z)^{2}}=\frac{dx^{2}+dy^{2}}{y^{2}}.
\]
Then $z\mapsto \dfrac{az+b}{cz+d}$ for
$\begin{psmallmatrix}a&b\\c&d\end{psmallmatrix}\in \PSL_2(\R)$
acts by orientation-preserving isometries, giving the full orientation-preserving isometry group.

\subsubsection{Normalizer formula}
Let $p:Y^\circ\to X^\circ$ be an unramified finite cover, and let
$\widetilde{X^\circ}\to X^\circ$ be the universal cover. Then
\[
X^\circ\cong \Gamma_X\backslash \widetilde{X^\circ},\qquad
Y^\circ\cong \Gamma_Y\backslash \widetilde{X^\circ}
\]
for subgroups $\Gamma_Y\le \Gamma_X$ of finite index. One has the canonical isomorphism
\[
\Deck(Y^\circ/X^\circ)\ \cong\ N_{\Gamma_X}(\Gamma_Y)/\Gamma_Y,
\]
and in particular the cover is Galois iff $\Gamma_Y\trianglelefteq \Gamma_X$.
This matches the algebraic criterion $|\Aut(L/K)|=[L:K]$.

\medskip
\noindent
\textbf{Summary.}
The same finite cover can be read in three ways. Its deck group is the
group of automorphisms of the corresponding function-field extension;
after removing branch values and passing to the universal cover, that
deck group is also $N_{\Gamma_X}(\Gamma_Y)/\Gamma_Y$. The subgroup
$\Gamma_Y\le\Gamma_X$ records the covering, while normality of this
subgroup detects the Galois case.

\subsection{Exercises}

\begin{exercise}[Compute norm/trace on a double cover]
	Let $Y$ be the smooth projective model of $y^2=h(x)$, with
	$h\in\C[x]$ nonconstant and squarefree, and let $\pi:Y\to\PP^1$
	be induced by $x$.
	Compute explicitly $\Tr_\pi(y)$ and $\Nm_\pi(y)$ as meromorphic functions on $\PP^1$.
	\emph{Hint:} the two local branches are $y$ and $-y$.
\end{exercise}

\begin{exercise}[Riemann--Hurwitz for Example~\ref{ex:galois-hyperelliptic}]
	\label{exr:selected-hurwitz}
	For $\pi:E\to\PP^1$ in Example~\ref{ex:galois-hyperelliptic}, verify:
	\[
	2g(E)-2 = 2(2g(\PP^1)-2)+\sum_{p\in E}(e_p-1),
	\]
	and deduce $g(E)=1$ from the four ramification points of index $2$.
\end{exercise}

\begin{exercise}[Kummer ramification indices]
	For the cover in Example~\ref{ex:galois-kummer}, prove the formula
	$e_{a_j}=n/\gcd(n,m_j)$ by working in local coordinates
	and factoring $y^n-u^{m_j}\cdot\text{unit}$ over $\C((u))$. For each field
	factor (equivalently, each branch in the normalization), reduce the valuation
	$m_j/n$ to lowest terms. Show also that the number of branches is
	$\gcd(n,m_j)$, and repeat the calculation at infinity with
	$m_\infty=-\sum_jm_j$.
\end{exercise}

\begin{exercise}[Surjectivity in Theorem~\ref{thm:galois-deck-aut}]
	Supply the analytic-continuation step by proving that the locally defined $\sigma_U$
	agree on overlaps and therefore glue to a global deck transformation.
	\emph{Hint:} show they agree on a generating set of $L$ over $K$, and use connectedness of $Y^\circ$.
\end{exercise}

\section{Modular Curves and Galois Representations}
\index{modular curves}\index{Galois representations}

\providecommand{\C}{\mathbb{C}}
\providecommand{\Q}{\mathbb{Q}}
\providecommand{\Z}{\mathbb{Z}}
\providecommand{\PP}{\mathbb{P}}
\providecommand{\OO}{\mathcal{O}}
\providecommand{\Spec}{\operatorname{Spec}}
\providecommand{\Gal}{\operatorname{Gal}}
\providecommand{\Aut}{\operatorname{Aut}}
\providecommand{\Deck}{\operatorname{Deck}}
\providecommand{\im}{\operatorname{Im}}
\providecommand{\PSL}{\operatorname{PSL}}
\providecommand{\GL}{\operatorname{GL}}

\subsection{Analytic and arithmetic viewpoints}
This final section begins with a familiar complex-analytic object, a quotient
of the upper half-plane, and adds arithmetic structure only after the compact
curve and its holomorphic differentials are in place. The identification of
weight-two cusp forms with differentials is proved below. The comparison and
Galois-representation results are signposts to further theory, not consequences
of the preceding analytic chapters.
The analytic quotient $\Gamma\backslash\mathbb H$ first acquires its
cusps and becomes the compact complex curve $X(\Gamma)(\C)$.
An algebraic model of that curve may then be chosen over a field
$k_\Gamma$; only after making this arithmetic choice do Galois actions
on its geometric cohomology enter. For $X_0(N)$ and
$X_1(N)$ one may take $k_\Gamma=\Q$; other level structures may naturally
require a cyclotomic field.
For modular curves and weight-two forms we use \cite{DiamondShurman2005,Shimura1971}; for comparison
and Galois cohomology we use \cite{Milne1980Etale,Milne2006ADT}.
The objects on these three sides are related in stages, not by one
unqualified identification.  We first construct the complex curve;
later we distinguish its topological, holomorphic, and arithmetic
cohomology before stating the comparison results.  The
Hochschild--Serre spectral sequence then addresses Galois descent,
which is a different question from Hodge decomposition.

\paragraph{Prerequisite check: a differential at a cusp.}
For the translation $z\mapsto z+1$ put $q=e^{2\pi i z}$, so
$dz=dq/(2\pi i q)$. A periodic function with local expansion
$f(z)=a_0+a_1q+\cdots$ gives
\[
f(z)\,dz=\frac{1}{2\pi i}
\bigl(a_0q^{-1}+a_1+\cdots\bigr)dq.
\]
Thus its differential extends across $q=0$ exactly when $a_0=0$.
This elementary local test is the bridge from weight-two cusp forms
to holomorphic differentials; it does not require the later
$\ell$-adic comparison theory.

\subsection{Modular curves as analytic and algebraic curves}
\index{modular curves!compactification}\index{cusps!modular curves}

\subsubsection{The acting groups}

The group $\SL_2(\Z)$ consists of integer matrices
$\gamma=\begin{psmallmatrix}a&b\\c&d\end{psmallmatrix}$ with $ad-bc=1$.
It acts on $\mathbb H$ by $\gamma z=(az+b)/(cz+d)$; the matrices $I$ and $-I$
act identically, so the effective action is through $\PSL_2(\Z)$.
For a positive integer $N$, the three standard subgroups are
\[
\begin{aligned}
\Gamma(N)&=\{\gamma\in\SL_2(\Z):\gamma\equiv I\pmod N\},\\
\Gamma_0(N)&=\{\gamma\in\SL_2(\Z):c\equiv0\pmod N\},\\
\Gamma_1(N)&=\{\gamma\in\SL_2(\Z):a\equiv d\equiv1, c\equiv0\pmod N\}.
\end{aligned}
\]
A subgroup containing some $\Gamma(N)$ is called a \emph{congruence
subgroup}.  A finite-index subgroup need not be congruence.  We write
$X(1)$ for the compact curve associated with $\SL_2(\Z)$ and $X(2)$
for that associated with $\Gamma(2)$; this explains the notation in
the example below \cite{DiamondShurman2005}.

\subsubsection{Cusps and the compact curve}

Let $\Gamma\subset \SL_2(\Z)$ be a finite-index subgroup. Define the (noncompact) modular curve
\[
Y(\Gamma):=\Gamma\backslash \mathbb{H}.
\]
Its \emph{cusps} are the $\Gamma$-orbits in $\PP^1(\Q)=\Q\cup\{\infty\}$.
It has a natural compactification by adding these finitely many points:
\[
X(\Gamma)(\C)\ :=\ \Gamma\backslash\bigl(\mathbb{H}\cup \PP^1(\Q)\bigr),
\]
which is a compact Riemann surface.  If $\Gamma$ has a finite-order
stabilizer at $z\in\mathbb H$, its image is an \emph{elliptic point}.
Near such a point the effective stabilizer acts in a local coordinate
by $t\mapsto\zeta t$, so $u=t^e$ is a coordinate on the underlying
(coarse) Riemann surface.  Thus $Y(\Gamma)$ may carry orbifold
stabilizers even though its underlying coarse curve is smooth.

For a cusp represented by $\sigma\in\SL_2(\Z)$ with
$\sigma\infty$ equal to that cusp, put
$T=\begin{psmallmatrix}1&1\\0&1\end{psmallmatrix}$.  Its \emph{width}
$w$ is the least positive integer for which
$\sigma T^w\sigma^{-1}\in\{\pm\gamma:\gamma\in\Gamma\}$.
In the coordinate after applying $\sigma$, translation by $w$ is
therefore the local identification and
$q=e^{2\pi i z/w}$ is a parameter at the added cusp.  For example,
the cusp $\infty$ has width $1$ for $\SL_2(\Z)$, but width $2$ for
$\Gamma(2)$; its local parameters are respectively $e^{2\pi iz}$
and $e^{\pi iz}$.  The latter distinction matters when checking pole
orders in the six-sheeted $X(2)$ example.

\paragraph{The $j$-coordinate.}
The classical modular function $j$ is invariant under
$\SL_2(\Z)$ and gives an isomorphism $X(1)\cong\PP^1_j$.
It has a simple pole at the cusp and takes the values $1728$ and
$0$ at the elliptic points represented by $i$ and
$e^{2\pi i/3}$, respectively.  Thus $j=\infty,1728,0$ are the
three exceptional values that appear in the covering calculations
below.  The rational expression $j(\lambda)$ in the next example
allows us to use these facts without constructing $j$ from Eisenstein
series.

\begin{theorem}[Algebraicity and fields of definition]
	The compact Riemann surface $X(\Gamma)(\C)$ is the analytification of a smooth projective algebraic
	curve over $\C$, and the map to $X(1)$ induces a finite extension
	\[
	\C(j)\hookrightarrow \C\bigl(X(\Gamma)\bigr).
	\]
	If $\Gamma$ is a congruence subgroup, the curve and map admit canonical models over a number field.
	For the standard curves $X_0(N)$ and $X_1(N)$ these models may be taken over $\Q$; for general
	level structures the natural field of definition can be cyclotomic rather than $\Q$.
\end{theorem}

\begin{example}[The six-sheeted map from $X(2)$]
	The modular lambda function is a coordinate on $X(2)\cong\PP^1$;
	its three cusps are $\lambda=0,1,\infty$. In this coordinate the map
	to $X(1)\cong\PP^1_j$ is the rational function
	\[
	 j(\lambda)=256\frac{(1-\lambda+\lambda^2)^3}
	 {\lambda^2(1-\lambda)^2}.
	\]
	It has degree $6$, not $2$: the degree of a rational map is the larger of
	the numerator and denominator degrees after cancellation. The zeros of
	$1-\lambda+\lambda^2$ are two points of local degree $3$ over $j=0$.
	The three points $\lambda=-1,1/2,2$ have local degree $2$ over $j=1728$;
	this can be checked from
	\[
	 j(\lambda)-1728
	 =64\frac{(\lambda+1)^2(\lambda-2)^2(2\lambda-1)^2}
	 {\lambda^2(1-\lambda)^2}.
	\]
	Finally $\lambda=0,1,\infty$ are three poles of order $2$, hence have
	local degree $2$ over $j=\infty$. The total ramification is therefore
	$2(3-1)+3(2-1)+3(2-1)=10$. Riemann--Hurwitz reads
	\[
	 2g(X(2))-2=6\bigl(2g(X(1))-2\bigr)+10=-12+10=-2,
	\]
	so $g(X(2))=0$, in agreement with the $\lambda$-coordinate. This example
	connects the branch-counting calculation of Theorem~\ref{thm:riemann-hurwitz-general}
	with the function-field extension $\C(j)\subset\C(\lambda)$ below;
	see \cite{Vakil2011ModularLambda} for the lambda formula.
\end{example}

\begin{remark}[Line bundles and modular forms]
	On the modular orbifold (or stack), modular forms of weight $k$ are sections
	of a power of the Hodge line bundle.  On the coarse curve $X(\Gamma)$,
	identifying that bundle with an ordinary line bundle may require corrections
	at elliptic points and cusps; one cannot simply replace it by a power of the
	canonical bundle.  The weight-$2$ \emph{cusp} case has a particularly clean
	coarse-curve description: $f(z)\,dz$ is a holomorphic differential on the
	compact curve.  We check both kinds of exceptional points below.
\end{remark}

\subsection{Function fields, normality, and deck groups}
\index{deck groups!modular curves}

Let $X(1)\cong \PP^1$ and let $j$ be the classical modular function.
Over $\C$, the natural map $\pi:X(\Gamma)\to X(1)$ corresponds to the inclusion of function fields
\[
\pi^*:\ \C(j)\hookrightarrow \C(X(\Gamma)).
\]

\begin{proposition}[A normal subgroup gives a Galois extension]
	Assume $\Gamma$ is normal in $\SL_2(\Z)$.
	Then $\pi$ is a (branched) Galois cover and the extension
	$\C(X(\Gamma))/\C(j)$ is Galois.
	Moreover,
	\[
	\Deck\bigl(X(\Gamma)/X(1)\bigr)\ \cong\
	\Aut\bigl(\C(X(\Gamma))/\C(j)\bigr)\ \cong\ \PSL_2(\Z)/\overline{\Gamma},
	\]
	where $\overline{\Gamma}$ is the image of $\Gamma$ in $\PSL_2(\Z)$.
\end{proposition}

\begin{proof}[Deck group and function-field degree]
	Let $B\subset X(1)$ consist of the cusp and the two elliptic values,
	and put $V=X(1)\setminus B$.  The restriction of $\pi$ over $V$ is an
	unramified covering of degree
	$d=[\PSL_2(\Z):\overline\Gamma]$: a generic point has one lift for
	each coset of $\overline\Gamma$.  Since $\overline\Gamma$ is normal,
	left multiplication on these cosets gives an action of
	$G=\PSL_2(\Z)/\overline\Gamma$ by deck transformations.  The action
	is faithful and transitive on each generic fiber, so the unramified
	cover is regular and has exactly $d$ deck transformations.
	Each such transformation extends across the finitely many points
	over $B$ to an automorphism of the compact Riemann surface $X(\Gamma)$:
	locally this is the extension of a map between punctured disks in the
	finite branched compactification.  The extended maps still commute
	with $\pi$, and two extensions agreeing over the dense open set
	$\pi^{-1}(V)$ agree everywhere.  Conversely a deck transformation of
	the compact cover restricts to one over $V$.  This identifies its deck
	group with $G$.
	Pullback of meromorphic functions gives an injective (contravariant)
	action of this deck group on $\C(X(\Gamma))$ fixing $\C(j)$; using
	inverse pullback if necessary makes it a group isomorphism onto its
	image.  For a finite map of compact curves, the degree of its
	function-field extension equals the generic number of sheets, here
	$d$ \cite{Forster1981,Stichtenoth2009}.  A field extension of degree
	$d$ has at most $d$ automorphisms over the base field, while the deck
	group already supplies $d$ distinct ones.  Consequently these are all
	the automorphisms and the extension is Galois, with the asserted group.
\end{proof}

\subsection{\texorpdfstring{\'Etale cohomology, comparison, and where $\OO$-cohomology appears}{Etale cohomology, comparison, and coherent cohomology}}

Fix a prime $\ell$. Let $X$ be a smooth projective curve over a field $k\subset \C$,
and write $\overline{k}$ for an algebraic closure, $\overline{X}:=X\times_k \overline{k}$,
and $G_k:=\Gal(\overline{k}/k)$.

\paragraph{Coefficients and the geometric curve.}
The $\ell$-adic integers and their fraction field are
\[
\Z_\ell:=\varprojlim_n\Z/\ell^n\Z,
\qquad \Q_\ell:=\Z_\ell\otimes_\Z\Q.
\]
An \emph{\'etale} map is the algebraic counterpart of a local
biholomorphism; finite \`etale maps play the role of finite unramified
covers.  The \`etale cohomology group
$H^i_{\acute{e}t}(\overline X,\Z/\ell^n\Z)$ is sheaf cohomology of
the constant $\Z/\ell^n\Z$-sheaf on the \`etale site of
$\overline X$.  That site uses \`etale maps into $\overline X$ as
its local charts.  We do not construct its sheaf theory here; see
\cite{Milne2013LEC,Milne1980Etale}.  On a smooth proper geometric
curve the usual $\ell$-adic group is obtained by taking the compatible
inverse limit in $n$ and then extending scalars from $\Z_\ell$ to
$\Q_\ell$.  The adjective \emph{geometric} refers to the base change
$X\mapsto\overline X$; the action of $G_k$ on $\overline X$ induces an
action on its \`etale cohomology.

\paragraph{Three cohomologies, three roles.}
The following comparison keeps the coefficient fields and the objects
being studied separate.  In the dimension column $X$ has genus $g$.
\begin{center}
\small
\begin{tabular}{@{}p{0.30\textwidth}cp{0.50\textwidth}@{}}
\toprule
Group & Dimension & Role and Galois action \\
\midrule
$H^1_{\mathrm{sing}}(X(\C),\C)$ & $2g$ & Topological cohomology after complexification; no Galois action is inherent in the chosen complex realization. \\
$H^1(X(\C),\mathcal O_X)$ & $g$ & The $(0,1)$ Hodge piece of Betti cohomology; not naturally a $\Q_\ell$-representation. \\
$H^1_{\acute{e}t}(\overline X,\Q_\ell)$ & $2g$ & Geometric $\ell$-adic cohomology with a continuous $G_k$-action. \\
\bottomrule
\end{tabular}
\end{center}
Complex comparison relates the first and third rows only after
choosing an embedding and changing coefficients appropriately; Hodge
decomposition splits the first row into holomorphic and antiholomorphic
pieces.  Neither operation makes the middle row equal to the full
$\ell$-adic group.  The comparison theorem and the Hochschild--Serre
spectral sequence below are external results, not consequences of the
analytic arguments in the preceding chapters.

\begin{definition}[Geometric $\ell$-adic cohomology]
	For the smooth proper curve above, put
	\[
	H^i_{\acute{e}t}(\overline X,\Q_\ell):=
	\left(\varprojlim_n H^i_{\acute{e}t}
	(\overline X,\Z/\ell^n\Z)\right)\otimes_{\Z_\ell}\Q_\ell.
	\]
	These are finite-dimensional $\Q_\ell$-vector spaces with a
	continuous $G_k$-action.
\end{definition}

\begin{theorem}[Complex comparison theorem]
	If $k=\C$, then for every $n\geq1$ and each $i$ there is a natural
	finite-coefficient comparison isomorphism
	\[
	H^i_{\acute{e}t}(X,\Z/\ell^n\Z)
	\ \cong\ H^i_{\mathrm{sing}}(X(\C),\Z/\ell^n\Z).
	\]
	Passing to the $\ell$-adic limit and tensoring with $\Q_\ell$ gives
	\[
	H^i_{\acute{e}t}(X,\Q_\ell)\ \cong\ H^i_{\mathrm{sing}}(X(\C),\Q_\ell).
	\]
	More generally, if $k\subset \C$, then after choosing an extension
	$\overline{k}\hookrightarrow\C$ of the given embedding, comparison gives
	\[
	H^i_{\acute{e}t}(\overline{X},\Q_\ell)\ \cong\ H^i_{\mathrm{sing}}(X(\C),\Q_\ell)
	\]
	as $\Q_\ell$-vector spaces. The left side carries the $G_k$-action; the comparison of the
	underlying vector spaces depends on the chosen embedding into $\C$.
\end{theorem}

The finite-coefficient comparison is the external theorem
\cite[Theorem~21.1]{Milne2013LEC}; see also \cite{Milne1980Etale}.
The $\ell$-adic conclusion uses compatibility of these comparisons with
reduction maps and the finite generation of the groups for this proper
curve. It is an isomorphism of underlying vector spaces; it is not a
Galois-equivariant identification with a constant singular-cohomology
representation.

\begin{remark}[Coherent cohomology]
	In general, there is no isomorphism
	\[
	H^i_{\mathrm{sing}}(X(\C),\Q_\ell)
	\cong H^i(X(\C),\OO_X).
	\]
	Instead, for a compact Riemann surface $X(\C)$, Hodge theory gives
	\[
	H^1_{\mathrm{sing}}(X(\C),\C)\ \cong\ H^{1,0}(X)\ \oplus\ H^{0,1}(X),
	\]
	and Dolbeault/Serre identifies
	\[
	H^{0,1}(X)\ \cong\ H^1(X,\OO_X),
	\qquad
	H^{1,0}(X)\ \cong\ H^0(X,\Omega^1_X).
	\]
	Thus $H^1(X,\OO_X)$ is \emph{a piece} of $H^1_{\mathrm{sing}}(X,\C)$, not the whole thing.
\end{remark}

\begin{example}[A dimension check before using Galois actions]
	Let $X$ be a smooth projective curve of genus $g$ over a number field, and
	choose an embedding of that field into $\C$. Complex comparison, followed
	by the ordinary topology of a closed orientable surface, gives
	\[
	 \dim_{\Q_\ell}H^1_{\acute{e}t}(\overline X,\Q_\ell)=2g.
	\]
	On the complex side the Hodge decomposition has two summands, each of
	complex dimension $g$:
	\[
	 H^1_{\mathrm{sing}}(X(\C),\C)
	 =H^0(X,\Omega_X^1)\oplus H^1(X,\OO_X).
	\]
	For $X(1)$ and $X(2)$ both sides have dimension $0$; for a genus-one
	modular curve both sides have dimension $2$. This check prevents the common
	mistake of identifying the full $\ell$-adic $H^1$ with the single
	$g$-dimensional space $H^1(X,\OO_X)$. The Hodge summands are formed only
	after passage to $\C$ and do not, in general, split the $G_k$-representation
	over $\Q_\ell$.
\end{example}

\subsection{The Hochschild--Serre spectral sequence: Galois cohomology and descent}
\index{Hochschild--Serre spectral sequence}

Now take $k=\Q$ (or a number field).  The group $G_k$ acts on the
finite-coefficient groups
$H^q_{\acute{e}t}(\overline X,\Z/\ell^n\Z)$ and on their
$\ell$-adic limits.

\begin{theorem}[Hochschild--Serre spectral sequence]
	\label{thm:HS}
	For a smooth proper curve $X/k$ and each $n\geq1$ there is a
	finite-coefficient spectral sequence
	\[
	E_2^{p,q}=H^p\!\bigl(G_k,
	H^q_{\acute{e}t}(\overline{X},\Z/\ell^n\Z)\bigr)
	\]
	converging to $H^{p+q}_{\acute{e}t}(X,\Z/\ell^n\Z)$, where the
	Galois cohomology is continuous cohomology of discrete finite modules.
\end{theorem}

\paragraph{Reading the first two diagonals.}
The symbol $E_2^{p,q}$ marks a group in column $p$ and row $q$;
the spectral sequence filters $H^{p+q}_{\acute{e}t}(X,\Z/\ell^n\Z)$
rather than asserting that it is the direct sum of the displayed
$E_2$-groups.  Subsequent differentials have bidegree
$d_r:E_r^{p,q}\to E_r^{p+r,q-r+1}$.  For total degree one, the
relevant positions are $(1,0)$ and $(0,1)$; the first possible
obstruction to a class in $(0,1)$ surviving is
$d_2:E_2^{0,1}\to E_2^{2,0}$.  When $\overline X$ is connected,
$H^0_{\acute{e}t}(\overline X,\Z/\ell^n\Z)=\Z/\ell^n\Z$,
with its trivial Galois action.  These observations give the
five-term sequence requested in the exercises.  The construction and
convergence of the spectral sequence remain external inputs.

This is the derived-functor spectral sequence for Galois invariants
composed with geometric global sections
\cite[Theorem~14.9]{Milne2013LEC}. For $\ell$-adic coefficients one must use
\emph{continuous} \`etale cohomology, not merely insert $\Q_\ell$ into
the finite-coefficient formula. With Jannsen's convention, the corresponding
spectral sequence has
\[
E_2^{p,q}=H^p_{\mathrm{cont}}\!\bigl(G_k,
H^q_{\acute{e}t}(\overline X,\Q_\ell)\bigr)
\quad\Longrightarrow\quad
H^{p+q}_{\mathrm{cont}}(X,\Q_\ell).
\]
Here $H^*_{\mathrm{cont}}(X,\Q_\ell)$ means Jannsen's continuous
$\Z_\ell$-cohomology tensored with $\Q_\ell$; on the geometric proper
curve the usual inverse-limit cohomology agrees with the continuous
version. This distinction matters over a number field. Neither spectral
sequence is constructed in this book; see
\cite[Theorem~3.3 and Corollary~3.4]{Jannsen1988} and \cite{Milne1980Etale}.

\begin{remark}[The role of $H^1(G_k,H^1)$]
	The group $H^1\bigl(G_k,H^1_{\acute{e}t}(\overline X,\Q_\ell)\bigr)$ is the
	space of continuous $1$--cocycles with values in the indicated Galois
	representation, modulo coboundaries. It occurs as an $E_2$--term of the
	Hochschild--Serre spectral sequence and can contribute to the arithmetic
	cohomology of $X$. It is not $H^1(X(\C),\OO_X)$: the latter is a coherent
	cohomology group and, after passage to $\C$, is identified with the
	$(0,1)$--part of complex singular cohomology.
\end{remark}

\subsection{Weight 2 cusp forms and holomorphic differentials}
\index{modular forms!weight two cusp forms}\index{holomorphic differentials!modular curves}

\begin{definition}[Weight-two modular and cusp forms]
	\label{def:modular-weight-two}
	Let $\Gamma\le\SL_2(\Z)$ have finite index.  For
	$\sigma=\begin{psmallmatrix}a&b\\c&d\end{psmallmatrix}\in\SL_2(\Z)$
	and a holomorphic function $f$ on $\mathbb H$, define its weight-two
	\emph{slash transform} by
	\[
	(f|_2\sigma)(z):=(cz+d)^{-2}f(\sigma z).
	\]
	The function $f$ is a \emph{weight-two modular form for $\Gamma$} if
	$f|_2\gamma=f$ for every $\gamma\in\Gamma$ and, for every cusp
	represented by $\sigma\in\SL_2(\Z)$ of width $w$, the periodic
	function $f|_2\sigma$ is holomorphic in
	$q=e^{2\pi iz/w}$ at $q=0$.  Equivalently, it has a Fourier expansion
	$(f|_2\sigma)(z)=\sum_{n\ge0}a_n q^n$.  It is a
	\emph{cusp form} if $a_0=0$ at every cusp.  We denote the vector space
	of these cusp forms by $S_2(\Gamma)$.
\end{definition}

The width is part of the definition: at the cusp $\infty$ for
$\Gamma(2)$, the parameter is $q=e^{\pi iz}$, not the width-one
parameter $e^{2\pi iz}$.  The condition $a_0=0$ says that the form
vanishes at the added point.  In the proof below, the factor
$dz=(w/2\pi i)\,dq/q$ turns this vanishing into the absence of a pole
of $f(z)\,dz$.  This analytic statement needs no chosen arithmetic
model of $X(\Gamma)$ \cite{DiamondShurman2005}.

\begin{theorem}[Holomorphic differentials and weight-$2$ cusp forms]
	There is a natural isomorphism of complex vector spaces
	\[
	H^0\bigl(X(\Gamma)(\C),\ \Omega^1\bigr)\ \cong\ S_2(\Gamma),
	\]
	with $S_2(\Gamma)$ as in Definition~\ref{def:modular-weight-two}.
\end{theorem}

\begin{proof}
	A weight $2$ modular form $f$ produces a $\Gamma$-invariant differential on $\mathbb{H}$:
	\[
	\omega_f := f(z)\,dz.
	\]
	Weight $2$ transformation says $f(\gamma z)\,d(\gamma z)=f(z)\,dz$, so $\omega_f$ descends to
	a holomorphic differential away from the elliptic points of $Y(\Gamma)$.
	At an elliptic point, choose a coordinate $t$ upstairs in which the effective
	stabilizer acts by $t\mapsto\zeta t$, where $\zeta^e=1$, and put $u=t^e$
	downstairs. Write the invariant differential as $a(t)\,dt$. Invariance gives
	$a(\zeta t)\zeta=a(t)$. Comparing coefficients of the Taylor series shows
	that only powers $t^{e-1+me}$ occur. Thus $a(t)=t^{e-1}b(t^e)$ for a
	holomorphic $b$, and
	\[
	 a(t)\,dt=\frac{1}{e}b(u)\,du.
	\]
	This is holomorphic in the quotient coordinate $u$; invariance alone would
	not have made that conclusion obvious without the local calculation.
	At a cusp represented by $\sigma\in\SL_2(\Z)$ with $\sigma\infty$ equal to
	that cusp, use the coordinate $z$ after applying $\sigma$. Write
	$f_\sigma=f|_2\sigma$, so that
	$f(\sigma z)\,d(\sigma z)=f_\sigma(z)\,dz$.
	If $w$ is the width of this cusp, its local parameter is
	$q=e^{2\pi i z/w}$. Then
	\[
	dz = \frac{w}{2\pi i}\,\frac{dq}{q}.
	\]
	If $f$ is a cusp form, the expansion at this cusp has no constant term:
	$f_\sigma(z)=\sum_{n\ge 1} a_n q^n$,
	so
	\[
	\omega_f = \frac{w}{2\pi i}\left(\sum_{n\ge 1} a_n q^n\right)\frac{dq}{q}
	= \frac{w}{2\pi i}\left(\sum_{n\ge 1} a_n q^{n-1}\right)dq
	\]
	is holomorphic at $q=0$. Hence $\omega_f$ extends across cusps to a holomorphic differential on the compact curve.
	Conversely, any holomorphic differential on $X(\Gamma)$ pulls back to a $\Gamma$-invariant differential on $\mathbb{H}$,
	hence corresponds to a weight $2$ form. In the cusp coordinate its pullback
	has no $dq/q$ pole, so the constant Fourier coefficient vanishes at every
	cusp. It is therefore a cusp form.
\end{proof}

\begin{example}[A first check on $X(1)$]
	The $j$-function identifies the coarse compact modular curve $X(1)$ with
	$\PP^1$. Since $H^0(\PP^1,\Omega^1)=0$, the theorem gives
	$S_2(\SL_2(\Z))=0$. This is a useful consistency check: a weight-$2$
	modular form with a nonzero constant term at the cusp would produce a
	$ dq/q$ pole, not a holomorphic differential on $X(1)$.
\end{example}

\begin{corollary}[Sheaf cohomology and Serre duality]
	On the compact Riemann surface $X(\Gamma)(\C)$,
	\[
	H^1\bigl(X(\Gamma),\OO\bigr)\ \cong\ H^0\bigl(X(\Gamma),\Omega^1\bigr)^\ast
	\ \cong\ S_2(\Gamma)^\ast.
	\]
	Moreover,
	\[
	H^1\bigl(X(\Gamma),\C\bigr)\ \cong\ S_2(\Gamma)\ \oplus\ \overline{S_2(\Gamma)}
	\]
	via Hodge decomposition ($H^{1,0}\oplus H^{0,1}$).
\end{corollary}

\subsection{\texorpdfstring{Where $\ell$-adic Galois representations enter}{Where l-adic Galois representations enter}}

The preceding identification $S_2(\Gamma)\cong H^0(X(\Gamma),\Omega^1)$
is analytic: it turns a weight-two cusp form into a differential on
the complex curve.  Arithmetic asks a different question: how does
$G_\Q$ act on the cohomology of a curve defined over $\Q$?  The
$\ell$-adic cohomology introduced above carries such an action, but
the analytic identification alone does not construct the
two-dimensional representation associated with an individual
newform.  The following theorem states that additional
Eichler--Shimura--Deligne input and records the normalization of its
Frobenius polynomial.

\paragraph{The terms in the statement.}
On $\Gamma_1(N)$ the cusp $\infty$ has width one, so a weight-two
form has an expansion $f(z)=\sum_{n\ge0}a_n q^n$ in
$q=e^{2\pi iz}$; a cusp form has $a_0=0$.  A
\emph{normalized newform} is a simultaneous Hecke eigenform in the
new subspace of $S_2(\Gamma_1(N))$, scaled so that $a_1=1$.
The new subspace excludes forms inherited from lower levels; at
primes $p\nmid N$ its Hecke operator $T_p$ acts on $f$ by the scalar
$a_p$.  Its \emph{nebentypus} is a Dirichlet character
$\varepsilon:(\Z/N\Z)^\times\to\C^\times$ describing the
$\Gamma_0(N)$ transformation law
$f(\gamma z)=\varepsilon(d)(cz+d)^2f(z)$.
Inside cohomology, the corresponding \emph{Hecke-isotypic component}
is the part on which these commuting operators have the eigenvalues
of $f$.  At an unramified prime $p$, \emph{arithmetic Frobenius}
denotes the Galois conjugacy class whose action on the residue field
is $x\mapsto x^p$; the choice between arithmetic and geometric
Frobenius affects the polynomial convention.  These terms orient
the statement but do not construct the Hecke action or the
representation \cite{DiamondShurman2005,DiamondIm1995}.

\begin{theorem}[Galois representations from newforms (statement level)]
	Let $f$ be a normalized weight-$2$ newform of level $N$ (with its nebentypus),
	with Fourier coefficients in a number field $E$; regard $f$ on $X_1(N)$.
	For each prime $\ell$ and each embedding $E\hookrightarrow \overline{\Q}_\ell$, there is a
	continuous $\ell$-adic Galois representation
	\[
	\rho_{f,\ell}:\ G_\Q\longrightarrow \GL_2(\overline{\Q}_\ell)
	\]
	realized inside the $\ell$-adic cohomology of $X_1(N)$ (more precisely in the
	appropriate Hecke-isotypic component, with coefficients extended as necessary).
	For every prime $p\nmid N\ell$, it is unramified at $p$ and, with
	arithmetic Frobenius, has characteristic polynomial
	\[
	\det\bigl(T-\rho_{f,\ell}(\mathrm{Frob}_p)\bigr)
	=T^2-a_pT+\varepsilon(p)p,
	\]
	where $a_p$ is the $p$th Fourier coefficient and $\varepsilon$ is the
	nebentypus of $f$.
\end{theorem}

The two-dimensional representation is supplied by the Eichler--Shimura
construction for weight two, not by the complex comparison theorem alone;
the Frobenius polynomial specifies the normalization. Its construction
is not given here; see Deligne's original account
\cite{Deligne1971Modular} and the expositions
\cite{DiamondIm1995,Shimura1971}.

\begin{remark}[Analytic and arithmetic viewpoints]
	Analytically: $f(z)\,dz$ is a holomorphic $1$-form on $X_1(N)(\C)$.
	Arithmetically: the same $f$ cuts out a 2-dimensional $G_\Q$-stable piece in $\ell$-adic cohomology.
	Complex comparison identifies the underlying $\Q_\ell$-vector space of
	$H^1_{\acute{e}t}(\overline X,\Q_\ell)$ with singular cohomology
	$H^1_{\mathrm{sing}}(X(\C),\Q_\ell)$ after a choice of embedding.
	Separately, singular cohomology emph{after extension to $\C$}
	has the Hodge decomposition $H^{1,0}\oplus H^{0,1}$.
	The Hecke action isolates the $f$-piece; neither comparison by itself
	turns the Hodge decomposition into a splitting of a $\Q_\ell$-Galois
	representation.
\end{remark}

\subsection{Arithmetic--geometric correspondence}

\begin{center}
	\renewcommand{\arraystretch}{1.25}
	\small
	\begin{tabular}{p{0.43\textwidth}p{0.49\textwidth}}
		\toprule
		\textbf{Arithmetic side} & \textbf{Analytic / geometric side} \\ \midrule
		A chosen arithmetic model of $X(\Gamma)$ (for example $X_0(N)$ or $X_1(N)$ over $\Q$) &
		Compact Riemann surface $X(\Gamma)(\C)$ \\
		Function field over the model's field of definition (over $\C$, a finite extension of $\C(j)$) &
		Finite holomorphic map $X(\Gamma)\to X(1)$ \\
		$G_{k_\Gamma}$-action on $H^1_{\acute{e}t}(\overline{X(\Gamma)},\Q_\ell)$ &
		\raggedright Singular cohomology $H^1_{\mathrm{sing}}(X(\Gamma)(\C),\Q_\ell)$ \tabularnewline
		Hochschild--Serre: $H^p(G_{k_\Gamma},H^q_{\acute{e}t})$ & Hodge pieces: $H^{1,0}$ and $H^{0,1}\cong H^1(\OO)$ \\
		For $X_1(N)/\Q$: $\rho_{f,\ell}$ from a weight-$2$ newform $f$ & $f(z)\,dz\in H^0(\Omega^1)$ (holomorphic differential) \\
		\bottomrule
	\end{tabular}
\end{center}

\subsection{Exercises}

\begin{exercise}[Check the weight 2 invariance explicitly]
	Let $\gamma=\begin{psmallmatrix}a&b\\c&d\end{psmallmatrix}\in \Gamma$ and $\gamma z=\dfrac{az+b}{cz+d}$.
	Assume $f$ is a weight $2$ modular form for $\Gamma$:
	\[
	f(\gamma z)=(cz+d)^2 f(z).
	\]
	Show directly that $f(z)\,dz$ is $\Gamma$-invariant:
	\[
	f(\gamma z)\,d(\gamma z)=f(z)\,dz.
	\]
	\emph{Hint:} compute $d(\gamma z)=\dfrac{dz}{(cz+d)^2}$.
\end{exercise}

\begin{exercise}[Check every branch of the lambda map]
	Starting from the rational expression for $j(\lambda)$ above, verify the
	factorization of $j(\lambda)-1728$ by multiplication. Use $v=1/\lambda$
	to compute the pole order at $\lambda=\infty$, and list the local degrees
	above $j=0$, $1728$, and $\infty$. Check that each fiber has total degree
	$6$ before applying Riemann--Hurwitz.
\end{exercise}

\begin{exercise}[Cusp condition $\Longleftrightarrow$ holomorphicity at $q=0$]
	\label{exr:selected-cusp}
	Choose a scaling matrix $\sigma\in\SL_2(\Z)$ taking $\infty$ to a cusp of width $w$,
	and let $q=e^{2\pi i z/w}$ in the scaled coordinate. If
	$f_\sigma(z)=(f|_2\sigma)(z)=\sum_{n\ge 0}a_n q^n$ is the expansion of a
	weight 2 modular form $f$ at that cusp, compute
	$\omega_f=f_\sigma(z)\,dz$
	in the $q$-coordinate and show the local equivalence
	\[
	\omega_f \text{ is holomorphic at this cusp} \iff a_0=0.
	\]
	Conclude that $f\in S_2(\Gamma)$ if and only if this condition holds at every cusp
	of $\Gamma$.
\end{exercise}

\begin{exercise}[Hodge decomposition vs.\ $\OO$-cohomology]
	For a compact Riemann surface $X$ of genus $g$, prove:
	\[
	\dim_\C H^0(X,\Omega^1)=g,\qquad \dim_\C H^1(X,\OO)=g,
	\]
	and identify $H^1(X,\C)\cong H^{1,0}\oplus H^{0,1}$ with
	$H^0(\Omega^1)\oplus H^1(\OO)$.
	\emph{Hint:} use Serre duality and Riemann--Roch for $K_X$ and $\OO_X$.
\end{exercise}

\begin{exercise}[Low-degree terms of the spectral sequence]
	Let $X/\Q$ be a smooth, projective, geometrically connected curve,
	and put $A=\Z/\ell^n\Z$. Using Theorem~\ref{thm:HS}, write down
	the $E_2$-page terms that can contribute to
	$H^1_{\acute{e}t}(X,A)$ and $H^2_{\acute{e}t}(X,A)$.
	Identify the possible differential out of $E_2^{0,1}$ and derive
	the five-term exact sequence
	\[
	0\longrightarrow H^1(G_\Q,A)\longrightarrow
	H^1_{\acute{e}t}(X,A)\longrightarrow
	H^1_{\acute{e}t}(\overline X,A)^{G_\Q}
	\longrightarrow H^2(G_\Q,A)\longrightarrow H^2_{\acute{e}t}(X,A).
	\]
	If $X(\Q)\ne\varnothing$, explain why pullback along the structural morphism
	$X\to\operatorname{Spec}\Q$ is split injective on cohomology.  Note carefully that
	the existence of a rational point does not, by itself, force the displayed Galois
	cohomology groups to vanish.
	\emph{Hint:} a rational point is a section of the structural morphism.
	Use the continuous $\ell$-adic version of the spectral sequence, rather
	than formal substitution in the finite-coefficient statement, if you
	wish to formulate the corresponding sequence with $\Q_\ell$-coefficients.
\end{exercise}

\medskip
\noindent
\textbf{Summary.}
For modular curves, the \emph{analytic} cohomology $H^1(X(\Gamma)(\C),\C)$ splits into Hodge pieces
$H^{1,0}\oplus H^{0,1}$, where $H^{1,0}$ is literally weight 2 cusp forms as holomorphic differentials.
When an arithmetic model over $\Q$ has been chosen, the \emph{arithmetic} side packages
the same geometric $H^1$ into a $G_\Q$-representation
$H^1_{\acute{e}t}(\overline{X(\Gamma)},\Q_\ell)$, and the Hochschild--Serre spectral sequence explains
how Galois cohomology measures descent and obstructions relative to $\Q$.

\clearpage
\renewcommand{\refname}{References and Further Reading}

\begingroup
\small
\setlength{\columnsep}{3em}
\printindex
\endgroup

\end{document}